\documentclass[12pt]{amsart}
\usepackage{etex}
\usepackage{amscd}
\usepackage{amssymb}
\usepackage[noadjust]{cite}
\usepackage{booktabs}
\usepackage{url}
\usepackage{hyphenat}
\usepackage{mathtools}
\usepackage{cancel} % Allows cancelling terms
\usepackage{ifpdf}
\usepackage[all]{xypic}
\usepackage{mathrsfs}
\usepackage[T1]{fontenc}
\usepackage[utf8]{inputenc}
\usepackage{chngcntr}
\usepackage{enumitem}
\usepackage{xr-hyper}
\ifpdf
  \usepackage[pdftex,final]{graphicx}
  \usepackage[pdftex,lmargin=1in,rmargin=1in,tmargin=1in,bmargin=1in]{geometry}
  \usepackage[bookmarks=true, bookmarksopen=true,%
    bookmarksdepth=3,bookmarksopenlevel=2,%
    colorlinks=true,%
    linkcolor=blue,%
    citecolor=blue,%
    filecolor=blue,%
    menucolor=blue,%
    urlcolor=blue]{hyperref}
  \hypersetup{pdftitle={Bimodules in bordered Heegaard Floer homology}}
  \hypersetup{pdfauthor={Robert Lipshitz, Peter S. Ozsváth, and Dylan
      P. Thurston}}
\else
  \usepackage[dvips,final]{graphicx}
  \usepackage[dvips,lmargin=1in,rmargin=1in,tmargin=1in,bmargin=1in]{geometry}
  \usepackage[draft]{hyperref}
\fi
\usepackage{tikz}
\usetikzlibrary{matrix,arrows}

\tikzset{cdlabel/.style={above,sloped,
    execute at begin node=$\scriptstyle,execute at end node=$}}
\tikzset{algarrow/.style={->, thick}}
\tikzset{blgarrow/.style={->, thick}}
\tikzset{clgarrow/.style={->, thick}}
\tikzset{tensoralgarrow/.style={double, double equal sign distance, -implies}}
\tikzset{tensorblgarrow/.style={double, double equal sign distance, -implies}}
\tikzset{tensorclgarrow/.style={double, double equal sign distance, -implies}}
\tikzset{tensorelgarrow/.style={double, double equal sign distance, -implies}}
\tikzset{modarrow/.style={->, dashed}}
\tikzset{Amodar/.style={->, dashed}}
\tikzset{Dmodar/.style={->, dashed}}
\tikzset{DAmodar/.style={->, dashed}}

\newread\testin

\graphicspath{{draws/}{mpdraws/}}
\makeatletter
\def\input@path{{}{draws/}}
\makeatother

\def\mathcenter#1{%
  \vcenter{\hbox{$#1$}}%
}

\def\mfig#1{
        \mathcenter{\includegraphics{#1}}
}

\def\mfigb#1{
        \mathcenter{\includegraphics[trim=-1 -1 -1 -1]{#1}}
}

\makeatletter
\newcommand\mi@kern[1]{%
  \settowidth\@tempdima{$\mi@obj^{#1}$}
  \kern-\@tempdima
  #1
  \settowidth\@tempdima{$\mi@obj$}
  \kern\@tempdima
}

\newtoks\mi@toksp
\newtoks\mi@toksb
\DeclareRobustCommand{\manyindices}[5]{
  \def\mi@obj{#5}
  \mi@toksp\expandafter{\mi@kern{#2}}
  \mi@toksb\expandafter{\mi@kern{#1}}
  \@mathmeasure4\textstyle{#5_{#1}^{#2}}
  \@mathmeasure6\textstyle{#5_{#3}^{#4}}
  \dimen0-\wd6 \advance\dimen0\wd4
  \@mathmeasure8\textstyle{\hphantom{{}_{#1}^{#2}}#5^{\the\mi@toksp#4}_{\the\mi@toksb#3}}
  \hbox to \dimen0{}{\kern-\dimen0\box8}
}
\makeatother 

\usepackage{ifpdf}
\ifpdf
  \let\textalt\texorpdfstring
\else
  \newcommand{\textalt}[2]{#1}
\fi

\newcommand{\RR}{\mathbb R}

\newcommand{\ZZ}{\mathbb Z}
\newcommand{\QQ}{\mathbb Q}

\newcommand{\FF}{\mathbb F}

\newcommand{\bD}{\mathbb D}

\newcommand{\connectsum}{\mathbin{\#}}

\newcommand{\co}{\colon}
\newcommand{\semico}{;\penalty 300}

\newcommand{\abs}[1]{{\lvert #1 \rvert}}
\newcommand{\OneHalf}{{\textstyle\frac{1}{2}}}
\newcommand{\OneQuart}{{\textstyle\frac{1}{4}}}

\newcommand{\isom}{\simeq}

\newcommand{\tensor}{\otimes}
\newcommand{\bdy}{\partial}

\newcommand{\lbracket}{[}
\newcommand{\rbracket}{]}

\DeclareMathOperator{\PD}{PD}
\DeclareMathOperator{\Out}{Out}
\newcommand{\spinc}{\mathfrak s}
\newcommand{\relspinc}{\underline\spinc}

\DeclareMathOperator{\divis}{div}

\DeclareMathOperator{\End}{End}

\DeclareMathOperator{\Sym}{Sym}

\DeclareMathOperator{\Hom}{Hom}
\DeclareMathOperator{\Extop}{Ext}
\DeclareMathOperator{\Tor}{Tor}

\DeclareMathOperator{\spin}{spin}
\newcommand{\SpinC}{\spin^c}
\newcommand{\RelSpinC}{{\underline\spin}^c}

\DeclareMathOperator{\ind}{ind}

\DeclareMathOperator{\inv}{inv}
\DeclareMathOperator{\Inv}{Inv}
\DeclareMathOperator{\ev}{ev}
\DeclareMathOperator{\gr}{gr}

\DeclareMathOperator{\supp}{supp}
\DeclareMathOperator{\Int}{Int}
\newcommand{\emb}{{\mathrm{emb}}} % embedded

\newcommand\merid{m}
\newcommand\longitude{\ell}
\newcommand\Lmerid{\widetilde\merid}
\newcommand\Llongitude{\widetilde\longitude}

\theoremstyle{plain}
\newtheorem{theorem}{Theorem}
\numberwithin{equation}{section}
\newtheorem{proposition}[equation]{Proposition}
\newtheorem{lemma}[equation]{Lemma}
\newtheorem{corollary}[equation]{Corollary}

\newtheorem{definition}[equation]{Definition}

\theoremstyle{definition}

\newtheorem{construction}[equation]{Construction}
\newtheorem{convention}[equation]{Convention}

\theoremstyle{remark}
\newtheorem{example}[equation]{Example}
\newtheorem{remark}[equation]{Remark}
\newcommand{\HF}{\mathit{HF}}
\newcommand{\HFa}{\widehat {\HF}}

\newcommand{\CFa}{\widehat {\mathit{CF}}}
\newcommand{\HFKa}{\widehat {\mathit{HFK}}}
\newcommand{\x}{\mathbf x}
\newcommand{\y}{\mathbf y}
\newcommand{\z}{\mathbf z}
\newcommand{\w}{\mathbf w}

\newcommand\CH{\mathit{CH}}
\newcommand\tCH{\widetilde{\CH}}
\newcommand\HH{\mathit{HH}}
\newcommand\Hochschild\HH

\newcommand{\MCG}{\mathit{MCG}}

\newcommand{\Ainf}{\mathcal A_\infty}
\newcommand{\Ainfty}{\Ainf}

\newcommand{\SetI}{\mathbf i}
\newcommand{\SetS}{\mathbf s}
\newcommand{\SetT}{\mathbf t}
\newcommand{\SetSS}{\mathbf S}
\newcommand{\SetTT}{\mathbf T}

\newcommand{\Alg}{\mathcal{A}}
\newcommand{\AlgAS}[2]{\Alg(#1,#2)}
\newcommand{\AlgA}[1]{\Alg(#1)}
\newcommand{\Idem}{\mathcal{I}}
\newcommand{\IdemAS}[2]{\Idem(#1,#2)}
\newcommand{\IdemA}[1]{\Idem(#1)}
\newcommand\Gen{\mathfrak{S}}

\renewcommand{\S}{\Gen}
\newcommand{\SA}[1]{\S(#1)}

\newcommand{\Blg}{\mathcal{B}}
\newcommand{\Clg}{\mathcal{C}}
\newcommand{\Elg}{\mathcal{E}}
\newcommand{\HAlg}{\mathcal{H}}
\newcommand\Halg{\HAlg}

\newcommand{\alphas}{{\boldsymbol{\alpha}}}
\newcommand{\betas}{{\boldsymbol{\beta}}}

\newcommand{\rhos}{{\boldsymbol{\rho}}}

\newcommand{\cM}{\mathcal{M}}
\newcommand{\Mod}{\cM}

\newcommand{\CFD}{\mathit{CFD}}
\newcommand{\CFDD}{\mathit{CFDD}}
\newcommand{\CFA}{\mathit{CFA}}

\newcommand{\CFDA}{\mathit{CFDA}}
\newcommand{\CFDAa}{\widehat{\CFDA}}
\newcommand{\CFAA}{\mathit{CFAA}}
\newcommand{\CFAAa}{\widehat{\CFAA}}
\newcommand{\CFDa}{\widehat{\CFD}}

\newcommand{\CFK}{\mathit{CFK}}
\newcommand{\CFKa}{\widehat{\CFK}}

\newcommand{\CFDDa}{\widehat{\CFDD}}

\newcommand{\CFAa}{\widehat{\CFA}}

\newcommand{\s}{\mathfrak{s}}
\newcommand{\Source}{{S^{\mspace{1mu}\triangleright}}}

\newcommand{\cZ}{\mathcal{Z}}
\newcommand{\PtdMatchCirc}{\cZ}
\newcommand{\PMC}{\PtdMatchCirc}

\newcommand{\CircPts}{{\mathbf{a}}}

\newcommand\DGA{A}
\newcommand\DGB{B}

\newcommand{\dg}{\textit{dg} }
\newcommand{\DD}{\textit{DD}}
\newcommand{\DA}{\textit{DA}}
\newcommand{\AD}{\textit{AD}}
\newcommand{\AAm}{\textit{AA}} %conflicts with a Danish letter

\newcommand{\cDGA}{\mathcal{A}}

\newcommand\Id{\mathbb{I}}
\newcommand\Ground{\mathbf k}
\newcommand\Groundk{\mathbf k}
\newcommand\Groundl{\mathbf j}

\newcommand\DTP{\mathbin{\widetilde\otimes}}
\newcommand\DT{\boxtimes}

\newcommand\Tensor{T}

\newcommand{\Field}{\FF_2}
\DeclareMathOperator{\nbd}{nbd}

\newcommand{\Heegaard}{\mathcal{H}}
\newcommand{\HD}{\Heegaard}

\renewcommand{\th}{^\text{th}}

\newcommand{\bigGroup}{G'}
\newcommand{\smallGroup}{G}
\newcommand{\smallGrSet}{S}
\newcommand{\bigGrSet}{S'}

\newcommand{\sG}{\smallGroup}

\newcommand{\grb}{\gr'}

\newcommand{\Torus}{\mathbb{T}}
\DeclareMathOperator{\Fix}{Fix}
\DeclareMathOperator{\Rest}{Rest}
\DeclareMathOperator{\Induct}{Induct}
\newcommand{\tR}{\widetilde{R}}
\newcommand{\tP}{\widetilde{P}}

\newcommand{\Hyph}{\text{-}}

\newcommand{\ModCat}{\mathsf{Mod}}
\newcommand{\setModCat}{\widetilde{\ModCat}}
\newcommand{\DuModCat}[1]{\manyindices{\mathfrak u}{#1}{}{}\ModCat}
\newcommand{\DModCatD}[2]{\manyindices{}{#1}{}{#2}\ModCat}
\newcommand{\DuModCatDu}[2]{\manyindices{\mathfrak u}{#1}{\mathfrak u}{#2}\ModCat}

\newcommand{\DBModCatDB}[2]{\manyindices{\mathfrak B}{#1}{}{#2}\ModCat}
\newcommand{\DBModCatA}[2]{\manyindices{\mathfrak B}{#1}{#2}{}\ModCat}

\newcommand{\DbModCat}[1]{\manyindices{\mathfrak b}{#1}{}{}\ModCat}
\newcommand{\DModCat}[1]{\lsupv{#1}\ModCat}
\newcommand{\DuModCatA}[2]{\manyindices{\mathfrak u}{#1}{#2}{}\ModCat}
\newcommand{\DModCatA}[2]{\lsupv{#1}\ModCat_{#2}}
\newcommand{\kbModCat}{\ModCat^{\mathfrak{b}}}

\newcommand{\HMod}{\mathsf{H}}
\newcommand{\ZMod}{\mathsf{Z}}

\newcommand{\SM}{\mathit{SM}}
\newcommand{\Cat}{\mathscr{C}}
\newcommand{\Dat}{\mathscr{D}}
\newcommand{\CatEnd}{\mathscr{E}\!\textit{nd}}

\newcommand{\CatCat}{\mathscr{C}\!\textit{at}}
\DeclareMathOperator{\TMor}{2Mor}
\DeclareMathOperator{\Mor}{Mor}
\DeclareMathOperator{\Barop}{Bar}
\DeclareMathOperator{\Cobarop}{Cob}
\newcommand{\CCobarop}{\widehat{\Cobarop}}

\DeclareMathOperator{\ob}{Ob}

\newcommand{\sos}[3]{\mathbin{{}_{#1}\mathord#2_{#3}}}
\newcommand{\lsub}[2]{{}_{#1}\mathord{#2}}
\newcommand{\lsup}[2]{{}^{#1}\mskip-.6\thinmuskip\mathord{#2}}
\newcommand{\lsupv}[2]{{}^{#1}\mskip-.2\thinmuskip\mathord{#2}}
\newcommand{\lsubsup}[3]{\manyindices{#1}{\mskip.6\thinmuskip#2\mskip-.6\thinmuskip}{}{}{\mathord{#3}}}

\newcommand{\CatHonQi}{\mathcal{D}_{H,qi}}
\newcommand{\CatHonAinf}{\mathcal{D}_{H,\infty}}
\newcommand{\CatHonAinfQi}{\mathcal{D}_{H,\infty qi}}
\newcommand{\CatAinfAinf}{\mathcal{D}_{\infty,\infty}}
\newcommand{\CatAinfQi}{\mathcal{D}_{\infty,\infty qi}}

\newcommand\Lgen{\ell}
\newcommand\Rgen{r}
\newcommand\drHD{{\mathcal H}_{dr}}
\newcommand{\drtHD}{{\mathcal H}_{dr2}}
\newcommand{\drY}{Y_{dr}}
\newcommand{\dr}{{dr}}
\newcommand\fHD{{\mathcal H}_{f}}
\newcommand\sbY{\mathcal Y}

\newcommand\op{{\mathrm{op}}}
\newcommand\md{{\mathrm{mid}}}
\newcommand\PunctF{F^\circ}
\newcommand\Interior{\mathrm{int}}

\makeatletter
\newcommand\honestalg[3]{\bigl\lbracket
\begin{smallmatrix} #1\@ifempty{#3}{}{&#3} \\ #2 \end{smallmatrix}
\bigr\rbracket}
\makeatother
\newcommand{\lab}[1]{$\scriptstyle #1$}

\begin{document}
\title[Bimodules in bordered Heegaard Floer]{Bimodules in bordered
  Heegaard Floer homology}

\author[Lipshitz]{Robert Lipshitz}
\thanks{RL was supported by an NSF Mathematical Sciences
  Postdoctoral Fellowship, NSF Grant DMS-0905796, and a Sloan Research
  Fellowship.}
\address{Department of Mathematics, Columbia University\\
  New York, NY 10027}
\email{lipshitz@math.columbia.edu}

\author[Ozsv\'ath]{Peter S. Ozsv\'ath}
\thanks{PSO was supported by NSF Grant DMS-0505811.}
\address {Department of Mathematics, Princeton University\\ 
Princeton, NJ 08544}
\email {petero@math.princeton.edu}

\author[Thurston]{Dylan~P.~Thurston}
\thanks {DPT was supported by NSF
  Grant DMS-1008049 and a Sloan Research Fellowship.}
\address{Department of Mathematics,
         Indiana University\\
         Bloomington, IN 47405}
\email{dpthurst@indiana.edu}

\begin{abstract}
  Bordered Heegaard Floer homology is a three-manifold invariant
  which associates to a surface $F$ an algebra $\Alg(F)$ and to a
  three-manifold $Y$ with boundary identified with $F$ a module
  over $\Alg(F)$. In this paper, we
  establish naturality properties of this invariant.  Changing the
  diffeomorphism between $F$ and the boundary of~$Y$ tensors the
  bordered
  invariant with a suitable bimodule over $\Alg(F)$.  These bimodules
  give an action of a suitably based mapping class group on the
  category of modules over $\Alg(F)$. The
  Hochschild homology of such a bimodule is identified with the knot
  Floer homology of the associated open book decomposition. In the
  course of establishing these results, we also calculate the homology
  of $\Alg(F)$. We also prove a duality theorem relating the two
  versions of the $3$-manifold invariant. Finally, in the case of a
  genus one surface, we
  calculate the mapping class group action explicitly. This completes
  the description of bordered Heegaard Floer homology for knot
  complements in terms of the knot Floer homology.
\end{abstract}

\maketitle

\setcounter{tocdepth}{3}
\tableofcontents

\section{Introduction}
Bordered Heegaard Floer homology is an invariant associated to a
three-manifold with boundary~\cite{LOT1}, depending on some additional
data.  More specifically, let $F$ be a closed, oriented surface of
genus $k$.  A bordered three-manifold with boundary $F$ is a
compact, oriented three-manifold $Y$ equipped with an
orientation-preserving diffeomorphism $\phi\colon F \to
\partial Y$. Bordered Heegaard
Floer homology associates to $F$ (and some extra data, see below)
a differential graded (\textit{dg}) algebra
$\mathcal{A}(F)$.  If $Y$ is a bordered three-manifold with boundary
$F$, the theory associates to $Y$ a right $\Ainf$-module $\CFAa(Y)$ over
$\Alg(F)$, the {\em type $A$ module of $Y$},  whose
quasi-isomorphism type depends only on the underlying diffeomorphism
type of~$Y$.  
To a
three-manifold with boundary the
theory also associates a  left
$\Ainf$-module over $\mathcal{A}(-F)$, $\CFDa(Y)$, the {\em type $D$ module
  of $Y$}.

Bordered Heegaard Floer homology is related to Heegaard Floer homology
$\HFa(Y)$ via a pairing theorem: if $Y$ is a three-manifold
which is divided into $Y_1$ and $Y_2$ by a separating surface $F$, 
then $\CFa(Y)$, a chain complex whose homology calculates $\HFa(Y)$,
is obtained as the $\Ainf$-tensor product of $\CFAa(Y_1)$ with
$\CFDa(Y_2)$. In other words, $\HFa(Y)\cong
\Tor_{\Alg(F)}(\CFAa(Y_1),\CFDa(Y_2))$. 

\subsection{Reparametrization and the bordered Floer invariants}

A key goal of this paper is to study how the bordered Heegaard Floer
invariants changes under reparametrization of the boundary.
More precisely, we fix a closed surface $F$, and also a preferred disk
$D\subset F$, together with a point $z\in \partial D$. Consider the
space of diffeomorphisms of $F$ which preserve the disk $D$ and
the point $z\in\partial D$.  This topological group will be called the
\emph{strongly based diffeomorphism group of $(F,D,z)$}.  Its group of
path components is called the \emph{strongly based
  mapping class group} of $F$, and two diffeomorphisms in the same
path component are called \emph{strongly isotopic}. This agrees with the usual mapping
class group of $F\setminus\Interior(D)$ (fixing the
boundary).

Consider next a handle decomposition of $F$ with one zero-handle
and where $D$ is the unique two-handle. We will mark in addition a basepoint
$z$ on the boundary of $F\setminus D$. This data can be combinatorially
encoded in the form of a {\em pointed matched circle} (see Definition~\ref{def:PMC} below). Bordered Floer homology
associates to a pointed matched circle~$\PMC$
a differential-graded
algebra~$\AlgA{\PtdMatchCirc}$.
Modules over these algebras are independent of the decomposition
$\PMC$ in the
following sense:
\begin{theorem}
  \label{thm:AlgebraDependsOnSurface}
  If $\PMC_1$ and $\PMC_2$ are two pointed matched
  circles representing the same underlying surface, then
  the derived categories of \dg $\AlgA{\PMC_1}$- and
  $\AlgA{\PMC_2}$-modules are equivalent.
\end{theorem}

For the purpose of this introduction, we will typically suppress the
pointed matched circle $\PMC$ from the notation, referring somewhat imprecisely
to $\AlgA{F}$. Theorem~\ref{thm:AlgebraDependsOnSurface} provides some justification for this practice.

A \emph{bordered three-manifold} is a quadruple
$(Y,\Delta,z_1,\psi)$, where $Y$ is an oriented three-manifold-with-boundary,
$\Delta$ is a disk in $\partial Y$, $z_1$ is a point on $\partial \Delta$,
and
$$\psi\colon (F,D,z) \to (\partial Y,\Delta,z)$$
is a
diffeomorphism from $F$ to $\partial Y$ sending $D$ to $\Delta$ and $z$ to~$z_1$.

The strongly based diffeomorphism group of $F$ acts on the set of
bordered three-manifolds by composition:
$$\phi\cdot(Y,\Delta,z_1,\psi) =(Y,\Delta,z_1,\psi\circ \phi^{-1}).$$
There are bimodules which encode this action, using the $\Ainf$-tensor product.
Specifically, let $M$ be a right $\Ainf$-module over the \dg algebra
$\Alg$ and $N$ be an $\Ainf$-bimodule over $\Alg$ and $\Blg$, where
$\Blg$ is another \dg algebra. Then we
can form the derived (or $\Ainf$) tensor product $M\DTP_{\Alg} N$, to obtain a
right $\Ainf$-module over $\Blg$.
We have bimodules associated to reparameterizing the boundary, as given
in the following:

\begin{theorem} 
  \label{thm:Reparameterization}
  Given a strongly based diffeomorphism $\phi\colon (F_1, D, z)
  \to (F_2, D, z)$ between surfaces $F_1$ and $F_2$ (corresponding to possibly
  different pointed matched circles), there are associated bimodules:
\[\CFAAa(\phi)_{\AlgA{-F_1},\AlgA{F_2}},\qquad
  \lsub{\AlgA{F_1}}\CFDAa(\phi)_{\AlgA{F_2}}, \qquad
  \lsub{\AlgA{F_1},\AlgA{-F_2}}\CFDDa(\phi).\footnote{Given algebras
    $A$ and $B$, we write ${}_AM_B$ to denote a module with a left
    action of $A$ and a right action of $B$. We write ${}_{A,B}M$ to
    denote a module with commuting left actions of $A$ and $B$, and
    $M_{A,B}$ to denote a module with commuting right actions of $A$
    and $B$.}\] 
  If $(Y_1, \Delta_1, z_1, \psi_1\co F_1\to \bdy Y_1)$
  and $(Y_2, \Delta_2, z_2,(-\psi_2)\co -{F_2}\to\bdy Y_2)$ (so
  $\psi_2\co F_2\to -\bdy Y_2$) are bordered 3-manifolds then:
  \begin{align*}
    \CFAa(Y_1,\psi_1)\DTP_{\AlgA{F_1}} \CFDAa(\phi) & \simeq \CFAa(Y_1,\psi_1\circ\phi^{-1})\\
    \CFAAa(\phi)\DTP_{\AlgA{F_2}}\CFDa(Y_2,\psi_2) & \simeq
    \CFAa(Y_2,-(\psi_2\circ\phi)) \\
    \CFAa(Y_1,\psi_1)\DTP_{\AlgA{F_1}} \CFDDa(\phi) & \simeq \CFDa(Y_1,-(\psi_1\circ\phi^{-1})) \\
    \CFDAa(\phi)\DTP_{\AlgA{F_2}}\CFDa(Y_2,\psi_2) & \simeq \CFDa(Y_2,\psi_2\circ\phi).
  \end{align*}
\end{theorem}
(This is proved in Section~\ref{sec:connected-pairing}. See
particularly Figure~\ref{fig:ModBimodPairingSchematic} for a schematic
illustrating why the parametrizations are as given.)

The bimodules satisfy the following invariance property:

\begin{theorem} 
  \label{thm:IsotopyInvariance}
  If $\phi$ and $\phi'$ are strongly isotopic diffeomorphisms of $F$
  then their associated bimodules are quasi-isomorphic.
\end{theorem}
(This is proved in Section~\ref{sec:AutBimodules}.)

The bimodules also behave functorially under composition, according to the following two results:

\begin{theorem}
  \label{thm:Id-is-Id}
  The type \DA\ bimodule associated to identity map from $F$ to itself,
  $\CFDAa(\Id_F)$, is quasi-isomorphic to $\AlgA{F}$ as an
$(\AlgA{F},\AlgA{F})$-bimodule.
\end{theorem}
(This is proved in Section~\ref{sec:id-bim}.) Note that $\AlgA{F}$ is
the identity for the tensor product operation.

\begin{theorem}
  \label{thm:Composition}
   Given two strongly based diffeomorphisms $\phi_1\co F_1\to F_2$ and
   $\phi_2\co F_2\to F_3$, we have that 
    $$\CFDAa(\phi_1)\DTP_{\AlgA{F_2}}\CFDAa(\phi_2)\simeq\CFDAa(\phi_2\circ\phi_1).$$
\end{theorem}
(This is proved in Section~\ref{sec:connected-pairing}.)

Together, Theorems~\ref{thm:IsotopyInvariance}, \ref{thm:Id-is-Id},
and~\ref{thm:Composition} can be summarized by saying that the
bimodules induce an action of the based mapping class group of a
surface $F$ on the module category of $\AlgA{F}$; see
Theorem~\ref{thm:MCG-2-functoriality} in Section~\ref{sec:mcg} for a
more precise statement. Actions of mapping class
groups---particularly, braid groups---on categories have arisen in
other contexts; see for example~\cite{KhTomas07:CobordismsCategories}
and the references therein.

Theorem~\ref{thm:Reparameterization} is also interesting in the case
where $\phi$ is the identity map. In that case, the theorem allows one
to convert a type $D$ module into a type $A$ module, and vice versa.
The statement is given in the following corollary, which can be
thought of as exhibiting a kind of Koszul
duality~\cite{Priddy70:Koszul} between the algebra of a
surface and its orientation reverse:

\begin{corollary}
  \label{cor:ConvertDintoA}
  There are dualizing modules $\CFAAa(\Id)$ and $\CFDDa(\Id)$ (bimodules
  associated to the identity map on $F$), which can be used to convert
  $\CFDa$-modules to $\CFAa$-modules and vice versa, in the sense that
  \begin{align*}
    \CFAAa(\Id)\DTP_{\AlgA{-F}}\CFDa(Y,-\psi)&\simeq \CFAa(Y,\psi) \\
    \CFAa(Y,\psi)\DTP_{\AlgA{F}}\CFDDa(\Id)&\simeq \CFDa(Y,-\psi).
  \end{align*}
\end{corollary}

The modules $\CFAAa(\Id)$ and $\CFDDa(\Id)$ uniquely determine one
another, as explained in Section~\ref{sec:Duality}. Indeed, this
description, along with the above corollary, quickly leads to the 
following result which allows one to express type $A$ modules entirely in 
terms of type $D$ modules:

\begin{theorem}
  \label{thm:Duality}
  Let $Y$ be a bordered three-manifold with boundary $F$.  Then
  $\CFAa(Y)$ is quasi-isomorphic, as a right $\Ainf$ module over
  $\AlgA{F}$, to the chain complex of maps from $\CFDDa(\Id)$
  to $\CFDa(Y)$.
\end{theorem}
(This theorem is stated more precisely and proved in Section~\ref{sec:Duality}, as
Theorem~\ref{thm:Duality-precise}.)

In fact, $\CFDDa(\Id)$, and hence also $\CFAAa(\Id)$, can be
calculated explicitly; see~\cite{LOT4}.

The bimodule ${}_{\AlgA{F}}\CFDAa(\phi)_{\AlgA{F}}$ of a strongly
based diffeomorphism 
$$\phi\colon(F,D,z)\to (F,D,z)$$ is a
bimodule over a single algebra, $\AlgA{F}$, on both sides; and hence,
it is natural to
study its Hochschild homology. We give this operation a
topological interpretation. To state this interpretation, recall that a strongly based
diffeomorphism $\phi\colon(F,D,z)\to (F,D,z)$ naturally
gives rise to a three-manifold with an open book decomposition. More
precisely, consider the three-manifold with torus boundary, defined as
the quotient of $[0,1]\times (F\setminus D)$ by the equivalence relation
$(0,\phi(x))\sim (1,x)$. This manifold is equipped with an embedded,
closed curve on the boundary, $([0,1]\times\{z\})/(0,z)\sim (1,z)$. By
filling along the curve on the boundary, we get a closed
three-manifold, equipped with a knot $K$ induced from $\{0\}\times
\partial D$. We denote the resulting three-manifold by $Y(\phi)$, and
let $K\subset Y(\phi)$ be the canonical knot in it.  This presentation
of the three-manifold $Y$ underlying $Y(\phi)$ is called an {\em open
  book decomposition} of $Y$, and $K$ is called the {\em binding}.

\begin{theorem}
  \label{thm:Hochschild}
  Let $\phi$ be a strongly-based diffeomorphism, and consider its
  associated $\AlgA{F}$-bimodule $\CFDAa(\phi)$. The Hochschild homology of the
  bimodule $\CFDAa(\phi)$ is
  the knot Floer
  homology of $Y(\phi)$ with respect to its binding $K$,
  $\HFKa(Y(\phi),K)$.
\end{theorem}
(This is proved in Section~\ref{sec:self-pairing}.)

\begin{remark}
  Let $K$ be a fibered knot in $S^3$, $F$ a fiber surface for $K$ and
  $\phi\co F\to F$ the monodromy. A classical theorem states that the
  Alexander polynomial of $K$ is the characteristic polynomial of
  $\phi_*\co H_1(F)\to H_1(F)$. Theorem~\ref{thm:Hochschild}
  categorifies this in the following sense. On the one hand,
  by~\cite[Equation (1)]{OS04:Knots}, the Euler characteristic of
  $\HFKa(Y(\phi),K)$ is $\Delta_K(t)$.  On the other hand, the
  (modulo~$2$) Grothendieck group of the category of $\Alg(F)$-modules
  is
  $\Lambda^*H_1(F;\Field)$ and the functor $\CFDAa(\phi)\DTP\cdot$
  decategorifies to $\phi_*\co \Lambda^*H_1(F;\Field)\to
  \Lambda^*H_1(F;\Field)$; compare~\cite{LOT13:faith}. Hochschild
  homology decategorifies to the graded trace; and the graded trace of
  $\phi_*\co \Lambda^*H_1(F;\Field)\to \Lambda^*H_1(F;\Field)$ is (by
  definition) the characteristic polynomial of $\phi_*\co
  H_1(F;\Field)\to H_1(F;\Field)$. (To obtain the decategorification
  with $\ZZ$-coefficients would require an absolute $\ZZ/2$-grading on
  the bimodules, which we do not construct.)
\end{remark}

\subsection{A more general framework}

The bimodules associated to mapping classes come from a more general
construction, which gives an invariant for a bordered three-manifold
$Y_{12}$ with two boundary components $F_1$ and $F_2$, and extra data
specified by a {\em strong framing}, which is a parameterization of the
boundary components of $Y_{12}$ and a framed arc connecting
the two boundary components of $Y_{12}$.  More precisely, we have the
following.

\begin{definition}
\label{def:IntroStronglyBordered}
Fix two connected surfaces $(F_1,D_1,z_1)$ and $(F_2,D_2,z_2)$ equipped with
preferred disks and basepoints on the boundaries of the
disks. A {\em strongly bordered three-manifold with boundary
  $F_1\amalg F_2$} is an oriented three manifold $Y_{12}$ with two
boundary components,
equipped with
\begin{itemize}
  \item preferred disks $\Delta_1$ and $\Delta_2$ on its two boundary
    components,
  \item basepoints $z_i'\in \partial \Delta_i$,
  \item diffeomorphisms $\psi\colon 
    (F_1\amalg F_2,D_1\amalg D_2,z_1\amalg z_2)
    \to (\partial Y_{12},\Delta_1\amalg \Delta_2,z_1'\amalg z_2')$,
  \item an arc $\gamma$ connecting $z_1'$ to $z_2'$, and
  \item a framing of $\gamma$, pointing into $\Delta_i$ at 
    $z_i'$ for $i=1,2$.
\end{itemize}
\end{definition}
(See also Definition~\ref{def:StronglyBordered}.)
Note that we can remove a neighborhood of $\gamma$ from $Y_{12}$ to
obtain a three-manifold $M$ with boundary $F_1\connectsum F_2$. The
trivialization of the normal bundle of $\gamma$ is the additional data needed to
construct a parameterization of the boundary of $M$ from the
parameterization of $\bdy Y$.

To obtain a bimodule, we must fix further data: we mark
each boundary component of $Y$ with an $A$ or a $D$. This
determines whether the corresponding boundary component is treated
as a type $D$ module or a type $A$ module (and hence, whether the underlying
algebra is associated to the boundary component with its induced orientation,
or the opposite of its induced orientation). More explicitly, we have the following:

\begin{theorem}
  \label{thm:InvarianceOfBimodules}
  Let $Y_{12}$ be a bordered three-manifold with two boundary
  components $F_1$ and $F_2$ and a strong framing.  Then, we can
  associate the following bimodules to $Y_{12}$:
 \[\CFAAa(Y_{12})_{\AlgA{F_1},\AlgA{F_2}}, \qquad
  {}_{\AlgA{-F_1}}\CFDAa(Y_{12})_{\AlgA{F_2}},\qquad
  {}_{\AlgA{-F_1},\AlgA{-F_2}}\CFDDa(Y_{12}).\]
  The
  quasi-isomorphism types of these bimodules are diffeomorphism
  invariants of the bordered three-manifold $Y_{12}$ with its strong
  boundary framing.
\end{theorem}
(A graded refinement of this theorem is proved in
Section~\ref{sec:CFBimodules}, as
Theorem~\ref{thm:gradedInvarianceOfBimodules}.)

The bimodules associated to an automorphism of $F$ are gotten as a
special case of the above construction, where $Y_{12}=[0,1]\times F$,
with the identity parametrization on one component, and $\phi$ on the other
component. In Section~\ref{sec:PairingTheorems}, we shall prove
general versions of the pairing theorem, from which
Theorems~\ref{thm:Reparameterization} and~\ref{thm:Composition} follow
as a corollaries.

There are other versions of the pairing theorem, including a
``double-pairing theorem'', where we glue two three-manifolds with two
boundary components together, Theorem~\ref{thm:DoublePairing}. Theorem~\ref{thm:Hochschild} is the
specialization of this to the case where we are self-gluing $[0,1]\times F$.

After building the basic background, we turn to the particular case of
a surface of genus one.  In this case, we calculate the bimodules
associated to an arbitrary mapping class, hence giving an explicit
description of the dependence of the bordered Heegaard Floer homology
of a three-manifold with torus boundary on the parameterization of the
boundary.  This also completes the description of $\CFDa$ in terms
of knot Floer homology for a knot in $S^3$.  Specifically, we showed
in~\cite[Chapter~\ref*{LOT:chap:TorusBoundary}]{LOT1} how to calculate
the type $D$ module of a knot
complement in terms of the knot Floer homology of the knot when the
framing on the knot is sufficiently large.  With the help of the torus
mapping class group calculations and
Theorem~\ref{thm:Reparameterization}, we are now able to calculate the
type $D$ module for arbitrary framings, as promised in~\cite{LOT1}.

\subsection{Organization}

This paper is organized as follows.  In
Section~\ref{sec:algebra-modules} we recall the algebraic background
on $\Ainf$-algebras and modules used throughout this paper. These include the
familiar notions of $\Ainf$ modules, which we call here {\em type $A$
  structures} to distinguish them from {\em type $D$ structures},
which are a variant of projective modules (see Corollary~\ref{cor:tensor-projective};
compare also
Remarks~\ref{rmk:Comodule} and~\ref{rmk:TwistedComplex}).
These
are then combined to
give various notions of bimodules. We discuss various operations on
modules and bimodules---in particular, the tensor product, $\Hom$ and
Hochschild homology functors---and review some category
theory. Section~\ref{sec:algebra-modules} concludes with a discussion
of group-valued gradings.

In Section~\ref{sec:PointedMatchedCircles} we recall the notion of a {\em
  pointed matched circle} $\PtdMatchCirc$, which is effectively the
handle decomposition of a surface used to define its associated
algebra. We then recall the definition of the algebra $\AlgA{\PMC}$
(introduced in~\cite{LOT1}). In Section~\ref{sec:homology-algebra}, we
turn to the calculation of the homology of the algebra
$\AlgA{\PMC}$. This calculation is used significantly in the
subsequent proof of Theorem~\ref{thm:Id-is-Id} in
Section~\ref{sec:id-bim}.  It is interesting to note that, as a consequence of
this calculation, we obtain a smaller differential graded algebra
${\mathcal A}'$ (a quotient of ${\mathcal A}$ by a
differential ideal), which can be used in place of ${\mathcal A}$ for
the purposes of the invariant; see
Proposition~\ref{prop:SmallerModel}.  
(Although we do not pursue this further in the present work, it does considerably simplify computations in practice.)

In Section~\ref{sec:Diagrams}, we study Heegaard diagrams associated to strongly framed
three-manifolds, and recall their existence and uniqueness
properties. We also turn our attention to a case of particular
importance: strongly framed three-manifolds associated to strongly
based surface automorphisms. We show how to construct explicitly the
corresponding Heegaard diagrams in this case.  In
Section~\ref{sec:CFBimodules}, we turn to the construction of bordered
bimodules for strongly based three-manifolds, and verify their
invariance properties, verifying
Theorems~\ref{thm:InvarianceOfBimodules}
and~\ref{thm:IsotopyInvariance}.

In Section~\ref{sec:PairingTheorems}, we turn to the various pairing
theorems which these bimodules enjoy.
Theorems~\ref{thm:Reparameterization}, \ref{thm:Composition}
and~\ref{thm:Hochschild} are deduced from these pairing theorems.  In
Section~\ref{sec:mcg}, we employ the calculations from
Section~\ref{sec:homology-algebra}, together with a suitable variant
of the pairing theorem, to deduce that the type \DA\ bimodule
associated to the identity map is $\Alg(F)$,
Theorem~\ref{thm:Id-is-Id}.  Thus armed,
we complete the proof that the
derived category of modules over the algebra
associated to a pointed matched circle depends only on the
homeomorphism type of the underlying surface,
Theorem~\ref{thm:AlgebraDependsOnSurface}. This information also
allows us to construct an action of the
mapping class group on the derived category of
$\Alg(F)$-modules.
Theorem~\ref{thm:Id-is-Id} is also the main ingredient we use in
Section~\ref{sec:Duality} to prove the duality theorem, Theorem~\ref{thm:Duality}.

Finally, in
Section~\ref{sec:torus-calc} we compute the bimodules associated to
torus automorphisms.

\subsection{Acknowledgements}
We are grateful for the supportive and stimulating mathematical
environment provided by the Columbia mathematics department, where
most of this work was done.  We also wish to thank MSRI for its hospitality
during the completion of this project. In addition, the third author
thanks the Center for the Topology and Quantization of Moduli Spaces
at Aarhus University for its hospitality.

The authors thank M.~Khovanov and
A.~Lauda for substantial help with the categorical aspects of this
paper, particularly with formulating
Theorem~\ref{thm:MCG-2-functoriality}.  The authors would also like to
thank R.~Zarev for sharing with us his insights during the course of
this work.  We thank C.~Manolescu
for a correction to a previous version. Finally, we thank the referee
for many helpful comments.

%%% Local Variables: 
%%% mode: latex
%%% TeX-master: "Bimodules"
%%% End: 

\counterwithin{equation}{subsection}
\section{\textalt{$\Ainf$}{A-infty} algebras and modules}
\label{sec:algebra-modules}
In this section we recall various notions from the theory of
$\Ainf$-algebras and derived categories. Most of these results are
standard (see
\cite{AinftyAlg,BernsteinLunts94:EquivariantSheaves,LefevreAInfinity,SeidelBook}),
and are collected here for the reader's convenience.  Our treatment is
slightly nonstandard in that we use extensively a certain algebraic
object defined  over $\Ainf$-algebras, which we call ``type $D$
structures.'' The reader is encouraged to think of these as
projective modules over the algebra. 
These arise naturally in the context of bordered Floer
theory---the bordered invariant $\CFDa(Y)$ of~\cite{LOT1} is a type
$D$ structure---and are convenient for various algebraic
constructions. (In fact, type $D$ structures have appeared under
various guises elsewhere, see Remarks~\ref{rmk:Comodule}
and~\ref{rmk:TwistedComplex}.)

\begin{convention}\label{convention:GroundandFD}
  Throughout, ($\Ainf$-) algebras will be algebras over the a ring
  $\Ground$, which is either $\FF_2$ or, more generally,
  $\oplus_{i=1}^N\FF_2$. Unless otherwise stated, tensor products
  denote tensor products over $\Ground$. (When we need to refer to
  a second such ground ring we will denote it $\Groundl$.)
\end{convention}

Note that for most of this paper it is enough to consider \dg algebras
rather than more general $\Ainf$-algebras, so the reader could skip
Section~\ref{sec:A-infty-algs} if desired. More general
$\Ainf$-algebras do appear in
Section~\ref{sec:short-chords-massey-gen}.

In one unusual twist, our $\Ainf$-algebras and modules are graded by
non-commutative groups, as explained in
Section~\ref{sec:algebras-gradings}.

\subsection{\textalt{$\Ainf$}{A-infty}-algebras}
\label{sec:A-infty-algs}
\subsubsection{Definition of \textalt{$\Ainf$}{A-infty}-algebras}
\label{sec:defin-Ainfty-alg}
\begin{definition}\label{def:Ainf-algebra}
  An $\Ainf$-algebra $\cDGA$ over~$\Ground$
  is a $\ZZ$-graded $\Ground$-bimodule~$A$, equipped with degree~$0$
  $\Ground$-linear multiplication maps
  \begin{equation}
    \label{eq:a-inf-ops}
    \mu_i\co \DGA^{\otimes i} \to \DGA[2-i]
  \end{equation}
  defined for all $i\geq 1$, satisfying the compatibility conditions
  that, for each~$n \ge 1$ and elements $a_1,\dots,a_n$,
  \begin{equation}\label{eq:a-inf-alg-rel}
    \sum_{i+j=n+1}\sum_{\ell=1}^{n-j+1}\mu_i(a_1\otimes \dots\otimes a_{\ell-1}\otimes \mu_j(a_\ell\otimes \dots\otimes a_{\ell+j-1})\otimes a_{\ell+j}\otimes \dots\otimes a_{n}) 
        =0.
  \end{equation}
  Here, $\DGA^{\otimes i}$ denotes the $\Ground$-bimodule
  $\overbrace{\DGA\otimes_{\Ground}\dots\otimes_{\Ground}\DGA}^{i}$
  and $[2-i]$ denotes a degree shift.
  We use $\cDGA$ for the $\Ainf$ algebra and $\DGA$ for its
  underlying $\Ground$-bimodule.

  An $\Ainf$ algebra is \emph{strictly unital} (or just \emph{unital})
  if there is an element
  $1\in A$ with the property that $\mu_2(a,1)=\mu_2(1,a)=a$ and
  $\mu_i(a_1,\dots,a_i)=0$ if $i\ne 2$ and $a_j=1$ for some $j$.
   For a unital $\Ainf$-algebra, the unit gives a preferred map
  $\iota\co \Ground \to \cDGA$.  

  An \emph{augmentation} of an $\Ainf$-algebra is a map $\epsilon\co
  \DGA\to\Ground$, satisfying the conditions that
  \begin{equation}\label{eq:augmentation}
  \begin{aligned}
    \epsilon(1) &= 1\\
    \epsilon(\mu_2(a_1,a_2)) &= \epsilon(a_1)\epsilon(a_2)\\
    \epsilon(\mu_k(a_1,\ldots,a_k)) &= 0&\text{for $k \ne 2$}.
  \end{aligned}
  \end{equation}
  This gives an \emph{augmentation ideal}
  $\DGA_+=\ker \epsilon$.
\end{definition}

One could consider a more general notion of augmentation, where
$\epsilon$ is an $\Ainf$ homomorphism in the sense of
Section~\ref{sec:defin-Ainf-alg-maps}.  We do not do this here, as we
do not need this level of generality for our present purposes, and
indeed, it would cause undue complication, especially in
Section~\ref{sec:hochschild-homology}.

Note that in particular $\mu_1$ gives $\DGA$ the structure of a chain
complex.
A differential graded algebra over $\Ground$ is an $\Ainf$ algebra
with $\mu_i=0$ for all $i>2$.

\begin{convention}\label{Convention:AlgsUniConFd}
  Throughout this paper, $\Ainf$-algebras will be assumed strictly
  unital and augmented.
\end{convention}

We can think about the $\Ainf$-relation
(Equation~(\ref{eq:a-inf-alg-rel})) graphically.  First, we associate
operations to graphs. 
\begin{definition} An \emph{$\Ainf$-operation tree} $\Gamma$ is a
  finite, directed tree embedded in the plane 
  such that each non-leaf vertex of $\Gamma$ has exactly one outgoing
  edge.
\end{definition}
In every $\Ainf$-operation tree $\Gamma$ there is a unique leaf that is a sink, along with
$n$ source leaves.  Then, given an $\Ainf$-algebra $\Alg$ we can
associate to $\Gamma$ an operation
\[
\mu_\Gamma\co A^{\otimes n}\to A
\]
(with some grading shift) as follows. To compute
$\mu_\Gamma(a_1\otimes\cdots\otimes a_n)$,
start with $a_1,\dots,a_n$ at the source leaves (labelled in order,
clockwise around the boundary of~$\Gamma$). Flow these
elements along $\Gamma$, and
when a string of elements enter a vertex of valence $k>1$, apply~$\mu_k$.
The element at the sink is the output.

In these terms, the basic $\Ainf$ relation states that the sum of
$\mu_\Gamma$ over $\Ainf$-operation trees~$\Gamma$ with exactly
two non-leaf vertices is zero.

Given an $\Ainf$-algebra $\cDGA=(\DGA,\{\mu_i\})$ we can form the tensor algebra
\[
\Tensor^*(\DGA[1])\coloneqq\bigoplus_{n=0}^{\infty}A^{\otimes n}[n].
\]
Setting $\mu_0=0,$ we can combine all the $\mu_i$ to form a single map
\begin{equation}
\mu\co \Tensor^*(\DGA[1])\to \DGA[2].\label{eq:def-mu}
\end{equation}
Defining
${\overline D}^{\Alg}\co \Tensor^*(\DGA[1])\to \Tensor^*(\DGA[1])$ by
\[{\overline D}^{\Alg}(a_1\otimes\dots\otimes a_n) = \sum_{j=1}^n
\sum_{\ell=1}^{n-j+1} a_1\otimes\dots \otimes
\mu_j(a_\ell\otimes\dots\otimes a_{\ell+j-1})\otimes\dots\otimes
a_n,\] the $\Ainf$ compatibility relations are encoded in the
relation $\mu\circ {\overline D}^{\Alg}=0$ or, equivalently,
${\overline D}^{\Alg}\circ{\overline D}^{\Alg}=0$.

Our $\Ainf$-algebras also need to be appropriately bounded.  Before
giving the general
case, we start with the version for \dg algebras:
\begin{definition}\label{def:Alg-bounded-dga}
  We say that an augmented \dg algebra $\Alg$ has \emph{nilpotent
    augmentation ideal} if there exists
  an $n$ so that $(A_+)^n$ = 0.
  We will also abuse terminology and say
  that $\Alg$ itself is \emph{nilpotent}.
\end{definition}
Note that this is stronger than saying that every element of $A_+$ is
nilpotent.  Also, a unital algebra can never be
nilpotent in the strict sense.

\begin{definition}\label{def:Alg-bounded}
  An augmented $\Ainf$-algebra $\cDGA=(\DGA,\{\mu_i\})$ is called
  \emph{nilpotent}
  if there exists an $n$ so that for any $i>n$, any elements $a_1,\dots,a_i\in
  A_+$  and any $\Ainf$-operation tree~$\Gamma$,
  $\mu_\Gamma(a_1\otimes\cdots\otimes a_i)=0$.
\end{definition}
(An equivalent definition of nilpotent would be to require that there
are only finitely many $\Gamma$ for which $\mu_\Gamma$ is not
identically zero on inputs from $A_+$. The definitions of
``operationally bounded'' for modules and bimodules, below, can be
reformulated similarly. Also, the fact that the two definitions of
nilpotent agree in the case of \dg algebras involves using the Leibniz
rule to push any instances of $\mu_1$ onto the inputs and then using
$(\mu_1)^2 = 0$.)
\begin{remark}\label{rem:ainf-bounded}
  In~\cite{LOT1} we assumed a weaker condition on our algebra: that
  $\mu_i=0$ for $i$ sufficiently large.  The stronger condition of
  Definition~\ref{def:Alg-bounded} is used to ensure that our smaller
  model $\DT$ of the $\Ainf$-tensor product of bimodules (rather than
  just modules) is
  well-defined in all cases (Proposition~\ref{prop:DT-many-cases}) and
  for some of the
  categorical aspects of this paper (in particular,
  Proposition~\ref{prop:D-to-A-and-back}). See also
  Remark~\ref{rmk:weaker-stronger-boundedness}.
\end{remark}
\begin{convention}
  With the exception of Section~\ref{sec:induced-Ainf-alg}, all of the
  $\Ainf$-algebras that show up in this paper will be assumed (or
  proved) to be nilpotent.
\end{convention}
\subsubsection{Definition of \textalt{$\Ainf$}{A-infty}-algebra maps}
\label{sec:defin-Ainf-alg-maps}
\begin{definition}\label{def:Ainf-alg-morph} Let $\Alg$ and $\Blg$ be
  $\Ainf$-algebras. An \emph{$\Ainf$-homomorphism from $\Alg$ to
    $\Blg$} is a collection of degree~0 maps
  $\phi=\{\phi_i\co\Alg[1]^{\otimes i}\to\Blg[1]\}$, $i\geq 1$,
  satisfying a compatibility
  condition which we state in terms of an auxiliary map
  $F^\phi=\Tensor^*(\Alg[1])\to\Tensor^*(\Blg[1])$, defined by
  \begin{multline*}
  F^\phi(a_1\otimes\dots\otimes
  a_n)=
  \sum_{i_1+\dots+i_k=n}\phi_{i_1}(a_1\otimes\dots\otimes
  a_{i_1})\otimes\phi_{i_2}(a_{i_1+1}\otimes\dots\otimes
  a_{i_1+i_2})\otimes\cdots \\
   \cdots\otimes\phi_{i_k}(a_{n-i_k+1}\otimes\dots\otimes
  a_n).
  \end{multline*}
  The compatibility condition is that
  \begin{equation}
    \label{eq:Ainf-alg-morph}
    \overline{D}^\Blg\circ F^\phi=F^\phi\circ\overline{D}^\Alg.
  \end{equation}
\end{definition}
Note that if $\phi=\{\phi_i\}$ is an $\Ainf$-homomorphism then
$\phi_1$ is a chain map.

Composition of $\Ainf$-homomorphisms is characterized by the property that
$F^{\phi\circ\psi}=F^\phi\circ F^\psi$; we leave it to the reader to
verify that such a composition exists.

It turns out that $\Ainf$-algebra \emph{isomorphisms} are just
homomorphisms with $\phi_1$ invertible:
\begin{lemma}\label{lemma:1isos-are-isos}
  Let $\Alg$ be an $\Ainf$-algebra and $\phi\co\Alg\to\Blg$ an
  $\Ainf$-algebra homomorphism such that $\phi_1$ is an
  isomorphism. Then $\phi$ is
  invertible, i.e., there is an $\Ainf$-algebra homomorphism $\psi\co
  \Blg\to\Alg$ such that $\phi\circ\psi=\Id_\Blg$ and $\psi\circ\phi=\Id_\Alg$.
\end{lemma}
\begin{proof}
  It suffices to show that $\phi$ has both left and right inverses; we
  will show that $\phi$ has a left inverse.
  Set $\psi_1=\phi_1^{-1}$. Observe that $\psi$ satisfies the first
  relation for $\Ainf$-homomorphisms ($\mu_1^\Alg \circ \psi_1+\psi_1\circ \mu_1^\Blg=0$); moreover, for
  any way of completing $\psi$ to
  an $\Ainf$-homomorphism, $(\psi\circ\phi)_1=\Id_\Alg$.

  Now, assume inductively that we have found $\psi_i$ for $i<n$ so
  that for $1<i<n$,
  $(\psi\circ\phi)_i=0.$
  Observe that
  \begin{align*}
    (\psi\circ\phi)_n(a_1\otimes\dots\otimes a_n)
    &=\psi_n(\phi_1(a_1)\otimes\dots\otimes\phi_1(a_n))+\sum_{i=1}^{n-1}\psi_i(F^\phi_i(a_1\otimes\dots\otimes a_n)).
  \end{align*}
  (Here $F^\phi_i$ is the component of $F^\phi$ that lands in
  $\Blg^{\otimes i}\subset \Tensor^*\Blg$.)
  So, if we set
  \begin{align*}
    \psi_n(b_1\otimes\dots\otimes
    b_n)&=\sum_{i=1}^{n-1}\psi_i(F^\phi_i(\phi_1^{-1}(b_1)\otimes\dots\otimes
    \phi_1^{-1}(b_n))),
  \end{align*}
  we have $(\psi\circ\phi)_n=0.$ Continuing in this way, we construct a
  map $\psi$ so that $\psi\circ\phi=\Id_\Alg$. We can also find a map
  $\psi'$ so that $\phi\circ\psi'=\Id_\Blg$ similarly; it follows, of
  course, that $\psi=\psi'$.
  
  It remains to check that $\psi$ satisfies the relation for
  $\Ainf$-homomorphisms, or
  equivalently that $F^\psi$ is a chain map. But we already know that
  $F^\psi=(F^\phi)^{-1}$, so the result follows from the fact that the
  inverse of a chain isomorphism is a chain map.
\end{proof}

\begin{definition}
\label{def:quasiIsomorphism}
An $\Ainf$-algebra homomorphism $\phi:\Alg \to \Blg$ is called a
\emph{quasi-isomorphism} if $\phi_1$ induces an isomorphism from the
homology of~$A$ with respect to $\mu_1^\Alg$ to the homology of~$B$
with respect to $\mu_1^\Blg$.
\end{definition}

\subsubsection{Induced \textalt{$\Ainf$}{A-infty}-algebra structures}
\label{sec:induced-Ainf-alg}
The following is the main lemma of homological perturbation theory;
see \cite[Section
 (1i)]{SeidelBook},~\cite[Section
6.4]{KontsevichSoibelman01:HMS}, \cite[Section
 3.3]{AinftyAlg},~\cite[p. 4]{Keller:OtherAinfAlg}.
For the proof we refer the reader to the references.
\begin{proposition}\label{prop:homol-pert-thy} Let $\Alg$ be an $\Ainf$-algebra and $B_*$ a chain
  complex over $\Ground$.
  \begin{enumerate}
  \item If $f\co B_*\to\Alg$ is a homotopy equivalence of chain
    complexes then there is an $\Ainf$-algebra structure $\{\mu_i\}$
    on $B_*$ and maps $f_i\co B^{\otimes i}\to \Alg[1-i]$ ($i\geq 1$) so
    that
    \begin{itemize}
    \item $\mu_1$ is the differential on $B_*$,
    \item $f_1$ is the given chain map $f$, and
    \item $\{f_i\}\co (B_*,\{\mu_i\})\to \Alg$ is a quasi-isomorphism.
    \end{itemize}
    Moreover, if $B_*$ is a \dg algebra and $f_1$ is an algebra map
    then we can choose $\{\mu_2\}$ to be the multiplication on $B_*$.
  \item If $g\co \Alg\to B_*$ is a homotopy equivalence of chain
    complexes then there is an $\Ainf$-algebra structure $\{\mu_i\}$
    on $B_*$ and maps $g_i\co \Alg^{\otimes i}\to B[1-i]$ so
    that
    \begin{itemize}
    \item $\mu_1$ is the differential on $B_*$,
    \item $g_1$ is the given chain map $g$, and
    \item $\{g_i\}\co \Alg\to (B_*,\{\mu_i\})$ is a quasi-isomorphism.
    \end{itemize}
  \end{enumerate}
\end{proposition}

\begin{remark}
  Note that over $\Ground$ there is no distinction between homotopy
  equivalences and quasi-isomorphisms of chain complexes. 
\end{remark}

\begin{corollary}\label{cor:AinfOnHomology}
  Let $\Alg$ be an $\Ainf$-algebra and $\HAlg$ the homology of $\Alg$,
  which inherits an associative algebra structure from $\Alg$. Let
  $i\co \HAlg\to\Alg$ denote an inclusion choosing a representative for
  each homology class and $p\co \Alg\to \HAlg$ any projection that sends each
  cycle to its homology class. Then there is an $\Ainf$-algebra
  structure $\HAlg$ consisting of maps $\{\mu_i\}$ on $\HAlg$ such that
  there are
  $\Ainf$-quasi-isomorphisms $f\co\HAlg\to\Alg$ and
  $g\co\Alg\to\HAlg$ extending $i$ and $p$ respectively.
\end{corollary}
\begin{proof}
  The existence of $\{\mu_i\}$ and either $f$ or $g$ is immediate from
  Proposition~\ref{prop:homol-pert-thy}; it remains to show that the
  $\Ainf$-algebra structures on $\HAlg$ given by the two parts of
  Proposition~\ref{prop:homol-pert-thy} can be chosen to be the
  same.

  Proposition~\ref{prop:homol-pert-thy} gives us two $\Ainf$-algebra
  structures $\HAlg$ and $\HAlg'$ on $H$ and $\Ainf$-quasi-isomorphisms
  $f\co \HAlg\to \Alg$, $g\co\Alg\to\HAlg'$. Observe that $g_1\circ
  f_1$ is the identity map on $H$. So, by
  Lemma~\ref{lemma:1isos-are-isos}, $g\circ f\co\HAlg\to\HAlg'$ is an
  isomorphism. Then $f\co\HAlg\to\Alg$ and $(g\circ f)^{-1}\circ
  g\co\Alg\to\HAlg$ are the desired maps. 
\end{proof}

\begin{remark}
  Even if $\Alg$ is nilpotent, the induced
  $\Ainf$-structure on the homology $\HAlg$ of $\Alg$ may not be nilpotent (or even satisfy
  the weaker condition that $\mu_i = 0$ for sufficiently
  large~$i$). However, for the algebras $\Alg(\PMC)$ of interest in
  this
  paper, the induced $\Ainf$-structures on $\HAlg$ will be
  nilpotent; see also Remark~\ref{rmk:our-H-is-nilpotent}.
\end{remark}

\begin{definition}
  We call the quasi-isomorphisms $f\co\HAlg\to\Alg$ given by
  Corollary~\ref{cor:AinfOnHomology} \emph{standard
    quasi-isomorphisms}.
\end{definition}

Note that the higher products on $\Halg$ depended on some choices. In
certain situations, however, they are canonically defined; we discuss
one instance of this (which will play a role in
Section~\ref{sec:homology-algebra}).  Let $\Alg$ be an $\Ainf$-algebra,
with products $\mu_i$, and $\HAlg$ be its homology algebra, with
products $\overline{\mu}_i$.

\begin{definition}
  \label{def:MasseyAdmissible}
  A sequence $\alpha_1,\dots,\alpha_m\in\HAlg$ is said to be {\em Massey admissible} 
  if for any $1\leq i<j \leq m$ with $(i,j)\neq(1,m)$, we have
  ${\overline \mu}_{j-i+1}(\alpha_i,\alpha_{i+1},\dots,\alpha_{j})=0$ and
  $\HAlg_{g(i,j)+1}=0$, where $g(i,j)$ is the grading of
  ${\overline \mu}_{j-i+1}(\alpha_i,\dots,\alpha_{j})$, 
  i.e., $g(i,j)=j-i-1+\gr(\alpha_i)+\cdots+\gr(\alpha_j)$.
\end{definition}

\begin{lemma}
  \label{lem:MasseyAdmissible}
  Let $\Alg$ be a $\dg$ algebra, and let $\Halg$ denote  its homology.
  Suppose $\alpha_1,\dots,\alpha_m\in\Halg$ is Massey admissible. Then, there are 
  elements $\xi_{i,j}\in\Alg$ for $0\leq i< j\leq m$ and $(i,j)\neq (0,m)$
  such that
  $$d\xi_{i,j}=\sum_{i<k<j}\xi_{i,k}\cdot \xi_{k,j},$$
  and where, for $i=1,\dots,m$, 
  $\xi_{i-1,i}$ is a cycle representing the homology class $\alpha_i$.
  Moreover, 
  ${\overline \mu}_m(\alpha_1,\dots,\alpha_m)$ is represented by the cycle
  \begin{equation}\label{eq:Massey-prod}
    \sum_{0<k<m}\xi_{0,k}\cdot \xi_{k,m}.
  \end{equation}
  The homology class 
  of this cycle is independent of the choices of the $\xi_{i,j}$.
\end{lemma}

\begin{proof}
  Fix a standard quasi-isomorphism $f\co \Halg\to \Alg$
  and consider the $\Ainf$
  relation for the map $f$, with inputs $\alpha_1,\dots,\alpha_m$. 
  Since the target is a $\dg$ algebra, this relation contains no
  trees with more than two nodes labelled by $f$.
  Indeed, there are the following four types of trees:
  \begin{enumerate}
    \item \label{case:Multiply}
      trees where there are two nodes labelled by $f$,
      whose two outputs get multiplied in $\Alg$;
    \item \label{case:Vanishes}
      trees with one node which is a multiplication ${\overline\mu}_i$ of 
      some proper (consecutive) subsequence of $\alpha_1,\dots,\alpha_m$,
      followed by a node labelled $f_{m-i+1}$;
    \item 
      \label{case:HigherMult}
      the tree representing $f_1({\overline
        \mu}_m(\alpha_1,\dots,\alpha_m))$; and
    \item 
      \label{case:NullHomologous} the tree representing $\mu_1(f_m(\alpha_1,\dots,\alpha_m))$.
    \end{enumerate}
    Now set $\xi_{i,j}=f_{j-i}(\alpha_{i+1},\dots,\alpha_{j})$.  We
    can interpret the sum of trees of Type~\eqref{case:Multiply} as
    the sum appearing in Equation~\eqref{eq:Massey-prod}.  Massey
    admissibility guarantees that higher multiplication
    ${\overline\mu}_i$ on a proper consecutive subsequence vanishes,
    and hence that terms of Type~\eqref{case:Vanishes} vanish.  The
    term of Type~\eqref{case:HigherMult} gives a cycle representing
    $\mu_m(\alpha_1,\dots,\alpha_m)$. The term of
    Type~\eqref{case:NullHomologous} evidently gives a coboundary.  It
    follows that for one choice of the $\xi_{i,j}$ the lemma holds.

    It remains to show that the homology class is independent of all
    the choices made, and hence that, for a Massey admissible
    sequence, ${\overline\mu}_m(\alpha_1,\dots,\alpha_m)$ is
    independent of the choice of compatible $\Ainf$ algebra structure
    on $H_*(\Alg)$. To this end, we show that if for some $c$, we
    exchange exactly one of the $\xi_{a,b}$ for $b-a\leq c$ by
    $\xi_{a,b}'$, then we can complete this to a system of
    $\xi_{i,j}'$ so that the final result changes by a coboundary.
    Specifically, note that $\xi_{a,b}-\xi_{a,b}'$ is a cycle, and
    it is supported in grading $g(a,b)+1$. Hence, by Massey
    admissibility, in fact $\xi_{a,b}-\xi_{a,b}'=d\eta_{a,b}$ for some
    choice of $\eta_{a,b}$. Now, define, for all $i\leq a$ and $j\geq
    b$,
    \begin{align*}
      \xi'_{a,j}&=\xi_{a,j}+\eta_{a,b}\cdot \xi_{b,j} \\
      \xi'_{i,b}&=\xi_{i,b}+\xi_{i,a}\cdot \eta_{a,b},
    \end{align*} and all other $i< j$, $\xi'_{i,j}=\xi_{i,j}$.  It is
    straightforward to check that the $\xi'_{i,j}$ satisfy the same
    equations as the original $\xi_{i,j}$, and also that
    $$\sum_{0<k<m} \xi'_{i,k}\cdot \xi_{k,j}'
    - \sum_{0<k<m} \xi_{i,k}\cdot \xi_{k,j}
    = \left\{\begin{array}{ll}
        d(\eta_{0,{m}}) & {\text{if $1=a$, $b=m$}} \\
        d(\eta_{0,b}\cdot \xi_{b,j}) & {\text{if $1=a$, $b<m$}} \\
        d(\xi_{0,a}\cdot \eta_{a,m}) & {\text{if $1<a$, $b=m$}} \\
        0 & {\text{otherwise.}}
      \end{array}\right.
      $$
      It is easy to see that we can go between any two solutions
      $\{\xi_{i,j}\}$ and $\{\xi'_{i,j}\}$ by a sequence of moves of the 
      above type, and each step leaves the homology class of the expression
      from Equation~\eqref{eq:Massey-prod} unchanged.
\end{proof}

\begin{remark}
  Equation~\eqref{eq:Massey-prod} is the traditional definition of the
  Massey product, which is typically defined only up to some
  indeterminacy.  The Massey admissibility condition guarantees that
  there is no ambiguity in its definition.
\end{remark}

\subsection{Modules over \textalt{$\Ainf$}{A-infty}-algebras}
\label{sec:infty-modules}
In this section we define various notions of $\Ainf$-modules and
bimodules which are used throughout the paper. As noted earlier, our
treatment is nonstandard in that we introduce an object which we call
a ``type $D$ structure,'' which the reader can think of
as a type of  projective
module, see Corollary~\ref{cor:tensor-projective}. 

Another slightly unusual feature of our treatment is that we will
define our categories of modules as \dg categories. We start by
reviewing this notion.
\subsubsection{Background on \dg categories}
\label{sec:dg-categories}
The material in this section is standard, but perhaps unfamiliar to
the low-dimensional topology community. Our treatment is drawn
from~\cite{Keller06:DGCategories}, to which we refer the reader for
more details and further results.
\begin{definition}
  A \emph{differential graded category} is a category
  $\Cat$ such that the morphism spaces are chain complexes
  and composition of morphisms is bilinear and commutes with
  the differential, i.e., such that composition of functions gives
  chain maps $\circ\co\Mor(y,z)\otimes_\Ground\Mor(x,y)\to \Mor(x,z)$.
\end{definition}

The prototypical example is the category of chain complexes:
\begin{example}
  \label{ex:MorTypeD}
  The \emph{\dg category of chain complexes} has objects chain
  complexes $C_*$ and morphism spaces
  \[
  \Mor(C_*,D_*)_n=\{f=(f_i\co C_i\to D_{n+i})_{i=-\infty}^\infty\}
  \]
  with differential defined by
  \[
  (\bdy f)(x)=\bdy_D(f(x))+f(\bdy_C x).
  \]

  Note that the $0$-cycles in $(\Mor(C_*,D_*),\bdy)$ are exactly the
  degree $0$ chain maps between $C_*$ and $D_*$ and the boundaries are
  the null-homotopic chain maps. The homology of
  $(\Mor(C_*,D_*),\bdy)$ is the group of chain maps modulo homotopy.
\end{example}

\begin{definition}\label{def:homomorphisms}
  Given a \dg category $\Cat$, let $\ZMod(\Cat)$ (respectively
  $\ZMod_*(\Cat)$) denote the category with the same objects as $\Cat$
  and morphisms $\Hom_{\ZMod(\Cat)}(x,y)=Z_0(\Mor_{\Cat}(x,y))$
  (respectively $\Hom_{\ZMod_*(\Cat)}(x,y)=Z_*(\Mor_{\Cat}(x,y))$) the
  degree $0$ cycles (respectively cycles of any degree) in $\Cat$'s
  morphism space. We call the morphisms in $\ZMod_*(\Cat)$ the
  \emph{homomorphisms in $\Cat$}, and sometimes denote the set of
  homomorphisms simply by $\Hom$ (as distinct from the set of all
  morphisms $\Mor$).

  Let $\HMod(\Cat)$ (respectively $\HMod_*(\Cat)$) denote the category
  with the same objects as $\Cat$ and morphisms
  $\Hom_{\HMod(\Cat)}(x,y)=H_0(\Mor_{\Cat}(x,y))$ (respectively
  $\Hom_{\HMod_*(\Cat)}(x,y)=H_*(\Mor_{\Cat}(x,y))$), the degree $0$
  homology (respectively total homology) of $\Cat$'s morphism space.
\end{definition}
\begin{example}
  For $\Cat$ the \dg category of chain complexes, $\ZMod(\Cat)$ is the
  usual category of chain complexes, in which the morphism spaces are
  the degree $0$ chain maps. The category $\HMod(\Cat)$ is the
  homotopy category of chain complexes.
\end{example}
The category $\HMod(\Cat)$ is naturally a triangulated category; see
\cite[Section 3.4]{Keller06:DGCategories}. One could go further and
invert quasi-isomorphisms, but for our purposes this will not be necessary.

\begin{definition}
  Let $\Cat$ be a \dg category. Morphisms $f, g\in\Mor_\Cat(x,y)$ are
  called \emph{homotopic} if there is a morphism $h\in\Mor_\Cat(x,y)$
  so that $(\bdy h)=f-g$; in this case we write $f\sim g$. A cycle
  $f\in\Mor_\Cat(x,y)$ is a \emph{homotopy equivalence} if there is a
  cycle $g\in\Mor_\Cat(y,x)$ so that $g\circ f\sim\Id_x$ and $f\circ
  g\sim\Id_y$.
\end{definition}

\begin{definition}
  For $f \in \Mor(x,y)$, let 
  \begin{align*}
  f_{*w} \co \Mor(w,x) \to \Mor(w,y)\\
  f^{*z} \co \Mor(y,z) \to \Mor(x,z)    
  \end{align*}
  be the maps obtained by pre-
  and post-composing with~$f$.  We will sometimes write $f_*$ or $f^*$
  if the spaces are clear from context.
\end{definition}

\begin{lemma}
  \label{lem:HomotopyEquivalence}
  If $f,g\in\Mor_\Cat(x,y)$ are homotopic morphisms then
  $f_{*w} \sim g_{*w}$ and
  $f^{*z} \sim g^{*z}$.

  Similarly, if $f\in\Mor_\Cat(x,y)$ is a homotopy equivalence then
  the maps $f_{*w}$ and $f^{*z}$ are chain homotopy
  equivalences.
\end{lemma}
The proof is straightforward.

\begin{definition}
  Let $\Cat$ and $\Dat$ be \dg categories. A \emph{\dg functor
    ${\mathcal F}\co\Cat\to\Dat$} is a functor $\Cat\to\Dat$ such that for any
  objects $x$ and $y$ of $\Cat,$ ${\mathcal F}(\Mor(x,y))$ is a degree $0$ \dg
  module homomorphism.
\end{definition}

\begin{lemma}\label{lemma:dg-funct-blah}
  If ${\mathcal F}$ is a \dg functor and $f,g\in\Mor(x,y)$ are homotopic then
  ${\mathcal F}(f)$ and ${\mathcal F}(g)$ are homotopic. If $f$ is a homotopy equivalence
  then ${\mathcal F}(f)$ is a homotopy equivalence.
\end{lemma}
The proof is immediate from the definitions.

\begin{definition}\label{def:homotopic-functors}
  Let $\Cat$ and $\Dat$ be \dg categories and ${\mathcal F}, {\mathcal G}\co \Cat\to \Dat$
  \dg functors. We say ${\mathcal F}$ is \emph{homotopic to ${\mathcal G}$} if for each
  $x\in\ob(\Cat)$ there are homotopy equivalences 
  $\eta_x\co {\mathcal F}(x)\to {\mathcal G}(x)$ so that
  for any $x,y\in\ob(\Cat)$ the following diagram commutes up to
  homotopy:
  \[
  \begin{tikzpicture}
    \node at (0,0) (xy) {$\Mor_\Cat(x,y)$};
    \node at (5,0) (FxFy) {$\Mor_\Dat({\mathcal F}(x),{\mathcal F}(y))$};
    \node at (0,-1.5) (GxGy) {$\Mor_\Dat({\mathcal G}(x),{\mathcal G}(y))$};
    \node at (5,-1.5) (FxGy) {$\Mor_\Dat({\mathcal F}(x),{\mathcal G}(y)).$};
    \draw[->] (xy) to node[above]{\lab{{\mathcal F}}} (FxFy);
    \draw[->] (xy) to node[left]{\lab{{\mathcal G}}} (GxGy);
    \draw[->] (FxFy) to node[right]{\lab{(\eta_y)_*}} (FxGy);
    \draw[->] (GxGy) to node[below]{\lab{(\eta_x)^*}} (FxGy);
  \end{tikzpicture}
  \]

  A \dg functor ${\mathcal F}\co\Cat\to\Dat$ is a \emph{homotopy equivalence} if
  there is a functor ${\mathcal G}\co\Dat\to\Cat$ so that ${\mathcal G}\circ {\mathcal F}$ is homotopic
  to $\Id_{\Cat}$ and ${\mathcal F}\circ {\mathcal G}$ is homotopic to $\Id_\Dat$.
\end{definition}
(In Definition~\ref{def:homotopic-functors} we have not required any
coherence for the homotopies in the diagram; one could formulate
stronger notions with such coherence built in.)

\begin{definition}\label{def:quasi-equiv}
  Let $\Cat$ and $\Dat$ be \dg categories. A functor
  ${\mathcal F}\co\Cat\to\Dat$ is a \emph{quasi-equivalence} if
  \begin{itemize}
  \item for all $x,y\in\ob(\Cat)$, the map ${\mathcal F}(x,y)\co \Mor_\Cat(x,y)\to
    \Mor_\Dat({\mathcal F}(x),{\mathcal F}(y))$ is a quasi\hyp isomorphism and
  \item the induced map $\mathsf{H}({\mathcal F})\co
    \mathsf{H}(\Cat)\to\mathsf{H}(\Dat)$ is an equivalence of
    categories.
  \end{itemize}
\end{definition}
(See~\cite[Section 2.3]{Keller06:DGCategories} for more details.)

\begin{proposition}\label{prop:dg-homotopy-equiv-quasi-equiv}
  If ${\mathcal F}\co \Cat\to\Dat$ is a homotopy equivalence then ${\mathcal F}$ is a
  quasi-equivalence.
\end{proposition}
Again, the proof is straightforward.

Just as one can generalize the notion of \dg algebras to
$\Ainf$-algebras, one can generalize the notion of \dg categories and
\dg functors to $\Ainf$-categories and $\Ainf$-functors. 
Some of the categories studied in this paper (in particular,
the category of type~$D$ modules) are $\Ainf$-categories, and some of
the functors (in
particular, $\DT$) are $\Ainf$-functors. Most of the additional complications
are,
however, not important for the applications in this paper: when
working with \dg algebras (rather than $\Ainf$-algebras), the categories
we consider are honest \dg categories (see Remark~\ref{rmk:TypeDsimple}).
So, we will not
spell out the notions of $\Ainf$-categories and $\Ainf$-functors,
trusting the reader to provide them if desired, or to consult,
e.g.,~\cite{Seidel02:FukayaDef}.

\begin{remark}
  The reader might find the following analogy helpful for
  understanding the role of these \dg categories. It has been
  understood for some time that when working with complexes, rather
  than taking homology of complexes it is often better to pass to the
  derived category, i.e., to invert morphisms which induce
  isomorphisms on homology. For example, operations like tensor
  product and $\Hom$ are better behaved with respect to inverting
  quasi-isomorphisms than with respect to taking homology,
  giving rise to the derived functors $\Tor$ and
  $\Extop$. The language of $\Ainf$-algebras allows one to view the
  derived category of $R$-modules itself as the homology of a \dg
  category. Then, if one is interested in studying categories of
  modules, it is better to work with \dg categories and invert
  quasi-equivalences of categories rather than take homology and work
  with derived categories.
\end{remark}

\subsubsection{The category of \textalt{$\Ainf$}{A-infty}-modules}
\label{sec:categ-Ainf-modules}
\begin{definition}\label{def:ainf-mod}
  A \emph{(right) $\Ainf$-module $M_\Alg$ over $\Alg=(\DGA,\{\mu_i\})$} is a
  $\ZZ$-graded $\Ground$-module $M$ together with degree $0$ $\Ground$-linear maps
  $m_{j+1}\co
  M\otimes \DGA[1]^{\otimes j}\to M[1]$ ($j=0,\dots,\infty$) such that for
  each $i=0,\dots,\infty$, $x\in M$ and $a_1,\dots,a_i\in A$,
  \begin{multline}
    \label{eq:inf-mod-rel}
    0=\sum_{j=0}^{i-1}
    m_{i-j}(m_{j+1}(x,a_1,\dots,a_j),a_{j+1},\dots,a_i)\\
    +\sum_{j=1}^i\sum_{k=1}^{i-j+1}m_{i-j+1}(x,a_1,\dots,a_{k-1},\mu_j(a_k,\dots,a_{k+j}),a_{k+j+1},\dots,a_i).
  \end{multline}
  An $\Ainf$-module is \emph{strictly unital} if $m_2(x,1) = x$ and
  $m_{i+1}(x,a_1,\dots,a_i)=0$ if $i\neq 1$ and one of the $a_i\in\Ground$.
\end{definition}

We will sometimes refer to $\Ainf$-modules as \emph{type $A$ modules},
to place them on equal footing with the type $D$ modules which will
appear later (Definition~\ref{def:TypeD}). The bordered
invariant $\CFAa(Y)$ is an $\Ainf$-module.

\begin{convention}
  All $\Ainf$-modules will be assumed strictly unital.
\end{convention}

As with the $\Ainf$-relation for algebras
(Equation~(\ref{eq:a-inf-alg-rel})), Equation~(\ref{eq:inf-mod-rel})
has an interpretation in terms of trees. In this interpretation, 
there are two types of strands: those corresponding to algebra elements,
and distinguished (dotted) ones corresponding to module elements. 
Precisely:
\begin{definition}
  A (right) \emph{$\Ainf$-module operation tree} is a finite directed
  tree~$\Gamma$ embedded in the plane so that all edges point
  downwards, with the edges of two types, the algebra edges (labelled
  by `$A$') or
  module edges (labelled by `$M$') and
  non-leaf vertices marked either `$\mu$' or `$m$', and so that:
  \begin{itemize}
  \item Each non-leaf vertex has exactly one outgoing edge.
  \item Each `$\mu$' vertex touches only algebra edges.
  \item At each `$m$' vertex, the leftmost incoming edge and the
    outgoing edge are module edges, and the other incoming edges are
    algebra edges.
  \item There is at least one module edge.
  \end{itemize}
\end{definition}

To an $\Ainf$-module operation tree we can associate an operation $m_\Gamma \co M \otimes
A^{\otimes n} \to M$ by flowing along the edges as before, applying
$\mu_i$ or $m_i$ depending on the label at the vertex.  Then
Equation~\eqref{eq:inf-mod-rel} says that the sum
of $m_\Gamma$ over all $\Ainf$-module operation trees~$\Gamma$ with
two vertices and a fixed number
of inputs vanishes.

An $\Ainf$-module operation tree is called \emph{spinal} if each node
is labelled `$m$', i.e., the line of module edges goes through each
node.

\begin{definition}\label{def:A-bounded}
  We say that an $\Ainf$-module $M_\Alg$ is \emph{operationally
    bounded} (or just
  \emph{bounded}) if, for all $x \in M$ there exists an $n$ so that
  for any $i>n$ and any spinal $\Ainf$-module operation tree~$\Gamma$ with
  $i+1$ input edges,
  $m_\Gamma(x\otimes \cdot)$ vanishes on $(A_+)^{\otimes i}$.
\end{definition}
Graphically, this says that for each $x$ there exists an $n$ so that
if $i_1+\cdots+i_k>n$ then for any $a_{1,1},\dots,a_{k,i_k}\in A_+$,
\[
\mathcenter{
\begin{tikzpicture}
  \node at (0,0) (x) {$x$};
  \node at (0,-1.5) (mu1) {$\mu_{i_1+1}$};
  \node at (0,-3) (dots) {$\vdots$};
  \node at (0,-4.5) (muk) {$\mu_{i_k+1}$};
  \node at (0,-6) (bottom) {$0\mathrlap{.}$};
  \node at (2.5,0) (a1) {$a_{1,1}\otimes\cdots\otimes a_{1,i_1}$};
  \node at (5,0) (cdots) {$\cdots$};
  \node at (7.5,0) (ak) {$a_{k,1}\otimes\cdots\otimes a_{k,i_k}$};
  \draw[Amodar] (x) to (mu1);
  \draw[Amodar] (mu1) to (dots);
  \draw[Amodar] (dots) to (muk);
  \draw[Amodar] (muk) to (bottom);
  \draw[tensoralgarrow] (a1) to (mu1);
%  \draw[tensoralgarrow] (cdots) to (dots);
  \draw[tensoralgarrow] (ak) to (muk);
\end{tikzpicture}}
\] 

The notion of operationally bounded modules
is different from more traditional definitions of boundedness; for
instance, it does not imply any bound on degrees which have non-zero
homology groups.
\begin{remark}\label{rem:mod-bounded}
  As in Remark~\ref{rem:ainf-bounded}, this is a stronger condition
  than the condition called ``operationally bounded'' in our previous
  paper~\cite{LOT1}, where we just assumed that $m_i = 0$ for $i$
  sufficiently large.  But note that a \dg module over a nilpotent \dg
  algebra is automatically operationally bounded in the stronger sense
  above. The reason we want the stronger condition relates to
  subtleties of the box tensor product of bimodules, not modules.
\end{remark}

We can combine the multiplications $m_i$ on an $\Ainf$-module $M_\Alg$
into a single map $m\co M \otimes \Tensor^*(A[1]) \to M[1]$. 
Because we assume $M_\Alg$ is
strictly unital, all the information is contained in a map $M \otimes
\Tensor^*(A_+[1])\to M[1]$, also denoted~$m$.
If $M_\Alg$ is operationally bounded then the map $m$ extends to a map
from the completed tensor product: $m\co M \otimes
\overline{\Tensor^*}(A[1])\to M[1]$.  (Here, 
$\overline{\Tensor^*}(A[1])=\prod_{i=0}^\infty A[1]^{\otimes i}$.)

\begin{definition}\label{def:Ainf-mod-homs}
  The category $\ModCat_\Alg$ of (right) $\Ainf$-modules over $\Alg$
  is the \dg category whose objects are $\Ainf$-modules $M_\Alg$ and
  whose morphism spaces are defined as follows.
  The $\ZZ$-graded vector spaces underlying the morphism spaces are 
  the vector spaces of $\Ground$-module homomorphisms
  $$\Mor_\Alg(M_{\Alg},N_{\Alg})\coloneqq\Hom_{\Ground}(M\otimes \Tensor^*(A_+[1]),N),$$
  where $A_+$ is the augmentation ideal of $A$.
  Here we think of the
  homomorphisms between two $\ZZ$-graded $\Ground$-modules $V$ and $W$ as
  a graded $\Field$ vector space, with
  \begin{equation*}
    \Hom_\Ground(V, W)_i \coloneqq \bigoplus_j \Hom_\Ground(V_j, W_{i+j}).
  \end{equation*}
  If $h \in \Mor_\Alg(M_{\Alg}, N_{\Alg})$, we denote the component
  parts by $h_i: M \otimes (A_+[1])^{\otimes (i-1)} \to N$.

  The differential on morphisms is given by
  \[
  \begin{split}
    (\bdy h)(x,a_1,\dots,a_n)&\coloneqq\sum_{i+j=n}h_{i+1}(m_{j+1}(x,a_1,\dots,a_j),a_{j+1},\dots,a_n)\\
      &\quad+\sum_{i+j=n}m_{i+1}(h_{j+1}(x,a_1,\dots,a_j),a_{j+1},\dots,a_n)\\
      &\quad+\sum_{i+j=n}\sum_{k=1}^{n-j}h_{i+1}(x,a_1,\dots,a_{k-1},\mu_{j+1}(a_k,\dots,a_{k+j}),a_{k+j+1},\dots,a_n).
  \end{split}
  \]

  The \emph{$\Ainf$-homomorphisms} from $M_\Alg$ to $N_\Alg$ in the
  usual sense (see for example~\cite{AinftyAlg}) are the cycles in the morphism complex;
  cf.~Definition~\ref{def:homomorphisms}.
\end{definition}
\begin{remark}
  Strict unitality is built into the definition 
  of $\Ainf$-morphisms, since we use the augmentation ideal in the
  definition of the morphism spaces.
\end{remark}

We can again give a graphical notation for operations built from
morphisms, as follows.  Let $h$ be an $\Ainf$-morphism from $M_\Alg$
to $N_\Alg$.  An \emph{$\Ainf$-module morphism tree}~$\Gamma$ is a planar
directed tree satisfying all the conditions of an $\Ainf$-module
operation tree, except that the module edges may be marked either
`$M$' or `$N$', and there is one node labelled by `$h$' with an `$M$' input and
`$N$' output which is otherwise like an `$m$' node.  Again, we can
define a map
$h_\Gamma: M \otimes
(A_+[1])^{\otimes n} \to N$ by applying at each vertex $\mu_i$, $h_i$, or
$m_i$ in either $M$ or $N$, as appropriate.

\begin{definition}
  An $\Ainf$-morphism $h$ from $M_\Alg$ to $N_\Alg$ is
  \emph{operationally bounded} if for each $x\in M$ there is an $n$ so that
  $h_\Gamma(x\otimes \cdot)$ vanishes on $(A_+)^{\otimes i}$ for all spinal
  $\Ainf$-module morphism trees
  $\Gamma$ with $i > n$ inputs.
\end{definition}

We can also view a morphism $(h_j)_{j=1}^\infty$ as a single map
$h\co M\otimes\Tensor^*(A[1])\to N$. 
Writing $\Delta\co \Tensor^*(V)\to \Tensor^*(V)\otimes
\Tensor^*(V)$ for the natural comultiplication for any vector
space~$V$, we can also draw the
differential of a morphism $h$ as
\[
\bdy h=
\mathcenter{\begin{tikzpicture}
    \node at (0,0) (tlblank) {};
    \node at (2,0) (trblank) {};
    \node at (1,-1) (Delta) {$\Delta$};
    \node at (0,-2) (m) {$m$};
    \node at (0,-3) (h) {$h$};
    \node at (0,-4) (blblank) {};
    \draw[Amodar] (tlblank) to (m);
    \draw[Amodar] (m) to (h);
    \draw[Amodar] (h) to (blblank);
    \draw[tensoralgarrow] (trblank) to (Delta);
    \draw[tensoralgarrow] (Delta) to (m);
    \draw[tensoralgarrow] (Delta) to (h);
\end{tikzpicture}}
+
\mathcenter{\begin{tikzpicture}
    \node at (0,0) (tlblank) {};
    \node at (2,0) (trblank) {};
    \node at (1,-1) (Delta) {$\Delta$};
    \node at (0,-3) (m) {$m$};
    \node at (0,-2) (h) {$h$};
    \node at (0,-4) (blblank) {};
    \draw[Amodar] (tlblank) to (h);
    \draw[Amodar] (h) to (m);
    \draw[Amodar] (m) to (blblank);
    \draw[tensoralgarrow] (trblank) to (Delta);
    \draw[tensoralgarrow] (Delta) to (m);
    \draw[tensoralgarrow] (Delta) to (h);  
\end{tikzpicture}}
+
\mathcenter{
  \begin{tikzpicture}
    \node at (0,0) (tlblank) {};
    \node at (2,0) (trblank) {};
    \node at (0,-2) (h) {$h$};
    \node at (1,-1) (D) {$\overline{D}$};
    \node at (0,-3) (blblank) {};
    \draw[Amodar] (tlblank) to (h);
    \draw[Amodar] (h) to (blblank);
    \draw[tensoralgarrow] (trblank) to (D);
    \draw[tensoralgarrow] (D) to (h);
  \end{tikzpicture}
}.
\]
We have used dashed lines for module elements, and solid lines for
algebra elements, either in $\Alg$ in $\Alg_+$.  A doubled arrow
denotes elements of $T^*(A[1])$ or $T^*(A_+[1])$.
    
Given $f\in\Mor_{\Alg}(M,N)$ and $g\in\Mor_{\Alg}(N,L)$, we
define the composite morphism $g\circ f\in \Mor_{\Alg}(M,L)$ by
$$(g\circ f)(x\otimes\alpha)=g\circ 
(f\otimes \Id_{\Tensor^*(\Alg[1])}) (x\otimes \Delta(\alpha)).$$
This induces a chain map
$$\Mor_{\Alg}(N,L)\otimes \Mor_{\Alg}(M,N) \to \Mor_{\Alg}(M,L).$$
We can draw composition as
\[
g\circ f =
\mathcenter{
\begin{tikzpicture}
  \node at (0,0) (tlblank) {};
  \node at (2,0) (trrblank) {};
  \node at (1,-1) (Delta) {$\Delta$};
  \node at (0,-2) (phi) {$f$};
  \node at (0,-3) (psi) {$g$};
  \node at (0,-4) (bblank) {};
  \draw[Amodar] (tlblank) to (phi);
  \draw[Amodar] (phi) to (psi);
  \draw[Amodar] (psi) to (bblank);
  \draw[tensoralgarrow] (trrblank) to (Delta);
  \draw[tensoralgarrow] (Delta) to (phi);
  \draw[tensoralgarrow] (Delta) to (psi);
\end{tikzpicture}}.
\]

\subsubsection{The category of type \textalt{$D$}{D} structures}
\label{sec:cat-type-d-str}

Throughout this subsubsection, let $\Alg=(\DGA,\{\mu_i\})$ be an
$\Ainf$-algebra satisfying Convention~\ref{Convention:AlgsUniConFd}.

\begin{definition}
  \label{def:TypeD}
  A (left) {\em{type $D$ structure}} over $\Alg$ is an object $\lsup{\Alg}{N}$
  consisting of a graded $\Ground$-module $N$ equipped with a
  degree~$0$ linear map
  $\delta^1 \co N \to \DGA[1] \otimes N$, satisfying a 
  compatibility condition which is best described after introducing
  an auxiliary construction.

Define maps $\delta^i \co N \to (\DGA[1])^{\otimes i}
\otimes N$ ($i\geq 2$) inductively by
\[
\delta^{i+1} \coloneqq (\Id_{\Alg^{\otimes i}} \otimes \delta^1) \circ
\delta^i
\]
and define $\delta \co N \to \overline{\Tensor^*}(\DGA[1]) \otimes N$
by 
\begin{equation}
  \label{eq:DefineDelta}
  \delta(x)\coloneqq \sum_{i=0}^\infty \delta^i.
\end{equation} (Here,
$\overline{\Tensor^*}(\DGA[1])=\prod_{i=0}^\infty\Alg[1]^{\otimes
  i}$.)
Then the compatibility condition is that for any $x \in N$,
\begin{equation}\label{eq:typeD-def}
  (\mu \otimes \Id_N) \circ \delta(x) = 0.
\end{equation}
Here, $\mu$ is the sum of the structure maps of $\Alg$, as in Formula~(\ref{eq:def-mu}).

  We say that $\lsup{\Alg}N$ is \emph{operationally bounded} if 
  for each $x\in \lsup{\Alg}N$, there is a constant $n=n(x)$
  with the property that for all $i>n$, $\delta^i(x)=0$.
  This is equivalent to saying that the above ${\delta}$ 
  factors thorough the inclusion 
  of $\Tensor^*(\DGA[1])$ in ${\overline{\Tensor^*}}(\DGA[1])$.

  When more than one module is present, we will often write
  $\delta^{N,i}$ for the operation $\delta^i$ on the type $D$
  structure $N$ and $\delta^N$ to denote the map $\delta$ on $N$. This
  conflict of notation with $\delta^i$ should not cause confusion.
\end{definition}
Sometimes we refer to type $D$ structures as \emph{type $D$ modules}.

Note that the defining equation~(\ref{eq:typeD-def}) makes sense since
$\Alg$ is nilpotent. We can represent Definition~\ref{def:TypeD}
graphically by
\[
\delta\,=
\mathcenter{
  \begin{tikzpicture}
    \node at (6,0) (blank0N) {};
    \node at (6,-1.5) (d0N) {$\delta$};
    \node at (6,-3) (blank2N) {};
    \node at (5,-3) (blankh1N) {};
    \draw[modarrow] (blank0N) to (d0N);
    \draw[modarrow] (d0N) to (blank2N); 
    \draw[tensoralgarrow, bend right=15] (d0N) to (blankh1N);
  \end{tikzpicture}}
\coloneqq
\,\,\mathcenter{
  \begin{tikzpicture}
    \node at (0,0) (tblank) {};
    \node at (0,-3) (bblank) {};
    \draw[modarrow] (tblank) to (bblank);    
  \end{tikzpicture}
}\,\,
+
\mathcenter{
  \begin{tikzpicture}
      \node at (-3,0) (blank01) {};
      \node at (-3,-1) (d01) {$\delta^1$};
      \node at (-3,-3) (blank21) {};
      \node at (-4,-3) (blankh11) {};
      \draw[modarrow] (blank01) to (d01);
      \draw[modarrow] (d01) to (blank21);
      \draw[algarrow, bend right=15] (d01) to (blankh11);    
  \end{tikzpicture}
}
+
  \mathcenter{
    \begin{tikzpicture}[x=1cm,y=32pt]
      \node at (0,0) (blank02) {};
      \node at (0,-1) (d02) {$\delta^1$};
      \node at (0,-2) (d12) {$\delta^1$};
      \node at (0,-3) (blank22) {};
      \node at (-1,-3) (blankh12) {};
      \node at (-.9,-3) (blankh22) {};
      \draw[modarrow] (blank02) to (d02);
      \draw[modarrow] (d02) to (d12);
      \draw[modarrow] (d12) to (blank22);
      \draw[algarrow, bend right=15] (d02) to (blankh12);
      \draw[algarrow, bend right=15] (d12) to (blankh22);
\end{tikzpicture}}
+\ \cdots.
\]

Then, Equation~(\ref{eq:typeD-def}) takes the form:
\[
  \mathcenter{
    \begin{tikzpicture}[x=1cm,y=32pt]
      \node at (1,0) (blank0) {};
      \node at (1,-1) (d0) {$\delta$};
      \node at (1,-3) (blank2) {};
      \node at (0,-2) (mu) {$\mu$};
      \node at (0,-3) (blank3) {};
      \draw[modarrow] (blank0) to (d0);
      \draw[modarrow] (d0) to (blank2);
      \draw[tensoralgarrow, bend right=15] (d0) to (mu);
      \draw[algarrow] (mu) to (blank3);
    \end{tikzpicture}}
  =0.
\]

We can also form type $D$ structures into a \dg or $\Ainf$ category
$\DuModCat{\Alg}$ where $\Mor(\lsup{\Alg}M, \lsup{\Alg}N)$ is the chain
complex whose underlying space consists of maps
\begin{align*}
  h^1 \co M \to \DGA \otimes N.
\end{align*}
Such a map can be upgraded to a map to the tensor algebra: define
$h^i \co M \to \DGA^{\otimes i} \otimes N$ by
\[
h^i \coloneqq \sum_{j=0}^{i-1} (\Id_{\DGA^{\otimes j+1}} \otimes
\delta^{N,i-j-1}) \circ (\Id_{\DGA^{\otimes j}} \otimes h^1) \circ
\delta^{M,j}
\]
and $h \co M \to \overline{\Tensor^*}\DGA \otimes N$ by $h \coloneqq
\sum_{i=1}^\infty h^i.$ The symbol $u$ in the notation
$\DuModCat{\Alg}$ indicates that objects in the category of type
$D$ structures are not required to be bounded (i.e., they are possibly
unbounded; see Definition~\ref{def:bounded-D}).

The degree of an element of $\Mor(\lsup{\Alg}M,\lsup{\Alg}N)$ is the
degree of the corresponding map of vector spaces $M\to \DGA\otimes
N$.

The boundary operator on morphisms is defined by
\[
(\partial h)^1 \coloneqq (\mu \otimes \Id_N) \circ h,
\]
as depicted on the left in Figure~\ref{fig:op-type-d}.
\begin{figure}
  \[
  \mathcenter{
    \begin{tikzpicture}[x=1cm,y=32pt]
      \node at (1,0) (blank0) {};
      \node at (1,-1) (dM) {$\delta^M$};
      \node at (1,-2) (h) {$h^1$};
      \node at (1,-3) (dN) {$\delta^N$};
      \node at (0,-4) (mu) {$\mu$};
      \node at (1,-5) (blank1) {};
      \node at (0,-5) (blank2) {};
      \draw[modarrow] (blank0) to (dM);
      \draw[modarrow] (dM) to (h);
      \draw[modarrow] (h) to (dN);
      \draw[modarrow] (dN) to (blank1);
      \draw[algarrow] (mu) to (blank2);
      \draw[tensoralgarrow, bend right=15] (dM) to (mu);
      \draw[tensoralgarrow, bend right=15] (dN) to (mu);
      \draw[algarrow, bend right=15] (h) to (mu);
    \end{tikzpicture}}
  \qquad\qquad
  \mathcenter{
    \begin{tikzpicture}[x=1cm,y=32pt]
      \node at (1,0) (blank0) {};
      \node at (1,-1) (dM) {$\delta^M$};
      \node at (1,-2) (h) {$h^1$};
      \node at (1,-3) (dN) {$\delta^N$};
      \node at (1,-4) (hp) {$(h')^1$};
      \node at (1,-5) (dP) {$\delta^P$};
      \node at (0,-3) (place1) {};
      \node at (0,-4) (place2) {};
      \node at (0,-5) (place3) {};
      \node at (0,-6) (mu) {$\mu$};
      \node at (1,-7) (blank1) {};
      \node at (0,-7) (blank2) {};
      \draw[modarrow] (blank0) to (dM);
      \draw[modarrow] (dM) to (h);
      \draw[modarrow] (h) to (dN);
      \draw[modarrow] (dN) to (hp);
      \draw[modarrow] (hp) to (dP);
      \draw[modarrow] (dP) to (blank1);
      \draw[algarrow] (mu) to (blank2);
      \draw[tensoralgarrow, bend right=15] (dM) to (place1);
      \draw[tensoralgarrow] (dN) to (place2);
      \draw[tensoralgarrow] (dP) to (mu);
      \draw[tensoralgarrow] (place1) to (place2);
      \draw[tensoralgarrow] (place2) to (place3);
      \draw[tensoralgarrow] (place3) to (mu);
      \draw[algarrow] (h) to (place1);
      \draw[algarrow] (hp) to (place3);
    \end{tikzpicture}}
  \qquad\qquad
  \mathcenter{\begin{tikzpicture}[x=1cm,y=32pt]
      \node at (0,1) (blank0) {};
      \node at (0,0) (d0) {$\delta^{M_0}$} ;
      \node at (0,-1) (h1) {$(h_1)^1$};
      \node at (0,-2) (d1) {$\delta^{M_1}$};
      \node at (0,-3) (h2) {$(h_2)^1$};
      \node at (0,-4) (dots) {$\vdots$};
      \node at (0,-5) (hk) {$(h_k)^1$};
      \node at (0,-6) (dk) {$\delta^{M_k}$};
      \node at (-1,-2) (place2) {};
      \node at (-1,-3) (place3) {};
      \node at (-1,-4) (place4) {};
      \node at (-1,-5) (place5) {$\vdots$};
      \node at (-1,-6) (place6) {};
      \node at (-1,-7) (mu) {$\mu$};
      \node at (-1,-8) (blank) {};
      \node at (0,-8) (blank2) {};
      \draw[modarrow] (d0) to (h1);
      \draw[modarrow] (blank0) to (d0);
      \draw[modarrow] (h1) to (d1);
      \draw[modarrow] (d1) to (h2);
      \draw[modarrow] (h2) to (dots);
      \draw[modarrow] (dots) to (hk);
      \draw[modarrow] (hk) to (dk);
      \draw[modarrow] (dk) to (blank2);
      \draw[tensoralgarrow, bend right=15] (d0) to (place2);
      \draw[algarrow] (h1) to (place2);
      \draw[tensoralgarrow] (d1) to (place3);
      \draw[algarrow] (h2) to (place4);
      \draw[algarrow] (hk) to (place6);
      \draw[tensoralgarrow] (dk) to (mu);
      \draw[tensoralgarrow] (place2) to (place3);
      \draw[tensoralgarrow] (place3) to (place4);
      \draw[tensoralgarrow] (place4) to (place5);
      \draw[tensoralgarrow] (place5) to (place6);
      \draw[tensoralgarrow] (place6) to (mu);
      \draw[algarrow] (mu) to (blank);
    \end{tikzpicture}}
  \]
  \caption{Operations on type~$D$ morphisms.  From left to right, the
    differential on $\Mor(\lsup{\Alg}M,\lsup{\Alg}N)$, the composition
  of morphisms, and the higher ($\Ainf$) compositions.}
  \label{fig:op-type-d}
\end{figure}

The composition of two morphisms $h^1 \in
\Mor(\lsup{\Alg}M,\lsup{\Alg}N)$ and
$(h')^1\in\Mor(\lsup{\Alg}N,\lsup{\Alg}P)$ is defined by
\begin{equation}\label{eq:type-d-composition}
(h' \circ h)^1 \coloneqq
(\mu \otimes \Id_N) \circ
(\Id_{\Tensor^*} \otimes \delta^P) \circ
(\Id_{\Tensor^*} \otimes (h')^1) \circ
(\Id_{\Tensor^*} \otimes \delta^N) \circ
(\Id_{\Tensor^*} \otimes h^1) \circ
\delta^M;
\end{equation}
this is illustrated in the middle in Figure~\ref{fig:op-type-d}.

It is straightforward to verify that this is a map of chain
complexes, i.e., $\partial(h\circ h') = \partial h \circ h' + h
\circ \partial h'$.

The attentive reader will notice that the composition of morphisms is
not associative but associative only up to homotopy. In fact:
\begin{lemma}\label{lem:D-higher-comp}
  There are higher composition maps making $\DuModCat{\Alg}$ into
  an $\Ainf$-category.
\end{lemma}
\begin{proof}
  Let $\lsup{\Alg}M_0$,\dots,$\lsup{\Alg}M_k$ be type $D$ structures
  over $\Alg$. Define a higher composition map
  \[
  \circ_k\co \Mor(\lsup{\Alg}M_0, \lsup{\Alg}M_1)\otimes\cdots\otimes
  \Mor(\lsup{\Alg}M_{k-1}, \lsup{\Alg}M_k)\to \Mor(\lsup{\Alg}M_0, \lsup{\Alg}M_k)
  \]
  by the diagram on the right in Figure~\ref{fig:op-type-d}.

  It is easy to check that the higher composition maps $\circ_k$
  satisfy the $\Ainf$-relation.
\end{proof}

\begin{remark}
  \label{rmk:TypeDsimple}
  In the case where $\Alg$ is a \dg algebra, the higher composition
  maps of Lemma~\ref{lem:D-higher-comp} vanish, so
  $\DuModCat{\Alg}$ is an honest \dg category.
\end{remark}

We call the cycles in $\Mor(\lsup{\Alg}M, \lsup{\Alg}N)$ \emph{type
  $D$ module homomorphisms}.

Note that, since composition of morphisms in
$\Mor(\DuModCat{\Alg})$ is associative up to homotopy, the
homotopy category $\HMod(\DuModCat{\Alg})$ is an honest category.

Finally, we can also form the category of bounded type $D$ structures:
\begin{definition}\label{def:bounded-D}
  We call a type $D$ morphism $h\co \lsup{\Alg}M\to \lsup{\Alg}N$
  \emph{bounded} if for any $x\in \lsup{\Alg}M$,
  there is a constant $n=n(x)$ with the property that
  for all $i\geq n$, 
  $(\Id_{\DGA^{\otimes i}}\otimes h)\circ
  \delta^{M,i}(x)=0$ and $(\Id_{\DGA}\otimes \delta^{N,i})\circ
  h(x)=0$.

  Let $\DbModCat{\Alg}$ denote the category of bounded type $D$
  structures and bounded morphisms.
\end{definition}

In fact, it is clear from the definition that:
\begin{lemma}
  \label{lem:BoundedSubcategory}
  Any morphism between bounded type $D$ structures is bounded.
\end{lemma}

One can, of course, define right type $D$ structures similarly; we
denote the category of right type $D$ structures over $\Alg$ by
$\ModCat^\Alg_{\mathfrak{u}}$.

\begin{definition}\label{def:opposite-type-D}
  Given a type $D$ structure $(\lsup{\Alg}M,\delta^1)$, define the
  \emph{opposite} type~$D$ structure
  $(\overline{M}{}^\Alg,\overline{\delta}{}^1)$ as follows.  First
  suppose that $\Alg$ and $M$ are finite-dimensional.  As a
  $\Ground$-module, $\overline{M}$ is just
  $M^*=\Hom_\Ground(M,\Ground)$. The map $\delta^1$ on~$M$ is an element of
  $\Hom_\Ground(M,A\otimes M)\isom M^*\otimes A\otimes M \isom
  \Hom_\Ground(M^*,M^*\otimes A)$, and
  $\overline{\delta}{}^1$ is obtained by viewing $\delta^1$ as lying in
  $\Hom_\Ground(M^*,M^*\otimes A)$.  If $\Alg$ or $M$ are not finite-dimensional, $\Hom_\Ground(M,
  A \otimes M)$ is a subspace of $\Hom_\Ground(M^*, M^* \otimes
  A)$, and $\overline{\delta}{}^1$ is the image of $\delta^1$ under
  the inclusion.
\end{definition}
\begin{lemma}\label{lem:opposite-is-D}
  Given a type $D$ structure $(M,\delta^1)$, the opposite structure
  $((\overline{M}){}^\Alg,\overline{\delta}{}^1)$ satisfies the type $D$
  structure equation (Equation~(\ref{eq:typeD-def})).
\end{lemma}
\begin{proof}
  The two structure equations are the same, as can be seen most easily
  graphically in the finite-dimensional case:
  \[
  0\,=
  \mathcenter{
    \begin{tikzpicture}[x=1cm,y=32pt]
      \node at (1,0) (blank0) {};
      \node at (1,-1) (dM) {$\delta^M$};
      \node at (1,-2) (h) {$h^1$};
      \node at (1,-3) (dN) {$\delta^N$};
      \node at (0,-4) (mu) {$\mu$};
      \node at (1,-5) (blank1) {};
      \node at (0,-5) (blank2) {};
      \draw[modarrow] (blank0) to (dM);
      \draw[modarrow] (dM) to (h);
      \draw[modarrow] (h) to (dN);
      \draw[modarrow] (dN) to (blank1);
      \draw[algarrow] (mu) to (blank2);
      \draw[tensoralgarrow, bend right=15] (dM) to (mu);
      \draw[tensoralgarrow, bend right=15] (dN) to (mu);
      \draw[algarrow, bend right=15] (h) to (mu);
    \end{tikzpicture}}
  \ =\ 
  \mathcenter{
    \begin{tikzpicture}[x=1cm,y=32pt]
      \node at (1,0) (blank0) {};
      \node at (1,-1) (dN) {$\delta^N$};
      \node at (1,-2) (h) {$h^1$};
      \node at (1,-3) (dM) {$\delta^M$};
      \node at (2,-4) (mu) {$\mu$};
      \node at (1,-5) (blank1) {};
      \node at (2,-5) (blank2) {};
      \draw[modarrow] (blank1) to (dM);
      \draw[modarrow] (dM) to (h);
      \draw[modarrow] (h) to (dN);
      \draw[modarrow] (dN) to (blank0);
      \draw[algarrow] (mu) to (blank2);
      \draw[tensoralgarrow, bend left=15] (dM) to (mu);
      \draw[tensoralgarrow, bend left=15] (dN) to (mu);
      \draw[algarrow, bend left=15] (h) to (mu);
    \end{tikzpicture}}
  \ =\ 
  \mathcenter{
    \begin{tikzpicture}[x=1cm,y=32pt]
      \node at (1,0) (blank0) {};
      \node at (1,-1) (dN) {$\overline{\delta}{}^N$};
      \node at (1,-2) (h) {$h^1$};
      \node at (1,-3) (dM) {$\overline{\delta}{}^M$};
      \node at (2,-4) (mu) {$\mu$};
      \node at (1,-5) (blank1) {};
      \node at (2,-5) (blank2) {};
      \draw[modarrow] (blank0) to (dN);
      \draw[modarrow] (dN) to (h);
      \draw[modarrow] (h) to (dM);
      \draw[modarrow] (dM) to (blank1);
      \draw[algarrow] (mu) to (blank2);
      \draw[tensoralgarrow, bend left=15] (dM) to (mu);
      \draw[tensoralgarrow, bend left=15] (dN) to (mu);
      \draw[algarrow, bend left=15] (h) to (mu);
    \end{tikzpicture}}.
  \]
  In the last diagram, the dotted arrow represents the dual space
  $M^*$, for which the maps travel in the opposite direction.
\end{proof}

Lemma~\ref{lem:BoundedSubcategory} of course says that
$\DbModCat{\Alg}$ is a full subcategory of $\DuModCat{\Alg}$.  We will
find it most convenient to work in a category which is in-between the
two:

\begin{definition}
  \label{def:DModCat}
  Let $\DModCat{\Alg}\subset \DuModCat{\Alg}$ denote the
  full subcategory whose objects are type $D$ structures $N$
  which are homotopy equivalent to bounded type $D$ structures.
\end{definition}

See also Proposition~\ref{prop:CharacterizeDMod} for an alternate
characterization, in terms of bar resolutions.  Our reason for
preferring the above category is that it is quasi-equivalent to the
category of type $A$ modules;
see~Proposition~\ref{prop:D-to-A-and-back}.

A type $D$ structure $\lsup{\Alg}M$ can be seen as a way to generalize
the notion of ``projective module'' to modules over an
$\Ainf$-algebra. In particular, we will see in
Section~\ref{sec:tensor-products} how to make $\Alg
\otimes_\Ground M$ into an $\Ainf$-module over $\Alg$, which, when
$\Alg$ is a differential graded algebra and $\lsup{\Alg}M$ is bounded, is
a projective module (Corollary~\ref{cor:tensor-projective}).

\begin{remark}\label{rmk:bounded-not-iso-invt} The notion of
  boundedness is not invariant under isomorphisms of type $D$
  structures. For example, consider the algebra $\Alg=\Alg(\Torus,0)$
  discussed in Section~\ref{subsec:GenusOneAlgebra} and the following
  type $D$ structures over $\Alg$:
  \begin{itemize}
  \item $\lsup{\Alg}M$ given by $M=\Field\langle a,b\rangle$ with 
      $\delta^1(a)=\rho_{12}\otimes a+(\rho_{1}+\rho_3)\otimes b$ and
      $\delta^1(b)=\rho_{23}\otimes b.$
  \item  $\lsup{\Alg}N$ given by $M=\Field\langle c,d\rangle$ with 
      $\delta^1(c)=(\rho_{1}+\rho_3)\otimes d$ and  $\delta^1(d)=0.$
  \end{itemize}
  (These are both models for $\CFDa$ of the $-1$-framed solid torus;
  compare~\cite[Section~\ref*{LOT:sec:surg-exact-triangle}]{LOT1}.)
  Clearly, $\lsup{\Alg}N$ is bounded but $\lsup{\Alg}M$ is not. The
  map $f\co\lsup{\Alg}N\to\lsup{\Alg}M$ given by
  $f^1(c)=a$ and $f^1(d)=b+\rho_2\otimes a$
  is an isomorphism.
\end{remark}

\begin{remark}\label{rmk:A-is-D-over-cobar}
  Let $\Cobarop(\Alg)=T^*(A_+[1]^*)$ be the tensor algebra on the dual
  of $A_+$ (the ``cobar resolution'' of $\Alg$) and $\CCobarop(\Alg)$
  its completion with respect to the
  length filtration. Then the operation $\overline{D}{}^\Alg$ dualizes
  to endow $\Cobarop(\Alg)$ and $\CCobarop(\Alg)$ with a
  differential. One can show that an $\Ainf$-module over $\Alg$ is
  exactly a type $D$ structure over $\CCobarop(\Alg)$.
\end{remark}

\begin{remark}
  \label{rmk:Comodule}
  In a similar vein to Remark~\ref{rmk:A-is-D-over-cobar}, we have the
  following reformulation of bounded type $D$ structures in terms of
  differential comodules.
    Consider the bar resolution
    ${\mathrm{Bar}}(\Alg)={\Tensor}^*(\Alg_+[1])$, endowed with its
    differential, which sums over all the
    ways of grouping together $k$ consecutive
    elements and applying $\mu_k$ to them.  This can be thought of
    as an associative coalgebra equipped with the comultiplication
    $\Delta$, which sums over all the ways of splitting up an element
    $a_1\otimes\dots\otimes a_n$ as a tensor product of two elements
    of ${\mathrm{Bar}}(\Alg)$, $(a_1\otimes\dots\otimes a_i)\otimes
    (a_{i+1}\otimes\dots\otimes a_n)$. Bounded type $D$ structures over $\Alg$
    are precisely differential comodules over
    ${\mathrm{Bar}}(\Alg)$. Specifically, fix a differential comodule 
    $$\Delta\co N \to {\mathrm{Bar}}(\Alg)\otimes N,$$
    equipped with compatible differential $D\co N
    \to N$.  Letting $\Pi^i$ denote the natural projection map
    $\Pi^i\co \Tensor^*(\Alg_+[1])\to
    \Alg_+[1]^{\otimes i}$, we can define the associated type $D$
    structure by $\delta^1=D+(\Pi^1\otimes \Id_N)\circ \Delta$. Conversely, given a type $D$
    structure $\delta^1$, the comodule structure $\Delta$ is given by 
    $\Delta=[\Tensor^*(\Id_A-\epsilon)\otimes \Id_M]\circ \delta$, where
    $\delta$ is given in Equation~\eqref{eq:DefineDelta}, and
    the differential is given by $(\epsilon\otimes \Id_M)\circ \delta^1$.
    Indeed, it is interesting to compare $\DT$ with the twisted tensor product
    of~\cite{LefevreAInfinity}.  (See also~\cite{KellerLefevre}.)
\end{remark}

\begin{remark}
  \label{rmk:TwistedComplex}
  Given the \dg algebra (respectively $\Ainf$
  algebra) $\Alg$, we can form the \dg category (respectively
  $\Ainf$-category) $C_0(\Alg)$ whose objects correspond to elementary
  idempotents of
  $\Alg$. Given idempotents $I$ and $J$, let $\Mor(I,J)$ be
  the algebra elements with $J\cdot a\cdot I=a$ (endowed with the
  natural differential), and let
  composition correspond to multiplication in the algebra
  (respectively higher composition correspond to higher multiplications in the
  algebra).  A type $D$ structure over $\Alg$ can be thought of as a
  twisted complex over the additive closure of $C_0(\Alg)$. (For the
  definition of twisted
  complexes, see \cite{BondalKapranov,Kontsevich}.)
\end{remark}

\subsubsection{Bimodules of various types}
\label{sec:bimod-var-types}
Just as there are two notions of module over an
$\Ainf$-algebra---those of an $\Ainf$-module (or type $A$ module) and
a type $D$ structure---there are four notions of bimodule: type \DD,
\AAm, \DA\ and \AD. Actually, there are subtleties about
defining type \DD\ modules over general $\Ainf$-algebras, so for these
we restrict to \dg algebras.

\begin{definition}
  Let $\Alg$ and $\Blg$ be  (strictly unital, augmented) $\Ainf$-algebras over $\Groundk$ and
  $\Groundl$ respectively. Then an
  \emph{$\Ainf$-bimodule} or \emph{type
    \AAm\ bimodule $\lsub{\Alg}M_\Blg$ over $\Alg$ and $\Blg$}
  consists of a graded $(\Groundk,\Groundl)$-bimodule $M$ and degree $0$ maps
  \[
   m_{i,1,j}\co A[1]^{\otimes i}\otimes M\otimes B[1]^{\otimes j}\to M[1]
  \]
  such that, for $m=\sum_{i,j}m_{i,1,j}$, the following
  analogue of Equation~(\ref{eq:inf-mod-rel}) is satisfied:
  \begin{equation}
    \label{eq:AAbimodule-relation}
  \mathcenter{
    \begin{tikzpicture}
      \node at (0,0) (topblank0) {};
      \node at (0,-2) (m1) {$m$};
      \node at (0,-3) (m2) {$m$};
      \node at (0,-4) (botblank) {};
      \node at (-1,0) (topblank1) {};
      \node at (1,0) (topblank2) {};
      \node at (-1,-1) (leftdelta) {$\Delta$};
      \node at (1,-1) (rightdelta) {$\Delta$};
      \draw[modarrow] (topblank0) to (m1);
      \draw[tensoralgarrow] (topblank1) to (leftdelta);
      \draw[tensoralgarrow] (topblank2) to (rightdelta);
      \draw[tensoralgarrow] (leftdelta) to (m1);
      \draw[tensoralgarrow] (rightdelta) to (m1);
      \draw[tensoralgarrow] (leftdelta) to (m2);
      \draw[tensoralgarrow] (rightdelta) to (m2);
      \draw[modarrow] (m1) to (m2);
      \draw[modarrow] (m2) to (botblank);
    \end{tikzpicture}}
  +
  \mathcenter{
    \begin{tikzpicture}
      \node at (-1, 0) (tlblank) {};
      \node at (0,0) (tcblank) {};
      \node at (1,0) (trblank) {};
      \node at (-1,-1) (DAlg) {$\overline{D}{}^{\mathrlap{\Alg}}$};
      \node at (0,-2) (m) {$m$};
      \node at (0,-3) (bblank) {};
      \draw[tensoralgarrow] (tlblank) to (DAlg);
      \draw[tensoralgarrow] (trblank) to (m);
      \draw[tensoralgarrow] (DAlg) to (m);
      \draw[modarrow] (tcblank) to (m);
      \draw[modarrow] (m) to (bblank);
    \end{tikzpicture}}
  +
  \mathcenter{
    \begin{tikzpicture}
      \node at (-1, 0) (tlblank) {};
      \node at (0,0) (tcblank) {};
      \node at (1,0) (trblank) {};
      \node at (1,-1) (DBlg) {$\overline{D}{}^{\mathrlap{\Blg}}$};
      \node at (0,-2) (m) {$m$};
      \node at (0,-3) (bblank) {};
      \draw[tensoralgarrow] (trblank) to (DBlg);
      \draw[tensoralgarrow] (tlblank) to (m);
      \draw[tensoralgarrow] (DBlg) to (m);
      \draw[modarrow] (tcblank) to (m);
      \draw[modarrow] (m) to (bblank);      
    \end{tikzpicture}}
  =0.
  \end{equation}

  An $\Ainf$-bimodule is \emph{strictly unital} if
  $m_{1,1,0}(1,x)=x=m_{0,1,1}(x,1)$ for any $x$ and if
  $m_{i,1,j}(a_1,\dots,a_i,x,b_1,\dots,b_j)$ vanishes if $i+j>1$ and
  one of the $a_k$ or $b_\ell$ lies in $\Groundk$ or $\Groundl$. We shall
  always work with strictly unital $\Ainf$-bimodules.

  A morphism $\lsub{\Alg}M_\Blg\to \lsub{\Alg}N_\Blg$ is a collection
  of maps $f_{i,1,j}\co A_+[1]^{\otimes i}\otimes M\otimes B_+[1]^{\otimes
    j}\to N$.  The set of morphisms is naturally a graded
  $\Field$ vector space; compare
  Definition~\ref{def:Ainf-mod-homs}. Moreover, the set of morphisms
  forms a chain complex: writing $f$ to denote the total map $f\co
  T^*(A_+[1])\otimes M\otimes T^*(B_+[1])\to N$, the differential of such a
  morphism $f$ is
\[
\bdy f=
\mathcenter{
    \begin{tikzpicture}[x=1cm,y=32pt]
      \node at (0,0) (topblank0) {};
      \node at (0,-2) (m1) {$f$};
      \node at (0,-3) (m2) {$m$};
      \node at (0,-4) (botblank) {};
      \node at (-1,0) (topblank1) {};
      \node at (1,0) (topblank2) {};
      \node at (-1,-1) (leftdelta) {$\Delta$};
      \node at (1,-1) (rightdelta) {$\Delta$};
      \draw[modarrow] (topblank0) to (m1);
      \draw[tensoralgarrow] (topblank1) to (leftdelta);
      \draw[tensoralgarrow] (topblank2) to (rightdelta);
      \draw[tensoralgarrow] (leftdelta) to (m1);
      \draw[tensoralgarrow] (rightdelta) to (m1);
      \draw[tensoralgarrow] (leftdelta) to (m2);
      \draw[tensoralgarrow] (rightdelta) to (m2);
      \draw[modarrow] (m1) to (m2);
      \draw[modarrow] (m2) to (botblank);
    \end{tikzpicture}}
+
\mathcenter{
    \begin{tikzpicture}[x=1cm,y=32pt]
      \node at (0,0) (topblank0) {};
      \node at (0,-2) (m1) {$m$};
      \node at (0,-3) (m2) {$f$};
      \node at (0,-4) (botblank) {};
      \node at (-1,0) (topblank1) {};
      \node at (1,0) (topblank2) {};
      \node at (-1,-1) (leftdelta) {$\Delta$};
      \node at (1,-1) (rightdelta) {$\Delta$};
      \draw[modarrow] (topblank0) to (m1);
      \draw[tensoralgarrow] (topblank1) to (leftdelta);
      \draw[tensoralgarrow] (topblank2) to (rightdelta);
      \draw[tensoralgarrow] (leftdelta) to (m1);
      \draw[tensoralgarrow] (rightdelta) to (m1);
      \draw[tensoralgarrow] (leftdelta) to (m2);
      \draw[tensoralgarrow] (rightdelta) to (m2);
      \draw[modarrow] (m1) to (m2);
      \draw[modarrow] (m2) to (botblank);
    \end{tikzpicture}}
+
\mathcenter{
    \begin{tikzpicture}
      \node at (-1, 0) (tlblank) {};
      \node at (0,0) (tcblank) {};
      \node at (1,0) (trblank) {};
      \node at (-1,-1) (DAlg) {$\overline{D}{}^{\mathrlap{\Alg}}$};
      \node at (0,-2) (m) {$f$};
      \node at (0,-3) (bblank) {};
      \draw[tensoralgarrow] (tlblank) to (DAlg);
      \draw[tensoralgarrow] (trblank) to (m);
      \draw[tensoralgarrow] (DAlg) to (m);
      \draw[modarrow] (tcblank) to (m);
      \draw[modarrow] (m) to (bblank);
    \end{tikzpicture}}
  +
  \mathcenter{
    \begin{tikzpicture}
      \node at (-1, 0) (tlblank) {};
      \node at (0,0) (tcblank) {};
      \node at (1,0) (trblank) {};
      \node at (1,-1) (DBlg) {$\overline{D}{}^{\mathrlap{\Blg}}$};
      \node at (0,-2) (m) {$f$};
      \node at (0,-3) (bblank) {};
      \draw[tensoralgarrow] (trblank) to (DBlg);
      \draw[tensoralgarrow] (tlblank) to (m);
      \draw[tensoralgarrow] (DBlg) to (m);
      \draw[modarrow] (tcblank) to (m);
      \draw[modarrow] (m) to (bblank);      
    \end{tikzpicture}}
  .
  \]
  Given another morphism $g\co \lsub{\Alg}N_\Blg\to \lsub{\Alg}P_\Blg$
  define
  \[
  g\circ f=
  \mathcenter{
    \begin{tikzpicture}[x=1cm,y=32pt]
      \node at (0,0) (topblank0) {};
      \node at (0,-2) (m1) {$f$};
      \node at (0,-3) (m2) {$g$};
      \node at (0,-4) (botblank) {};
      \node at (-1,0) (topblank1) {};
      \node at (1,0) (topblank2) {};
      \node at (-1,-1) (leftdelta) {$\Delta$};
      \node at (1,-1) (rightdelta) {$\Delta$};
      \draw[modarrow] (topblank0) to (m1);
      \draw[tensoralgarrow] (topblank1) to (leftdelta);
      \draw[tensoralgarrow] (topblank2) to (rightdelta);
      \draw[tensoralgarrow] (leftdelta) to (m1);
      \draw[tensoralgarrow] (rightdelta) to (m1);
      \draw[tensoralgarrow] (leftdelta) to (m2);
      \draw[tensoralgarrow] (rightdelta) to (m2);
      \draw[modarrow] (m1) to (m2);
      \draw[modarrow] (m2) to (botblank);
    \end{tikzpicture}}.
  \]

  We let $\lsub{\Alg}\ModCat_\Blg$ denote the \dg category of type
  \AAm\  modules over $\Alg$ and $\Blg$.  The cycles in
  $\Mor(\lsub{\Alg}M_\Blg,\lsub{\Alg}N_\Blg)$ are the \emph{type \AAm\ 
    bimodule homomorphisms}.
\end{definition}

\begin{example}
  If $\Alg$ and $\Blg$ are \dg algebras and $M$ is a type \AAm\ 
  bimodule such that $m_{i,1,j}$ vanishes whenever $i+j>1$ then $M$ is
  an ordinary \dg bimodule.
\end{example}

An \emph{\AAm\ module operation tree} is like an $\Ainf$-module
operation tree, except that the algebra edges are now labelled either
`$A$' or `$B$', and
nodes labelled `$m$' the incoming module edge
need not be the leftmost edge, but edges to the left of the module
edge are labelled `$A$' and edges to the right are labelled `$B$'.
Similarly for \AAm\
morphism trees.

\begin{definition}\label{def:AA-bounded}
  An $\Ainf$-bimodule $\lsub{\Alg}M_\Blg$ is \emph{[operationally] bounded} if for
  each $x \in M$ there is an $n$ so that for $i+j>n$ and any spinal \AAm\
  module operation tree~$\Gamma$ with $i+1+j$
  total inputs, $i$ labelled `$A$' (to the left of the module edge)
  and $j$ labelled `$B$' (to the right of the module edge), $m_\Gamma(\cdot \otimes x \otimes \cdot)$ vanishes on $(A_+[1])^{\otimes i} \otimes
  (B_+[1])^{\otimes j}$.  It is \emph{left} (respectively
  \emph{right}) \emph{[operationally] bounded} if for each $x\in M$
  and each $i$ there
  exists an $n$ so that for all spinal \AAm\ module operation trees~$\Gamma$
  with $i$ right (respectively left) inputs and $j
  > n$ left (respectively right) inputs, $m_\Gamma(\cdot \otimes x \otimes
  \cdot)$ vanishes on $(A_+[1])^{\otimes j} \otimes (B_+[1])^{\otimes i}$
  (respectively $(A_+[1])^{\otimes i} \otimes (B_+[1])^{\otimes j}$).

  Similarly, a morphism $f\co
  \lsub{\Alg}M_\Blg\to \lsub{\Alg}N_\Blg$ is called \emph{bounded}
  if for each $x\in M$ there is an $n$ so that for any spinal \AAm\
  morphism tree~$\Gamma$ with $i+1+j>n$ total inputs, $f_\Gamma(\cdot \otimes x \otimes \cdot)$ vanishes on $(A_+[1])^{\otimes i} \otimes
  (B_+[1])^{\otimes j}$.  It is \emph{left} (respectively
  \emph{right}) \emph{bounded} if for each $x\in M$
  and each $i$ there
  exists an $n$ so that for any spinal \AAm\ morphism tree~$\Gamma$ with
  $i$ right (respectively left) inputs and $j
  > n$ left (respectively right) inputs, $f_\Gamma(\cdot \otimes x \otimes
  \cdot)$ vanishes on $(A_+[1])^{\otimes j} \otimes (B_+[1])^{\otimes i}$
  (respectively $(A_+[1])^{\otimes i} \otimes (B_+[1])^{\otimes j}$).
\end{definition}
(See also Lemma~\ref{lem:bounded-char}, for a mild reformulation of
these notions, and unification with the notion of boundedness for type
\DA\ and \DD\ structures, defined below.)

\begin{example}\label{eg:Id-AA-mod}
  An $\Ainf$ algebra $\Alg$ can be viewed as an $\Ainf$ bimodule
  $\lsub{\Alg}\Alg_{\Alg}$ over itself, with
$$m_{i,1,j}(a_1,\dots,a_m,x,b_1,\dots,b_n)=\mu_{i+j+1}(a_1,\dots,a_i,x,b_1,\dots,b_j).$$
  This bimodule is operationally bounded if $\Alg$ is nilpotent (as we
  are assuming throughout), but not usually otherwise.
\end{example}

\begin{definition}\label{def:DA-structure}
  Let $\Alg$ and $\Blg$ be $\Ainf$-algebras over $\Groundk$ and
  $\Groundl$ respectively.  Then a \emph{type \DA\ bimodule
    $\lsup{\Alg}N_\Blg$ over $\Alg$ and $\Blg$} consists of a
  graded $(\Groundk,\Groundl)$-bimodule $N$ and degree $0$,
  $(\Groundk,\Groundl)$-linear maps
  $$
  \delta_{1+j}^1\co N \otimes \DGB[1]^{\otimes j} \to \DGA[1]
  \otimes N.
  $$
  The compatibility condition is as follows. Let $\delta^1=\sum_j
  \delta^1_j$.
  Define maps $\delta^i\co N\otimes T^*(\DGB[1])\to \DGA[1]^{\otimes i}\otimes
  N$ inductively by
  \begin{align*}
    \delta^0&= \Id_N  \\
    \delta^{i+1}&=(\Id_{\DGA^{\otimes i}}\otimes\delta^1)\circ(\delta^{i}\otimes\Id_{T^*\DGB})\circ(\Id_N\otimes\Delta)
  \end{align*}
  where $\Delta\co T^*(\DGB)\to T^*(\DGB)\otimes T^*(\DGB)$ is the
  canonical comultiplication.  Let
  \[\delta^N\co N\otimes \Tensor^*(B)\to \overline{\Tensor^*}(A[1])\otimes N\]
  be the map defined by
  \[\delta^N=\sum_{i=0}^\infty\delta^i.\] That is, graphically,
  \[
  \delta^N=
  \mathcenter{\begin{tikzpicture}
      \node at (0,0) (tcblank) {};
      \node at (1,0) (trblank) {};
      \node at (0,-1) (delta) {$\delta^N$};
      \node at (-1,-2) (blblank) {};
      \node at (0,-2) (bcblank) {};
      \draw[modarrow] (tcblank) to (delta);
      \draw[tensoralgarrow] (trblank) to (delta);
      \draw[tensoralgarrow] (delta) to (blblank);
      \draw[modarrow] (delta) to (bcblank);
    \end{tikzpicture}}
  \coloneqq
  \,\,\mathcenter{\begin{tikzpicture}
      \node at (0,0) (tcblank) {\hskip 1cm};
      \node at (0,-1) (bcblank) {~};
      \draw[modarrow] (tcblank) to (bcblank);
    \end{tikzpicture}}\,\,
  \oplus
  \mathcenter{\begin{tikzpicture}
      \node at (0,0) (tcblank) {};
      \node at (1,0) (trblank) {};
      \node at (0,-1) (delta) {$\delta^1$};
      \node at (-1,-2) (blblank) {};
      \node at (0,-2) (bcblank) {};
      \draw[modarrow] (tcblank) to (delta);
      \draw[tensoralgarrow] (trblank) to (delta);
      \draw[algarrow] (delta) to (blblank);
      \draw[modarrow] (delta) to (bcblank);
    \end{tikzpicture}}
  \oplus
  \mathcenter{
    \begin{tikzpicture}
      \node at (0,0) (tcblank) {};
      \node at (1,0) (trblank) {};
      \node at (1,-1) (Delt) {$\Delta$};
      \node at (0,-2) (deltai) {$\delta^1$};
      \node at (0,-3) (deltaj) {$\delta^1$};
      \node at (0,-4) (bcblank) {};
      \node at (-1,-4) (br2blank) {};
      \node at (-.8,-4) (br1blank) {};
      \draw[modarrow] (tcblank) to (deltai);
      \draw[modarrow] (deltai) to (deltaj);
      \draw[modarrow] (deltaj) to (bcblank);
      \draw[algarrow] (deltai) to (br2blank);
      \draw[algarrow] (deltaj) to (br1blank);
      \draw[tensoralgarrow] (trblank) to (Delt);
      \draw[tensoralgarrow] (Delt) to (deltai);
      \draw[tensoralgarrow] (Delt) to (deltaj);
    \end{tikzpicture}}
  \oplus
  \mathcenter{
    \begin{tikzpicture}
      \node at (0,0) (tcblank) {};
      \node at (1,0) (trblank) {};
      \node at (1,-1) (Delt) {$\Delta$};
      \node at (0,-2) (deltai) {$\delta^1$};
      \node at (0,-3) (deltaj) {$\delta^1$};
      \node at (0,-4) (deltak) {$\delta^1$};
      \node at (0,-5) (bcblank) {};
      \node at (-1,-5) (br2blank) {};
      \node at (-.8,-5) (br1blank) {};
      \node at (-.6,-5) (br0blank) {};
      \draw[modarrow] (tcblank) to (deltai);
      \draw[modarrow] (deltai) to (deltaj);
      \draw[modarrow] (deltaj) to (deltak);
      \draw[modarrow] (deltak) to (bcblank);
      \draw[algarrow] (deltai) to (br2blank);
      \draw[algarrow] (deltaj) to (br1blank);
      \draw[algarrow] (deltak) to (br0blank);
      \draw[tensoralgarrow] (trblank) to (Delt);
      \draw[tensoralgarrow] (Delt) to (deltai);
      \draw[tensoralgarrow] (Delt) to (deltaj);
      \draw[tensoralgarrow] (Delt) to (deltak);
    \end{tikzpicture}}
  \cdots .
  \]
  Then, the compatibility condition is given graphically by
  \[
  \mathcenter{
    \begin{tikzpicture}
      \node at (0,0) (tcblank) {};
      \node at (1,0) (trblank) {};
      \node at (1,-1) (DBlg) {$\overline{D}^{\mathrlap{\Blg}}$};
      \node at (0,-2) (delta) {$\delta^N$};
      \node at (0,-3) (bcblank) {};
      \node at (-1,-3) (blblank) {};
      \draw[modarrow] (tcblank) to (delta);
      \draw[modarrow] (delta) to (bcblank);
      \draw[tensoralgarrow] (delta) to (blblank);
      \draw[tensoralgarrow] (trblank) to (DBlg);
      \draw[tensoralgarrow] (DBlg) to (delta);
    \end{tikzpicture}
  }+
  \mathcenter{
    \begin{tikzpicture}
      \node at (0,0) (tcblank) {};
      \node at (1,0) (trblank) {};
      \node at (-1,-2) (DAlg) {$\overline{D}^{\mathrlap{\Alg}}$};
      \node at (0,-1) (delta) {$\delta^N$};
      \node at (0,-3) (bcblank) {};
      \node at (-1,-3) (blblank) {};
      \draw[modarrow] (tcblank) to (delta);
      \draw[modarrow] (delta) to (bcblank);
      \draw[tensoralgarrow] (trblank) to (delta);
      \draw[tensoralgarrow] (delta) to (DAlg);
      \draw[tensoralgarrow] (DAlg) to (blblank);
    \end{tikzpicture}
  }
  =0,
  \]
  or symbolically by
  \begin{equation}\label{eq:DA-def}
    \delta^N \circ (\Id_N \otimes \overline{D}^\Blg) + (\overline{D}^\Alg \otimes \Id_N) \circ \delta^N = 0.
  \end{equation}
  (Compare Formula~(\ref{eq:typeD-def}).) 

  A type \DA\ structure $\lsup{\Alg}M_\Blg$ is called \emph{strictly
    unital} if $\delta^1_2(x,1)=1\otimes x$ for any $x\in M$ and
  $\delta^1_{1+i}(x,b_1,\dots,b_i)=0$ if $i>1$ and some $b_\ell\in\Groundl$,
  so $\delta^1_{1+i}$ is induced by a map from $N\otimes
  B_+[1]^{\otimes i}$ to  $A[1]\otimes N$, which we also denote
  $\delta^1_{1+i}$. We will assume our type \DA\ 
  structures are strictly unital.

  A morphism of type \DA\ structures $f^1\co \lsup{\Alg}M_\Blg\to
  \lsup{\Alg}N_\Blg$ is a collection of maps $f^1_{1+j}\co M\otimes
  B_+[1]^{\otimes j}\to A \otimes N$.  The set of morphisms is
  naturally a graded
  vector space; compare
  Definition~\ref{def:Ainf-mod-homs}. Moreover, the set of morphisms
  forms a chain complex: the differential of a morphism $f^1$ is shown
  in Figure~\ref{fig:DA-diff-compose} (left).
  \begin{figure}
    \centering
    \[
    \bdy(f^1)= \!\mathcenter{\begin{tikzpicture} 
        \node at (0,0) (tcblank) {}; 
        \node at (1.5,0) (trblank) {}; 
        \node at (1,-1) (Delta) {$\Delta$}; 
        \node at (0,-2) (deltaM) {$\delta^M$};
        \node at (0,-3) (f) {$f^1$}; 
        \node at (0,-4) (deltaN) {$\delta^N$};
        \node at (-1,-5) (mu) {$\mu^\Alg$};
        \node at (0,-6) (bcblank) {}; 
        \node at (-1.5,-6) (blblank) {};
        \draw[DAmodar] (tcblank) to (deltaM); 
        \draw[DAmodar] (deltaM) to (f);
        \draw[DAmodar] (f) to (deltaN); 
        \draw[DAmodar] (deltaN) to (bcblank);
        \draw[tensorblgarrow, bend left=5] (trblank) to (Delta);
        \draw[tensorblgarrow] (Delta) to (deltaM);
        \draw[tensorblgarrow, bend left=5] (Delta) to (f); 
        \draw[tensorblgarrow, bend left=15] (Delta) to (deltaN);
        \draw[tensoralgarrow, bend right=15] (deltaM) to (mu);
        \draw[algarrow, bend right=5] (f) to (mu); 
        \draw[tensoralgarrow] (deltaN) to (mu);
        \draw[algarrow, bend right=5] (mu) to (blblank);
      \end{tikzpicture}}\!
    + \!\mathcenter{
      \begin{tikzpicture}
        \node at (0,0) (tcblank) {}; 
        \node at (2,0) (trblank) {};
        \node at (1,-1) (D) {$\overline{D}^{\mathrlap{\Blg}}$}; 
        \node at (0,-2) (f) {$f^1$}; 
        \node at (0,-3) (bcblank) {}; 
        \node at (-1,-3) (blblank) {}; 
        \draw[DAmodar] (tcblank) to (f); 
        \draw[DAmodar] (f) to (bcblank); 
        \draw[tensorblgarrow] (trblank) to (D);
        \draw[tensorblgarrow] (D) to (f); 
        \draw[algarrow] (f) to (blblank);
      \end{tikzpicture}
    }
    \qquad\qquad
     \mathcenter{\begin{tikzpicture} 
        \node at (0,0) (tcblank) {}; 
        \node at (2.5,0) (trblank) {}; 
        \node at (2,-1) (Delta) {$\Delta$}; 
        \node at (0,-2) (deltaM1) {$\delta^{M_{k+1}}$};
        \node at (0,-3) (f1) {$f^1_{k}$}; 
        \node at (0,-4) (deltaM2) {$\delta^{M_k}$};
        \node at (0,-5) (f2) {$f^1_{k-1}$};
        \node at (0,-6) (vdots) {$\vdots$};
        \node at (0,-7) (deltaM3) {$\delta^{M_2}$};
        \node at (0,-8) (f3) {$f^1_1$};
        \node at (0,-9) (deltaM4) {$\delta^{M_1}$};
        \node at (-2,-10) (mu) {$\mu^\Alg$};
        \node at (0,-11) (bcblank) {}; 
        \node at (-2.5,-11) (blblank) {};
        \draw[DAmodar] (tcblank) to (deltaM1); 
        \draw[DAmodar] (deltaM1) to (f1);
        \draw[DAmodar] (f1) to (deltaM2); 
        \draw[DAmodar] (deltaM2) to (f2);
        \draw[DAmodar] (f2) to (vdots);
        \draw[DAmodar] (vdots) to (deltaM3);
        \draw[DAmodar] (deltaM3) to (f3);
        \draw[DAmodar] (f3) to (deltaM4);
        \draw[DAmodar] (deltaM4) to (bcblank);
        \draw[tensorblgarrow, bend left=5] (trblank) to (Delta);
        \draw[tensorblgarrow, bend left=5] (Delta) to (deltaM1);
        \draw[tensorblgarrow, bend left=5] (Delta) to (f1); 
        \draw[tensorblgarrow, bend left=5] (Delta) to (deltaM2);
        \draw[tensorblgarrow, bend left=15] (Delta) to (f2); 
        \draw[tensorblgarrow, bend left=15] (Delta) to (deltaM3);
        \draw[tensorblgarrow, bend left=15] (Delta) to (f3); 
        \draw[tensorblgarrow, bend left=15] (Delta) to (deltaM4);
        \draw[algarrow, bend right=15] (f1) to (mu); 
        \draw[algarrow, bend right=15] (f2) to (mu); 
        \draw[algarrow, bend right=5] (f3) to (mu); 
        \draw[tensoralgarrow, bend right=15] (deltaM1) to (mu);
        \draw[tensoralgarrow, bend right=15] (deltaM2) to (mu);
        \draw[tensoralgarrow, bend right=5] (deltaM3) to (mu);
        \draw[tensoralgarrow, bend right=5] (deltaM4) to (mu);
        \draw[algarrow, bend right=5] (mu) to (blblank);
      \end{tikzpicture}}.
    \]
    \caption{\textbf{The differential and composition of \DA\ morphisms.} Left: the differential of a type \DA\ morphism, which can also be thought of as the composition map $\circ_1$. Right: the composition $\circ_k$ of morphisms $f_1,\dots,f_k$.}
\label{fig:DA-diff-compose}
\end{figure}
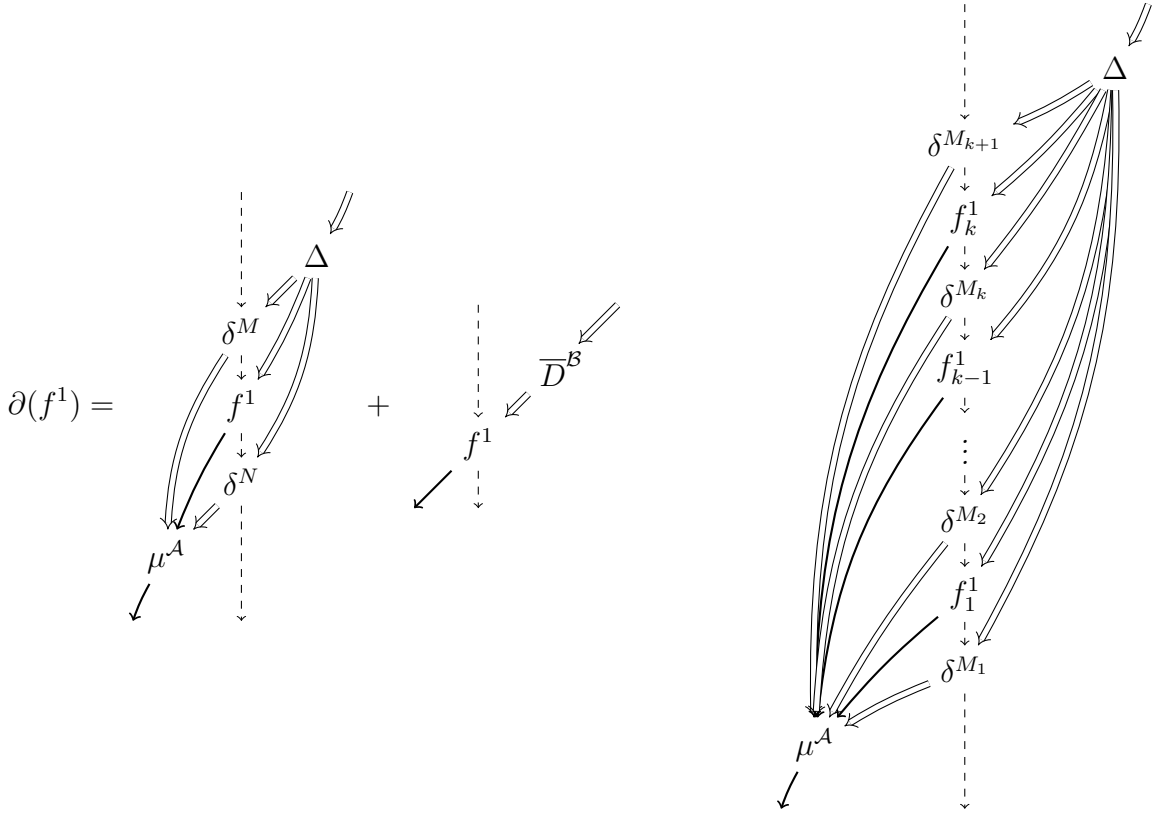
  Define
  higher composition maps $\circ_k(f_1,\dots,f_k)$ as in
  Figure~\ref{fig:DA-diff-compose} (right).
  These composition maps make the collection of type \DA\ structures
  over $\Alg$ and $\Blg$ into 
  an $\Ainf$ category, which we denote 
  $\DuModCatA{\Alg}{\Blg}$. 
  The cycles in
  $\Mor(\lsup{\Alg}M_\Blg,\lsup{\Alg}N_\Blg)$ are called \emph{type
    \DA\ structure homomorphisms.}
\end{definition}

Let 
\[
  \delta^i_j=\delta^i|_{N\otimes B_+^{\otimes (j-1)}}\qquad\text{and}\qquad
  \delta^N_j=\delta^N|_{N\otimes B_+^{\otimes (j-1)}}
\]
be the parts of $\delta^i$ and $\delta^N$ taking $j$ inputs.

\begin{definition}
  A \emph{\DA\ module operation graph} $\Gamma$ consists of:
  \begin{itemize}
  \item A connected, directed graph~$G$.
  \item An embedding $\iota\co G\to\overline{\bD}^2$ of $G$ in the
    disk, so that all edges point downwards.
  \item A labeling of the vertices of $G$ mapped by $\iota$ to the
    interior of the disk (the \emph{interior vertices}) by `$\mu$',
    `$\delta$', or `$\epsilon$'.
  \item A marking of each edge of $G$ as either a module edge
    (labelled `$M$') or an algebra edge (labelled `$A$' or `$B$').
  \end{itemize}
  This data is required to satisfy:
  \begin{itemize}
  \item The directed graph $G$ has no oriented cycles.
  \item The vertices mapped by $\iota$ to $\bdy\overline{\bD}^2$ (the
    \emph{exterior vertices}) are leaves.
  \item At each `$\mu$' vertex, there is at least one incoming and
    exactly one outgoing edge, all algebra edges with the same label.
  \item At each `$\epsilon$' vertex, there is one incoming algebra
    edge and no outgoing edges.
  \item At each `$\delta$' vertex, there is at least one incoming and
    exactly two outgoing edges,
    such that the leftmost incoming and right outgoing edge are
    module edges, and with the other incoming edges
    marked `$B$ and the other outgoing edge marked `$A$'.
  \item There is a module edge.
  \end{itemize}
\end{definition}
Call an exterior vertex an `in' (respectively `out') vertex if it is a
source (respectively sink). Note that all `in' (respectively `out')
vertices are consecutive with respect to the cyclic order on
$\bdy\overline{\bD}^2$. It also follows from the conditions that the
edge from the leftmost `in' (respectively rightmost `out') vertex is
a module edge, and the edges from all other `in' (respectively `out')
vertices are algebra edges.

Associated to a \DA\ module operation graph with $(i+1)$ `in' and $(j+1)$ `out' vertices
is a map
\[
\delta_\Gamma \co M\otimes B_+^{\otimes i}\to A^{\otimes j}\otimes M,
\]
defined in the obvious way.

A \DA\ module operation graph is \emph{spinal} if it has no `$\mu$' nodes. 

\begin{definition}\label{def:DA-bounded}
  A type \DA\ structure $\lsupv{\Alg}M_\Blg$ is called
  \emph{[operationally] bounded} if for each $x \in M$ there is an $n$
  so that for all $i+j > n$ and spinal \DA\ module operation
  graphs~$\Gamma$ with $i$ (right) algebra inputs and $j$
  (left) algebra outputs, $\delta_\Gamma(x\otimes\cdot)$ vanishes on
  $(B_+[1])^{\otimes i}$.
  It is
  \emph{right [operationally] bounded} if for each $x$ and each $j$ there is an $n$ so
  that for spinal \DA\ module operation graphs~$\Gamma$ with $i>n$
  algebra inputs and $j$ algebra outputs,
  $\delta_\Gamma(x \otimes \cdot)$ vanishes on $(B_+[1])^{\otimes i}$.
  It is \emph{left [operationally] bounded}
  if for each $x\in \lsupv{\Alg}M_{\Blg}$ and each $i$
  there is an $n$ so that for all \DA\ module operation
  graphs~$\Gamma$ with $i$ algebra inputs and $j > n$ algebra outputs,
  $\delta_\Gamma(x\otimes \cdot)$ vanishes on $(B_+[1])^{\otimes i}$.
  Boundedness for morphisms is defined similarly.

  We denote the category of bounded type \DA\ structures and bounded
  type \DA\ morphisms by 
  $\DBModCatA{\Alg}{\Blg}$.
\end{definition}

\begin{remark}
  One reason that $\epsilon$ appears in the boundedness conditions for
  type \DA\ structures is that we want the forgetful functor from type
  \DA\ structures to type $A$ modules (see
  Section~\ref{sec:forgetful-functors}) to take right bounded
  structures to bounded modules.
\end{remark}

\begin{definition}\label{def:rank-1-DA-mods}
Let $\Alg$ and $\Blg$ be $\Ainf$-algebras over $\Ground$. 
Given an $\Ainf$-morphism~$\phi\co\Blg\to\Alg$, defined by
maps $\phi_k \co B^{\otimes k} \to A$, define a bimodule
$\lsupv{\Alg}[\phi]_\Blg$ with underlying space a free rank-1 module over
$\Ground$ and structure maps given by the $\phi_k$. That is, let
$\iota\in [\phi]$ be the generator and define
\[
\delta^1_{1+k}(\iota,b_1,\dots,b_k)=\phi_k(b_1,\dots,b_k)\otimes\iota.
\]
\end{definition}

As a special case, given an $\Ainf$-algebra $\Alg$ we have the module
$\lsupv{\Alg}[\Id]_\Alg$. As a $\Ground$-module,
$\lsupv{\Alg}[\Id]_\Alg$ is isomorphic to $\Ground$. For $k\neq 2$,
$\delta^1_k=0$, while 
\[
\delta^1_2(\iota,a)=a\otimes\iota
\]
where $\iota$ is the generator of $\lsupv{\Alg}[\Id]_\Alg$. We call
$\lsupv{\Alg}[\Id]_\Alg$ the \emph{identity bimodule.}

\begin{remark}
  The identity bimodule $\lsupv{\Alg}[\Id]_\Alg$ is typically not
  [operationally] bounded: in fact, it is
  bounded if and only if $\Tensor^* A_+$ is finite dimensional.
  The module $\lsupv{\Alg}[\Id]_{\Alg}$ is always
  left and right [operationally] bounded, however.
\end{remark}

In fact, the modules from Definition~\ref{def:rank-1-DA-mods} have the following
easy characterization:
\begin{lemma}\label{lem:characterize-induced} Let $\Alg$ and $\Blg$ be
  $\Ainf$-algebras over $\Ground$.
Let $\lsup{\Alg}M_{\Blg}$ be a type \DA\ bimodule whose underlying $\Ground$-bimodule
is $\Ground$. Assume furthermore that $\delta^1_1=0$. Then there is an $\Ainf$-algebra morphism
$\phi\co\Blg\to \Alg$ with the property that
$\lsup{\Alg}M_\Blg\cong \lsupv{\Alg}[\phi]_{\Blg}$.
\end{lemma}
\begin{proof}
Given $\lsup{\Alg}M_\Blg$ with generator $\iota$ as a $\Ground$-module, the homomorphism~$\phi$ is uniquely characterized by
\[
\delta^1_{n+1}(\iota,b_1,\dots,b_n)=\phi_n(b_1,\dots,b_n)\otimes \iota.
\]
The hypothesis that $\delta^1_1=0$ ensures that $\phi_0=0$, so we drop it.
We must verify that $\phi=\{\phi_i\}_{i=1}^{\infty}$ satisfies the $\Ainf$-relation (Equation~\eqref{eq:Ainf-alg-morph}). Define
$\tau\co \Ground\otimes T^*(\Blg)\simeq T^*(\Blg)\otimes\Ground$ to be
the canonical identification (as both are isomorphic to
$T^*(\Blg)$). Then we have
$\delta=(F^\phi\otimes\Id_\Ground)\circ\tau$. So, by the
$\Ainf$-relation for $\delta$,
\begin{align*}
  0&=\delta\circ(\Id_\Ground\otimes
  \overline{D}^\Blg)+(\overline{D}^\Alg\otimes \Id_\Ground)\circ\delta\\
  &=(F^\phi\otimes\Id_\Ground)\circ\tau\circ(\Id_\Ground\otimes
  \overline{D}^\Blg)+(\overline{D}^\Alg\otimes
  \Id_\Ground)\circ(F^\phi\otimes\Id_\Ground)\circ\tau\\
  &=((F^\phi\circ \overline{D}^\Blg)\otimes\Id_\Ground)\circ\tau
+((\overline{D}^\Alg\circ F^\phi)\otimes\Id_\Ground)\circ\tau
\end{align*}
hence $F^\phi\circ \overline{D}^\Blg+\overline{D}^\Alg\circ
F^\phi=0$, as desired.
\end{proof}

\begin{remark}
  The hypothesis that $\delta^1_1=0$ can be dropped, if we allow for
  more general types of $\Ainf$-algebra morphisms, i.e., those which
  contain a term~$\phi_0$. The hypotheses of
  Lemma~\ref{lem:characterize-induced} are satisfied in the
  case we use it (Theorem~\ref{thm:Id-is-Id}).
\end{remark}

\begin{definition}
  Let $\DModCatA{\Alg}{\Blg}$ denote the full subcategory
  of $\DuModCatA{\Alg}{\Blg}$ consisting of type \DA\ bimodules
  which are homotopy equivalent to bounded type \DA\ bimodules.
\end{definition}

This is equivalent to the category of type \DA\ bimodules which
are homotopy equivalent to {\em left} bounded type \DA\ bimodules;
see Proposition~\ref{prop:CharacterizeDMod} below.

One can define \emph{type \AD\  modules} similarly, by reflecting all
of the pictures.  For instance, given $\phi\co \Blg\to\Alg$, one
can define a module $\lsub{\Blg}[\phi]^{\Alg}$ analogous to the one
from Definition~\ref{def:rank-1-DA-mods}.

Like for the $D$ structures, type \DA\ modules have opposite type \AD\
modules. We will explain this operation only under some finiteness
assumptions:
\begin{definition}\label{def:opposite-type-DA}
  Suppose that $(\lsup{\Alg}M_\Blg,\delta^1)$ is a type \DA\ structure
  and $\Alg$, $\Blg$ and $M$ are finite-dimensional.  Define the
  \emph{opposite} type~\AD\ structure
  $(\lsub{\Blg}{\overline{M}}{}^\Alg,\overline{\delta}{}^1)$ as follows.
  As a $(\Ground,\Groundl)$-bimodule, $\overline{M}$ is just
  $M^*=\Hom_{\Ground\otimes\Groundl}(M,\Ground\otimes\Groundl)\cong \Hom_{\Field}(M,\Field)$. The map $\delta^1_k$ on~$M$ is an
  element of $\Hom_\Ground(M\otimes B_+^{\otimes k},A\otimes M)\isom
  \Hom_\Ground(B_+^{\otimes k}\otimes M^*,M^*\otimes A),$ and
  $\overline{\delta}{}^1$ is obtained by viewing $\delta^1$ as lying
  in the right hand side.
  \end{definition}
\begin{lemma}
  Given a type \DA\ structure $\lsup{\Alg}M_\Blg$, the opposite type
  \AD\ structure $\lsub{\Blg}{\overline{M}}^\Alg$ satisfies the type
  \AD\ structure equation.
\end{lemma}
\begin{proof}
  We leave the verification, which is similar to the proof of
  Lemma~\ref{lem:opposite-is-D}, to the reader.
\end{proof}

Next, we turn to our final notion of bimodule, a type \DD\ 
structure. We will only define these when $\Alg$ and $\Blg$ are \dg
algebras.
\begin{definition}\label{def:DD-str-over-dgas}
  Let $\Alg$ and $\Blg$ be \dg algebras over $\Groundk$ and
  $\Groundl$. We define the \emph{category
    of type \DD\ structures over $\Alg$ and $\Blg$},
  $\DuModCatDu{\Alg}{\Blg}$, to be the category of type $D$
  structures over $\Alg\otimes_{\FF_2}\Blg^\op$, that is,
  $\DuModCat{\Alg\otimes_{\FF_2}\Blg^\op}$.
\end{definition}
We denote a type \DD\ structure $M$ by $\lsup{\Alg}M^\Blg$. The cycles
in $\Mor(\lsup{\Alg}M^\Blg,\lsup{\Alg}N^\Blg)$ are the \emph{type \DD\
  structure homomorphisms.}

We think of the data of a type \DD\ structure as a graded
$(\Groundk,\Groundl)$-bimodule
$M$ and a degree $0$ map $\delta^1\co M\to \DGA\otimes
M\otimes\DGB[1]$, such that the
following compatibility condition holds:
\[
\mathcenter{
  \begin{tikzpicture}
    \node at (0,0) (tblank) {};
    \node at (0,-1) (delta1) {$\delta^1$};
    \node at (0,-2) (delta2) {$\delta^1$};
    \node at (-1,-3) (mul) {$\mu_2$};
    \node at (1,-3) (mur) {$\mu_2$};
    \node at (-1,-4) (blblank) {};
    \node at (0,-4) (bcblank) {};
    \node at (1,-4) (brblank) {};
    \draw[modarrow] (tblank) to (delta1);
    \draw[modarrow] (delta1) to (delta2);
    \draw[modarrow] (delta2) to (bcblank);
    \draw[algarrow, bend right=15] (delta1) to (mul);
    \draw[algarrow, bend right=15] (delta2) to (mul);
    \draw[algarrow, bend left=15] (delta1) to (mur);
    \draw[algarrow, bend left=15] (delta2) to (mur);
    \draw[algarrow] (mul) to (blblank);
    \draw[algarrow] (mur) to (brblank);
  \end{tikzpicture}
}
+
\mathcenter{
  \begin{tikzpicture}
    \node at (0,0) (tblank) {};
    \node at (0,-1) (delta) {$\delta^1$};
    \node at (-1,-2) (mu) {$\mu_1$};
    \node at (-1,-3) (blblank) {};
    \node at (0,-3) (bcblank) {};
    \node at (1,-3) (brblank) {};
    \draw[modarrow] (tblank) to (delta);
    \draw[modarrow] (delta) to (bcblank);
    \draw[algarrow, bend right=15] (delta) to (mu);
    \draw[algarrow] (mu) to (blblank);
    \draw[algarrow, bend left=15] (delta) to (brblank);
  \end{tikzpicture}
}
+
\mathcenter{
  \begin{tikzpicture}
    \node at (0,0) (tblank) {};
    \node at (0,-1) (delta) {$\delta^1$};
    \node at (1,-2) (mu) {$\mu_1$};
    \node at (-1,-3) (blblank) {};
    \node at (0,-3) (bcblank) {};
    \node at (1,-3) (brblank) {};
    \draw[modarrow] (tblank) to (delta);
    \draw[modarrow] (delta) to (bcblank);
    \draw[algarrow, bend left=15] (delta) to (mu);
    \draw[algarrow] (mu) to (brblank);
    \draw[algarrow, bend right=15] (delta) to (blblank);
  \end{tikzpicture}
}
=0.
\]

A morphism $g^1\co \lsup{\Alg}M^\Blg\to \lsup{\Alg}N^\Blg$ of type
\DD\ structures is a map $g^1\co M\to A\otimes N\otimes B$.  The set
of morphisms is naturally a graded vector space, as
in Definition~\ref{def:Ainf-mod-homs}. Moreover, the set of morphisms
forms a chain complex: the
differential of a morphism $g^1$ is:
\[
\bdy(g^1)= \mathcenter{
  \begin{tikzpicture}
    \node at (0,0) (tblank) {}; 
    \node at (0,-1) (delta) {$\delta^1$}; 
    \node at (0,-2) (g) {$g^1$};
    \node at (-1,-3) (muA) {$\mu_2^\Alg$}; 
    \node at (1,-3) (muB) {$\mu_2^\Blg$};
    \node at (-1,-4) (blblank) {}; 
    \node at (0,-4) (bcblank) {};
    \node at (1,-4) (brblank) {};
    \draw[modarrow] (tblank) to (delta); 
    \draw[modarrow] (delta) to (g);
    \draw[modarrow] (g) to (bcblank);
    \draw[algarrow, bend right=15] (g) to (muA); 
    \draw[blgarrow, bend left=15] (g) to (muB); 
    \draw[algarrow, bend right=15] (delta) to (muA); 
    \draw[blgarrow, bend left=15] (delta) to (muB); 
    \draw[algarrow] (muA) to (blblank);
    \draw[blgarrow] (muB) to (brblank);
  \end{tikzpicture}
} 
+
\mathcenter{
  \begin{tikzpicture}
    \node at (0,0) (tblank) {}; 
    \node at (0,-2) (delta) {$\delta^1$}; 
    \node at (0,-1) (g) {$g^1$};
    \node at (-1,-3) (muA) {$\mu_2^\Alg$}; 
    \node at (1,-3) (muB) {$\mu_2^\Blg$};
    \node at (-1,-4) (blblank) {}; 
    \node at (0,-4) (bcblank) {};
    \node at (1,-4) (brblank) {};
    \draw[modarrow] (tblank) to (g); 
    \draw[modarrow] (g) to (delta);
    \draw[modarrow] (delta) to (bcblank);
    \draw[algarrow, bend right=15] (g) to (muA); 
    \draw[blgarrow, bend left=15] (g) to (muB); 
    \draw[algarrow, bend right=15] (delta) to (muA); 
    \draw[blgarrow, bend left=15] (delta) to (muB); 
    \draw[algarrow] (muA) to (blblank);
    \draw[blgarrow] (muB) to (brblank);
  \end{tikzpicture}
} 
+
\mathcenter{
  \begin{tikzpicture}
    \node at (0,0) (tblank) {}; 
    \node at (0,-1) (delta) {$g^1$};
    \node at (1,-2) (mu) {$\mu_1^\Blg$}; 
    \node at (-1,-3) (blblank) {}; 
    \node at (0,-3) (bcblank) {}; 
    \node at (1,-3) (brblank) {};
    \draw[modarrow] (tblank) to (delta); 
    \draw[modarrow] (delta) to (bcblank); 
    \draw[algarrow, bend right=15] (delta) to (blblank); 
    \draw[blgarrow, bend left=15] (delta) to (mu); 
    \draw[algarrow] (mu) to (brblank);
  \end{tikzpicture}
}
+ \mathcenter{
  \begin{tikzpicture}
    \node at (0,0) (tblank) {}; 
    \node at (0,-1) (delta) {$g^1$};
    \node at (-1,-2) (mu) {$\mu_1^\Blg$}; 
    \node at (1,-3) (blblank) {}; 
    \node at (0,-3) (bcblank) {}; 
    \node at (-1,-3) (brblank) {};
    \draw[modarrow] (tblank) to (delta); 
    \draw[modarrow] (delta) to (bcblank); 
    \draw[algarrow, bend left=15] (delta) to (blblank); 
    \draw[blgarrow, bend right=15] (delta) to (mu); 
    \draw[algarrow] (mu) to (brblank);
  \end{tikzpicture}
}.
\]

Composition in $\DuModCatDu{\Alg}{\Blg}$ is given as follows: for
$g^1\co \lsup{\Alg}M^\DGB\to \lsup{\Alg}N^\DGB$ and $h^1\co
\lsup{\Alg}N^\DGB\to \lsup{\Alg}P^\DGB$,
\[
h^1\circ g^1= \mathcenter{
  \begin{tikzpicture}
    \node at (0,0) (tblank) {}; 
    \node at (0,-2) (delta) {$h^1$}; 
    \node at (0,-1) (g) {$g^1$};
    \node at (-1,-3) (muA) {$\mu_2^\Alg$}; 
    \node at (1,-3) (muB) {$\mu_2^\Blg$};
    \node at (-1,-4) (blblank) {}; 
    \node at (0,-4) (bcblank) {};
    \node at (1,-4) (brblank) {};
    \draw[modarrow] (tblank) to (g); 
    \draw[modarrow] (g) to (delta);
    \draw[modarrow] (delta) to (bcblank);
    \draw[algarrow, bend right=15] (g) to (muA); 
    \draw[blgarrow, bend left=15] (g) to (muB); 
    \draw[algarrow, bend right=15] (delta) to (muA); 
    \draw[blgarrow, bend left=15] (delta) to (muB); 
    \draw[algarrow] (muA) to (blblank);
    \draw[blgarrow] (muB) to (brblank);
  \end{tikzpicture}
}.
\]

Given a type \DD\ structure $(M,\delta^1)$ define maps
$\delta^i$ inductively by
\[
\delta^i=(\Id_{\Alg^{\otimes (i-1)}}\otimes \delta^1\otimes
\Id_{\Blg^{\otimes (i-1)}})\circ \delta^{i-1}\co M\to \Alg^{\otimes
  i}\otimes M\otimes \Blg^{\otimes i}
\]
and set $\delta=\bigoplus_{i=0}^\infty\delta^i$. 

Similarly, given a morphism $g^1\co \lsup{\Alg}M^\Blg\to
\lsup{\Alg}N^\Blg$, define $g^i\co M\to \Alg^{\otimes i}\otimes
N\otimes \DGB^{\otimes i}$ by
\[
g^i=\sum_{j=1}^{i} (\Id_{\Alg^{\otimes j}}\otimes
\delta^{N,i-j}\otimes\Id_{\DGB^{\otimes j}})\circ (\Id_{\Alg^{\otimes
    (j-1)}}\otimes g^1\otimes \Id_{\DGB^{\otimes (j-1)}})\circ
\delta^{M,j-1}
\]
and let $g=\sum_{i=1}^\infty g^i$.

Note that the augmentation $\epsilon\co\DGA\to\Ground$ of $\Alg$
extends to a map $\epsilon\co T^*\DGA\to\Ground$ by
$\epsilon(a_1\otimes\cdots\otimes
a_k)=\epsilon(a_1)\cdots\epsilon(a_k)$.

\begin{definition}\label{def:DD-bounded}
  We call a type \DD\ structure $\lsup{\Alg}M^\Blg$ \emph{left}
  (respectively \emph{right}) \emph{[operationally] bounded} 
  if for each $x\in \lsup{\Alg}M^{\Blg}$, there is
  a constant $n$ with the property that for
  all $i>n$, $(\Id_{A^{\otimes i}}\otimes \Id_{M}\otimes
  \epsilon_\Blg)\circ \delta^i=0$ (respectively $(\epsilon_\Alg\otimes
  \Id_{M}\otimes\Id_{B^{\otimes i}})\circ\delta^i =0$). We call $\lsup{\Alg}M^\Blg$
  \emph{[operationally] bounded} if for each $x$, $n$ can be
  chosen so that $\delta^i(x)=0$ for
  $i>n$. Boundedness for morphisms of type \DD\ structures is defined
  similarly. 
\end{definition}

As usual, the condition of being operationally bounded is stronger
than the condition of being both left and right bounded.

\begin{definition}\label{def:DD-sep}
    Call a type \DD\ structure \emph{separated} if the map
    $\delta^1$ can be written as $\delta^{1L}+\delta^{1R}$ where
    $\delta^{1L}\co M\to \Alg\otimes M\otimes\Groundl$ and
    $\delta^{1R}\co M\to \Groundk\otimes M\otimes \Blg$. 
\end{definition}

\begin{remark}\label{rem:DD-difficult}
  In general, over $\Ainf$-algebras, a type \DD\ bimodule
  ${}^{\Alg}M^{\Blg}$ should be a type $D$ module over $\Alg
  \otimes_{\FF_2} \Blg^\op$.  The difficulty is in defining the tensor
  product of $\Ainf$-algebras. This has been
  done~\cite{SU04:Diagonals,MS06:AssociahedraProdAinf,Loday11:DiagonalStasheff},
  but is somewhat complicated and is unnecessary for this paper.
\end{remark}

\begin{definition}
  Given $\dg$ algebras $\Alg$ and $\Blg$, let
  $\DBModCatDB{\Alg}{\Blg}$ denote the category
  whose objects consist of bounded type \DD\ bimodules.
  Similarly, we define
  $\DModCatD{\Alg}{\Blg}$ to be the full subcategory of
  $\DuModCatDu{\Alg}{\Blg}$ consisting of type \DD\ bimodules
  which are homotopy equivalent to bounded ones.
\end{definition}

See also Proposition~\ref{prop:CharacterizeDMod}.

So far, we have discussed bimodules with a single left and a single
right action. One can also consider bimodules with two left actions
or two right actions---and, indeed, it is most natural to define the
invariant $\CFAAa$ (respectively $\CFDDa$) of
Section~\ref{sec:CFBimodules} as a bimodule with two right
(respectively left) actions. Obviously, there are no new
mathematical difficulties in this theory. Moreover, the notation
extends easily; for example, $M_{\Alg,\Blg}$ denotes a type \AAm\ 
bimodule with two right actions and $\ModCat_{\Alg,\Blg}$ denotes
the \dg category of such bimodules.

\subsection{Operations on bimodules}
\subsubsection{Forgetful functors}\label{sec:forgetful-functors}
In defining the tensor product and one sided $\Mor$ operations, it will be
convenient to invoke forgetful functors between certain of our
categories. In particular, there are forgetful functors
\begin{align*}
  \mathcal{F}&\co \DModCatA{\Alg}{\Blg}\to
  \DModCat{\Alg}\qquad\text{and}\\
  \mathcal{F}&\co \lsub{\Alg}\ModCat_\Blg\to 
  \lsub{\Alg}\ModCat
\end{align*}
gotten by
$\mathcal{F}(\lsup{\Alg}M_\Blg,\{\delta^1_i\})=(\lsup{\Alg}M,\overline{\delta}^1)$
where $\lsup{\Alg}M$ is isomorphic to $M$ as a $(\Groundk,\Groundl)$-bimodule and
$\overline{\delta}^1=\delta^1_1$; and similarly
$\mathcal{F}(\lsub{\Alg}M_\Blg,\{m_{i,1,j}\})=(\lsub{\Alg}M,\overline{m}_i)$
where $\overline{m}_{i+1}=m_{i,1,0}$. (The forgetful functor is
defined similarly on morphisms.)

Similarly, there are forgetful functors 
\begin{align*}
  \mathcal{F}&\co \lsupv{\Alg}\ModCat^\Blg\to
  \lsupv{\Alg}\ModCat\qquad\text{and}\\
  \mathcal{F}&\co \lsub{\Alg}\ModCat^\Blg\to 
  \lsub{\Alg}\ModCat
\end{align*}
gotten by
$\mathcal{F}(\lsup{\Alg}M^\Blg,\{\delta^1\})=(\lsup{\Alg}M,
\overline{\delta}^1)$ where $\lsup{\Alg}M$ is isomorphic to $M$ as a $(\Groundk,\Groundl)$-bimodule and
$\overline{\delta}^1=(\Id_\Alg\otimes M\otimes \epsilon_\Blg)\circ
\delta^1$; and similarly $\mathcal{F}(\lsub{\Alg}M^\Blg,\{\delta_i^1\})=(\lsub{\Alg}M,\overline{m}_i)$
where $\overline{m}_{i+1}=(M\otimes
\epsilon_\Blg)\circ \delta_i^1$. As before, $\epsilon_\Blg\co
\Blg\to\Groundl$ denotes the augmentation.

These forgetful functors interact well with our definition of
boundedness:
\begin{lemma}
  Let $M$ be a bimodule of any type. If $M$ is left bounded then
  $\mathcal{F}(M)$ is bounded.
\end{lemma}
We leave the proof to the reader; it is not hard, but involves several
cases. Note that this lemma is actually implicit in writing, say, that
$\mathcal{F}\co \DModCatA{\Alg}{\Blg}\to \DModCat{\Alg}$, since
$\DModCat{\Alg}$ consists of type $D$ structures homotopy equivalent
to bounded ones.

There are, of course, also functors which forget the left action; we
will denote these $\mathcal{F}$ as well.

\subsubsection{Tensor products}
\label{sec:tensor-products}

In the present section, we define a pairing between type $A$ modules and type $D$
modules (and their generalizations to bimodules), which gives a model for the derived tensor product of $\Ainf$ modules.
This model for the derived product comes up naturally when one studies the gluing
problems for pseudoholomorphic curves (see Section~\ref{sec:PairingTheorems}). Indeed,
this model typically has smaller rank than the usual derived tensor product
(though its differentials are correspondingly more complicated).

We start by considering modules with a single action, and then proceed
to bimodules.
\begin{definition}
  For $\Alg$ an $\Ainf$-algebra, $M_\Alg\in \ModCat_\Alg$, and
  $\lsup{\Alg}N \in \lsupv{\Alg}\ModCat$, with at least one of $M_\Alg$ or
  $\lsup{\Alg}N$ bounded, define $M_\Alg \DT \lsup{\Alg}N$ to be chain
  complex with underlying space $M \otimes_\Ground N$ and boundary operator
  \[
  \partial \coloneqq (m_M \otimes \Id_N) \circ (\Id_M \otimes \delta^N).
  \]
\end{definition}
The boundedness hypothesis implies that the tensor product is
well-defined, as follows. If $M_\Alg$ is bounded then the operations
$m_i$ vanish for sufficiently large $i$, so they sum to give a map
$M\otimes \overline{T}^*(A)\to M$. If $\lsup{\Alg}N$ is bounded then the image of
$\delta^N$ lies in $T^*A\otimes M$, to which we can apply $\sum m_i$.

Graphically, the differential on $M_\Alg\DT\lsup{\Alg}N$ is given by
\[
\mathcenter{
  \begin{tikzpicture}
    \node at (0,0) (tlblank) {};
    \node at (1,0) (trblank) {};
    \node at (1,-1) (delta) {$\delta^N$};
    \node at (0,-2) (m) {$m_M$};
    \node at (0,-3) (blblank) {};
    \node at (1,-3) (brblank) {};
    \draw[modarrow] (trblank) to (delta);
    \draw[modarrow] (tlblank) to (m);
    \draw[tensoralgarrow] (delta) to (m);
    \draw[modarrow] (delta) to (brblank);
    \draw[modarrow] (m) to (blblank);    
  \end{tikzpicture}
}.
\]

The fact that $\partial^2=0$ is verified in
\cite[Lemma~\ref*{LOT:lem:ThetaComplex}]{LOT1}.

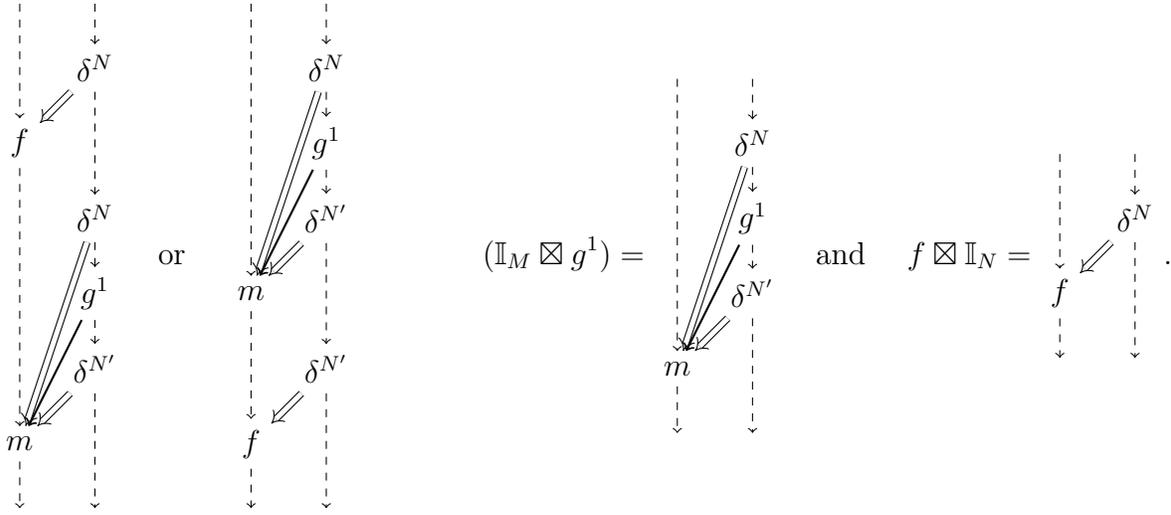
\begin{figure}
  \centering
\[
\mathcenter{
  \begin{tikzpicture}
    \node at (0,0) (tlblank) {};
    \node at (1,0) (trblank) {};
    \node at (1,-1) (delta1) {$\delta^N$};
    \node at (0,-2) (f) {$f$};
    \node at (1,-3) (delta2) {$\delta^N$};
    \node at (1,-4) (g) {$g^1$};
    \node at (1,-5) (delta3) {$\delta^{N'}$};
    \node at (0,-6) (m) {$m$};
    \node at (0,-7) (blblank) {};
    \node at (1,-7) (brblank) {};
    \draw[Amodar] (tlblank) to (f);
    \draw[Amodar] (f) to (m);
    \draw[Amodar] (m) to (blblank);
    \draw[Dmodar] (trblank) to (delta1);
    \draw[Dmodar] (delta1) to (delta2);
    \draw[Dmodar] (delta2) to (g);
    \draw[Dmodar] (g) to (delta3);
    \draw[Dmodar] (delta3) to (brblank);
    \draw[tensoralgarrow] (delta1) to (f);
    \draw[tensoralgarrow] (delta2) to (m);
    \draw[tensoralgarrow] (delta3) to (m);
    \draw[algarrow] (g) to (m);
  \end{tikzpicture}
}
\text{\quad or \quad}
\mathcenter{
  \begin{tikzpicture}
    \node at (0,0) (tlblank) {};
    \node at (1,0) (trblank) {};
    \node at (1,-1) (delta2) {$\delta^N$};
    \node at (1,-2) (g) {$g^1$};
    \node at (1,-3) (delta3) {$\delta^{N'}$};
    \node at (0,-4) (m) {$m$};
    \node at (1,-5) (delta1) {$\delta^{N'}$};
    \node at (0,-6) (f) {$f$};
    \node at (0,-7) (blblank) {};
    \node at (1,-7) (brblank) {};
    \draw[Amodar] (tlblank) to (m);
    \draw[Amodar] (m) to (f);
    \draw[Amodar] (f) to (blblank);
    \draw[Dmodar] (trblank) to (delta2);
    \draw[Dmodar] (delta3) to (delta1);
    \draw[Dmodar] (delta2) to (g);
    \draw[Dmodar] (g) to (delta3);
    \draw[Dmodar] (delta1) to (brblank);
    \draw[tensoralgarrow] (delta1) to (f);
    \draw[tensoralgarrow] (delta2) to (m);
    \draw[tensoralgarrow] (delta3) to (m);
    \draw[algarrow] (g) to (m);
  \end{tikzpicture}
}\qquad\qquad
(\Id_M\DT g^1)=
\mathcenter{
  \begin{tikzpicture}
    \node at (0,0) (tlblank) {};
    \node at (1,0) (trblank) {};
    \node at (1,-1) (delta1) {$\delta^N$};
    \node at (1,-2) (g) {$g^1$};
    \node at (1,-3) (delta2) {$\delta^{N'}$};
    \node at (0,-4) (m) {$m$};
    \node at (0,-5) (blblank) {};
    \node at (1,-5) (brblank) {};
    \draw[Amodar] (tlblank) to (m);
    \draw[Amodar] (m) to (blblank);
    \draw[Dmodar] (trblank) to (delta1);
    \draw[Dmodar] (delta1) to (g);
    \draw[Dmodar] (g) to (delta2);
    \draw[Dmodar] (delta2) to (brblank);
    \draw[tensoralgarrow] (delta1) to (m);
    \draw[algarrow] (g) to (m);
    \draw[tensoralgarrow] (delta2) to (m);
  \end{tikzpicture}
}
\text{\quad and \quad} f\DT \Id_N =
\mathcenter{
  \begin{tikzpicture}
    \node at (0,0) (tlblank) {};
    \node at (1,0) (trblank) {};
    \node at (1,-1) (delta) {$\delta^N$};
    \node at (0,-2) (f) {$f$};
    \node at (0,-3) (blblank) {};
    \node at (1,-3) (brblank) {};
    \draw[Amodar] (tlblank) to (f);
    \draw[Amodar] (f) to (blblank);
    \draw[Dmodar] (trblank) to (delta);
    \draw[Dmodar] (delta) to (brblank);
    \draw[tensoralgarrow] (delta) to (f);    
  \end{tikzpicture}
}.
\]  
  \caption{\textbf{Diagrams for defining the box product of
      morphisms.} Left: two options for defining $f\DT g$. Right: box
    products with the identity morphisms.}
  \label{fig:DT-module-maps}
\end{figure}

We will see presently that $\DT$ induces a bifunctor on the level of
derived categories. Functoriality for $\DT$ on the \dg level is
somewhat subtle, however. Given a $f: M_\Alg \to M'_\Alg$
and $g^1: \lsup{\Alg}N \to \lsup{\Alg}N'$, there are two natural
diagrams that one
might use to define $f\DT g^1$, shown in the left of 
Figure~\ref{fig:DT-module-maps}.
In other words, if we define $(\Id\DT g^1)$ and $(f\DT \Id)$ as in the
right of Figure~\ref{fig:DT-module-maps}
then the two different choices correspond to $(\Id_{M'}\DT g^1)\circ (f\DT
\Id_N)$ and $(f\DT \Id_{N'})\circ (\Id_M\DT g^1)$ respectively. (Of course, for
the diagram defining $\Id_M \DT g^1$ to make sense we either need
$M_\Alg$ or both $\lsup{\Alg}N$ and $\lsup{\Alg}N'$ to be
operationally bounded, and for the diagram defining
$f\DT \Id_N$ to be defined we need either the map $f$ or the module
$\lsup{\Alg}N$ to be operationally bounded.)

The two choices above are homotopic. Indeed:
\begin{lemma}
  \label{lem:DTfunct}
	Fix $M_\Alg, N_\Alg, L_\Alg\in\ModCat_{\Alg}$ and $\lsup{\Alg}P, \lsup{\Alg}Q, \lsup{\Alg}R\in
        \lsupv{\Alg}\ModCat$ with either $\lsup{\Alg}P, \lsup{\Alg}Q$ and
        $\lsup{\Alg}R$ bounded or $M_\Alg, N_\Alg, L_\Alg$ and the
        morphisms between them bounded.
	Then, the maps $$(\cdot\DT\Id_P)\co
	\Mor_{\Alg}(M_\Alg,N_\Alg)\to \Mor_\Ground(M_\Alg\DT \lsup{\Alg}P,N_\Alg\DT \lsup{\Alg}P)$$ and
	$$(\Id_M\DT\cdot)\co \Mor^{\Alg}(\lsup{\Alg}P,\lsup{\Alg}Q)\to
	\Mor_\Ground(M_\Alg\DT \lsup{\Alg}P,M_\Alg\DT \lsup{\Alg}Q)$$
        defined above are chain maps. Further:
        \begin{enumerate}
        \item\label{item:id-DT-1} The maps $(\Id_M\DT\cdot)$ and
          $(\cdot\DT\Id_P)$ are functorial under composition in sense
          that the following square commutes:
          \[
          \begin{tikzpicture}
            \node at (0,0) (tlcorn) {$\Mor_{\Alg}(M,N) \otimes
              \Mor_{\Alg}(N,L)$};
            \node at (9,0) (trcorn) {$\begin{array}{l} \\
                \Mor_{\Ground}(M\DT P,N \DT P) \\
                \qquad\otimes \Mor_{\Ground}(N\DT P,L\DT P)
              \end{array}$};
            \node at (0,-2) (blcorn) {$\Mor_{\Alg}(M,L)$};
            \node at (9,-2) (brcorn) {$\Mor_{\Ground}(M\DT P, L\DT P);$};
            \draw[->] (tlcorn) to node[above]{\lab{f\otimes g \mapsto (f\DT\Id_P)\otimes(g\DT \Id_P)}} (trcorn);
            \draw[->] (tlcorn) to node[left]{\lab{(f\otimes g)\mapsto
                g\circ f}} (blcorn);
            \draw[->] (trcorn) to node[right]{\lab{(k\otimes l)\mapsto l\circ k}} (brcorn);
            \draw[->] (blcorn) to node[below]{\lab{h\mapsto (h\DT \Id_P)}} (brcorn);
          \end{tikzpicture}
          \]
        while the following square commutes up to homotopy:
          \[
          \begin{tikzpicture}
            \node at (0,0) (tlcorn) {$\Mor^{\Alg}(P,Q)\otimes \Mor^{\Alg}(Q,R)$};
            \node at (9,0) (trcorn) {${\begin{array}{l}
              \Mor_{\Ground}(M\DT P,M \DT Q) \\
              \qquad \otimes \Mor_{\Ground}(M \DT Q,M\DT R)
            \end{array}}$};
            \node at (0,-2) (blcorn) {$\Mor^{\Alg}(P,R)$};
            \node at (9,-2) (brcorn) {$\Mor_{\Ground}(M\DT P, M\DT R).$};
            \draw[->] (tlcorn) to node[above]{\lab{f^1\otimes g^1 \mapsto (\Id_M\DT f^1)\otimes(\Id_M\DT
            g^1)}} (trcorn);
            \draw[->] (tlcorn) to node[left]{\lab{(f^1\otimes g^1)\mapsto g^1\circ f^1}} (blcorn);
            \draw[->] (trcorn) to node[right]{\lab{(k\otimes l)\mapsto
                l\circ k}} (brcorn);
            \draw[->] (blcorn) to node[below]{\lab{h^1\mapsto (\Id_M\DT h^1)}} (brcorn);
          \end{tikzpicture}
          \]
      \item\label{item:id-DT-2} The maps $(\cdot\DT\Id)$ and $(\Id\DT\cdot)$ commute with
        each other in the sense that the square
          \[
          \begin{tikzpicture}
            \node at (0,0) (tlcorn) {$\Mor_\Alg(M,N)\otimes
          \Mor^{\Alg}(P,Q)$};
            \node at (9,0) (trcorn) {${\begin{array}{l}\Mor_\Ground(M\DT P, M\DT Q)\\ \qquad\otimes \Mor_\Ground(M\DT Q, N\DT Q)\end{array}}$};
            \node at (0,-2) (blcorn) {$\begin{array}{l}\Mor_\Ground(M\DT P, N\DT P) \\
              \qquad \otimes \Mor_\Ground(N\DT P, N\DT Q)
            \end{array}$};
            \node at (9,-2) (brcorn) {$\Mor_\Ground(M\DT P, N\DT Q)$};
            \draw[->] (tlcorn) to node[above]{\lab{f \otimes g^1\mapsto (\Id_M\DT
            g^1)\otimes (f\DT \Id_Q)}} (trcorn);
            \draw[->] (tlcorn) to node[left]{
              \lab{\substack{(f\otimes g^1)\mapsto(f\DT \Id_P)\hfill\\
                  \qquad\otimes (\Id_N\DT g^1)}}} (blcorn);
            \draw[->] (trcorn) to node[right]{\lab{(k\otimes
                l)\mapsto l\circ k}} (brcorn);
            \draw[->] (blcorn) to node[below]{\lab{(k\otimes
                l)\mapsto l\circ k}} (brcorn);
          \end{tikzpicture}
          \]
	commutes up to homotopy.
      \end{enumerate}
\end{lemma}
\begin{proof}
	We leave as exercises that the $(\cdot\DT \Id_P)$ and $(\Id_M\DT \cdot)$ are
	chain maps. The fact that the first square commutes is straightforward.
	The homotopy $H$ for the second square
	$$H\co \Mor(P,Q) \otimes \Mor(Q,R) \to \Mor(M \DT P, M \DT R)$$
	is defined by
        \begin{multline}
	H(f^1\otimes g^1)=
	(\mu\otimes \Id_R)\circ (\Id_M \otimes (
        (\Id_{\overline{\Tensor^*}(\Alg[1])}\otimes \delta^R)\circ
        (\Id_{\overline{\Tensor^*}(\Alg[1])}\otimes g^1)\circ\\
        (\Id_{\overline{\Tensor^*}(\Alg[1])}\otimes \delta^Q)\circ
        (\Id_{\overline{\Tensor^*}(\Alg[1])}\otimes f^1)\circ
        \delta^P)).
        \end{multline}
         Note the similarity of the right hand side to the definition
         of the composition $f^1 \circ g^1$ in
         Formula~\eqref{eq:type-d-composition}.

        The homotopy $K$ for the third square
        $$K\co \Mor(M,N) \otimes \Mor(P,Q) \to \Mor(M\DT P, N\DT Q)$$
        is furnished by
        $$K(f,g^1)(x\otimes p) = (f\otimes \Id_Q)(x\otimes
        g(p)),$$
        or pictorially:
  \[
    \mathcenter{
      \begin{tikzpicture}
        \node at (0,0) (tlblank) {};
        \node at (1,0) (trblank) {};
        \node at (1,-1) (delta) {$\delta^N$};
        \node at (1,-2) (g) {$g^1$};
        \node at (1,-3) (delta2) {$\delta^N$};
        \node at (0,-4) (phi) {$f$};
        \node at (0,-5) (blblank) {};
        \node at (1,-5) (brblank) {};        
        \draw[modarrow] (trblank) to (delta);
        \draw[modarrow] (delta) to (g);
        \draw[modarrow] (g) to (delta2);
        \draw[modarrow] (delta2) to (brblank);
        \draw[modarrow] (tlblank) to (phi);
        \draw[modarrow] (phi) to (blblank);    
        \draw[tensoralgarrow] (delta) to (phi);
        \draw[algarrow] (g) to (phi);
        \draw[tensoralgarrow] (delta2) to (phi);
      \end{tikzpicture}
    }
  \]        
\end{proof}

Let us choose, arbitrarily, to define $f\DT g^1=(f\DT\Id)\circ (\Id\DT
g^1)$. Then:
\begin{corollary}\label{cor:DT-descends}
  \mbox{}
  \begin{enumerate}
  \item\label{item:dtfunc0}
    The operation $\DT$ induces a chain map
    $$\Mor_{\Alg}(M,N)\otimes \Mor^{\Alg}(P,Q)\to \Mor(M\DT P,N\DT Q).$$
  \item\label{item:dtfunc1} The operation $\DT$ is functorial up to
    homotopy. That is, $(f\DT g^1)\circ (f'\DT (g')^1)$ is homotopic to
    $(f\circ f')\DT (g^1\circ (g')^1)$.
  \item\label{item:dtfunc2} If $f\in\Mor(M_\Alg,N_\Alg)$ and
    $g^1\in\Mor(\lsup{\Alg}M,\lsup{\Alg}N)$ are cycles then $f\DT g^1$ is
    a cycle.
  \item\label{item:dtfunc3} If $f\in\Mor(M_\Alg,N_\Alg)$ and
    $g^1\in\Mor(\lsup{\Alg}M,\lsup{\Alg}N)$ are cycles and either $f$ or
    $g^1$ is nullhomotopic then $f\DT g$ is nullhomotopic.
  \item\label{item:dtfunc4} The operation $\DT$ descends to bifunctors
    \begin{align*}
      \DT&\co \HMod(\kbModCat_\Alg) \times \HMod({}^{\Alg}\ModCat)\to
      \HMod(\ModCat_\Ground)\\
      \DT&\co \HMod(\ModCat_\Alg) \times \HMod({}^{\Alg}\kbModCat)\to
      \HMod(\ModCat_\Ground).
    \end{align*}
  \end{enumerate}
\end{corollary}
\begin{proof}
  The fact that $\DT$ induces a chain map on morphism spaces
  follows from Lemma~\ref{lem:DTfunct}: it is defined as a composite
  of two chain maps.
  To verify part~(\ref{item:dtfunc1}), note that
  \begin{align*}
    (f\DT g^1)\circ (f'\DT (g')^1)&=[(f\DT \Id)\circ (\Id\DT g^1)]\circ[(f'\DT
    \Id)\circ (\Id\DT (g'^1))]\\
    &\sim [(f\DT \Id)\circ (f'\DT
    \Id)]\circ[\Id\DT g^1\circ\Id\DT (g')^1]\\
    &\sim (f\circ f'\DT \Id)\circ(\Id\DT g^1\circ (g')^1)\\
    &=(f\circ f')\DT (g^1\circ (g')^1)
  \end{align*}
  where the first homotopy uses part~(\ref{item:id-DT-2})
  of Lemma~\ref{lem:DTfunct} while the second homotopy
  uses part~(\ref{item:id-DT-1}) of Lemma~\ref{lem:DTfunct}.

  Parts~(\ref{item:dtfunc2}) and~(\ref{item:dtfunc3}) are easy to
  verify. Part~(\ref{item:dtfunc4}) then follow formally.
\end{proof}

\begin{corollary}
  \label{cor:DTChainHomotopyEquivalence}
  If $f^1\in\Mor^{\Alg}(P,Q)$ is a chain homotopy equivalence of type
  $D$ structures, then $\Id_M\DT f^1\in \Mor(M\DT P,M\DT Q)$ is a chain
  homotopy equivalence of complexes. Similarly, if $\phi\in
  \Mor_{\Alg}(M,N)$ is a chain homotopy equivalence of $\Alg$-modules,
  then $\phi\DT \Id_P\co \Mor(M\DT P, N\DT P)$ is a chain homotopy
  equivalence of complexes.
\end{corollary}

While commuting $f\DT \Id$ and $\Id\DT g$ was somewhat
subtle, ($\Ainf$) functoriality of $\DT$ in each factor is more
straightforward:
\begin{lemma}\label{lem:DT-ainf-each-factor} Let $M_\Alg$ be an $\Ainf$-module and $\lsup{\Alg}N$ a
  type $D$ structure. Then:
  \begin{enumerate}
  \item The operation $\cdot_\Alg\DT\lsup{\Alg}N$ gives a \dg functor
    $\ModCat_\Alg\to \ModCat_{\Ground}$ and
  \item The operation $M_\Alg\DT\lsup{\Alg}\cdot$ extends to an
    $\Ainf$-functor $\lsupv{\Alg}\ModCat\to \ModCat_{\Ground}$.
  \end{enumerate}
\end{lemma}
\begin{proof}
  The first part is straightforward. For the second, if
  $\lsup{\Alg}N_i$, $i=0,\dots,n$ are type $D$ structures and $f_i\co
  \lsup{\Alg}N_i\to\lsup{\Alg}N_{i+1}$ are morphisms, define
  \[
  \mathcenter{(M_\Alg \DT\cdot)_{1,n}(\Id_M,f_1,\dots,f_n)=}
  \mathcenter{
  \begin{tikzpicture}
    \node at (-1,0) (tlblank) {};
    \node at (1,0) (trblank) {};
    \node at (1,-1) (delta1) {$\delta^{N_{0}}$};
    \node at (1,-2) (f1) {$f_1$};
    \node at (1,-3) (delta2) {$\delta^{N_{1}}$};
    \node at (1,-4) (rdots) {$\vdots$};
    \node at (1,-5) (delta3) {$\delta^{N_{n-1}}$};
    \node at (1,-6) (f2) {$f_n$};
    \node at (1,-7) (delta4) {$\delta^{N_n}$};
    \node at (-1,-8) (m) {$m_M$};
    \node at (1,-9) (brblank) {};
    \node at (-1,-9) (blblank) {};
    \draw[modarrow] (tlblank) to (m);
    \draw[modarrow] (m) to (blblank);
    \draw[modarrow] (trblank) to (delta1);
    \draw[modarrow] (delta1) to (f1);
    \draw[modarrow] (f1) to (delta2);
    \draw[modarrow] (delta2) to (rdots);
    \draw[modarrow] (rdots) to (delta3);
    \draw[modarrow] (delta3) to (f2);
    \draw[modarrow] (f2) to (delta4);
    \draw[modarrow] (delta4) to (brblank);
    \draw[tensoralgarrow, bend right=10] (delta1) to (m);
    \draw[tensoralgarrow, bend right=8] (delta2) to (m);
    \draw[tensoralgarrow, bend right=5] (delta3) to (m);
    \draw[tensoralgarrow] (delta4) to (m);
    \draw[algarrow, bend right=9] (f1) to (m);
    \draw[algarrow, bend right=2] (f2) to (m);
  \end{tikzpicture}}.
  \]
  It is straightforward to verify that this makes $M_\Alg\DT\cdot$
  into an $\Ainf$-functor.
\end{proof}
\begin{remark}
  Even when $\Alg$ is a \dg algebra, so $\lsupv{\Alg}\ModCat$ is an
  honest \dg category, the operation $M_\Alg\DT\lsup{\Alg}\cdot$ is
  still only an $\Ainf$-functor.
\end{remark}

Next we turn to the behavior of $\DT$ for bimodules. Since we have
various kinds of bimodules, there are various cases of the tensor
product:
\begin{itemize}
\item $\DA \DT\DD$ is a type \DD\ module.
\item $\AAm\DT\DD$ is a type \AD\ module.
\item $\DA\DT\DA$ is a type \DA\ module.
\item $\AAm\DT\DA$ is a type \AAm\  module.
\end{itemize}
(In each case, we assume one of the factors in the tensor product is
appropriately bounded; see Proposition~\ref{prop:DT-many-cases}.)

\begin{figure}
  \centering
  \[
  \mathcenter{
    \begin{tikzpicture}
      \node at (-1,0) (tllblank) {};
      \node at (0,0) (tlblank) {}; 
      \node at (1,0) (trblank) {};
      \node at (2,0) (trrblank) {};
      \node at (1,-1) (delta) {$\delta^N$}; 
      \node at (0,-2) (m) {$m$}; 
      \node at (0,-3) (blblank) {};
      \node at (1,-3) (brblank) {};
      \draw[modarrow] (trblank) to (delta); 
      \draw[modarrow] (delta) to (brblank); 
      \draw[modarrow] (tlblank) to (m); 
      \draw[modarrow] (m) to (blblank);
      \draw[tensoralgarrow] (tllblank) to (m);
      \draw[tensoralgarrow] (trrblank) to (delta);
      \draw[tensoralgarrow] (delta) to (m);
    \end{tikzpicture}
  }
  \hspace{.25in}
  \mathcenter{
    \begin{tikzpicture}
      \node at (0,0) (tlblank) {}; 
      \node at (1,0) (trblank) {};
      \node at (2,0) (trrblank) {};
      \node at (1,-1) (deltar) {$\delta^N$}; 
      \node at (0,-2) (deltal) {$\delta^1$}; 
      \node at (-1,-3) (bllblank) {}; 
      \node at (0,-3) (blblank) {};
      \node at (1,-3) (brblank) {};
      \draw[modarrow] (trblank) to (deltar); 
      \draw[modarrow] (deltar) to (brblank); 
      \draw[modarrow] (tlblank) to (deltal); 
      \draw[modarrow] (deltal) to (blblank);
      \draw[algarrow] (deltal) to (bllblank);
      \draw[tensoralgarrow] (trrblank) to (deltar);
      \draw[tensoralgarrow] (deltar) to (deltal);
    \end{tikzpicture}
  }
  \hspace{.25in}  
  \mathcenter{
    \begin{tikzpicture}
      \node at (0,0) (tlblank) {}; 
      \node at (1,0) (trblank) {};
      \node at (-1,0) (tllblank) {};
      \node at (1,-1) (delta) {$\delta$}; 
      \node at (0,-2) (m) {$m$}; 
      \node at (2,-2) (prod) {$\Pi$};
      \node at (2,-3) (brrblank) {}; 
      \node at (0,-3) (blblank) {};
      \node at (1,-3) (brblank) {};
      \draw[modarrow] (trblank) to (delta); 
      \draw[modarrow] (delta) to (brblank); 
      \draw[modarrow] (tlblank) to (m); 
      \draw[modarrow] (m) to (blblank);
      \draw[tensoralgarrow] (tllblank) to (m);
      \draw[tensoralgarrow] (delta) to (prod);
      \draw[tensoralgarrow] (delta) to (m);
      \draw[algarrow] (prod) to (brrblank);
    \end{tikzpicture}
  }
  \hspace{.25in}  
  \mathcenter{
    \begin{tikzpicture}
      \node at (0,0) (tlblank) {}; 
      \node at (1,0) (trblank) {};
      \node at (1,-1) (delta) {$\delta$}; 
      \node at (0,-2) (deltal) {$\delta^1$}; 
      \node at (2,-2) (prod) {$\Pi$};
      \node at (2,-3) (brrblank) {}; 
      \node at (1,-3) (brblank) {};
      \node at (0,-3) (blblank) {};
      \node at (-1,-3) (bllblank) {};
      \draw[modarrow] (trblank) to (delta); 
      \draw[modarrow] (delta) to (brblank); 
      \draw[modarrow] (tlblank) to (m); 
      \draw[modarrow] (deltal) to (blblank);
      \draw[algarrow] (deltal) to (bllblank);
      \draw[tensoralgarrow] (delta) to (prod);
      \draw[tensoralgarrow] (delta) to (deltal);
      \draw[algarrow] (prod) to (brrblank);
    \end{tikzpicture}
  }
  \]
  \caption{\textbf{Bimodule structure on $\DT$-tensor products of
      bimodules.} Far left: $\AAm \DT \DA$. Center left: $\DA\DT
    \DA$. Center right: $\AAm\DT \DD$. Far right: $\DA\DT \DD$. Here, for
    elements $a_1,\dots,a_n\in\DGA$ of a \dg algebra we let
    $\Pi(a_1\otimes\dots\otimes a_n)$ denote the product of the
    elements, $a_1\cdots a_n$.}
\label{fig:bimod-on-DT}
\end{figure}
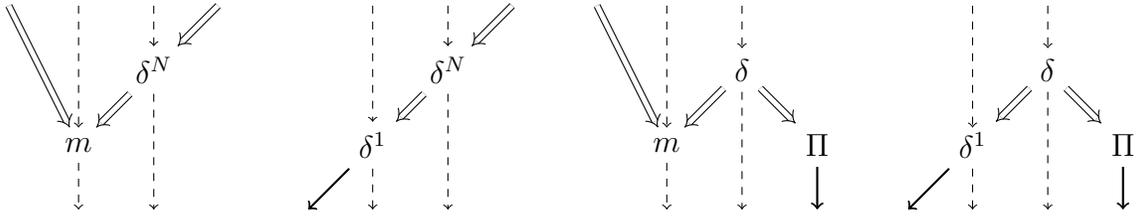

\begin{definition}
  Let $\lsub{\Alg}M_\Blg$ and $\lsup{\Blg}N_\Clg$ (respectively
  $\lsup{\Alg}M_\Blg$ and $\lsup{\Blg}N_\Clg$, $\lsub{\Alg}M_\Blg$ and
  $\lsup{\Blg}N^\Clg$,  $\lsup{\Alg}M_\Blg$ and $\lsup{\Blg}N^\Clg$)
  by a type \AAm\  and \DA\ (respectively \DA\ and \DA, \AAm\  and \DD,
  \DA\ and \DD)\footnote{In the last two cases we assume that $\Blg$
    and $\Clg$ are \dg algebras. In general, we implicitly add this hypothesis
    any time a type \DD\ structure is mentioned.} bimodules with either $M$ right bounded or $N$ left
  bounded. As a chain complex, define
  \[
  \lsub{\Alg}M_\Blg\DT
  \lsup{\Blg}N_\Clg=\mathcal{F}(\lsub{\Alg}M_\Blg)_\Blg\DT \lsup{\Blg}{\mathcal{F}}(\lsup{\Blg}N_\Clg)
  \]
  with \AAm\  (respectively \DA, \AD, \DD) structure map given as in
  Figure~\ref{fig:bimod-on-DT} far left (respectively center left,
  center right, far right).
\end{definition}
(In the figure, we use the following notation: given elements
$b_1,\dots, b_k$ in a \dg algebra $\Blg$,
$\Pi(b_1,\dots,b_k)=b_1\cdots b_k$ denotes their product.  Note that
\DD\ bimodules by assumption involve \dg algebras.)

Since most of the results in all of these cases are quite similar, we
will often use the ambiguous notation $M$ or $N$ to refer to any
consistent way of placing superscripts and subscripts.
\begin{proposition}\label{prop:DT-many-cases}
  The condition that $M$ be right bounded or $N$ be left bounded
  guarantees that the sums defining the structure maps for $M\DT N$ in
  Figure~\ref{fig:bimod-on-DT} are finite. In this case they satisfy
  the corresponding structure equations. Moreover:
  \begin{enumerate}
  \item If $M$ and $N$ are both left bounded (respectively right
    bounded) then $M\DT N$ is left bounded (respectively right
    bounded).
  \item If $M$ (respectively $N$) is bounded then $M\DT N$ is left
    bounded (respectively right bounded).
  \item If $M$ (respectively $N$) is bounded and $N$ (respectively
    $M$) is right bounded (respectively left bounded) then $M \DT N$ is
    bounded.
 \end{enumerate}
\end{proposition}
The reader might find it
interesting to compare Proposition~\ref{prop:DT-many-cases} with
Lemma~\ref{lem:admiss-glues}.

In order to prove Proposition~\ref{prop:DT-many-cases} without
checking all the cases individually, we will reformulate and unify
the various definitions of boundedness.  Given an operation graph of one of the types
considered above
(planar, directed graphs with labelled nodes and edges, obeying
certain restrictions), we can restrict the inputs and outputs to lie
in $A_+$ and $B_+$. (For the outputs,
this means applying
$(1-\epsilon)$ to each output edge.)  This gives a map
\[
m_\Gamma^+: (A_+)^{\otimes k_1} \otimes M \otimes (B_+)^{\otimes k_2}
 \to (A_+)^{\otimes l_1} \otimes M \otimes (B_+)^{\otimes l_2}
\]
for appropriate values of $k_1$, $k_2$, $l_1$ and $l_2$.  (Some of $k_1$, $k_2$,
$l_1$ and $l_2$ will necessarily be~$0$, depending on the type of
bimodule.)

\begin{lemma}\label{lem:bounded-char}
  A bimodule~$M$ (of any type) is bounded if and only if for each~$x \in M$ there
  is a bound on the number of leaves of bimodule operation
  trees~$\Gamma$ for which the corresponding operation $m_\Gamma^+$ is
  non-zero when applied to~$x$.  It is left bounded if and only if for any $x$ and
  any bound on the number of right inputs/outputs there is a bound on
  the number of left inputs/outputs.  Similar statements hold for
  right bounded bimodules and for bimodule morphisms.
\end{lemma}
\begin{proof}[Proof sketch]
  This is very close to the definition of boundedness or left/right
  boundedness in each case, with the exception of the restriction to
  spinal graphs (without `$\mu$' nodes).  Let $\Gamma'$ be the graph
  obtained from $\Gamma$ by pushing all $\epsilon$ nodes upstream through
  $\mu$ nodes, using the relations in
  Equation~\eqref{eq:augmentation}, as far as possible, and then
  taking the subgraph formed by the `$m$' or `$\delta$' nodes
  and any adjacent `$\epsilon$' nodes.
  Because $A_+$ and $B_+$ are nilpotent, there
  is a bound on the number of inputs to $\Gamma$ that can
  contribute to an input to~$\Gamma'$.  Similarly,
  there is a bound on the number of outputs from~$\Gamma'$ that
  can contribute to an output of~$\Gamma$.  Thus, bounds on the
  inputs/outputs of~$\Gamma$ give bounds on the inputs/outputs of
  $\Gamma'$, and vice versa.
\end{proof}
\begin{proof}[Proof sketch of Proposition~\ref{prop:DT-many-cases}]
  All of the bimodule structure operations in
  Figure~\ref{fig:bimod-on-DT} can be expanded out so that each
  algebra edge carries an element of $A_+$, $B_+$ or $C_+$, rather than just
  $A$, $B$ or $C$, simply by
  taking a sum of terms where we apply $\epsilon$ or $(1-\epsilon)$ on
  each edge.  Now suppose $N$ is left bounded and we wish to compute a
  structure map on $M \DT N$ with some fixed number of outputs/inputs
  on the right and left.  Then there is a bound on the number of $B_+$
  edges leaving the dotted line corresponding to~$N$, which
  immediately gives a bound on the number of terms contributing to the
  definition of the structure map.

  A similar argument works if $M$ is right bounded.

  If $M$ and $N$ are both left bounded and we have a bound on the
  number of terms on the right of the entire diagram, then
  left-boundedness of~$N$ gives us a bound on the number of algebra
  edges communicating between the two dotted lines, and left-boundedness
  of~$M$ then gives a bound on the number on the number of algebra
  edges on the left of the diagram, as desired to show that $M \DT N$
  is left bounded.

  The other cases are similar.
\end{proof}

\begin{remark}\label{rmk:weaker-stronger-boundedness}
  If we drop the assumption that the algebras involved are nilpotent,
  Lemma \ref{lem:bounded-char} becomes false, but most cases of
  Proposition~\ref{prop:DT-many-cases} remain true.
  However, in a tensor product where the right factor is a \DD\ module
  (i.e., $\lsub{\Alg}M_\Blg \DT \lsup{\Blg}N^\Clg$ or $\lsup{\Alg}M_\Blg \DT
  \lsup{\Blg}N^\Clg$), if we only assume that $N$ is left bounded it
  does not follow that the sums in Figure~\ref{fig:bimod-on-DT} are
  finite.  A strengthening of the definition of left/right boundedness
  for \DD\ bimodules fixes this case.  However, if $\Alg$ is not
  nilpotent, $\lsub{\Alg}\Alg_\Alg$ is not usually left or right bounded
  which, for instance, breaks Proposition~\ref{prop:D-to-A-and-back}.
\end{remark}

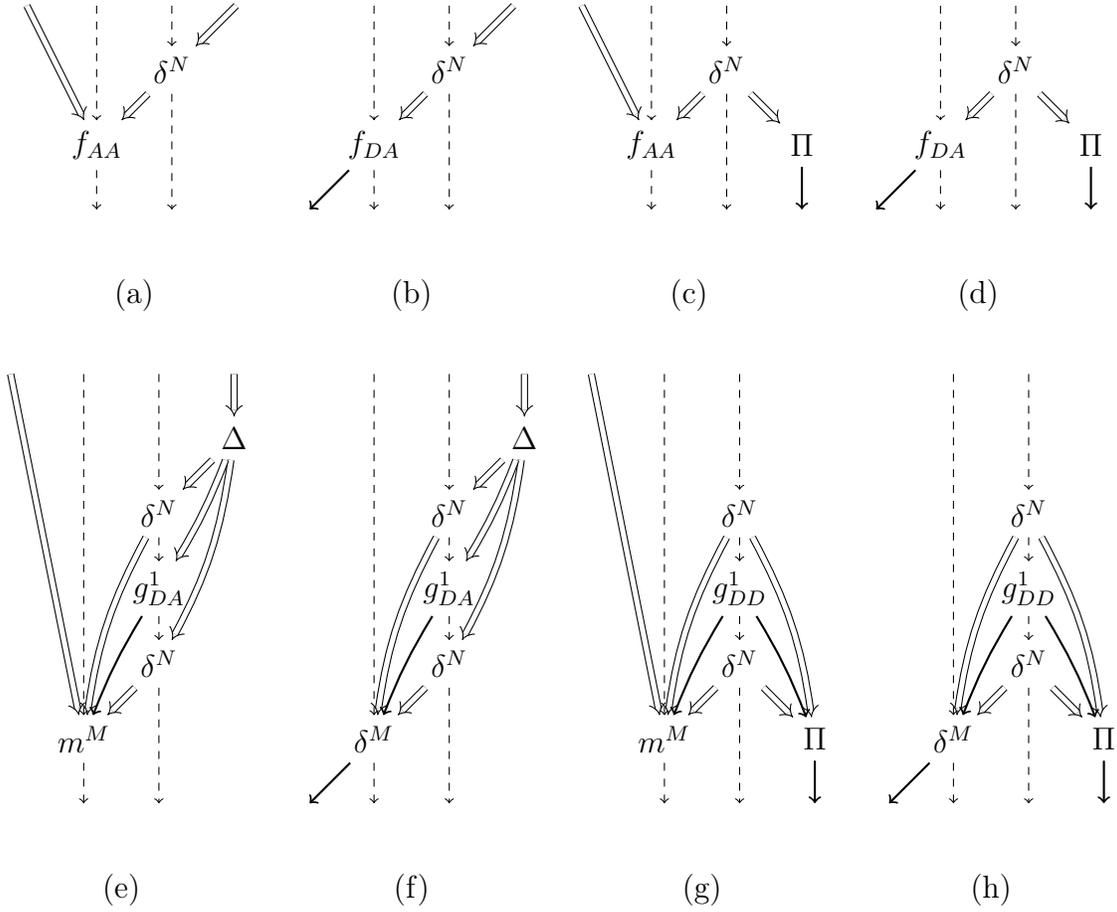
\begin{figure}
  \centering
  \[ 
  \mathcenter{
    \begin{tikzpicture}
      \node at (0,0) (tlblank) {};
      \node at (1,0) (trblank) {};
      \node at (1,-1) (delta) {$\delta^N$};
      \node at (0,-2) (f) {$f_{AA}$};
      \node at (0,-3) (blblank) {};
      \node at (1,-3) (brblank) {};
      \node at (2,0) (trrblank) {};
      \node at (-1,0) (tllblank) {};
      \node at (.5,-4) (caption) {(a)};
      \draw[Amodar] (tlblank) to (f);
      \draw[Amodar] (f) to (blblank);
      \draw[Dmodar] (trblank) to (delta);
      \draw[Dmodar] (delta) to (brblank);
      \draw[tensorblgarrow] (delta) to (f);    
      \draw[tensorclgarrow] (trrblank) to (delta);
      \draw[tensoralgarrow] (tllblank) to (f);
    \end{tikzpicture}
  }
  \quad
  \mathcenter{
    \begin{tikzpicture}
      \node at (0,0) (tlblank) {};
      \node at (1,0) (trblank) {};
      \node at (1,-1) (delta) {$\delta^N$};
      \node at (0,-2) (f) {$f_{DA}$};
      \node at (0,-3) (blblank) {};
      \node at (1,-3) (brblank) {};
      \node at (2,0) (trrblank) {};
      \node at (-1,-3) (bllblank) {};
      \node at (.5,-4) (caption) {(b)};
      \draw[Amodar] (tlblank) to (f);
      \draw[Amodar] (f) to (blblank);
      \draw[Dmodar] (trblank) to (delta);
      \draw[Dmodar] (delta) to (brblank);
      \draw[tensorblgarrow] (delta) to (f);    
      \draw[tensorclgarrow] (trrblank) to (delta);
      \draw[algarrow] (f) to (bllblank);
    \end{tikzpicture}
  }
  \quad
  \mathcenter{
    \begin{tikzpicture}
      \node at (0,0) (tlblank) {};
      \node at (1,0) (trblank) {};
      \node at (1,-1) (delta) {$\delta^N$};
      \node at (0,-2) (f) {$f_{AA}$};
      \node at (2,-2) (Pi) {$\Pi$};
      \node at (0,-3) (blblank) {};
      \node at (1,-3) (brblank) {};
      \node at (2,-3) (brrblank) {};
      \node at (-1,0) (tllblank) {};
      \node at (.5,-4) (caption) {(c)};
      \draw[Amodar] (tlblank) to (f);
      \draw[Amodar] (f) to (blblank);
      \draw[Dmodar] (trblank) to (delta);
      \draw[Dmodar] (delta) to (brblank);
      \draw[tensorblgarrow] (delta) to (f);    
      \draw[tensorclgarrow] (delta) to (Pi);
      \draw[clgarrow] (Pi) to (brrblank);
      \draw[tensoralgarrow] (tllblank) to (f);
    \end{tikzpicture}
  }
  \quad
  \mathcenter{
    \begin{tikzpicture}
      \node at (0,0) (tlblank) {};
      \node at (1,0) (trblank) {};
      \node at (1,-1) (delta) {$\delta^N$};
      \node at (0,-2) (f) {$f_{DA}$};
      \node at (2,-2) (Pi) {$\Pi$};
      \node at (0,-3) (blblank) {};
      \node at (1,-3) (brblank) {};
      \node at (2,-3) (brrblank) {};
      \node at (-1,-3) (bllblank) {};
      \node at (.5,-4) (caption) {(d)};
      \draw[Amodar] (tlblank) to (f);
      \draw[Amodar] (f) to (blblank);
      \draw[Dmodar] (trblank) to (delta);
      \draw[Dmodar] (delta) to (brblank);
      \draw[tensorblgarrow] (delta) to (f);    
      \draw[tensorclgarrow] (delta) to (Pi);
      \draw[clgarrow] (Pi) to (brrblank);
      \draw[algarrow] (f) to (bllblank);
    \end{tikzpicture}
  }
\]

\[
\mathcenter{
  \begin{tikzpicture}
    \node at (0,0) (tlblank) {};
    \node at (1,0) (trblank) {};
    \node at (1,-2) (delta1) {$\delta^N$};
    \node at (1,-3) (g) {$g_{DA}^1$};
    \node at (1,-4) (delta2) {$\delta^N$};
    \node at (0,-5) (deltaM) {$m^M$};
    \node at (2, -1) (Delta) {$\Delta$};
    \node at (0,-6) (blblank) {};
    \node at (1,-6) (brblank) {};
    \node at (2,0) (trrblank) {};
    \node at (-1,0) (tllblank) {};
    \node at (.5,-7) (caption) {(e)};
    \draw[Amodar] (tlblank) to (deltaM);
    \draw[Amodar] (deltaM) to (blblank);
    \draw[Dmodar] (trblank) to (delta1);
    \draw[Dmodar] (delta1) to (g);
    \draw[Dmodar] (g) to (delta2);
    \draw[Dmodar] (delta2) to (brblank);
    \draw[tensorblgarrow, bend right=10] (delta1) to (deltaM);
    \draw[blgarrow, bend right=5] (g) to (deltaM);
    \draw[tensorblgarrow] (delta2) to (deltaM);
    \draw[tensoralgarrow] (tllblank) to (deltaM);
    \draw[tensorclgarrow] (Delta) to (delta1);
    \draw[tensorclgarrow, bend left=10] (Delta) to (delta2);
    \draw[tensorclgarrow, bend left=5] (Delta) to (g);
    \draw[tensorclgarrow] (trrblank) to (Delta);
  \end{tikzpicture}
}  
\quad
\mathcenter{
  \begin{tikzpicture}
    \node at (0,0) (tlblank) {};
    \node at (1,0) (trblank) {};
    \node at (1,-2) (delta1) {$\delta^N$};
    \node at (1,-3) (g) {$g_{DA}^1$};
    \node at (1,-4) (delta2) {$\delta^N$};
    \node at (0,-5) (deltaM) {$\delta^M$};
    \node at (2, -1) (Delta) {$\Delta$};
    \node at (0,-6) (blblank) {};
    \node at (1,-6) (brblank) {};
    \node at (2,0) (trrblank) {};
    \node at (-1,-6) (bllblank) {};
    \node at (.5,-7) (caption) {(f)};
    \draw[Amodar] (tlblank) to (deltaM);
    \draw[Amodar] (deltaM) to (blblank);
    \draw[Dmodar] (trblank) to (delta1);
    \draw[Dmodar] (delta1) to (g);
    \draw[Dmodar] (g) to (delta2);
    \draw[Dmodar] (delta2) to (brblank);
    \draw[tensorblgarrow, bend right=10] (delta1) to (deltaM);
    \draw[blgarrow, bend right=5] (g) to (deltaM);
    \draw[tensorblgarrow] (delta2) to (deltaM);
    \draw[algarrow] (deltaM) to (bllblank);
    \draw[tensorclgarrow] (Delta) to (delta1);
    \draw[tensorclgarrow, bend left=10] (Delta) to (delta2);
    \draw[tensorclgarrow, bend left=5] (Delta) to (g);
    \draw[tensorclgarrow] (trrblank) to (Delta);
  \end{tikzpicture}
}  
\quad
\mathcenter{
  \begin{tikzpicture}
    \node at (0,0) (tlblank) {};
    \node at (1,0) (trblank) {};
    \node at (1,-2) (delta1) {$\delta^N$};
    \node at (1,-3) (g) {$g_{DD}^1$};
    \node at (1,-4) (delta2) {$\delta^N$};
    \node at (0,-5) (deltaM) {$m^M$};
    \node at (2, -5) (Pi) {$\Pi$};
    \node at (0,-6) (blblank) {};
    \node at (1,-6) (brblank) {};
    \node at (2,-6) (brrblank) {};
    \node at (-1,0) (tllblank) {};
    \node at (.5,-7) (caption) {(g)};
    \draw[Amodar] (tlblank) to (deltaM);
    \draw[Amodar] (deltaM) to (blblank);
    \draw[Dmodar] (trblank) to (delta1);
    \draw[Dmodar] (delta1) to (g);
    \draw[Dmodar] (g) to (delta2);
    \draw[Dmodar] (delta2) to (brblank);
    \draw[tensorblgarrow, bend right=10] (delta1) to (deltaM);
    \draw[blgarrow, bend right=5] (g) to (deltaM);
    \draw[tensorblgarrow] (delta2) to (deltaM);
    \draw[tensoralgarrow] (tllblank) to (deltaM);
    \draw[tensorclgarrow, bend left=10] (delta1) to (Pi);
    \draw[tensorclgarrow] (delta2) to (Pi);
    \draw[clgarrow, bend left=5] (g) to (Pi);
    \draw[clgarrow] (Pi) to (brrblank);
  \end{tikzpicture}
}  
\quad
\mathcenter{
  \begin{tikzpicture}
    \node at (0,0) (tlblank) {};
    \node at (1,0) (trblank) {};
    \node at (1,-2) (delta1) {$\delta^N$};
    \node at (1,-3) (g) {$g_{DD}^1$};
    \node at (1,-4) (delta2) {$\delta^N$};
    \node at (0,-5) (deltaM) {$\delta^M$};
    \node at (2, -5) (Pi) {$\Pi$};
    \node at (0,-6) (blblank) {};
    \node at (1,-6) (brblank) {};
    \node at (2,-6) (brrblank) {};
    \node at (-1,-6) (bllblank) {};
    \node at (.5,-7) (caption) {(h)};
    \draw[Amodar] (tlblank) to (deltaM);
    \draw[Amodar] (deltaM) to (blblank);
    \draw[Dmodar] (trblank) to (delta1);
    \draw[Dmodar] (delta1) to (g);
    \draw[Dmodar] (g) to (delta2);
    \draw[Dmodar] (delta2) to (brblank);
    \draw[tensorblgarrow, bend right=10] (delta1) to (deltaM);
    \draw[blgarrow, bend right=5] (g) to (deltaM);
    \draw[tensorblgarrow] (delta2) to (deltaM);
    \draw[algarrow] (deltaM) to (bllblank);
    \draw[tensorclgarrow, bend left=10] (delta1) to (Pi);
    \draw[tensorclgarrow] (delta2) to (Pi);
    \draw[clgarrow, bend left=5] (g) to (Pi);
    \draw[clgarrow] (Pi) to (brrblank);
  \end{tikzpicture}
}  
  \]  
  \caption{\textbf{The box product of a bimodule morphism with the
      identity morphism.}}
  \label{fig:f-box-id-bimod}
\end{figure}

As was the case for tensoring type $D$ and $A$ modules, the tensor
product for bimodules is not strictly functorial. Again, we define the
box product of two morphisms in terms of the box product of a morphism
with the identity morphism. There are now eight cases:
\begin{itemize}
\item Given $f_{AA}\co \lsub{\Alg}M_\Blg\to \lsub{\Alg}M'_\Blg$ and
  $\lsup{\Blg}N_\Clg$ define $f_{AA}\DT \Id_{N}$ as in
  Figure~\ref{fig:f-box-id-bimod} part (a).
\item Given $f_{DA}\co \lsup{\Alg}M_\Blg\to \lsup{\Alg}M'_\Blg$ and
  $\lsup{\Blg}N_\Clg$ define $f_{DA}\DT \Id_N$ as in
  Figure~\ref{fig:f-box-id-bimod} part (b).
\item Given $f_{AA}\co \lsub{\Alg}M_\Blg\to \lsub{\Alg}M'_\Blg$ and
  $\lsup{\Blg}N^\Clg$ define $f_{AA}\DT \Id_N$ as in
  Figure~\ref{fig:f-box-id-bimod} part (c).
\item Given $f_{DA}\co \lsub{\Alg}M_\Blg\to \lsub{\Alg}M'_\Blg$ and
  $\lsup{\Blg}N^\Clg$ define $f_{DA}\DT \Id_N$ as in
  Figure~\ref{fig:f-box-id-bimod} part (d).
\item Given $\lsub{\Alg}M_\Blg$ and $g_{DA}\co \lsup{\Blg}N_\Clg\to
  \lsup{\Blg}N'_\Clg$ define $\Id_M\DT g_{DA}$ as in 
  Figure~\ref{fig:f-box-id-bimod} part (e).
\item Given $\lsup{\Alg}M_\Blg$ and $g_{DA}\co \lsup{\Blg}N_\Clg\to
  \lsup{\Blg}N'_\Clg$ define $\Id_M\DT g_{DA}$ as in 
  Figure~\ref{fig:f-box-id-bimod} part (f).
\item Given $\lsub{\Alg}M_\Blg$ and $g_{DD}\co \lsup{\Blg}N^\Clg\to
  \lsup{\Blg}(N')^\Clg$ define $\Id_M\DT g_{DD}$ as in 
  Figure~\ref{fig:f-box-id-bimod} part (g).
\item Given $\lsup{\Alg}M_\Blg$ and $g_{DD}\co \lsup{\Blg}N^\Clg\to
  \lsup{\Blg}(N')^\Clg$ define $\Id_M\DT g_{DD}$ as in 
  Figure~\ref{fig:f-box-id-bimod} part (h).
\end{itemize}
In all cases, define $f\DT g$ to be $(f\DT\Id)\circ (\Id\DT g^1)$.

With these definitions, the obvious analogue of
Lemma~\ref{lem:DTfunct} holds. Moreover, similarly to
Corollary~\ref{cor:DT-descends} and
Lemma~\ref{lem:DT-ainf-each-factor} we have:
\begin{lemma}  \label{lem:DT-descends-bimod}
\mbox{}
  \begin{enumerate}
  \item The map $\DT\co \Mor(M,M')\otimes \Mor(N,N')\to \Mor(M\DT
    N,M'\DT N')$ is a chain map.
  \item The operation $\DT$ is functorial up to
    homotopy. That is, $(f\DT g)\circ (f'\DT g')$ is homotopic to
    $(f\circ f')\DT (g\circ g')$.
  \item The operations $M_\Alg\DT\cdot $ and $\cdot\DT\lsup{\Alg}N$
    extend to $\Ainf$-functors, and so
  \item the operations $\DT$ descend to bifunctors of homotopy
    categories.
  \end{enumerate}
\end{lemma}

Next we turn to the question of associativity of tensor product. Like
functoriality on the \dg level, this is somewhat subtle. There are several
cases when it is straightforward, however:

\begin{lemma}\label{lem:tens-assoc-1-strict}
\mbox{}
  \begin{enumerate}
  \item Let $\lsub{\Alg}N_\Blg$ be a type \AAm\ module and $M^\Alg$
    and $\lsup{\Blg}P$ type $D$ structures. Then there is a canonical
    isomorphism
    \[
    (M^\Alg\DT \lsub{\Alg}N_\Blg)\DT\lsup{\Blg}P\cong M^\Alg\DT
    (\lsub{\Alg}N_\Blg\DT\lsup{\Blg}P).
    \]
  \item Let $\lsup{\Alg}N_\Blg$ be a type \DA\ structure, $M_\Alg$ a
    type $A$ module and $\lsup{\Blg}P$ a type $D$ structure. Then
    there is a canonical isomorphism
    \[
    (M_\Alg\DT \lsup{\Alg}N_\Blg)\DT\lsup{\Blg}P\cong M_\Alg\DT
    (\lsup{\Alg}N_\Blg\DT\lsup{\Blg}P).
    \]
  \item 
    \label{case:DAintheMiddle}
    Let $\lsup{\Alg}N^\Blg$ be a separated type \DD\ structure as in
    Definition~\ref{def:DD-sep}. Then
    for any $M_\Alg$ and $\lsub{\Blg}P$ there is a canonical isomorphism
    \[
    (M_\Alg\DT \lsup{\Alg}N^\Blg)^\Blg\DT\lsub{\Blg}P\cong
    M_\Alg\DT \lsup{\Blg}(\lsup{\Alg}N^\Blg\DT\lsub{\Blg}P).
    \]
  \end{enumerate}
  Similar statements hold if $M$ and / or $P$ is a bimodule (of any
  type compatible with the $\DT$ products).
\end{lemma}
\begin{proof}
  The differential on the triple box product in the three cases is given
  by the diagrams
  \[
  \text{(1)}\mathcenter{
    \begin{tikzpicture}
      \node at (-1,0) (tlblank) {};
      \node at (0,0) (tcblank) {};
      \node at (1,0) (trblank) {};
      \node at (-1,-1) (deltaM) {$\delta^M$};
      \node at (1,-1) (deltaP) {$\delta^P$};
      \node at (0,-2) (mN) {$m^N$};
      \node at (-1,-3) (blblank) {};
      \node at (0,-3) (bcblank) {};
      \node at (1,-3) (brblank) {};
      \draw[Dmodar] (tlblank) to (deltaM);
      \draw[Dmodar] (deltaM) to (blblank);
      \draw[Amodar] (tcblank) to (mN);
      \draw[Amodar] (mN) to (bcblank);
      \draw[Dmodar] (trblank) to (deltaP);
      \draw[Dmodar] (deltaP) to (brblank);
      \draw[tensoralgarrow] (deltaM) to (mN);
      \draw[tensorblgarrow] (deltaP) to (mN);
    \end{tikzpicture}
  }
  \qquad
  \text{(2)}\mathcenter{
    \begin{tikzpicture}
      \node at (-1,0) (tlblank) {};
      \node at (0,0) (tcblank) {};
      \node at (1,0) (trblank) {};
      \node at (-1,-3) (mM) {$m^M$};
      \node at (1,-1) (deltaP) {$\delta^P$};
      \node at (0,-2) (deltaN) {$\delta^N$};
      \node at (-1,-4) (blblank) {};
      \node at (0,-4) (bcblank) {};
      \node at (1,-4) (brblank) {};
      \draw[Dmodar] (tlblank) to (mM);
      \draw[Dmodar] (mM) to (blblank);
      \draw[Amodar] (tcblank) to (deltaN);
      \draw[Amodar] (deltaN) to (bcblank);
      \draw[Dmodar] (trblank) to (deltaP);
      \draw[Dmodar] (deltaP) to (brblank);
      \draw[tensoralgarrow] (deltaN) to (mM);
      \draw[tensorblgarrow] (deltaP) to (deltaN);      
    \end{tikzpicture}
  }
  \qquad
  \text{(3)}\mathcenter{
    \begin{tikzpicture}
      \node at (-1,0) (tlblank) {};
      \node at (0,0) (tcblank) {};
      \node at (1,0) (trblank) {};
      \node at (-1,-2) (mM) {$m^M$};
      \node at (0,-1) (deltaN) {$\delta^{N,L}$};
      \node at (-1,-3) (blblank) {};
      \node at (0,-3) (bcblank) {};
      \node at (1,-3) (brblank) {};
      \draw[Dmodar] (tlblank) to (mM);
      \draw[Dmodar] (mM) to (blblank);
      \draw[Amodar] (tcblank) to (deltaN);
      \draw[Amodar] (deltaN) to (bcblank);
      \draw[Dmodar] (trblank) to (brblank);
      \draw[tensoralgarrow] (deltaN) to (mM);
    \end{tikzpicture}
  }\,\,
  +
  \,\,\mathcenter{
    \begin{tikzpicture}
      \node at (1,0) (tlblank) {};
      \node at (0,0) (tcblank) {};
      \node at (-1,0) (trblank) {};
      \node at (1,-2) (mM) {$m^N$};
      \node at (0,-1) (deltaN) {$\delta^{N,L}$};
      \node at (1,-3) (blblank) {};
      \node at (0,-3) (bcblank) {};
      \node at (-1,-3) (brblank) {};
      \draw[Dmodar] (tlblank) to (mM);
      \draw[Dmodar] (mM) to (blblank);
      \draw[Amodar] (tcblank) to (deltaN);
      \draw[Amodar] (deltaN) to (bcblank);
      \draw[Dmodar] (trblank) to (brblank);
      \draw[tensoralgarrow] (deltaN) to (mM);
    \end{tikzpicture}
  }
  \]
  respectively, independently of which way one associates. (In the
  separated type \DD\ module, the map $\delta^L$ is defined from
  $\delta^{1L}$ in the same way $\delta$ is defined from
  $\delta^1$. The fact that these are the only two terms in the
  differential follows from strict unitality of the modules.)  The
  same holds for bimodules, with slightly extended diagrams.
\end{proof}

For $N$ a non-separated type \DD, the box tensor product is not strictly
associative as in Lemma~\ref{lem:tens-assoc-1-strict}. However,
associativity does hold up to homotopy equivalence:
\begin{proposition}\label{prop:tens-assoc-2-weak} Let $\Alg$ and $\Blg$ be
  \dg algebras and $M_\Alg$, $\lsup{\Alg}N^\Blg$ and $\lsub{\Blg}P$ be
  right type $A$, type \DD, and left type $A$ structures respectively. 
  Suppose moreover that $\lsup{\Alg}N^\Blg$ is homotopy equivalent to a 
  bounded type \DD\ structure.
  Then
  $(M_\Alg\DT \lsup{\Alg}N^\Blg) \DT \lsub{\Blg}P$ is
  homotopy equivalent to $M_\Alg\DT (\lsup{\Alg}N^\Blg \DT
  \lsub{\Blg}P)$. The analogous statements hold if $M$ is a type \AAm\ or
  \DA\ module and/or $P$ is a type \AAm\ or \AD\ module.
\end{proposition}
We will prove this in Section~\ref{sec:bar-cobar-modules}, after
introducing the bar resolution.

\subsubsection{Bar resolutions of modules}
\label{sec:bar-cobar-modules}
\begin{definition}
  For $\Alg$ a \dg algebra, $\lsupv{\Alg}\Barop(\Alg)^\Alg$ is the type \DD\
  bimodule with underlying $\Groundk$-module
  $\Tensor^*(\Alg[1])$, with basis written $[a_1|\cdots|a_k]$ for $k
  \ge 0$, and structure maps
  \begin{multline*}
    \delta^1[a_1|\cdots|a_k] \coloneqq a_1\otimes[a_2|\cdots|a_k]\otimes 1
      + 1\otimes [a_1|\cdots|a_{k-1}]\otimes a_k\\
      +\sum_{1\le i\le k}1\otimes
      [a_1|\cdots|\mu_1(a_i)|\cdots|a_k]\otimes 1 +
      \sum_{1\le i\le
        k-1}1\otimes [a_1|\cdots|\mu_2(a_i,a_{i+1})|\cdots|a_k]\otimes 1.
  \end{multline*}
\end{definition}

Note that $\lsupv{\Alg}\Barop(\Alg)^\Alg$ is bounded as a type \DD\ structure.

We can use the bar resolution to define the tensor product of
$\Ainf$-modules:
\begin{definition}\label{def:DTP} Given $\Ainf$-modules $M_\Alg$ and
  $\lsub{\Alg}N$ over a \dg algebra~$\Alg$,
  define the \emph{$\Ainf$-tensor product of $M$ and $N$} to be
  \[
  M_\Alg \DTP \lsub{\Alg} N \coloneqq M_\Alg \DT
  \lsupv{\Alg}\Barop(\Alg)^\Alg \DT \lsub{\Alg} N.
  \]
\end{definition}
The module $\lsupv{\Alg}\Barop(\Alg)^\Alg$ is separated, so by
Lemma~\ref{lem:tens-assoc-1-strict} it is okay that we have not
parenthesized the triple box product.  This definition agrees with the
standard definition (e.g., \cite[Section 6.3]{AinftyAlg}).

\begin{proposition}
  \label{prop:D-to-A-and-back}
  Let $\Alg$ be a \dg algebra.  The map $\lsupv{\Alg}N \mapsto
  \lsub{\Alg}{\Alg}_{\Alg}\DT \lsupv{\Alg}N$ induces an $\Ainf$-functor
  $\lsupv{\Alg}\ModCat\to \lsub{\Alg}\ModCat$.  Similarly, the map
  $M_{\Alg}\mapsto M_{\Alg}\DT \lsupv{\Alg}\Barop(\Alg)^\Alg$ induces
  an $\Ainf$-functor $\lsub{\Alg}\ModCat\to
  \lsupv{\Alg}\ModCat$. These two functors are homotopy inverses to
  one another (and the homotopy is canonical). They entwine the tensor
  products $\DT$ and $\DTP$ in the sense that there is a canonical
  homotopy equivalence
  \[
  M_\Alg\DT\lsup{\Alg}N\simeq (M_\Alg)\DTP(\lsub{\Alg}\Alg_\Alg \DT
  \lsup{\Alg}N).
  \]

  In particular, the categories $\lsupv{\Alg}\ModCat$ and
  $\lsub{\Alg}\ModCat$ are quasi-equivalent, and hence their derived
  categories are equivalent.

  Corresponding statements hold for the categories of type \DD, \DA,
  and \AAm\ modules.
\end{proposition}

We will prove Proposition~\ref{prop:D-to-A-and-back} presently. The
proposition justifies the following abuse of notation: given a type
$D$ structure $\lsup{\Alg}M$, let
$\lsub{\Alg}M=\lsub{\Alg}\Alg_\Alg\DT \lsup{\Alg}M$. Similarly, given
a type $A$ module $\lsub{\Alg}N$, let
$\lsup{\Alg}N=\lsupv{\Alg}\Barop(\Alg)^\Alg\DT \lsub{\Alg}N$. This notation
extends in an obvious way to bimodules. The statement
in Proposition~\ref{prop:D-to-A-and-back} about tensor products
becomes
\[
M_\Alg\DT\lsup{\Alg}N\simeq M_\Alg\DTP\lsub{\Alg}N.
\]

The proof of Proposition~\ref{prop:D-to-A-and-back}
is based on the following key lemma:
\begin{lemma}\label{lem:bar-res}
  For any \dg algebra $\Alg$, the type \DA\ module
  $\lsupv{\Alg}\Barop(\Alg)^\Alg \DT \lsub{\Alg}\Alg_\Alg$ is homotopy
  equivalent to $\lsupv{\Alg}[\Id]_\Alg$.
  Moreover, the homotopy equivalence  
  $\kappa\co \lsupv{\Alg}\Barop(\Alg)^\Alg \DT \lsub{\Alg}\Alg_\Alg
  \rightarrow\lsupv{\Alg}[\Id]_\Alg$
  is bounded.
\end{lemma}
(See Example~\ref{eg:Id-AA-mod} and Definition~\ref{def:rank-1-DA-mods} for
definitions of $\lsub{\Alg}\Alg_\Alg$ and $\lsupv{\Alg}[\Id]_\Alg$.
Also, the homotopy inverse to $\kappa$ is not necessarily bounded.)

\begin{proof}
  This lemma is a version of the standard fact that the bar
  resolution is a resolution.  (See \cite[Proposition~\ref*{LOT:prop:BarResolution}]{LOT1} for a
  version that is not far from the one we give below.)  We
  translate the proof
  into our language.  For convenience, write $\lsup{\Alg}M_\Alg$ or
  just $M$ for
  $\lsupv{\Alg}\Barop(\Alg)^\Alg \DT \lsub{\Alg}\Alg_\Alg$.  Now
  define $\kappa \co M \to \Id$ by
  \[
  \kappa^1_{1+l}([a_1|\cdots|a_k]b, c_1, \ldots, c_l) \coloneqq
  \begin{cases}
    0 & k > 0 \text{ or } l > 0\\
    b\otimes 1 & k = l = 0.
  \end{cases}
  \]
  Define $\psi \co \Id \to M$ by
  \[
  \psi^1_{1+l} (1, a_1, \ldots, a_l) \coloneqq 1 \otimes  [a_1|\cdots|a_l]1.
  \]
  It is elementary to check that $\kappa$ and $\psi$ are homomorphisms of
  type \DA\ structures (i.e., cycles in their respective morphism spaces)
  and that $\kappa \circ \psi$ is the identity.
  The other composition $\psi \circ \kappa$ is not the identity, but it
  is homotopic to the identity by the map $h \co M \to
  M$ defined by
  \[
  h^1_{1+l}([a_1|\cdots|a_k]b, c_1,\ldots, c_l) \coloneqq
    1\otimes [a_1|\cdots|a_k|b|c_1|\cdots|c_l]1. \qedhere
  \]
\end{proof}

\begin{lemma}
  \label{lem:IdentityDAIsIdentity}
  Let $\Alg$ be a \dg algebra and $M_\Alg$ an $\Ainf$-module
  over $\Alg$. Then $M_\Alg\DT\lsupv{\Alg}[\Id]_\Alg$ is canonically
  isomorphic to
  $M_\Alg$. Similarly, if $\lsup{\Alg}N$ is a type $D$ structure then
  $\lsupv{\Alg}[\Id]_\Alg\DT\lsup{\Alg}N$ is canonically isomorphic to
  $\lsup{\Alg}N$. 
  Similar statements hold when $M$ is a type \AAm\ or \DA\
  module and when $N$ is a type \DA\ or \DD\ module.
\end{lemma}
\begin{proof}
  This is immediate from the definitions.
\end{proof}
\begin{proof}[Proof of Proposition~\ref{prop:D-to-A-and-back}]
  The fact that $\lsub{\Alg}{\Alg}_{\Alg}\DT \cdot$ and
  $\lsupv{\Alg}\Barop(\Alg)^\Alg\DT\cdot$ are functorial is part of
  Lemma~\ref{lem:DT-descends-bimod}. To see that these two
  functors are homotopy inverses, note that  
  \[
  \lsub{\Alg}{\Alg}_{\Alg}\DT(\lsupv{\Alg}\Barop(\Alg)^\Alg\DT\cdot)\cong
  (\lsub{\Alg}{\Alg}_{\Alg}\DT\lsupv{\Alg}\Barop(\Alg)^\Alg)\DT\cdot\simeq \lsub{\Alg}[\Id]^{\Alg}\DT\cdot
  \]
  where the isomorphism uses Proposition~\ref{lem:tens-assoc-1-strict}
  and the homotopy equivalence uses Lemmas~\ref{lem:bar-res}
  and~\ref{lem:DT-descends-bimod}. But
  $\lsub{\Alg}[\Id]^{\Alg}\DT\cdot$ is an equivalence of
  categories. Similar reasoning applies to the composition
  $\lsupv{\Alg}\Barop(\Alg)^\Alg\DT(\lsub{\Alg}{\Alg}_{\Alg}\DT\cdot)$, proving
  the result.

  The corresponding statements about bimodules follow similarly. The
  fact that the functors intertwine $\DT$ and $\DTP$ is obvious from
  the definitions.
\end{proof}

\begin{definition}
  Given an $\Ainf$ module $M_\Alg$ over a \dg algebra define its
  \emph{bar resolution} to be
  \[
  \Barop(M_\Alg)=M_\Alg\DT\lsupv{\Alg}\Barop(\Alg)^\Alg\DT\lsub{\Alg}\Alg_\Alg.
  \]
  Similarly, if $\lsub{\Alg}M_\Blg$ is a type \AAm\  structure then
  define its \emph{bar resolution} to be
  \[
  \Barop(\lsub{\Alg}M_\Blg)=\lsub{\Alg}\Alg_\Alg\DT\lsupv{\Alg}\Barop(\Alg)^\Alg
  \DT\lsub{\Alg}M_\Blg\DT\lsupv{\Blg}\Barop(\Blg)^\Blg\DT\lsub{\Blg}\Blg_\Blg
  \]
\end{definition}

By taking bar
resolutions, over a \dg algebra, every $\Ainf$-module is
$\Ainf$-homotopy equivalent to an honest \dg module:
\begin{proposition}\label{prop:Bar-htpy-equiv}
  Let $\Alg$ be a \dg algebra and $M_\Alg$ an $\Ainf$-module
  over $\Alg$. Then:
  \begin{enumerate}
  \item $\Barop(M_\Alg)$ is $\Ainf$-homotopy equivalent to $M_\Alg$ and
  \item $\Barop(M_\Alg)$ is an honest \dg module.
  \end{enumerate}
  Similarly, if $\Blg$ is another \dg algebra and $\lsub{\Alg}M_\Blg$
  is a type \AAm\  structure then:
  \begin{enumerate}
  \item $\Barop(\lsub{\Alg}M_\Blg)$ is $\Ainf$-homotopy equivalent to
    $\lsub{\Alg}M_\Blg$ and
  \item $\Barop(\lsub{\Alg}M_\Blg)$ is an honest \dg module.
  \end{enumerate}
\end{proposition}

\begin{proof}
  This is a combination of Lemmas~\ref{lem:bar-res}
  and~\ref{lem:IdentityDAIsIdentity}.  Lemma~\ref{lem:bar-res}
  furnishes a homotopy equivalence between left-bounded type \DA\ 
  modules $\kappa\co\lsupv{\Alg}\Barop(\Alg)^\Alg \DT \lsub{\Alg}\Alg_\Alg\rightarrow
  \lsupv{\Alg}[\Id]_\Alg$. This induces a homotopy equivalence
  $$\Id_M\DT\kappa\co M_{\Alg}\DT \lsupv{\Alg}\Barop(\Alg)^\Alg \DT \lsub{\Alg}\Alg_\Alg\rightarrow
  M_{\Alg}\DT \lsupv{\Alg}[\Id]_{\Alg},$$ while
  Lemma~\ref{lem:IdentityDAIsIdentity} furnishes the isomorphism
  $M_\Alg\cong M_\Alg\DT\lsupv{\Alg}[\Id]_\Alg$.
\end{proof}

It is immediate from
Lemma~\ref{lem:bar-res} that:
\begin{corollary} If $\Alg$ is a \dg algebra, then
  \begin{itemize}
  \item $M_\Alg \DTP \lsub{\Alg} \Alg_\Alg$ is chain-homotopy
    equivalent to $M_\Alg$;
  \item $(M_\Alg \DTP \lsub{\Alg}N_\Blg) \DTP \lsub{\Blg}P = M_\Alg
    \DTP (\lsub{\Alg}N_\Blg \DTP \lsub{\Blg}P)$; and
  \item $M_\Alg \DTP (\lsub{\Alg}\Alg_\Alg \DT \lsup{\Alg}N) \simeq M_\Alg
    \DT \lsup{\Alg}N$.
  \end{itemize}
\end{corollary}

Thus, for instance, $\DTP$ turns $\HMod(\lsub{\Alg}\ModCat_\Alg)$ into a
monoidal category.

The bar resolution can be used to give an alternate characterization
of the subcategory $\DModCat{\Alg}\subset \DuModCat{\Alg}$:

\begin{proposition}
  \label{prop:CharacterizeDMod}
  The following conditions on a type $D$ structure are equivalent:
  \begin{enumerate}[label=(D-\arabic*),ref=D-\arabic*]
    \item\label{item:D-homotopy-bdd} $\lsupv{\Alg}N$ is homotopy
      equivalent to a bounded
      type $D$ structure.
    \item\label{item:D-homotopy-bar} The canonical map 
      $\lsupv{\Alg}\Barop(\Alg)^\Alg\DT\lsub{\Alg}{\Alg}_{\Alg}\DT\lsupv{\Alg}N
      \rightarrow \lsupv{\Alg}N$ (induced by combining maps from
      Lemma~\ref{lem:bar-res} and~\ref{lem:IdentityDAIsIdentity})
      is a homotopy equivalence.
  \end{enumerate}
  Analogously, the following are equivalent for a type \DA\ structure:
  \begin{enumerate}[label=(DA-\arabic*),ref=DA-\arabic*]
    \item $\lsupv{\Alg}N_{\Blg}$ is homotopy equivalent to a left bounded
      type \DA\ structure.
    \item $\lsupv{\Alg}N_{\Blg}$ is homotopy equivalent to a bounded
      type \DA\ structure.
    \item The canonical map 
      $\lsupv{\Alg}\Barop(\Alg)^{\Alg}\DT \lsub{\Alg}\Alg_{\Alg}\DT 
      \lsupv{\Alg}N_{\Blg}\rightarrow 
      \lsupv{\Alg}N_{\Blg}$
      is a homotopy equivalence.
    \item The canonical map 
      $\lsupv{\Alg}\Barop(\Alg)^{\Alg}\DT \lsub{\Alg}\Alg_{\Alg}\DT 
      \lsupv{\Alg}N_{\Blg}\DT \lsupv{\Blg}\Barop(\Blg)^{\Blg}
      \DT \lsub{\Blg}\Blg_{\Blg}\
      \rightarrow 
      \lsupv{\Alg}N_{\Blg}$
      is a homotopy equivalence.
  \end{enumerate}
  The following conditions on a type \DD\ bimodule are equivalent
  \begin{enumerate}[label=(DD-\arabic*),ref=DD-\arabic*]
    \item $\lsupv{\Alg}N^{\Blg}$ is homotopy equivalent to a bounded
      type \DD\ structure.
    \item The canonical map 
      \[\lsupv{\Alg}\Barop(\Alg)^\Alg\DT\lsub{\Alg}{\Alg}_{\Alg}\DT\lsupv{\Alg}N^{\Blg}
      \DT \lsub{\Blg}{\Blg}_{\Blg}\DT\lsupv{\Blg}\Barop(\Blg)^\Blg
      \rightarrow \lsupv{\Alg}N^{\Blg}\]
      is a homotopy equivalence.
  \end{enumerate}
\end{proposition}

\begin{proof}
  We start with the case of a type $D$ module $\lsupv{\Alg}N$.

  $\eqref{item:D-homotopy-bdd}\Rightarrow\eqref{item:D-homotopy-bar}$. Let
  $\lsupv{\Alg}N'$ be a bounded type $D$
  structure, and $\phi\co N\longrightarrow N'$ be a homotopy
  equivalence.  According to Lemma~\ref{lem:DT-descends-bimod}, we
  have a homotopy commutative diagram. 
  
  \[
  \begin{tikzpicture}
    \node at (0,0) (tlcorn) {$\lsupv{\Alg}\Barop(\Alg)^{\Alg}\DT \lsub{\Alg}\Alg_{\Alg}\DT \lsupv{\Alg}N$};
    \node at (5,0) (tmid) {$\lsupv{\Alg}[\Id]_{\Alg}\DT \lsupv{\Alg}N$};
    \node at (10,0) (trcorn) {$\lsupv{\Alg}N$};
    \node at (0,-2) (blcorn)  {$\lsupv{\Alg}\Barop(\Alg)^{\Alg}\DT \lsub{\Alg}\Alg_{\Alg}\DT \lsupv{\Alg}N'$};
    \node at (5,-2) (bmid) {$\lsupv{\Alg}[\Id]_{\Alg}\DT \lsupv{\Alg}N'$};
    \node at (10,-2) (brcorn) {$\lsupv{\Alg}N'$};
    \draw[->] (tlcorn) to node[above]{$\kappa\DT\Id_N$} (tmid);
    \draw[->] (tlcorn) to node[left]{$\Id_{\Barop(\Alg)\DT \Alg}\DT \phi$} (blcorn);
    \draw[->] (trcorn) to node[right]{$\phi$} (brcorn);
    \draw[->] (blcorn) to node[below]{$\kappa\DT\Id_{N'}$} (bmid);
    \draw[->] (tmid) to node[above]{$\cong$} (trcorn);
    \draw[->] (bmid) to node[below]{$\cong$} (brcorn);
  \end{tikzpicture}
  \]
  All arrows with the possible exception of $\kappa\DT\Id_N$ are homotopy equivalences.
  It follows that $\kappa\DT\Id_N$ is a homotopy equivalence, as well.

  $\eqref{item:D-homotopy-bar}\Rightarrow\eqref{item:D-homotopy-bdd}$. Since
  $\lsub{\Alg}\Alg_{\Alg}$ is a bounded
  type \AAm\ structure, we can form $\lsupv{\Alg}\Barop(\Alg)^{\Alg}\DT
  \lsub{\Alg}\Alg_{\Alg}\DT \lsupv{\Alg}N$.  Moreover, since
  $\lsupv{\Alg}\Barop(\Alg)^{\Alg}$ is a bounded type \DD\ structure,
  $\lsupv{\Alg}\Barop(\Alg)^{\Alg}\DT \lsub{\Alg}\Alg_{\Alg}\DT
  \lsupv{\Alg}N$ is a bounded type $D$ structure.

  The type \DA\ case is similar. The type \DD\ case follows from the type $D$ case.
\end{proof}

\begin{corollary}\label{cor:tensor-projective}
  Let $\Alg$ be a \dg algebra. Suppose $\lsup{\Alg}M$ is homotopy
  equivalent to a bounded type~$D$ structure (i.e., $\lsup{\Alg}M$ is
  an object in $\lsupv{\Alg}\ModCat$). Then $\lsub{\Alg}\Alg_\Alg \DT
  \lsup{\Alg}M$ is a projective $\Alg$-module.
\end{corollary}

\begin{proof}
  By Lemma~\ref{lem:DT-descends-bimod} and
  Proposition~\ref{prop:CharacterizeDMod},
  $\lsub{\Alg}\Alg_\Alg\DT\lsup{\Alg}M$ is homotopy equivalent to its
  bar resolution. But the bar resolution of any module is
  projective~\cite[Proposition 10.12.2.6]{BernsteinLunts94:EquivariantSheaves}, and
  projectivity is preserved by homotopy equivalences.
\end{proof}

\begin{remark}
  \label{rmk:NeedBoundedness}
  The condition that $\lsup{\Alg}M$ is homotopy equivalent to a
  bounded type $D$ structure is essential.  For
  instance, let $\Alg=\Field[t]/t^2$, and $\lsupv{\Alg}M$ have one
  generator $x$ with $\delta^1 x=t\otimes x$, then
  $\lsub{\Alg}M=\lsub{\Alg}\Alg_\Alg \DT \lsup{\Alg}M$ is acyclic,
  but, since $(\Alg/t)\otimes \lsub{\Alg}M$ has homology, we can
  conclude that $\lsub{\Alg}M$ is not homotopy equivalent to the
  trivial module. Thus, $\lsub{\Alg}M$ is not projective.

  Note, in particular, that this implies that not every type $D$
  structure is homotopy equivalent to a bounded type $D$ structure.
\end{remark}

Finally, we turn to the proof of Proposition~\ref{prop:tens-assoc-2-weak},
which is based on the following observation:
\begin{lemma}
  If $\lsupv{\Alg}N^{\Blg}$ is homotopy equivalent 
  to a bounded type \DD\ structure, then it is 
  homotopy equivalent to a separated one.
\end{lemma}
\begin{proof}
  Let $\lsup{\Alg}M^\Blg$ be a type \DD\ module. Then
  $\lsupv{\Alg}\Barop(\Alg)^\Alg\DT\lsub{\Alg}\Alg_\Alg\DT\lsup{\Alg}M^\Blg\DT
  \lsub{\Blg}{\Blg}_{\Blg}\DT
  \lsupv{\Blg}\Barop(\Blg)^{\Blg}$
  is a separated type \DD\ module, and is homotopy equivalent to
  $\lsup{\Alg}M^\Blg$ by Proposition~\ref{prop:CharacterizeDMod}.
\end{proof}
\begin{proof}[Proof of Proposition~\ref{prop:tens-assoc-2-weak}]
  Let $\lsup{\Alg}{\widetilde{N}}^\Blg$ be a separated type \DD\ module
  homotopy equivalent to $\lsup{\Alg}N^\Blg$. Then
  \[
  (M_\Alg\DT \lsup{\Alg}N^\Blg) \DT \lsub{\Blg}P\simeq (M_\Alg\DT
  \lsup{\Alg}{\widetilde{N}}^\Blg) \DT \lsub{\Blg}P \simeq M_\Alg\DT
  (\lsup{\Alg}{\widetilde{N}}^\Blg \DT \lsub{\Blg}P) \simeq M_\Alg\DT
  (\lsup{\Alg}N^\Blg \DT \lsub{\Blg}P)
  \]
  where the outside equivalences use Lemma~\ref{lem:DT-descends-bimod}
  and the middle equivalence uses Lemma \ref{lem:tens-assoc-1-strict}.
\end{proof}

\begin{remark}\label{remark:bar-in-general}
  The results of this section generalize in a transparent way to the
  case of $\Ainf$-algebras; the only obstruction is that we have not
  defined the bar resolution except over \dg algebras---an obstruction
  which is more terminological than mathematical.  The key observation
  is that the difficulties mentioned in Remark~\ref{rem:DD-difficult}
  do not arise for separated \DD\ modules like $\lsup{\Alg}M^\Alg$.
\end{remark}

\subsubsection{More \textalt{$\Mor$}{Mor}}
\label{sec:mod-hom}

Recall from Section~\ref{sec:forgetful-functors} that there are
functors $\mathcal{F}$ which forget one of the actions on a bimodule.
\begin{definition}\label{def:bimod-on-Hom}
  Define $\Mor_{\Blg}(\lsub{\Alg}M_\Blg,\lsub{\Clg}N_\Blg)$ to be the
  complex $\Mor(\mathcal{F}(M)_\Blg,\mathcal{F}(N)_\Blg)$. We define a
  $(\Clg,\Alg)$-bimodule module structure on
  $\Mor_{\Blg}(\lsub{\Alg}M_\Blg,\lsub{\Clg}N_\Blg)$ as follows. The
  operation $m_{i,1,j}=0$ unless either $i=0$ or $j=0$. The operation
  $m_{0,1,0}$ is, of course, the differential on the $\Mor$ complex.
  Finally, for exactly one of $i$ and $j$ nonzero, define
  \[
  \begin{split}
    m_{i,1,0}&(c_1,\dots,c_i,f)_{k+1}(x,b_1,\dots,b_k)\\
    &= \sum_{p+q=k}
    m^N_{i,1,q}\left(c_1,\dots,c_i,f_{p+1}\left(x,b_{1},\dots,b_{p}\right),b_{p+1},\dots,b_k\right)
  \end{split}
  \]
  \[
  \begin{split}
    m_{0,1,j}&(f,a_c,\dots,a_j)_{k+1}(x,b_1,\dots,b_k)\\
    &= \sum_{p+q=k}
    f_{q+1}\left(m_{j,1,p}^M(a_1,\dots,a_j,x,b_1,\dots,b_p),b_{p+1},\dots,b_{k}\right).
  \end{split}
  \]
  
  Graphically, if we draw $f$ as
  \[
  \begin{tikzpicture}
    \node at (0,0) (tcblank) {};
    \node at (1,0) (trblank) {};
    \node at (0,-1) (f) {$f$};
    \node at (0,-2) (bcblank) {};
    \draw[Amodar] (tcblank) to node[below,sloped]{\lab{M}} (f);
    \draw[Amodar] (f) to node[below,sloped]{\lab{N}} (bcblank);
    \draw[tensorblgarrow] (trblank) to node[below,sloped]{\lab{T^*B}} (f);
  \end{tikzpicture}
  \]
  then the module structure on
  $\Mor_{\Blg}(\lsub{\Alg}M_\Blg,\lsub{\Clg}N_\Blg)$ is given as in Figure~\ref{fig:mor-module-structure}.
  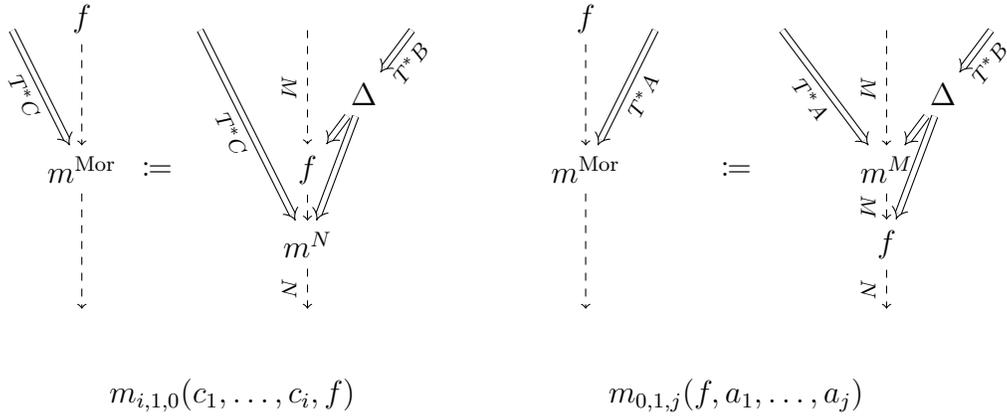
\begin{figure}
    \centering
    \[
    \mathcenter{
      \begin{tikzpicture}
        \node at (-3,0) (ftcblank) {$f$};
        \node at (-3,-2) (action) {$m^{\Mor}$};
        \node at (-3,-4) (fbcblank) {};
        \node at (-4,0) (cfs)  {}; %{$c_1\otimes\cdots\otimes c_i$};
        \node at (-2,-2) (equals) {$\coloneqq$};
%        \draw[Amodar] (ftcblank) to (fbcblank);
        \draw[Amodar] (ftcblank) to (action);
        \draw[Amodar] (action) to (fbcblank);
        \draw[tensorclgarrow] (cfs) to node[below,sloped]{\lab{T^*C}} (action);
        \node at (0,0) (tlblank) {}; 
        \node at (1.5,0) (trblank) {};
        \node at (-1.5,0) (tllblank) {};% {$c_1\otimes\cdots\otimes c_i$};
        \node at (0,-2) (f) {$f$}; 
        \node at (0,-3) (mN) {$m^N$}; 
        \node at (0,-4) (blblank) {};
        \node at (.75,-1) (Delta) {$\Delta$};
        \draw[tensorclgarrow] (tllblank) to node[below,sloped]{\lab{T^*C}} (mN); 
        \draw[tensorblgarrow] (trblank) to node[below,sloped]{\lab{T^*B}} (Delta);
        \draw[tensorblgarrow] (Delta) to (f); 
        \draw[tensorblgarrow] (Delta) to (mN); 
        \draw[Amodar] (tlblank) to node[below,sloped]{\lab{M}} (f); 
        \draw[Amodar] (f) to (mN); 
        \draw[Amodar] (mN) to node[below,sloped]{\lab{N}} (blblank);
        \node at (-1,-5) (label) {$m_{i,1,0}(c_1,\dots, c_i,f)$};
      \end{tikzpicture}
    } 
    \qquad \quad\mathcenter{
      \begin{tikzpicture}
        \node at (-4,0) (ftcblank) {$f$};
        \node at (-4,-2) (action) {$m^{\Mor}$};
        \node at (-4,-4) (fbcblank) {};
        \node at (-3,0) (cfs)  {}; %{$a_1\otimes\cdots\otimes a_i$};
        \node at (-2,-2) (equals) {$\coloneqq$};
        \draw[Amodar] (ftcblank) to (action);
        \draw[Amodar] (action) to (fbcblank);
%        \draw[Amodar] (ftcblank) to (fbcblank);
        \draw[tensorclgarrow] (cfs) to node[below,sloped]{\lab{T^*A}} (action);
        \node at (0,0) (tlblank) {};
        \node at (-1.5,0) (tcblank) {};% {$a_1\otimes\cdots\otimes a_j$};
        \node at (1.5,0) (trblank) {};
        \node at (0,-2) (mM) {$m^M$};
        \node at (0,-3) (f) {$f$};
        \node at (0,-4) (blblank) {};
        \node at (.75,-1) (Delta) {$\Delta$}; 
        \draw[tensorblgarrow] (trblank) to node[below,sloped]{\lab{T^*B}} (Delta);
        \draw[tensorblgarrow] (Delta) to (mM); 
        \draw[tensorblgarrow] (Delta) to (f);
        \draw[tensoralgarrow] (tcblank) to
        node[below,sloped]{\lab{T^*A}} (mM); 
        \draw[Amodar] (tlblank) to node[below,sloped]{\lab{M}} (mM); 
        \draw[Amodar] (mM) to node[below,sloped]{\lab{M}} (f);
        \draw[Amodar] (f) to node[below,sloped]{\lab{N}} (blblank);
        \node at (-2,-5) (label) {$m_{0,1,j}(f,a_1,\dots, a_j)$};
      \end{tikzpicture}
    }
    \]
    
    \caption{\textbf{Module structure on the $\Mor$-complex.}}
\label{fig:mor-module-structure}
\end{figure}
\end{definition}
\begin{lemma}
  The structure defined in Definition~\ref{def:bimod-on-Hom} satisfies
  the $\Ainf$-bimodule relation (Equation~\eqref{eq:AAbimodule-relation}).
\end{lemma}
\begin{proof}
  This is a straightforward, if tedious, verification. It is clear
  from the diagrams that the left and right actions commute. The
  $\Ainf$-relation for the right action involves the following terms:
  \[
  \mathcenter{
    \begin{tikzpicture}
      \node at (0,0) (tlblank) {};
      \node at (-1,0) (tcblank) {};
      \node at (2,0) (trblank) {};
      \node at (0,-2) (mM) {$m^M$};
      \node at (0,-3) (mM2) {$m^M$};
      \node at (0,-4) (f) {$f$};
      \node at (0,-5) (blblank) {};
      \node at (-1,-1) (turn) {$\Delta$};
      \node at (1,-1) (Delta) {$\Delta$};
      \draw[tensorblgarrow] (trblank) to node[below,sloped]{\lab{T^*B}} (Delta);
      \draw[tensorblgarrow] (Delta) to (mM);
      \draw[tensorblgarrow] (Delta) to (mM2);
      \draw[tensorblgarrow, bend left=5] (Delta) to (f);
      \draw[tensoralgarrow] (tcblank) to node[above,sloped]{\lab{T^*A}} (turn);
      \draw[tensoralgarrow] (turn) to (mM);
      \draw[tensoralgarrow] (turn) to (mM2);
      \draw[Amodar] (tlblank) to node[above,sloped]{\lab{M}} (mM);
      \draw[Amodar] (mM) to node[below,sloped]{} (mM2);
      \draw[Amodar] (mM2) to node[below,sloped]{} (f);
      \draw[Amodar] (f) to node[below,sloped]{\lab{N}} (blblank);
    \end{tikzpicture}
  }\qquad
  \mathcenter{
    \begin{tikzpicture}
      \node at (0,0) (tlblank) {};
      \node at (-1,0) (tllblank) {};
      \node at (3,0) (trblank) {};
      \node at (0,-2) (mM) {$m^M$};
      \node at (0,-3) (f) {$f$};
      \node at (0,-4) (blblank) {};
      \node at (2,-1) (Delta) {$\Delta$};
      \node at (1,-2) (DB) {$\overline{D}^\Blg$};
      \draw[tensorblgarrow] (trblank) to node[below,sloped]{\lab{T^*B}} (Delta);
      \draw[tensorblgarrow] (Delta) to (mM);
      \draw[tensorblgarrow] (Delta) to (DB);
      \draw[tensorblgarrow] (DB) to (f);
      \draw[tensoralgarrow] (tllblank) to node[above,sloped]{\lab{T^*A}} (mM);
      \draw[Amodar] (tlblank) to node[above,sloped]{\lab{M}} (mM);
      \draw[Amodar] (mM) to node[below,sloped]{} (f);
      \draw[Amodar] (f) to node[below,sloped]{\lab{N}} (blblank);
    \end{tikzpicture}
  }\qquad
  \mathcenter{
    \begin{tikzpicture}
      \node at (0,0) (tlblank) {};
      \node at (-1,0) (tllblank) {};
      \node at (3,0) (trblank) {};
      \node at (2,-.5) (DB) {$\overline{D}^\Blg$};
      \node at (0,-2) (mM) {$m^M$};
      \node at (0,-3) (f) {$f$};
      \node at (0,-4) (blblank) {};
      \node at (1,-1) (Delta) {$\Delta$};
      \draw[tensorblgarrow] (trblank) to node[below,sloped]{\lab{T^*B}} (DB);
      \draw[tensorblgarrow] (DB) to (Delta);
      \draw[tensorblgarrow] (Delta) to (mM);
      \draw[tensorblgarrow] (Delta) to (f);
      \draw[tensoralgarrow] (tcblank) to node[above,sloped]{\lab{T^*A}} (mM);
      \draw[Amodar] (tlblank) to node[above,sloped]{\lab{M}} (mM);
      \draw[Amodar] (mM) to node[below,sloped]{} (f);
      \draw[Amodar] (f) to node[below,sloped]{\lab{N}} (blblank);
    \end{tikzpicture}
  }
  \]
  \[
  2\cdot 
  \mathcenter{
    \begin{tikzpicture}
      \node at (0,0) (tlblank) {};
      \node at (-1,0) (tllblank) {};
      \node at (2,0) (trblank) {};
      \node at (0,-2) (mM) {$m^M$};
      \node at (0,-3) (f) {$f$};
      \node at (0,-4) (mN) {$m^N$};
      \node at (0,-5) (blblank) {};
      \node at (1,-1) (Delta) {$\Delta$};
      \draw[tensorblgarrow] (trblank) to node[below,sloped]{\lab{T^*B}} (Delta);
      \draw[tensorblgarrow] (Delta) to (mM);
      \draw[tensorblgarrow] (Delta) to (f);
      \draw[tensorblgarrow, bend left=5] (Delta) to (mN);
      \draw[tensoralgarrow] (tllblank) to node[above,sloped]{\lab{T^*A}} (mM);
      \draw[Amodar] (tlblank) to node[above,sloped]{\lab{M}} (mM);
      \draw[Amodar] (mM) to node[below,sloped]{} (f);
      \draw[Amodar] (f) to node[below,sloped]{\lab{N}} (mN);
      \draw[Amodar] (mN) to (blblank);
    \end{tikzpicture}
  }\qquad
  \mathcenter{
    \begin{tikzpicture}
      \node at (0,0) (tlblank) {};
      \node at (-1,0) (tllblank) {};
      \node at (2,0) (trblank) {};
      \node at (0,-2) (mM) {$m^M$};
      \node at (0,-3) (f) {$f$};
      \node at (0,-4) (blblank) {};
      \node at (-1,-1) (turn) {$\overline{D}^\Alg$};
      \node at (1,-1) (Delta) {$\Delta$};
      \draw[tensorblgarrow] (trblank) to node[below,sloped]{\lab{T^*B}} (Delta);
      \draw[tensorblgarrow] (Delta) to (mM);
      \draw[tensorblgarrow] (Delta) to (f);
      \draw[tensoralgarrow] (tllblank) to
      node[above,sloped]{\lab{T^*A}} (turn);
      \draw[tensoralgarrow] (turn) to (mM);
      \draw[Amodar] (tlblank) to node[below,sloped]{\lab{M}} (mM);
      \draw[Amodar] (mM) to node[below,sloped]{} (f);
      \draw[Amodar] (f) to node[below,sloped]{\lab{N}} (blblank);
    \end{tikzpicture}
  }.
  \]
  Here, the first four terms come from $m(m(f,a_1,\dots,
  a_j),a_{j+1},\dots, a_k)$; the first one is the generic
  case, the second is part of the $j=0$ case, the third is part of the $j=k$
  case, and the fourth occurs in both the $j=0$ and $j=k$ cases. The
  fifth term comes from $m(f,\overline{D}^\Alg(a_1,\dots,
  a_k))$.  Applying the $\Ainf$-relation for $M$, these terms
  cancel. 

  A similar argument applies to the left action.
\end{proof}

Bimodule morphisms $f\co\lsub{\Alg}M_\Blg\to \lsub{\Alg}M'_\Blg$ and
$g\co \lsub{\Clg}N_\Blg\to\lsub{\Clg}N'_\Blg$ induce $\Ainf$-maps
\begin{align*}
  f^*&\co \Mor_\Blg(\lsub{\Alg}M'_\Blg,\lsub{\Clg}N_\Blg)\to
  \Mor_\Blg(\lsub{\Alg}M_\Blg,\lsub{\Clg}N_\Blg)\qquad\text{ and}\\
  g_*&\co \Mor_\Blg(\lsub{\Alg}M_\Blg,\lsub{\Clg}N_\Blg)\to
  \Mor_\Blg(\lsub{\Alg}M_\Blg,\lsub{\Clg}N'_\Blg)
\end{align*}
as follows. Define $f^*_{i,1,j}$ to be zero if $i>0$ and define
\[f^*\co \Mor_\Blg(\lsub{\Alg}M'_\Blg,\lsub{\Clg}N_\Blg)\otimes T^*(A)\to
\Mor_\Blg(\lsub{\Alg}M_\Blg,\lsub{\Clg}N_\Blg)
\]
by
\[
\mathcenter{
  \begin{tikzpicture}
    \node at (0,0) (tcblank) {};
    \node at (1,0) (trblank) {};
    \node at (0,-1) (h) {$h$};
    \node at (0,-2) (bcblank) {};
    \draw[Amodar] (tcblank) to node[below,sloped]{\lab{M'}} (h);
    \draw[Amodar] (h) to node[below,sloped]{\lab{N}} (bcblank);
    \draw[tensorblgarrow] (trblank) to node[below,sloped]{\lab{T^*B}} (h);
  \end{tikzpicture}}\,\, , \,\,
\mathcenter{
  \begin{tikzpicture}
    \node at (0,0) (tblank) {};
    \node at (0,-2) (bblank) {};
    \draw[tensoralgarrow] (tblank) to node[below,sloped]{\lab{T^*A}} (bblank);
  \end{tikzpicture}}\,\,
\quad
\stackrel{f^*}{\longrightarrow}
\quad
\mathcenter{
  \begin{tikzpicture}
    \node at (0,0) (tcblank) {};
    \node at (1,0) (trblank) {};
    \node at (-1,0) (tlblank) {};
    \node at (0,-2) (f) {$f$};
    \node at (1,-1) (Delta) {$\Delta$};
    \node at (0,-3) (h) {$h$};
    \node at (0,-4) (bcblank) {};
    \draw[Amodar] (tcblank) to node[below,sloped]{\lab{M}} (f);
    \draw[Amodar] (f) to node[below,sloped]{\lab{M'}} (h);
    \draw[Amodar] (h) to node[below,sloped]{\lab{N}} (bcblank);
    \draw[tensorblgarrow] (trblank) to node[below,sloped]{\lab{T^*B}} (Delta);    
    \draw[tensorblgarrow] (Delta) to (f);
    \draw[tensorblgarrow] (Delta) to (h);
    \draw[tensoralgarrow] (tlblank) to node[below,sloped]{\lab{T^*A}} (f);
  \end{tikzpicture}}.
\]
Similarly, define $g_{*,i,1,j}$ to be zero if $j>0$
and \[
g_*\co T^*(C)\otimes\Mor_\Blg(\lsub{\Alg}M_\Blg,\lsub{\Clg}N_\Blg)\to
  \Mor_\Blg(\lsub{\Alg}M_\Blg,\lsub{\Clg}N'_\Blg)
\]
by
\[
\mathcenter{
\begin{tikzpicture}
  \node at (0,0) (tblank) {};
  \node at (0,-2) (bblank) {};
  \draw[tensoralgarrow] (tblank) to node[below,sloped]{\lab{T^*C}} (bblank);
\end{tikzpicture}}
\,\, , \,\,
\mathcenter{
\begin{tikzpicture}
  \node at (0,0) (tcblank) {};
  \node at (1,0) (trblank) {};
  \node at (0,-1) (h) {$h$};
  \node at (0,-2) (bcblank) {};
  \draw[Amodar] (tcblank) to node[below,sloped]{\lab{M}} (h);
  \draw[Amodar] (h) to node[below,sloped]{\lab{N}} (bcblank);
  \draw[tensorblgarrow] (trblank) to node[below,sloped]{\lab{T^*B}} (h);
\end{tikzpicture}}
\quad
\stackrel{g_*}{\longrightarrow}
\quad
\mathcenter{
\begin{tikzpicture}
  \node at (0,0) (tcblank) {};
  \node at (1,0) (trblank) {};
  \node at (-1,0) (tlblank) {};
  \node at (0,-3) (g) {$g$};
  \node at (1,-1) (Delta) {$\Delta$};
  \node at (0,-2) (h) {$h$};
  \node at (0,-4) (bcblank) {};
  \draw[Amodar] (tcblank) to node[below,sloped]{\lab{M}} (h);
  \draw[Amodar] (h) to (g); %node[below,sloped]{\lab{N}}
  \draw[Amodar] (g) to node[below,sloped]{\lab{N'}} (bcblank);
  \draw[tensorblgarrow] (trblank) to node[below,sloped]{\lab{T^*B}} (Delta);    
  \draw[tensorblgarrow] (Delta) to (g);
  \draw[tensorblgarrow] (Delta) to (h);
  \draw[tensoralgarrow] (tlblank) to node[below,sloped]{\lab{T^*C}} (g);
\end{tikzpicture}}.
\]

\begin{proposition}\label{Prop:Hom-funct}
  The assignments $f\mapsto f^*$ and $g\mapsto g_*$ define chain maps
  \begin{align*}
    \Mor_\Blg(\lsub{\Alg}M_\Blg,\lsub{\Alg}M'_\Blg)&\to
    \lsub{\Clg}{\Mor}_\Alg(\Mor_\Blg(\lsub{\Alg}M'_\Blg,\lsub{\Clg}N_\Blg),\Mor_{\Blg}(\lsub{\Alg}M_\Blg,\lsub{\Clg}N_\Blg))\\ \shortintertext{and}
    \Mor_\Blg(\lsub{\Clg}N_\Blg,\lsub{\Clg}N'_\Blg)&\to \lsub{\Clg}\Mor_\Alg(\Mor_\Blg(\lsub{\Alg}M_\Blg,\lsub{\Clg}N_\Blg),\Mor_{\Blg}(\lsub{\Alg}M_\Blg,\lsub{\Clg}N'_\Blg)).
  \end{align*}
  In particular,
  \begin{enumerate}
  \item If $f$ and $g$ are $\Ainf$ bimodule
    homomorphisms then $f^*$ and $g_*$ respect the bimodule structures
    on $\Mor_\Blg$.
  \item If $f'\co\lsub{\Alg}M_\Blg\to \lsub{\Alg}M'_\Blg$ is
    homotopic to $f$ then $(f')^*$ is homotopic to $f^*$. If
    $g'\co\lsub{\Clg}N_\Blg\to\lsub{\Clg}N'_\Blg$ is homotopic to $g$
    then $g'_*$ is homotopic to $g_*$.
  \item If $f$ is a homotopy equivalence then
    $f^*$ is a homotopy equivalence. If $g$ is a homotopy equivalence
    then $g_*$ is a homotopy equivalence.
  \end{enumerate}
  Moreover, these maps are functorial: given $f'\co
  \lsub{\Alg}M'_\Blg\to \lsub{\Alg}M''_\Blg$ and $g'\co
  \lsub{\Clg}N'_\Blg\to\lsub{\Clg}N''_\Blg$, $(f'\circ f)^*=f^*\circ
  (f')^*$ and $(g'\circ g)_*=g'_*\circ g_*$.
\end{proposition}
\begin{proof}
  The statement that the assignments $f\mapsto f^*$ and $g\mapsto g_*$
  are chain maps follows from a straightforward
  computation once you sort through the definitions.
  The
  functoriality statement is
  also verified directly; this verification is less painful.
\end{proof}

Recall that the operation $\DT$ is functorial (up to homotopy). In
particular, given a morphism $f\in\Mor(M_\Blg,N_\Blg)$ and a type $D$
structure $\lsup{\Blg}P$ there is an associated morphism
$(f\DT\Id_P)\in\Mor_\Ground(M_\Blg\DT\lsup{\Blg}P,N_\Blg\DT\lsup{\Blg}P)$,
defined in Figure~\ref{fig:DT-module-maps}. If $M$, $N$ and $P$ are
bimodules $\lsub{\Alg}M_\Blg$, $\lsub{\Clg}N_\Blg$ and
$\lsup{\Blg}P_\Elg$ then we may view $f\DT\Id$ as an element of
$\Mor_{\Elg}(\lsub{\Alg}M_\Blg\DT\lsup{\Blg}P_\Elg,\lsub{\Clg}N_\Blg\DT\lsup{\Blg}P_\Elg).$
\begin{proposition}\label{prop:Hom-tensor-id} Let $\lsub{\Alg}M_\Blg$, $\lsub{\Clg}N_\Blg$ and
  $\lsup{\Blg}P_\Elg$ be bimodules. Then ``tensoring with the
  identity map''
  \[
  \cdot\DT\Id_P\co \Mor_{\Blg}(\lsub{\Alg}M_\Blg,\lsub{\Clg}N_\Blg)\to \Mor_{\Elg}(\lsub{\Alg}M_\Blg\DT\lsup{\Blg}P_\Elg,\lsub{\Clg}N_\Blg\DT\lsup{\Blg}P_\Elg)
  \]
  is a map of $(\Clg,\Alg)$-bimodules.

  An analogous result holds if $P$ is a type \AAm\ module, with $\DTP$
  in place of $\DT$.
\end{proposition}
\begin{proof}
  We will discuss the left action by $T^*\Clg$; the other cases are
  similar. On the one hand,
  \[
  [(c_1\otimes\dots\otimes c_k)\cdot f]\DT\Id_P=
  \mathcenter{
  \begin{tikzpicture}
    \node at (0,0.5) (tcblank) {};
    \node at (-2,0.5) (tlblank) {$c_1\otimes\dots\otimes c_k$};
    \node at (2,0.5) (trblank) {};
    \node at (3,0.5) (trrblank) {};
    \node at (2,-.5) (deltaP) {$\delta^P$};
    \node at (1,-1) (Delta) {$\Delta$};
    \node at (0,-1.5) (f) {$f$};
    \node at (0,-2.5) (mN) {$m^N$};
    \node at (0,-3.5) (bcblank) {};
    \node at (2,-3.5) (brblank) {};
    \draw[Amodar] (tcblank) to (f);
    \draw[Amodar] (f) to (mN);
    \draw[Amodar] (mN) to (bcblank);
    \draw[DAmodar] (trblank) to (deltaP);
    \draw[DAmodar] (deltaP) to (brblank);
    \draw[tensorelgarrow] (trrblank) to (deltaP);
    \draw[tensorblgarrow] (deltaP) to (Delta);
    \draw[tensorblgarrow] (Delta) to (f);
    \draw[tensorblgarrow] (Delta) to (mN);
    \draw[tensorclgarrow] (tlblank) to (mN);
  \end{tikzpicture}}
  \]
  while on the other hand,
  \[
  (c_1\otimes\dots\otimes c_k)\cdot (f\DT\Id_P)=
\mathcenter{
  \begin{tikzpicture}
    \node at (0,0) (tcblank) {};
    \node at (-2,0) (tlblank) {$c_1\otimes\dots\otimes c_k$};
    \node at (1,0) (trblank) {};
    \node at (3,0) (trrblank) {};
    \node at (2,-.5) (Delta) {$\Delta$};
    \node at (1,-1) (deltaP1) {$\delta^P$};
    \node at (1,-2) (deltaP2) {$\delta^P$};
    \node at (0,-1.5) (f1) {$f$};
    \node at (0,-2.5) (mN) {$m^N$};
    \node at (0,-3.5) (bcblank) {};
    \node at (1,-3.5) (brblank) {};
    \draw[Amodar] (tcblank) to (f1);
    \draw[Amodar] (f1) to (mN);
    \draw[Amodar] (mN) to (bcblank);
    \draw[DAmodar] (trblank) to (deltaP1);
    \draw[DAmodar] (deltaP1) to (deltaP2);
    \draw[DAmodar] (deltaP2) to (brblank);
    \draw[tensorelgarrow] (trrblank) to (Delta);
    \draw[tensorelgarrow] (Delta) to (deltaP1);
    \draw[tensorelgarrow] (Delta) to (deltaP2);
    \draw[tensorblgarrow] (deltaP1) to (f1);
    \draw[tensorblgarrow] (deltaP2) to (mN);
    \draw[tensorclgarrow] (tlblank) to (mN);
  \end{tikzpicture}}.
  \]
  But, by the definition of $\delta^P$, these two diagrams are exactly
  the same.

  The analogue if $P$ is a type \AAm\ module follows from the
  definition of $\DTP$ in terms of $\DT$.
\end{proof}

\begin{lemma}\label{lem:Tautology}
  For any strongly unital $\Ainf$-algebras $\Alg$ and $\Blg$ and
  $\Ainf$-bimodule $\lsub{\Blg}M_{\Alg}$, we have that
  $\lsub{\Blg}M_{\Alg}$ is quasi-isomorphic to
  $\Mor_{\Alg}(\lsub{\Alg}\Alg_{\Alg},\lsub{\Blg}M_{\Alg})$ as a
  $(\Blg,\Alg)$-bimodule.
\end{lemma}
\begin{proof}
  For non-negative integers $i$ and $j$, 
  define maps
  $$\phi_{i,1,j}\co
  (B_+)^{\otimes i}\otimes
  \Mor_{\Alg}(\lsub{\Alg}\Alg_{\Alg},\lsub{\Blg}M_{\Alg})
  \otimes (A_+)^{\otimes j}\to
  \lsub{\Blg}M_{\Alg}$$ by 
  \[
  \phi_{i,1,j}(b_1,\dots,b_i,f,a_1,\dots,a_j)=
  \begin{cases}
      f_{j+1}(1,a_1,\dots,a_j) & {\text{if $i=0$}} \\
      0 & {\text{if $i>0$.}}
  \end{cases}
  \]
  One can check that the $\phi_{i,1,j}$ piece
  together to give a chain map from the complex
  $
  \Tensor^*(B_+[1])
  \otimes 
  \Mor_{\Alg}(\lsub{\Alg}\Alg_{\Alg},\lsub{\Blg}M_{\Alg})
  \otimes 
  \Tensor^*(A_+[1])
  $ to $\lsub{\Blg}M_{\Alg}$; i.e., $\phi$ gives a bimodule morphism.

  We claim that the component $\phi_{0,1,0}\co
  \Mor_{\Alg}(\lsub{\Alg}\Alg_{\Alg},\lsub{\Blg}M_{\Alg})\to
  \lsub{\Blg}M_{\Alg}$ is an isomorphism in homology.
  To this end, observe that 
  $$\Mor_{\Alg}(\lsub{\Alg}\Alg_{\Alg},\lsub{\Blg}M_{\Alg})
  \cong \Hom_{\Ground}(A\otimes \Tensor^*(A_+[1]),M),$$
  where the isomorphism is as chain complexes and $\Hom_\Ground$
  means the full chain complex of morphisms, not just the subset of
  chain maps.  The chain complex structure on
  $A\otimes\Tensor^*(A_+[1])$ that makes this isomorphism true is
  simply the bar resolution of $\Ground$
  (i.e., it is $\lsub{\Alg}\Alg_{\Alg}\DT \lsupv{\Alg}\Barop(\Alg)^{\Alg}\DT \Ground$,
  where here $\Ground$ is thought of as a $\Alg$-module using the
  augmentation $\epsilon$). Thus, by Proposition~\ref{prop:Bar-htpy-equiv},
  it follows that $\Hom_{\Ground}(A\otimes \Tensor^*(\Alg_+[1]),M)$
  is homotopy equivalent to $\Hom_{\Ground}(\Ground,M)\cong M$. Indeed, it 
  is straightforward to verify that the homotopy equivalence is furnished
  by our map $\phi_{0,1,0}$. 
\end{proof}

\begin{lemma}\label{lem:Mor-tensor-relation}
  Let $\lsub{\Alg}M_\Blg$ and $\lsub{\Alg}N_\Clg$ be bimodules and
  $\lsup{\Clg}P_\Elg$ a type \DA\ structure. Then
  \[
  \Mor_{\Alg}(\lsub{\Alg}M_\Blg,\lsub{\Alg}N_\Clg\DT\lsup{\Clg}P_\Elg)
  \cong
  \Mor_{\Alg}(\lsub{\Alg}M_\Blg,\lsub{\Alg}N_\Clg)\DT\lsup{\Clg}P_\Elg
  \]
  as $(\Blg,\Elg)$-bimodules.
\end{lemma}
\begin{proof}
  This is immediate from the definitions.
\end{proof}

For computations, it is often more convenient to work with $\Mor$
complexes of type $D$ structures, which tend to be much smaller. We
will outline how this theory works, leaving the reader to supply most
of the proofs.
\begin{definition}
  Let $\lsup{\Alg}M_\Blg$ and $\lsup{\Alg}N_\Clg$ be type \DA\ 
  bimodules. Let
  $\Mor^\Alg(\lsup{\Alg}M_\Blg,\lsup{\Alg}N_\Clg)=\Mor(\lsup{\Alg}{\mathcal{F}}(M),
  \lsup{\Alg}{\mathcal{F}}(N))$ where $\mathcal{F}$ is the functor
  forgetting the right action. Endow
  $\Mor^\Alg(\lsup{\Alg}M_\Blg,\allowbreak\lsup{\Alg}N_\Clg)$ with an $\Ainf$
  $(\Clg,\Blg)$-bimodule structure by setting $m^{\Mor}_{i,1,j}=0$ if
  $i$ and $j$ are both nonzero, and defining the product
  $T^*B\otimes \Mor^\Alg(\lsup{\Alg}M_\Blg,\lsup{\Alg}N_\Clg)\to
  \Mor^\Alg(\lsup{\Alg}M_\Blg,\lsup{\Alg}N_\Clg)$ by
  \[
  \mathcenter{
    \begin{tikzpicture}
      \node at (0,0) (tblank) {};
      \node at (0,-2) (bblank) {};
      \draw[tensoralgarrow] (tblank) to node[above,sloped]{\lab{T^*\Blg}} (bblank);
    \end{tikzpicture}
  }
  \,\, ,
  \mathcenter{
    \begin{tikzpicture}
      \node at (0,0) (tcblank) {};
      \node at (0,-1) (f) {$f^1$};
      \node at (0,-2) (bcblank) {};
      \node at (-1,-2) (blblank) {};
      \draw[Dmodar] (tcblank) to node[above,sloped]{\lab{M}} (f);
      \draw[Dmodar] (f) to node[above,sloped]{\lab{M}} (bcblank);
      \draw[algarrow] (f) to node[above,sloped]{\lab{\Alg}} (blblank);
    \end{tikzpicture}
  }
  \quad\longrightarrow\quad
  \mathcenter{
    \begin{tikzpicture}
      \node at (0,0) (tcblank) {};
      \node at (1,0) (tlblank) {};
      \node at (0,-1) (delta) {$\delta^M$};
      \node at (0,-2) (f) {$f^1$};
      \node at (-1,-3) (mu) {$\mu^\Alg$};
      \node at (-1,-4) (blblank) {};
      \node at (0,-4) (bcblank) {};
      \draw[Dmodar] (tcblank) to node[above,sloped]{\lab{M}} (delta);
      \draw[Dmodar] (delta) to (f);
      \draw[Dmodar] (f) to node[above,sloped]{\lab{N}} (bcblank);
      \draw[tensorblgarrow] (tlblank) to node[below,sloped]{\lab{T^*\Blg}} (delta);
      \draw[tensoralgarrow] (delta) to node[above,sloped]{\lab{T^*\Alg}} (mu);
      \draw[algarrow] (f) to (mu);
      \draw[algarrow] (mu) to (blblank);
    \end{tikzpicture}
  }
  \]
  and the product
  $\Mor^\Alg(\lsup{\Alg}M_\Blg,\lsup{\Alg}N_\Clg)\otimes T^*\Clg\to
  \Mor^\Alg(\lsup{\Alg}M_\Blg,\lsup{\Alg}N_\Clg)$ by
  \[
  \mathcenter{
    \begin{tikzpicture}
      \node at (0,0) (tcblank) {};
      \node at (0,-1) (f) {$f^1$};
      \node at (0,-2) (bcblank) {};
      \node at (-1,-2) (blblank) {};
      \draw[Dmodar] (tcblank) to node[above,sloped]{\lab{M}} (f);
      \draw[Dmodar] (f) to node[above,sloped]{\lab{M}} (bcblank);
      \draw[algarrow] (f) to node[above,sloped]{\lab{\Alg}} (blblank);
    \end{tikzpicture}
  }
  \,\, , \,\,
  \mathcenter{
    \begin{tikzpicture}
      \node at (0,0) (tblank) {};
      \node at (0,-2) (bblank) {};
      \draw[tensoralgarrow] (tblank) to node[above,sloped]{\lab{T^*\Clg}} (bblank);
    \end{tikzpicture}
  }
  \quad\longrightarrow\quad
  \mathcenter{
    \begin{tikzpicture}
      \node at (0,0) (tcblank) {};
      \node at (1,0) (trblank) {};
      \node at (0,-1) (f) {$f^1$};
      \node at (0,-2) (delta) {$\delta^N$};
      \node at (-1,-3) (mu) {$\mu^\Alg$};
      \node at (-1,-4) (blblank) {};
      \node at (0,-4) (bcblank) {};
      \draw[Dmodar] (tcblank) to node[above,sloped]{\lab{M}} (f);
      \draw[Dmodar] (f) to (delta);
      \draw[Dmodar] (delta) to (bcblank);
      \draw[tensorclgarrow] (trblank) to node[below,sloped]{\lab{T^*\Clg}} (delta);
      \draw[tensoralgarrow] (delta) to node[below,sloped]{\lab{T^*\Alg}} (mu);
      \draw[algarrow] (f) to (mu);
      \draw[algarrow] (mu) to (blblank);
    \end{tikzpicture}
  }.
  \]
\end{definition}

The obvious analogue of Proposition~\ref{Prop:Hom-funct} holds in this
context. Moreover, the following analogue of
Proposition~\ref{prop:Hom-tensor-id} is true:
\begin{proposition}\label{prop:Hom-tensor-id-D}
  Let $\lsup{\Alg}M_\Blg$ and $\lsup{\Alg}N_\Clg$ be type \DA\ 
  bimodules and $\lsub{\Elg}P_\Alg$ a type \AAm\  bimodule. Then the
  ``tensoring with the identity map''
  \[
  \Id_P\DT\cdot\co \Mor^\Alg(\lsup{\Alg}M_\Blg,\lsup{\Alg}N_\Clg)\to \Mor_\Elg(\lsub{\Elg}P_\Alg\DT\lsup{\Alg}M_\Blg,\lsub{\Elg}P_\Alg\DT\lsup{\Alg}N_\Clg)
  \]
  is a map of $(\Clg,\Blg)$-bimodules.
\end{proposition}
\begin{corollary}
  \label{cor:InterpretMor}
  If $\Alg$ is a \dg algebra then
  $\Mor^\Alg(\lsup{\Alg}M_\Blg,\lsup{\Alg}N_\Clg)$ is canonically
  isomorphic to
  the chain complex of maps
  \[\psi\co \lsub{\Alg}\Alg_{\Alg}\DT \lsup{\Alg}M_\Blg\to \lsub{\Alg}\Alg_\Alg\DT \lsup{\Alg}N_\Clg\]
  which commute with the $\Alg$ action, equipped with the differential
  \[\partial\psi=\partial_{\Alg\DT N}\circ\psi
  +\psi\circ\partial_{\Alg\DT M}.\]
  If $\lsup{\Alg}M_\Blg$ and $\lsup{\Alg}N_\Clg$ are homotopy
  equivalent to bounded type \DA\ structures then
  the inclusion of this subcomplex into
  $\Mor_\Alg(\lsub{\Alg}\Alg_\Alg\DT
  \lsup{\Alg}M_\Blg,\lsub{\Alg}\Alg_\Alg\DT \lsup{\Alg}N_\Clg)$ is a
  quasi-isomorphism.
\end{corollary}
\begin{proof}
  Taking $P=\lsub{\Alg}\Alg_\Alg$ in
  Proposition~\ref{prop:Hom-tensor-id-D}, we get a map 
  \[
  \Id_\Alg\DT\cdot\co \Mor^\Alg(\lsup{\Alg}M_\Blg,\lsup{\Alg}N_\Clg)\to \Mor_\Alg(\lsub{\Alg}\Alg_\Alg\DT\lsup{\Alg}M_\Blg,\lsub{\Alg}\Alg_\Alg\DT\lsup{\Alg}N_\Clg).
  \]
  It is clear that this map is injective, and so identifies
  $\Mor^\Alg(\lsup{\Alg}M_\Blg,\lsup{\Alg}N_\Clg)$ with some
  subcomplex of the $\Ainf$-maps
  $\lsub{\Alg}\Alg_\Alg\DT\lsup{\Alg}M_\Blg\to \lsub{\Alg}\Alg_\Alg\DT
  \lsup{\Alg}N_\Clg$. It is straightforward to see that it is the
  stated subcomplex. The fact that this subcomplex is homotopy
  equivalent to the entire $\Mor_\Alg$-complex follows from
  Proposition~\ref{prop:D-to-A-and-back}.
\end{proof}

There is one more case that we will consider: that of the type \DA\ structure on the
type $D$ structure morphisms from a type \DA\ module to a type \DD\
module.
Before discussing this
morphism space, we pause to note another interpretation of the chain
complex of morphisms between two type $D$ structures:
\begin{lemma}
  \label{lem:MorInTermsOfDTP}
  Let $\lsup{\Alg}M$ and $\lsup{\Alg}N$ be type $D$ structures. Then 
  \[
  \Mor^\Alg(\lsup{\Alg}M,\lsup{\Alg}N)\cong \overline{M}^\Alg\DT\lsub{\Alg}\Alg_\Alg\DT\lsup{\Alg}N.
  \]
  Here, $\overline{M}^\Alg$ denotes the opposite type $D$ structure to
  $\lsup{\Alg}M$, as in Definition~\ref{def:opposite-type-D}.
\end{lemma}
\begin{proof}
  This is immediate from the definitions, as follows. The chain
  complex for the morphism complex is:
  \[
  \mathcenter{
  \begin{tikzpicture}
    \node at (0,0) (tblank) {};
    \node at (0,-1) (h) {$h^1$};
    \node at (0,-2) (bcblank) {};
    \node at (-1,-2) (blblank) {};
    \draw[Dmodar] (tblank) to node[right]{\lab{M}} (h);
    \draw[Dmodar] (h) to node[right]{\lab{N}} (bcblank);
    \draw[algarrow] (h) to node[above]{\lab{\Alg}} (blblank);
  \end{tikzpicture}}
  \mathcenter{\stackrel{\bdy}{\longrightarrow}}
  \mathcenter{
    \begin{tikzpicture}
    \node at (0,0) (tblank) {};
    \node at (0,-1) (delta1) {$\delta^M$};
    \node at (0,-2) (h) {$h^1$};
    \node at (0,-3) (delta2) {$\delta^N$};
    \node at (-1,-4) (mu) {$\mu$};
    \node at (0,-5) (bcblank) {};
    \node at (-1,-5) (blblank) {};
    \draw[Dmodar] (tblank) to node[right]{\lab{M}} (delta1);
    \draw[Dmodar] (delta1) to node[right]{\lab{M}} (h);
    \draw[Dmodar] (h) to node[right]{\lab{N}} (delta2);
    \draw[Dmodar] (delta2) to node[right]{\lab{N}} (bcblank);
    \draw[tensoralgarrow, bend right=15] (delta1) to (mu);
    \draw[tensoralgarrow, bend right=5] (delta2) to (mu);
    \draw[algarrow, bend right=10] (h) to (mu);
    \draw[algarrow] (mu) to (blblank);
  \end{tikzpicture}
  }
  \]
  which corresponds to the differential on the box complex,
  \[
  \mathcenter{
  \begin{tikzpicture}
    \node at (-1,-3) (tblank) {};
    \node at (0,-1) (h) {$h^1$};
    \node at (1,-3) (bcblank) {};
    \node at (0,-3) (blblank) {};
    \draw[Dmodar] (tblank) to node[above, sloped]{\lab{M}} (h);
    \draw[Dmodar] (h) to node[above, sloped]{\lab{N}} (bcblank);
    \draw[algarrow] (h) to (blblank);
  \end{tikzpicture}}
  \mathcenter{\stackrel{\bdy}{\longrightarrow}}
  \mathcenter{
    \begin{tikzpicture}
    \node at (-2,-5) (tblank) {};
    \node at (-1,-3) (delta1) {$\delta^M$};
    \node at (0,-2) (h) {$h^1$};
    \node at (1,-3) (delta2) {$\delta^N$};
    \node at (0,-4) (mu) {$\mu$};
    \node at (2,-5) (bcblank) {};
    \node at (0,-5) (blblank) {};
    \draw[Dmodar] (tblank) to node[above, sloped]{\lab{M}} (delta1);
    \draw[Dmodar] (delta1) to node[above, sloped]{\lab{M}} (h);
    \draw[Dmodar] (h) to node[above, sloped]{\lab{N}} (delta2);
    \draw[Dmodar] (delta2) to node[above, sloped]{\lab{N}} (bcblank);
    \draw[tensoralgarrow] (delta1) to (mu);
    \draw[tensoralgarrow] (delta2) to (mu);
    \draw[algarrow] (h) to (mu);
    \draw[algarrow] (mu) to (blblank);
  \end{tikzpicture}
  }.
  \]
\end{proof}

Now, suppose that $\lsup{\Alg}M_\Blg$ is a type \DA\ structure and
$\lsup{\Alg}N^\Clg$ is a type \DD\ module. Then we can give the
morphism space $\Mor^\Alg(\lsup{\Alg}M_\Blg,\lsup{\Alg}N^\Clg)$ the
structure of a type \DA\ module via the isomorphism
\begin{equation}
  \Mor^\Alg(\lsup{\Alg}M_\Blg,\lsup{\Alg}N^\Clg)\cong \lsub{\Blg}{\overline{M}}^\Alg\DT\lsub{\Alg}\Alg_\Alg\DT\lsup{\Alg}N^\Clg.\label{eq:MorDT}
\end{equation}
(The opposite type \DA\ structure $\lsub{\Blg}{\overline{M}}^\Alg$ is
defined in Definition~\ref{def:opposite-type-DA}.)
This is related to the module structure on the space of type \AAm\
morphisms as follows:
\begin{proposition}\label{prop:mor-iso-equiv}
  For $\lsup{\Alg}M_\Blg \in \lsupv{\Alg}{\ModCat}_\Blg$ and
  $\lsup{\Alg}N^\Clg \in \lsupv{Alg}{\ModCat}^\Clg$,  
  there is a canonical inclusion 
  \[
  \Mor^\Alg(\lsup{\Alg}M_\Blg,\lsup{\Alg}N^\Clg)\to \Mor_\Alg(\lsub{\Alg}\Alg_\Alg\DT\lsup{\Alg}M_\Blg,\lsub{\Alg}\Alg_\Alg\DT\lsup{\Alg}N^\Clg)
  \]
  respecting the bimodule structure and inducing an isomorphism on
  homology.
\end{proposition}
\begin{proof}
  That there is such an inclusion inducing an isomorphism on homology
  follows from Proposition~\ref{prop:D-to-A-and-back}. We leave
  verification that it respects the bimodule structure to the reader.
\end{proof}

\subsubsection{Hochschild homology}
\label{sec:hochschild-homology}
We next turn to the Hochschild homology, or self tensor product, of a
bimodule with two actions of the same algebra. We first introduce the
classical Hochschild complex of an $\Ainf$-bimodule and then give a
version for type \DA\ structures analogous to the $\DT$ tensor
product; we also prove that the two definitions are equivalent in the
obvious sense (Proposition~\ref{prop:HochschildDA}). The Hochschild
complex of a type \DA\ structure arises naturally when studying the
knot Floer homology of open books; see Section~\ref{sec:self-pairing}.
\begin{definition}
  Let $\Alg$ be an $\Ainf$-algebra over a commutative ground ring
  $\Ground$, and $\lsub{\Alg}M_{\Alg}$ be an $\Ainf$-bimodule.
  The {\em Hochschild complex $\CH(\lsub{\Alg}M_\Alg)$ of $\lsub{\Alg}M_\Alg$} is defined as follows. 
  Let $\CH_n(M)$ be the $\Field$ vector space which is the quotient of
  $M\otimes_\Ground \overbrace{A[1] \otimes_\Ground \dots \otimes_\Ground A[1]}^{n}$ by the relations
  $$e\cdot x \otimes a_1\otimes\dots\otimes a_n=
  x \otimes a_1\otimes\dots\otimes a_n\cdot e,$$ where $e$ ranges over
  $\Ground$. As a vector space,
  $$\CH(\lsub{\Alg}M_\Alg)=\bigoplus_{n=0}^{\infty} \CH_n(M).$$ The differential on
  $\CH(\lsub{\Alg}M_\Alg)$ is given by
  \begin{multline*}
    D(x\otimes a_1\otimes\dots\otimes a_\ell)\\
    = 
    \begin{aligned}[t]
      &\sum_{m+n\leq \ell}
    m_{m,1,n}(a_{\ell-m+1},\dots,a_\ell,x,a_1,\dots,a_n)
    \otimes
    a_{n+1}\otimes\dots\otimes a_{\ell-m} \\
    &+ \sum_{1\leq m< n \leq \ell} x\otimes a_1\dots \otimes
    \mu_{n-m}(a_{m},\dots,a_{n-1}) \otimes a_{n}\otimes\dots\otimes
    a_\ell.
    \end{aligned}
  \end{multline*}
\end{definition}
It is perhaps more instructive to think of $\CH_n(M)$ as generated by
equivalence classes of collections of $n$ elements arranged on a
circle (the equivalence relation coming from a circular tensor
product). From this description, then, $D$ is a sum of
maps, which take any collection of $i\leq n$ consecutive terms and
applies
whichever higher multiplication map is available to this
collection. The result of this component of $D$
lies in $\CH_{n-i+1}(M)$. Graphically, $D$ is
\[
\mathcenter{
  \begin{tikzpicture}
    \node at (0,0) (tcblank) {};
    \node at (2,0) (trblank) {};
    \node at (2,-1) (Delta) {$\Delta$};
    \node at (0,-2) (m) {$m^M$};
    \node at (3,-1.5) (rwrap) {};
    \node at (-1,-1.5) (lwrap) {};
    \node at (0,-3) (bcblank) {};
    \node at (2,-3) (brblank) {};
    \draw[Amodar] (tcblank) to (m);
    \draw[Amodar] (m) to (bcblank);
    \draw[tensoralgarrow] (trblank) to (Delta);
    \draw[tensoralgarrow] (Delta) to (brblank);
    \draw[tensoralgarrow] (Delta) to (m);
    \draw[tensoralgarrow] (Delta) to (rwrap);
    \draw[tensoralgarrow] (lwrap) to (m);
  \end{tikzpicture}
}
+
\,\,\mathcenter{
  \begin{tikzpicture}
    \node at (0,0) (tcblank) {};
    \node at (0,-3) (bcblank) {};
    \node at (1,0) (trblank) {};
    \node at (1,-2) (D) {$\overline{D}^\Alg$};
    \node at (1,-3) (brblank) {};
    \draw[Amodar] (tcblank) to (bcblank);
    \draw[tensoralgarrow] (trblank) to (D);
    \draw[tensoralgarrow] (D) to (brblank);
  \end{tikzpicture}
},
\]
where in the first diagram a bundle of strands in $\tensor^*\Alg$
runs off the right edge and comes back on the left.

Note that the subscript $n$ on $\CH_n$ is not a grading (though it is
a filtration). The grading on $\CH_n$ is given by
\[
\gr(x\otimes a_1\otimes\dots\otimes
a_n)=\gr(x)+\gr(a_1)+\dots+\gr(a_n)+n.\]

\begin{lemma}
  \label{lem:HochschildWellDefined}
  The endomorphism $D$ of $\CH(\lsub{\Alg}M_\Alg)$ is a differential.
  If $\lsub{\Alg}M_\Alg$ and $\lsub{\Alg}N_\Alg$ are $\Ainf$-homotopy
  equivalent bimodules, then the Hochschild complexes
  $\CH(\lsub{\Alg}M_\Alg)$ and $\CH(\lsub{\Alg}N_\Alg)$ are homotopy
  equivalent chain complexes.
\end{lemma}
\begin{proof}
  The fact that $D^2=0$ follows easily from the fact that $m$ and
  $\mu$ satisfy the $\Ainf$-relations.
  For the second statement,
  for $f\in\Mor(\lsub{\Alg}M_\Alg,\lsub{\Alg}N_\Alg)$ define
  \[\CH(f)\co \CH(\lsub{\Alg}M_\Alg)\to\CH(\lsub{\Alg}N_\Alg)\]
  by 
  \[
  \CH(f)=
  \mathcenter{
    \begin{tikzpicture}
      \node at (0,0) (tcblank) {};
      \node at (2,0) (trblank) {};
      \node at (2,-1) (Delta) {$\Delta$};
      \node at (0,-2) (f) {$f^M$};
      \node at (3,-1.5) (rwrap) {};
      \node at (-1,-1.5) (lwrap) {};
      \node at (0,-3) (bcblank) {};
      \node at (2,-3) (brblank) {};
      \draw[Amodar] (tcblank) to (f);
      \draw[Amodar] (m) to (bcblank);
      \draw[tensoralgarrow] (trblank) to (Delta);
      \draw[tensoralgarrow] (Delta) to (brblank);
      \draw[tensoralgarrow] (Delta) to (f);
      \draw[tensoralgarrow] (Delta) to (rwrap);
      \draw[tensoralgarrow] (lwrap) to (f);
    \end{tikzpicture}
  }.
  \]
  (As before, a bundle of strands running off the right edge of a
  diagram comes back on the right.)
  It is easy to verify that this definition makes $\CH$ into a \dg
  functor, and that $\CH(\Id)=\Id$. The result follows.
\end{proof}

\begin{definition}
        The homology of $\CH(M)$, denoted $\HH(M)$, is called the
        {\em Hochschild homology} of the bimodule $M$.
\end{definition}

We next show that, in certain cases, the Hochschild homology can be
computed from a much smaller complex.

First, some terminology. Given a $\Ground$-bimodule $N$, let
$[N,\Ground]$ be the submodule of $N$ generated by all
elements of the form $nk-kn$ where $n\in N$ and $k\in\Ground$, and let
$N^\circ=N/[N,\Ground]$ denote the vector space quotient of $N$ by
$[N,\Ground]$. We call $N^\circ$ the \emph{cyclicization} of $N$. Note
that for bimodules $M$ and $N$, any $(\Ground,\Ground)$-bilinear map
$N\to M$ descends to a linear map $N^\circ\to M^\circ$.

Now, let $\lsup{\Alg}N_\Alg$ be a
type \DA\ structure, and
$\lsub{\Alg}M_\Alg=\lsub{\Alg}\Alg_\Alg\DT\lsup{\Alg}N_\Alg$ the
associated type \AAm\ bimodule. Let
$$\delta^1_{j}\co N\otimes A_+^{\otimes j}
\to A\otimes N$$ denote the structure maps for $\lsup{\Alg}N_\Alg$. We
will assume that $\lsup{\Alg}N_\Alg$ is bounded in the sense of
Definition~\ref{def:DA-bounded}. The vector space $N^\circ$ will be the
underlying vector space for the smaller model for the Hochschild complex of $\lsub{\Alg}M_\Alg$.

The maps $\delta^1_j$ fit together to give a degree $-1$ map
$$\tilde{\delta}\co N \otimes \Tensor^*(\DGA_+[1])\to
\DGA[1]\otimes N\otimes \Tensor^*(\DGA_+[1]),$$
defined by
$$\tilde{\delta}(x\otimes a_1\otimes\dots\otimes a_m)
=\sum_{0\leq j\leq m} \delta^1_j(x\otimes\dots a_1\otimes\dots\otimes a_j)\otimes
a_{j+1}\otimes\dots\otimes a_m.$$

There is an $\Field$-linear, degree $1$ cyclic rotation map
$$R\co \left(\DGA[1]\otimes N \otimes \Tensor^*(\DGA_+[1])\right)^\circ\to
\left(N\otimes \Tensor^*(\DGA_+[1])\right)^\circ$$
defined by
$$R(a_0\otimes x \otimes a_1 \otimes \dots \otimes a_m)
=x\otimes a_1\otimes\dots\otimes a_m\otimes [(\Id-\epsilon)(a_0)],$$
where $\epsilon$ denotes the augmentation on $\Alg$. Note that the
map $R$ would not make sense without cyclicizing: if
$\iota_1,\iota_2\in\Ground$ are orthogonal idempotents, and $a\in
\DGA$ and $n\in N$ are such that $a=\iota_1a\iota_1$ and
$n=\iota_2n\iota_1$ then $a\otimes n=0$ but $R(a\otimes n)=n\otimes
a\neq 0.$
We can similarly define a map
\[ {\overline R}\co (\Tensor^*(A_+[1])\otimes N\otimes 
\Tensor^*(A_+[1])) \to (N\otimes \Tensor^*(A_+[1]))^\circ\]
by
\[ {\overline R}(a_0\otimes\dots\otimes a_n\otimes x 
\otimes a_{n+1}\otimes\dots\otimes a_{m})
=x\otimes a_{n+1}\otimes\dots\otimes a_{m}\otimes a_0
\otimes\dots \otimes a_n. \]

The map $\Id-\epsilon$ and its tensor powers will come up frequently,
so we let $\pi\co A[1]^{\otimes k}\to A_+[1]^{\otimes k}$ (respectively
$\pi\co T^*(A[1])\to T^*(A_+[1])$) denote
$(\Id-\epsilon)^{\otimes k}$ (respectively
$\bigoplus_k(\Id-\epsilon)^{\otimes k}$). We will need also the obvious
inclusion map $\iota\co N\to N\otimes \Tensor^*(\DGA_+[1]),$ and
its cyclicization, which we also denote by $\iota$.

Provided that $\lsup{\Alg}N_\Alg$ is bounded, these ingredients can be
assembled to form a linear map $\widetilde{\partial} \co N^\circ \to N^\circ$
defined by
\begin{equation}
  \label{eq:HochschildDiff}
  {\widetilde \partial} = \sum_{n=1}^{\infty} \epsilon\circ\tilde{\delta}\circ
  \overbrace{(R\circ\tilde{\delta})\circ\dots\circ (R\circ\tilde{\delta})}^{n-1} {}\circ \iota.
\end{equation}
We can draw the map
$\widetilde{\partial}$ as
\begin{equation}\label{eq:DA-hoch-diff-V1}
  \widetilde{\partial} = \sum\quad
\mathcenter{
\begin{tikzpicture}
  \node at (0,0) (tcblank) {};
  \node at (0,-1) (delta1) {$\delta^1$};
  \node at (-1,-1.5) (lwrap1) {};
  \node at (3,-1.5) (rwrap1) {};
  \node at (2,-2) (IdE1) {$\pi$};
  \node at (1,-2.5) (Delta1) {$\Delta$};
  \node at (0,-3) (delta2) {$\delta^1$};
  \node at (-1,-3.5) (lwrap2) {};
  \node at (3,-3.5) (rwrap2) {};
  \node at (2,-4) (IdE2) {$\pi$};
  \node at (1,-4.5) (merge1) {};
  \node at (1,-5.5) (Delta2) {$\Delta$};
  \node at (0,-6) (delta3) {$\delta^1$};
  \node at (-1,-6.5) (lwrap3) {};
  \node at (3,-6.5) (rwrap3) {};
  \node at (2,-7) (IdE3) {$\pi$};
  \node at (1,-7.5) (merge2) {};
  \node at (1,-8.5) (rdots) {$\vdots$};
  \node at (0,-8.5) (cdots) {$\vdots$};
  \node at (1,-9.5) (blank) {};
  \node at (0,-10) (deltak) {$\delta^1$};
  \node at (-1,-10.5) (epsilon) {$\epsilon$};
  \node at (0,-11) (bcblank) {};
  \draw[DAmodar] (tcblank) to (delta1);
  \draw[DAmodar] (delta1) to (delta2);
  \draw[DAmodar] (delta2) to (delta3);
  \draw[DAmodar] (delta3) to (cdots);
  \draw[DAmodar] (cdots) to (deltak);
  \draw[DAmodar] (deltak) to (bcblank);
  \draw[algarrow] (delta1) to (lwrap1);
  \draw[algarrow] (rwrap1) to (IdE1);
  \draw[algarrow] (IdE1) to (Delta1);
  \draw[tensoralgarrow] (Delta1) to (delta2);
  \draw[tensoralgarrow] (Delta1) to (merge1);
  \draw[algarrow] (delta2) to (lwrap2);
  \draw[algarrow] (rwrap2) to (IdE2);
  \draw[algarrow] (IdE2) to (merge1);
  \draw[tensoralgarrow] (merge1) to (Delta2);
  \draw[tensoralgarrow] (Delta2) to (delta3);
  \draw[tensoralgarrow] (Delta2) to (merge2);
  \draw[algarrow] (delta3) to (lwrap3);
  \draw[algarrow] (rwrap3) to (IdE3);
  \draw[algarrow] (IdE3) to (merge2);
  \draw[tensoralgarrow] (merge2) to (rdots);
  \draw[tensoralgarrow] (rdots) to (blank) to (deltak);
  \draw[algarrow] (deltak) to (epsilon);
\end{tikzpicture}}.
\end{equation}

\begin{lemma}
  Assume $\lsupv{\Alg}N_{\Alg}$ is bounded. For any $x\in
  \lsupv{\Alg}N_{\Alg}$ the sum defining ${\widetilde\partial}(x)$ is
  finite.
\end{lemma}

\begin{proof}
  This is immediate from the definitions.
\end{proof}

Our next goal is to show that $\widetilde{\partial}^2=0$. In order to
do this, we develop some more notation.

For $N$ and $\Alg$ as above, define maps $f_k\co (N\otimes
T^*(A_+[1]))^\circ\to (N\otimes T^*(A_+[1]))^\circ$ by defining $f_0$ to be the
identity map and defining $f_k$ by the diagram
\[
\mathcenter{
  \begin{tikzpicture}
    \node at (0,0) (tcblank) {};
    \node at (1,0) (trblank) {};
    \node at (0,-1) (fk) {$f_{k+1}$};
    \node at (0,-2) (bcblank) {};
    \node at (1,-2) (brblank) {};
    \draw[DAmodar] (tcblank) to (fk);
    \draw[DAmodar] (fk) to (bcblank);
    \draw[tensoralgarrow] (trblank) to (fk);
    \draw[tensoralgarrow] (fk) to (brblank);
  \end{tikzpicture}
}
=
\mathcenter{
\begin{tikzpicture}
  \node at (0,.5) (tcblank) {};
  \node at (1,.5) (trblank) {};
  \node at (0,-.5) (delta) {$\delta$};
  \node at (-1,-1) (IdE) {$\pi$};
  \node at (-2,-1.5) (lwrap) {};
  \node at (1.5,-1.5) (rwrap) {};
  \node at (0,-2) (fk) {$f_k$};
  \node at (0,-3) (bcblank) {};
  \node at (1,-3) (brblank) {};
  \draw[DAmodar] (tcblank) to (delta);
  \draw[DAmodar] (delta) to (fk);
  \draw[DAmodar] (fk) to (bcblank);
  \draw[tensoralgarrow] (trblank) to (delta);
  \draw[tensoralgarrow] (delta) to (IdE);
  \draw[tensoralgarrow] (IdE) to (lwrap);
  \draw[tensoralgarrow] (rwrap) to (fk);
  \draw[tensoralgarrow] (fk) to (brblank);
\end{tikzpicture}}
=\!
\mathcenter{
\begin{tikzpicture}
  \node at (0,.5) (tcblank) {};
  \node at (1,.5) (trblank) {};
  \node at (0,-.5) (fk) {$f_k$};
  \node at (0,-1.5) (delta) {$\delta$};
  \node at (-1,-2) (IdE) {$\pi$};
  \node at (-2,-2.5) (lwrap) {};
  \node at (2,-2.5) (rwrap) {};
  \node at (0,-3) (bcblank) {};
  \node at (1,-3) (brblank) {};
  \draw[DAmodar] (tcblank) to (fk);
  \draw[DAmodar] (fk) to (delta);
  \draw[DAmodar] (delta) to (bcblank);
  \draw[tensoralgarrow] (trblank) to (fk);
  \draw[tensoralgarrow, bend left=60] (fk) to (delta);
  \draw[tensoralgarrow] (delta) to (IdE);
  \draw[tensoralgarrow] (IdE) to (lwrap);
  \draw[tensoralgarrow] (rwrap) to (brblank);
\end{tikzpicture}}.
\]

Similarly, define maps $f_k^{(2)}\co (N\otimes T^*(A_+[1])\otimes T^*(A_+[1]))^\circ\to
(N\otimes T^*(A_+[1])\otimes T^*(A_+[1]))^\circ$ by letting $f_0^{(2)}$ be the identity
map and 
\[
\mathcenter{
  \begin{tikzpicture}
    \node at (0,0) (tcblank) {};
    \node at (1,0) (trblank) {};
    \node at (2,0) (trrblank) {};
    \node at (0,-1) (fk) {$f_{k+1}^{(2)}$};
    \node at (0,-2) (bcblank) {};
    \node at (1,-2) (brblank) {};
    \node at (2,-2) (brrblank) {};
    \draw[DAmodar] (tcblank) to (fk);
    \draw[DAmodar] (fk) to (bcblank);
    \draw[tensoralgarrow] (trblank) to (fk);
    \draw[tensoralgarrow] (fk) to (brblank);
    \draw[tensoralgarrow] (trrblank) to (fk);
    \draw[tensoralgarrow] (fk) to (brrblank);
  \end{tikzpicture}
}
=
\mathcenter{
\begin{tikzpicture}
  \node at (0,.5) (tcblank) {};
  \node at (1,.5) (trblank) {};
  \node at (2,.5) (trrblank) {};
  \node at (0,-.5) (delta) {$\delta$};
  \node at (-1,-1) (IdE) {$\pi$};
  \node at (-2,-1.5) (lwrap) {};
  \node at (2.5,-1.5) (rwrap) {};
  \node at (0,-2) (fk) {$f_k^{(2)}$};
  \node at (0,-3) (bcblank) {};
  \node at (1,-3) (brblank) {};
  \node at (2,-3) (brrblank) {};
  \draw[DAmodar] (tcblank) to (delta);
  \draw[DAmodar] (delta) to (fk);
  \draw[DAmodar] (fk) to (bcblank);
  \draw[tensoralgarrow] (trblank) to (delta);
  \draw[tensoralgarrow] (delta) to (IdE);
  \draw[tensoralgarrow] (IdE) to (lwrap);
  \draw[tensoralgarrow] (rwrap) to (fk);
  \draw[tensoralgarrow] (fk) to (brblank);
  \draw[tensoralgarrow] (trrblank) to (fk);
  \draw[tensoralgarrow] (fk) to (brrblank);
\end{tikzpicture}}
=
\mathcenter{
\begin{tikzpicture}
  \node at (0,0) (tcblank) {};
  \node at (1,0) (trblank) {};
  \node at (2,0) (trrblank) {};
  \node at (0,-1) (fk) {$f_k^{(2)}$};
  \node at (0,-2) (delta) {$\delta$};
  \node at (-1,-2.5) (IdE) {$\pi$};
  \node at (-2,-3) (lwrap) {};
  \node at (3,-3) (rwrap) {};
  \node at (0,-3.5) (bcblank) {};
  \node at (1,-3.5) (brblank) {};
  \node at (2,-3.5) (brrblank) {};
  \draw[DAmodar] (tcblank) to (fk);
  \draw[DAmodar] (fk) to (delta);
  \draw[DAmodar] (delta) to (bcblank);
  \draw[tensoralgarrow] (trblank) to (fk);
  \draw[tensoralgarrow] (trrblank) to (fk);
  \draw[tensoralgarrow, bend left=60] (fk) to (delta);
  \draw[tensoralgarrow, bend left=45] (fk) to (brblank);
  \draw[tensoralgarrow] (delta) to (IdE);
  \draw[tensoralgarrow] (IdE) to (lwrap);
  \draw[tensoralgarrow] (rwrap) to (brrblank);
 \end{tikzpicture}}.
\]

\begin{lemma}\label{lem:hoch-stabilizes}
  Assume that $\lsup{\Alg}N_\Alg$ is 
  a bounded type \DA\ structure. Then 
  given any $x\in \lsup{\Alg}N_{\Alg}$, 
  there is a constant $C=C(x)$ with the property that
  for all $k\geq C$,  $f_k(x,\cdot)=f_{k+1}(x,\cdot)$ and moreover $f_k(x\otimes a)=0$
  if $a\in (A_+[1])^{\otimes i}$ for some $i>0$. Consequently, for $k$ sufficiently
  large the $f_k$ are induced by a map
  \[
  f_\infty\co N^\circ\to (N\otimes T^*(A_+[1]))^\circ
  \]
  in the sense that for $k > C(x)$,
  \[
  f_k(x,a) =
  \begin{cases}
    f_\infty^{(2)}(xa) & a \in \Ground\\
    0 & \text{otherwise}.
  \end{cases}
  \]
  Similarly, given $x\in \lsup{\Alg}N_{\Alg}$,
  there is a constant $C=C(x)$ so that for $k\geq C$, $f_k^{(2)}(x,\cdot,\cdot)=f_{k+1}^{(2)}(x,\cdot,\cdot)$;
  $f_k(x\otimes a\otimes b)=0$ if $a$ or $b$ is in $(A_+[1])^{\otimes i}$ for some
  $i>0$; and so the $f_k^{(2)}$ are induced by a map
  \[
  f_\infty^{(2)}\co N^\circ\to (N\otimes T^*(A_+[1])\otimes T^*(A_+[1]))^\circ
  \]
  in the same sense.
\end{lemma}
\begin{proof}
  Recall that $\delta^N=\Id_N+\delta^1+\cdots$. By the definition of
  admissibility there is a $C$ so that $\delta^n=0$ for
  $n>C$.  In $f_k$ the operator $\delta$ occurs $k$ times in a
  row, and consequently all but $C$ of these terms must be $\Id_N$. This
  implies that $f_k(x, \cdot) = f_{k+1}(x,\cdot)$ for sufficient
  large~$k$. To see that $f_k$ is induced from $f_\infty$, note that
  $\Id_N(x\otimes 
  a)=0$ if $a\notin T^0(A_+[1])=\Ground$.

  The arguments for the corresponding statements about $f_k^{(2)}$ are
  similar.
\end{proof}

With this notation, we are now able to reinterpret the
Hochschild differential:
\begin{lemma}\label{lem:hoch-diff-reinterp}
  The operator $\widetilde{\partial}$ is given by
  \begin{equation}\label{eq:DA-hoch-diff-V2}
  \widetilde{\partial}=
  \mathcenter{
    \begin{tikzpicture}
      \node at (0,0) (tcblank) {};
      \node at (0,-1) (finfty) {$f_\infty$};
%      \node at (1,-1.5) (turn) {};
      \node at (0,-2) (delta) {$\delta^1$};
      \node at (-1,-2.5) (epsilon) {$\epsilon$};
      \node at (0,-3) (bcblank) {};
      \draw[DAmodar] (tcblank) to (finfty);
      \draw[DAmodar] (finfty) to (delta);
      \draw[DAmodar] (delta) to (bcblank);
      \draw[tensoralgarrow, bend left=60] (finfty) to (delta);
 %     \draw[tensoralgarrow] (finfty) to (turn) to (delta);
      \draw[algarrow] (delta) to (epsilon);
    \end{tikzpicture}
  }.
  \end{equation}
\end{lemma}
\begin{proof}
  Each term in Formula~(\ref{eq:DA-hoch-diff-V2}) can
  be expanded to give a term in Formula~(\ref{eq:DA-hoch-diff-V1}) (by
  expanding $\delta$ into copies of $\delta^1$). We
  must show that each term in Formula~(\ref{eq:DA-hoch-diff-V1})
  occurs exactly once this way. The idea is to read the expression in
  Formula~(\ref{eq:DA-hoch-diff-V1}) from bottom to top. More
  precisely, consider a term in Formula~(\ref{eq:DA-hoch-diff-V1}); we
  want to write this term in the form of
  Formula~(\ref{eq:DA-hoch-diff-V2}). Label the $\delta^1$'s occurring
  in order as $\delta^1_{(1)},\dots,\delta^1_{(k)}$. The operation
  $\delta^1_{(k)}$ has a sequence of inputs $a_1,\dots,a_\ell$. The input
  $a_1$ came from some $\delta^1_{(i)}$ for some $i<k$. Let
  $\delta_{(1)}= \delta^1_{(k-1)}\circ R\circ\dots\circ R\circ
  \delta^1_{(i)}$. Similarly, the operation $\delta^1_{(i-1)}$ has
  inputs $a_1',\dots,a_l'$, where $a_1'$ is produced by
  $\delta^1_{(j)}$ for some $j<i-1$. Let $\delta_{(2)}=
  \delta^1_{(i-1)}\circ R\circ\dots\circ R\circ
  \delta^1_{(j-1)}$. Repeat, producing operations
  $\delta_{(3)},\dots,\delta_{(m)}$. Then the term in
  Formula~(\ref{eq:DA-hoch-diff-V2}) with operations
  $\delta_{(m)},\dots,\delta_{(1)},\delta^1_{(k)}$ corresponds to the
  given term in Formula~(\ref{eq:DA-hoch-diff-V1}), and moreover it is
  clear from the construction that this is the unique sequence of
  operations corresponding to the given term.
\end{proof}

Next, we summarize the properties of the operators we have introduced:
\begin{lemma}\label{lem:Hoch-op-props}
  With notation as above,
  \begin{enumerate}
  \item\label{item:HLem5}
    $f_\infty^{(2)}=(\Id\otimes\Delta)\circ f_\infty=[({\overline R}\circ
    \delta)\otimes \Id]\circ (\Id\otimes\Delta)\circ f_\infty$.  That
    is,
    \[
    \mathcenter{
      \begin{tikzpicture}
        \node at (0,0) (tcblank) {}; \node at (0,-1) (finfty2)
        {$f_\infty^{(2)}$}; \node at (0,-2) (bcblank) {};
        \node at (1,-2) (brblank) {}; \node at (2,-2) (brrblank) {};
        \draw[DAmodar] (tcblank) to (finfty2); \draw[DAmodar]
        (finfty2) to (bcblank);
        \draw[tensoralgarrow] (finfty2) to (brblank);
        \draw[tensoralgarrow] (finfty2) to (brrblank);
      \end{tikzpicture}} = \mathcenter{
      \begin{tikzpicture}
        \node at (0,0) (tcblank) {}; \node at (0,-1) (finfty)
        {$f_\infty$}; \node at (1,-2) (Delta) {$\Delta$}; \node at
        (0,-3) (bcblank) {}; \node at (1,-3) (brblank) {}; \node at
        (2,-3) (brrblank) {}; \draw[DAmodar] (tcblank) to (finfty);
        \draw[DAmodar] (finfty) to (bcblank); \draw[tensoralgarrow]
        (finfty) to (Delta); \draw[tensoralgarrow] (Delta) to
        (brblank); \draw[tensoralgarrow] (Delta) to (brrblank);
      \end{tikzpicture}
    } =
    \mathcenter{
      \begin{tikzpicture}
        \node at (0,0) (tcblank) {}; \node at (0,-1) (finfty)
        {$f_\infty$}; \node at (1,-1.5) (Delta) {$\Delta$}; \node at
        (0,-2) (delta) {$\delta$}; \node at (-1,-2.5) (IdE)
        {$\pi$}; \node at (-2,-3) (lwrap) {}; \node at (3,-3)
        (rwrap) {}; \node at (0,-3.5) (bcblank) {}; \node at (1,-3.5)
        (brblank) {}; \node at (2,-3.5) (brrblank) {}; \draw[DAmodar]
        (tcblank) to (finfty); \draw[DAmodar] (finfty) to (delta);
        \draw[DAmodar] (delta) to (bcblank); \draw[tensoralgarrow]
        (finfty) to (Delta); \draw[tensoralgarrow] (Delta) to
        (brblank); \draw[tensoralgarrow] (Delta) to (delta);
        \draw[tensoralgarrow] (delta) to (IdE); \draw[tensoralgarrow]
        (IdE) to (lwrap); \draw[tensoralgarrow] (rwrap) to (brrblank);
      \end{tikzpicture}}.
    \]
    \item\label{item:HLem4} The operator
  \[
  \mathcenter{
    \begin{tikzpicture}
      \node at (0,0) (tcblank) {};
      \node at (0,-1) (finfty) {$f_\infty$};
      \node at (1.5,-1.5) (D) {$\overline{D}^\Alg$};
      \node at (0,-2) (fk) {$f_k$};
      \node at (0,-3) (bcblank) {};
      \node at (1,-3) (brblank) {};
      \draw[DAmodar] (tcblank) to (finfty);
      \draw[DAmodar] (finfty) to (fk);
      \draw[DAmodar] (fk) to (bcblank);
      \draw[tensoralgarrow] (finfty) to (D);
      \draw[tensoralgarrow] (D) to (fk);
      \draw[tensoralgarrow] (fk) to (brblank);
    \end{tikzpicture}}
  +
  \mathcenter{
    \begin{tikzpicture}
      \node at (0,0) (tcblank) {};
      \node at (0,-1) (finfty) {$f_\infty$};
      \node at (1.5,-1.5) (Delta) {$\Delta$};
      \node at (0,-2) (delta) {$\delta^1$};
      \node at (-1,-2.5) (epsilon) {$\epsilon$};
      \node at (0,-3) (fk) {$f_k$};
      \node at (0,-4) (bcblank) {};
      \node at (1,-4) (brblank) {};
      \draw[DAmodar] (tcblank) to (finfty);
      \draw[DAmodar] (finfty) to (delta);
      \draw[DAmodar] (delta) to (fk);
      \draw[DAmodar] (fk) to (bcblank);
      \draw[tensoralgarrow] (finfty) to (Delta);
      \draw[tensoralgarrow] (Delta) to (fk);
      \draw[tensoralgarrow] (Delta) to (delta);
      \draw[tensoralgarrow] (fk) to (brblank);
      \draw[algarrow] (delta) to (epsilon);
    \end{tikzpicture}}
  \]
  is independent of $k$.
  \end{enumerate}
\end{lemma}
\begin{proof}
  We prove Part~(\ref{item:HLem5}) by induction on $k$, showing more
  generally that
  \begin{align*}
  \mathcenter{
    \begin{tikzpicture}
      \node at (0,0) (tcblank) {};
      \node at (0,-1) (fk2) {$f_{2k}^{(2)}$};
      \node at (0,-2) (bcblank) {};
      \node at (1,-2) (brblank) {};
      \node at (2,-2) (brrblank) {};
      \draw[DAmodar] (tcblank) to (fk2);
      \draw[DAmodar] (fk2) to (bcblank);
      \draw[tensoralgarrow] (fk2) to (brblank);
      \draw[tensoralgarrow] (fk2) to (brrblank);
    \end{tikzpicture}}
  &=
  \,\,\mathcenter{
    \begin{tikzpicture}
      \node at (0,0) (tcblank) {};
      \node at (0,-1) (fk) {$f_k$};
      \node at (1,-2) (Delta) {$\Delta$};
      \node at (0,-3) (bcblank) {};
      \node at (1,-3) (brblank) {};
      \node at (2,-3) (brrblank) {};
      \draw[DAmodar] (tcblank) to (fk);
      \draw[DAmodar] (fk) to (bcblank);
      \draw[tensoralgarrow] (fk) to (Delta);
      \draw[tensoralgarrow] (Delta) to (brblank);      
      \draw[tensoralgarrow] (Delta) to (brrblank);      
    \end{tikzpicture}
  } \displaybreak[1] \\
  \mathcenter{
    \begin{tikzpicture}
      \node at (0,0) (tcblank) {};
      \node at (0,-1) (fk2) {$f_{2k+1}^{(2)}$};
      \node at (0,-2) (bcblank) {};
      \node at (1,-2) (brblank) {};
      \node at (2,-2) (brrblank) {};
      \draw[DAmodar] (tcblank) to (fk2);
      \draw[DAmodar] (fk2) to (bcblank);
      \draw[tensoralgarrow] (fk2) to (brblank);
      \draw[tensoralgarrow] (fk2) to (brrblank);
    \end{tikzpicture}}
  &=
    \mathcenter{
      \begin{tikzpicture}
        \node at (0,0) (tcblank) {};
        \node at (0,-1) (finfty) {$f_k$};
        \node at (1,-1.5) (Delta) {$\Delta$}; 
        \node at (0,-2) (delta) {$\delta$}; 
        \node at (-1,-2.5) (IdE) {$\pi$}; 
        \node at (-2,-3) (lwrap) {};
        \node at (3,-3) (rwrap) {}; 
        \node at (0,-3.5) (bcblank) {}; 
        \node at (1,-3.5) (brblank) {};
        \node at (2,-3.5) (brrblank) {}; 
        \draw[DAmodar] (tcblank) to (finfty); 
        \draw[DAmodar] (finfty) to (delta);
        \draw[DAmodar] (delta) to (bcblank); 
        \draw[tensoralgarrow] (finfty) to (Delta); 
        \draw[tensoralgarrow] (Delta) to (brblank); 
        \draw[tensoralgarrow] (Delta) to (delta);
        \draw[tensoralgarrow] (delta) to (IdE); 
        \draw[tensoralgarrow] (IdE) to (lwrap); 
        \draw[tensoralgarrow] (rwrap) to (brrblank);
      \end{tikzpicture}}.
  \end{align*}
  (These identities relate to the portion of $f_k$, etc., with no
  algebra inputs.)
  Indeed, the second equality follows from the first and the
  definition of $f_k^{(2)}$. For the first equality, the $k=0$ case is
  trivial: both sides reduce to the identity map on $N$. For the
  inductive step,
  \begin{align*}
  \mathcenter{
    \begin{tikzpicture}
      \node at (0,0) (tcblank) {};
      \node at (0,-1) (fk) {$f_{k+1}$};
      \node at (1,-2) (Delta) {$\Delta$};
      \node at (0,-3) (bcblank) {};
      \node at (1,-3) (brblank) {};
      \node at (2,-3) (brrblank) {};
      \draw[DAmodar] (tcblank) to (fk);
      \draw[DAmodar] (fk) to (bcblank);
      \draw[tensoralgarrow] (fk) to (Delta);
      \draw[tensoralgarrow] (Delta) to (brblank);      
      \draw[tensoralgarrow] (Delta) to (brrblank);      
    \end{tikzpicture}
  }
  &=
  \mathcenter{
    \begin{tikzpicture}
      \node at (0,0) (tcblank) {};
      \node at (0,-1) (fk) {$f_{k}$};
      \node at (0,-2) (delta) {$\delta$};
      \node at (-1,-2.5) (IdE) {$\pi$};
      \node at (-2,-3) (lwrap) {};
      \node at (2,-3) (rwrap) {};
      \node at (1,-3.5) (Delta) {$\Delta$};
      \node at (0,-4.5) (bcblank) {};
      \node at (1,-4.5) (brblank) {};
      \node at (2,-4.5) (brrblank) {};
      \draw[DAmodar] (tcblank) to (fk);
      \draw[DAmodar] (fk) to (delta);
      \draw[DAmodar] (delta) to (bcblank);
      \draw[tensoralgarrow, bend left=60] (fk) to (delta);
      \draw[tensoralgarrow] (delta) to (IdE);
      \draw[tensoralgarrow] (IdE) to (lwrap);
      \draw[tensoralgarrow] (rwrap) to (Delta);      
      \draw[tensoralgarrow] (Delta) to (brblank);      
      \draw[tensoralgarrow] (Delta) to (brrblank);      
    \end{tikzpicture}}
  =
  \mathcenter{
    \begin{tikzpicture}
      \node at (0,0) (tcblank) {};
      \node at (0,-1) (fk) {$f_{k}$};
      \node at (1,-1.5) (Delta) {$\Delta$};
      \node at (0,-2) (delta1) {$\delta$};
      \node at (-1,-2.5) (IdE1) {$\pi$};
      \node at (0,-3) (delta2) {$\delta$};
      \node at (-1,-3.5) (IdE2) {$\pi$};
      \node at (-2,-3) (lwrap1) {};
      \node at (3,-3) (rwrap1) {};
      \node at (-2,-4) (lwrap2) {};
      \node at (3,-4) (rwrap2) {};
      \node at (0,-5) (bcblank) {};
      \node at (1,-5) (brblank) {};
      \node at (2,-5) (brrblank) {};
      \draw[DAmodar] (tcblank) to (fk);
      \draw[DAmodar] (fk) to (delta1);
      \draw[DAmodar] (delta1) to (delta2);
      \draw[DAmodar] (delta2) to (bcblank);
      \draw[tensoralgarrow] (fk) to (Delta);
      \draw[tensoralgarrow] (Delta) to (delta1);
      \draw[tensoralgarrow] (Delta) to (delta2);
      \draw[tensoralgarrow] (delta1) to (IdE1);
      \draw[tensoralgarrow] (IdE1) to (lwrap1);
      \draw[tensoralgarrow] (rwrap1) to (brblank);
      \draw[tensoralgarrow] (delta2) to (IdE2);
      \draw[tensoralgarrow] (IdE2) to (lwrap2);
      \draw[tensoralgarrow] (rwrap2) to (brrblank);
    \end{tikzpicture}}\displaybreak[1] \\
  &=
  \mathcenter{
    \begin{tikzpicture}
      \node at (0,0) (tcblank) {};
      \node at (0,-1) (fk) {$f_{2k}^{(2)}$};
      \node at (0,-2) (delta1) {$\delta$};
      \node at (-1,-2.5) (IdE1) {$\pi$};
      \node at (0,-3) (delta2) {$\delta$};
      \node at (-1,-3.5) (IdE2) {$\pi$};
      \node at (-2,-3) (lwrap1) {};
      \node at (3,-3) (rwrap1) {};
      \node at (-2,-4) (lwrap2) {};
      \node at (3,-4) (rwrap2) {};
      \node at (0,-5) (bcblank) {};
      \node at (1,-5) (brblank) {};
      \node at (2,-5) (brrblank) {};
      \draw[DAmodar] (tcblank) to (fk);
      \draw[DAmodar] (fk) to (delta1);
      \draw[DAmodar] (delta1) to (delta2);
      \draw[DAmodar] (delta2) to (bcblank);
      \draw[tensoralgarrow, bend left=60] (fk) to (delta1);
      \draw[tensoralgarrow, bend left=60] (fk) to (delta2);
      \draw[tensoralgarrow] (delta1) to (IdE1);
      \draw[tensoralgarrow] (IdE1) to (lwrap1);
      \draw[tensoralgarrow] (rwrap1) to (brblank);
      \draw[tensoralgarrow] (delta2) to (IdE2);
      \draw[tensoralgarrow] (IdE2) to (lwrap2);
      \draw[tensoralgarrow] (rwrap2) to (brrblank);
    \end{tikzpicture}}
  =\,
  \mathcenter{
    \begin{tikzpicture}
      \node at (0,0) (tcblank) {};
      \node at (0,-1) (fk2) {$f_{2k+2}^{(2)}$};
      \node at (0,-2) (bcblank) {};
      \node at (1,-2) (brblank) {};
      \node at (2,-2) (brrblank) {};
      \draw[DAmodar] (tcblank) to (fk2);
      \draw[DAmodar] (fk2) to (bcblank);
      \draw[tensoralgarrow] (fk2) to (brblank);
      \draw[tensoralgarrow] (fk2) to (brrblank);
    \end{tikzpicture}}.
  \end{align*}
  where the first equality uses the inductive definition of $f_{k+1}$; the
  second equality uses the fact that $\delta$ and $\pi$ respect the
  coalgebra structure of $T^*A_+$; the third equality uses the
  inductive hypothesis; and the fourth uses the inductive definition
  of $f_{2k+2}^{(2)}$ (twice).  

  For Part~(\ref{item:HLem4}), observe that
  \begin{align*}
  \mathcenter{
    \begin{tikzpicture}
      \node at (0,0) (tcblank) {};
      \node at (0,-1) (finfty) {$f_\infty$};
      \node at (1.5,-1.5) (D) {$\overline{D}^\Alg$};
      \node at (0,-2) (fk) {$f_k$};
      \node at (0,-3) (bcblank) {};
      \node at (1,-3) (brblank) {};
      \draw[DAmodar] (tcblank) to (finfty);
      \draw[DAmodar] (finfty) to (fk);
      \draw[DAmodar] (fk) to (bcblank);
      \draw[tensoralgarrow] (finfty) to (D);
      \draw[tensoralgarrow] (D) to (fk);
      \draw[tensoralgarrow] (fk) to (brblank);
    \end{tikzpicture}}
  +
  \mathcenter{
    \begin{tikzpicture}
      \node at (0,0) (tcblank) {};
      \node at (0,-1) (finfty) {$f_\infty$};
      \node at (1.5,-1.5) (Delta) {$\Delta$};
      \node at (0,-2) (delta) {$\delta^1$};
      \node at (-1,-2.5) (epsilon) {$\epsilon$};
      \node at (0,-3) (fk) {$f_k$};
      \node at (0,-4) (bcblank) {};
      \node at (1,-4) (brblank) {};
      \draw[DAmodar] (tcblank) to (finfty);
      \draw[DAmodar] (finfty) to (delta);
      \draw[DAmodar] (delta) to (fk);
      \draw[DAmodar] (fk) to (bcblank);
      \draw[tensoralgarrow] (finfty) to (Delta);
      \draw[tensoralgarrow] (Delta) to (fk);
      \draw[tensoralgarrow] (Delta) to (delta);
      \draw[tensoralgarrow] (fk) to (brblank);
      \draw[algarrow] (delta) to (epsilon);      
    \end{tikzpicture}}
  &=
  \mathcenter{
    \begin{tikzpicture}
      \node at (0,0) (tcblank) {};
      \node at (0,-1) (finfty) {$f_\infty$};
      \node at (1.5,-1.5) (D) {$\overline{D}^\Alg$};
      \node at (0,-2) (delta) {$\delta$};
      \node at (-1,-2.5) (IdE) {$\pi$};
      \node at (-2, -3) (lwrap) {};
      \node at (1,-3) (rwrap) {};
      \node at (0,-3.5) (fk) {$f_{k-1}$};
      \node at (0,-4.5) (bcblank) {};
      \node at (1,-4.5) (brblank) {};
      \draw[DAmodar] (tcblank) to (finfty);
      \draw[DAmodar] (finfty) to (delta);
      \draw[DAmodar] (delta) to (fk);
      \draw[DAmodar] (fk) to (bcblank);
      \draw[tensoralgarrow] (finfty) to (D);
      \draw[tensoralgarrow] (D) to (delta);
      \draw[tensoralgarrow] (delta) to (IdE);
      \draw[tensoralgarrow] (IdE) to (lwrap);
      \draw[tensoralgarrow] (rwrap) to (fk);
      \draw[tensoralgarrow] (fk) to (brblank);
    \end{tikzpicture}}
  +
  \mathcenter{
    \begin{tikzpicture}
      \node at (0,0) (tcblank) {};
      \node at (0,-1) (finfty) {$f_\infty$};
      \node at (1.5,-1.5) (Delta) {$\Delta$};
      \node at (0,-2) (delta) {$\delta^1$};
      \node at (-1,-2.5) (epsilon) {$\epsilon$};
      \node at (0,-3) (delta2) {$\delta$};
      \node at (-1,-3.5) (IdE) {$\pi$};
      \node at (-2,-4) (lwrap) {};
      \node at (1.5,-4) (rwrap) {};
      \node at (0,-4.5) (fk) {$f_{k-1}$};
      \node at (0,-5.5) (bcblank) {};
      \node at (1,-5.5) (brblank) {};
      \draw[DAmodar] (tcblank) to (finfty);
      \draw[DAmodar] (finfty) to (delta);
      \draw[DAmodar] (delta) to (delta2);
      \draw[DAmodar] (delta2) to (fk);
      \draw[DAmodar] (fk) to (bcblank);
      \draw[tensoralgarrow] (finfty) to (Delta);
      \draw[tensoralgarrow] (Delta) to (delta2);
      \draw[tensoralgarrow] (Delta) to (delta);
      \draw[tensoralgarrow] (delta2) to (IdE);
      \draw[tensoralgarrow] (IdE) to (lwrap);
      \draw[tensoralgarrow] (rwrap) to (fk);
      \draw[tensoralgarrow] (fk) to (brblank);
      \draw[algarrow] (delta) to (epsilon);      
    \end{tikzpicture}}\displaybreak[1] \\
  &=
    \mathcenter{
    \begin{tikzpicture}
      \node at (0,0) (tcblank) {};
      \node at (0,-1) (finfty) {$f_\infty$};
%      \node at (1.5,-1.5) (D) {};
      \node at (0,-2) (delta) {$\delta$};
      \node at (-1,-2.5) (IdE) {$\pi$};
      \node at (-2, -3) (lwrap) {};
      \node at (2.5,-3) (rwrap) {};
      \node at (1.5,-3.5) (Dnew) {$\overline{D}^\Alg$};
      \node at (0,-4) (fk) {$f_{k-1}$};
      \node at (0,-5) (bcblank) {};
      \node at (1,-5) (brblank) {};
      \draw[DAmodar] (tcblank) to (finfty);
      \draw[DAmodar] (finfty) to (delta);
      \draw[DAmodar] (delta) to (fk);
      \draw[DAmodar] (fk) to (bcblank);
%      \draw[tensoralgarrow] (finfty) to (D) to (delta);
      \draw[tensoralgarrow, bend left=60] (finfty) to (delta);
      \draw[tensoralgarrow] (delta) to (IdE);
      \draw[tensoralgarrow] (IdE) to (lwrap);
      \draw[tensoralgarrow] (rwrap) to (Dnew);
      \draw[tensoralgarrow] (Dnew) to (fk);
      \draw[tensoralgarrow] (fk) to (brblank);
    \end{tikzpicture}}
  +
  \mathcenter{
    \begin{tikzpicture}
      \node at (0,0) (tcblank) {};
      \node at (0,-1) (finfty) {$f_\infty$};
      \node at (1.5,-1.5) (Delta) {$\Delta$};
      \node at (0,-3) (delta) {$\delta^1$};
      \node at (-1,-3.5) (epsilon) {$\epsilon$};
      \node at (0,-2) (delta2) {$\delta$};
      \node at (-1,-2.5) (IdE) {$\pi$};
      \node at (-2,-3) (lwrap) {};
      \node at (1.5,-3) (rwrap) {};
      \node at (0,-4) (fk) {$f_{k-1}$};
      \node at (0,-5) (bcblank) {};
      \node at (1,-5) (brblank) {};
      \draw[DAmodar] (tcblank) to (finfty);
      \draw[DAmodar] (finfty) to (delta2);
      \draw[DAmodar] (delta2) to (delta);
      \draw[DAmodar] (delta) to (fk);
      \draw[DAmodar] (fk) to (bcblank);
      \draw[tensoralgarrow] (finfty) to (Delta);
      \draw[tensoralgarrow] (Delta) to (delta2);
      \draw[tensoralgarrow] (Delta) to (delta);
      \draw[tensoralgarrow] (delta2) to (IdE);
      \draw[tensoralgarrow] (IdE) to (lwrap);
      \draw[tensoralgarrow] (rwrap) to (fk);
      \draw[tensoralgarrow] (fk) to (brblank);
      \draw[algarrow] (delta) to (epsilon);      
    \end{tikzpicture}}\displaybreak[1] \\
  &=\,
  \mathcenter{
    \begin{tikzpicture}
      \node at (0,0) (tcblank) {};
      \node at (0,-1) (finfty) {$f_\infty$};
      \node at (1.5,-1.5) (D) {$\overline{D}^\Alg$};
      \node at (0,-2) (fk) {$f_{k-1}$};
      \node at (0,-3) (bcblank) {};
      \node at (1,-3) (brblank) {};
      \draw[DAmodar] (tcblank) to (finfty);
      \draw[DAmodar] (finfty) to (fk);
      \draw[DAmodar] (fk) to (bcblank);
      \draw[tensoralgarrow] (finfty) to (D);
      \draw[tensoralgarrow] (D) to (fk);
      \draw[tensoralgarrow] (fk) to (brblank);
    \end{tikzpicture}}
  +
  \mathcenter{
    \begin{tikzpicture}
      \node at (0,0) (tcblank) {};
      \node at (0,-1) (finfty) {$f_\infty$};
      \node at (1.5,-1.5) (Delta) {$\Delta$};
      \node at (0,-2) (delta) {$\delta^1$};
      \node at (-1,-2.5) (epsilon) {$\epsilon$};
      \node at (0,-3) (fk) {$f_{k-1}$};
      \node at (0,-4) (bcblank) {};
      \node at (1,-4) (brblank) {};
      \draw[DAmodar] (tcblank) to (finfty);
      \draw[DAmodar] (finfty) to (delta);
      \draw[DAmodar] (delta) to (fk);
      \draw[DAmodar] (fk) to (bcblank);
      \draw[tensoralgarrow] (finfty) to (Delta);
      \draw[tensoralgarrow] (Delta) to (fk);
      \draw[tensoralgarrow] (Delta) to (delta);
      \draw[tensoralgarrow] (fk) to (brblank);
      \draw[algarrow] (delta) to (epsilon);      
    \end{tikzpicture}},
  \end{align*} 
  where the first equality uses the definition of $f_{k}$, the second
  uses the type \DA\ structure relation, combined with the fact that
  $$\overline{D}^\Alg\circ \pi+\pi\circ \overline{D}^\Alg
  = \epsilon\otimes \pi+\pi\otimes\epsilon$$
  (which follows from the assumption that $\Alg$ is strictly unital),
  and the third uses
  Part~(\ref{item:HLem5}) and the definition of $f_\infty$.
\end{proof}

\begin{proposition}\label{prop:Hoch-diff-sq-0}
  If $(\lsup{\Alg}N_{\Alg},\delta)$ is a bounded type \DA\ structure,
  then $\widetilde{\partial}^2=0$.
\end{proposition}
\begin{proof}
  To see that $\widetilde{\partial}^2=0$, we prove that for any $k$,
  \begin{equation}\label{eq:inductive-hoch-diff}
  \widetilde{\partial}^2=
    \mathcenter{
    \begin{tikzpicture}
      \node at (0,0) (tcblank) {};
      \node at (0,-1) (finfty) {$f_\infty$};
      \node at (1.5,-1.5) (D) {$\overline{D}^\Alg$};
      \node at (0,-2) (fk) {$f_k$};
%      \node at (1,-4) (turn) {};
      \node at (0,-3) (delta1) {$\delta^1$};
      \node at (-1,-3.5) (epsilon) {$\epsilon$};
      \node at (0,-4) (bcblank) {};
      \draw[DAmodar] (tcblank) to (finfty);
      \draw[DAmodar] (finfty) to (fk);
      \draw[DAmodar] (fk) to (delta1);
      \draw[DAmodar] (delta1) to (bcblank);
      \draw[tensoralgarrow] (finfty) to (D);
      \draw[tensoralgarrow] (D) to (fk);
      \draw[tensoralgarrow, bend left=60] (fk) to (delta1);
%      \draw[tensoralgarrow] (fk) to (turn) to (delta1);
      \draw[algarrow] (delta1) to (epsilon);
    \end{tikzpicture}
  }
  +\,
  \mathcenter{
    \begin{tikzpicture}
      \node at (0,0) (tcblank) {};
      \node at (0,-1) (finfty) {$f_\infty$};
      \node at (1.5,-1.5) (Delta) {$\Delta$};
      \node at (0,-2) (delta1) {$\delta^1$};
      \node at (-1,-2.5) (epsilon) {$\epsilon$};
      \node at (0,-3) (fk) {$f_k$};
%      \node at (1,-5) (turn) {};
      \node at (0,-4) (delta2) {$\delta^1$};
      \node at (-1,-4.5) (epsilon2) {$\epsilon$};
      \node at (0,-5) (bcblank) {};
      \draw[DAmodar] (tcblank) to (finfty);
      \draw[DAmodar] (finfty) to (delta1);
      \draw[DAmodar] (delta1) to (fk);
      \draw[DAmodar] (fk) to (delta2);
      \draw[DAmodar] (delta2) to (bcblank);
      \draw[tensoralgarrow] (finfty) to (Delta);
      \draw[tensoralgarrow] (Delta) to (delta1);
      \draw[tensoralgarrow] (Delta) to (fk);
      \draw[tensoralgarrow, bend left=60] (fk) to (delta2);
%      \draw[tensoralgarrow] (fk) to (turn) to (delta2);
      \draw[algarrow] (delta1) to (epsilon);
      \draw[algarrow] (delta2) to (epsilon2);
    \end{tikzpicture}
  }  
  \end{equation}
  Indeed, for $k\gg0$, the first term of
  Equation~(\ref{eq:inductive-hoch-diff}) vanishes (because of the
  $\overline{D}^\Alg$ followed by the $f_k$; see
  Lemma~\ref{lem:hoch-stabilizes}), while, in light of
  Lemma~\ref{lem:hoch-stabilizes}, the second term reduces to the
  reinterpretation of $\widetilde{\partial}$ from
  Lemma~\ref{lem:hoch-diff-reinterp}.

  On the other hand, it is immediate from Part~(\ref{item:HLem4}) of
  Lemma~\ref{lem:Hoch-op-props} that the expression on the right of
  Equation~(\ref{eq:inductive-hoch-diff}) is independent of $k$. For
  $k=0$, Equation~(\ref{eq:inductive-hoch-diff}) reduces to
  \[
  \widetilde{\partial}^2=
  \mathcenter{
    \begin{tikzpicture}
      \node at (0,0) (tcblank) {};
      \node at (0,-1) (finfty) {$f_\infty$};
      \node at (1.5,-1.5) (D) {$\overline{D}^\Alg$};
      \node at (0,-2) (delta1) {$\delta^1$};
      \node at (-1,-2.5) (epsilon) {$\epsilon$};
      \node at (0,-3) (bcblank) {};
      \draw[DAmodar] (tcblank) to (finfty);
      \draw[DAmodar] (finfty) to (delta1);
      \draw[DAmodar] (delta1) to (bcblank);
      \draw[tensoralgarrow] (finfty) to (D);
      \draw[tensoralgarrow] (D) to (delta1);
      \draw[algarrow] (delta1) to (epsilon);
    \end{tikzpicture}
  }
  \,+\,\,
  \mathcenter{
    \begin{tikzpicture}
      \node at (0,0) (tcblank) {};
      \node at (0,-1) (finfty) {$f_\infty$};
      \node at (1.5,-1.5) (Delta) {$\Delta$};
      \node at (0,-2) (delta1) {$\delta^1$};
      \node at (-1,-2.5) (epsilon) {$\epsilon$};
      \node at (0,-3) (delta2) {$\delta^1$};
      \node at (-1,-3.5) (epsilon2) {$\epsilon$};
      \node at (0,-4) (bcblank) {};
      \draw[DAmodar] (tcblank) to (finfty);
      \draw[DAmodar] (finfty) to (delta1);
      \draw[DAmodar] (delta1) to (delta2);
      \draw[DAmodar] (delta2) to (bcblank);
      \draw[tensoralgarrow] (finfty) to (Delta);
      \draw[tensoralgarrow] (Delta) to (delta1);
      \draw[tensoralgarrow] (Delta) to (delta2);
      \draw[algarrow] (delta1) to (epsilon);
      \draw[algarrow] (delta2) to (epsilon2);
    \end{tikzpicture}    
  }
  \]
  which is zero by the type \DA\ structure relation (with one algebra
  element output).
\end{proof}

We let $\tCH(\lsup{\Alg}N_\Alg)=(N/[\Ground,N],\tilde{\bdy})$ denote
the chain complex associated as above to the type \DA\ structure
$\lsup{\Alg}N_\Alg$.

Similarly, given a bounded type \DA\ morphism between
two bounded type \DA\ structures $f\co \lsup{\Alg}N_\Alg\to
\lsup{\Alg}N'_\Alg$ define a map $\tCH(f)\co N/[\Ground,N]\to
N'/[\Ground,N']$ by
\[
  \tCH(f) = \sum_{n,m=0}^{\infty} \epsilon\circ
  \mathord{\overbrace{(\tilde{\delta}\circ R)\circ\dots\circ
    (\tilde{\delta}\circ R)}^{n}} \circ f\circ
  \mathord{\overbrace{(R\circ\tilde{\delta})\circ\dots\circ (R\circ\tilde{\delta})}^{m}}\circ \iota.
\]

\begin{proposition}
  \label{prop:HochschildPrimeWellDefined}The assignment $\tCH$ is an $\Ainf$
  functor from the category of bounded type \DA\ bimodules to
  $\ModCat_{\Ground}$, the category of chain complexes over~$\Ground$.
  In particular, if
  $f\co \lsup{\Alg}N_\Alg\to \lsup{\Alg}N'_\Alg$ is a bounded type
  \DA\ homomorphism, then $\tCH(f)$ is a
  chain map, and if $f$ and $f'$ are homotopic morphisms, then the
  maps $\tCH(f)$ and $\tCH(f')$ are chain homotopic.
\end{proposition}
\begin{proof}
  The proof is essentially the same as the proof that
  $\widetilde{\partial}^2=0$ (Proposition~\ref{prop:Hoch-diff-sq-0}),
  and we leave it to the reader.
\end{proof}

The relation between $\tCH(\lsup{\Alg}N_\Alg)$ and
$\CH(\lsub{\Alg}M_\Alg)$ is similar to the relation between $\DT$ and
$\DTP$ (Proposition~\ref{prop:D-to-A-and-back}). In particular:
\begin{proposition} Let $\Alg$ be a \dg algebra.
  \label{prop:HochschildDA}
  \begin{enumerate}
  \item\label{item:HHisHH1} Suppose that $\lsub{\Alg}M_\Alg$ is a type
    \AAm\ module. Then $\lsup{\Alg}\Barop(\Alg)^\Alg\DT\lsub{\Alg}M_\Alg$ is
    a bounded type \DA\ structure and the complex
    $\CH(\lsub{\Alg}M_\Alg)$ is isomorphic to the complex
    $\tCH(\lsup{\Alg}\Barop(\Alg)^\Alg\DT\lsub{\Alg}M_\Alg)$.
  \item\label{item:HHisHH2} Suppose that $\lsup{\Alg}N_\Alg$ is a
    bounded type \DA\ structure. Then the complex
    $\tCH(\lsup{\Alg}N_\Alg)$ is homotopy equivalent to the
    complex $\CH(\lsub{\Alg}\Alg_\Alg\DT \lsup{\Alg}N_\Alg)$.
  \end{enumerate}
\end{proposition}
\begin{proof}
  Part~(\ref{item:HHisHH1}) is immediate from the
  definitions. Part~(\ref{item:HHisHH2}) follows from
  Part~(\ref{item:HHisHH1}), which gives the first isomorphism in the
  string
  \[
  \CH(\lsub{\Alg}\Alg_\Alg\DT \lsup{\Alg}N_\Alg)\cong
  \tCH(\lsup{\Alg}\Barop(\Alg)^\Alg\DT\lsub{\Alg}\Alg_\Alg\DT
  \lsup{\Alg}N_\Alg)\simeq\tCH(\lsup{\Alg}N_\Alg);
  \]
  together with the following observation: since $\lsupv{\Alg}N_{\Alg}$
  is bounded, the natural map 
  $$\kappa\otimes \Id_N \lsupv{\Alg}\Barop(\Alg)^{\Alg}\DT \lsub{\Alg}\Alg_{\Alg}
  \DT \lsup{\Alg}N_{\Alg}\longrightarrow
  \lsup{\Alg}N_{\Alg}$$
  is a {\em bounded} homotopy equivalence (a fact which can be seen by looking 
  at the maps from Proposition~\ref{prop:Bar-htpy-equiv}). In view of this fact, 
  Proposition~\ref{prop:HochschildPrimeWellDefined} supplies the second homotopy equivalence.
\end{proof}
\begin{remark}
  The extension of Proposition~\ref{prop:HochschildDA} to
  $\Ainf$-algebras is straightforward; the reason that we restrict to
  the \dg case is that we have not defined the bar resolution more
  generally.
\end{remark}

\subsection{Equivalences of categories}
\label{sec:equiv-categ}
In Section~\ref{sec:algebra-modules} we introduced many different
categories of modules and bimodules. In
Section~\ref{sec:bar-cobar-modules} we showed that the \dg categories
$\ModCat^\Alg$ and $\ModCat_{\Alg}$ are quasi-equivalent, and also
corresponding statements for bimodules; in particular, their
homological categories $\HMod(\ModCat_{\Alg})$ and
$\HMod(\ModCat^{\Alg})$ are equivalent triangulated categories. In
this section we continue to tame the multitude, showing that:
\begin{itemize}
\item If $\Alg$ is a \dg algebra then $\HMod(\ModCat_{\Alg})$ is
  triangle-equivalent to the derived category of \dg $\Alg$-modules.
\item If $\Alg$ and $\Blg$ are quasi-isomorphic \dg (or, more
  generally, $\Ainf$-) algebras then $\ModCat_\Alg$ and $\ModCat_\Blg$
  are quasi-equivalent.
\end{itemize}

The results in this section are not purely of aesthetic interest:
various of them will be used in Sections~\ref{sec:mcg}
and~\ref{sec:Duality}, with
consequences that are useful for computation.

\subsubsection{Homotopy equivalence and quasi-isomorphism}
\label{sec:models-der-cat}
Let $\Alg$ be a differential graded algebra. Throughout this section, the
word \emph{honest} is used to distinguish ordinary (``honest'')
differential graded modules from $\Ainf$-modules. We can consider the
following different
models for the derived category of $\Alg$-modules:
\begin{itemize}
\item The category $\CatHonQi$ with objects honest \dg $\Alg$-modules and
  morphisms obtained by localizing the homotopy category of honest
  module maps with respect to quasi-isomorphisms.  (Recall that a
  quasi-isomorphism is a chain map inducing an isomorphism on
  homology.)
\item The category $\CatHonAinf$ with objects honest \dg $\Alg$-modules and
  morphisms $\Ainf$\hyp homotopy classes of $\Ainf$-morphisms.
\item The category $\CatHonAinfQi$ obtained from $\CatHonAinf$ by
  localizing with respect to $\Ainf$-quasi-isomorphisms.
\item The category $\CatAinfAinf$ with objects $\Ainf$-modules and
  morphisms $\Ainf$-homotopy classes of $\Ainf$-morphisms. (This has
  been denoted $\HMod(\ModCat_\Alg)$ elsewhere in this section.)
\item The category $\CatAinfQi$ obtained from $\CatAinfAinf$ by localizing
  with respect to $\Ainf$\hyp quasi\hyp isomorphisms.
\end{itemize}

\begin{proposition}\label{prop:derived-is-derived-is-derived-is-derived}
  The categories $\CatHonQi$, $\CatHonAinf$, $\CatAinfAinf$ and
  $\CatAinfQi$ are all equivalent triangulated
  categories. Corresponding statements hold for categories of bimodules. 
\end{proposition}
\begin{proof}
  The proof is standard so we will only
  sketch it. The main point is that all of these categories are
  equivalent to the full subcategory of honest, projective modules. To
  see this, one observes that:
  \begin{enumerate}
  \item The bar resolution functor $M\mapsto \Barop(M)$ maps the
    category of ($\Ainf)$ modules into the subcategory of projective
    modules, and
    takes $\Ainf$-module homomorphisms to honest \dg module homomorphisms.
  \item Any map between projective \dg modules inducing an isomorphism on
    homology has a homotopy inverse (see, e.g., \cite[Lemma
    10.12.2.2]{BernsteinLunts94:EquivariantSheaves}).
  \item\label{item:derived-pf-1} The canonical map $\Barop(M)\to M$ is an isomorphism in any of
    $\CatHonQi$, $\CatHonAinf$, $\CatAinfAinf$ and $\CatAinfQi$, and
    so the bar resolution functor is naturally isomorphic to the
    identity functor.
  \end{enumerate}
  These, together, imply the result.
\end{proof}
Note that, concretely, Proposition~\ref{prop:derived-is-derived-is-derived-is-derived} implies that every quasi-isomorphism of $\Ainf$-modules is a homotopy equivalence. The analogue for (bounded) type $D$ structures, and for bimodules, is given in Corollary~\ref{cor:quasi-iso-is-htpy-equiv-D-DA-DD}. 

The honest homotopy category of \dg $\Alg$-modules is not in
general equivalent to the derived category; it is
part~(\ref{item:derived-pf-1}) of the proof that breaks down.

\begin{proposition}\label{prop:DTP-is-derived-tensor-product}
  With respect to the identification from
  Proposition~\ref{prop:derived-is-derived-is-derived-is-derived} the
  tensor product $\DTP$ of Definition~\ref{def:DTP} is identified with
  the usual derived tensor product.
\end{proposition}
\begin{proof}
  This is clear from the definitions and the fact that the bar
  resolution of a module is projective.
\end{proof}

\begin{definition}
  A map $f\co \lsup{\Alg}M\to \lsup{\Alg}N$ of type $D$ modules is
  called a \emph{quasi\hyp isomorphism} if the induced map $\Id_{\Alg}\DT
  f\co \lsub{\Alg}\Alg_\Alg\DT\lsup{\Alg}M\to
  \lsub{\Alg}\Alg_\Alg\DT\lsup{\Alg}N$ is a quasi-isomorphism. Similarly, a map $f\co
  \lsup{\Alg}M_\Blg \to \lsup{\Alg}N_\Blg$ of type \DA\ modules is
  called a \emph{quasi-isomorphism} if the induced map $\Id_{\Alg}\DT
  f\co \lsub{\Alg}\Alg_\Alg\DT\lsup{\Alg}M_\Blg\to
  \lsub{\Alg}\Alg_\Alg\DT\lsup{\Alg}N_\Blg$ is a quasi-isomorphism,
  and a map $f\co \lsup{\Alg}M^\Blg \to \lsup{\Alg}N^\Blg$ of type
  \DD\ modules is called a \emph{quasi-isomorphism} if the induced map
  $\Id_{\Alg}\DT f\DT\Id_\Blg\co
  \lsub{\Alg}\Alg_\Alg\DT\lsup{\Alg}M^\Blg\DT\lsub{\Blg}\Blg_\Blg\to
  \lsub{\Alg}\Alg_\Alg\DT\lsup{\Alg}N^\Blg\DT\lsub{\Blg}\Blg_\Blg$ is a
  quasi-isomorphism.
\end{definition}
\begin{corollary}\label{cor:quasi-iso-is-htpy-equiv-D-DA-DD}
  Let $\lsup{\Alg}M,\lsup{\Alg}N\in\lsupv{\Alg}\ModCat$. A map $f\co
  \lsup{\Alg}N\to\lsup{\Alg}N$ is a quasi-isomorphism if and only if
  it is a homotopy equivalence. The same holds if $M$ and $N$ are
  instead in $\lsub{\Alg}\ModCat_\Blg$, $\lsupv{\Alg}\ModCat_\Blg$ or
  $\lsupv{\Alg}\ModCat^\Blg$.
\end{corollary}
\begin{proof}
  For type $D$ structures, this is immediate from
  Proposition~\ref{prop:D-to-A-and-back} and
  Proposition~\ref{prop:derived-is-derived-is-derived-is-derived}. The
  bimodule analogues follow from these results together with
  Proposition~\ref{prop:bimod-vs-ABop-mod}, below.
\end{proof}

\subsubsection{Induction and restriction} 
\label{sec:induct-restr}
Consider a map of $\Ainf$ algebras $\phi\co \Alg\to
\Blg$. Associated to $\phi$ are restriction and induction functors
\begin{align*}
\Rest_{\phi}&\co \ModCat_\Blg\to \ModCat_\Alg\\
\lsupv{\phi}\Induct&\co \lsupv{\Alg}\ModCat\to 
\lsupv{\Blg}\ModCat
\end{align*}
defined by $\cdot\DT \lsupv{\Blg}[\phi]_{\Alg}$ and
$\lsupv{\Blg}[\phi]_{\Alg}\DT\cdot$ respectively, where
$\lsupv{\Blg}[\phi]_\Alg$ is as defined in
Definition~\ref{def:rank-1-DA-mods}.\footnote{Since
$\lsupv{\Alg}\ModCat$ and $\lsupv{\Blg}\ModCat$ are
$\Ainf$-categories, the induction functor
is, of course, an $\Ainf$-functor; if $\Alg$ and $\Blg$ are \dg
algebras then this complication disappears.} One can define
restriction functors of left modules and induction functors of right
type $D$ structures similarly, using $\lsub{\Alg}[\phi]^\Blg$ instead
of $\lsupv{\Blg}[\phi]_{\Alg}$.

It is obvious from their definitions that these functors behave well
with respect to composition of algebra homomorphisms:
\begin{lemma}
  If $\phi\co \Alg\to \Blg$ and $\psi\co
  \Blg\to {\mathcal C}$ are $\Ainf$ algebra homomorphisms then
  \begin{align*}
    \Rest_{\phi}\circ\Rest_{\psi} &= \Rest_{\psi\circ \phi}\\
    \lsupv{\psi}\Induct\circ\lsupv{\phi}\Induct &=\lsupv{\psi\circ\phi}\Induct.
  \end{align*}
\end{lemma}

For most of the rest of this section we restrict to the case that
$\Alg$ and $\Blg$ are \dg algebras.
\begin{lemma}\label{lem:ind-rest-fund-lem}
  Let $\phi\co\Alg\to\Blg$ be an $\Ainf$ morphism.
  Then there is a natural map of $\Alg$-$\Alg$ type \AAm\ bimodules
  \[
  \lsub{\Alg}\Alg_{\Alg}\to \lsub{\Alg}[\phi]^{\Blg}\DT 
  \lsub{\Blg}\Blg_{\Blg}\DT \lsupv{\Blg}[\phi]_{\Alg}.
  \]
  Similarly, if $\Alg$ and $\Blg$ are \dg algebras then there is a natural map of
  $\Blg$-$\Blg$ type \DD\ bimodules
  $$\lsupv{\Blg}[\phi]_{\Alg}\DT \lsup{\Alg}\Barop(\Alg)^{\Alg}\DT \lsub{\Alg}[\phi]^\Blg\to \lsup{\Blg}\Barop(\Blg)^{\Blg}.$$

  When $\phi$ is a quasi-isomorphism, then these two natural maps
  are quasi-isomorphisms.
\end{lemma}
\begin{proof}
  Define a map $f\co \lsub{\Alg}\Alg_\Alg\to\lsub{\Alg}[\phi]^{\Blg}\DT 
  \lsub{\Blg}\Blg_{\Blg}\DT \lsupv{\Blg}[\phi]_{\Alg}$ by
  \[
  f_{i,1,j}(a_1,\dots,a_i,a',a''_1,\dots,a_j'')=1_\Ground\otimes \phi_{i+1+j}(a_1,\dots,a_i,a',a''_1,\dots,a_j'')\otimes 1_\Ground.
  \]
  Using the fact that
  \[
  m_{[\phi]\DT\Blg\DT[\phi]}(a_1,\dots,a_i,1_\Ground\otimes b\otimes
  1_\Ground,a''_1,\dots,a''_j)=1_\Ground\otimes
  \mu_\Blg(F^\phi(a_1,\dots,a_i),b,F^\phi(a''_1,\dots,a''_j))\otimes 1_\Ground
  \]
  it is easy to verify that $f$ is an $\Ainf$-bimodule homomorphism
  (i.e., $d(f)=0$). Furthermore, if $\phi_1$ induces an isomorphism on homology,
  $f_{0,1,0}$ induces an isomorphism on homology, i.e., $f$ is a
  quasi-isomorphism.

  For the second part, define a map $g\co
  \lsupv{\Blg}[\phi]_{\Alg}\DT
  \lsup{\Alg}\Barop(\Alg)^{\Alg}\DT
  \lsub{\Alg}[\phi]^\Blg\to
  \lsup{\Blg}\Barop^{\Blg}$ of separated type \DD\ modules by
  \[
  g(1_\Ground\otimes[a_1|\cdots|a_k]\otimes
  1_\Ground)=1_\Ground\otimes [F^\phi(a_1\otimes\cdots\otimes
  a_k)]\otimes 1_\Ground.
  \]
  Then it is easy to see that $g$ is a type \DD\ homomorphism (i.e.,
  $d(g)=0$). Moreover, if $\phi$ is a quasi-isomorphism then $F^\phi$
  is a quasi-isomorphism, and hence $g$ is a quasi-isomorphism.
\end{proof}

\begin{definition}
  \label{def:QuasiInvertible}
  Fix $\Ainfty$ algebras $\Alg$ and $\Blg$. A \emph{quasi-inverse} to
  a type \DA\ bimodule 
  $\lsupv{\Alg}P_{\Blg}$ is a type \DA\ bimodule 
  $\lsupv{\Blg}Q_{\Alg}$ with the property that
  $\lsupv{\Alg}P_{\Blg}\DT\lsupv{\Blg}Q_{\Alg}\simeq
  \lsupv{\Alg}\Id_{\Alg}$  and
  $\lsupv{\Blg}Q_{\Alg}\DT\lsupv{\Alg}P_{\Blg}\simeq
  \lsupv{\Blg}\Id_{\Blg}$. A type \DA\ bimodule is
  \emph{quasi-invertible} if it has a quasi-inverse.

  There are obvious analogous definitions if $P$ is a type \AAm\ and
  $Q$ is a type \DD\ module, and if $P$ and $Q$ are both type \AAm\
  modules (using $\DTP$ instead of $\DT$).
\end{definition}

Although a quasi-isomorphism between $\Ainf$-algebras need not be
invertible, it does have a quasi-inverse:

\begin{proposition}
  \label{prop:quasi-isom-invert-bimodule}
  If $\phi\co \Alg\to \Blg$ is a quasi-isomorphism,
  then $\lsupv{\Alg}[\phi]_{\Blg}$ is quasi-invertible.
\end{proposition}
\begin{proof}
  We will see that the quasi-inverse to $M=\lsupv{\Alg}[\phi]_{\Blg}$ is
  $N=\lsup{\Blg}\Barop(\Blg)^\Blg \DT\lsub{\Blg}[\phi]^\Alg
  \DT\lsub{\Alg}\Alg_\Alg$. Indeed,
  \[
  N\DT M\cong \lsup{\Blg}\Barop(\Blg)^\Blg \DT\lsub{\Blg}[\phi]^\Alg
  \DT\lsub{\Alg}\Alg_\Alg\DT\lsupv{\Alg}[\phi]_\Blg\simeq
  \lsup{\Blg}\Barop(\Blg)^\Blg \DT\lsub{\Blg}\Blg_\Blg
  \simeq \lsupv{\Blg}[\Id]_\Blg
  \]
  where the first quasi-isomorphism uses
  Lemma~\ref{lem:ind-rest-fund-lem} and the second uses
  Lemma~\ref{lem:bar-res}.

  Similarly,
  \[
  M\DT N\cong \lsupv{\Alg}[\phi]_\Blg\DT\lsup{\Blg}\Barop(\Blg)^\Blg
  \DT\lsub{\Blg}[\phi]^\Alg \DT\lsub{\Alg}\Alg_\Alg\simeq
  \lsup{\Alg}\Barop(\Alg)^\Alg\DT\lsub{\Alg}\Alg_\Alg\simeq \lsupv{\Alg}[\Id]_\Alg
  \]
  where again the first quasi-isomorphism uses
  Lemma~\ref{lem:ind-rest-fund-lem} and the second uses
  Lemma~\ref{lem:bar-res}.
\end{proof}

One reason to be interested in quasi-invertible bimodules is the
following:
\begin{lemma}\label{lem:quasi-invert-gives-equiv-DA}
  If $\lsup{\Blg}M_\Alg$ is quasi-invertible then the functors
  $\lsup{\Blg}M_\Alg\DT\cdot$ and $\cdot\DT\lsup{\Blg}M_\Alg$ are
  quasi-equivalences of categories. Analogous statements hold for type
  \DD\ and \AAm\ modules.  An analogous statement holds for type \AAm\
  modules with $\DTP$ in place of $\DT$.
\end{lemma}
The proof is straightforward.

We have now proved the \dg case of the following proposition,
which we state merely for the reader's edification (compare
\cite[Theorem 10.12.5.1]{BernsteinLunts94:EquivariantSheaves}):
\begin{proposition}\label{prop:quasi-iso-equivalence}
  Let $\Alg$ and $\Blg$ be $\Ainf$-algebras, and
  $\phi\co\Alg\to\Blg$ a quasi-isomorphism. 
  Let
  $\Induct_\phi: \ModCat_\Alg \to \ModCat_\Blg$ denote the functor
  \[
  M_\Alg\mapsto \Induct^\phi(M_\Alg\DT\lsup{\Alg}\Barop(\Alg)^\Alg)^\Blg\DT\lsub{\Blg}\Blg_\Blg.
  \]
  Then $\Rest_\phi$ and $\Induct_\phi$ are inverse quasi-equivalences
  of \dg categories. In particular, the derived categories
  $\mathsf{H}(\ModCat_\Alg)$ and $\mathsf{H}(\ModCat_\Blg)$ are
  equivalent.  Corresponding statements apply to categories of type
  $D$ structures and bimodules of all kinds.
\end{proposition}
The argument we have given works in general, with the only obstruction
being the fact that we have not defined type \DD\ modules over
$\Ainf$-algebras. This is not a serious difficulty for the purpose of
this result; cf.\ Remark~\ref{remark:bar-in-general}.

\subsubsection{Bimodules and \textalt{$(\Alg^{\op}\otimes \Blg)$}{Aop tensor B}-modules.}
\label{sec:A-Bop-vs-bimods}
Let $\Alg$ and $\Blg$ be differential graded algebras.
An honest $(\Alg,\Blg)$-bimodule is exactly the same as a right $(\Alg^\op\otimes
\Blg)$-module. For $\Ainf$-bimodules, this is not quite
true: a $(\Alg^\op\otimes \Blg)$-module has distinct operations $m_3(x,\overline{\gamma}\otimes
1,1\otimes\eta)$ and
$m_3(x,1\otimes\eta,\overline{\gamma}\otimes 1)$, both of which
intuitively correspond to the operation $m_3(\eta,x,\gamma)$ on an
$\Ainf$-bimodule.

Nevertheless, the two categories are equivalent.  More precisely, let
$M$ be a right $\Ainf$ $(\Alg^\op\otimes
\Blg)$-module. We define an $\Ainf$ $(\Alg,\Blg)$-bimodule structure on
$M$ as follows. Given sequences $S_\Alg=(a_1,\dots,a_m)$ of elements of $\Alg$
and $S_\Blg=(b_1,\dots,b_n)$ of elements of $\Blg$, we say that a sequence
$S_{\Alg^\op\otimes \Blg}=(c_1,\dots,c_{m+n})$ of elements of $\Alg^\op\otimes \Blg$
\emph{interleaves} $S_\Alg$ and $S_\Blg$ if:
\begin{itemize}
\item each $c_k$ is either $\overline{a}_i\otimes 1$ or $1\otimes b_i$,
\item the sequence obtained from $\{c_i\}$ by forgetting the
  $1\otimes b_i$'s is exactly
  $(\overline{a}_1\otimes1,\dots,\overline{a}_m\otimes1)$, and
\item the sequence obtained from $\{c_i\}$ by forgetting the
  $\overline{a}_i\otimes 1$'s is exactly
  $(1\otimes b_1,\dots,1\otimes b_n)$.
\end{itemize}

The $(\Alg,\Blg)$-bimodule structure on $M$ is defined by
\[
m_{m,1,n}(a_1,\dots,a_m,x,b_1,\dots,b_n)=\qquad \mathclap{\sum_{\substack{(c_1,\dots,c_{m+n})\text{
    interleaves} \\(a_1,\dots,a_m)\text{ and } (b_1,\dots,b_n)}}}\qquad  m_{m+n+1}(x,c_1,\dots,c_{m+n}).
\]

It is routine to verify that these higher products do, in fact, make
$M$ into an $(\Alg,\Blg)$-bimodule. Moreover, this construction extends in
an obvious way to maps, leading to a functor from $\Mod_{\Alg^\op\otimes
  \Blg}$ to $\lsub{\Alg}\Mod_\Blg$.

\begin{proposition}\label{prop:bimod-vs-ABop-mod}
  The categories $\HMod(\ModCat_{\Alg^\op\otimes \Blg})$ and
  $\HMod(\lsub{\Alg}\ModCat_\Blg)$ are equivalent triangulated categories.
\end{proposition}
\begin{proof}
  This follows from
  Proposition~\ref{prop:derived-is-derived-is-derived-is-derived} for
  $(\Alg^\op\otimes \Blg)$ modules and for $(\Alg, \Blg)$ bimodules.
\end{proof}

\subsection{Group-valued gradings}
\label{sec:algebras-gradings}
The gradings on Floer homology theories differ from gradings in
classical homology or homological algebra in (at least) three
important ways:
\begin{enumerate}
\item\label{item:Floer-gr-1} It is often easiest to consider
  Floer complexes as relatively graded groups, rather than absolutely graded ones.
\item\label{item:Floer-gr-2} The relative grading is usually only
  partially defined. That is, Floer chain complexes break up
  as direct sums in which there is a relative grading on each summand
  but no way to compare the gradings across summands.
\item\label{item:Floer-gr-3} The relative gradings on Floer complexes are often
  cyclic, i.e., by $\ZZ/n$ rather than by $\ZZ$.
\end{enumerate}
All three points lead to difficulties with homological algebra. In
particular, there is no notion of a degree $0$ morphism between two different
relatively graded modules, nor is it clear how to define the cone of a
morphism between relatively graded modules.

In the literature, points~(\ref{item:Floer-gr-1})
and~(\ref{item:Floer-gr-2}) are usually treated by fixing a lift of
the (partially defined) relative grading to an absolute grading, or
quantifying over all such
lifts~\cite{Seidel00:Graded}. (Alternatively, for Heegaard Floer
theory of rational homology spheres there is a natural absolute
$\QQ$-grading lifting the relative
$\ZZ$-grading~\cite{OS06:HolDiskFour}.) In papers where
point~(\ref{item:Floer-gr-3}) leads to difficulties, authors usually
either restrict the generality of their results or work with a
periodic $\ZZ$-graded lift of the $\ZZ/n$-graded module (although
this, too, leads to difficulties with homological algebra).

None of these approaches seem satisfactory for bordered Floer
theory. Indeed, we observed in \cite{LOT1} that there is a natural grading on
the algebras involved in bordered Floer theory by non-commutative
groups $G$, and on the modules involved by $G$-sets. This section is
devoted to reviewing and expanding the algebraic framework of such
gradings. While group-valued gradings have occurred before in the
literature~\cite{GroupGradings, Khovanov05:Hopfological}, our
perspective and goals seem somewhat different.

Partially-defined relative $\ZZ$- or
$\ZZ/n$-gradings are a special case of the construction here. More
dramatically, cyclic gradings arise naturally when taking tensor
products of non-cyclically-graded modules in the $G$-set graded
context.

This section is organized as follows. In
Section~\ref{sec:G-graded-basics} we review the basic definitions of
$G$-valued gradings from~\cite{LOT1}, and extend the theory to
bimodules. Section~\ref{sec:G-graded-basics} should provide enough
background on non-commutative gradings for the reader with a little
faith to understand most of the rest of this paper, except for
Section~\ref{sec:Duality}. Section~\ref{sec:set-graded-dg-cats}
extends the notion of \dg categories as appropriate for categories of
$G$-set graded modules; the generalization is called a $\ZZ$-set
graded \dg category. Section~\ref{sec:G-set-mod-cats} then
organizes the $G$-set graded modules into $\ZZ$-set-graded \dg
categories. This then allows one to extend easily the homological
algebra introduced earlier in this section to $G$-set graded modules.

\subsubsection{Basics of group-valued gradings}\label{sec:G-graded-basics}
We start by recalling some notions of group-valued gradings of
algebras and modules from~\cite{LOT1}, and then generalize them
somewhat. Note that all of our $\Ainf$-algebras and modules have
underlying vector spaces (or $\Ground$-modules), and all of the
structure maps of $\Ainf$-algebras and modules are maps between tensor
products of vector spaces. So, to grade these modules and speak about
the degrees of structure maps it suffices to explain gradings of
vector spaces and tensor products of vector spaces. We do so as
follows:
\begin{definition}\label{def:G-set-graded-k-mod}
  Let $(G,\lambda)$ be a pair of a group~$G$ and a distinguished
  element~$\lambda$ in the center of~$G$.  A
  \emph{$G$-graded $\Ground$-bimodule} is a $\Ground$-bimodule $V$
  which is decomposed as a direct sum
  \[
  V=\bigoplus_{g\in G} V_g.
  \]
  We say that an element $v \in V_g$ is \emph{homogeneous of
    degree~$g$} and write $\gr(v) = g$.

  If $V$ and $W$ are two $G$-graded $\Ground$-bimodules then $V
  \otimes W$ is itself $G$-graded by $\gr(v \otimes w) =
  \gr(v)\gr(w)$.

  For a $G$-graded $\Ground$-bimodule~$V$, the space $V[n]$ is a
  $G$-graded $\Ground$-bimodule with gradings shifted by $\lambda^n\co
  V[n]_g = V_{\lambda^{-n} g}$.

  A homomorphism $f\co V\to W$ between $G$-graded $\Ground$-bimodules
  is \emph{homogeneous of degree $k\in\ZZ$} if for all $g\in G$,
  $f(V_g)\subset W_{g\lambda^k}$.
\end{definition}

Using Definition~\ref{def:G-set-graded-k-mod}, the definitions of
$\Ainf$-algebras and 
$\Ainf$-algebra homomorphisms (Sections~\ref{sec:defin-Ainfty-alg}
and~\ref{sec:defin-Ainf-alg-maps}) carry over to the $G$-graded case
without change. For example:
\begin{definition} \label{def:ainf-nc-grading}
  Let $(G,\lambda)$ be a group with a distinguished central element.
  An \emph{$\Ainf$-algebra graded by~$(G,\lambda)$} is an
  $\Ainf$-algebra~$\Alg$ with a grading $\gr$ by $G$ (as a
  $\Ground$-bimodule),
  i.e., a decomposition $\Alg=\bigoplus_{g\in G} \Alg_g$, satisfying
  the following condition:
  for homogeneous elements $a_i$, we require that
  \begin{equation}
    \label{eq:ainf-nc-grading}
    \gr(\mu_j(a_1,\dots,a_j)) = \gr(a_1)\cdots\gr(a_j)\lambda^{j-2}.
  \end{equation}
\end{definition}
\begin{example}
  In the case that $G=\ZZ$ and $\lambda=1$,
  Definition~\ref{def:ainf-nc-grading} reduces to a classical
  ($\ZZ$-graded) $\Ainf$-algebra.
\end{example}

If $\Alg$ is $G$-graded, we could consider modules that are also
$G$-graded.  However, we prefer to consider modules which are
graded by $G$-sets. Again, it suffices to explain the notion of
$\Ground$-modules graded by $G$-sets:
\begin{definition}\label{def:G-set-gr-k-bimod}
  Let $(G,\lambda)$ be a group with a distinguished central element,
  and let $S$ be a right $G$-set.  An \emph{$S$-graded
    $\Ground$-module} is a $\Ground$-module $V$ which is decomposed as
  a direct sum
  \[
  M=\bigoplus_{s\in S}V_s.
  \]

  If $V\!$ is an $S$-graded
  $\Ground$-module and $W\!$ is a $G$-graded $\Ground$-bimodule, then $V
  \otimes W$ is itself $S$-graded by $\gr(v \otimes w) =
  \gr(v)\gr(w)$.  Similarly if $T$ is a left $G$-set and $V\!$ is a
  $T$-graded $\Ground$-module, then $W \otimes V$ is $T$-graded by
  $\gr(w \otimes v) = \gr(w) \gr(v)$.

  For any module $V\!$ graded by a set with an action of $\lambda$, the
  space $V[n]$ is $V\!$ with shifted grading: $V[n]_x = V_{\lambda^{-n}
    x}$.  Note that we need not distinguish between left and right
  actions of~$\lambda$ since $\lambda$ is central.
\end{definition}

Via Definition~\ref{def:G-set-gr-k-bimod} the definitions of
$\Ainf$-modules and type $D$ structures
(Sections~\ref{sec:categ-Ainf-modules} and~\ref{sec:cat-type-d-str})
carry over to the $G$-set graded case without change. For example:
\begin{definition}\label{def:ainf-mod-nc-grading}
  For $(G,\lambda)$ a group with a distinguished central element,
  $\Alg$ a $G$-graded $\Ainf$-algebra, and $S$ a right $G$-set,
  a \emph{right $S$-graded $\Ainf$-module} is an $\Ainf$-module
  $M_\Alg$ whose underlying $\Ground$-module is
  graded by~$S$, such that for homogeneous elements $x\in M$ and $a_i \in A$,
\[
\gr(m_{j+1}(x, a_1, \dots, a_j)) = \gr(x)\cdot\lambda^{j-1}\gr(a_1)\cdots\gr(a_j).
\]
\end{definition}

\begin{example}
  The case that $G=\ZZ$, $\lambda=1$, and $S$ is a freely transitive $G$-set
  is equivalent to relatively $\ZZ$-graded modules. In
  particular, given a transitive $\ZZ$-set graded module $M$ one can define a
  relative grading on $M$ by $\gr(x,y)=n$ if $x$ and $y$ are
  homogeneous and $\gr(x)\lambda^n=\gr(y)$. Conversely, given a
  relatively $\ZZ$-graded module $M$ it is not hard to construct a
  canonical $\ZZ$-set~$S$ and an associated $S$-graded module.
\end{example}

\begin{definition}
  \label{def:EssentialGradingSet}
  For $S$ a $G$-set and $M$ an $S$-graded module, we say that a
  $G$-invariant subset $T\subset S$ is
  \emph{essential} if
  each $G$-orbit in~$T$ contains an element~$s$ so that $M_s$ is
  nontrivial.
\end{definition}

The theory of $G$-set gradings has a simple extension to bimodules:
\begin{definition}
  \label{def:GradingBimodules}
  If $(G_1, \lambda_1)$ and $(G_2, \lambda_2)$ are two groups with
  distinguished central elements, define the group
  \[
  G_1 \times_\lambda G_2 \coloneqq
  G_1 \times G_2 / (\lambda_1 = \lambda_2)
  \]
  with the distinguished central element $\lambda = [\lambda_1] =
  [\lambda_2]$.  So, a $(G_1 \times_\lambda G_2)$-set is a set
  with commuting actions of $G_1$ and $G_2$, where the actions of
  $\lambda_1$ and $\lambda_2$ agree.

  A \emph{left-right $(G_1,G_2)$-set graded $\Ground$-module} is a left
  $(G_1\times_\lambda G_2^{op})$-set graded $\Ground$-module. (Left-left
  and right-right graded modules are defined similarly.)
\end{definition}
Via Definition~\ref{def:GradingBimodules} the definitions of bimodules
of various types (Section~\ref{sec:bimod-var-types}) carry over to the
$G$-set graded case without change: if $\Alg$ is $G_1$-graded and
$\Blg$ is $G_2$-graded, then we consider bimodules that are left-right
$(G_1,G_2)$-set graded.

We next turn to tensor products of $G$-set graded modules.

\begin{definition}
  \label{def:ProductGSet}
  If $(G,\lambda)$ is a group with a distinguished central element,
  $S$ is a right $G$-set, and $T$ is a left $G$-set, define their 
  {\em twisted product}
  \[
  S \times_G T \coloneqq (S \times T)/(s\times gt) \sim (sg\times t).
  \]
  The set $S \times_G T$ has an action of $\lambda$ defined by $\lambda
  \cdot [s\times t] \coloneqq [s\lambda\times  t] = [s\times \lambda t]$, but in
  general has no further structure.
  For $V$ and $W$ two $\Ground$-modules graded by $S$ and $T$,
  respectively, $V \otimes W$ is graded by the $\ZZ$-set $S \times_G T$, where if
  $v\in V_s$ and $w\in W_t$, then $v \otimes w \in (V \otimes
  W)_{[s\times t]}$.

  Thus, for instance, if $M_\Alg$ is graded by~$S$ and $\lsup{\Alg}N$
  is graded by~$T$, $M_\Alg \DT \lsup{\Alg}N$ is graded by $S \times_G
  T$.

  More generally, for the tensor product of bimodules, if $G_1, G_2,
  G_3$ are all groups with distinguished central elements, $S$ is a
  left-right $(G_1, G_2)$-set, and $T$ is a left-right $(G_2,
  G_3)$-set, then $S\times_{G_2} T$ is a left-right $(G_1, G_3)$-set in an
  obvious way.  As before, if $V$ and $W$ are two spaces graded by $S$
  and~$T$, respectively, $V\otimes W$ is graded by $S \times_{G_2} T$.
  In particular, if $\lsub{\Alg} M_\Blg$ and $\lsup{\Blg}N_\Clg$ are two
  left-right set-graded $\Ainf$-bimodules, then $\lsub{\Alg}M_\Blg \DT
  \lsup{\Blg}N_\Clg$ is also a left-right set-graded $\Ainf$-bimodule.
\end{definition}

Note that the complex $M_\Alg \DT \lsup{\Alg}N$ is not in general
$\ZZ$-graded; rather, the action of $\lambda$ on $S\times_G T$ breaks
up the morphism space into a sum of chain complexes according to the
orbits of the action, with possibly cyclic grading on each orbit.  (See
\cite[Example~\ref*{LOT:ex:order-tensor}]{LOT1} for an example where the
grading ends up being a finite cyclic group, even though the action of $\lambda$ on
each side has infinite order.)

The grading on the Hochschild complex of a $G$-set graded bimodule
behaves as one would expect:
\begin{definition}\label{def:Hoch-grading}
  For $S$ a left-right $(G,G)$-set, define
  \[
  S^\circ \coloneqq S/(s\cdot g)\sim (g\cdot s)
  \]
  (where $s \in S$, $g \in G$).  $S^\circ$ has a $\ZZ$-action given by
  multiplication by~$\lambda$ on the left or right. 

  There is a quotient map from $S$ to $S^\circ$, which we denote
  $s\mapsto s^\circ$.
\end{definition}

\begin{lemma}\label{lem:Hoch-grading}
  Let $\lsup{\Alg}M_\Alg$ be a type \DA\ bimodule graded by a
  left-right $(G,G)$-set $S$. Then $\tCH(M)$ is graded by $S^\circ$. The
  analogous statement holds for the Hochschild complex of a type
  \AAm\ bimodule.
\end{lemma}

There are two ways to form $G$-set graded modules into categories. The
easier of the two, which will suffice for most of the paper, is to
define a category for each transitive $G$-set $S$, and to restrict to
morphisms which only shift the grading by a power of~$\lambda$:
\begin{definition}\label{def:cat-fixed-gr-set}
  Fix $(G,\lambda)$ a group with a distinguished central element and
  $S$ a right $G$-set where the action of $G$ is transitive. Let
  $n$ be the order of $\lambda$ in~$S$.  (This order is
  well-defined since the $G$ action is transitive.)  For $V$ and $W$
  two $S$-graded $\Ground$-modules, let $\widetilde\Hom_S(V,W)$ be the
  $(\ZZ/n\ZZ)$-graded $\Field$ vector space with
  \[
  \widetilde\Hom_S(V,W)_i \coloneqq
    \bigoplus_{s\in S} \Hom(V_s, W_{s \lambda^i}).
  \]
  These can be organized into a category $\widetilde\ModCat_{\Ground,S}$, whose
  objects are $S$-graded $\Ground$-modules, and whose morphism spaces are
  $\widetilde\Hom_S(V,W)$. 

  For $\Alg$ a $G$-graded $\Ainf$-algebra, the
  category $\widetilde\ModCat_{\Alg,S}$ of right $S$-graded $\Ainf$-modules over
  $\Alg$ is the \dg category whose objects are right $S$-graded
  $\Ainf$-modules $M_\Alg$ and whose morphism spaces are
  $$\widetilde\Mor_{\Alg,S}(M_{\Alg},N_{\Alg})=\widetilde\Hom_S(M\otimes \Tensor^*(A_+[1]),N),$$
  with the differential as in
  Definition~\ref{def:Ainf-mod-homs}.
  (Here we have extended the notion of \dg categories in the obvious way
  to allow cyclic gradings.)

  The definitions in Section~\ref{sec:dg-categories} of homotopic
  morphisms and homotopy equivalences carry over unchanged.  We
  extend the definitions of $\ZMod(\widetilde\ModCat_{\Alg,S})$ and
  $\HMod(\widetilde\ModCat_{\Alg,S})$ in the obvious ways, and
  $\HMod(\widetilde\ModCat_{\Alg,S})$
  is still a triangulated category.

  We define categories of other kinds of modules in the $G$-set graded case
  analogously. For instance we have a category 
  $\lsupv{\Alg}{\widetilde\ModCat}_S$ of left $S$-graded type $D$
  structures over~$\Alg$, as well as categories of bimodules of
  various types.
\end{definition}

The invariants $\CFDDa(Y,\spinc)$, $\CFAAa(Y,\spinc)$ and
$\CFDAa(Y,\spinc)$ defined in
Section~\ref{sec:CFBimodules} will be well-defined up to isomorphism
in $\HMod(\widetilde\ModCat_{S})$ (for appropriate
grading sets $S$).

\subsubsection{Set-graded \dg
  categories}\label{sec:set-graded-dg-cats}
It is natural to talk about morphisms between modules
graded by different $G$-sets, or in other words to collect the
$S$-graded modules for all $G$-sets $S$ into a single category. In
doing so, one encounters the following issues:
\begin{itemize}
\item There is no natural way to define a degree $0$ morphism from an
  $S$-graded module to a $T$-graded module if $S$ is different from
  $T$.
\item  More generally, if $M$ and $N$ are graded by $G$-sets $S$ and $T$
  respectively then, like $M\otimes N$, the collection of morphisms
  $\Mor(M,N)$ is graded by a $\ZZ$-set constructed from $S$ and
  $T$ rather than simply by $\ZZ$.
\item  Further, with this construction, there are natural morphisms,
  like the identity morphism, which
  are not homogeneous (not supported in a single grading).
\end{itemize}

We will formalize these notions as follows. We will define a category
of $\ZZ$-set graded chain complexes. If $M$ and $N$ are graded by
$G$-sets $S$ and $T$ respectively then the morphisms from $M$ to $N$
will be a $\ZZ$-set graded chain complex.  (That is, the categories of
$G$-set graded modules of various kinds are enriched over the category
of $\ZZ$-set graded chain complexes.)

\begin{definition}
  \label{def:compose-rels}
  Let $S$ and $T$ be two sets, each endowed with a $\ZZ$-action
  (written as multiplication by $\lambda$).  A {\em relation} between
  $S$ and $T$ is a subset $R\subset S\times T$. We write $sRt$ to mean
  that $(s,t)\in R$.  We say that $R$ is {\em $\lambda$-invariant} if it
  satisfies the property that $sRt$ if and
  only if $(\lambda s) R (\lambda t)$. The {\em composite} of two
  relations $R_1\subset S \times T$ and $R_2 \subset T \times U$,
  written $R_1\circ R_2\subset S \times U$, is the relation
  \[ R_1 \circ R_2 = \{\,(s,u) \mid s \in S, u \in U,
  \exists t \in T: sR_1t
  \text{ and } tR_2u\}. \]
\end{definition}

\begin{definition}\label{def:set-graded-cx}
  A \emph{$\ZZ$-set graded chain complex} over $\Ground$ is a triple $(S,
  C, \partial)$, where $S$ is a set with a $\ZZ$ action (written as
  multiplication by~$\lambda$), $C$ is a $\Ground$-module graded
  by~$S$, and $\partial$ is an operator on $C$ satisfying $\partial^2 = 0$
  and such that for $x \in C_s$, $\partial x \in C_{\lambda^{-1}s}$.  We form
  $\ZZ$-set graded chain complexes into a
  category as follows. A \emph{morphism} from $(S,C,\partial_1)$ to $(T,C',\partial_2)$ is
  \begin{itemize}
  \item a $\lambda$-invariant relation $R$ between $S$ and $T$, and
  \item a $\Ground$-module map $\phi \co C \to C'$ that is compatible
    with~$R$, in the sense that
    \[
    \phi(C_s) \subset
      \bigoplus_{\substack{t \in T\\ s R t}}
      C'_t.
    \]
  \end{itemize}
  Let $\Mor((S,C,\partial_1),(T,C',\partial_2))$ denote this space of
  morphisms, with its natural differential.

  If in addition $\phi$ intertwines the actions of $\partial_1$ and
  $\partial_2$, we say that $\phi$ is a \emph{homomorphism}.  
  To compose two morphisms, we compose the relations in the
  sense of Definition~\ref{def:compose-rels} and compose the
  $\Ground$-module maps.

  The category of $\ZZ$-set graded chain complexes has objects triples
  $(S,C,\bdy)$ as above and $\Hom\left((S,C,\bdy),(T,D,\bdy')\right)$
  the set of homomorphisms (not just morphisms) from $(S,C,\bdy)$
  to $(T,D,\bdy')$.

  Given two $\ZZ$-set graded complexes $(S,C,\partial_1)$ and $(T,D,\partial_2)$,
  we can form their {\em tensor product}:
  \[
  (S, C, \partial_1) \otimes (T, D, \partial_2) \coloneqq
  (S \times_\ZZ T, C \otimes_\Ground D, \partial_1 \otimes \Id_D + \Id_C
  \otimes \partial_2),
  \]
  where $\gr(x \otimes y) = [\gr(x)\times \gr(y)]$ for homogeneous $x,
  y$.  This tensor product extends in an obvious way to
  (homo)morphisms, giving a monoidal structure on this category.
\end{definition}

\begin{example}
  Consider a $3$-manifold $Y$. The Heegaard Floer complex $\CFa(Y)$ is
  a set-graded chain complex where the grading set $S$ is
  (non-canonically) given by
  \[
  \bigcup_{\spinc\in\SpinC(Y)}\ZZ/\divis(c_1(\spinc)).
  \]
  (In fact, the grading set can be defined canonically.  First one
  fixes a Heegaard diagram and base generator and uses these to define
  a $G$-set grading as in Section~\ref{sec:cf-gradings}.  Then one
  observes that different base generator or Heegaard diagram lead to
  canonically isomorphic $G$-sets.)
\end{example}

The most useful notion of an element of a $\ZZ$-set graded chain
complex $(S, C, \partial)$ is not the naive notion of an element of
$C$: elements carry with them some additional grading information, as
follows.
\begin{definition}\label{def:elements}
  Let $\underline{\Ground}$ denote the $\ZZ$-set graded chain complex
  $(\ZZ,\Ground_0,0)$, where $\Ground_0$ denotes a copy of
  $\Ground$ lying in grading $0$.
  An \emph{element} of a $\ZZ$-set graded chain complex $(S,C,\partial)$ is a
  morphism of $\ZZ$-set graded chain complexes $\underline{\Ground}\to
  (S,C,\partial)$.  (Note that $\underline{\Ground}$ is the identity
  for the tensor product, and so this is natural from the point of
  view of category theory.)
\end{definition}

\begin{lemma}
  \label{lem:WhatIsAnElement}
  An element of $(S,C,\partial)$ (in the sense of Definition~\ref{def:elements})
  is equivalent to a pair $(T,x)$ where
  $T\subset S$ and $x\in \bigoplus_{s\in T} C_s$.
\end{lemma}

\begin{proof}
  A morphism $\underline{\Ground} \to (S,C,\partial)$ is given by a
  $\lambda$-invariant relation~$R$ between $\ZZ$ and $S$, 
  and a map~$\phi$.  Let $T$ be
  $\{\,t \in S \mid 0Rt\,\}$, the set of elements that are related
  to~$0$, and let $x=\phi(1)$. In the opposite direction, for any
  subset $T$ of $S$, there is a unique $\lambda$-invariant relation $R$
  between $\ZZ$ and $S$ with the property that $T$ is the set of elements related to~$0$;
  and similarly, any element $x\in C$ can be viewed as $\phi(1)$ for a uniquely determined
  $\phi\co\underline{\Ground}\to C$.
\end{proof}

For $T\subset S$ we let $C_T$ denote the set of elements of $C$ with
grading set $T$, i.e., the elements of the form $(T,x)$. For $s\in S$, then,
$C_{\{s\}}$ is isomorphic to $C_s$ as defined in
Definition~\ref{def:set-graded-cx}.

We can use the tensor product to define an internal $\Mor$ for
(appropriately finite) $\ZZ$-set graded chain complexes, just as for ordinary
chain complexes. We start by explaining $G$-set gradings of dual spaces:

\begin{definition}
  For $G$ any group and $S$ a right $G$-set, let $S^*$ be a set with
  elements $s^*$ in bijection with elements $s$ of~$S$,
  and with a left $G$-action defined by $g\cdot s^* \coloneqq (s \cdot
  g^{-1})^*$.  If $V$ is an $S$-graded
  $\Ground$-module, the \emph{graded dual} of $V$, which we will
  denote $V^{\gr *}$, is
  an $S^*$-graded module with
  \[
  (V^{\gr *})_{s^*} \coloneqq (V_s)^*.
  \]
  With these definitions, if $\phi: V \to V$
  is homogeneous with respect to the right action of $g$ on~$S$, then
  $\phi^*: V^* \to V^*$ is homogeneous with respect to the left action
  of $g$ on~$S^*$.
\end{definition}

\begin{definition}\label{def:set-graded-internal-hom}
  Suppose that $(S, C, \partial_1)$ and $(T, D, \partial_2)$ are
  set-graded chain complexes such that $C$ is supported on finitely
  many orbits of the action of~$\ZZ$ on~$S$.  Then we define a
  $\ZZ$-set graded chain complex $\Mor((S, C, \partial_1), (T,
  D, \partial_2))$ to be a chain complex graded by $T\times_\lambda S^*$ with
  \[
  \Mor((S, C, \partial_1), (T, D, \partial_2))_u \coloneqq
    \prod_{s^*\in S^*} \bigoplus_{\substack{t\in T\\u = [t\times s^*]}}
    \Hom(C_s, D_t).
  \]
  The differential on $\Mor((S, C, \partial_1), (T, D, \partial_2))$
  is given as usual by
  \[
  (\partial\phi)(x) = \partial(\phi(x)) + \phi(\partial x).
  \]
\end{definition}

Note that if $C$ is finite dimensional, then
\[
\Mor((S, C, \partial_1), (T, D, \partial_2)) \isom
(T, D, \partial_2) \otimes (S^*, C^{\gr *}, \partial_1^*).
\]

We have defined two different notions of ``a morphism'' between
$\ZZ$-set graded chain complexes. Fortunately, they agree:
\begin{lemma}
  The set of elements (Definition~\ref{def:elements}) of the internal 
  morphism complex $\Mor((S,C,\partial_1),(T,D,\partial_2))$
  (Definition~\ref{def:set-graded-internal-hom}) are in natural
  bijection with the morphisms of $\ZZ$-set graded complexes defined
  in Definition~\ref{def:set-graded-cx}.
\end{lemma}
\begin{proof}
  An element of the internal morphism complex gives a morphism as follows.
  Thanks to Lemma~\ref{lem:WhatIsAnElement}, we can think of
  an element of the internal morphism complex as a subset
  $R'$ of $T\times_{\lambda} S^*$, together with a choice of 
  $$\phi\in\bigoplus_{u\in R'\subset T\times_{\lambda} S^*} \Mor((S,C,\partial_1),(T,D,\partial_2))_{u}
  =
  \bigoplus_{u\in R'\subset T\times_{\lambda} S^*}
  \prod_{s^*\in S^*} \bigoplus_{\substack{t\in T\\u = [t\times s^*]}}
  \Hom(C_s, D_t).
  $$  
  This $\phi$, of course, can be thought of as a $\Ground$-module map
  from $C$ to $D$.  Now, subsets of $T\times_{\lambda} S^*$ are in
  one-to-one correspondence with subsets of $(T\times S)/\ZZ$, which
  in turn are in one-to-one correspondence with $\lambda$-invariant
  relations $R$ between $S$ and $T$. Under this correspondence, we think of
  $\phi$ as a $\Ground$-module map from $C$ to
  $D$ which is compatible with $R$ (in the sense of
  Definition~\ref{def:set-graded-cx}). 

  Conversely, given a morphism, an element of the morphism complex can
  be constructed by reversing the above process.
\end{proof}

We next turn to issues of injectivity.
The notion of injectivity for morphisms of $\ZZ$-set graded chain complexes 
takes into account also grading information, as follows:
\begin{definition}
  \label{def:Injective}
  A relation $R \subset S \times T$ is said to be \emph{injective} if
  for every $s \in S$, there is a $t \in T$ so that $s R t$
  and $s' R t$ does not hold for any other $s' \in S$.

  A morphism $(R,\phi)\in \Hom((S,C,\partial),(T,D,\partial'))$
  (in the sense of Definition~\ref{def:set-graded-cx})
  is said to be {\em homology injective} if $R$ is injective in the
  sense above and
  the element $\phi$, thought of as a homomorphism from $C$ to $D$,
  induces an injective map on homology.
\end{definition}

If $R$ is a relation between $S$ and $T$, and $\Sigma\subset S$, let
$\Sigma R\subset T$ denote the subset of all $t\in T$ with the
property that $sRt$ holds for some $s\in \Sigma$.  

\begin{lemma}\label{lem:injective-rel}
  A relation $R$ is injective in the sense of
  Definition~\ref{def:Injective}
  iff for any two subsets $\Sigma, \Sigma'\subset S$, $\Sigma R=
  \Sigma'R$
 implies that $\Sigma=\Sigma'$ (i.e., $R$ induces an injective
  function from subsets of~$S$ to subsets of~$T$).
\end{lemma}

\begin{proof}
  Suppose first that $R$ is injective in the sense of
  Definition~\ref{def:Injective}.  For $s\in S$, let $f(s) \in T$ be
  the element promised by that definition.  Then $R$ induces a
  bijection between $S$ and $f(S)$, and if $\Sigma R =
  \Sigma' R$, then
  $f(\Sigma) = \Sigma R \cap f(S) = \Sigma' R \cap
  f(S) = f(\Sigma')$, so $\Sigma = \Sigma'$.

  Conversely, suppose that $R$ induces an injective function on subsets
  of~$S$.  Then in particular, for each element $s\in S$, $SR \supsetneq (S
  \setminus \{s\})R$.  Let $f(s) \in T$ be an element of $SR$ that is
  not in $(S \setminus \{s\})R$.  Then $f(s)$ satisfies the conditions
  of Definition~\ref{def:Injective}.
\end{proof}

\begin{lemma}
  A homomorphism $(R,\phi)$ of $\ZZ$-set graded complexes is
  injective in the sense of Definition~\ref{def:Injective} if and 
  only if composition with $\phi$ induces an injection on the set of
  elements of $(S,C)$ (in the sense of Definition~\ref{def:elements}) 
  to the set of elements of $(T,D)$.
\end{lemma}
\begin{proof}  
  Let $(R,\phi)$ be a morphism from $(S,C)$ to $(T,D)$.  For an
  element $(\Sigma, x)$ of $(S,C)$ (where $\Sigma\subset S$ and $x\in
  C_{\Sigma}$), the image (element of $(T,D)$) under $(R,\phi)$ is the
  pair $(\Sigma R, \phi(x))$.

  From Lemma~\ref{lem:injective-rel},
  we see $(R,\phi)$ is injective on elements of $(S,C)$ of the form
  $(\Sigma,0)$ (with $\Sigma\subset S$) 
  if and only if $R$ is injective; and it is injective on more general
  elements if and only if $\phi$, thought of as a homomorphism from $C$ to 
  $D$, is injective.
\end{proof}

Since the set of morphisms between $\ZZ$-set graded chain complex
$(S,C,\partial_1)$ and $(T,D,\partial_2)$ (with $C$ supported on
finitely many $\ZZ$-orbits) is itself a chain complex, there is an
obvious notion of when two morphisms are homotopic.  In fact, the
definition of homotopic morphisms can be extended without difficulty
to the general case without the finiteness restriction on~$C$.

\begin{definition}
  A \emph{$\ZZ$-set graded \dg category} is a category~$\Cat$
  where the morphism spaces are $\ZZ$-set graded chain complexes over
  $\Field$ and composition of morphisms gives a $\ZZ$-set graded chain
  maps $\circ\co\Mor(y,z) \otimes \Mor(x,y) \to \Mor(x,z)$.
\end{definition}
\begin{example}
  \label{ex:dgcat}
  An ordinary \dg category $\Cat$ gives a $\ZZ$-set graded \dg
  category $\widetilde{\Cat}$ as
  follows. Let $\ob(\widetilde{\Cat})=\ob(\Cat)$.  For objects
  $M,N\in\ob(\widetilde{\Cat})$ define
  $\Mor_{\widetilde{\Cat}}(M,N)=(\ZZ,\Mor_{\Cat}(M,N),\partial)$ where
  $\partial$ is the differential on $\Mor_{\Cat}(M,N)$.
\end{example}

\begin{example}
  As a special case of Example~\ref{ex:dgcat}, 
  consider a $\ZZ$-graded $\Ainf$-algebra $\Alg$ and its
  $\ZZ$-set graded \dg category of right modules $\Mod_\Alg$. Given
  modules $M_\Alg$ and $N_\Alg$, the space $\Mor(M_\Alg,N_\Alg)$ is a
  $\ZZ$-set graded chain complex. As in Definition~\ref{def:elements},
  an element of $\Mor(M_\Alg,N_\Alg)$---i.e., a morphism from $M_\Alg$
  to $N_\Alg$---has some grading information built into it. 

  More precisely, an ``element'' of the morphism space $\Mor(M_\Alg,N_\Alg)$ consists of a
  pair $(S,\phi)$ where $\phi\co M\otimes T^*\Alg\to N$ is a map of
  $\Ainf$-modules and $S\subset \ZZ$ is such that if $x\in M_i$ then
  $\phi_1(x)\subset \bigoplus_{j\in S+i}N_j$, and similarly for higher
  products.
\end{example}

The definitions in Section~\ref{sec:dg-categories} of homotopic
morphisms and homotopy equivalences carry over unchanged to an
arbitrary $\ZZ$-set graded \dg category~$\Cat$. The
definitions of $\ZMod_*(\Cat)$ and $\HMod_*(\Cat)$ also carry over.
However, it no longer makes sense to talk about morphisms of
degree~$0$, so we no longer have $\ZMod(\Cat)$ or the triangulated category
$\HMod(\Cat)$.

\begin{definition}\label{def:set-gr-dg-functor}
  A \emph{$\ZZ$-set graded \dg functor} $F \co \Cat\to\Dat$ is a
  functor enriched so that for $x,y\in\Cat$, the map $F\co
  \Mor(x,y)\to\Mor(F(x),F(y))$ is a $\ZZ$-set graded chain complex
  homomorphism (Definition~\ref{def:set-graded-cx}).
\end{definition}

Lemma~\ref{lemma:dg-funct-blah},
Definitions~\ref{def:homotopic-functors} and~\ref{def:quasi-equiv},
and Proposition~\ref{prop:dg-homotopy-equiv-quasi-equiv} carry over almost
unchanged to the context of $\ZZ$-set graded \dg functors.  (In
Definition~\ref{def:quasi-equiv}, consider $\HMod_*(F)$ instead of $\HMod(F)$.)
The basic
example is provided by tensor products: see
Example~\ref{ex:DTPasDGfunctor}.

\subsubsection{Categories of \textalt{$G$}{G}-set graded modules}\label{sec:G-set-mod-cats}

We now define a category of $G$-set graded modules where the sets may
vary. 

We start by defining the grading set for the $\Hom$-space:
\begin{definition}
  If $(G,\lambda)$ is a group with a distinguished central element,
  $S$ and $T$ are two right $G$-sets, and $V$ and $W$ two
  $\Ground$-modules graded by $S$ and $T$, respectively, then 
  $\Hom_\Ground(V,W)$ is defined to be $V^{\gr *} \otimes W$ (as a
  $\Field$ vector space) with its $T \times_G S^*$
  grading.
\end{definition}
\begin{remark}
  One might imagine that the grading set for the $\Hom$ space between
  $S$- and $T$-graded modules should be the set of $G$-equivariant
  maps $S\to T$. However, in many situations (e.g., if $G=\ZZ$,
  $S=\ZZ/2$ and $T=\ZZ$) this prevents there being any
  homogeneous homomorphisms at all. Our philosophy is that any module map
  (between finite-dimensional modules, say) should be decomposable
  as a sum of homogeneous maps.
\end{remark}

\begin{definition}\label{def:cat-vary-gr-set}
  For $\Alg$ a $(G,\lambda)$-graded $\Ainf$-algebra over $\Ground$, let
  $\ModCat_\Alg$ be the set-graded \dg category whose objects are
  strictly unital, set-graded $\Ainf$-modules (i.e., pairs $(S, M_\Alg)$ where $M_\Alg$
  is graded by~$S$), and whose morphism spaces are
  \begin{equation}
  \Mor_\Alg((S,M_\Alg), (T,N_\Alg)) \coloneqq \Hom_{\Ground}(M \otimes
    \Tensor^*(A_+[1]), N),\label{eq:G-mod-mor-def}
  \end{equation}
  with a grading by $T \times_G S^*$ and differential as in
  Definition~\ref{def:Ainf-mod-homs}. (The right hand side
  of Equation~\ref{eq:G-mod-mor-def} is to be interpreted as a
  $\ZZ$-set graded chain complex.) To define composition, we need to
  give $\ZZ$-set graded chain maps
  \[
  \Mor_\Alg((S, M_\Alg), (T, N_\Alg)) \otimes
  \Mor_\Alg((T, N_\Alg), (U, P_\Alg)) \to
  \Mor_\Alg((S, M_\Alg), (U, P_\Alg)).
  \]
  On the chain level, this map is defined as in
  Section~\ref{sec:categ-Ainf-modules}.  It preserves the grading
  relation~$R$ generated by
  \[
  \bigl((s^* \times_G t) \times_\lambda (t^* \times_G u)\bigr) \mathrel{R}
    (s^* \times_G u)
  \]
  for $s \in S$, $t\in T$, and $u \in U$.
\end{definition}

\begin{example}
  Specializing Definition~\ref{def:cat-vary-gr-set} to the case where
  $\Alg=\Ground=\Field$, we see that $\ModCat_{\Field}$ is the category of
  $\ZZ$-set graded chain complexes. Indeed, the above definition allows us
  to consider the category of set-graded chain complexes as a
  set-graded \dg category (compare Example~\ref{ex:MorTypeD}).
\end{example}

For a fixed transitive $G$-set $S$, the category
$\widetilde\ModCat_{\Alg,S}$ with grading set~$S$
(Definition~\ref{def:cat-fixed-gr-set}) is a subcategory of
$\ModCat_\Alg$, but not a full subcategory: morphisms in
$\widetilde\ModCat_{\Alg,S}$ are only allowed to shift the grading by
a power of~$\lambda$, while in $\ModCat_\Alg$ there is no such
restriction.

We define $\ModCat^\Alg$ similarly, using type~$D$ structures instead
of $\Ainf$-modules.  When $\Alg$ is an $\Ainf$-algebra rather than a
\dg algebra, $\Mod^\Alg$ is a $\ZZ$-set graded $\Ainf$-category, which
is defined analogously to $\ZZ$-set graded \dg categories. The
variants $\lsub{\Alg}\ModCat$ and
$\lsupv{\Alg}\ModCat$ are defined symmetrically, using left actions
instead of right actions. Categories of bimodules are defined
similarly, cf.~Definition~\ref{def:GradingBimodules}.

We can consider $\ZZ$-set graded functors from the set-graded category
$\Mod_{\Alg}$ (from Definition~\ref{def:cat-vary-gr-set}) to the category
of $\ZZ$-set graded chain complexes.  Although the definition of this
notion can be pieced together by what has been written so far, we
spell it out for the reader's convenience:

\begin{example}
  \label{examp:SetGradedFunctor}
Let $\Alg$ be a $(G,\lambda)$-graded $\Ainf$ algebra, let $\Mod_{\Alg}$
(whose objects are pairs $(S,M_{\Alg})$, where $S$ is a $G$-set and
$M_{\Alg}$ is graded by $S$) be its category of $G$-set graded
$\Ainf$-modules, and let $\ModCat_{\Field}$ denote the category of $\ZZ$-set
graded chain complexes. A {\em differential set-graded functor}
(or, less precisely, a {\em \dg functor})
${\mathcal F}\co \Mod_{\Alg}\longrightarrow \ModCat_{\Field}$ consists of the following data:
\begin{itemize}
\item for each $(S,M_{\Alg})\in \Mod_{\Alg}$, a pair 
  ${\mathcal F}(S,M_{\Alg})=(S',M')$, 
  where
  $S'$ is a $\ZZ$-set and $M'$ is a chain complex graded by $S'$, and
\item for each pair $(S,M_{\Alg}), (T,N_{\Alg})\in\Mod_{\Alg}$,
  an element (in the sense of Definition~\ref{def:elements}) 
  ${\mathcal F}_{(S,M),(T,N)}$
  of 
  the internal hom set
  $$\Hom_{\mathcal C}(\Mor_{\Alg}((S,M),(T,N)),
  \Mor_{\mathcal C}((S',M'),(T',N')),$$
  which in turn consists of a pair $(R_{(S,M),(T,N)},\Phi_{(S,M),(T,N)})$ where $R$ is a $\lambda$-invariant relation
  $$R_{(S,M),(T,N)} \subset (T\times_G S^*)^*\times (T'\times_\lambda (S')^*)$$
  and $\Phi_{(S,M),(T,N)}$ is a cycle 
  in  
  $$\Mor_{\mathcal C}(\Mor_{\Alg}((S,M),(T,N)),
  \Mor_{\mathcal C}((S',M'),(T',N'))$$
  which is supported in $R_{(S,M),(T,N)}$.
\end{itemize}
Moreover, we demand that the 
${\mathcal F}_{(S,M),(T,N)}$ respect composition, in the sense that 
the following two conditions are satisfied:
\begin{enumerate}
\item $R_{(S,M),(U,L)}=R_{(S,M),(T,N)}\circ R_{(T,N),(U,L)}$
and 
\item the following diagram commutes:
\[
  \begin{tikzpicture}
    \node at (0,0) (xyTyz) {$\Mor_{\mathcal A}((S,M),(T,N))\otimes
      \Mor_{\mathcal A}((T,N),(U,L))$}; 
    \node at (8,0) (xz) {$\Mor_{\mathcal A}((S,M),(U,L))$}; 
    \node at (0,-2)(FxyTFyz) {$\Mor_{\mathcal C}({\mathcal F}(S,M),{\mathcal F}(T,N))
      \otimes      \Mor_{\mathcal C}({\mathcal F}(T,N),{\mathcal F}(U,L))$};
    \node at (8,-2) (Fxz) {$\Mor_{\mathcal C}({\mathcal F}(S,M),{\mathcal F}(U,L))$}; 
    \draw[->] (xyTyz) to node[above]{\lab{\circ_{\Alg}}} (xz);
    \draw[->] (FxyTFyz) to node[above]{\lab{\circ_{\mathcal C}}} (Fxz);
    \draw[->] (xyTyz) to node[left]{\lab{\Phi_{(S,M),(T,N)}\otimes \Phi_{(T,N),(U,L)}}} (FxyTFyz);
    \draw[->] (xz) to node[right]{\lab{\Phi_{(S,M),(U,L)}}} (Fxz);
  \end{tikzpicture}
\]
\end{enumerate}
\end{example}

The following example (and its generalizations to bimodules) will be
of importance to us:

\begin{example}
  \label{ex:DTPasDGfunctor}
  Let $\Alg$ be a $G$-graded $\Ainf$ algebra, and fix a $G$-set graded
  type $D$ structure $(U,\lsupv{\Alg}P)$. The operation
  $(S,M_{\Alg})\mapsto (S\times_G U,M_{\Alg}\DT \lsupv{\Alg}P)$
  induces a \dg functor~${\mathcal F}$ from the category of $G$-set graded
  $\Ainf$-modules, $\ModCat_{\Alg}$, to the category of $\ZZ$-set
  graded chain complexes. In particular, given two $G$-set graded 
  $\Ainf$ modules $(S,M_\Alg)$ and $(T,N_{\Alg})$, we define a
  homomorphism
  $${\mathcal F}_{(S,M),(T,N)}=(R,\Phi)\in\Hom(\Mor(M,N),\Mor(M\DT P,N\DT
  P))$$ as follows.  $R$ is the tautological relation in 
  $$(T\times_G S^*)\times
  ((T\times_G U)\times_{\lambda} (S\times_G U)^*))
  =
  (T\times_G S^*)\times
  ((T\times_G U)\times_{\lambda} (U^* \times_G S^*))$$
  given by
  $$[t\times_{G} s^*]R[[t\times_{G} u]\times_{\lambda}[u^*\times_{G} s^*]],$$ 
  and $\Phi$ is the map $\cdot\DT \Id_{N}$
  from Section~\ref{sec:tensor-products}.
\end{example}

\begin{definition}
  A \dg functor ${\mathcal F}\co \Mod_{\Alg}\to \ModCat_{\Field}$
  (as in Example~\ref{examp:SetGradedFunctor}) is said to 
  be {\em homology faithful} if the homomorphisms ${\mathcal F}_{(S,M),(T,N)}$
  are homology injective in the sense of
  Definition~\ref{def:Injective} (that is, $\mathcal{F}$ induces an
  injective functor $\HMod(\cM_\Alg) \to \HMod(\ModCat_{\Field})$).
\end{definition}

It is necessary to take morphisms of bimodules over one of the two
actions, and to consider the bimodule structure on the result (see
especially Sections~\ref{sec:mod-hom} and~\ref{sec:Duality}). This
operation interacts with the gradings as follows:
\begin{definition}
  If $S$ is left $(G_1 \times_\lambda G_2^\op)$-set and $T$ is a left
  $(G_3 \times_\lambda G_2^\op)$-set, define
  \begin{align*}
  \Hom_{G_2}(S,T) &\coloneqq \{\,\langle s, t\rangle \mid s \in S, t \in T\,\}
   / \{\,\langle sg, tg \rangle \sim \langle s, t \rangle\mid g \in G_2\,\}\\
   &= T \times_{G_2} S^*.
  \end{align*}
  We view the result as a left $(G_3 \times_\lambda G_1^\op)$-set, as in the
  description as a product over $G_2$.  

  If $V$\! and $W$\! are graded by $S$ and~$T$, then $\Hom_\Ground(V, W)$
  may be graded by $\Hom_{G_2}(S,T)$.  In particular, all of the one-sided
  $\Mor$-spaces between various bimodules are graded by sets like
  this.  For instance, if $\Alg$, $\Blg$, $\Clg$ are $\Ainf$-algebras
  graded by $G_1$, $G_2$, $G_3$, respectively, and $\lsub{\Alg}M_\Blg$
  and $\lsub{\Clg}N_\Blg$ are graded by $S$ and~$T$ as above, then
  $\Mor_\Blg(\lsub{\Alg}M_\Blg, \lsub{\Clg}N_\Blg)$ has underlying
  space $\Hom_\Ground(M \otimes \Tensor^*(\Blg[1]), N)$, which is
  graded by $\Hom_{G_2}(S, T)$.  With this grading on
  $\Mor_\Blg(\lsub{\Alg}M_\Blg, \lsub{\Clg}N_\Blg)$, the differential
  is a graded map.
\end{definition}

%%% Local Variables: 
%%% mode: latex
%%% TeX-master: "Bimodules"
%%% End: 

\counterwithout{equation}{subsection}
\counterwithin{equation}{section}
\section{Pointed matched circles}
\label{sec:PointedMatchedCircles}

\subsection{The algebra of a pointed matched circle}

We start by recalling the following definition~\cite{LOT1}:
\begin{definition}\label{def:PMC}
 A \emph{pointed matched circle} is an oriented circle $Z$, equipped
with a basepoint $z$, and also an additional 
$4k$ points ${\mathbf a}=\{a_1,\dots,a_{4k}\}$ (all distinct from $z$), which are
partitioned into pairs, in such a manner that the one-manifold obtained by
performing surgeries on the $2k$ pairs of points gives a circle. 
\end{definition}

The pairing of points can be thought of as map~$M$ taking an element
of ${\mathbf a}$ to its equivalence class
(consisting of two elements).  Alternately, we may think of the
pairing as an involution~$x \mapsto x'$ on $\CircPts$, taking a point to
the other element of its pair.
We often abbreviate the data  $(Z,\{a_1,\dots,a_{4k}\},M,z)$ by $\PMC$.

\begin{construction}\label{construct:PMC-gives-surf}
  Given a pointed matched circle $\PMC$, we can associate a surface
  whose boundary is a circle, containing a marked point $z$. We
  denote this surface $\PunctF(\PMC)$, and let $F(\PMC)$ denote the
  result of filling in the boundary component of $\PunctF(\PMC)$ with a disk
  $D$.  Note that any two surfaces specified by the same pointed
  matched circle are homeomorphic, via a homeomorphism which is
  uniquely determined up to isotopy.
\end{construction}

\begin{figure}
  \centering
  %Font is 12 point.
  \includegraphics{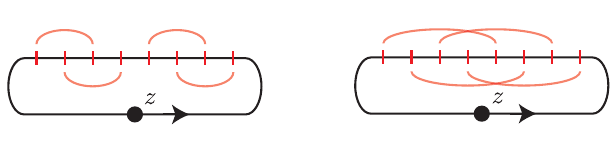}
  \caption{\textbf{Two convenient pointed matched circles for the genus $2$ surface.} Left: the split pointed matched circle. Right: the antipodal pointed matched circle.}
  \label{fig:split-antipodal-pointed}
\end{figure}

Any surface of genus $k>1$ can be represented by more than
one pointed matched circle. There are two convenient families of
pointed matched circles which can be used to describe an arbitrary
oriented surface. One is the {\em split} pointed matched circle, which
is obtained as the $k$-fold connected sum of pointed matched circles
representing genus one surfaces (and dropping extra basepoints). The
other we call the {\em antipodal} pointed matched circle, where $x
\mapsto x'$ is
the map exchanging antipodal points. See
Figure~\ref{fig:split-antipodal-pointed}. One can also reverse the
orientation of
a pointed matched circle $\PMC$ to get a new pointed matched circle
$-\PMC$. 

Recall from~\cite{LOT1} the algebra $\Alg(n,k)$. An $\Field$-basis for
$\Alg(n,k)$ is given by triples $(\SetSS,\SetTT,\phi)$ where $\SetSS$ and $\SetTT$ are
$k$-element subsets of $[n]=\{1,\dots,n\}$ and $\phi\co \SetSS\to \SetTT$ is a
bijection with $\phi(x) \ge x$ for each $x\in\SetSS$. Associated to any
such triple is a number $\inv(\phi)$, the
number of inversions of $\phi$. We define
\[
(\SetSS_1,\SetTT_1,\phi_1)\cdot(\SetSS_2,\SetTT_2,\phi_2)=
\begin{cases}
  (\SetSS_1,\SetTT_2,\phi_2\circ\phi_1) & \text{$\SetTT_1=\SetSS_2$,
    $\inv(\phi_2\circ\phi_1)=\inv(\phi_1)+\inv(\phi_2)$ }\\
  0 & \text{otherwise.}
\end{cases}
\]
Let $\Inv(\phi)$ denote the set of inversions of $\phi$, so
$\inv(\phi)=\#\Inv(\phi)$. For $\sigma=(i,j)\in\Inv(\phi)$ define
$\phi_\sigma$ by $\phi_\sigma(i)=\phi(j)$, $\phi_\sigma(j)=\phi(i)$, and
$\phi_\sigma(k)=\phi(k)$ for $k\neq i,j$. There is a differential on 
$\AlgAS{n}{k}$
given by
\[
\bdy(\SetSS,\SetTT,\phi)=\,\,\,\sum_{\mathclap{\substack{\sigma\in\Inv(\phi)\\ \inv(\phi_\sigma)=\inv(\phi)-1}}}\,\,\,(\SetSS,\SetTT,\phi_\sigma).
\]
We draw basis elements of $\AlgAS{n}{k}$ as strand diagrams with
upward-veering strands, like this:
\[
\mfigb{strands-0}
\]
The product becomes concatenation (with the convention that double
crossings are set to zero) and the differential corresponds to
smoothing crossings.

Now, let $\PMC$ be a pointed matched circle. For each
$i=-k,\dots,k$ we define a subalgebra $\AlgAS{\PMC}{i}\subset
\mathcal{A}(4k,k+i)$ as follows. Cutting $Z$ at the basepoint $z$, the
orientation of $Z$ identifies $\mathbf{a}$ with $[4k]$, so we can view
the matching $M$ as inducing an involution $x \mapsto x'$ of $[4k]$,
which we extend to an involution on the set of subsets of $[4k]$.
For $a=(\SetSS,\SetTT,\phi)$ a basis element of
$\mathcal{A}(4k,k+i)$ and $\SetI\subset\Fix(\phi)$ a set of fixed points of $\phi$
define a map 
$$\phi_\SetI\co (\SetSS\setminus \SetI)\cup (\SetI')\to
(\SetTT\setminus \SetI)\cup (\SetI')$$ by
\[
\phi_\SetI(j)=
\begin{cases}
  \phi(j) & \text{if $j\notin \SetI'$}\\
  j & \text{if $j\in\SetI'$.}
\end{cases}
\]
Define
\[
E(\SetSS,\SetTT,\phi)=\sum_{\SetI\subset \Fix(\phi)}((\SetSS\setminus \SetI)\cup
\SetI',(\SetTT\setminus \SetI)\cup \SetI',\phi_\SetI).
\] 
(Loosely speaking, $E(\SetSS, \SetTT, \phi)$ is obtained by
``smearing'' the locations of the horizontal strands in~$\phi$
according to the matching.)
Then, $\AlgAS{\PMC}{i}$ is the subalgebra of
$\mathcal{A}(4k,k+i)$ generated by 
\[
\{E(\SetSS,\SetTT,\phi)\mid \SetSS\cap \SetSS'=\SetTT\cap \SetTT'=\emptyset\}.
\] 
(In fact, these elements span $\AlgAS{\PMC}{i}$ as an
$\Field$-vector space.)

Note in particular that the primitive idempotents in
$\AlgAS{\PMC}{i}$ correspond naturally to $(k+i)$-element subsets of
$[2k]=[4k]/M$. Given a $(k+i)$-element subset $\SetS$ of $[2k]$ we let $I(\SetS)$
denote the corresponding idempotent. 

The subalgebra $\AlgAS{\PMC}{0}\subset \AlgA{\PMC}$ is the summand
directly relevant to the case of $3$-manifolds with connected
boundary. The other summands become necessary when considering
three-manifolds with disconnected boundary (as in, for instance,
Theorem~\ref{thm:Hochschild}).  Consider the extreme cases
$\AlgAS{\PMC}{-k}$ and $\AlgAS{\PMC}{k}$. The algebra
$\AlgAS{\PMC}{-k}$ is isomorphic to $\Field$; the generator is the
identity permutation of the empty set. By contrast, the algebra
$\AlgAS{\PMC}{k}$ is quite large, but it follows from
Theorem~\ref{thm:halg-support} in Section~\ref{sec:homology-algebra} that
$H_*(\AlgAS{\PMC}{k})\cong\Field$.

We can put the algebras together to form
\[
\AlgA{\PMC}=\bigoplus_{i=-k}^{k}\AlgAS{\PMC}{i},
\]
which we think of as an algebra with an extra integer grading
(represented by the integer $i$), so that elements with different
grading multiply to~$0$. We call this grading the {\em
  strands grading} on $\AlgA{\PMC}$.

Let $\IdemAS{\PMC}{i}$ denote the ring of idempotents in
$\AlgAS{\PMC}{i}$ and $\IdemA{\PMC}$ the ring of
idempotents in $\AlgA{\PMC}$. Of course, we can think of
$\AlgAS{\PMC}{i}=\AlgA{\PMC}\cdot \IdemAS{\PMC}{i} = 
\IdemAS{\PMC}{i}\cdot \AlgA{\PMC}$.
We think of the ring of idempotents as the ground ring. The algebra
$\Alg(\PMC)$ has a natural augmentation $\epsilon\co \Alg(\PMC)
\to \Idem(\PMC)$, sending any strand diagram with a non-horizontal
strand to~$0$.

The algebra $\AlgA{\PMC}$ can also be defined in terms of
Reeb chords in $(Z,\mathbf{a})$; see
\cite[Section~\ref*{LOT:sec:reeb-chords-def}]{LOT1}. In particular,
given a set of Reeb chords $\rhos$ there
is an associated algebra element $a(\rhos)$ (which may be zero if
$\rhos$ does not satisfy some compatibility conditions).

\begin{remark}
  It is not hard to show that if $\PMC$ represents a surface~$F$ of
  genus~$k$ then
  the Grothendieck group of projective $\AlgAS{\PMC}{i}$-modules
  has rank $\binom{2k}{k+i}$. It is interesting to compare with the
  homology of $\Sym^{k+i}(\PunctF)$.
\end{remark}

\subsection{Gradings}
\label{sec:Gradings}

As discussed in \cite[Section~\ref*{LOT:sec:gradings-algebra}]{LOT1},
the algebra $\AlgA{\PMC}$ is
not graded by~$\ZZ$, but rather by a non-commutative group
$\smallGroup$, a $\ZZ$-central extension of $H_1(F(\PMC))$. In fact,
it is more naturally graded by
a larger group $\bigGroup$. We recall that construction first (and discuss
some of its key properties) and then
we turn to the refinement to gradings in 
$\smallGroup$ (which depends on some auxiliary  choices).

Let $Z'=Z\setminus z$.
Fix $p\in {\mathbf a}$ and $\alpha\in H_1(Z',{\mathbf a})$. We define the
{\em multiplicity $m(\alpha,p)$} of $\alpha$ at $p$ 
to be the average of the local multiplicity of $\alpha$ just above $p$ and 
just below $p$. Extend this to a bilinear pairing
$$m\co H_1(Z',{\mathbf a})\times H_0({\mathbf a})\to \OneHalf \ZZ.$$

Let $\bigGroup(\PMC)$ denote the $\ZZ$-central extension of
$H_1(Z',{\mathbf a})$ with the commutation relation
$$\tilde\alpha\cdot\tilde\beta=\tilde\beta\cdot\tilde\alpha\cdot \lambda^{2m(\beta,\bdy\alpha)},$$
where $\alpha,\beta\in H_1(Z',{\mathbf a})$, $\tilde\alpha$ and
$\tilde\beta$ are lifts of $\alpha$ and $\beta$ to $G'$, $\lambda$ is a
generator for the center, and $\bdy \colon H_1(Z',{\mathbf
  a})\to H_0({\mathbf a})$ denotes the connecting homomorphism.
The group $\bigGroup(\PMC)$ can be realized explicitly as an index two
subgroup of the group
$\OneHalf \ZZ\times H_1(Z',{\mathbf a})$ endowed with the
multiplication map
$$(\ell_1,\alpha_1)\cdot (\ell_2,\alpha_2)=(\ell_1+\ell_2+m(\alpha_2,\bdy \alpha_1),\alpha_1+\alpha_2),$$
This index two subgroup is generated by elements $\lambda=(1,0)$ and $(-\OneHalf,[i,i+1])$
where $i=1,\dots,4k-1$. We call $\ell$ (respectively $\alpha$) the
\emph{Maslov component} (respectively \emph{$\SpinC$ component}) of
$(\ell,\alpha)$.
Concretely, $\bigGroup(\PMC)$ consists of pairs $(\ell,\alpha)$ where
\begin{align*}
  \ell &\equiv \OneQuart\#(\text{parity changes in $\alpha$}) \pmod{1}.
\end{align*}

A basis element $a$ of the algebra $\Alg(\PMC)$ has an associated
one-chain in
$H_1(Z,{\mathbf a})$, $[a]$. Recall also that $a$ is a linear
combination of basic elements of $\Alg(4k,k+i)$; let $a_0$ be any of
the terms appearing in this linear combination. Using $a_0$ and $[a]$ we can construct
the $\bigGroup(\PMC)$-grading by the formulas
\begin{equation}
  \label{eq:DefGrading}
  \begin{split}
    \iota(a) &= \inv(a_0)-m([a_0],S)\\
    \gr'(a)&=(\iota(a),[a]),
  \end{split}
\end{equation}
where $S$ is the initial idempotent of $a_0$ (i.e., 
$I(S)\cdot a_0=a_0$).  A short argument shows that the quantity
$\iota(a)$ is independent of the choice of $a_0$; for more details, see
\cite[Proposition~\ref*{LOT:prop:grading-descends}]{LOT1}.
It is also not hard to show that for any $a$, $\gr'(a)$ is an element of
$\bigGroup(\PMC)$, and the map $\gr'$
gives $\AlgA{\PMC}$ a
$\bigGroup(\PMC)$-grading, in the sense that $\gr'(a\cdot
b)=\gr'(a)\cdot \gr'(b)$
and $\gr'(\partial a) = \lambda^{-1}\gr'(a)$;
see~\cite[Proposition~\ref*{LOT:prop:grb-is-grading}]{LOT1}.

As an immediate application of the grading, our algebras satisfy the
condition of Definition~\ref{def:Alg-bounded}.
\begin{lemma}\label{lem:Alg-nilpotent}
  For any pointed matched circle $\PMC$, the algebra $\Alg(\PMC)$
  has nilpotent augmentation ideal.
\end{lemma}
\begin{proof}
  This is immediate from the facts that $\Alg(\PMC)$ is
  finite-dimensional and elements of the augmentation ideal have
  positive $H_1(Z',\mathbf{a})$ gradings.
\end{proof}

In Section~\ref{sec:id-bim}, we will use a stronger grading-positivity
result for our algebras:
\begin{lemma}
  \label{lem:NegativeGradingsOnAlgebra}
  For any pointed matched circle $\PMC$, $\Alg(\PMC)$ is $\grb$-graded in
  non-positive Maslov degrees.  More precisely,
  if $a$ is a generator of $\Alg(\PMC)$, then the Maslov degree of $a$
  is less than or equal to zero,
  with equality if and only if $a$ is an idempotent.
\end{lemma}

\begin{proof}
  We will show that $\iota(a)$ is less than or equal to
  $-1/2$ times the number of moving (non-horizontal) strands in $a$.
    Letting $\Sigma$ denote the set of moving strands in $a$, observe that
  $$\iota(a)=
  -\frac{\#|\Sigma|}{2}+\sum_{\substack{s_i,s_j\in\Sigma\\s_i \ne s_j}}
  \bigl(\#(s_i\cap s_j) - m(I(s_i),[s_j])-m(I(s_j),[s_i])\bigr),$$
  where $I(s_i)$ denotes the initial point of the strand $s_i$,
  and $[s_i]$ is its associated interval. (Horizontal strands do not
  contribute to $\iota(a)$.)  Evidently, the only positive
  contributions here come from crossings of $s_i$ and $s_j$.  For each,
  the contribution is at most $1$ for each pair $\{s_i,s_j\}$.
  However, if $s_i$
  and $s_j$ cross, then either the initial point of $s_i$  is
  contained in the interior of $[s_j]$ or the initial point of $s_j$
  is contained  in the interior of $[s_i]$. Thus, we can cancel
  off each positive contribution with a corresponding~$-1$.
\end{proof}

\subsubsection{Refined gradings}
\label{subsec:SmallGroup}

The \emph{Heisenberg group $\smallGroup(\PMC)$ of $H_1(F)$} is the
central extension of $H_1(F;\ZZ)$ by a subgroup $\ZZ$ generated by
$\lambda$, with commutation relation
$$g\cdot h = h\cdot g \cdot \lambda^{2\#([g]\cap[h])}$$
for any $g,h\in \smallGroup(\PMC)$. (Here, $[g]$ and $[h]$ are the
images of $g$ and $h$ in $H_1(F)$.)

As in \cite[Section~\ref*{LOT:sec:refined-grading}]{LOT1}, there
is a natural inclusion $i_*\co H_1(F)\to H_1(Z',{\mathbf a})$, with
image the kernel of $M_* \circ \bdy$.  Since
$[g]\cap[h]=m(i_*[h],\bdy i_*[g])$, it follows that this inclusion 
$\smallGroup(\PMC)\hookrightarrow \bigGroup(\PMC)$ is 
a group homomorphism.  As explained
in \cite[Section~\ref*{LOT:sec:refined-grading}]{LOT1}, the
$\bigGroup$ grading on the algebra can be refined to a
$\smallGroup$-valued grading. (The refined grading leads to cleaner
statements of the pairing theorems; see, for instance,
Theorem~\ref{thm:GradedPairing} in Section~\ref{sec:PairingTheorems}.) That construction involves certain
choices, as we now elaborate.

\begin{definition}
  \label{def:Compatible}
  Fix a pointed matched circle.
  Given an element $\alpha\in H_1(Z,{\mathbf a})$ and subsets $\SetS$
  and $\SetT$ of $[2k]$ we say that $\alpha$ is {\em compatible
    with the idempotents $I(\SetS)$ and $I(\SetT)$} if
  \[
  M_*(\bdy\alpha) = \SetT - \SetS.
  \]
\end{definition}

In particular, for a generator $a$ of $\Alg(\PMC)$ with
$I(\SetS)\cdot a\cdot I(\SetT)=a$, the homology class $[a]$ is
compatible with $\SetS$ and~$\SetT$.

\begin{definition}\label{def:grading-refinement}
  {\em Grading refinement data} for the algebra $\Alg(\PMC)$ consists of 
  a function
  \[\psi\colon \{\,\SetS\mid\SetS\subset [2k]\,\}\to \bigGroup(\PMC)\]
  satisfying the condition that if
  $g'\in\bigGroup(\PMC)$ is a group element so that
  $[g']$ is compatible with $I(\SetS)$ and $I(\SetT)$,
  then $\psi(\SetS)\cdot g'\cdot \psi(\SetT)^{-1}$ lies in
  $\smallGroup(\PMC)\subset\bigGroup(\PMC)$.
\end{definition}

Grading refinement data for $\Alg(\PMC,i)$ can be specified by
choosing, for each $i=0,\dots,2k$, a base idempotent
$\SetT_i\subset[2k]$ with $|\SetT_i|=i$ and maps
$$\psi_i\colon \{\,\SetS\mid\SetS\subset [2k],\ |\SetS|=i\,\} \to G'(\PMC)$$
satisfying
$$M_*\bdy[\psi_i(\SetS)]=\SetS-\SetT_i. $$ 

\begin{definition}\label{def:refined-grading}
    Given grading refinement data $\psi$
    as above, we can define a corresponding $\smallGroup(\PMC)$-valued
    grading $\gr_\psi$ on $\Alg(\PMC)$
    as follows. For any generator $a$
    of $\Alg(\PMC,i)\subset\Alg(\PMC)$ with idempotents $\SetS$
    and~$\SetT$, define
    $$\gr_{\psi}(a)=\psi(\SetS)\cdot \gr'(a)\cdot \psi(\SetT)^{-1}.$$
    It is straightforward to verify that this is indeed a grading
    with values in $\smallGroup(\PMC)\subset\bigGroup(\PMC)$.  
\end{definition}

For fixed refinement data $\psi$, we can consider the categories of
$\smallGroup(\PMC)$-graded $\Ainf$ modules and type $D$ structures
(and corresponding bimodules).  For the few times we wish to call
attention to all this information, we will decorate the relevant
categories, writing, for example,
$\ModCat_{\Alg,\psi}^{\smallGroup(\PMC)}$ for the category of
$\smallGroup(\PMC)$-graded $\Ainf$-modules over $\Alg$ with its
$\smallGroup(\PMC)$-grading induced by the refinement data $\psi$.

The reader is warned that if $\psi$ and $\psi'$ are two different
refinement data for $\PMC$ then the induced $\smallGroup(\PMC)$-graded
algebras are typically not graded quasi-isomorphic. They do, however,
have isomorphic module categories, according to the following:

\begin{proposition}
  \label{prop:ChangeOfRefinementData}
  If $\psi$ and $\psi'$ are two different refinement
  data for $\Alg(\PMC)$, the identity type \DA\ bimodule
  $\lsupv{\Alg}[\Id]_{\Alg}$ can be given a grading,
  $\lsupv{\Alg,\psi}[\Id]_{\Alg,\psi'}$
  so that the functors 
  \begin{align*}
    (\cdot\DT \lsupv{\Alg,\psi}[\Id]_{\Alg,\psi'})
    \co \ModCat_{\Alg,\psi}^{\smallGroup(\PMC)}& \to 
    \ModCat_{\Alg,\psi'}^{\smallGroup(\PMC)}\\
    (\lsupv{\Alg,\psi}[\Id]_{\Alg,\psi'}
    \DT\cdot)\co \ModCat^{\Alg,\psi'}_{\smallGroup(\PMC)}& \to 
    \ModCat^{\Alg,\psi}_{\smallGroup(\PMC)}
  \end{align*}
  induce equivalences of categories.
\end{proposition}
(See Definition~\ref{def:rank-1-DA-mods} for the definition of
$\lsupv{\Alg}[\Id]_{\Alg}$.)
\begin{proof}
  Consider the left-right
  $(\smallGroup(\PMC),\smallGroup(\PMC))$-set
  $T=\smallGroup(\PMC)$, with action given by
  $(g_1\times g_2)\cdot g=g_1\cdot g \cdot g_2$. 
  We endow $\lsupv{\Alg}[\Id]_{\Alg}$ 
  with a grading by $T$, as follows.
  Recall that generators of $\lsupv{\Alg}[\Id]_{\Alg}$ are of the form
  $I(\SetS)$ where $\SetS$ is any subset of $[2k]$.
  We endow the generator 
  $I(\SetS)$ of $\lsupv{\Alg}[\Id]_{\Alg}$ with the grading
  \begin{equation}
    \label{eq:grading-change}
    \gr^{\psi}_{\psi'}(I(\SetS))=\psi(\SetS)\cdot \psi'(\SetS)^{-1}\in T.
  \end{equation}
  We denote the resulting $T$-graded module by
  $\lsupv{\Alg,\psi}[\Id]_{\Alg,\psi'}$. 
  We must check that,
  for any homogeneous algebra element $a=I(\SetS)\cdot a \cdot
  I(\SetT)$, the gradings of $I(s) \cdot a$ and $a \cdot I(t)$ are
  consistent.  We compute (using
  Definition~\ref{def:refined-grading}):
  \begin{equation}
    \label{eq:refined-grading-check}
  \begin{aligned}
    \gr^{\psi}_{\psi'}(I(\SetS))\cdot\gr_{\psi'}(a)&=
      \bigl(\psi(\SetS)\psi'(\SetS)^{-1}\bigr)
      \bigl(\psi'(\SetS)\gr'(a)\psi'(\SetT)^{-1}\bigr)\\
      &= \psi(\SetS)\gr'(a)\psi'(\SetT)^{-1}\\
      &=\gr_{\psi}(a)\cdot \gr^{\psi}_{\psi'}(I(\SetT)),
  \end{aligned}
  \end{equation}
  showing that $\lsupv{\Alg,\psi}[\Id]_{\Alg,\psi'}$ is indeed a $T$-graded
  $\Alg(\PMC,\psi)$-$\Alg(\PMC,\psi')$ bimodule. An inverse to this
  bimodule is supplied by the analogously-defined bimodule
  $\lsupv{\Alg,\psi'}[\Id]_{\Alg,\psi}$, where we use
  the canonical identification $T\times_{\smallGroup(\PMC)}T\to T$
  induced by multiplication to grade the isomorphism from the tensor product
  to the identity.
\end{proof}

\begin{definition}
  \label{def:Refinable}
  Let $S$ be a right $\bigGroup$-space and $M_{\Alg}$ an
  $\Ainf$-module over $\Alg$ graded by~$S$. Consider choices of
  $\x,\y\in M_{\Alg}$
  and $\SetS,\SetT\subset [2k]$ with the following properties:
  \begin{itemize}
  \item $\x\cdot I(\SetS)=\x$,
  \item $\y\cdot I(\SetT)=\y$, and
  \item $\x$ and $\y$ are non-zero elements which are 
    homogeneous with respect to the $S$-grading.
  \end{itemize}
  We say $(S,M_{\Alg})$ is {\em $\smallGroup$-refinable} if for
  all choices of $\x$, $\y$, $\SetS$, and $\SetT$ as above,
  any group element $g\in\bigGroup$ such that
  $\gr'(\x)\cdot g=\gr'(\y)$ is compatible with the idempotents
  $I(\SetS)$ and $I(\SetT)$.
\end{definition}

\begin{remark}
  \label{rmk:Homogeneity}
  Let $(S,M_{\Alg})$ be $\smallGroup$-refinable. If $\x\in
  M_{\Alg}$ is a homogeneous element, then there is an elementary
  idempotent $I(\SetS)$ with the property that $\x\cdot
  I(\SetS)=\x$. To see this, note that if $\x$ is homogeneous of some
  degree (in $S$), then $\x\cdot I(\SetS)$ and $\x\cdot I(\SetT)$ are
  homogeneous of the same degree; thus, if both elements are non-zero,
  then the identity element
  $e\in\bigGroup$ is compatible with the idempotents $I(\SetS)$ and
  $I(\SetT)$, which in turn forces $\SetS=\SetT$.
\end{remark}

\begin{lemma}
  \label{lem:RefineGrading}
  For any grading refinement data $\psi$, there is a natural
  functor
  $${\mathcal F}^{\psi}\co \ModCat^{\smallGroup(\PMC)}_{\Alg,\psi}
  \to \ModCat^{\bigGroup(\PMC)}_{\Alg}$$
  which is homology faithful
  (i.e., the induced functor on homology is injective on morphisms).  
  These functors are compatible with changing refinement data,
  in the sense that there is a natural isomorphism of functors
  $${\mathcal F}^{\psi}\cong{\mathcal F}^{\psi'}\circ (\cdot\DT \lsupv{\Alg,\psi}\Id_{\Alg,\psi'}).$$
  Let
  $M_{\Alg}\in\ModCat^{\bigGroup(\PMC)}_{\Alg}$ be an object with
  grading set $S$ that is essential in the sense of
  Definition~\ref{def:EssentialGradingSet}. Then, $M_\Alg$
  is isomorphic to ${\mathcal F}^{\psi}(N_{\Alg,\psi})$ for
  some $N_{\Alg,\psi}\in\ModCat^{\smallGroup(\PMC)}_{\Alg,\psi}$ if and only 
  if $M_{\Alg}$ is $\smallGroup$-refinable.
\end{lemma}

\begin{proof}
  First we define the functor $\mathcal{F}^\psi$.
  Recall that objects in $\ModCat^{\smallGroup(\PMC)}_{\Alg,\psi}$ consist of pairs
  $(S,M_{\Alg})$, where 
  \begin{itemize}
  \item $S$ is a right $\smallGroup(\PMC)$-set and
  \item $M_{\Alg}$ is a right $\Alg$ $\Ainf$-infinity module graded by~$S$
  \end{itemize}
  (satisfying the necessary compatibility conditions). 
  The functor is defined on objects by
  $${\mathcal F}^{\psi}(S,M_{\Alg})=
  (S\times_{\smallGroup(\PMC)}\bigGroup(\PMC),N_{\Alg}),$$
  where on the right hand side, the module $M_{\Alg}$ is given the grading
  $\gr(\x)=\grb(\x)\cdot\psi^{-1}(\SetS)$,
  where $\x\cdot I(\SetS)=\x$.

  We now must define the functor on morphisms.  (See
  Example~\ref{examp:SetGradedFunctor}.)
  Given two objects $(M_{\Alg},S)$ and $(N_{\Alg},T)$ in
  $\ModCat^{\smallGroup(\PMC)}_{\Alg,\psi}$, we must define a
  homomorphism of set-graded complexes
  $$\Mor_{\Alg,\psi}^{\smallGroup(\PMC)}((S,M_{\Alg}),(T,N_{\Alg}))
  \rightarrow \Mor_{\Alg}^{\bigGroup(\PMC)}
  ((S\times_{\smallGroup(\PMC)}\bigGroup(\PMC),M_{\Alg}),(T\times_{\smallGroup(\PMC)}\bigGroup(\PMC),N_{\Alg})).$$
  This consists of two pieces of data:
  a chain map $\Phi$ of morphism spaces, which in this case we take
  to be the identity map (as a graded linear map), and a
  $\lambda$-invariant relation (on
  which it is supported)
  \begin{align*}
    R\subset &(T\times_{\smallGroup(\PMC)} S^*)\times_{\lambda}
    \left((T\times_{\smallGroup(\PMC)}\bigGroup(\PMC))\times_{\bigGroup(\PMC)}
      (S\times_{\smallGroup(\PMC)}\bigGroup(\PMC))^*\right)
  \end{align*}
  which in this case we take to be the relation
  $$(t\times_{\smallGroup} s^*) R 
  (t\times_{\smallGroup} g')\times_{\bigGroup} (g'\times_{\smallGroup} s^*),$$
  for arbitrary $s\in S$, $t\in T$, $g'\in\bigGroup$.

  The identity map $\Phi$ is clearly injective, both as a chain map
  and on homology.
  To check that $R$ is injective (as in
  Definition~\ref{def:Injective}), note that
  \[
  (t \times_\smallGroup s^*) R (t \times_\smallGroup 1)
  \times_{\bigGroup} (1 \times_\smallGroup s^*)
  \]
  and that
  $$(t_1\times_\smallGroup 1)\times_{\bigGroup} (1\times_\smallGroup s_1^*)
  =(t_2\times_\smallGroup g')\times_{\bigGroup} (g'\times_{\smallGroup} s_2^*)$$
  if and only if there is an element $k\in \bigGroup$ so that
  \begin{align*}
    (t_1\times_\smallGroup 1)\times(1 \times_{\smallGroup} s_1^*)
    =(t_2\times_\smallGroup g' k)\times k^{-1}\cdot
    (g'\times_\smallGroup s_2^*)
    =(t_2\times_\smallGroup g' k)\times (g'k\times_\smallGroup s_2^*);
  \end{align*}
  i.e., there are $g,h \in\smallGroup$
  so that the following equations hold:
  \begin{align*}
    t_1\times 1&=t_2\cdot g\times g^{-1} g' k \\
    s_1\times 1&= s_2\cdot h\times  h^{-1} g'k.
  \end{align*}
  We conclude that $g=h$, and hence $t_1=t_2\cdot g$ and $s_1=s_2\cdot g$; i.e.,
  $(s_1\times_{\smallGroup} t_1^*)=(s_2\times_{\smallGroup} t_2^*)$, as needed.

  One can think of the functor  ${\mathcal F}^{\psi}$ as given by
  $(\cdot\DT \lsupv{\Alg,\psi}\Id_\Alg)$, where we think of
  $\lsupv{\Alg,\psi}\Id_\Alg$ as the identity bimodule with its
  $\smallGroup(\PMC)$-$\bigGroup(\PMC)$ grading by the set
  $\bigGroup(\PMC)$, endowed with the left $\smallGroup(\PMC)$ action
  by left translation, and right $\bigGroup(\PMC)$-grading by
  $\bigGroup(\PMC)$. The stated isomorphism of functors (gotten by varying $\psi$)
  follows from the fact that
  $$\lsupv{\Alg_\psi'}\Id_{\Alg,\psi}\DT \lsupv{\Alg,\psi}\Id_{\Alg}\cong 
  \lsupv{\Alg,\psi'}\Id_{\Alg}.$$

  Finally, we prove the claim about the image of
  $\mathcal{F}^\psi$. Given a $\smallGroup(\PMC)$-refinable module
  $(M_{\Alg},S)$, we proceed as follows.
  For each non-zero, $\gr'$-homogeneous element $\x$ we define
  \begin{equation}
    \label{eq:RefineTheGrading}
    \gr(\x)=\gr'(\x)\cdot \psi(\SetS)^{-1},
  \end{equation}
  where $I(\SetS)$ is the idempotent with
  $\x\cdot I(\SetS)=\x$ (see Remark~\ref{rmk:Homogeneity}).
  Let $\gr(M_{\Alg})$ denote the
  image in $S$ of all homogeneous elements. Let $T\subset S$
  denote the $\smallGroup(\PMC)$-orbits of $\gr(M_{\Alg})$. 
  The grading $\gr$ can be viewed as giving a $T$-grading
  on $M_{\Alg}$, for the $\smallGroup(\PMC)$-set $T$. Moreover,
  applying ${\mathcal F}^{\psi}$ to this object, we obtain an object
  which is isomorphic to $M_{\Alg}$ with its original $S$-grading, via
  the $\bigGroup(\PMC)$-space isomorphism
  $T\times_{\smallGroup(\PMC)}\bigGroup(\PMC)\cong S$ (given by
  $(t,g')\mapsto t\cdot g'$). (Note that we are using here the fact that
  the grading set $S$ is essential; otherwise the canonical map
  $T\times_{\smallGroup(\PMC)}\bigGroup(\PMC)\to  S$ would
  fail to surject.)
\end{proof}

\begin{remark}
  \label{rmk:RefineOtherModules}
  We have phrased Lemma~\ref{lem:RefineGrading} in terms of right $\Ainf$-modules
  over $\Alg$; but the same arguments work for type $D$ structures,
  and indeed bimodules of various types. Note, however, when refining
  left, rather than right, modules,
  Equation~\eqref{eq:RefineTheGrading} gets replaced by 
  $\gr(\x)=\psi(\SetS)\cdot\gr'(\x)$.
  Similarly, if $M$ is a left-right $\Alg$-bimodule,
  we define
  $\gr(\x)=\psi(\SetS)\cdot\gr'(\x)\cdot \psi(\SetT)^{-1}$
  when $I(\SetS)\cdot\x\cdot I(\SetT)=\x$.
\end{remark}

When working with type $D$ modules, we also need to relate the grading
on $\Alg(\PMC)$ and $\Alg(-\PMC)$.
Recall from \cite[Equation~(\ref*{LOT:eq:def-R})]{LOT1} that, if $r: Z
\to -Z$ is the (orientation-reversing) identity map, then
\begin{equation}\label{eq:def-R}
R(j,\alpha) = (j, r_*(\alpha))
\end{equation}
defines an group anti-homomorphism from $G(\PMC)$ to $G(-\PMC)$.  For
$\SetS \subset [2k]$ corresponding to an idempotent in $\Alg(\PMC)$,
let $\SetS' = [2k] \setminus \SetS$ correspond to the complementary
idempotent in $\Alg(-\PMC)$.

\begin{definition}\label{def:reverse-grading-refine}
  Given grading refinement data~$\psi$ for $\Alg(\PMC)$, the
  \emph{reverse} of $\psi$ is defined on
  the idempotents of $\Alg(-\PMC)$ by
  \[
  \psi'(\SetS') = R(\psi(\SetS))^{-1}.
  \]
\end{definition}

\begin{lemma}\label{lem:reverse-grading-refine}
  The reverse grading refinement data $\psi'$ of
  Definition~\ref{def:reverse-grading-refine} is grading refinement
  data for $\Alg(-\PMC)$ (i.e., satisfies Definition~\ref{def:grading-refinement}).
\end{lemma}
(Compare \cite[Equation~(\ref*{LOT:eq:R-small-g})]{LOT1}.)
\begin{proof}
  For $g' \in G'(-\PMC)$ compatible with $I(\SetS')$ and $I(\SetT')$,
  we have
  \begin{align*}
    \psi'(\SetS') \cdot g' \cdot \psi'(\SetT')^{-1}
      &= R\bigl(\psi(\SetS)\bigr)^{-1} \cdot g'
           \cdot R\bigl(\psi(\SetT)\bigr)\\
      &= R\Bigl(\psi(\SetT) \cdot R(g') \cdot
           \psi(\SetS)^{-1}\Bigr).
  \end{align*}
  Now observe that $R(g')$ is compatible with $I(\SetT)$
  and $I(\SetS)$, so $\psi(\SetT) \cdot R(g') \cdot
  \psi(\SetS)^{-1}$ is in $G(\PMC)$ as $\psi$ is grading refinement
  data for $\Alg(\PMC)$.  Thus the last line in
  the displayed equation is in $G(-\PMC)$, as desired.
\end{proof}

\subsection{An example: the algebra of the surface of genus one}
\label{subsec:GenusOneAlgebra}

We recall the algebra for a genus one surface, as described in
\cite[Section~\ref*{LOT:sec:torus-algebra}]{LOT1}.

Consider the surface $F$ of genus $1$.  This can be represented
by a unique pointed matched circle $\PMC$. The corresponding
algebra $\Alg(\PMC)$ has three non-trivial summands
\[\Alg(\PMC)=\Alg(\PMC,-1)\oplus \Alg(\PMC,0)\oplus\Alg(\PMC,1).\]
The two outermost summands are not very interesting.  The summand
$\Alg(\PMC,-1)$ is isomorphic to $\Field$ (there are zero moving strands). Although
$\Alg(\PMC,1)$ is not one-dimensional, its homology is, and indeed
$\Alg(\PMC,1)$ is quasi-isomorphic to $\Field$.

Thus, in the genus one case, the interesting algebra is the summand
$\Alg=\Alg(\PMC,0)$.

That algebra is generated (over $\Field$) by two idempotents $\iota_0$
and $\iota_1$, and six non-trivial elements
$\rho_1$, $\rho_2$, $\rho_3$, $\rho_{12}$, $\rho_{23}$, and $\rho_{123}$.

The differential is zero, and the non-zero products are
\[
  \rho_1\rho_2 = \rho_{12} \qquad \rho_2\rho_3 = \rho_{23} \qquad
  \rho_1\rho_{23} = \rho_{123} \qquad \rho_{12}\rho_{3} = \rho_{123}.
\]
(All other products of two $\rho$'s vanish.)  There
are also compatibility conditions with the idempotents:
\begin{align*}
  \rho_1&=\iota_0\rho_1 \iota_1&
  \rho_2&=\iota_1\rho_2 \iota_0&
  \rho_3&=\iota_0\rho_3 \iota_1 \\
  \rho_{12}&=\iota_0\rho_{12} \iota_0&
  \rho_{23}&=\iota_1\rho_{23} \iota_1&
  \rho_{123}&=\iota_0\rho_{123} \iota_1.
\end{align*}

The unrefined grading takes values in the group $\bigGroup$ generated by
quadruples $(m;a,b,c)$ where $j\in\OneHalf\ZZ$, $a,b,c\in\ZZ$ and $j$
is an integer if $a,b,c\in2\ZZ$ or if $b\in2\ZZ$ and $a,c\not\in2\ZZ$,
and a half-integer otherwise.
The multiplication on $\bigGroup$ is given by
\begin{multline*}
(m_1; a_1,b_1,c_1) \cdot (m_2; a_2,b_2,c_2)
= \biggl(m_1+m_2+
\frac{1}{2}
\begin{vmatrix}
a_1 & b_1 \\
a_2 & b_2 
\end{vmatrix}
+ 
\frac{1}{2}
\begin{vmatrix}
b_1 & c_1 \\
b_2 & c_2 
\end{vmatrix};\\ a_1+a_2, b_1+b_2, c_1+c_2\biggr).
\end{multline*}
Here, $\lambda$ is the element $(1;0,0,0)$.
Gradings on the algebra are specified by:
\[
  \grb(\rho_1)=\left(-\textstyle\frac{1}{2};1,0,0\right) \qquad
  \grb(\rho_2)=\left(-\textstyle\frac{1}{2};0,1,0\right) \qquad
  \grb(\rho_3)=\left(-\textstyle\frac{1}{2};0,0,1\right).
\]

The group $\smallGroup\subset \bigGroup$ is generated by elements
$(1;0,0,0)$, $(\OneHalf;0,1,1)$ and $(\OneHalf;1,1,0)$.

One choice of grading refinement data $\{\psi\}$ is given by the function
$\psi\colon \{\iota_0,\iota_1\}\to \bigGroup$
\begin{equation}
  \label{eq:RefineTorus}
\psi(\iota_0)=(0;0,0,0) \qquad \psi(\iota_1)=(-\OneHalf;1,0,0).
\end{equation}
With respect to this refinement, then, the induced $\smallGroup$-grading
on the algebra is specified by
\[
  \gr(\rho_1)=\left(0;0,0,0\right) \qquad
  \gr(\rho_2)=\left(-\textstyle\frac{1}{2};1,1,0\right) \qquad
  \gr(\rho_3)=\left(0;-1,0,1\right).
\]

\subsection{Induction and restriction functors}\label{sec:induction-restriction}
\begin{figure}
  %Font is 12 point.
  \centering
  \includegraphics{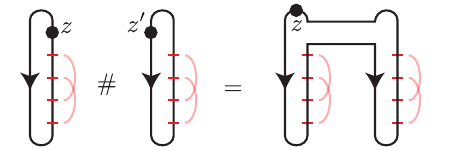}
  \caption{\textbf{The connect sum of pointed matched circles.}}
  \label{fig:ptd-connect-sum}
\end{figure}

Consider two pointed matched circles $\PMC=(Z,\mathbf{a},M,z)$ and
$\PMC'=(Z',\mathbf{a'},M',z')$. Taking the connect sum of $Z$ and $Z'$
at the points $z$ and
$z'$ we obtain a new matched circle
$(Z\connectsum Z',\mathbf{a}\cup\mathbf{a'},M\cup M')$. There are two natural choices of
where to put the basepoint for $Z\connectsum Z'$; for definiteness, we will put
the basepoint in $Z\cap(Z\connectsum Z')$, slightly counter-clockwise of where
the original $z$ was. Thus, we obtain a new pointed matched circle
$\PMC\connectsum \PMC'$. See
Figure~\ref{fig:ptd-connect-sum}. If $\PMC$ specifies a surface $F$
and $\PMC'$ specifies $F'$ then $\PMC\connectsum\PMC'$ specifies
$F\connectsum F'$.

There is an obvious inclusion map
\[
i_{\PMC,\PMC'}\co \AlgA{\PMC}\otimes\AlgA{\PMC'}\to\AlgA{\PMC\connectsum\PMC'}.
\]
This is a morphism of differential algebras. If we fix grading
refinement data for $\psi_1$ for $\Alg(\PMC)$, $\psi_2$ for $\Alg(\PMC')$, and
$\psi_{12}$ for $\AlgA{\PMC\connectsum\PMC'}$ independently, then
the bimodule $\lsupv{\Alg(\PMC \connectsum
  \PMC')}[i_{\PMC,\PMC'}]_{\Alg(\PMC)\otimes \Alg(\PMC')}$ is graded by
$\smallGroup(\PMC)\times_{\ZZ} \smallGroup(\PMC') \cong
\smallGroup(\PMC\connectsum \PMC')$ as follows.
For $\iota_{\SetS_1 \times
\SetS_2}$ an idempotent in $\Alg(\PMC) \otimes \Alg(\PMC')$, set the
grading of the generator $I(\SetS_1 \times \SetS_2)$ to be
\begin{equation}
  \label{eq:grading-inclusion}
\gr(I(\SetS_1 \times \SetS_2)) = \psi_{12}(\SetS_1 \times\SetS_2)
  \psi_1(\SetS_1)^{-1} \psi_2(\SetS_2)^{-1}.
\end{equation}
(Compare Formula~\eqref{eq:grading-change}.)
We must check that, for $a \in \Alg(\PMC) \otimes \Alg(\PMC')$ from
idempotent $\SetS_1 \times \SetS_2$ to $\SetT_1 \times \SetT_2$, that
$\iota_{\SetS_1 \times \SetS_2}\cdot a$ and $a\cdot\iota_{\SetT_1 \times
  \SetT_2}$ have consistent grading, i.e.,
\begin{equation*}
  \gr(I(\SetS_1 \times \SetS_2)) \cdot \gr_{\psi_1\times\psi_2}(a) =
     \gr_{\psi_{12}}(a) \cdot \gr(I(\SetT_1 \times\SetT_2)).
\end{equation*}
This straightforward computation is
analogous to Formula~\eqref{eq:refined-grading-check}.
Thus, these inclusion maps induce 
(graded) restriction functors of module categories
\begin{align}\label{eq:graded-rest}
  \Rest_{\PMC,\PMC'}\co \ModCat_{\AlgA{\PMC\connectsum \PMC'}}&\to
  \ModCat_{\AlgA{\PMC},\AlgA{\PMC'}}
\end{align}
see Section~\ref{sec:induct-restr}.

There are also projection maps
\[
p_{\PMC,\PMC'}\co \AlgA{\PMC\connectsum \PMC'}\to\AlgA{\PMC}\otimes\AlgA{\PMC'}
\]
obtained by setting to zero any basis element which crosses between
$Z$ and $Z'$. (This can also be thought of as adding a second
basepoint to $\PMC\connectsum \PMC'$.)  For arbitrary grading
refinement data for $\Alg(\PMC)$, $\Alg(\PMC')$, and
$\AlgA{\PMC\connectsum\PMC'}$, the bimodule
$\lsupv{\Alg(\PMC)\otimes\Alg(\PMC')}[p_{\PMC,\PMC'}]_{\Alg(\PMC
  \connectsum \PMC')}$ can be graded by Formula~\eqref{eq:grading-inclusion} as
before.  Thus, we get an induction
functor
\begin{align*}
\Induct^{\PMC,\PMC'}\co \lsupv{\AlgA{\PMC\connectsum \PMC'}}\ModCat&
\to{}\lsupv{\AlgA{\PMC},\AlgA{\PMC'}}\ModCat;
\end{align*}
again, see Section~\ref{sec:induct-restr}.

These restriction and induction functors will be used in defining the
type \AAm\ and \DD\ modules, respectively, in
Section~\ref{sec:CFBimodules}.

\subsection{Mapping class groupoid}
Fix an integer $k$.  Let $\PMC$ be a pointed matched circle
on $4k$ points, and let $\PunctF(\PMC)$ be the associated
surface-with-boundary. 

Given pointed matched circles $\PMC_1$ and $\PMC_2$ let
\[
\MCG_0(\PMC_1,\PMC_2)=\{\phi\co\PunctF(\PMC_1)\stackrel{\cong}{\to}
\PunctF(\PMC_2)\mid \phi(z_1)=z_2\}/\text{isotopy}
\] 
denote the set of basepoint-respecting isotopy class of
homeomorphisms $\phi\co\PunctF(\PMC_1)\to
\PunctF(\PMC_2)$ carrying $z_1$ to $z_2$, where $z_i\in
\partial \PunctF(\PMC_i)$ is the basepoint. (The subscript
$0$ is to indicate that the maps respect the boundary and the basepoint.)

The {\em{genus $k$ bordered mapping class groupoid}} $\MCG_0(k)$ is
the category whose objects are pointed matched circles with $4k$
points and with morphism set between $\PMC_1$ and $\PMC_2$ given by
$\MCG_0(\PMC_1,\PMC_2)$. This is clearly a groupoid: each morphism is
invertible, and any two objects can be connected by a
morphism. Moreover, the endomorphisms of a given pointed matched
circle is identified with the strongly based mapping class group of
the surface of genus~$k$.

\begin{remark}
  Our mapping class groupoid is closely related to the \emph{Ptolemy
    groupoid} studied by Penner~\cite{Penner87:MCGroupoid}. In particular,
  our pointed matched circles are called \emph{chord
    diagrams} in the fat graph
  literature~\cite{AndersenBenePenner09:MCGroupoid},~\cite{Bene08:ChordDiagrams}.
\end{remark}

%%% Local Variables: 
%%% mode: latex
%%% TeX-master: "Bimodules"
%%% End: 

\section{Homology of the algebra}
\label{sec:homology-algebra}

In this section we compute the homology $H(\Alg(\PtdMatchCirc))$ of
the algebras $\Alg(\PtdMatchCirc)$. We will denote $H(\Alg(\PMC))$ by
$\HAlg(\PtdMatchCirc)$.  Like $\Alg(\PtdMatchCirc)$,
$\HAlg(\PtdMatchCirc)$ is graded by $\bigGroup(\PMC)$.
Specifically, the degree $(\ell,h)$ part of $\HAlg(\PtdMatchCirc)$ is
represented by cycles whose $\SpinC$ part is $h$, and whose
number of crossings is determined, up to an additive constant
depending on the one-chain and initial idempotent, by $\ell$
(Equation~\eqref{eq:DefGrading}).

This calculation will be used later in proving functorial properties
of bordered Heegaard Floer homology (especially in showing that the
identity diffeomorphism induces the identity bimodule,
Theorem~\ref{thm:Id-is-Id}). For the purpose of the following
statement, recall that if $p\in [4k]$, we let $M(p)$
denote its corresponding equivalence class with respect to the matching.

\begin{figure}
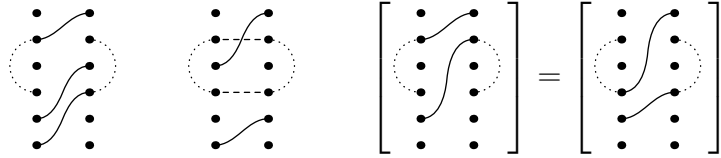

  \[
  \mfigb{strands-70}\qquad
  \mfigb{strands-71}\qquad
  \left[\mfigb{strands-72}\right] = \left[\mfigb{strands-73}\right]
  \]
  \caption{\textbf{Illustration of the statement of
      Theorem~\ref{thm:halg-support}.}
    From left to right, we have an algebra element whose homology
    class is ruled out by condition~\ref{item:halg-multiplicity}), an
    algebra element whose homology class is ruled out by
    condition~\ref{item:halg-half-match}, and two algebra elements
    representing the same homology class.  In each case only a portion of the
    matched circle is drawn.}
  \label{fig:halg-support}
\end{figure}

\begin{theorem}\label{thm:halg-support}
  Let $\SetS, \SetT$ be a pair of subsets of $[2k]$.  The degree
  $(\ell,h)$ part of $I(\SetS)\cdot \HAlg(\PtdMatchCirc)\cdot
  I(\SetT)$ is trivial unless all of the following conditions hold:
  \begin{enumerate}
  \item\label{item:halg-compatible} The homology class $h$ is
    compatible with $I(\SetS)$ and $I(\SetT)$, in the sense
    of Definition~\ref{def:Compatible}.
  \item\label{item:halg-multiplicity} The local multiplicities
    of $h\in H_1(Z,\CircPts)$ are $0$ or $1$.
  \item\label{item:halg-half-match} If $p_1,p_2\in\CircPts$ are
    matched points (so $M(p_1) = M(p_2)$), $p_1 \in
    \Int(\supp(h))$, and $p_2 \not\in \Int(\supp(h))$, then
    $M(p_1)\not\in \SetS\cap\SetT$.
  \item\label{item:halg-minimal-deg} The Maslov degree $\ell$ is minimal
    among all algebra elements
    with associated one-chain~$h$.
  \end{enumerate}
  Furthermore, for every $(\ell,h)$, $I(\SetS)$, and $I(\SetT)$ satisfying the
  conditions above, the degree $(\ell,h)$ part of $I(\SetS)\cdot \HAlg\cdot  I(\SetT)$ is
  1-dimensional, and is represented by any crossingless diagram of the
  correct grading.
\end{theorem}

The above result can be used to calculate the ranks of the homology.
For example, by counting the $\SpinC$ parts allowed in the above theorem,
one finds that for the split pointed
matched circle for the surface of genus two 
$\PtdMatchCirc_{\text{spl}}$,
$\HAlg(\PtdMatchCirc_{\text{spl}})$ has total rank $164$, divided up
according to the number of strands as follows:
\[
\sum_{i}\dim(\HAlg(\PtdMatchCirc_{\text{spl}},i)) \cdot T^i=
T^{-2}+32\cdot T^{-1} + 98 + 32\cdot T^1+T^2. \]
On the other hand, if
$\PtdMatchCirc$ is the antipodal pointed matched circle for the
surface of genus two then $\HAlg(\PtdMatchCirc,0)$ has rank $70$, though the other ranks
are the same; i.e.,
\[
\sum_{i}\dim(\HAlg(\PtdMatchCirc,i))\cdot T^i=
T^{-2}+32 \cdot T^{-1} + 70 + 32 \cdot T+T^2 \]
(so the total rank is $136$).
The different ranks here
underscore the fact that the algebra is really associated to a pointed matched
circle, rather than its underlying surface. We shall, however, see
later that the derived categories of modules for different representatives
of the same surface are equivalent,
cf.\ Theorem~\ref{thm:AlgebraDependsOnSurface}.

\subsection{Computing the homology}
\label{sec:computing-homology}

Before proving Theorem~\ref{thm:halg-support}, we first study the
(rather simple) homology of $\Alg(n,k)$.

\begin{lemma}\label{lem:noneg-mult-exist}
  For $\SetSS$ and $\SetTT$ two subsets of $\{1,\dots,n\}$, if the unique
  homology class $h\in H_1(Z,\CircPts)$ satisfying $\bdy h = \SetTT - \SetSS$ has
  all non-negative local multiplicities then $I(\SetSS)\cdot \Alg(n,k)\cdot I(\SetTT)$ is
  non-empty.  If some local multiplicity is at least two, then
  $I(\SetSS)\cdot \Alg(n,k)\cdot I(\SetTT)$ has dimension at least~$2$ (i.e., there is some
  element with a crossing in it).
\end{lemma}

\begin{proof}
  If $h=0$, then $\SetSS=\SetTT$ and we can take the element $I(\SetSS)$.  If $h$ has
  non-negative local multiplicities, pick some maximal interval $[i,j]
  \subset Z$ with positive multiplicity, connect $i$ to~$j$, and
  subtract~$1$ from all multiplicities in $[i,j]$.  The result still
  has non-negative multiplicity and by induction we can find an
  element of $\Alg(n,k-1)$ with these multiplicities.  When we add the
  first strand we get our desired element.  If $h$ has multiplicity at
  least~$2$ anywhere, we will consider nested intervals in this
  procedure and therefore construct an element of the strands algebra
  with a crossing.
\end{proof}

\begin{proposition}\label{prop:a-n-k-homology}
  For $\SetSS$ and $\SetTT$ two subsets of $\{1,\dots,n\}$, if
  $I(\SetSS)\cdot \Alg(n,k)\cdot I(\SetTT)$ is 1-dimensional, so is its homology; otherwise
  the homology is~0.
\end{proposition}

\begin{proof}
  Each summand $I(\SetSS)\cdot\Alg(n,k)\cdot I(\SetTT)$ has a well-defined
  induced one-chain $h$ with $\partial h = \SetTT-\SetSS$.
  If $I(\SetSS)\cdot\Alg(n,k)\cdot I(\SetTT)$ is more than one dimensional,
  it is easy to see that at least one of following two conditions is satisfied:
  \begin{enumerate}
  \item 
    \label{case:HorizontalCrossing}
    some element of $\SetSS\cap\SetTT$ is in the interior of the
    support of $h$, or
  \item 
    \label{case:MultBiggerThanTwo}
    $h$ has local multiplicity at least $2$ at some point.
  \end{enumerate}
  Thus, our goal is to show that if $h$ satisfies either of the above
  two conditions then the corresponding summand of $\Alg(n,k)$
  is acyclic.

  We start with Case~\eqref{case:HorizontalCrossing}: suppose that
  $i\in\SetSS\cap\SetTT$ is contained in a part of the support of $h$
  where the local multiplicity (both just above and just below $i$)
  is~$1$.  (The other possibilities in Case~\eqref{case:HorizontalCrossing}
  will be handled in Case~\eqref{case:MultBiggerThanTwo} below.) Then,
  we define a map
  $$H\colon I(\SetSS)\cdot \Alg(n,k)\cdot I(\SetTT)\to 
  I(\SetSS)\cdot \Alg(n,k)\cdot I(\SetTT)$$
  as follows. Recall that the part of the algebra
  $I(\SetSS)\cdot \Alg(n,k)\cdot I(\SetTT)$ is
  spanned by upward-veering maps $\phi\colon \SetSS\to \SetTT$.
  Given such a map $\phi$, define $H(\phi)$ by
  \[
  H(\phi) =
  \begin{cases}
    0 & \phi(i) = i\\
    \phi_{(i,j)} & \phi(i) \ne i, \phi(j) = i.
  \end{cases}
  \]
  (Recall that for $\sigma$ a transposition, $\phi_\sigma$ was defined
  in Section~\ref{sec:PointedMatchedCircles} and switches the roles of
  the inputs $i$ and~$j$.)

  We claim that this map $H$ is a null-homotopy of the identity on
  this part of the algebra:
  \begin{equation}
    \label{eq:NullHomotopy}
    \partial\circ H + H\circ \partial = \Id.
  \end{equation}
  It suffices to verify
  this for generators $\phi\co \SetSS\to\SetTT$, with the
  same two cases as in the definition of~$H$:
  \begin{itemize}
  \item Suppose $\phi(i)=i$. Then, $H(\phi)=0$. Moreover, there is
    exactly one term in $\partial \phi$ for which $H$ does not vanish:
    it is the term corresponding to the resolution of the horizontal
    strand at $i$ with a unique strand corresponding to $\phi(j)=k$, where
    $j<i<k$. (Note that this term is non-zero: since the local
    multiplicity of $h$ at $i$ is one, the resolution of this crossing
    cannot introduce a new double-crossing.)  Of course, $H$ applied
    to this resolution is $\phi$; thus we have verified
    Equation~\eqref{eq:NullHomotopy} in this case.
  \item Suppose there are $k<i<j$ with $\phi(i)=j$ and $\phi(k)=i$. 
    Now, terms in $\partial H(\phi)$ and those in $H\partial (\phi)$
    pair off except for the single term in $\partial H(\phi)$
    where we resolve the crossing in $H(\phi)$ which did not already
    appear in $\phi$. This term, of course, gives $\phi$.
  \end{itemize}
  We have verified Equation~\eqref{eq:NullHomotopy}, and hence
  the homology of $I(\SetSS)\cdot\Alg(n,k)\cdot I(\SetTT)$ is trivial
  in this part of Case~\eqref{case:HorizontalCrossing}.

  For Case~\eqref{case:MultBiggerThanTwo}, we proceed similarly.  Now,
  choose $i$ so that the local multiplicity of $h$ just below $i$ is
  $1$ and the local multiplicity just above $i$ is $2$. 
  Indeed, we can
  choose $i$ so that the local multiplicity of $h$ at any point less
  than $i$ is strictly less than $2$.  We then let $j$ denote the
  element of $\SetSS$ just below $i$. Since we have already
  considered Case~\eqref{case:HorizontalCrossing}, we can also assume
  without loss of generality that the local multiplicity of $h$ just 
  below $j$ is zero.
  Once again, we define a map
  $H$ satisfying Equation~\eqref{eq:NullHomotopy}, by defining
  $H(\phi)$ for $\phi\colon \SetSS\to \SetTT$; only now,
  the cases are slightly different:
  \[
  H(\phi) =
  \begin{cases}
    0 & \phi(i) < \phi(j)\\
    \phi_{(i,j)} & \phi(j) < \phi(i).
  \end{cases}
  \]
  We verify Equation~\eqref{eq:NullHomotopy} for generators
  $\phi\co\SetSS\to\SetTT$.  In either case, there is a strand~$s_i$
  starting at~$i$ and a strand~$s_j$ starting at~$j$.
  \begin{itemize}
  \item Suppose that $\phi(i)<\phi(j)$. Then $H(\phi)=0$, and there is
    exactly one term in $\partial \phi$ for which $H(\phi)$ is
    non-zero, corresponding to the resolution of the crossing of $s_i$
    and~$s_j$.  (Again, the assumptions on the multiplicities
    guarantee that this term is non-zero.)
  \item Suppose that $\phi(i)>\phi(j)$. We divide
    the crossings in $\phi$ into two kinds: a {\em special} crossing
    is a crossing of a strand $s$ with $s_i$, which has the property
    that $s$ also crosses $s_j$. Any other crossing is said to be {\em
      generic}. If $\psi$ is a resolution of $\phi$ at a generic
    crossing, then it is easy to see that $\psi(i)>\psi(j)$, and hence
    that the terms in $H(\partial\phi)$ which come from generic
    resolutions of $\phi$ cancel corresponding terms in $\partial
    H(\phi)$.  If $\psi$ is a resolution of $\phi$ at a special
    crossing, however, $\psi(i)<\psi(j)$ and hence
    $H(\psi)=0$. Likewise, the corresponding terms in $\partial
    H(\phi)$ which are gotten by resolving the special crossings are
    also zero, as these resolutions all introduce double-crossings in
    $H(\phi)$. Finally, there is one left-over term in $\partial
    H(\phi)$ which does not appear in $H(\partial \phi)$, namely the term gotten by resolving the crossing between $s_i$
    and $s_j$; this gives $\phi$ back, as desired.
  \end{itemize}
  This completes the proof of the proposition.
\end{proof}

\begin{remark}
  We could alternately prove Proposition~\ref{prop:a-n-k-homology} by
  identifying $I(\SetSS)\cdot \Alg(n,k)\cdot I(\SetTT)$ with an
  interval in Hasse diagram of the Bruhat order of the symmetric group
  on $k$ letters, and applying the results of~\cite{BW82:BruhatOrder},
  which
  imply that such an interval is acyclic if it has more than one
  element.  The argument above can be viewed as an explicit proof of
  this fact.
\end{remark}

\begin{proof}[Proof of Theorem~\ref{thm:halg-support}]
  We first prove the necessary hypotheses as stated.
  \begin{enumerate}
  \item Without the compatibility condition,
    the degree $(\ell,h)$ part of $I(\SetS)\cdot \Alg(\PMC)\cdot I(\SetT)$
    is trivial.
  \item If the local multiplicity of $h$ is negative anywhere, then
    the degree $(*,h)$ part of $I(\SetS)\cdot \Alg(\PtdMatchCirc)\cdot
    I(\SetT)$ is empty.  Otherwise, suppose that the multiplicity of
    $h$ is bigger than~$1$ somewhere.  Consider the filtration on
    $I(\SetS)\cdot\Alg(\PtdMatchCirc)\cdot I(\SetT)$ given by the
    number of horizontal strands.  The differential on
    $\Alg(\PtdMatchCirc)$ can decrease the number of horizontal
    strands by at most one, in the case where we smooth a crossing
    with a horizontal strand.  We will show that the homology of the
    associated graded complex $C(\ell,h,*)=\bigoplus_{m\in\ZZ} C(\ell,h,m)$ 
    is zero.  Here, $\ell$ and $h$ are the gradings inherited from
    before, and $m$ is the newly-introduced grading, which we think
    of as the number of {\em non-horizontal} strands.
    We claim that $C(*,h,m)$ is 
    nearly identified with the corresponding subcomplex of $\Alg(4k,
    m)$.  The
    complex is not quite the same since we need to forbid horizontal
    strands; however, we can shift down the endpoints on the right by
    one half unit, to obtain an embedding of
    $C(*,h,m)$ into $I(\SetSS)\cdot\Alg(8k,m)\cdot I(\SetTT')\subset \Alg(8k,m)$
    (where here $\SetTT'$ is gotten from
    $\SetTT$ by shifting down by one half). If $h$ has local multiplicity
    $\geq 2$ somewhere, then so do the generators of
    $I(\SetSS)\cdot\Alg(8k,m)\cdot I(\SetTT')$. It follows from
    Lemma~\ref{lem:noneg-mult-exist} that the dimension of
    $I(\SetSS)\cdot \Alg(8k,m)\cdot I(\SetTT')$ is at least two, and hence, by
    Proposition~\ref{prop:a-n-k-homology}, that its homology is trivial.
  \item Suppose not; then $M(p_1)$ is in both the initial and final
    idempotent.  Modify the filtration from
    case~\eqref{item:halg-multiplicity} by considering the filtration
    on $I(\SetS)\cdot \Alg(\PtdMatchCirc)\cdot I(\SetT)$ given by the
    number of horizontal strands \emph{other than the one at $p_1$}
    (if there is one).  Again we can identify the associated graded
    pieces with an appropriate part of $\Alg(8k,m)$, where we shift all
    right endpoints other than $p_1$ down by half a unit.  
    (To see this identification, note that
    since $p_2\notin\Int(\supp(h))$, there are no crossings
    with the horizontal strand at $p_2$.) If we delete
    $p_1$ from the initial and final idempotent, by
    Lemma~\ref{lem:noneg-mult-exist} we can construct an algebra
    element with these endpoints in the corresponding 
    summand of $\Alg(8k,m-1)$. If we then add a horizontal strand
    at~$p_1$, we introduce an element with a
    crossing (since $p_1\in\Int(\supp(h))$). It follows from 
    Proposition~\ref{prop:a-n-k-homology} that the homology of this piece
    is trivial.
  \end{enumerate}
  If none of the other conditions apply, our homology class~$h$
    consists of a disjoint union of intervals with multiplicity~$1$,
    so that for every matched pair $\{p_1,p_2\}$ that is contained in
    both $\SetS$ and~$\SetT$, either both of $p_1$ and~$p_2$ are in
    $\Int(\supp(h))$ or neither is.  Let $m$ be the number where both
    are.  As before, consider the filtration by the number of
    horizontal strands. In this case the filtration is actually a
    grading (which the differential drops by one), since the only
    possible crossings are between non-horizontal strands and
    horizontal strands.  Then the
    complex $I(\SetS)\cdot \Alg(\PtdMatchCirc)\cdot I(\SetT)$ is
    isomorphic to the
    standard complex for the $m$-dimensional hypercube, with the
    gradings matching (up to an overall shift).  In particular, there
    is a unique element in
    $H_0$ (which is the lowest Maslov grading), represented by any
    generator in this grading.  Such a generator corresponds to a
    crossingless matching.
\end{proof}

\subsection{Massey products}
\label{sec:short-chords-massey-gen}

Recall (Corollary \ref{cor:AinfOnHomology}) that one can endow
$\HAlg(\PMC)$ with an $\Ainf$-structure so that $\Alg(\PMC)$ is
quasi-isomorphic to $\HAlg(\PMC)$. As discussed in
Section~\ref{sec:induced-Ainf-alg}, while many of the induced higher
products on $\HAlg(\PMC)$  are not entirely canonical, some of them are.
The aim of the present section is to show that the homology
$\HAlg(\PMC)$ is generated by canonically determined higher products
of chords of length one, in a suitable sense.
This will have as a corollary a certain rigidity of the algebra
(Proposition~\ref{prop:Rigidity}),
which will be useful for us in Section~\ref{sec:id-bim}.

\begin{proposition}\label{prop:massey-generate}
  Let $\zeta\in \HAlg(\PtdMatchCirc)$ 
  be a non-trivial, homogeneous homology class
  whose support has length greater than one.
  Then there is some $m>1$ and a Massey admissible sequence (in the
  sense of
  Definition~\ref{def:MasseyAdmissible}) of homogeneous 
  elements $\alpha_1,\dots,\alpha_m\in\HAlg(\PMC)$ with non-zero support
  with the property that
  $\zeta={\overline \mu}_{m}(\alpha_1,\dots,\alpha_m)$.
\end{proposition}

We will use the following technical lemma to ensure Massey admissibility.

\begin{lemma}
  \label{lemma:GradedMasseyProducts}
  Let $\Alg=\Alg(\PMC)$.  Let $\alpha_1,\dots,\alpha_m$ be
  a collection of homogeneous homology classes in $\HAlg$, and choose
  homogeneous representing cycles $a_1,\dots,a_m$. 
  Suppose that there are elements $\xi_{i,j}\in \Alg$ defined for $1\leq i<j\leq m$ 
  so that
  \begin{itemize}
  \item for $i=2,\dots,m$, 
    $\xi_{i-1,i}=a_i$;
  \item $d\xi_{i,j}=\sum\limits_{i<k<j}\xi_{i,k}\cdot \xi_{k,j};$
  \item for each $j>i+1$, we have $d\xi_{i,j}\neq 0$; and
  \item for $1 < i < m$,
    there is no algebra element whose support is
    $[\alpha_1]+\dots + [\alpha_i]$
    and whose initial idempotent agrees with the initial idempotent
    of $\alpha_1$.
  \end{itemize}
  Then the associated homology classes 
  $\alpha_1,\dots,\alpha_m$ form a Massey admissible sequence, and moreover
  $a_1\cdot \xi_{1,m}$ represents
  ${\overline\mu}_{m}(\alpha_1,\dots,\alpha_{m})$ (for any choice of compatible ${\overline \mu}_i$).
\end{lemma}

\begin{proof}
  We wish to apply Lemma~\ref{lem:MasseyAdmissible}.  When $1\leq
  i<j\leq m$, $\xi_{i,j}$ is a chain whose grading is $\lambda\cdot
  g(i+1,j)$, where $g(i,j)$ is the grading of
  $\overline{\mu}_{j-i+1}(\alpha_{i},\dots,\alpha_{j})$.  (We write
  now $\lambda\cdot g(i,j)$ in place of the $g(i,j)+1$ appearing in
  Definition~\ref{def:MasseyAdmissible}, where the grading set was
  $\ZZ$.) Since for $1<i<i+1<j\leq m$ the differential of $d\xi_{i,j}$
  is non-zero (by hypothesis), it follows from
  Theorem~\ref{thm:halg-support} that the homology group $\Halg_{\lambda\cdot
    g(i+1,j)}$ is $0$ for all $1\leq i<i+1 <j\leq m$.  Moreover, the
  final condition ensures that for $1<j<m$, ${\overline
    \mu}_i(\alpha_1,\dots,\alpha_j)=0$ since in fact there is no
  algebra element in the appropriate grading.  Similarly,
  $\Halg_{\lambda\cdot g(1,j)}$ is $0$.  Thus,
  $\alpha_1,\dots,\alpha_m$ is a Massey admissible sequence, and in
  fact Lemma~\ref{lem:MasseyAdmissible} applies (after we extend the
  $\xi_{i,j}$'s above by setting $\xi_{0,1}=a_1$ and $\xi_{0,j} = 0$
  for $1 < j < m$).
\end{proof}

Before proving Proposition~\ref{prop:massey-generate}, we introduce
some more terminology, and then some further lemmas.

Let $\zeta$ be an element of $\HAlg(\PtdMatchCirc)$ supported in
grading $(\ell,h)$.
A point $p$ in $\CircPts$ is called {\em fully unoccupied} if
$M(p)$ does not appear in either the initial or the final idempotent of
$\zeta$.
A point $p_1$ in $\CircPts$ is called {\em fully internal} if both $p_1$ and
its mate $p_2$ are contained in $\Int(\supp(h))$. It is
called {\em fully internal and unoccupied} if $p_1$ is fully internal
and $M(p_1)$ is not contained in the initial (and hence also the terminal)
idempotent. 

\begin{lemma}
  \label{lem:MasseyUnoccupied}
  Suppose that $\zeta$ is a non-trivial homology class in $\HAlg(\PtdMatchCirc)$
  supported in degrees $(\ell,h)$, and suppose that there is a point $p$
  which is fully internal and unoccupied.
  Then, we have a Massey admissible sequence of $m>1$ homology classes $\alpha_1,\dots,\alpha_m$, each of
  which has non-trivial support, with
  $$\zeta={\overline\mu}_{m}(\alpha_1,\dots,\alpha_m)$$
\end{lemma}

\begin{proof}
  According to Theorem~\ref{thm:halg-support}, $\zeta$ can be
  represented by a single diagram (rather than a formal sum of
  diagrams).  Now, any representative for $\zeta$ has some strand
  which crosses~$p$. In fact, we can find a representative for $\zeta$
  with the property that on the strand $s_1$ through $p$, there are no
  points which are fully internal and occupied.  We find a
  representative $a\in\Alg$ for $\zeta$ with the property that all the
  other strands of $a$ become stationary after the moment when $s_1$
  hits $p_1$, after which point only $s_1$ moves. See
  Figure~\ref{fig:MasseyFactoring} for an illustration.

\begin{figure}
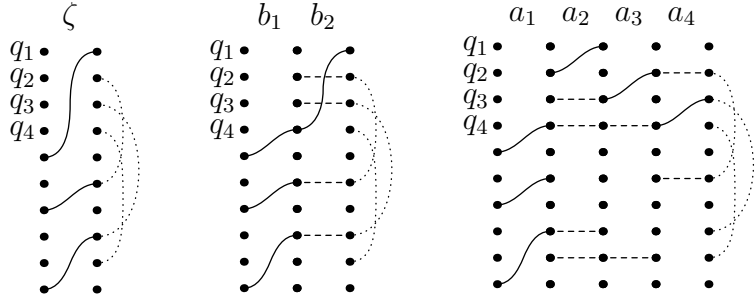

  \centering
  \[
  \mfigb{strands-60}\qquad
  \mfigb{strands-61}\qquad
  \mfigb{strands-62}
  \]
  \caption{\textbf{Illustration of the proof of
      Lemma~\ref{lem:MasseyUnoccupied}.}
    We start with the cycle
    $\zeta$, which has a strand $s_1$ which crosses the position
    $q_4$, which is fully internal and unoccupied.  The strand $s_1$
    then crosses two other positions ($q_3$ and $q_2$) which
    are matched with terminal points in $\zeta$. Thus, factoring the
    portion of the strand starting at $q_4$ to the right, we obtain a
    factorization of $\zeta$ into $b_1\cdot b_2$, where $b_2$ is not a
    cycle. The Massey factorization into $a_1,\dots,a_4$ is shown on the right.}
  \label{fig:MasseyFactoring}
\end{figure}

This provides a factorization (in $\Alg$) $a=b_1\cdot b_2$, where
$b_1$ is the part of $\zeta$ before $s_1$ crosses $p$, and $b_2$
corresponds to the portion of $s_1$ starting at $p$ (and whose initial
idempotent coincides with the terminal idempotent of $b_1$). Our
condition on the strand $s_1$ ensures that $b_1$ is a cycle. By
contrast, $b_2$ need not be a cycle: we can break $s_1$ into $m-1$
segments by positions $q_2,\dots,q_{m-1}$, where the $q_i$ are those
points in the interior of the support of $s_1$ but above $p$ which are
occupied but only partially internal in~$b_2$.  In order for $\zeta$ to be homologically
non-trivial, each $q_i$ must be matched with some
corresponding $q_i'$ contained in the support of $b_1$.  The $q_i'$
are necessarily in the terminal but not initial idempotent of $b_1$,
so that there
is a strand in $b_1$ entering $q_i'$. (The
$q_i'$ are not in the interior of $h$, since the $q_i$ were not fully
internal and occupied; we will use this observation towards the end of
the proof.)  We label the $q_i$ with the following conventions:
\begin{itemize}
\item  $q_1$ is the terminal point of the
  strand~$s_1$,
\item  $\{q_i\}$ are ordered in the opposite order to the order they are encountered
  along~$s_1$, and
\item $q_m=p$. 
\end{itemize}
We construct cycles $\{a_i\}_{i=1}^m$ as follows. Let $a_1=b_1$.
To define the other $a_i$, note that 
the terminal idempotent of $b_1$ 
has the form $I(\SetS)$, where $\SetS=\SetS_1 \setminus \{M(q_1)\}$ and
$\SetS_1$ includes all the $M(q_i)$. 
Now, for $i=2,\dots,m$, the initial idempotent
of $a_i$ is $\SetS_1\setminus \{M(q_{i-1})\}$ and
$a_i$ consists of a single moving
strand, which is the portion of $s_1$ from $q_{i}$ to $q_{i-1}$.
Clearly, each $a_i$ is a cycle.

Now, define $\xi_{i,j}$ for $1\leq i<j\leq m$ to be the algebra element
obtained from the substrand of $s_1$ which goes from $q_{j}$ to
$q_{i}$, and with initial idempotent $I(\SetS_1\setminus \{M(q_j)\})$.
In particular for $2 \le i \leq m$, we have $\xi_{i-1,i}=a_{i}$.  Any
algebra element with support $[a_1]+\dots+[a_i]$ with $i<m$ has
some strand which starts at $q_i$; but $M(q_i)$ is not in the
initial idempotent of $\alpha_1$, so there are no such elements. (Here
we are using the fact
that none of the $q_i$ was fully internal and occupied, as promised). 
Lemma~\ref{lemma:GradedMasseyProducts} applies, showing 
$\alpha_1,\dots,\alpha_m$ is Massey admissible, and further that
${\overline\mu}_{m}(\alpha_1,\dots,\alpha_m)=a_1\cdot
\xi_{1,m}=b_1\cdot b_2=[\zeta]$.
\end{proof}

\begin{lemma}
  \label{lem:MasseyBoundaryMatched}
  Suppose that $\zeta$ is a non-trivial homology class in
  $\HAlg(\PtdMatchCirc)$ supported in degrees $(\ell,h)$, and suppose
  that there is a point $p$ in the interior of the support of $h$
  which is matched with another point $p'$ on the boundary of the support.  Then, we have a Massey
  admissible sequence of $m>1$ homology classes
  $\alpha_1,\dots,\alpha_m$, each with non-trivial support,
  with
  $$\zeta={\overline\mu}_{m}(\alpha_1,\dots,\alpha_m).$$
\end{lemma}

\begin{proof}
  By hypothesis, $M(p)$ occurs in either the initial or the terminal
  idempotent of~$\zeta$, but not both.  We focus first on the case
  where $M(p)$ is in the initial idempotent.

  As in the proof of Lemma~\ref{lem:MasseyUnoccupied}, we find a
  strand $s_1$ which moves across $p$, and in fact, we can find a
  representative of $\zeta$ so that no points which are
  fully internal and occupied are encountered on that strand.  We then
  find a representative $a$ for $\zeta$ so that all the other strands
  of $a$ are stationary after the moment where $s_1$ hits~$p$.

  As before, this gives a factorization $a=b_1\cdot b_2$, where $b_2$
  corresponds to the substrand of $s_1$ starting at $p$.
  Once again, $b_1$ is a cycle. The strands picture for
  $b_2$ might have some crossings, which are in one-to-one
  correspondence with points in the support of $s_1$ above $p$ and
  which are matched with terminal points of~$b_1$.
  (Figure~\ref{fig:MasseyFactoring} can be modified to give a picture
  of the present case: simply cut off the bottom-most length one interval in the pointed matched
  circle, so that now $q_4$ is matched with an initial point of~$b_1$.)
  Our sequence $a_1,\dots, a_m$
  is now obtained by the same mechanism as in the proof of 
  Lemma~\ref{lem:MasseyUnoccupied}.

  The case where $M(p)$ is in the terminal idempotent of $\zeta$
  follows similarly. In this case, we find a representative for our
  strands so that all the strands of $z$ are stationary until the moment where
  $s_1$ hits $p$. This now gives a factorization $z=b_1\cdot b_2$, where 
  $b_1$ corresponds to the substrand of $s_1$ starting at its initial point and going
  until $p$, and $b_2$ consists of the rest of $s_1$ and all other strands. Now, $b_2$
  is a cycle. The strands picture for
  $b_1$ might have some crossings, which are in one-to-one
  correspondence with the points in the support of $s_1$ after its
  initial point but before $p$. A slight modification to Lemma~\ref{lemma:GradedMasseyProducts}
  then applies.
\end{proof}

\begin{lemma}
  \label{lem:NoTermInitMatcheds}
  Suppose that $\zeta$ is a non-trivial homology class in $\HAlg(\PtdMatchCirc)$
  supported in degrees $(\ell,h)$, and suppose that there is no point in
  the interior of the support of $h$ which is matched with either an
  initial or terminal point of $\zeta$, or which is fully
  unoccupied. Then, there is some terminal point which is not matched
  with any other point in the support of~$h$.
\end{lemma}

\begin{proof}
  Suppose on the contrary that every terminal point is matched with another
  point in $h$. Since no point in the interior is matched with an initial
  or terminal point, it follows that each terminal point is matched with an 
  initial point. Moreover, since there are no fully unoccupied points,
  and $\zeta$ is non-trivial,
  Criterion~\ref{item:halg-half-match} in 
  Theorem~\ref{thm:halg-support} ensures
  every point in the interior of $h$ is matched with another point
  in the interior of $h$. It is now easy to see that the support of $h$
  disconnects the one-manifold obtained by doing surgery on the
  matched pairs in $\PMC$, contrary to the definition of a pointed
  match circle.
\end{proof}

\begin{proof}[Proof of Proposition~\ref{prop:massey-generate}]
  Suppose $\zeta\in I(\SetS)\cdot\HAlg(\PMC)\cdot I(\SetT)$ is a non-trivial 
  homology class whose support has length greater than
  one.  Let $p_1$ be any point in the interior of the support of $\zeta$,
  and let $p_2$ be the other point with $M(p_1)=M(p_2)$.
  We have the following cases:
  \begin{itemize}
    \item If $M(p_1)\not\in\SetS\cup \SetT$, then 
      Lemma~\ref{lem:MasseyUnoccupied} provides the desired factorization.
    \item If $p_2$ is in the boundary of the support of $\zeta$,
      then  Lemma~\ref{lem:MasseyBoundaryMatched}
      provides the desired factorization.
    \item If $p_2$ is not contained in the closure of the support of $\zeta$,
      then but $M(p_1)\in\SetS\cup\SetT$, then this also forces
      $M(p_1)\in \SetS\cap\SetT$ (since $p_1$ is in the interior of the support of $\zeta$), 
      which, in view of
      Theorem~\ref{thm:halg-support} (Case~\ref{item:halg-half-match}),
      contradicts
      the assumption that $\zeta$ is homologically non-trivial.
    \end{itemize}
  Thus, we can assume that every interior point is occupied in both initial and terminal
  idempotents and is equivalent to some
  other point in the interior (i.e., it is fully internal and
  occupied). Moreover,
  Lemma~\ref{lem:NoTermInitMatcheds} ensures that there is a point $p$ on
  the boundary of the support which is not equivalent to any other
  point in the support of~$h$. We assume that $p$ is an initial
  point; the case when it is a terminal point is similar.

  Consider the point $q$ just above $p$. 
  We have shown that it is in the
  terminal idempotent, whether or not it is in the interior of the
  support; in any
  case, there is a representing cycle $a$ which contains a length one
  strand from $p$ to $q$. We can factor this strand off the right,
  so as to factor $a$ as a product of two cycles.
\end{proof}

As a digression, the following is a corollary of the above proofs,
showing that for the {\em split} genus $k$ matched circle,
Proposition~\ref{prop:massey-generate} can be strengthened.

\begin{corollary}\label{cor:split-generate}
  For the split genus~$k$ matched circle~$\PtdMatchCirc_{\text{spl}}$,
  $\HAlg(\PtdMatchCirc_{\text{spl}})$ is generated as an algebra by homology
  classes of Reeb chords of length~1.
\end{corollary}

\begin{proof}
  Suppose first there is a point~$p_1$ on the boundary of the support
  which is matched with a point~$p_2$ in the interior.  Then, since
  the since the distance between $p_1$ and $p_2$ is two, there
  is a unique point~$q$ in between them.  We can factor off the
  chord between $p_1$ and $q$ on one side or the other (depending
  on whether the idempotent of~$q$ is occupied or not), so that the
  two factors are both cycles.

  Otherwise, if Lemma~\ref{lem:MasseyUnoccupied} applies, the
  factorization of $a$ into
  algebra elements $b_1\cdot b_2$ is a factorization into cycles.

  In the remaining case of Proposition~\ref{prop:massey-generate}
  (when neither Lemma~\ref{lem:MasseyUnoccupied}
  nor~\ref{lem:MasseyBoundaryMatched} applies), we obtain a
  factorization into cycles as in the proof of
  Proposition~\ref{prop:massey-generate}.
\end{proof}

\begin{figure}[ht]
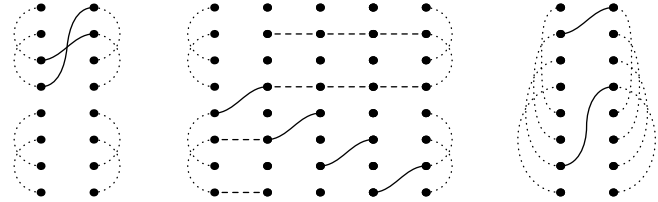

  \begin{gather*}
    \mfig{strands-50}\qquad\mfig{strands-55}\qquad\mfig{strands-51}
  \end{gather*}
\caption{\label{fig:Factor}
{\bf{Factorization in the algebra.}}
On the left, we have exhibited an element of the algebra
for the split genus two matched circle
which cannot be factored into length one chords; however, this
is not a cycle. (Note also that its differential
cannot be written as a product of length one chords.)
In the middle, there are 4 algebra elements in the split genus two
case with a well-defined, non-zero Massey product, showing that the
algebra is not formal.
On the right, we
have exhibited a homology class for the antipodal pointed matched
circle which can be written as a Massey product,
but not an ordinary product, of Reeb chords of length~$1$.
}
\end{figure}

Aspects of Proposition~\ref{prop:massey-generate}
are illustrated in Figure~\ref{fig:Factor}.
In particular, the second example shows that
$\Alg(\PtdMatchCirc_{\text{spl}})$ is not formal when the genus is
bigger than one, and the last
example shows that Corollary~\ref{cor:split-generate} is
not true for the antipodal matched circle of genus~2.

Proposition~\ref{prop:massey-generate} has the following consequence which will be used in 
the proof of Theorem~\ref{thm:Id-is-Id}:

\begin{proposition}
  \label{prop:Rigidity}
  Suppose that $\phi\colon \Alg(\PMC)\to \Alg(\PMC)$
  is a $G$-graded $\Ainf$ morphism such that $\phi_1(\xi)=\xi$ where $\xi$ is
  any Reeb chord of length one. Then
  $\phi_1$ induces the identity map on homology.
\end{proposition}

\begin{proof}
  Pre-composing $\phi$ with the standard quasi-isomorphism $f\co\Halg\to\Alg$,
  we reduce to the following statement:
  Let $\psi\co\Halg\longrightarrow \Alg$ be a morphism of $\Ainf$-algebras
  with the property that for each length $1$ chord $\xi$,
  $\psi(\xi)$ is a cycle representing the homology class of $\xi$, then
  $\psi_1$ induces the identity map on homology; i.e., if $\zeta$ is any homology class,
  then $\psi_1(\zeta)$ is a cycle representing $\zeta$.

  We prove this by induction on the total size of the support of
  $\zeta$.  Let $\zeta$ be a non-trivial homology class with support
  of size bigger than~$1$, and find a
  corresponding Massey admissible sequence $\alpha_1,\dots,\alpha_m$
  so that $\zeta=\overline{\mu}_m(\alpha_1,\dots,\alpha_m)$, where the
  support of each $\alpha_i$ is smaller than the support of $\zeta$.
  Define $\xi_{i,j}=\psi_{j-i}(\alpha_{i+1},\dots,\alpha_j)$, so that
  $\xi_{i-1,i}=\psi_1(\alpha_i)$. By the inductive hypothesis, 
  $\xi_{i-1,i}$ is a cycle representing the homology
  class $\alpha_i$. The $\Alg_\infty$ relation for $\psi$,
  together with Massey admissibility, ensures that 
  for $(i,j)\neq
  (0,m)$, 
  $$d\xi_{i,j}=\sum_{i<k<j}\xi_{i,k}\cdot \xi_{k,j}.$$
  Thus, Lemma~\ref{lem:MasseyAdmissible} ensures that
  $\sum_{0<k<m}\xi_{0,k}\cdot \xi_{k,m}$ is a cycle representing
  ${\overline \mu}_m(\alpha_1,\dots,\alpha_m)$.  On the other hand,
  another application of the $\Ainf$ relation (and Massey admissibility)
  gives 
  $$d\xi_{0,m}=\sum_{0<k<m}\xi_{0,k}\cdot \xi_{k,m}+\psi_1({\overline \mu}_m(\alpha_1,\dots,\alpha_m));$$
  i.e., $\psi_1(\zeta)=\psi_1({\overline \mu}_m(\alpha_1,\dots,\alpha_m))$ 
  represents the homology class ${\overline \mu}_m(\alpha_1,\dots,\alpha_m)=\zeta$.
\end{proof}

\begin{remark}
  The above proof works, provided $\phi$ preserves the relevant
  notions of homogeneity: i.e., it works whether $\phi$ is
  $\smallGroup$- or $\bigGroup$-graded. In the application (see
  Section~\ref{sec:id-bim}), we are interested in the case where
  $\phi$ is $\smallGroup$-graded.
\end{remark}

\begin{remark}\label{rmk:our-H-is-nilpotent}
  The $\Ainf$-structure on $\HAlg(\PMC)$ is nilpotent
  (Definition~\ref{def:Alg-bounded}), by the argument of
  Lemma~\ref{lem:Alg-nilpotent}.
\end{remark}

\subsection{A smaller model for \textalt{$\Alg(\PMC)$}{A(Z)}.}
\label{subsec:SmallerAlgebra}

Let $\PMC$ be a pointed matched circle.
Of course, $\HAlg(\PMC)$ is derived equivalent to $\Alg(\PMC)$.  Thus, for our
purposes, we could always work with $\Ainf$ modules over this
homology. This has the advantage that the underlying algebra has
smaller rank, but the disadvantage that now one must always keep track
of $\Ainf$ operations. There is, however, a natural intermediate
level: there is a differential graded algebra $\Alg'(\PMC)$ which is
a quotient of $\Alg$, but which is quasi-isomorphic to $\Alg$.

\begin{definition}
  Let ${\mathcal I}\subset \Alg(\PMC)$ denote the differential ideal
  generated by all algebra elements which have local multiplicity
  greater than $1$ somewhere.
\end{definition}

\begin{proposition}
  \label{prop:SmallerModel}
  The quotient map 
  $$\Alg(\PMC)\to \Alg(\PMC)/{\mathcal I}=\Alg'(\PMC)$$
  is a quasi-isomorphism.
  Moreover, the map
  sending $M\mapsto M\DTP_{\Alg(\PMC)}\Alg'(\PMC)$
  induces an equivalence of derived categories.
\end{proposition}

\begin{proof}
  The quasi-isomorphism statement is a direct consequence of
  Theorem~\ref{thm:halg-support}, Part~\eqref{item:halg-multiplicity}.
  The equivalence of derived categories statement follows, see
  Proposition~\ref{prop:quasi-iso-equivalence} (or indeed
  \cite[Theorem 10.12.5.1]{BernsteinLunts94:EquivariantSheaves}).
\end{proof}

In view of Proposition~\ref{prop:SmallerModel}, we could with use
$\Alg'(\PMC)$ in place of $\Alg(\PMC)$ throughout the present
paper. We chose not to do this for aesthetic reasons; but note that,
for practical calculations, it is indeed preferable to work in
$\Alg'(\PMC)$.

%%% Local Variables: 
%%% mode: latex
%%% TeX-master: "Bimodules"
%%% End: 

\section{Bordered Heegaard diagrams}
\label{sec:Diagrams}
In this section, we extend the notion of a bordered Heegaard diagram
of \cite[Chapter~\ref*{LOT:chap:heegaard-diagrams-boundary}]{LOT1} to
$3$-manifolds with two boundary
components. (The generalization to manifolds with more than two
boundary components is straightforward, and we mostly leave it to the
interested reader; see also Remark~\ref{remark:more-bdy-comps}.) This
generalization was first sketched in the
appendix to~\cite{LOT1}, which the reader may want to consult for a
condensed treatment. 

First, we recall the notion of strongly bordered three-manifolds with
two boundary components, introduced in
Definition~\ref{def:IntroStronglyBordered}. We adapt the definition slightly,
so that borderings are specified by pointed matched circles, and 
give notation which will be used later.

\begin{definition}
  \label{def:StronglyBordered}
  A {\em strongly bordered three-manifold with two boundary components}~$\sbY$
  is a tuple 
  $${\mathcal Y}=(Y,\PMC_L,\Delta_L,z_L,\phi_L,\PMC_R,\Delta_R,z_R,\phi_R,
  \gamma_z,\nu_z)$$ where:
  \begin{itemize}
    \item $Y$ is a compact, oriented three-manifold with two 
      boundary components $\partial_L Y$ and $\partial_R Y$.
    \item $\Delta_L\subset \partial_L Y$ and $\Delta_R \subset \partial_R Y$
      are preferred disks.
    \item $z_L\in\partial\Delta_L$ and $z_R\in\partial \Delta_R$
      are basepoints on the boundaries of the disks.
    \item $\phi_L$ and $\phi_R$
      are homeomorphisms
      \begin{align*}
        \phi_L&\co (F(\PMC_L),D_L,z_L)\to
        (\partial_L Y,\Delta_L,z_L)\\
        \phi_R&\co (F(\PMC_R),D_R,z_R)\to
        (\partial_R Y,\Delta_R,z_R).
      \end{align*}
      (Here, $D_L$ and $D_R$ are the preferred disks in $F(\PMC_L)$
      and $F(\PMC_R)$, equipped with the basepoints $z_L$ and $z_R$
      coming from the pointed matched circles; see
      Construction~\ref{construct:PMC-gives-surf}. In the interest of
      notational simplicity, we do not distinguish the notation for
      the basepoints in the model surface from the preferred
      basepoints on $\partial Y$.)
    \item $\gamma_z$ is a path in $Y$ connecting $z_L$ to $z_R$.
    \item $\nu_z$ is an isotopy class of nowhere vanishing normal
      vector fields to $\gamma_z$ pointing into $\Delta_L$ at $z_L$
      and into $\Delta_R$ at $z_R$.
    \end{itemize}
\end{definition}

We wish to describe Heegaard diagrams which specify the above data.
Before doing this, we recall how we specify bordered three-manifolds
(with one boundary component) by diagrams.
\begin{definition}\label{sec:bordered-HD-one-bdy-component}
  A \emph{pointed bordered Heegaard diagram with one boundary
    component} is a quadruple
  $\HD=(\overline{\Sigma},\overline{\alphas},\betas,z)$ where
  ${\overline\Sigma}$ is a compact surface of genus $g$ with one
  boundary component; $\betas$ is a $g$-tuple of pairwise disjoint
  curves in the interior $\Sigma$ of $\overline{\Sigma}$;
  \[
  \overline{\alphas}=\{\overbrace{\overline{\alpha}_1^a,\dots,\overline{\alpha}_{2k}^a}^{\overline{\alphas}^a},\overbrace{\alpha_1^c,\dots,\alpha_{g-k}^c}^{\alphas^c}\}
  \]
  is a collection of pairwise disjoint embedded arcs (the
  $\overline{\alpha}_i^a$) with boundary on $\bdy\overline{\Sigma}$
  and circles (the $\alpha_i^c$) in the interior $\Sigma$ of
  $\overline{\Sigma}$; and $z$ is a basepoint in
  $\bdy\overline\Sigma\setminus\alphas^a$. We require that
  $\overline\Sigma\setminus\overline\alphas$ and
  $\overline\Sigma\setminus\betas$ both be connected; this translates
  to the condition that the $\alpha$- (respectively $\beta$-) curves
  be linearly independent in
  $H_1(\overline{\Sigma},\bdy\overline{\Sigma})$.
\end{definition}

Given a bordered Heegaard diagram $\HD=(\Sigma,\alphas,\betas,z)$, the
boundary $\bdy\overline{\Sigma}$ is a pointed matched circle in a
natural way, with
$\mathbf{a}=\overline\alphas^a\cap\bdy\overline{\Sigma}$, $M$ the
matching pairing up the endpoints of each $\overline\alpha_i^a$, and
$z$ the basepoint $z$. We call this pointed matched circle $\PMC(\HD)$.

\begin{construction}\label{construct:bordered-HD-to-bordered-Y-one-bdy}
  A pointed bordered Heegaard diagram $\HD$ with one boundary
  component specifies a $3$-manifold $Y(\HD)$ with one boundary
  component as follows:
  \begin{enumerate}
  \item Thicken $\overline{\Sigma}$ to $\overline{\Sigma}\times[0,1]$.
  \item Attach three-dimensional two-handles along the
    $\alpha$-circles in ${\overline\Sigma}\times\{0\}$.
  \item Attach three-dimensional two-handles along the $\beta$-circles
    in ${\overline\Sigma}\times\{1\}$.
  \end{enumerate}
  \begin{figure}
    \centering
    \includegraphics{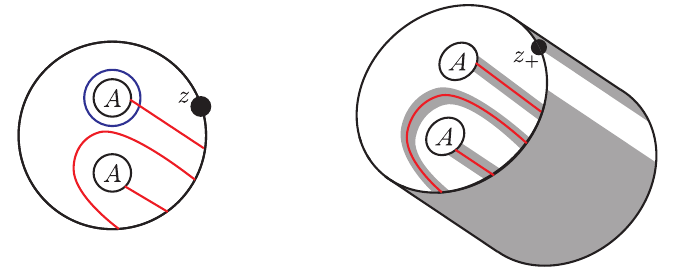}
    \caption{\textbf{A pointed bordered Heegaard diagram and the
        associated $3$-manifold.} The picture on the left is a
      Heegaard diagram for the bordered solid torus shown on the
      right.  The shaded part of the boundary is
      $\PunctF(\PMC(\HD))$.}
    \label{fig:PBHD-3-mfld}
  \end{figure}
  A parameterization of the boundary is specified as follows.
  Consider the graph
$$\left({\overline\alphas}^a\cup \left(\partial{\overline\Sigma}\setminus\nbd(z)\right)\right)\times \{0\}\subset {\overline\Sigma}\times\{0\},$$
thought of as a subset of $\partial Y$. The closure $\PunctF$ of a
neighborhood of this graph is naturally identified with the
surface-with-boundary $\PunctF(\PMC(\HD))$ associated to the pointed
matched circle $\PMC(\HD)$. (This identification \emph{reverses}
orientation with our orientation conventions.  Note that this
neighborhood of the $\alpha$-arcs is orientation-reversing
homeomorphic to a portion of $\bdy Y(\HD)$, so that the orientation of
$\PunctF(\PMC(\HD))$ agrees with the orientation of $\bdy Y(\HD)$.)
The deleted neighborhood of $z$ is an interval with endpoints $z_-$
and $z_+$, and we can take either of them (say, $z_+$) as
corresponding to the basepoint on the boundary.  The complement of
$\PunctF$ in $\partial Y$ is a disk. See Figure~\ref{fig:PBHD-3-mfld}.
\end{construction}

Thus fortified, we turn to the two boundary component case.
\begin{definition}
  \label{def:ArcedBordered}
  An \emph{arced bordered Heegaard diagram with two
    boundary components} is a quadruple
  $(\overline{\Sigma},\overline{\alphas},\betas,{\mathbf z})$ where
  \begin{itemize}
  \item ${\overline\Sigma}_g$ is a
    compact surface of genus $g$ with two boundary components,
    $\bdy_L\overline\Sigma$ and $\bdy_R\overline\Sigma$;
  \item $\betas$ is a $g$-tuple of pairwise
    disjoint curves in the interior $\Sigma$ of $\overline{\Sigma}$;
  \item 
    \(
    \overline{\alphas}=\{\overbrace{\overline{\alpha}_1^{a,L},\dots,\overline{\alpha}_{2\Lgen}^{a,L}}^{\overline{\alphas}^{a,L}},\overbrace{\overline{\alpha}_1^{a,R},\dots,\overline{\alpha}_{2\Rgen}^{a,R}}^{\overline{\alphas}^{a,R}},\overbrace{\alpha_1^c,\dots,\alpha_{g-\Lgen-\Rgen}^c}^{\alphas^c}\}
    \)
    is a collection of pairwise disjoint embedded arcs with boundary
    on $\bdy_L\overline{\Sigma}$ (the $\overline{\alpha}_i^{a,L}$\!), arcs with
    boundary on $\bdy_R\overline{\Sigma}$ (the $\overline{\alpha}_i^{a,R}$\!), and
    circles (the $\alpha_i^c$) in the interior $\Sigma$ of
    $\overline{\Sigma}$; and
  \item ${\mathbf z}$ is a path in $\overline{\Sigma}\setminus(\overline\alphas\cup\betas)$ between
    $\bdy_L\overline{\Sigma}$ and $\bdy_R\overline{\Sigma}$.
  \end{itemize}
  These are required to satisfy:
  \begin{itemize}
  \item $\overline{\Sigma}\setminus\overline{\alphas}$ and
    $\overline{\Sigma}\setminus\betas$ are connected and
  \item $\overline{\alphas}$ intersects $\betas$ transversely.
  \end{itemize}
\end{definition}

An arced bordered Heegaard diagram $\HD$ with two boundary components
specifies two pointed matched circles
\begin{align*}
  \PMC_L(\HD)&=(\bdy_L\overline{\Sigma},\overline\alphas^{a,L}\cap\bdy_L\overline{\Sigma},M_L,{\mathbf z}\cap\bdy_L\overline{\Sigma}) \quad\text{and}\\
  \PMC_R(\HD)&=(\bdy_R\overline{\Sigma},\overline\alphas^{a,R}\cap\bdy_R\overline{\Sigma},M_R,{\mathbf z}\cap\bdy_R\overline{\Sigma})
\end{align*}
where $M_L$ (respectively $M_R$) is the pairing matching the
endpoints of each arc in $\overline{\alpha}_i^{a,L}$ (respectively $\overline{\alpha}_i^{a,R}$).

\begin{figure}
  \centering
  \includegraphics{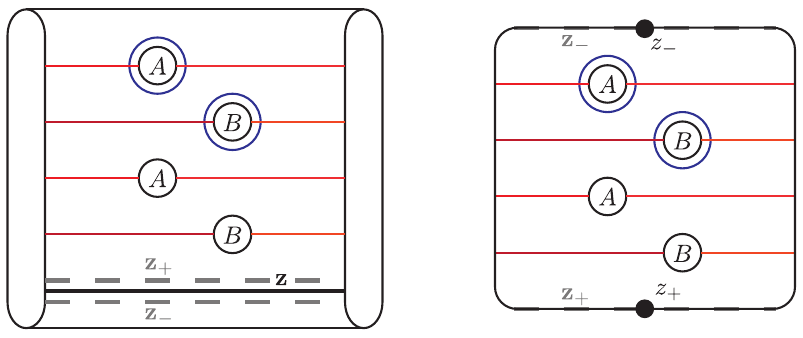}
  \caption{\textbf{Drilling Heegaard diagrams.} Left: an arced
    bordered Heegaard diagram $\HD$ for $T^2\times[0,1]$. Right: the
    result $\drHD$ of drilling the tunnel $\mathbf{z}$ from $\HD$, a
    bordered Heegaard diagram with a single boundary component.}
  \label{fig:drillingHD}
\end{figure}

\begin{definition}
  \label{def:Drilling}
  Given an arced bordered Heegaard diagram
  $\HD=(\overline{\Sigma},\overline{\alphas},\betas,\mathbf{z})$ with two
  boundary components, there is a bordered Heegaard diagram $\drHD$
  with one boundary component, obtained by deleting a neighborhood of
  the arc ${\mathbf z}$ from $\overline\Sigma$. The boundary of the
  deleted neighborhood of $\mathbf{z}$ consists of two disjoint
  push-offs ${\mathbf z}_+$ and ${\mathbf z}_-$ of ${\mathbf z}$, and
  $\partial {\overline \Sigma}_{dr}$ is
  \[
  \left(\partial_L{\overline \Sigma}\setminus \nbd(z_L)\right)
  \cup
  \left(\partial_R{\overline \Sigma}\setminus \nbd(z_R)\right)
  \cup {\mathbf z}_+
  \cup {\mathbf z}_-.
  \]
  
  We can equally well choose to put the basepoint of $\drHD$ on
  ${\mathbf z}_+$ or on ${\mathbf z}_-$. We will call these two
  choices of basepoint $z_+$ and $z_-$ respectively.
  For either of these choices, we
  call the pointed bordered Heegaard diagram $\drHD$ a diagram
  obtained from $\HD$ by {\em drilling}.  Note that $\PMC(\drHD)$ is
  the pointed matched circle $\PMC_L(\HD)\connectsum\PMC_R(\HD)$.

  We will use the notation $\drtHD$ to denote the doubly pointed
  Heegaard diagram obtained by viewing both $z_+$ and $z_-$ as
  basepoints of the drilled diagram.
\end{definition}
See Figure~\ref{fig:drillingHD} for an illustration of the drilling
construction on Heegaard diagrams.

\begin{figure}
  \centering
  \includegraphics{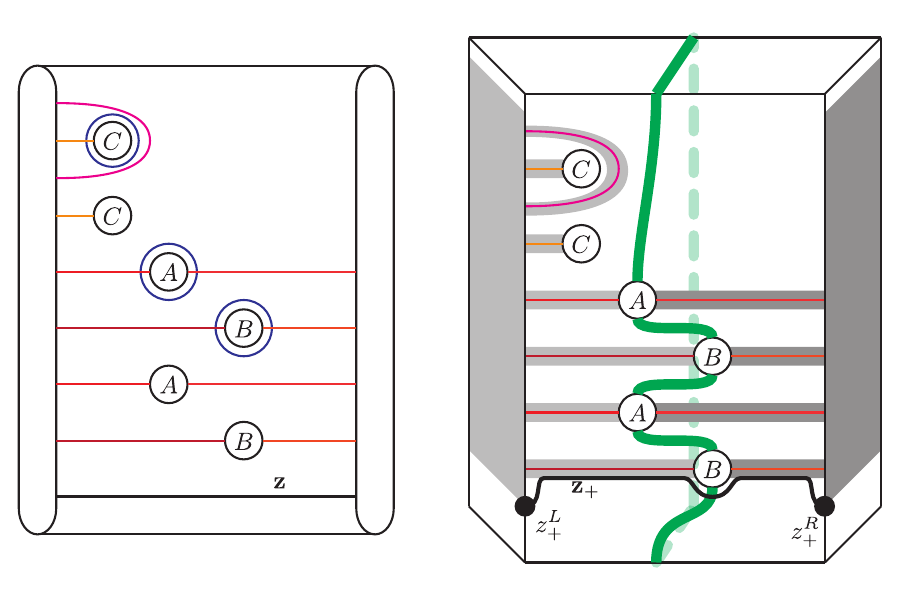}
  \caption{\textbf{Constructing a bordered 3-manifold with two
      boundary components from an arced bordered Heegaard diagram.}
    The Heegaard diagram on the left represents an elementary cobordism from the
    genus two surface to the genus one surface. The lightly
  shaded region on the right picture is
$\PunctF(\partial_L\HD)$ (a surface of genus two), while the darkly shaded one is
$\PunctF(\partial_R\HD)$ (a surface of genus one). There is a
$2$-handle attached along the thick (green) curve.}
  \label{fig:2bdy-construct}
\end{figure}

\begin{construction}\label{construct:framed-3-mfld}
  An arced bordered Heegaard diagram with two boundary components
  in the sense of Definition~\ref{def:ArcedBordered} gives rise
  to a strongly bordered three-manifold in the sense of
  Definition~\ref{def:StronglyBordered}, as follows.  Let $\drHD$ be
  the Heegaard diagram obtained from $\HD$ by drilling
  (Definition~\ref{def:Drilling}). The boundary of $Y(\drHD)$ is
  decomposed as a connect sum $F(\PMC_L)\connectsum F(\PMC_R)$. We attach a
  three-dimensional two-handle along the connect sum annulus to obtain
  $Y$. To see the other structure, we perform this construction with
  more care. The boundary of $Y(\drHD)$ consists of three pieces:
\begin{itemize}
  \item A neighborhood of the graph 
    $$\left({\overline\alphas}_L^a\cup 
      \left(\partial_L{\overline\Sigma}\setminus\nbd(z_L)\right)\right)\times 
    \{0\}\subset {\overline\Sigma}\times\{0\},$$
    whose closure
    is identified in an orientation-reversing way with
    $\PunctF(\PMC_L)$, which we denote $\PunctF_L$. (Note that
    $\PunctF_L$ contains the basepoint $z^+_L$ on its boundary.)
  \item A neighborhood of the graph 
    $$\left({\overline\alphas}_R^a\cup 
      \left(\partial_R{\overline\Sigma}\setminus\nbd(z_R)\right)\right)\times 
    \{0\}\subset {\overline\Sigma}\times\{0\},$$
    whose closure is identified in an orientation reversing way with
    $\PunctF_R \coloneqq \PunctF(\PMC_R)$.
    (Note that
    $\PunctF_L$ contains the basepoint $z^+_R$ on its boundary.)
  \item An annulus $A$, equipped with a path ${\mathbf z}^+$ connecting
    $z^+_L$ to $z^+_R$. (This is the path ${\mathbf z}_+$ on the
    boundary of
    ${\overline\Sigma}_{dr}={\overline\Sigma}_{dr}\times\{0\}$,
    thought of as a subset of the boundary of $Y(\drHD)$.)
\end{itemize}    
Now, we attach a three-dimensional two-handle to $Y(\drHD)$ along the
annulus $A$. More precisely, let $\Delta$ be a two-dimensional disk. We glue
$\Delta\times[0,1]$ to $Y(\drHD)$, identifying $(\partial \Delta)
\times[0,1]$ with the annulus $A$, so that
$(\partial\Delta)\times\{0\}$ is glued
to the boundary of $\PunctF_L$, while $\Delta\times\{1\}$ is glued to
the boundary of $\PunctF_R$. Let $\Delta_L=\Delta\times\{0\}$ and
$\Delta_R=\Delta\times\{1\}$. It is
easy to see that this gives a three-manifold homeomorphic to
$Y$. Moreover, in this model, the boundary of $Y$ consists of the
disjoint union of $\PunctF_L\cup\Delta_L$ and $\PunctF_R\cup\Delta_R$, and
$\Delta_L$ and $\Delta_R$ respectively contain
$z_L^+$ and $z_R^+$ on their boundary. Our preferred disks are
$\Delta_L$ and $\Delta_R$, and their basepoints are $z_L^+$ and
$z_R^+$. The homeomorphisms $\phi_L$ and
$\phi_R$ are supplied by
the pointed matched circles. The arc $\gamma_z$ is supplied by
${\mathbf z}_+$, which in turn gives a path on $(\partial \Delta)\times [0,1]$.
The framing is specified so as to point into $\Delta\times\{t\}$ for each 
$t\in [0,1]$. See Figure~\ref{fig:2bdy-construct}.
\end{construction}

We call the data 
$(Y,\Delta_L,z_L,\Delta_R,z_R, \phi_L,\phi_R,\gamma_z)$
from Construction~\ref{construct:framed-3-mfld} 
the strongly bordered
$3$-manifold associated to the arced bordered diagram $\HD$.
We will often abuse notation and use $Y$ or
$\sos{\phi_L}{Y}{\phi_R}$ or $\sos{\PMC_L}{Y}{\PMC_R}$ to denote
all the data of a strongly bordered $3$-manifold (depending on which
pieces we want to emphasize).

\begin{remark}\label{remark:more-bdy-comps}
  In the case of strongly bordered Heegaard diagrams with more than two
  boundary components, one replaces the arc ${\mathbf z}$ with a tree (or, as a
  special case, a sequence of arcs) connecting the various boundary
  components.
\end{remark}

There is an inverse to the drilling construction, {\em filling},
defined as follows.
\begin{definition}
  \label{def:Filling}
  Let $\drHD=(\overline\Sigma,\overline\alphas,\betas,z_+)$ be a
  pointed bordered Heegaard diagram with one boundary component. Let
  $z_-$ be another point in $\bdy\overline\Sigma$, and suppose that,
  writing $\bdy\overline\Sigma\setminus \{z_+,z_-\}=I_L\amalg I_R$,
  there are no $\alpha$-arcs running between $I_L$ and $I_R$. Then
  $\{z_+,z_-\}$ decomposes $\PMC(\drHD)$ as a connect sum,
  $\PMC(\drHD)=\PMC_L\connectsum\PMC_R$. Attaching a band (two-dimensional
  one-handle) to $\drHD$ between $z_+$ and $z_-$ gives a Heegaard diagram $\fHD$
  with two boundary components, with $\bdy_L{\fHD}=\PMC_L$ and
  $\bdy_R{\fHD}=\PMC_R$.
\end{definition}

The three-manifold $Y(\fHD)$ is obtained from $Y(\drHD)$ by attaching
a $3$-dimensional $2$-handle to $\bdy Y(\drHD)$ along the connect sum
curve.

\subsection{Arced bordered Heegaard diagrams}

We can use the drilling construction to re\-phrase questions about
arced bordered Heegaard diagrams with two boundary components in terms
of ordinary (one boundary component) bordered three-manifolds. For example,
we have the following:

\begin{definition}\label{def:diagrams-multidef}
  Let $\HD$ be an arced bordered Heegaard diagram.
  \begin{itemize}
  \item A \emph{generator} of $\HD$ is a generator of $\drHD$. We let
    $\Gen(\HD)$ denote the set of generators of $\HD$.
  \item Given generators $\x,\y\in\Gen(\HD)$, the set of \emph{domains
      connecting $\x$ and $\y$}, $\pi_2(\x,\y)$, is the set of domains
    in $\drHD$ connecting $\x$ to $\y$ that do not cross either $z_+$
    or $z_-$. We view domains as linear combinations of components of
    $\overline{\Sigma}\setminus (\alphas\cup\betas)$. Recall that
    $\bdy^\bdy(B)$ denotes the intersection of $\bdy B$ with
    $\bdy\overline{\Sigma}$.
  \item Let $\pi_2^\bdy(\x,\y)=\{B\in\pi_2(\x,\y)\mid
    \bdy^\bdy B=0\}$. These are the \emph{provincial} domains from
    $\x$ to $\y$.

    There are natural isomorphisms
    \begin{align*}
      \pi_2(\x,\x)&\cong H_2(Y(\HD),\bdy Y(\HD))\\
      \pi_2^\bdy(\x,\x)&\cong H_2(\drY(\HD))
    \end{align*}
    corresponding to \cite[Lemmas~\ref*{LOT:lem:pi2-h2}
    and~\ref*{LOT:lem:pi2-h2-bdy}]{LOT1}.
  \item We call elements of $\pi_2(\x,\x)$ \emph{periodic domains.}
  \item The arced bordered Heegaard diagram with two boundary
    components $\HD$ is called \emph{admissible} (respectively
    \emph{provincially admissible}) if the associated drilled Heegaard
    diagram $\drtHD$ is admissible (respectively provincially
    admissible) in the sense of
    \cite[Definition~\ref*{LOT:def:admissibility}]{LOT1}
    (respectively
    \cite[Definition~\ref*{LOT:def:provincial-admissibility}]{LOT1}),
    i.e., if every nontrivial
    periodic domain (respectively provincial periodic domain) of $\HD$
    has both positive and negative coefficients.  
  \end{itemize}
\end{definition}

\begin{proposition}\label{prop:diagrams-exist}
  Any strongly bordered three-manifold with two boundary components
  comes from an admissible arced bordered Heegaard diagram with two
  boundary components.
\end{proposition}
\begin{proof}
  Choose a bordered Heegaard diagram for the bordered three-manifold
  (with one boundary component) $Y\setminus \nbd(\gamma_z)$. This can
  be done according to
  \cite[Lemma~\ref*{LOT:lem:3mfld-heegaard}]{LOT1}; moreover, it can
  be made admissible by
  \cite[Proposition~\ref*{LOT:prop:heegaard-moves}]{LOT1}. (Note that
  admissibility with only one basepoint which is gotten 
  from \cite[Proposition~\ref*{LOT:prop:heegaard-moves}]{LOT1} is
  slightly stronger than the
  admissibility we require here.) The filling construction of
  Definition~\ref{def:Filling} then produces the desired diagram
  for~$Y$.
\end{proof}

\begin{proposition}\label{prop:heegaard-moves}
  If $\HD$ and $\HD'$ specify the same strongly bordered $3$-manifold
  then $\HD$ and $\HD'$ are related by a sequence of the following
  moves:
  \begin{itemize}
  \item Isotopies of the $\alpha$- and $\beta$-curves.
  \item Handleslides among the $\alpha$-circles and among the $\beta$-circles.
  \item Handleslides of an $\alpha$-arc over an $\alpha$-circle.
  \item Stabilizations of the diagram.
  \end{itemize}
  Moreover, if $\HD$ and $\HD'$ are admissible (respectively
  provincially admissible) then the moves can be chosen so that all
  intermediate diagrams are admissible (respectively provincially
  admissible).
\end{proposition}
\begin{proof} 
  Suppose $\HD$ and $\HD'$ are provincially admissible and specify the
  same strongly bordered $3$-manifold.
  Then $\drHD$ and $\drHD'$ specify the same bordered three-manifold (with one
  boundary component). Thus, $\drHD$ and $\drHD'$ can be connected
  by a sequence of provincially admissible
  Heegaard moves which do not cross either of the basepoints,
  as in \cite[Propositions~\ref*{LOT:prop:heegaard-moves}
  and~\ref*{LOT:prop:admis-achieve-maintain}]{LOT1}. Filling all the diagrams,
  we obtain the desired sequence connecting $\HD$ to $\HD'$.

  The case when $\HD$ and $\HD'$ are admissible is similar, except
  that $\drHD$ and $\drHD'$ are not necessarily admissible, as there
  may be periodic domains with positive coefficients crossing the
  extra basepoint (say $z_-$, if $z_+$ was the basepoint for $\HD$).  The
  doubly-pointed diagrams $\drtHD$ and $\drtHD'$ are admissible (in
  the obvious sense), and
  \cite[Propositions~\ref*{LOT:prop:heegaard-moves}
  and~\ref*{LOT:prop:admis-achieve-maintain}]{LOT1} adapt easily to the
  doubly-pointed case.
\end{proof}
We call two diagrams which are related by the moves of
Proposition~\ref{prop:heegaard-moves} \emph{equivalent.}

In the case of $3$-manifolds with two boundary components, we can
refine some of the notions related to domains.  Given a domain $B$,
$\bdy^\bdy_L B$ (respectively $\bdy^\bdy_RB$) denote the intersection
of $\bdy B$ with $\bdy_L\overline{\Sigma}$ (respectively
$\bdy_R\overline{\Sigma}$).  Let
$\pi_2^{\bdy_L}(\x,\y)=\{A\in\pi_2(\x,\y)\mid \bdy^{\bdy}_LA=0\}$
denote the set of \emph{left-provincial} domains connecting $\x$ and
$\y$, and $\pi_2^{\bdy_R}(\x,\y)=\{A\in\pi_2(\x,\y)\mid
\bdy^{\bdy}_RA=0\}$ the set of \emph{right-provincial} domains
connecting $\x$ and $\y$.  It is easy to show that
\begin{align*}
  \pi_2^{\bdy_L}(\x,\x)&\cong H_2(Y(\HD),\bdy_RY(\HD))\\
  \pi_2^{\bdy_R}(\x,\x)&\cong H_2(Y(\HD),\bdy_LY(\HD)).
\end{align*}

\begin{definition}\label{def:admissible-diag}
  An arced bordered Heegaard diagram $\HD$ with two boundary
  components is called \emph{left} (respectively \emph{right})
  \emph{admissible} if every nontrivial right-provincial (respectively
  left-provincial) periodic domain has both positive and negative
  coefficients.
\end{definition}
It is easy to show that Proposition~\ref{prop:heegaard-moves} still
holds if one replaces ``admissible'' by ``left admissible'' or ``right
admissible.''

\begin{lemma}\label{lemma:left-right-Area}
  The Heegaard diagram $\HD$ is left (respectively right) admissible
  if and only if there is an area form on $\overline{\Sigma}$ with
  respect to which every right-provincial (respectively
  left-provincial) periodic domain has signed area $0$. The diagram
  $\HD$ is admissible if and only if there is an area form on
  $\overline{\Sigma}$ with respect to which every periodic domain has
  signed area $0$.
\end{lemma}
\begin{proof}
  The proof is exactly the same as the proof of
  \cite[Lemma~\ref*{LOT:lemma:admiss-reform}]{LOT1}, which in turn is
  the same as the proof of \cite[Lemma~4.12]{OS04:HolomorphicDisks}.
\end{proof}

Note that
\begin{multline*}
\text{admissible}\implies\text{left or right
  admissible}\implies\\ \text{left and right
  admissible}\implies \text{provincially admissible}.
\end{multline*}
All of these implications are strict.

Finally we discuss how $\SpinC$-structures on manifolds with two
boundary components relate to arced, bordered Heegaard diagrams. 

Before doing this, we recall some generalities (see~\cite{OS04:Knots}).  Suppose that $M$ is a
three-manifold, equipped with an oriented, null-homologous knot
$C$. Then, there is a notion of {\em relative $\SpinC$ structures},
denoted $\RelSpinC(M,C)$. These are defined to be $\SpinC$ structures
on the zero-surgery manifold $M_0(C)$.  If $C$ is equipped with a
Seifert surface $F$, there is an identification 
$$\RelSpinC(M,C)\cong \SpinC(M)\oplus\ZZ,$$
The projection onto $\ZZ$ is  gotten by
$${\underline\spinc}\mapsto \OneHalf \langle c_1({\underline\spinc}), {\widehat F}\rangle,$$
where ${\widehat F}$ is gotten by closing off $F$ in $M_0(C)$.
This projection to $\ZZ$ is called the {\em Alexander grading} on relative $\SpinC$
structures.

Let $\sbY$ be a strongly bordered $3$-manifold with two boundary
components specified by an arced bordered Heegaard diagram $Y$. A
meridian $C$ of $\gamma_z\subset Y$ specifies a knot~$K$ in
$\drY=Y(\drHD)$. This knot is null-homologous. Indeed, the
surface $\PunctF(\PMC_R Y)$ provides a natural choice of Seifert
surface for $K$ (and in particular an orientation for~$K$).

\begin{lemma}
  \label{lem:RelAbsSpinC}
  There is a natural identification $\RelSpinC(\drY,K)\cong\SpinC(Y)$.
  Under this identification, the Alexander grading of $\relspinc$
  corresponds to the evaluation of the corresponding $\SpinC$ structure
  on $\partial_R Y$.
\end{lemma}

\begin{proof}
  Let $Y'$ denote zero-surgery on $\drY$ along $K$.
  It is easy to see that $Y'$ is naturally identified with
  the three-manifold obtained from $Y$ by attaching a one-handle
  to its boundary connecting the left and right basepoints. The
  identification
  $\RelSpinC(\drY,K)\cong\SpinC(Y)$
  follows at once. Under this identification, $\partial_R Y$ is clearly
  homologous to the capped-off Seifert surface.
\end{proof}

As in~\cite{OS04:Knots}, an oriented knot is specified by a Heegaard
diagram with two (ordered) basepoints $z_+$ and $z_-$ (denoted $w$ and
$z$ in~\cite{OS04:Knots}). The oriented knot in the three-manifold is
specified as follows: draw an arc from $z_-$ to $z_+$ which crosses
only $\beta$-circles, and then close this up by drawing an arc from
$z_+$ to $z_-$ which crosses only $\alpha$-circles. (It might be
necessary to push the two arcs into the two handlebodies to make the
knot be embedded.) 

\begin{lemma}
  \label{lem:DrilledDiagram}
  The Heegaard diagram $\drtHD$ is a doubly-pointed bordered
  Heegaard diagram for $(\drY,K)$, where $K$ is oriented as the
  boundary of $F^\circ(\PMC_RY)$.
\end{lemma}
\begin{proof}
  Let $\gamma$ be a path in $\Sigma_{dr}$ connecting $z_+$ to $z_-$ in
  the complement of the $\alpha$-curves and $\eta$ a path in
  $\Sigma_{dr}$ connecting $z_-$ to $z_+$ in the complement of the
  $\beta$-circles. Then the push-off of $\gamma\cup\eta$ is the knot specified by
  $\drtHD$.  We can choose $\gamma$ to lie in a neighborhood of the
  $\alpha^R$-arcs and $\eta$ to lie near $\bdy_R\overline{\Sigma_{dr}}$ and
  cross only $\alpha^R$-arcs. Then it is clear that $\gamma\cup\eta$
  is isotopic to $\pm K$.  With our orientation conventions, if we
  wish for $K$ to be oriented as the boundary of $\PunctF(\PMC_R)$,
  then we order $z_+$ and $z_-$ so that the arc from $z_+$ to $z_-$ in
  $\PMC(\drY)$ (with its induced orientation) contains the matched
  pairs for $\PMC_R$. (Equivalently, if we think of the arc $\z$ in
  $\HD$ connecting $\partial_L\HD$ to $\partial_R \HD$ as running left
  to right, then $\z_+$ is gotten by translating $\z$ upwards in
  $\HD$, and $\z_-$ is gotten by pushing it down.)
\end{proof}

\begin{definition}
  We can define a map $\spinc\co \Gen(\HD)\to \SpinC(Y)$
  as follows.  View $\x\in\Gen(\HD)$ as a generator for $\Gen(\drHD)$,
  let $\relspinc(\x)$ denote its corresponding relative $\SpinC$
  structure, and then let $\spinc(\x)$ be the corresponding $\SpinC$
  structure in $\SpinC(Y)$, according to the equivalence of
  Lemma~\ref{lem:RelAbsSpinC}. Let 
  \[
  \Gen(\HD,\spinc)=\{\x\in\Gen(\HD)\mid \spinc(\x)=\spinc\}.
  \]
\end{definition}

\begin{lemma}\label{lem:same-spinc}
  Given $\x,\y\in\Gen(\HD)$, the set
  $\pi_2(\x,\y)$ is non-empty if and only if $\spinc(\x)=\spinc(\y)$.
\end{lemma}

\begin{proof}
  In view of Lemma~\ref{lem:RelAbsSpinC}, this statement is equivalent
  to the corresponding statement for knot Floer homology (see
  \cite[Section 2.3]{OS04:Knots} and \cite[Sections 2.4 and
  2.6]{OS04:HolomorphicDisks}). 
  Recall that this is proved first by constructing a difference element
  $\epsilon(\x,\y)\in H_1(\drY\setminus C,\partial \drY)$ for $\x,\y\in\Gen(\HD)$
  (which vanishes if and only if $\pi_2(\x,\y)$ is non-empty), and then showing
  that $\relspinc(\y)=\relspinc(\x)+\PD(\epsilon(\x,\y))$.
\end{proof}

As noted in the discussion of the Alexander grading above, there is
a restriction map $\SpinC(Y)\to \SpinC(\bdy_RY)\cong\ZZ.$
We will see in the proof of Theorem~\ref{thm:DoublePairing} in
Section~\ref{sec:PairingTheorems} that this
restriction map is closely related to the strands grading on the
algebra.

\subsection{Gluing Heegaard diagrams}
\label{subsec:GluingDiagrams}
In Section~\ref{sec:PairingTheorems}, we will see how bordered Floer
homology groups transform under three gluing operations one can
perform on bordered three-manifolds.

The first of these gluing operations glues a bordered three-manifold
with one boundary component to one with two boundary components.  

\begin{construction}
Suppose Y is a bordered three-manifold with one boundary component
whose parameterization is specified by a homeomorphism $\phi\colon
F(\PMC) \to \partial Y$, and let $Y'$ be a strongly
bordered three-manifold with two boundary components $\partial_L Y'$
and $\partial_R Y'$ with parameterizations specified by
\[\xymatrix{
  \phi'_L\colon F(\PMC'_L) \to \partial_L Y' &
  \phi'_R\colon F(\PMC'_R) \to \partial_R Y'.}\]
Suppose moreover that $\PMC=-\PMC'_L$. Then we can form the bordered
three-manifold
$$Y\sos{\partial Y}\cup{\partial_L Y'} Y',$$
which is obtained by gluing $Y$ to $Y'$ via the identification
of $\partial Y$ with $-\partial_L Y'$ given by 
$\phi'_L\circ \phi^{-1}$. (Note that for the purpose of this definition, 
we do not need $Y'$ to be strongly bordered, just bordered; however,
in our applications $Y'$ will be equipped with this extra data.)
\end{construction}

The second of these gluing operations glues two strongly bordered
three-manifolds with two boundary components. 

\begin{construction}
Fix two strongly bordered three-manifolds with two
boundary components 
\begin{align*}
&(Y',\PMC_L',\Delta_L',z_L',\phi_L',\PMC_R',\Delta_R',z_R',\phi_R',
  \gamma_z',\nu_z')\text{ and}\\
&(Y'',\PMC_L'',\Delta_L'',z_L'',\phi_L'',\PMC_R'',\Delta_R'',z_R'',\phi_R'',
  \gamma_z'',\nu_z'').
\end{align*}
Suppose moreover that
$\PMC'_R=-\PMC''_L$. Then we can form a new strongly bordered
three-manifold with two boundary components. The underlying three-manifold
is
$$Y'''=Y'\sos{\partial_R Y'}{\cup}{\partial_L Y''} Y'',$$
gotten by gluing $\partial_R Y'$ to $\partial_L Y''$
via $\phi_L''\circ (\phi_R')^{-1}$.
The path $\gamma'''$ is gotten by connecting $\gamma_z'$ to $\gamma_z''$.
Framings are obtained similarly.
\end{construction}

The third gluing operation is a kind of self-gluing. 

\begin{construction}
\label{construct:GeneralizedOpenBooks}
Suppose that 
$$(Y,\PMC_L,\Delta_L,z_L,\phi_L,\PMC_R,\Delta_R,z_R,\phi_R,
  \gamma_z,\nu_z)$$
is a strongly bordered three-manifold with two boundary components;
and suppose moreover that $\PMC_L= -\PMC_R$.
Then, identifying the two boundary
components of $Y$ together, we obtain a new three-manifold which is
equipped with a framed knot, gotten by gluing up the framed arc
$\gamma_z$.  Performing
surgery on this framed knot in the self-glued three-manifold, we obtain a new
three-manifold denoted $(Y^{\circ},K)$, equipped with a knot $K$
gotten as the core of the surgery torus. 
We call $(Y^{\circ},K)$ the {\em generalized open book} associated
to the strongly bordered three-manifold $Y$.
\end{construction}
The justification for this terminology is the following.  If we
consider a strongly based mapping class $\phi\colon \PunctF(\PMC)
\to \PunctF(\PMC)$, there is an associated strongly
bordered three-manifold $M_\phi$ whose underlying topological space is
$[0,1]\times F(\PMC)$, parameterized by the identity on one boundary
and (the map on the closed surface induced by) $\phi$ on the
other; see Construction~\ref{construct:MCG-to-Bordered} and
Lemma~\ref{lem:open-book-is-open-book}. The
associated three-manifold $(Y^\circ,K)$ gotten as above is
classically known as the {\em open book} associated to~$\phi$.

On the level of Heegaard diagrams, the three gluing operations can be
described as follows.

Let $\HD=(\Sigma,\alphas,\betas,z)$ be a pointed bordered Heegaard
diagram with one boundary component, and
$\HD'=(\Sigma',\alphas',\betas',\mathbf{z'})$ an arced bordered Heegaard
diagram with two boundary components. Suppose that the pointed matched
circle $\PMC(\HD)$ associated to $\HD$ is the orientation reverse
$-\PMC_L(\HD')$ of the left pointed matched circle associated to
$\HD'$. Then gluing the boundary of $\HD$ to the left boundary of
$\HD'$ we obtain a new pointed bordered Heegaard diagram
\[
\HD\sos{\bdy}{\cup}{\bdy_L}\HD'=(\Sigma\sos{\bdy}{\cup}{\bdy_L}\Sigma',\alphas\sos{\bdy}{\cup}{\bdy_L}\alphas', \betas\cup\betas',\bdy_R\mathbf{z'}).
\]

Similarly, if $\HD''$ is another arced bordered Heegaard diagram with
two boundary components, such that $\PMC_L(\HD'')=-\PMC_R(\HD')$ then
we can glue $\HD''$ to $\HD'$ to get
\[
\HD'\sos{\bdy_R}{\cup}{\bdy_L}\HD''=(\Sigma'\sos{\bdy_R}{\cup}{\bdy_L}\Sigma'',\alphas'{}_{\bdy_R}\cup_{\bdy_L}\alphas'',
\betas'\cup\betas'',\mathbf{z'}\sos{\bdy_R}{\cup}{\bdy_L} \mathbf{z''}).
\]

\begin{lemma}
  \label{lem:GluingDiagrams}
With notation from above, gluing bordered Heegaard diagrams corresponds to gluing bordered
three-manifolds, as follows:
\begin{align*}
Y(\HD\sos{\bdy}{\cup}{\bdy_L}\HD')&=Y(\HD)\cup_{\phi'_L\circ\phi^{-1}}Y(\HD') \\
Y(\HD'\sos{\bdy}{\cup}{\bdy_L}\HD'')&=Y(\HD')\cup_{\phi''_L\circ(\phi'_R)^{-1}}Y(\HD'').
\end{align*}
\end{lemma}

\begin{proof}
  This is straightforward.
\end{proof}

Finally, we have the following construction  mirroring
Construction~\ref{construct:GeneralizedOpenBooks} on the level of 
Heegaard diagrams:
\begin{construction}
\label{construct:GenOBHeeg}
Suppose that $\HD$ is an arced bordered Heegaard diagram with two
boundary components such that $\PMC_L\HD=-\PMC_R\HD$. Then we can glue
$\bdy_L\HD$ to $\bdy_R\HD$. The result is a closed surface
$\Sigma^\Box$ of genus $g+1$ with $g$ $\alpha$-circles and $g$
$\beta$-circles, as well as a closed curve $\mathbf{z}^\circ$ corresponding to
the arc ${\mathbf z}$. Place basepoints $z_+$ and $z_-$ on the two
sides of $\mathbf{z}^\circ$, and then surger out the arc
$\mathbf{z}^\circ$ from $\Sigma^\Box$. The result
is a doubly-pointed Heegaard diagram
$\HD^\circ=(\Sigma^\circ,\alphas^\circ,\betas^\circ,z_+,z_-)$,
which we call the {\em self-glued diagram} associated to $\HD$.
\end{construction}

\begin{lemma}
  \label{lem:SelfGluingDiagrams}
  The self-glued diagram $\HD^{\circ}$ of
  Construction~\ref{construct:GenOBHeeg}  represents the generalized
  open book associated to the strongly bordered three-manifold
  $(Y^{\circ},K)$ of
  Construction~\ref{construct:GeneralizedOpenBooks}.
\end{lemma}

\begin{proof}
  This is straightforward; see also Lemma~\ref{lem:DrilledDiagram}.
\end{proof}

Next, we discuss how the gluing constructions interact with the
admissibility hypotheses.

\begin{lemma}
 \label{lem:admiss-glues}
 Let $\HD_1$ and $\HD_2$ be arced bordered Heegaard diagrams with two
 boundary components, such that $\PMC_R(\HD_1)=-\PMC_L(\HD_2)$. Let
 $\HD=\HD_1\sos{\bdy_R}{\cup}{\bdy_L}\HD_2$.  If $\HD_1$ is right
 admissible and $\HD_2$ is provincially admissible, or if $\HD_1$ is
 provincially admissible and $\HD_2$ is left admissible, then $\HD$ is
 provincially admissible. Moreover:
 \begin{enumerate}
 \item\label{item:Admis1} If $\HD_1$ and $\HD_2$ are both left
   admissible (respectively right admissible) then $\HD$ is left
   admissible (respectively right admissible).
 \item If $\HD_1$ (respectively $\HD_2$) is
   admissible then $\HD$ is left admissible (respectively right
   admissible).
 \item \label{item:Admis3}If $\HD_1$ (respectively $\HD_2$) is
   admissible and $\HD_2$ (respectively $\HD_1$) is right admissible
   (respectively left admissible) then $\HD$ is admissible.
 \end{enumerate}
 The obvious analogues hold in the case that $\HD_1$ has only one
 boundary component. (Many of the statements become the same in this
 case.)
\end{lemma}
\begin{proof}
  We will discuss the case that $\HD_1$ is right admissible and
  $\HD_2$ is provincially admissible, and the case that both $\HD_1$
  and $\HD_2$ are left admissible; the other cases are similar.

  Suppose $\HD_1$ is right admissible. If $P$ is a nontrivial
  provincial periodic domain in $\HD$ then $P\cap\HD_1$ is a
  left-provincial periodic domain in $\HD_1$. Hence either $P\cap
  \HD_1$ has both positive and negative coefficients or $P\cap \HD_1$
  is the trivial domain. In the latter case, $P\cap \HD_2$ is a
  provincial periodic domain, and hence has both positive and negative
  coefficients.  In either case, $P$ has both positive and negative
  coefficients.

  Similarly, suppose $\HD_1$ and $\HD_2$ are both left admissible. If
  $P$ is a nontrivial right-provincial periodic domain in $\HD$ then
  $P\cap \HD_2$ is a right-provincial periodic domain in $\HD_2$, and
  hence either has both positive and negative coefficients or is
  trivial. In the latter case, $ P\cap \HD_1$ is a nontrivial
  right-provincial periodic domain in $\HD_1$, and hence has both
  positive and negative coefficients.
\end{proof}

We say a doubly-pointed Heegaard diagram is {\em admissible} if all
periodic domains (i.e., domains which miss both basepoints) have both positive
and negative local multiplicities.  This is the condition required to
define knot Floer homology using the given Heegaard diagram. (It
corresponds to weakly admissibility for all $\SpinC$ structures for
singly-pointed Heegaard diagrams, in the sense of~\cite[Definition
4.10]{OS04:HolomorphicDisks}.)

\begin{lemma}\label{lem:dbl-bounded-self-glue} Let $\HD$ denote an
  arced bordered Heegaard diagram with two boundary components and
  such that $\PMC_L(\HD)=-\PMC_R(\HD)$. Suppose that $\HD$ is
  admissible.  Then the doubly-pointed Heegaard diagram $\HD^\circ$ is
  admissible.
\end{lemma}
\begin{proof}
  Note that the set of periodic domains in $\HD^\circ$ is a subset of
  the set of periodic domains in $\HD$. By
  Lemma~\ref{lemma:left-right-Area}, we can find an area form on
  $\widetilde{\Sigma}$ with respect to which any periodic domain has
  signed area zero. This induces an area form on $\Sigma^\circ$ with
  the corresponding property.
\end{proof}

\subsection{Bordered Heegaard diagrams for surface diffeomorphisms}
\label{sec:DiagramsForAutomorphisms}
\begin{definition}
  Given a bordered $3$-manifold $(Y,\psi)$ (where $\psi\co F(\PMC)\to \bdy
  Y$) and $\phi\in \MCG_0(\PMC,\PMC')$, let $\phi(Y,\psi)$ denote the bordered
  $3$-manifold $(Y,\psi\circ\phi^{-1})$, the result of twisting the
  parameterization of the boundary by $\phi$. This has a natural
  extension to strongly bordered three-manifolds as well:
  given a strongly bordered three-manifold $\sos{\psi_L}{Y}{\psi_R}$,
  where $\psi_R\co F(\PMC_R)\to\bdy_R Y$, and $\phi\in \MCG_0(\PMC_R,\PMC_R')$,
  let $\phi(\sos{\psi_L}{Y}{\psi_R})$ be the strongly bordered
  three-manifold $(\sos{\psi_L}{Y}{\psi_R\circ\phi^{-1}})$.
\end{definition}

When considering the above action of the mapping class group on bordered
three\hyp manifolds, the following strongly bordered 3-manifolds arise
naturally:

\begin{construction}
\label{construct:MCG-to-Bordered}
Fix pointed matched circles $\PMC_L$ and $\PMC_R$, and a strongly
based mapping class $\phi\colon (F(\PMC_L),D_L,z_L)\to
(F(\PMC_R),D_R,z_R)$. We can form a corresponding strongly bordered
three-manifold with $Y=[0,1]\times F(\PMC_R)$, 
$\partial_L Y = \{0\} \times -F(\PMC_R)$, $\partial_R Y = \{1\}\times F(\PMC_R)$,
$\Delta_L=\{0\}\times D_R$,
$\psi_L=\{0\}\times-\phi$, $\Delta_R=\{1\}\times D_R$,
$\psi_R=\{1\}\times\Id$,
$\gamma_z=[0,1]\times\{z_R\}$. We call it the \emph{mapping cylinder
  of~$\phi$}, and denote the resulting
arced bordered three-manifold by
$\sos{\phi}{{\left([0,1]\times F(\PMC_R)\right)}}{\Id}$ or
simply~$M_\phi$.  Note that $\bdy_L M_\phi = -F(\PMC_L)$.
\end{construction}

\begin{remark}
  In general, for two topological spaces $X$ and $Y$ and a map
  $\phi\co X \to Y$, the mapping cylinder of $\phi$ is the quotient space
  \[
  M_\phi = ([0,1] \times X \amalg Y) / ((1,x) \sim \phi(x)),
  \]
  equipped with maps
  \begin{align*}
    \psi_L\co X &\hookrightarrow M_\phi &
     \psi_R\co Y &\hookrightarrow M_\phi \\
    \psi_L(x) &= 0 \times x &
      \psi_R(y) &= y.
  \end{align*}
  This definition works for arbitrary maps $\phi$, and, for instance,
  gives a CW complex if $X$ and $Y$ are CW complexes and $\phi$ is a
  cellular map.  In the case when $\phi$ is a homeomorphism, the above
  space is equivalent to $[0,1] \times X$ with $\psi_L = \{0\}\times
  \Id_X$ and $\psi_R = \{1\}\times \phi^{-1}$, which in turn is
  equivalent to  $[0,1] \times Y$ with $\psi_L = \{0\}\times
  \phi$ and $\psi_R = \{1\}\times \Id_Y$.
\end{remark}

\begin{lemma}
  \label{lem:WellDefinedStronglyBordered}
  Any strongly bordered $3$-manifold~$Y$ whose underlying
  space can be identified with a product of a surface with an interval
  (so that $\gamma_z$ is identified with the product of a point with
  the interval, respecting the framing) is of the form
  $M_\phi$ for some choice of
  strongly based mapping class~$\phi$. Moreover, two such strongly
  bordered three-manifolds are isomorphic if and only if they
  represent the same strongly based mapping class.
\end{lemma}

\begin{proof}
  Suppose that $\sos{\psi_L}{Y}{\psi_R}$ is a strongly bordered three-manifold whose
  underlying space is homeomorphic to the product of a surface with an
  interval.  To keep orientations consistent, let $F(\PMC_L) =
  -\bdy_L(Y)$ and $F(\PMC_R) = \bdy_R(Y)$. We extract a strongly based
  mapping class
  $\phi\in\MCG_0(F(\PMC_L),F(\PMC_R))$ as follows.  Fix a
  diffeomorphism $\Phi\co [0,1]\times F(\PMC_R)\to Y$ so that the following hold:
  \begin{itemize}
  \item $\Phi|_{\{1\}\times F(\PMC_R)}=\psi_R$.
  \item $\Phi([0,1]\times z_R)=\gamma_z$.
  \item $\Phi(\{0\}\times D_R)=\Delta_L$.
  \item $\Phi(\{1\}\times D_R)=\Delta_R$.
  \item The normal vector $\nu_z$ is never tangent to $\Phi([0,1]
    \times \partial D_R)$.
  \end{itemize}
  We then define
  \[\phi=(\Phi|_{\{0\}\times F(\PMC_R)})^{-1} \circ (-\psi_L) \co
  F(\PMC_L(\HD))\to F(\PMC_R(\HD)).\]
  Then $\Phi$ provides an isomorphism between $Y$ and
  the strongly bordered three-manifold $\sos{\phi}{{\left([0,1]\times F(\PMC_R)\right)}}{\Id}$.
  
  Next, we claim that the strongly based mapping class of $\phi$ is
  independent of the choices made.  Indeed, if $\Phi'\co
  [0,1]\times F(\PMC_R)\to Y$ is an alternate choice of $\Phi$, then
  $\Phi^{-1}\circ\Phi'$ is a pseudoisotopy from
  $(\Phi|_{\{0\}\times F(\PMC_R)})^{-1}\circ\Phi'|_{\{0\}\times F(\PMC_R)}$
  to the identity map.
  (A \emph{pseudoisotopy} between two self-diffeomorphisms $f_0$ and $f_1$
  of a closed manifold $M$ is a self-diffeomorphism of
  $[0,1]\times M$ that restricts to $f_0$ on $\{0\}\times M$ and $f_1$
  on $\{1\} \times M$.)
  Since the equivalence relations induced by pseudoisotopy and isotopy agree in dimension
  $2$,\footnote{Proof: if a map $\phi$ is pseudoisotopic (concordant)
    to the identity then, in particular, $\phi$ is homotopic to the
    identity. But homotopic homeomorphisms of surfaces are isotopic.}  it
  follows that
  \[
  (\Phi|_{\{0\}\times F(\PMC_R)})^{-1} \circ (-\psi_L)  \text{\quad and\quad}
  (\Phi'|_{\{0\}\times F(\PMC_R)})^{-1} \circ  (-\psi_L)
  \]
  are isotopic, as desired.
\end{proof}

These three-manifolds encode the action of the mapping class group
on bordered three-manifolds, in the following sense.

\begin{lemma}
  \label{lem:BorderedForDiffeo}
  Fix a strongly based diffeomorphism
  $\phi\colon F(\PMC)\to F(\PMC')$. The associated
  strongly bordered three-manifold
  $M_\phi$ has the following properties:
  \begin{enumerate}
  \item
    \label{e:ActOnBordered}
    If $(Y,\psi)$ is a bordered three-manifold,
    where $\psi\co F(\PMC)\to \partial Y$ is
    a strongly based mapping class, then
    $\phi(Y,\psi)$ is obtained by gluing $Y$ and~$M_\phi$:
    $$\phi(Y,\psi) \cong Y \sos{\bdy Y}{\cup}{\bdy_L} M_\phi,$$
    canonically.
  \item 
    \label{e:ActOnStrongBordered}
    If $\sos{\psi_L}{Y}{\psi_R}$ is a strongly bordered
    three-manifold with two boundary components, where
    $\psi_R\colon F(\PMC_R)\to \partial_R Y$
    and $\PMC\cong \PMC_R$, then
    $\phi(\sos{\psi_L}{Y}{\psi_R})$ is obtained by 
    gluing $Y$ and $M_\phi$
    along $\partial_R Y$:
    $$\phi(\sos{\psi_L}{Y}{\psi_R}) \cong
    \left(\sos{\psi_L}{Y}{\psi_R}\right) \sos{\bdy_R}{\cup}{\bdy_L} 
  M_\phi,$$
  canonically.
  \item
    \label{e:Groupoid}
    Given another strongly based diffeomorphism $\phi'\colon
    F(\PMC')\to F(\PMC'')$, we have that
    $M_{\phi' \circ \phi}$
    is obtained 
    from gluing:
    \[
    M_{\phi' \circ \phi} \cong M_{\phi}
     \sos{\bdy_R}{\cup}{\bdy_L}
     M_{\phi'},
    \]
    canonically.
  \end{enumerate}
\end{lemma}

\begin{proof}
  It is straightforward to construct the isomorphism
  realizing Properties \eqref{e:ActOnBordered}
  and~\eqref{e:ActOnStrongBordered}.
  Property~\eqref{e:Groupoid} follows from
  Property~\eqref{e:ActOnStrongBordered}, as
  \begin{align*}
    \phi'(M_\phi) &= \phi'(\sos{\phi}{([0,1]\times F(\PMC'))}{\Id}) \\
      &= \sos{\phi}{([0,1]\times F(\PMC'))}{(\phi')^{-1}} \\
      &= \sos{\phi' \circ \phi}{([0,1]\times F(\PMC''))}{\Id} \\
      &= M_{\phi' \circ \phi}.\qedhere
  \end{align*}
\end{proof}

Similarly, for self gluing, we have:
\begin{lemma}\label{lem:open-book-is-open-book}
  Let $\phi\co F(\PMC)\to F(\PMC)$ be a strongly based
  diffeomorphism. Then the generalized open book associated to
  $M_\phi$ (Construction~\ref{construct:GeneralizedOpenBooks}) agrees
  with the open book associated to $\phi$ (with the orientation
  conventions from~\cite[Section 4.4.2]{Geiges08:ContactIntro}
  or~\cite[Section 2]{Etnyre06:OpenBookLecs}, say).
\end{lemma}
\begin{proof}
  The open book associated to
  $\phi$ is given by
  \[
  \Bigl([0,1]\times (F(\PMC)\setminus \Delta)\Bigr)\Big/\begin{pmatrix}
    (1,x)\sim(0,\phi(x))\\
  (t,x)\sim(t',x)\text{ for }x\in\bdy\Delta\end{pmatrix},
  \]
  which agrees with the conventions from
  Construction~\ref{construct:GeneralizedOpenBooks}.
\end{proof}

In the sequel, we will find it convenient to reformulate the above
properties in terms of Heegaard diagrams.

\begin{definition}
  \label{def:MCG-to-HD}
  Fix a strongly based mapping class $\phi$.  We say that a Heegaard
  diagram $\HD$ \emph{represents} $\phi$ if its underlying three-manifold
  $Y(\HD)$ is homeomorphic (respecting the marking) to the mapping
  cylinder $M_\phi$ of
  Construction~\ref{construct:MCG-to-Bordered}.
\end{definition}

\begin{lemma}\label{lemma:diagram-for-mcg-defined}
  If $\HD$ and $\HD'$ represent the same element
  $\phi\in\MCG_0(F(\PMC_L),F(\PMC_R))$ then $\HD$ and $\HD'$ are equivalent.
\end{lemma}
\begin{proof}
  By Proposition~\ref{prop:heegaard-moves}, it suffices to show that
  $Y(\HD)$ and $Y(\HD')$ are isomorphic strongly bordered
  $3$-manifolds. But this follows from
  Lemma~\ref{lem:WellDefinedStronglyBordered}.
\end{proof}

\begin{lemma}\label{lemma:compose-MCG-diagrams} 
  Fix strongly bordered mapping classes $\phi\in\MCG_0(F(\PMC),F(\PMC'))$ and 
  $\psi\in\MCG_0(F(\PMC'),F(\PMC''))$, and let $\HD_\phi$ and $\HD_\psi$
  be Heegaard diagrams representing $\phi$ and $\psi$ respectively.
  Then the union 
  $(\HD_\phi)\sos{\bdy_R}{\cup}{\bdy_L}(\HD_\psi)$
  is a Heegaard diagram which 
  represents the composite $\psi\circ\phi$.
\end{lemma}
\begin{proof}
  This follows from Lemma~\ref{lem:GluingDiagrams}
  and Part~\eqref{e:Groupoid} of Lemma~\ref{lem:BorderedForDiffeo}.
\end{proof}

It will be useful to have an explicit construction of a
Heegaard diagram associated to a strongly based mapping class $\phi$.

\begin{figure}
  \centering
  %Font is 12 point.
  \includegraphics{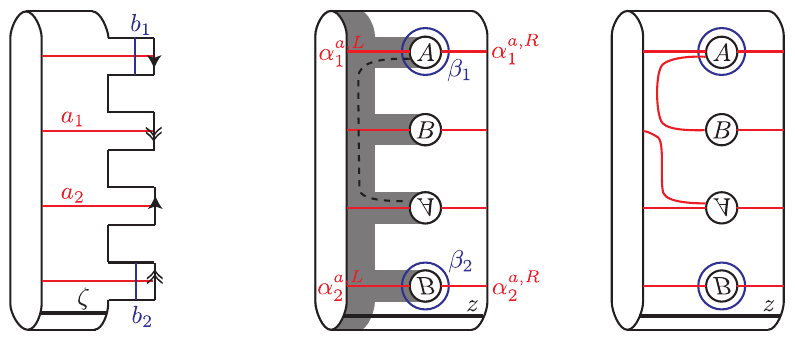}
  \caption{\textbf{Constructing a Heegaard diagram for a surface
      automorphism.} Left: the surface $F_0$ and the arcs $a_i$ and
    $b_i$ in it. Center: the surface $\Sigma$, and $\alpha$- and
    $\beta$-curves giving a diagram for the identity map. The
    subsurface $F_0$ of $\Sigma$ is shaded; note
    that $F_0\setminus \zeta$ is
    (orientation-preserving) homeomorphic to
    $\PunctF$. Right: the resulting diagram for a Dehn twist along the
    dashed curve in $\Sigma$.}
  \label{fig:Build-HD}
\end{figure}

Let $\PtdMatchCirc$ be a pointed matched circle.
Consider the product of the circle and an interval,
$[0,1]\times Z$.  Attach one-handles to $\{1\}\times Z$
as specified by the matching. For each
pair $p_i$ and $q_i$ on the pointed matched circle which are matched
(i.e., with $M(p_i)=q_i$)
we run an arc $a_i$ through the one-handle and extend $a_i$ as
$\{p_i,q_i\}\times [0,1]$ through the annulus $[0,1]\times Z$, so that its boundary
lies on $\{0\}\times \PtdMatchCirc$. This gives a
surface-with-boundary $F_0$ which is homeomorphic to the surface
$F(\PtdMatchCirc)$ with two disks removed, and which is
equipped with $2k$ arcs $\{a_1,\dots,a_{2k}\}$.  Equip $F_0$ with an additional $2k$ arcs
$\{b_1,\dots,b_{2k}\}$, which are chosen so that the arc $b_i$ is
contained in the $i\th$ one-handle attached to the original annulus,
and is dual to~$a_i$.  (That is, $b_i$ meets $a_i$ in a single, transverse
intersection point, and is disjoint from all the $a_j$ with $i\neq j$.)  The
boundary of $F_0$ has two components, one of which contains all the
endpoints of the $b_i$, which we denote $\partial_b F_0$, and the
other which contains all the endpoints of the $a_i$.  The basepoint in the
pointed matched circle equips $F_0$ with an arc $\zeta$ which connects
the two boundary components of $F_0$. See Figure~\ref{fig:Build-HD} on
the left.

Now, let ${\overline
F}_0$ be another copy of $F_0$ with orientation reversed,
equipped with curves $\{{\overline a}_i\}_{i=1}^{2k}$ and
$\{{\overline b}_i\}_{i=1}^{2k}$.
Let $\Sigma$ be the surface with two boundary components, obtained
from $F_0\amalg{\overline
F}_0$, by identifying $\partial_bF_0$ with $\partial_b {\overline
F}_0$ in such a manner that the boundary of $b_i$ is identified with
the boundary of ${\overline b}_i$, and the $\partial_bF_0$ boundary
point of $\zeta$ is identified with the corresponding boundary point
of ${\overline \zeta}$. The surface $\Sigma$ comes with $4k$
arcs $\alpha_i^{a,L}$ and $\alpha_i^{a,R}$ given by
$\alpha_i^{a,L} = a_i$ and
$\alpha_i^{a,R}={\overline a}_i$.  The surface
$\Sigma$ also comes with $2k$
circles $\beta_i=b_i\cup {\overline b}_i$ and one more arc  $\mathbf{z}=\zeta\cup{\overline \zeta}$. Then
$(\Sigma,\alpha_1^{i,L},\dots,\alpha_{2k}^{i,R},\beta_1,\dots,\beta_{2k},\mathbf{z})$
is a diagram for the identity map. Again, see
Figure~\ref{fig:Build-HD}, in the middle.

\begin{definition}
\label{def:ConstructHeegaardDiagram}
For $\phi$ a strongly based diffeomorphism $F(\PMC)
\to F(\PMC')$,
let $\Sigma$ be the surface obtained by gluing $F_0(\PMC')$ and ${\overline
  F}_0(\PMC)$ along $\bdy_b$, $\alpha_i^{a,L} = \phi(a_i)\subset F_0(\PMC')$,
$\alpha_i^{a,R}={\overline a}_i\subset \overline{F}_0(\PMC)$, $\beta_i=b_i\cup {\overline
  b}_i$, and $\mathbf{z}=\zeta\cup\overline{\zeta}$.
We call
$\HD(\phi)=(\Sigma,\alpha_1^{a,L},\dots,\alpha_{2k}^{a,L},\alpha_1^{a,R},\dots,\alpha_{2k}^{a,R},\beta_1,\dots,\beta_{2k},\mathbf{z})$
the \emph{canonical bordered Heegaard diagram associated to~$\phi$}.
\end{definition}
(Again, see Figure~\ref{fig:Build-HD}.) 

\begin{lemma}
  The canonical bordered Heegaard diagram $\HD(\phi)$ associated to $\phi$
  represents the map $\phi$ in the sense of
  Definition~\ref{def:MCG-to-HD}.
\end{lemma}
\begin{proof}
  We first verify the statement in the case where $\psi$ is the identity
  mapping class. Let $\HD$ be the corresponding Heegaard
  diagram. Clearly, the three-manifold $Y(\HD)$ is a surface times an
  interval. Thus, according to
  Lemma~\ref{lem:WellDefinedStronglyBordered}, $Y(\HD)$ represents
  some mapping class~$\psi$.  Moreover, after
  performing some handleslides and cancellations, one can see that
  $\HD\sos{\partial_R \HD}{\cup}{\partial_L \HD} \HD$ is equivalent to
  $\HD$. Thus, by Lemma~\ref{lemma:compose-MCG-diagrams},
  and~\ref{lem:WellDefinedStronglyBordered}, $\psi\circ\psi$
  represents the same mapping class as $\psi$. Thus, $\psi$ represents
  the identity mapping class.
  
  In general, the mapping cylinder $M_\phi$ is
  $\sos{\phi}{([0,1]\times F(\PMC'))}{\Id}$, which we can think of as
  obtained from the mapping cylinder for the identity by twisting the
  parametrization by $\phi$ on the left, which is just what we did
  above by changing the $\alpha^{a,R}$ arcs.
\end{proof}

%%% Local Variables: 
%%% mode: latex
%%% TeX-master: "Bimodules"
%%% End: 

\section{Bimodules for bordered manifolds}
\label{sec:CFBimodules}

Recall from~\cite{LOT1} that bordered Floer homology associates to a bordered $3$-manifold $Y$
with one boundary component modules $\CFDa(Y)$ and $\CFAa(Y)$.
The module $\CFDa(Y)$ is, more precisely, a type $D$ structure in the
sense of Definition~\ref{def:TypeD}, and 
encodes all the holomorphic curve counts in its differential.  The
module $\CFAa(Y)$ is an $\Ainf$ module, which encodes holomorphic
curve counts in all its actions.

For a $3$-manifold $Y$ with two boundary components, we can treat each
boundary component in either a type $A$ or a type $D$ manner. Treating
both as type $A$ boundaries leads to $\CFAAa(Y)$. Treating one as type
$A$ and the other as type $D$ leads to $\CFDAa(Y)$. Treating both as
type $D$ leads to $\CFDDa(Y)$.

It turns out that both $\CFAAa(Y)$ and $\CFDDa(Y)$ can be obtained
from the modules for $3$-manifolds with a single boundary component
via the drilling construction of Definition~\ref{def:Drilling}
and the restriction / induction functors of
Section~\ref{sec:induction-restriction}. Defining the module
$\CFDAa(Y)$ seems to require some new work, though the ideas (and
analytic machinery) are all present in the single boundary component
cases.

The reader is encouraged to consult Section~\ref{sec:torus-calc} for
examples of the bimodules. The reader may also want to refer
to \cite[Appendix~\ref*{LOT:app:Bimodules}]{LOT1} for an abbreviated (and, in the case of
$\CFAAa$ and $\CFDDa$, slightly different) account of this material.
\subsection{The type \textit{AA} bimodule}
\begin{definition}\label{def:AA-bimod}
  Let $Y$ be a strongly bordered $3$-manifold with $\bdy_LY=F(\PMC_L)$
  and $\bdy_RY=F(\PMC_R)$. Fix an arced bordered Heegaard diagram
  $\HD$ for $Y$, and assume that $\HD$ is provincially
  admissible. Then, define
  \[
  \CFAAa(\HD)=\Rest_{\PMC,\PMC'}(\CFAa(\drHD))
  \]
  which, in light of Section~\ref{sec:A-Bop-vs-bimods}, we can view as
  a bimodule with $\Ainf$-commuting right actions of
  $\AlgA{\PMC_L}$ and $\AlgA{\PMC_R}$.

  The module $\CFAAa(\HD)$ decomposes as a direct sum
  \[
  \CFAAa(\HD)=\bigoplus_{\spinc\in\SpinC(Y)}\CFAAa(Y,\spinc).
  \]
\end{definition}

Geometrically, Definition~\ref{def:AA-bimod} means that $\CFAAa(\HD)$
is generated over
$\FF_2$ by $g$-tuples of points in $\alphas\cap\betas$, one on each
$\alpha$- and $\beta$-circle and no two on the same $\alpha$-arc. The
differential counts provincial holomorphic curves, i.e., curves not
approaching $\bdy\overline{\Sigma}$. The right bimodule structure
come from counting curves with asymptotics at
$\bdy_L\overline{\Sigma}$ and $\bdy_R\overline{\Sigma}$, with
appropriate height constraints on the asymptotics. 

\begin{proposition}\label{prop:CFAA-inv}
  If $\HD$ and $\HD'$ are Heegaard diagrams for the same
  strongly bordered $3$-manifold $Y$ then $\CFAAa(\HD)_{\AlgA{\PMC_L},\AlgA{\PMC_R}}$ and
  $\CFAAa(\HD')_{\AlgA{\PMC_L},\AlgA{\PMC_R}}$ are $\Ainf$-homotopy equivalent bimodules. 
\end{proposition}
\begin{proof}
  As in the proof of Proposition~\ref{prop:heegaard-moves}, the
  drilled Heegaard diagrams $\drHD$ and $\drHD'$ are equivalent, so the
  result follows from invariance of $\CFAa$
  \cite[Theorem~\ref*{LOT:thm:A-invariance}]{LOT1}.
\end{proof}
Because of Proposition~\ref{prop:CFAA-inv}, we are justified in
writing $\CFAAa(Y)$ to denote the (homotopy equivalence class of)
$\CFAAa(\HD)$ for some (any) diagram $\HD$ for $Y$. 

\begin{lemma}\label{lem:admis-bd-AA}
  If $\HD$ is admissible (respectively left admissible, right
  admissible) in the sense of Definition~\ref{def:diagrams-multidef}
  (respectively Definition~\ref{def:admissible-diag}) then
  $\CFAAa(\HD)$ is bounded (respectively left bounded, right bounded)
  in the sense of Definition~\ref{def:AA-bounded}.
\end{lemma}
\begin{proof}
  This follows easily from Lemma~\ref{lemma:left-right-Area},
  similarly to \cite[Lemma~\ref*{LOT:lem:finite-typeA}]{LOT1}.
\end{proof}
\subsection{The type \DD\  bimodule}
\begin{definition}\label{def:DD-bimod}
  Let $Y$ be a strongly bordered $3$-manifold with $\bdy_LY=F(\PMC_L)$
  and $\bdy_RY=F(\PMC_R)$. Fix an arced bordered Heegaard diagram
  $\HD$ for $Y$, and assume that $\HD$ is provincially
  admissible. Then, define
  \[
  \CFDDa(\HD)=\Induct^{-\PMC_L,-\PMC_R}(\CFDa(\drHD)),
  \]
  a type \DD\ structure over $\AlgA{-\PMC_L}$ and
  $\AlgA{-\PMC_R}$ (where $-\PMC$ denotes $\PMC$ with its
  orientation reversed).

  The module $\CFDDa(\HD)$ decomposes as a direct sum
  \[
  \CFDDa(\HD)=\bigoplus_{\spinc\in\SpinC(Y)}\CFDDa(Y,\spinc).
  \]
\end{definition}

Geometrically, Definition~\ref{def:DD-bimod} means that $\CFDDa(\HD)$
is generated as a type \DD\ structure
over $\Alg(\PMC_L)\otimes\Alg(\PMC_R)$ by $g$-tuples of points in
$\alphas\cap\betas$, one on each $\alpha$- and $\beta$-circle and no
two on the same $\alpha$-arc. The differential counts holomorphic curves with
asymptotics at $\bdy\overline{\Sigma}$, without height constraints.
These curves contribute coefficients corresponding to their
asymptotics.

\begin{proposition}\label{prop:CFDD-inv}
  If $\HD$ and $\HD'$ are Heegaard diagrams for the same strongly
  bordered $3$-manifold $Y$ then
  ${}^{\AlgA{-\PMC_L},\AlgA{-\PMC_R}}\CFDDa(\HD)$ and
  ${}^{\AlgA{-\PMC_L},\AlgA{-\PMC_R}}\CFDDa(\HD')$ are homotopy
  equivalent bimodules.
\end{proposition}
\begin{proof}
  As in invariance of $\CFAAa$, the proof of
  Proposition~\ref{prop:heegaard-moves} implies that the drilled
  Heegaard diagrams $\drHD$ and $\drHD'$ are equivalent, so the result
  follows from invariance of $\CFDa$ \cite[Theorem~\ref*{LOT:thm:D-invariance}]{LOT1}.
\end{proof}
Because of Proposition~\ref{prop:CFDD-inv}, we are justified in
writing $\CFDDa(Y)$ to denote the (homotopy equivalence class of)
$\CFDDa(\HD)$ for some (any) diagram $\HD$ for $Y$. 

\begin{lemma}\label{lem:admis-bd-DD}
  If $\HD$ is admissible (respectively left admissible, right
  admissible) in the sense of Definition~\ref{def:diagrams-multidef}
  (respectively Definition~\ref{def:admissible-diag}) then
  $\CFDDa(\HD)$ is bounded (respectively left bounded, right bounded)
  in the sense of Definition~\ref{def:DD-bounded}.
\end{lemma}
\begin{proof}
  As for $\CFAAa$, this follows easily from
  Lemma~\ref{lemma:left-right-Area}.
\end{proof}
\subsection{The type \DA\  bimodule}
\label{subsec:TypeDAbimodule}
Fix a provincially admissible arced bordered Heegaard diagram
$\HD=(\overline\Sigma_g,\overline\alphas,\betas,\mathbf{z})$, with boundaries
$\PMC_L$ and $\PMC_R$ representing surfaces of genus $k_L$ and
$k_R$. We will associate to $\HD$ a bimodule over
$\AlgA{-\PMC_L}$ and $\AlgA{\PMC_R}$,
\[
\lsupv{\AlgA{-\PMC_L}}\CFDAa(\HD)_{\AlgA{\PMC_R}}.
\]
(As usual, $-\PMC_L$ denotes the orientation reverse of $\PMC_L$.)

Recall that $\SA{\HD}$ is the set of $g$-tuples of points $\x$ in
$\Sigma$ so that
\begin{itemize}
\item exactly one $x_i$ lies on each $\alpha$- and each $\beta$-circle and
\item no two $x_i$ lie on the same $\alpha$-arc.
\end{itemize}
Let $X(\HD)$ be the $\Field$--vector space spanned by
$\SA{\HD}$.

Given $\x\in\SA{\HD}$ let $o_L(\x)$ denote the indices of the $\alpha^L$-arcs
occupied by $\x$ and $o_R(\x)$ the indices of the $\alpha^R$-arcs occupied by
$\x$. We let $\SA{\HD,i}\subset \SA{\HD}$ denote the subset of
generators with $\#o_R(\x)=i$.  Note that $\#o_L(\x)+\#o_R(\x)=k_L+k_R$.
Let $I_{L,D}(\x)=I([2k_L]\setminus o^L(\x))$ and $I_{R,A}=I(o^R(\x))$.  We
define a left (respectively right) action of $\IdemA{\PMC_L}$
(respectively $\IdemA{\PMC_R}$) on $\SA{\HD}$ by
\[
I(\SetS)\cdot\x\cdot I(\SetT)\coloneqq
\begin{cases}
  \x & I(\SetS)=I_{L,D}(\x)\text{ and }I(\SetT)=I_{R,A}(\x)\\
  0 &\text{otherwise},
\end{cases}
\]
where $\SetS$ and $\SetT$ are subsets of $[2k]$.

As an $(\Idem(\PMC_L), \Idem(\PMC_R))$-bimodule, 
$\lsupv{\AlgA{-\PMC_L}}\CFDAa(\HD)_{\AlgA{\PMC_R}}$
 is $\SA{\HD}$, with the above action.

Our next task is to define the type \DA\ structure maps
on $\CFDAa(\HD)$, for which we resort to holomorphic
curves.

As in~\cite{LOT1}, we will count holomorphic curves in
\[
((\Sigma\setminus \mathbf{z})\times[0,1]\times\RR,(\alphas\times\{1\}\times\RR)\cup(\betas\times\{0\}\times\RR)).
\]
To avoid repeating the seemingly innumerable definitions and
propositions of \cite[Chapter~\ref*{LOT:chap:structure-moduli}]{LOT1},
we will use the drilling
construction of Section~\ref{sec:Diagrams} and simply use
moduli spaces defined in~\cite{LOT1}. (Since we are considering only curves
missing the region containing~$\mathbf{z}$, moduli spaces in the tunneled
diagram contain the moduli spaces in the original diagram.)

So, let
$\drHD=(\overline{\Sigma}_\dr,\overline{\alphas}_\dr,\betas_\dr,z_+)$
denote the bordered Heegaard diagram with one boundary component
obtained by drilling a tunnel from $\HD$. Reeb chords in
$(\bdy\overline{\Sigma}_\dr\setminus
z_+,\overline{\alphas}_\dr\cap\overline{\Sigma}_\dr)$
come in three kinds:
\begin{itemize}
\item Reeb chords connecting points in
  $\bdy\overline{\alphas}^{a,L}_\dr$. We refer to these as
  \emph{left Reeb chords}, and decorate them with an ``$L$''.
\item Reeb chords connecting points in
  $\bdy\overline{\alphas}^{a,R}_\dr$. We refer to these as
  \emph{right Reeb chords}, and decorate them with an ``$R$''.
\item Reeb chords connecting points in
  $\bdy\overline{\alphas}^{a,L}_\dr$ to points in
  $\bdy\overline{\alphas}^{a,R}_\dr$. We refer to these as
  \emph{mixed Reeb chords} and shall have no use for them in the
  present discussion.
\end{itemize}

Note that there is a one-to-one correspondence between $\SA{\HD}$ and $\S(\drHD)$.

Recall that a decorated source $\Source$ is a Riemann surface $S$ with
boundary and boundary punctures, where each puncture is either labelled
$+\infty$, $-\infty$ or $e\infty$, and the $e\infty$ punctures are
further labelled by Reeb chords at east infinity. Given generators $\x$
and $\y$, a homology class $B\in\pi_2(\x,\y)$ connecting $\x$ to $\y$,
a decorated source $\Source$, and an ordered partition $\vec{\rhos}$
of the Reeb chords labeling punctures of $\Source$ we have a moduli
space
\[
\cM^B(\x,\y;\Source;\vec{\rhos})
\]
of holomorphic curves $u\co
S\to\Sigma_\dr\times[0,1]\times\RR$ in the homology class $B$
with asymptotics specified by $\x$, $\y$ and $\vec{\rhos}$
\cite[Definition~\ref*{LOT:def:cM}]{LOT1}. The expected dimension of
$\cM^B(\x,\y;\Source;\vec{\rhos})$ is given by
\begin{equation}
g-\chi(S)+2e(B)+|\vec{\rhos}|-1=:\ind(B,\Source,\vec{\rhos})-1\label{eq:ind1}
\end{equation}
\cite[Proposition~\ref*{LOT:Prop:Index}]{LOT1}. If the curve $u$ is an embedding then
the Euler characteristic of $S$ is determined by
\begin{equation}
\chi(S)=\chi_{\emb}(B,\vec{\rhos}):=g+e(B)-n_\x(B)-n_\y(B)-\iota(\vec{\rhos})\label{eq:ind2}
\end{equation}
\cite[Proposition~\ref*{LOT:prop:asympt_gives_chi}]{LOT1}. In
particular, this leads us to
define
\begin{equation}
  \label{eq:ind3}
  \ind(B,\vec{\rhos})\coloneqq e(B)+n_\x(B)+n_\y(B)+\iota(\vec{\rhos})+\abs{\vec{\rhos}}.
\end{equation}

\begin{figure}
  \centering
  %Font is 12 point.
  \includegraphics{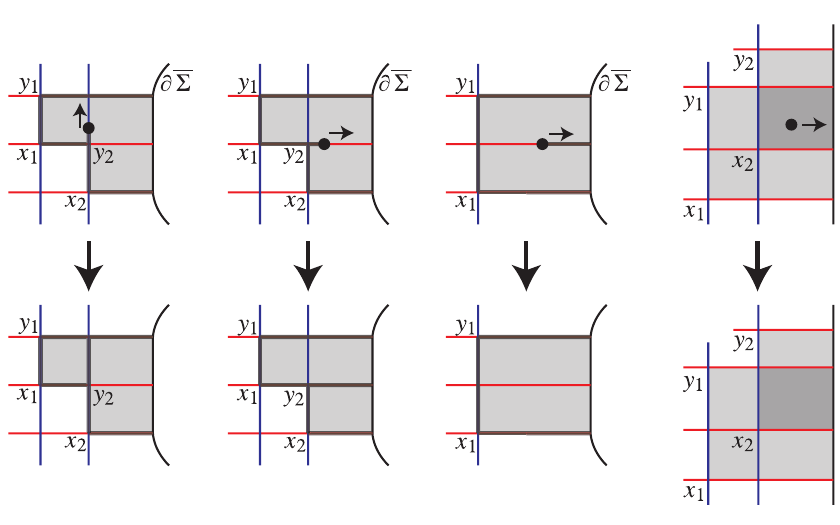}
  \caption{\textbf{Examples of the four codimension-one
      degenerations.} Far left: A two-level splitting. Center left:
    degenerating a join curve. Center right: degenerating a split
    curve (a collision of levels). Far right: degenerating a shuffle
    curve. The dark dots
    indicate branch points. This figure is
    adapted from \cite[Figures~\ref*{LOT:fig:degen_examples}
    and~\ref*{LOT:fig:shuffle2}]{LOT1}.}
  \label{fig:codim-one-degens}
\end{figure}

If the asymptotic data $(\x,\rhos)$ is such that
$u^{-1}(\widetilde\alpha_i^{a,L}\times(1,t))$ (respectively
$u^{-1}(\widetilde\alpha_i^{a,R}\times(1,t))$) consists of at most one
point for any given $t$ (i.e., $(\x,\rhos)$ is \emph{strongly boundary
  monotonic}) then the moduli space $\cM^B(\x,\y;\Source;\vec{\rhos})$
is well behaved, and we can understand its codimension-one boundary:
\begin{proposition} \cite[Theorem~\ref*{LOT:thm:master_equation}]{LOT1}
  \label{thm:master_equation}Suppose that $(\x,\vec{\rhos})$ is
  strongly boundary monotonic. Fix $\y$, $B\in\pi_2(\x,\y)$, and
  $\Source$ with $e\infty$ punctures labelled by $\vec{\rhos}$,
  such that $\ind(B,\Source,\vec\rhos)=2$. Let $\cM =
  \cM^B(\x,\y\semico\Source\semico\vec{\rhos})$.  Then the total number
  of all
  \begin{enumerate}
  \item two-story ends of~$\cM$,
  \item join curve ends of~$\cM$,
  \item odd shuffle curve ends of~$\cM$, and
  \item collision of levels $i$ and $i+1$ in $\cM$, where $\rhos_i$ and
    $\rhos_{i+1}$ are weakly composable
  \end{enumerate}
  is even.
\end{proposition}
Examples of the four types of degenerations are shown in
Figure~\ref{fig:codim-one-degens}.

To define the multiplications on $\CFDAa(\HD)$ we collect certain of
the $\cM(\x,\y;\Source;\vec{\rhos})$. Specifically, given a sequence
$\vec{\rho}^L$ of left Reeb chords and a sequence of sets of right
Reeb chords $(\rhos_1^R,\dots,\rhos_n^R)$, we say that a sequence
$\vec{\rhos}=(\rhos_1,\dots,\rhos_m)$ \emph{interleaves}
$(\vec{\rho}^L;\rhos_1^R,\dots,\rhos_n^R)$ if, as a multiset,
$\{\rhos_1,\dots,\rhos_m\}=\vec{\rho}^L\amalg\{\rhos_1^R,\dots,\rhos_n^R\}$,
and the orderings of $\vec{\rho}^L$ and $(\rhos_1,\dots,\rhos_n)$
agree with the orderings induced by $\vec{\rhos}$.

Now, let
\[
\Mod^B(\x,\y;\vec{\rho}^L;\rhos_1^R,\dots,\rhos_n^R)
=\quad\quad\bigcup_{\mathclap{\substack{\vec{\rhos}\text{ interleaves
    }(\vec{\rho}^L;\rhos_1^R,\dots,\rhos_n^R)\\
    \chi(S)=\chi_{\emb}(B,\vec{\rhos})}}}\quad\quad\cM^B(\x,\y;\Source;\vec{\rhos}).
\]

\begin{lemma}
  If $\vec{\rhos}$ and $\vec{\rhos}'$ both interleave
  $(\vec{\rho}^L;\rhos_1^R,\dots,\rhos_n^R)$ then $(\x,\vec{\rhos})$
  is strongly boundary monotonic if and only if $(\x,\vec{\rhos}')$
  is strongly boundary monotonic. Moreover, for any homology
  class $B$, $\ind(B,\vec{\rhos})=\ind(B,\vec{\rhos}')$.
\end{lemma}
\begin{proof}
  The boundary monotonicity statement is immediate from the
  definition.  It is also immediate from the definitions that if
  $\vec{\rhos}$ interleaves $(\vec{\rho}^L;\rhos_1^R,\dots,\rhos_n^R)$
  then
  \[
  \iota(\vec{\rhos})=\iota(\vec{\rho}^L)+\iota(\rhos_1^R,\dots,\rhos_n^R).
  \]
  So, the statement about $\ind$ follows from the definition,
  Formula~\eqref{eq:ind3}.
\end{proof}
Consequently, it makes sense to talk about a triple
$(\x;\vec{\rho}^L;\rhos_1^R,\dots,\rhos_n^R)$ being strongly boundary
monotonic: such a triple is strongly boundary monotonic if for some
(equivalently, any) sequence $\vec{\rhos}$ interleaving
$(\vec{\rho}^L;\rhos_1^R,\dots,\rhos_n^R)$, $(\x,\vec{\rhos})$ is
strongly boundary monotonic. Similarly, define
$\ind(B;\vec{\rho}^L;\rhos_1^R,\dots,\rhos_n^R)$ to be
$\ind(B,\vec{\rhos})$ for any $\vec{\rhos}$ which interleaves
$(\vec{\rho}^L;\rhos_1^R,\dots,\rhos_n^R)$.

With these moduli spaces in hand, define a type \DA\ structure
(Definition~\ref{def:DA-structure}) on $\CFDAa$ by
\begin{multline}
\label{eq:DefineTypeDA}
\delta^1_{n+1}(\x,a(\rhos_1^R),\dots,a(\rhos_n^R))  \\
\coloneqq
  \sum_{\y\in\SA{\HD}}\quad\quad
  \sum_{\mathclap{\substack{B\in\pi_2(\x,\y)\\
       \ind(B;\vec{\rho}^L;\rhos_1^R,\dots,\rhos_n^R)= 1}}}\quad
    \#\left(\Mod^B(\x,\y;\vec{\rho}^L;\rhos_1^R,\dots,\rhos_n^R)\right)a(-\rho_1^L)\cdots
    a(-\rho_m^L)\y 
\end{multline}
where $\vec{\rho}^L=(\rho_1^L,\dots,\rho_m^L)$ and $-\rho_i$ denotes
$\rho_i$ with its orientation reversed. The reader may find it helpful
to compare this definition with
\cite[Chapter~\ref*{LOT:chap:type-a-mod}]{LOT1}. Also, note
that the case $n=0$ is essentially the differential on $\CFDa$
from \cite[Chapter~\ref*{LOT:chap:type-d-mod}]{LOT1}. Finally, notice that if
$(\x;\vec{\rho}^L;\rhos_1^R,\dots,\rhos_n^R)$ is not strongly boundary
monotonic then either the left hand side of
Equation~(\ref{eq:DefineTypeDA}) is nonsensical or the element
$a(-\rho_1^L)\cdots a(-\rho_m^L)\y$ is automatically $0$.
\begin{lemma}
  Under the provincial admissibility hypothesis, the sum defining
  $\delta_n$ is finite. 
\end{lemma}
\begin{proof}
  This is a trivial adaptation of \cite[Lemma~\ref*{LOT:lem:finite-typeA}]{LOT1}.
\end{proof}

This completes the definition of $\CFDAa(\HD)$. It remains to check
that: 
\begin{itemize}
\item  The maps $\delta_n$ satisfy the compatibility conditions of
  Definition~\ref{def:DA-structure} (see Equation~\eqref{eq:DA-def}).
\item If $\HD$ and $\HD'$ define the same strongly bordered
  $3$-manifold then $\CFDAa(\HD)$ is $\Ainf$-homotopy-equivalent to $\CFDAa(\HD')$.
\end{itemize}

We start by refining Proposition~\ref{thm:master_equation}.
\begin{proposition}\label{prop:DA-master}
  Fix generators $\x$ and $\y$, $B\in\pi_2(\x,\y)$, a sequence of left
  Reeb chords~$\vec{\rho}^L$ and a sequence of sets of right Reeb
  chords $(\rhos_1^R,\dots,\rhos_n^R)$, so that
  $(\x;\vec{\rho}^L;\rhos_1^R,\dots,\rhos_n^R)$ is strongly boundary
  monotonic. Assume that
  $\ind(B;\vec{\rho}^L;\rhos_1^R,\dots,\rhos_n^R)=2$. Then the sum of
  the following numbers is even:
  \begin{enumerate}
  \item\label{item:2Story} The number of two-story ends of
    $\cM^B(\x,\y;\vec{\rho}^L;\rhos_1^R,\dots,\rhos_n^R)$, i.e.,
    \[
    \sum_{\mathclap{\substack{\ind(B_1;\vec{\rho}^L_1;\rhos^R_1,\dots,\rhos^R_i)=1\\\ind(B_2;\vec{\rho}^L_2;\rhos^R_{i+1},\dots,\rhos^R_n)=1}}}\quad\#\left(\cM^{B_1}(\x,\w;\vec{\rho}^L_1;\rhos^R_1,\dots,\rhos^R_i)\times\cM^{B_2}(\w,\y;\vec{\rho}^R_2;\rhos^R_{i+1},\dots,\rhos^R_{n})\right).
    \]
    where the sum is over $\w\in \SA{\HD}$, $B_1\in\pi_2(\x,\w)$,
    $B_2\in\pi_2(\w,\y)$, $B=B_1\ast B_2$, $i=0,\dots, n$ and
    $(\vec{\rho}^L_1,\vec{\rho}^L_2)=\vec{\rho}^L$.
  \item\label{item:RJoins} The number of join curve ends among right
    Reeb chords, i.e.,
    \[
    \sum_{\mathclap{\substack{i=1,\dots,n\\ \rho_{i,j}=\rho_a\uplus\rho_b}}}\#\cM^B(\x,\y;\vec{\rho}^L;\rhos_1^R,\dots,\rhos_i^{R,a,b},\dots,\rhos^R_n)
  \]
  where $\rhos_i^{R,a,b}$ is obtained from $\rhos_i^R$ by replacing
  $\rho_{i,j}\in\rhos_i^R$ by $\rho_a,\rho_b$.
\item\label{item:Shuffles} The number of odd shuffle curve ends among
  right Reeb chords, i.e.,
  \[
  \sum_{i=1}^n\#\cM^B(\x,\y;\vec{\rho}^L;\rhos_1^R,\dots,\rhos_i^{R,\prime},\dots,\rhos_n^R)
  \]
  where $\rhos_i^{R,\prime}$ is obtained from $\rhos_i^R$ by
  performing a weak shuffle.
\item\label{item:RCollisions} The number of collisions among right
  levels, i.e.,
  \[
  \sum_{i=1}^n\#\cM^B(\x,\y;\vec{\rho}^L;\rhos_1^R,\dots,\rhos_i^R\uplus\rhos_{i+1}^R,\dots,\rhos_n^R)
  \]
  where $\rhos_i$ and $\rhos_{i+1}$ are weakly composable.
\item\label{item:LJoins} The number of join curve ends among left Reeb
  chords, i.e.,
  \[
  \sum_{i=1}^n\#\cM^B(\x,\y;\vec{\rhos}^{L,\prime};\rhos^R_1,\dots,\rhos^R_n)
  \]
  where
  $\vec{\rhos}^{L,\prime}=(\rho_1^L,\dots,\rho_{i-1}^L,\{\rho_a^L,\rho_b^L\},\dots,\rho_m^L)$
  is obtained by replacing $\rho_i^L=\rho_a^L\uplus\rho_b^L$ in
  $\vec{\rho}^L$ with $\{\rho_a^L,\rho_b^L\}$.
\item\label{item:LSplits} The number of split curve ends among left
  Reeb chords, i.e.,
  \[
  \sum_{\rho_i^{L,+}=\rho_{i+1}^{L,-}}\#\cM^B(\x,\y;(\rho^L_1,\dots,\rho^L_{i-1},\rho^L_i\uplus\rho^L_{i+1},\dots,\rho^L_m);\rhos_1^R,\dots,\rhos_n^R).
  \]
\item\label{item:LCollisions} The number of other collisions of left levels
  $\rho_i,\rho_{i+1}$, i.e.,
  \[
  \sum\#\cM^B(\x,\y;\vec{\rhos}^{L,\prime};\rhos_1^R,\dots,\rhos_n^R)
  \]
  where
  $\vec{\rhos}^{L,\prime}=(\rho_1^L,\dots,\rho_{i-1}^L,\{\rho_i^L,\rho^L_{i+1}\},\dots,\rho_n^L)$
  and $\rho_i^{L,+}\neq\rho_{i+1}^{L,-}$. Moreover, $\rho_i^L$ and
  $\rho_{i+1}^L$ must satisfy: 
  \begin{itemize}
  \item $\rho_i^{L,+}$ and $\rho_{i+1}^{L,-}$ do not lie on the same
    $\alpha$-arc, and
  \item $(\rho_i^L,\rho_{i+1}^L)$ are not interleaved (in that order).
  \end{itemize}
\end{enumerate}
\end{proposition}
\begin{proof}
  Recall that embedded curves have maximal index and, in codimension one,
  families of embedded curves converge to embedded curves
  \cite[Proposition~\ref*{LOT:prop:asympt_gives_chi}, Lemma~\ref*{LOT:lem:splittings-embedded}
  and Lemma~\ref*{LOT:lemma:collision-is-composable}]{LOT1}.
  So, summing
  Proposition~\ref{thm:master_equation} over all $\Source$ with
  embedded Euler characteristic, the only possible ends not accounted
  for are:
  \begin{itemize}
  \item Collisions of levels between right and left Reeb chords. These
    cancel in pairs.
  \item Collisions of right levels which are not composable. These are
    prohibited
    by~\cite[Lemma~\ref*{LOT:lemma:collision-is-composable}]{LOT1}.
  \item Collisions of left levels not satisfying the conditions set
    out. The first condition comes from the fact that boundary
    degenerations are prohibited
    (see~\cite[Lemma~\ref*{LOT:lemma:collision-weakly-composable}]{LOT1}).
    (Note that collisions where $\rho_i^{L,+} = \rho_{i+1}^{L,-}$ are
    included in sum~\eqref{item:LSplits}.)
    The
    second comes from the fact that $\{\rho_i^L\}$ and $\{\rho_{i+1}^L\}$ must
    be composable~\cite[Lemma~\ref*{LOT:lemma:collision-is-composable}]{LOT1}.
  \item Shuffle curve ends among left Reeb chords. These are
    prohibited because each part of $\vec{\rho}^L$ has only a single
    Reeb chord.
  \end{itemize}
  The result follows.
\end{proof}

\begin{proposition}
  \label{prop:DA-module}
  The maps $\delta_n$ satisfy the compatibility conditions of a type
  \DA\  structure.
\end{proposition}
\begin{proof}
  The proof is a combination of the proofs of
  \cite[Propositions~\ref*{LOT:prop:typeD-d2}
  and~\ref*{LOT:prop:A-module-defined}]{LOT1}, and we shall be
  somewhat terse.  We must show that
  for any $\rhos_1,\dots,\rhos_n$,
\begin{align*}
  0&=(\bdy\otimes \Id_N)(\delta_{n+1} (\x\otimes a(\rhos_1)\otimes \dots\otimes a(\rhos_n))) \\
  &\quad +
  \!\!\!\sum_{i+j=n+2}\!\! (\mu_2\otimes \Id_N) \circ (\Id_A\otimes \delta_i) \left( 
  \delta_j(\x\otimes a(\rhos_1)\otimes \dots\otimes a(\rhos_{j-1}))\otimes a(\rhos_j)\otimes\dots\otimes a(\rhos_n)\right) \\
  &\quad +\sum_{i=1}^n\delta_{n+1}(\x\otimes a(\rhos_1)\otimes \dots\otimes \bdy a(\rhos_i)\otimes \dots\otimes \rhos_n)\\
  &\quad+\sum_{i=1}^n\delta_n(\x\otimes a(\rhos_1)\otimes \dots\otimes a(\rhos_i)a(\rhos_{i+1})\otimes \dots\otimes \rhos_n),
\end{align*}
cf.\ Definition~\ref{def:DA-structure}.

The second term corresponds to two level splittings,
sum~(\ref{item:2Story}) of Proposition~\ref{prop:DA-master}. The third
term corresponds to the right join and shuffle ends,
sums~(\ref{item:RJoins}) and~(\ref{item:Shuffles}) of
Proposition~\ref{prop:DA-master}; see also~\cite[Proof of
Proposition~\ref*{LOT:prop:A-module-defined}]{LOT1}. The fourth term
corresponds to the collisions of right
levels, sum~(\ref{item:RCollisions}) of
Proposition~\ref{prop:DA-master}; again, see also~\cite[Proof of
Proposition~\ref*{LOT:prop:A-module-defined}]{LOT1}.  The first term
corresponds to sum~(\ref{item:LSplits}) of
Proposition~\ref{prop:DA-master};
compare~\cite[Lemma~\ref*{LOT:lemma:third_term}]{LOT1}.

It remains to see that the sums~(\ref{item:LJoins})
and~(\ref{item:LCollisions}) of Proposition~\ref{prop:DA-master} cancel
in pairs as long as $a(-\rho_1^L)\cdots a(-\rho_m^L)\neq 0$. (See
also~\cite[Proof of Proposition~\ref*{LOT:prop:typeD-d2}]{LOT1} for
this part of the
proof.)  This product being nonzero imposes the following additional
conditions on collisions of left levels $\rho_i$ and $\rho_{i+1}$:
\begin{itemize}
\item $\rho_i^{L,-}$ and $\rho_{i+1}^{L,-}$ lie on different
  $\alpha$-arcs. Similarly, $\rho_i^{L,+}$ and $\rho_{i+1}^{L,+}$
  lie on different $\alpha$-arcs. See
  \cite[Lemma~\ref*{LOT:lemma:nonzero-implies-distinct}]{LOT1}.
\item If $\rho_i^{L,-}$ and $\rho_{i+1}^{L,+}$ lie on the same
  $\alpha$-arc then $\rho_i^{L,-}=\rho_{i+1}^{L,+}$. This is immediate
  from $a(-\rho_{i}^L)a(-\rho_{i+1}^L)\neq 0$.
\item $(\rho_{i+1}^L,\rho_i^L)$ are not interleaved (in that
  order). Again, this is immediate from the fact that
  $a(-\rho_{i}^L)a(-\rho_{i+1}^L)\neq 0$.
\end{itemize}
Thus, the two allowed kinds of left collisions which are not
algebraically $0$ are:
\begin{itemize}
\item Collisions with $\rho_i^{L,-}=\rho_{i+1}^{L,+}$. These moduli
  spaces cancel with the join curve ends of the factorization
  $a(-\rho_1^L)\cdots a(-(\rho_i^L\uplus\rho_{i+1}^{L}))\cdots
  a(-\rho_m^L).$
\item Collisions with the endpoints of $\rho_i^L$ and $\rho_{i+1}^L$
  lying on four different $\alpha$-arcs, and with $\rho_i^L$ and
  $\rho_{i+1}^L$ either nested or disjoint. In
  this case, the same degeneration also occurs
  for the factorization with $a(\rho_i^L)$ and $a(\rho_{i+1}^L)$
  switched. \qedhere
\end{itemize}
\end{proof}

We next turn to the issue of invariance.
\begin{proposition}\label{prop:CFDA-inv}
  If $\HD$ and $\HD'$ are provincially admissible arced bordered
  Heegaard diagrams defining the same strongly bordered $3$-manifold
  ${}_{\PMC_L}Y_{\PMC_R}$, then the corresponding bimodules
  $\lsupv{\AlgA{-\PMC_L}}\CFDAa(\HD)_{\AlgA{\PMC_R}}$ and
  $\lsupv{\AlgA{-\PMC_L}}\CFDAa(\HD')_{\AlgA{\PMC_R}}$ are
  $\Ainf$-homotopy equivalent.
\end{proposition}
\begin{proof}
  As in the case of a single boundary component, the invariance is
  proved by constructing homotopy equivalences corresponding to each
  of the Heegaard moves of Proposition~\ref{prop:heegaard-moves}; the
  reader should have no difficulty adapting the proof
  from \cite[Section~\ref*{LOT:sec:def-CFA}]{LOT1} to the present situation.
\end{proof}
Because of Proposition~\ref{prop:CFDA-inv}, we are justified in
writing $\CFDAa(Y)$ to denote the (homotopy equivalence class of)
$\CFDAa(\HD)$ for some (any) diagram $\HD$ for $Y$. 

We conclude this section with a lemma about admissibility.
\begin{lemma}\label{lem:admis-bd-DA}
  If $\HD$ is admissible (respectively left admissible, right
  admissible) in the sense of Definition~\ref{def:diagrams-multidef}
  (respectively Definition~\ref{def:admissible-diag}) then
  $\CFDAa(\HD)$ is bounded (respectively left bounded, right bounded)
  in the sense of Definition~\ref{def:DA-bounded}.
\end{lemma}
\begin{proof}
  As for $\CFAAa$ and $\CFDDa$, this follows easily from
  Lemma~\ref{lemma:left-right-Area}.  
\end{proof}

\begin{remark}
  It follows from the pairing theorems of
  Section~\ref{sec:PairingTheorems} that the module $\CFDAa(Y)$
  is determined by $\CFDDa(Y)$ (or, equally well,
  $\CFAAa(Y)$). Consequently, the bimodules associated to
  $3$-manifolds with two boundary components are completely determined
  by the invariants of $3$-manifolds with connected boundaries, via
  the induction / restriction functors.
\end{remark}

\subsection{Modules associated to surface automorphisms}
\label{sec:AutBimodules}
Given a strongly based diffeomorphism $\psi\colon
(F,D,z)\to (F,D,z)$, define
\begin{align*} 
  \CFAAa(\psi)_{\AlgA{F},\AlgA{F}}&:=\CFAAa(\HD)_{\AlgA{F},\AlgA{F}}\\
  {}^{\AlgA{-F},\AlgA{-F}}\CFDDa(\psi)&:= {}^{\AlgA{-F},\AlgA{-F}}\CFDDa(\HD)\\
  {}^{\AlgA{-F}}\CFDAa(\psi)_{\AlgA{F}}&:= {}^{\AlgA{-F}}\CFDAa(\HD)_{\AlgA{F}}
\end{align*}
where $\HD$ is any Heegaard diagram representing $\psi$ (in the sense
of Definition~\ref{def:MCG-to-HD}).

\begin{proof}[Proof of Theorem~\ref{thm:IsotopyInvariance}]
This is immediate from Proposition~\ref{lemma:diagram-for-mcg-defined}
and invariance of $\CFAAa$, $\CFDDa$ and $\CFDAa$,
Propositions~\ref{prop:CFAA-inv},~\ref{prop:CFDD-inv}
and~\ref{prop:CFDA-inv} respectively. 
\end{proof}

\subsection{Gradings}\label{sec:cf-gradings}

Suppose $Y$ is a strongly bordered three-manifold with boundary
parameterized by $\PMC_L$ and $\PMC_R$, and choose a compatible
provincially admissible arced bordered Heegaard diagram $\HD$.

The gradings on $\CFDDa(Y)$ and $\CFAAa(Y)$ are induced by the
induction and restriction functors, so we will focus on the grading of
$\CFDAa(Y)=\lsupv{\AlgA{-\PMC_L}}\CFDAa(Y)_{\AlgA{\PMC_R}}$.  (See
Remark~\ref{rem:DD-AA-grading} for an alternate approach to the
gradings on $\CFDDa(Y)$ and $\CFAAa(Y)$.)

The bimodule $\CFDAa(Y)$ interacts with the strands grading in the
following way.  The condition that $\#o_L(\x)+\#o_R(\x)=k_L+k_R$ for all
generators $\x$ gives
$$\IdemAS{-\PMC_L}{i}\cdot \CFDAa(Y)=\CFDAa(Y)\cdot
\IdemAS{\PMC_R}{i}.$$
Thus, defining
$\CFDAa(Y,i)\coloneqq\CFDAa(Y)\cdot \IdemAS{\PMC_R}{i}$, 
we have that $\CFDAa(Y,i)$ is a type \DA\ structure
over $\AlgA{-\PMC_L,i}$ and $\AlgA{\PMC_R,i}$, and
\begin{equation}
\label{eq:StrandSplittingDA}
\lsupv{\Alg(-\PMC_L)}\CFDAa(\HD)_{\Alg(\PMC_R)} =
\bigoplus_{i\in\ZZ}\lsupv{\Alg(-\PMC_L,i)}\CFDAa(\HD,i)_{\Alg(\PMC_R,i)}.
\end{equation}

Moreover, by Lemma~\ref{lem:same-spinc},
there is a natural splitting of $\CFDAa(Y)$ according to
$\SpinC$ structures:
$\CFDAa(Y)=\bigoplus_{\spinc\in\SpinC(Y)}\CFDAa(Y,\spinc)$.
In particular, to define the summand $\CFDAa(Y,\spinc)$, we repeat
the construction from Section~\ref{subsec:TypeDAbimodule}, 
using only the subset $\SA{\HD,\spinc}\subset \SA{\HD}$ of generators
representing $\spinc$.

We would like to endow $\CFDAa(Y)$ with the structure of a left-right
$(\bigGroup(-\PMC_L),\allowbreak\bigGroup(\PMC_R))$-set graded bimodule
(in the sense of Definition~\ref{def:GradingBimodules}); i.e.,
we would like to grade $\CFDAa(Y)$ by a set with a compatible right action by
\[
\bigGroup_{\DA}(\bdy \HD)\coloneqq
\bigGroup(-\PMC_L)^{\op}\times_{\lambda}
\bigGroup(\PMC_R).
\]
(When the Heegaard diagram is clear from the context, we will
sometimes write $\bigGroup_{\DA}$ to mean $\bigGroup_{\DA}(\bdy\HD)$.)
This is done by a suitable adaptation of the grading on $\CFDa$ and
$\CFAa$ from \cite{LOT1}.  

We will often work in the isomorphic group
\[
\bigGroup_\AAm(\bdy\HD) = \bigGroup(\PMC_L)\times_{\lambda}
\bigGroup(\PMC_R).
\]
Recall from Equation~\eqref{eq:def-R}
that, if $r: Z
\to -Z$ is the (orientation\hyp reversing) identity map, then
\[
R(j,\alpha) = (j, r_*(\alpha))
\]
defines a group anti-homomorphism from $G(\PMC)$ to $G(-\PMC)$, and
so an isomorphism $G(\PMC) \cong G(-\PMC)^\op$.  Then
\[
R \times_\lambda \Id \co \bigGroup_\AAm(\bdy\HD) \to \bigGroup_{DA}(\bdy\HD)
\]
is a canonical isomorphism, which we denote~$\tR$.

We will construct the grading one $\SpinC$ structure at a
time. Specifically, fix a Heegaard diagram $\HD$ for $Y$ and a
$\SpinC$ structure $\spinc$ over $Y$. Suppose that there is at least
one generator for $\Gen(\HD)$ which represents~$\spinc$; we will
return to the case that $\HD$ has no generator representing~$\spinc$
at the end of this subsection.

There is a map
$$g'\co \pi_2(\x,\y)\to \bigGroup_\AAm(\partial\HD)$$
defined by
\begin{equation}
  \label{eq:DefineGPrimed}
  g'(B)= (-e(B)-n_{\x}(B)-n_{\y}(B),\partial^{\partial_L} B,
\partial^{\partial_R} B),
\end{equation}
where $e(B)$ is the Euler measure of $B$
and $n_{\x}(B)$ is sum of the average local multiplicities of $B$ at
each coordinate of~$\x$. (Compare
\cite[Section~\ref*{LOT:sec:domains}]{LOT1}.)

Recall from Section~\ref{sec:Gradings} that $\bigGroup(\PMC)$ is an
index two subgroup of the group $\OneHalf \ZZ\times H_1(Z',{\mathbf
  a})$ (with a twisted multiplication), 
so we must show that
$g'(B)\in\bigGroup_\AAm(\bdy\HD)$. To this end, we have the following:
\begin{lemma}\label{lem:DA-grading-integrality}
  The tuple $g'(B)$ defined in Equation~\eqref{eq:DefineGPrimed} is an
  element of $\bigGroup_\AAm(\partial\HD)$.
\end{lemma}
\begin{proof}
  This follows from~\cite[Proposition~\ref*{LOT:lem:g-B-integral}]{LOT1} by drilling.
\end{proof}

\begin{lemma}
  \label{lem:GradingOfJuxatposition}
  If $B_1\in\pi_2(\x,\y)$ and $B_2\in\pi_2(\y,\w)$, then
  \begin{equation}
    \label{eq:GradingOfJuxtaposition}
    g'(B_1*B_2)=g'(B_1)\cdot g'(B_2).
  \end{equation}
\end{lemma}

\begin{proof}
  This follows from \cite[Lemma~\ref*{LOT:lem:gB-mult}]{LOT1} by drilling.
\end{proof}

For $\x \in \Gen(\HD)$, let $P'_\x \subset \bigGroup_\AAm(\bdy\HD)$ be
$g'(\pi_2(\x,\x))$.

\begin{corollary}
  \label{cor:GradingSet}
  For $\x\in\Gen(\HD)$, $P'_{\x}$ is a subgroup of
  $\bigGroup_\AAm(\partial \HD)$. Also, if $\y\in\Gen(\HD)$ is another generator and $C\in\pi_2(\x,\y)$, then
  $P'_{\x}=g'(C)\cdot P'_{\y}\cdot g'(C)^{-1}$.
\end{corollary}

\begin{proof}
  Both parts follow immediately from Lemma~\ref{lem:GradingOfJuxatposition}.
\end{proof}

\begin{definition}
  Fix $\x_0\in\Gen(\HD,\spinc)$. Let $S'_{\DA}(\HD,\x_0)$ denote the
  quotient~$\tR(P'_{\x_0})\backslash\bigGroup_\DA(\bdy\HD)$ as a set
  with a right action of $\bigGroup_{\DA}(\partial \HD)$, or equivalently as a left-right
  $(\bigGroup(-\PMC_L),\bigGroup(\PMC_R))$-set.
  There is a
  grading on $\CFDAa(\HD,\spinc)$ with values in
  $\bigGrSet_{\DA}(\HD,\x_0)$, defined by $\gr'_{\x_0}(\x)=[\tR(g'(B))]$
  for any $B\in\pi_2(\x_0,\x)$.
\end{definition}

With the above definition, Lemma~\ref{lem:GradingOfJuxatposition}
ensures that if $B\in\pi_2(\x,\y)$, then
\begin{equation}
  \label{eq:gr-x-y}
  \grb_{\x_0}(\y)=\grb_{\x_0}(\x)\cdot \tR(g'(B)),
\end{equation}
where the
multiplication on the right is right translation
in $\bigGroup_{\DA}(\partial\HD)$ of the right coset $\grb_{\x_0}(\x)$.

\begin{lemma}
  \label{lem:GradedModule}
  If $B\in\pi_2(\x,\y)$, and $(\vec{\rho}\,^L,\vec{\rhos}^R)$ is compatible with $B$,
  then
  \[\lambda^{\abs{\vec{\rhos}^R}-\ind(B;\vec{\rho}\,^L;\vec{\rhos}^R)}
  \cdot \gr'(\vec{\rhos}^R)
  = \tR(g'(B))\cdot \gr'(-\vec{\rho}\,^L)\]
  inside $\bigGroup_{\DA}(\bdy\HD)$.
\end{lemma}

\begin{proof}
  This is a combination of the following facts:
  \begin{align*}
    \gr'(\vec{\rhos}^R)&=(\iota(\vec{\rhos}^R),0,\partial_R^\partial B) \\
    \gr'(-{\vec\rho}\,^L)&=(-\abs{\vec{\rho}\,^L}-\iota({\vec\rho}\,^L),-\partial_L^{\partial} B,0) \\
    \ind(B;\vec\rho\,^L;\vec\rhos^R)&=
    e(B)+n_\x(B)+n_\y(B)+\iota(\vec\rhos^R)+\iota(\vec\rho\,^L)+\abs{\vec{\rho}\,^L}+\abs{\vec{\rhos}^R}.
  \end{align*}
  The first of these equations is verified
  in \cite[Lemma~\ref*{LOT:lem:iota-grading}]{LOT1};
  the second is verified
  in the proof of \cite[Lemma~\ref*{LOT:lem:index-vs-gB}]{LOT1};
% XXX: Use this sentence after InvPairing updated publicly next:
%  the second is verified in \cite[Lemma~\ref*{LOT:eq:gr-rho-seq}]{LOT1};
  and the third is
  the definition, Equation~\eqref{eq:ind3}.
  Note that the two terms
  on the right hand side, $\grb(-\vec{\rho}\,^L)$ and $\tR(g'(B))$, commute
  with each other, as the $\SpinC$ component of $\grb(-\vec\rho\,^L)$ is
  the negative of the portion of the $\SpinC$ component of $\tR(g'(B))$ that lies on the left
  boundary.
\end{proof}

\begin{proposition}\label{prop:gr-prime-is-grading}
  The map $\gr'_{\x_0}$ defines a grading of $\CFDAa$ as a \DA\ 
  structure with values in
  the right $\bigGroup_{\DA}(\bdy\HD)$-set
  $\bigGrSet_{\DA}(\HD,\x_0)$. Different
  choices of $\x_0\in\Gen(\HD,\spinc)$
  lead to canonically isomorphic $\bigGroup_{\DA}(\bdy\HD)$-set graded modules.
\end{proposition}

\begin{proof}
  For the first part, suppose that $a(-\vec{\rho}\,^L)\otimes \y$
  appears with non-zero multiplicity in
  $\delta^1_{n+1}(\x,a(\rhos_1),\dots,a(\rhos_n))$. Then there is a
  domain $B \in \pi_2(\x,\y)$ so that
  $({\vec\rho}\,^L,{\vec\rhos}^R)$ is compatible with $B$, and indeed
  $\ind(B,{\vec\rho}\,^L\semico{\vec\rhos}^R)=1$. We need to know
  $$\lambda^{n-1}\cdot \gr'_{\x_0}(\x)\cdot \prod_{i=1}^n\gr'(\rhos_i^R)
  =\gr'(-{\vec\rho}\,^L)\star \gr'_{\x_0}(\y)$$
  (where here $\star$ refers to the left action of 
  $\bigGroup(-\PMC_L)$, which in turn can be viewed
  as right translation by an element of 
  $\bigGroup(-\PMC_L)^{\op}\subset\bigGroup_{\DA}(\bdy\HD)$).
  But this follows from Lemma~\ref{lem:GradedModule} and
  Equation~\eqref{eq:gr-x-y}.

  Suppose now $\x_0$ and $\x_1$ are two different choices of generator
  both of which represent~$\spinc$. This means that there is a domain $C\in\pi_2(\x_0,\x_1)$.
  We define now an identification
  $$\Phi'^{\x_1}_{\x_0}\co \tR(P'_{\x_1})\backslash \bigGroup_\DA(\bdy\HD) \to \tR(P'_{\x_0})\backslash\bigGroup_\DA(\bdy\HD)$$
  by 
  $$\Phi'^{\x_1}_{\x_0}(\tR(P'_{\x_1})\cdot h)= \tR(P'_{\x_0})\cdot \tR(g'(C))\cdot h.$$
  This gives a well-defined map on coset spaces since
  $P'_{\x_0}=g'(C)\cdot P'_{\x_1}\cdot g'(C)^{-1}$, by Corollary~\ref{cor:GradingSet}.
  Now, for any other $\y$ representing $\spinc$ and any $B \in \pi_2(\x_1,
  \y)$,
  \begin{align*}
    \Phi'^{\x_1}_{\x_0}(\gr'_{\x_1}(\y)) &=
    \Phi'^{\x_1}_{\x_0}(\tR(P'_{\x_1}) \cdot \tR(g'(B))) 
    =\tR(P'_{\x_0}) \cdot \tR(g'(C))\cdot \tR(g'(B)) \\
    &=\tR(P'_{\x_0}) \cdot \tR(g'(C*B)) 
    =\gr'_{\x_0}(\y)
  \end{align*}
  by the various definitions and another application of Corollary~\ref{cor:GradingSet}.
  Thus, the desired isomorphism
  $$(\CFDAa(\HD,\spinc),\gr_{\x_1}',\bigGrSet_{\DA}(\HD,\x_1))\mapsto 
  (\CFDAa(\HD,\spinc),\gr_{\x_0}',\bigGrSet_{\DA}(\HD,\x_0))$$
  is supplied by the identity map on the modules, combined with
  the map $\Phi'^{\x_1}_{\x_0}$ on the $\bigGroup_{\DA}(\bdy\HD)$-sets.
\end{proof}

Finally, we comment briefly on the case that $\HD$ has no generators
representing the $\SpinC$-structure $\spinc$. In this case,
$\CFDAa(\HD,\spinc)$ is the trivial module, but (perhaps) we should
still specify its grading set. Choose another diagram $\HD'$ so that there
is a generator $\x_0\in\Gen(\HD')$ with $\spinc(\x_0)=\spinc$, and
define the grading set for $\CFDAa(\HD,\spinc)$ to be
$\bigGrSet_{\DA}(\HD',\x_0)$ and the grading on $\CFDAa(\HD,\spinc)$ to be
the unique map from $\Gen(\HD,\spinc)=\emptyset$ to
$\bigGrSet_{\DA}(\HD',\x_0)$.  A simplified version of the invariance proof
(keeping track only of the $G'$-set gradings, and not the modules
themselves) shows that, up to isomorphism (in the category of $G'$-set
graded bimodules), this is independent of the choice of $\HD'$; see
the proof of Proposition~\ref{prop:CFDA-inv-graded} for more details.

\subsubsection{Refined gradings}
\label{subsec:RefineCFDAGrading}

We now give $\CFDAa(Y,\spinc)$ a grading by a left-right
($\smallGroup(-\PMC_L),\allowbreak\smallGroup(\PMC_R)$)-set, using the smaller grading group from
Section~\ref{subsec:SmallGroup}.  Let $G_\AAm(\bdy\HD) = G(\PMC_L)
\times_\lambda G(\PMC_R)$ and $G_{\DA}(\bdy\HD) = G(-\PMC_L)^\op
\times_\lambda G(\PMC_R)$.

The existence of a refinement is a formal consequence of the
following:

\begin{lemma}
  \label{lem:Refinable}
  The image $P'_{\x}$ of $\pi_2(\x,\x)$ in $\bigGroup_\AAm(\partial \HD)$ is in
  fact contained in $\smallGroup_\AAm(\partial \HD)$. 
  Moreover, given two generators $\x$ and $\y$ representing $\spinc$,
  if $\gr'(\x)\cdot \tR(g)=\gr'(\y)$, with $g=g_L\times_{\lambda} g_R$,
  then $R(g_L)$ is compatible with
  the idempotents $I_{L,D}(\x)$ and $I_{L,D}(\y)$,
  and $g_R$ is compatible with the idempotents
  $I_{R,A}(\x)$ and $I_{R,A}(\y)$, in the sense of
  Definition~\ref{def:Compatible}.
\end{lemma}

\begin{proof}
  Suppose $\x$ and $\y$
  represent $\spinc$, and let $B\in\pi_2(\x,\y)$. It is clear
  that the homology class $\bdy^{\bdy}_L(B)$
  is compatible with the idempotents $I_{L,A}(\x)$ and $I_{L,A}(\y)$,
  or equivalently $r_*(\bdy^\bdy_L(B))$ is compatible with
  $I_{L,D}(\x)$ and $I_{L,D}(\y)$.  Similarly, $\bdy^{\bdy}_R(B)$ is
  compatible with $I_{R,A}(\x)$ and $I_{R,A}(\y)$.

  Specializing to the case where $\x=\y$, it follows that $P'_{\x}$ is
  contained in $\smallGroup_\AAm(\bdy\HD)\subset\bigGroup_\AAm(\bdy\HD)$.

We turn to the second condition. 
  For notational simplicity, let $I_\AAm(\x)$ denote the pair of
  idempotents $(I_{L,A}(\x),I_{R,A}(\x))$.  Since $\x$ and $\y$ both
  represent~$\spinc$, there is some $B\in\pi_2(\x,\y)$. Moreover,
  $\gr'(\y)=\gr'(\x)\cdot \tR(g'(B))$, which ensures that $g=h \cdot g'(B)$ for some $h\in
  P'_{\x}$ (since $\grb(\x)$ and $\grb(\y)$ are cosets
  of~$P'_\x$). Since $h\in \smallGroup_\AAm(\bdy\HD)$ and $g'(B)$ is compatible with $I(\x)$
  and $I(\y)$, it follows readily that
  $g$ is also compatible with the idempotents $I_\AAm(\x)$ and $I_\AAm(\y)$, or
  equivalently that $R(g_L)$ is compatible with $I_{L,D}(\x)$ and
  $I_{L,D}(\y)$ and $g_R$ is compatible with $I_{R,A}(\x)$ and
  $I_{R,A}(\y)$.
\end{proof}

Lemma~\ref{lem:Refinable} ensures that the type \DA\ bimodule
$\CFDAa(\HD,\spinc)$ is refinable, in the sense of
Definition~\ref{def:Refinable}, provided that our Heegaard diagram has
a generator which represents~$\spinc$ (except now we are using left-right
$(\smallGroup(-\PMC_L),\smallGroup(\PMC_R))$-sets, rather
than just right $\smallGroup$-sets). Thus, the analogue of
Lemma~\ref{lem:RefineGrading} (adapted to bimodules) applies,
allowing us to think of $\CFDAa(\HD,\spinc)$ as a left-right
$(\smallGroup(-\PMC_L),\smallGroup(\PMC_R))$-graded type \DA\ structure.

More concretely, fix a reference point $\x_0\in\Gen(\HD,\spinc)$, and
fix refinement data $\psi_{L,A}$ and $\psi_{R,A}$ for $\Alg(\PMC_L)$ and
$\Alg(\PMC_R)$, respectively.  Let $\psi_{L,D}$ be the reverse of
$\psi_{L,A}$ (see Definition~\ref{def:reverse-grading-refine}), which
is grading refinement data for $\Alg(-\PMC_L)$ by
Lemma~\ref{lem:reverse-grading-refine}.
Define $\psi_\AAm(\x) \in G'_\AAm(\bdy\HD)$ 
and $\psi_\DA(\x) \in G'_\DA(\bdy\HD)$ by
\begin{align*}
\psi_\AAm(\x)&=\bigl(\psi_{L,A}(I_{L,A}(\x)), \psi_{R,A}(I_{R,A}(\x))\bigr)\\
\psi_\DA(\x)&=\bigl(\psi_{L,D}(I_{L,D}(\x))^{-1},\psi_{R,A}(I_{R,A}(\x))\bigr)= \tR(\psi_\AAm(\x)).
\end{align*}
For $B \in \pi_2(\x,\y)$, let $g(B) = \psi_\AAm(\x)
\cdot g'(B) \cdot \psi_\AAm(\y)$; by Lemma~\ref{lem:Refinable}, $g(B)
\in G_\AAm(\bdy\HD)$.  Let $P_\x = g(\pi_2(\x,\x))$.
Let
$S_{\DA}(\HD,\x_0)$ be the quotient
$\tR(P_{\x_0})\backslash\smallGroup_\DA(\partial \HD)$ as a right
$G_{\DA}(\bdy\HD)$-set.  
For any $\x\in\Gen(\HD,\spinc)$, define
\begin{equation}\label{eq:bimod-grading-refine}
\gr_{\x_0}(\x)=\gr'_{\x_0}(\x) \cdot \tR(\psi_\AAm(\x)^{-1})
  = \gr'_{\x_0}(\x) \cdot \psi_\DA(\x)^{-1}
  = \psi_{L,D}(\x) \cdot \gr'_{\x_0}(\x) \cdot \psi_{R,A}(\x)^{-1}
\end{equation}
as an element of $S_{\DA}(\HD,\x_0)$.
(Compare Formula~\eqref{eq:RefineTheGrading}.)

For a different choice of initial point $\x_1$ representing $\spinc$,
we define a canonical identification of grading sets as in the proof
of Proposition~\ref{prop:gr-prime-is-grading}:
\begin{equation}
  \label{eq:IdentifyGradingSets}
  \Phi^{\x_1}_{\x_0}(\tR(P_{\x_1})\cdot h)= \tR(P_{\x_0})\cdot \tR(g(C))\cdot h,
\end{equation}
where $C\in\pi_2(\x_0,\x_1)$.

\begin{proposition}
  \label{prop:CFDA-inv-graded}
  Fix a $\SpinC$ structure $\spinc$ on some strongly bordered
  $3$-manifold $\lsub{\PMC_L}Y_{\PMC_R}$, and let $\HD$ and $\HD'$ be
  Heegaard diagrams for $Y$ in each of which there is at least one
  generator representing $\spinc$.  Fix refinement data $\psi_L$ and
  $\psi_R$ for $\Alg(-\PMC_L)$ and $\Alg(\PMC_R)$ respectively.  Then
  the induced bimodules
  $\lsupv{\AlgA{-\PMC_L}}\CFDAa(\HD,\spinc)_{\AlgA{\PMC_R}}$ and
  $\lsupv{\AlgA{-\PMC_L}}\CFDAa(\HD',\spinc)_{\AlgA{\PMC_R}}$ are
  $\Ainf$-homotopy equivalent as left-right
  $(\smallGroup(-\PMC_L),\smallGroup(\PMC_R))$-graded bimodules.
\end{proposition}

\begin{proof}
  Invariance in the ungraded sense was verified in
  Proposition~\ref{prop:CFDA-inv}. We will first produce a map of
  $\bigGroup_\DA(\bdy\HD)$-sets $\bigGrSet_{\DA}(\HD,\spinc)\to\bigGrSet_{\DA}(\HD',\spinc)$
  compatible with the continuation isomorphism of
  Proposition~\ref{prop:CFDA-inv}. The isomorphism of
  Proposition~\ref{prop:CFDA-inv} is a sequence of elementary
  isomorphisms corresponding to stabilizations, handleslides, and
  isotopies and changes of almost complex structure, and we will treat
  these separately. Note that, by
  Proposition~\ref{prop:gr-prime-is-grading}, in each case we are free
  to choose base generators for $\HD$ and $\HD'$ that are
  convenient.

  Suppose $\HD'$ is obtained from $\HD$ by stabilizing in the region
  containing the basepoint $z$. Then there is an obvious
  identification of generators in $\Gen(\HD,\spinc)$ and
  $\Gen(\HD',\spinc)$. Fix a base generator $\x_0\in\Gen(\HD,\spinc)$
  and let $\x_0'$ denote the corresponding generator of
  $\Gen(\HD',\spinc)$. Then the obvious identification between
  $P'_{\x_0}$ and $P'_{\x_0'}$ induces an identification of grading
  sets, which commutes with the stabilization isomorphism on
  $\CFDAa$.

  Next, suppose $\HD'$ is obtained from $\HD$ by a handleslide. There
  are several cases: a handleslide among the $\beta$-circles, among
  the $\alpha$-circles, or of an $\alpha$-arc over an $\alpha$-circle.
  The case of sliding an arc over a circle is the most complicated, so
  we will restrict to that one; for definiteness, suppose that
  $\alpha_1^{a,L}$ is slid over $\alpha_1^c$. As in
  \cite[Section~\ref*{LOT:sec:CFD-handleslide}]{LOT1}, arrange that
  each $\alpha$-arc of $\HD'$ is close to
  the corresponding $\alpha$-arc of $\HD$ and intersects it in a
  single point, denoted $\theta_i^L$ or $\theta_i^R$, and that for
  each $\theta_i^L$ (respectively $\theta_i^R$) there is a bigon in
  $(\Sigma,\alphas,\alphas')$ originating at $\theta_i^L$
  (respectively $\theta_i^R$) (i.e., that the $\theta_i^L$ would
  correspond to the top graded generator in a closed diagram). Let
  $f_{\HD,\HD'}$ denote the triangle map giving the isomorphism of
  Proposition~\ref{prop:CFDA-inv}.

  Then, fix a base generator $\x_0\in\Gen(\HD,\spinc)$. There is a
  corresponding generator $\x_0'\in\Gen(\HD',\spinc)$ so that there is
  a provincial domain in $(\Sigma,\alphas,\alphas',\betas)$ connecting
  $\x_0$, $\x_0'$ and a generator $\Theta_0$ composed of $\theta_i^L$
  and $\theta_i^R$; and so that this domain consists of a disjoint
  union of triangles supported in the isotopy region and possibly an
  annulus with boundary on $\alpha_1^{a,L}$, $\alpha_1^{a,L,\prime}$
  and~$\alpha_1^c$. Then there is an obvious identification between
  $P'_{\x_0}$ and $P'_{\x'_0}$, which gives an identification of grading
  sets $\bigGrSet_{\DA}(\HD,\spinc)\cong\bigGrSet_{\DA}(\HD',\spinc)$.

  To see that this identification is compatible with the isomorphism
  of Proposition~\ref{prop:CFDA-inv}, recall that to each domain $B'$
  in $(\Sigma,\alphas,\alphas',\betas)$ counted in the triangle map
  there is an associated domain $B$ in
  $(\Sigma,\alphas,\betas)$. Moreover, if $(B',\vec{\rhos})$
  contributes to $f_{\HD,\HD'}$ then $\ind(B,\vec{\rhos})=0$. It
  follows that from this and the fact that the map $f_{\HD,\HD'}$
  treats the $\vec{\rhos}$ in the same way that $\delta^1$ does that
  the map $f_{\HD,\HD'}$ preserves the relative $\bigGroup_{\DA}(\bdy\HD)$ gradings.

  The case that $\HD'$ differs from $\HD$ by an isotopy or a change of
  complex structure is similar to, but easier than, the case of a
  handleslide, so we leave it to the reader.

  Finally, we turn to the $\smallGroup_{\DA}(\bdy\HD)$-set gradings. In view of
  Lemmas~\ref{lem:RefineGrading} and~\ref{lem:Refinable}, the
  $\bigGroup$-set gradings can be lifted to $\smallGroup$-set
  gradings. Moreover, since we choose the same refinement data for the
  two Heegaard diagrams, compatibility of the map of $\bigGroup$-sets
  with the isomorphism of Proposition~\ref{prop:CFDA-inv} implies
  compatibility of the map of $\smallGroup$-sets with the isomorphism
  of Proposition~\ref{prop:CFDA-inv}.
\end{proof}

The above proposition allows us to use Heegaard-diagram-free notation for
$\CFDAa$: we write 
$\lsupv{\AlgA{-\PMC_L}}\CFDAa(Y,\spinc)_{\AlgA{\PMC_R}}$
for
$\lsupv{\AlgA{-\PMC_L}}\CFDAa(\HD,\spinc)_{\AlgA{\PMC_R}}$,
where $\HD$ is any Heegaard diagram which represents $Y$. We also
let $\smallGrSet_\DA(Y,\spinc)$ (or just $\smallGrSet(Y,\spinc)$)
denote the corresponding grading set.

Moreover, we can let
\begin{equation}
  \label{eq:DefGradingSet}
\smallGrSet(Y)=\cup_{\spinc\in\SpinC(Y)}\smallGrSet(Y,\spinc),
\end{equation}
and
define 
$$\lsupv{\AlgA{-\PMC_L}}\CFDAa(\HD)_{\AlgA{\PMC_R}}=
\bigoplus_{\spinc\in\SpinC(Y)}\lsupv{\AlgA{-\PMC_L}}\CFDAa(Y,\spinc)_{\AlgA{\PMC_R}},$$
and think of it as graded by $\smallGrSet(Y)$.

\begin{remark}
  \label{rem:DD-AA-grading}
  Instead of using induction and restriction, the gradings on the
  bimodules $\CFDDa(Y)$ and $\CFAAa(Y)$ can be treated similarly to
  the discussion above. The invariant $\CFAAa(\HD)$ is graded by the
  right $\bigGroup_{\AAm}(\bdy\HD)$-set
  $\bigGrSet_{\AAm}(\bdy\HD)=P'_\x\backslash\bigGroup_{\AAm}(\bdy\HD)$.
  The refined grading on $\CFAAa(\HD)$ is by the right
  $\smallGroup_{\AAm}=\smallGroup(\PMC_L)\times_\ZZ\smallGroup(\PMC_R)$-set
  $\smallGrSet_{\AAm}(\bdy\HD)=P_\x\backslash\smallGroup_{\AAm}$. Tracing
  through the definitions shows that this grading set (and the
  corresponding grading) is the same as that given by the restriction
  functor.

  Similarly, the invariant $\CFDDa(\HD)$ is graded by the left
  $\bigGroup_{\DD}(\bdy\HD)=\bigGroup(-\PMC_L)\times_\ZZ\bigGroup(-\PMC_R)$-set
  $\bigGrSet_{\DD}(\bdy\HD)=\bigGroup_{\AAm}(\bdy\HD)/\mathit{RR}(P'_\x)$ (where
  $\mathit{RR}$ denotes the map $
  \bigGroup(\PMC_L)^\op\times_\ZZ\bigGroup(\PMC_R)^\op\to
  \bigGroup(-\PMC_L)\times_\ZZ\bigGroup(-\PMC_R)$ gotten by applying
  $R$ to each factor). The refined grading on $\CFDDa(\HD)$ is by the
  left
  $\smallGroup_{\DD}(\bdy\HD)=\smallGroup(-\PMC_L)\times_\ZZ\smallGroup(-\PMC_R)$-set
  $\smallGrSet_{\DD}(\bdy\HD)=\smallGroup_{\AAm}(\bdy\HD)/RR(P_\x)$.
  Again, it follows from the definitions that this agrees with the
  grading given by the induction functor.
\end{remark}

\subsection{Invariance}
We collect the results from this section into the following:
\begin{theorem}
  \label{thm:gradedInvarianceOfBimodules}
  Let $Y_{12}$ be a three-manifold with boundary, strongly bordered by
  $\PMC_1$ and $\PMC_2$, and fix grading refinement data for
  $\Alg(\PMC_1)$ and $\Alg(\PMC_2)$. Then, we can associate the
  following $G$-set graded bimodules to $Y_{12}$:
  $\CFAAa(Y_{12})_{\AlgA{\PMC_1},\AlgA{\PMC_2}}$,
  ${}_{\AlgA{-\PMC_1}}\CFDAa(Y_{12})_{\AlgA{\PMC_2}}$ and
  ${}_{\AlgA{-\PMC_1},\AlgA{-\PMC_2}}\CFDDa(Y_{12})$.  The
  quasi-isomorphism types of these $G$-set graded bimodules are
  diffeomorphism invariants of the bordered three-manifold $Y_{12}$
  with its strong boundary framing.
\end{theorem}
\begin{proof}
  Without the gradings, this is immediate from
  Propositions~\ref{prop:CFAA-inv},~\ref{prop:CFDD-inv}
  and~\ref{prop:CFDA-inv}. The fact that the isomorphisms respect the
  grading on $\CFDAa$ is proved in
  Proposition~\ref{prop:CFDA-inv-graded}; the proofs that the
  isomorphisms respect the gradings on $\CFAAa$ and $\CFDDa$ are
  analogous.
\end{proof}

%%% Local Variables: 
%%% mode: latex
%%% TeX-master: "Bimodules"
%%% End: 

\section{Pairing theorems}
\label{sec:PairingTheorems}

Theorems~\ref{thm:Reparameterization}, \ref{thm:Composition}, and
\ref{thm:Hochschild} can all be seen as pairing theorems, which
express how the bordered Floer homology groups transform as bordered
three manifolds are glued in the three situations discussed in
Section~\ref{subsec:GluingDiagrams}.  The aim of the present section
is to study how the bordered invariants change under these three
gluing operations, to obtain proofs of the aforementioned three
theorems. (Indeed, we obtain three generalizations,
Theorems~\ref{thm:GenReparameterization}, \ref{thm:GenComposition} and
\ref{thm:DoublePairing} below.)

\subsection{Pairing along a connected surface}
\label{sec:connected-pairing}

Here is the promised generalization of Theorem~\ref{thm:Reparameterization}:
\begin{theorem}
  \label{thm:GenReparameterization}
  Let $Y_{12}$ be a strongly bordered three-manifold with boundary
  parameterized by $-\PMC_1$ and $\PMC_2$. 
  Let $Y_1$ be a three-manifold with
  boundary parameterized by $\PMC_1$. Then there are $\Ainf$-homotopy equivalences:
  \begin{align*}
    \CFAa(Y_1)\DT_{\AlgA{\PMC_1}}  \CFDAa(Y_{12}) & \simeq \CFAa(Y_{1}\cup_{F_1} Y_{12})\\
    \CFAAa(Y_{12}) \DT_{\AlgA{-\PMC_1}} \CFDa(Y_1)& \simeq \CFAa(Y_{1}\cup_{F_1} Y_{12}) \\
    \CFAa(Y_1)\DT_{\AlgA{\PMC_1}}  \CFDDa(Y_{12}) & \simeq \CFDa(Y_{1}\cup_{F_1} Y_{12}) \\
    \CFDAa(Y_{12}) \DT_{\AlgA{-\PMC_1}}  \CFDa(Y_1) & \simeq \CFDa(Y_{1}\cup_{F_1} Y_{12}),
  \end{align*}
  The first two are equivalences of type $A$ structures over $\Alg(\PMC_2)$,
  while the second two are equivalences of type $D$ structures
  over $\Alg(\PMC_1)$.
\end{theorem}

(Here and later, the $\Ainf$-(bi)modules like $\CFDAa(Y_{12})$ are
required to be appropriately bounded for the tensor product to exist.
This is always possible, as we just choose the Heegaard diagram to be
the appropriate variant of admissible.)

In a similar spirit, we have the following generalization of
Theorem~\ref{thm:Composition}:
\begin{theorem}
  \label{thm:GenComposition}
  Let $Y_{12}$ be a strongly bordered three-manifold with boundary
  parameterized by $-\PMC_1$ and $\PMC_2$. Let $Y_{23}$ be a strongly
  bordered three-manifold with
  boundary parameterized by $-\PMC_2$ and $\PMC_3$. Then there are $\Ainf$-quasi-isomorphisms:
  \begin{align*}
    \CFDAa(Y_{12})\DT_{\AlgA{\PMC_2}} \CFDAa(Y_{23}) & \simeq 
    \lsupv{\AlgA{\PMC_1}}\CFDAa(Y_{12}\cup_{F_2} Y_{23})_{\AlgA{\PMC_3}} \\
    \CFAAa(Y_{12})\DT_{\AlgA{\PMC_2}} \CFDAa(Y_{23}) & \simeq 
    \CFAAa(Y_{12}\cup_{F_2} Y_{23})_{\AlgA{-\PMC_1},\AlgA{\PMC_3}} \\
    \CFDAa(Y_{12})\DT_{\AlgA{\PMC_2}} \CFDDa(Y_{23}) & \simeq 
    \lsupv{\AlgA{\PMC_1},\AlgA{-\PMC_3}}\CFDDa(Y_{12}\cup_{F_2}
    Y_{23}) \\
    \CFAAa(Y_{12})\DT_{\AlgA{\PMC_2}} \CFDDa(Y_{23}) & \simeq 
    \lsupv{\AlgA{-\PMC_3}}\CFDAa(Y_{12}\cup_{F_2} Y_{23})_{\AlgA{-\PMC_1}}
  \end{align*}
\end{theorem}

\begin{proof}[Proof of Theorems~\ref{thm:GenReparameterization} and~\ref{thm:GenComposition}.]
  Both of the proofs of the pairing theorem in~\cite{LOT1} extend
  easily to these cases. To belabor the point, we will prove in detail
  the first equivalences of Theorem~\ref{thm:GenReparameterization}
  via nice diagrams; the proofs of the other parts of the theorems
  proceed similarly.

  So, let $\HD_1$ be a nice diagram for $Y_1$; existence of such is
  guaranteed by
  \cite[Proposition~\ref*{LOT:prop:nice-diagram}]{LOT1}. Let $\HD'$ be
  a Heegaard
  diagram for $Y_{12}$. Apply the algorithm from
  \cite[Proposition~\ref*{LOT:prop:nice-diagram}]{LOT1} to $\drHD'$
  and then fill in the tunnel; the result is a
  Heegaard diagram $\HD_{12}$ for $Y_{12}$ so that ${\drHD}_{12}$ is
  nice. (We will simply call $\HD_{12}$ nice in this case.)

  We turn to the first isomorphism of
  Theorem~\ref{thm:GenReparameterization}. Note that the fact that
  $\HD_1$ and $\HD_{12}$ are nice implies (by
  \cite[Lemma~\ref*{LOT:lem:nice-admissible}]{LOT1})
  that they are admissible; Proposition~\ref{prop:DT-many-cases}
  and Lemma~\ref{lem:admis-bd-DA} imply
  that the box product $\CFAa(\HD_1)\DT_{\AlgA{\PMC_2}}
  \CFDAa(\HD_{12})$ is well defined.  On the other side,
  Lemma~\ref{lem:admiss-glues}
  implies that the glued diagram
  $\HD=\HD_{1}\sos{\bdy}{\cup}{\bdy_L}\HD_{12}$ is admissible. (In
  fact, $\HD$ is nice, and hence admissible by
  \cite[Lemma~\ref*{LOT:lem:nice-admissible}]{LOT1}.)

  We claim that $\CFAa(\HD_1)\DT_{\AlgA{\PMC_2}} \CFDAa(\HD_{12})$ is
  exactly equal to $\CFAa(\HD)$. As a first step, note that there is
  an obvious correspondence between generators. (In particular, only
  $\CFDAa(\HD_{12},0)\subset \CFDAa(\HD_{12})$ contributes to the
  tensor product.) So, we need to check that this identification
  respects the differentials and right module structures.

  On the one hand, the diagram $\HD$ is obviously nice. So,
  by \cite[Proposition~\ref*{LOT:Prop:CurvesInNiceDiagrams}]{LOT1},
  the differential on $\CFAa(\HD)$
  counts (provincial) rectangles and bigons, while the only nontrivial
  (right) algebra actions correspond to unions of half-strips through
  $\bdy_R\HD$. Note that these half strips are entirely contained in
  $\HD_{12}\subset \HD$.

  On the other hand, in $\CFAa(\HD_1)\DT_{\AlgA{\PMC_2}}
  \CFDAa(\HD_{12})$, the differential comes from three different
  contributions:
  \begin{itemize}
  \item Provincial curves in $\HD_1$, which correspond to bigons and
    rectangles by \cite[Proposition~\ref*{LOT:Prop:CurvesInNiceDiagrams}]{LOT1},
  \item Provincial curves in $\HD_{12}$, which again correspond to bigons and
    rectangles by
    \cite[Proposition~\ref*{LOT:Prop:CurvesInNiceDiagrams}]{LOT1}, and
  \item Contributions of the form
    \[
    \x_1\otimes\x_{12}\stackrel{\Id\otimes\delta^1}{\longrightarrow}
    \x_1\otimes (\rho\otimes\y_{12})=(\x_1\otimes\rho)\otimes \y_{12}
    \stackrel{m_2\otimes\Id}{\longrightarrow} \y_1\otimes \y_2.
    \]
  \end{itemize}
  The third kind of contributions correspond exactly to rectangles
  crossing $\bdy\HD_1=\bdy_L\HD_{12}$. The other two correspond to
  bigons and rectangles in $\HD$ contained entirely in one of $\HD_1$
  or $\HD_{12}$. Consequently, the differentials on
  $\CFAa(\HD_1)\DT_{\AlgA{\PMC_2}} \CFDAa(\HD_{12})$ and $\CFAa(\HD)$
  agree.

  Since $\HD_{12}$ is nice, the right module structure comes entirely
  from juxtapositions of half-strips crossing $\bdy_R\HD_{12}$. These
  are exactly the same curves which define the module structure on
  $\CFAa(\HD)$, and they contribute in the same way.

  Thus, we have an isomorphism of right differential modules
  \[
  \CFAa(\HD_1)\DT_{\AlgA{\PMC_2}} \CFDAa(\HD_{12}) \cong \CFAa(\HD).
  \qedhere
  \]
\end{proof}

\begin{figure}
  \centering
  \includegraphics{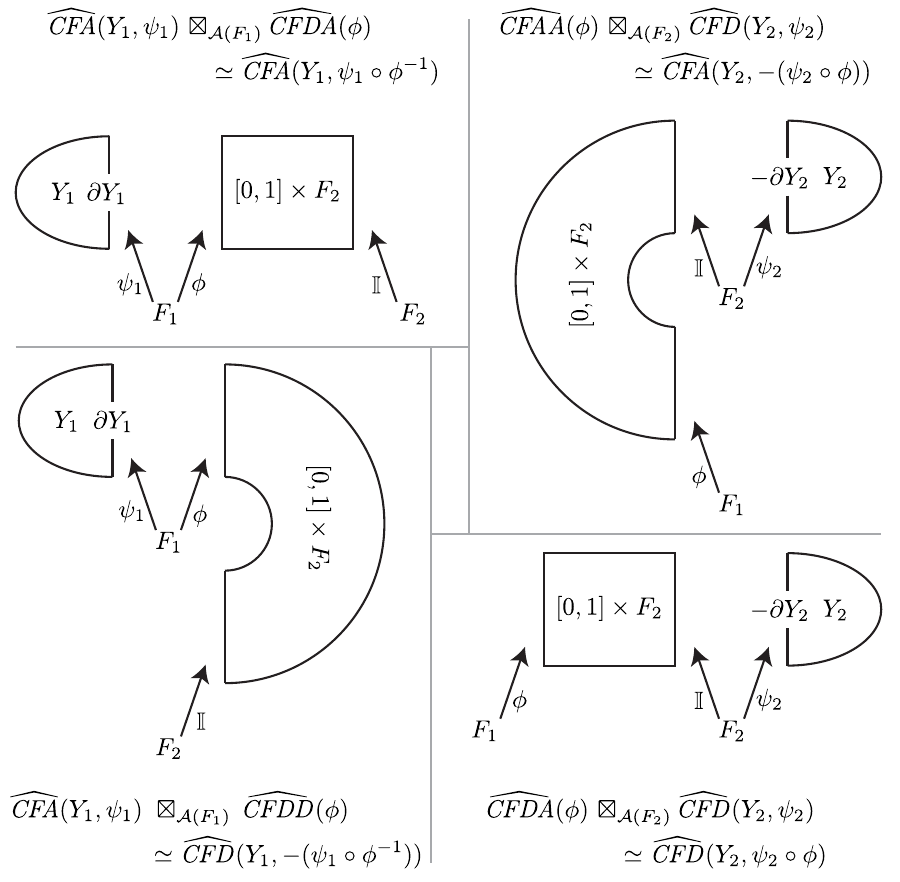}
  \caption{\textbf{Schematic illustration of
      Theorem~\ref{thm:Reparameterization}.} The four cases are shown
    in order left-to-right, top-to-bottom. To compute the
    parametrization of the boundary after
    gluing, start at the unglued boundary component and follow the
    arrows until you reach $\bdy Y_i$, composing the maps labeling the
    arrows (or their inverses).}
  \label{fig:ModBimodPairingSchematic}
\end{figure}

\begin{proof}[Proof of Theorems~\ref{thm:Reparameterization} and~\ref{thm:Composition}]
  Theorem~\ref{thm:Reparameterization} (respectively
  Theorem~\ref{thm:Composition}) is an immediate consequence of
  Theorem~\ref{thm:GenReparameterization} (respectively
  Theorem~\ref{thm:GenComposition}), together with the
  definition of the bimodule of a surface diffeomorphism, and the
  interpretation of the derived tensor product in terms of $\DT$,
  Propositions~\ref{prop:D-to-A-and-back}
  and~\ref{prop:DTP-is-derived-tensor-product}. (A schematic,
  illustrating one way to keep the compositions straight, is given in
  Figure~\ref{fig:ModBimodPairingSchematic}.)
\end{proof}

Similarly, we have the following:

\begin{proof}[Proof of Corollary~\ref{cor:ConvertDintoA}]
  This is a special case of Theorem~\ref{thm:Reparameterization}, using the identity map for $\psi$.
\end{proof}

\subsubsection{Gradings}\label{subsubsec:graded-pairing}
We discuss now how the pairing theorem intertwines the gradings on the
two sides. For definiteness, we will consider the $\DT$ product of two
type \DA\ modules; the other cases are similar.

Fix strongly bordered $3$-manifolds $Y_1$ and $Y_2$ with two boundary
components, where $Y_i$ is parameterized by $\PMC_L(Y_i)$
and $\PMC_R(
Y_i)$, with $\PMC_R( Y_1)=-\PMC_L( Y_2)$, so that we can form the
manifold $Y= Y_1 \sos{\partial_R Y_1}{\cup}{\partial_L Y_2} Y_2$.  For
brevity, let $\PMC_\md$ be $\PMC_R(Y_1) = -\PMC_L(Y_2)$.
As
in Equation~\eqref{eq:DefGradingSet}, we let $\smallGrSet_{\DA}(Y_i)$ and
$\smallGrSet_{\DA}(Y)$ denote the various grading sets for the bordered
Floer homology bimodules.  
Since
$\PMC_L Y=\PMC_L Y_1$ and $\PMC_R Y=\PMC_R Y_2$,
$\smallGrSet_{\DA}(Y_1)\times_{\smallGroup(\PMC_\md)}\smallGrSet_{\DA}(Y_2)$
is
naturally a left-right
$(\smallGroup(-\PMC_L(Y)),\smallGroup(\PMC_R(Y)))$-set.

As before, we will write, for instance, $\smallGroup_{\DA}(\partial Y)$ for
$\smallGroup(-\PMC_L(Y))^\op \times_\lambda \smallGroup(\PMC_R(Y))$.

\begin{theorem}\label{thm:GradedPairing}
  If $\PMC_R(Y_1) = -\PMC_L(Y_2) = \PMC_\md$ and $Y =
  Y_1\sos{\bdy_LY_1}{\cup}
  {\bdy_R Y_2}Y_2$ as above,
  there is an identification of $\smallGroup_{\DA}(\partial Y)$-sets
  \[
  \smallGrSet_{\DA}(Y_1)
  \times_{\smallGroup(\PMC_\md)}
  \smallGrSet_{\DA}(Y_2)\cong \smallGrSet_{\DA}(Y)
  \]
  so that the isomorphism in Theorem~\ref{thm:GenReparameterization}
  is a $\smallGroup_{\DA}(\partial Y)$-set
  graded isomorphism.
\end{theorem}
In fact, we can refine this statement slightly: there is a natural
identification between $\SpinC$-structures on $Y$ and
$\smallGroup(\partial Y)$-orbits in 
$\smallGrSet_{\DA}(Y_1)
\times_{\smallGroup(\PMC_\md)}
\smallGrSet_{\DA}(Y_2)\cong \smallGrSet_{\DA}(Y)$, which is refined by the identification in
Theorem~\ref{thm:GradedPairing}; the identification is given in the
proof.

\begin{proof}[Proof of Theorem~\ref{thm:GradedPairing}]
The  identification of $\smallGroup_{\DA}(\partial Y)$-sets
\[
\smallGrSet_{\DA}(Y_1)\times_{\smallGroup(\PMC_\md)}\smallGrSet_{\DA}(Y_2)\cong \smallGrSet_{\DA}(Y)
\]
is given one $\SpinC$ structure (over $Y$) at a time.
More precisely, fix $\spinc\in\SpinC(Y)$, and let $\spinc_i$
denote its restriction to $Y_i$. We will exhibit an identification
\[
\smallGrSet_{\DA}(Y_1,\spinc_1)\times_{\smallGroup(\PMC_\md)}\smallGrSet_{\DA}(Y_2,\spinc_2)\cong 
\bigcup_{h\in H_1(\partial_R Y_1;\ZZ)}\smallGrSet_{\DA}(Y,\spinc+\PD[h]).
\]
We work with a Heegaard diagram 
$\HD=\HD_1\sos{\PMC_R}{\cup}{\PMC_L}\HD_2$ for $Y$
which has a generator
$\x$ representing $\spinc$, so that the restrictions $\x_i$ of $\x$ to $\HD_i$
represent $\spinc_i$. Thus, our goal is to construct a map
\[
\smallGrSet_{\DA}(\HD_1,\x_1)\times_{\smallGroup(\PMC_\md)}\smallGrSet_{\DA}(\HD_2,\x_2)\cong 
\bigcup_{h\in H_1(\partial_R Y_1;\ZZ)}\smallGrSet_{\DA}(Y,\spinc_z(\x)+\PD[h]).
\]

First, however, we construct the identification on the level of orbit spaces,
i.e., $\SpinC$ structures.
And indeed, before this, we construct a map
\[
p_\x \co
\smallGrSet_{\DA}(\HD_1,\x_1)
\times_{\smallGroup(\PMC_\md)}
\smallGrSet_{\DA}(\HD_2,\x_2)
\to H_1(Y,\partial Y),
\]
as follows.
Recall from Section~\ref{sec:cf-gradings} that $\smallGrSet_{\DA}(\HD_1,\x_1)$ and $\smallGrSet_{\DA}(\HD_2,\x_2)$
are the coset spaces
$\tR(P_{\x_1})\backslash\smallGroup_\DA(\partial\HD_1)$
and $\tR(P_{\x_2})\backslash \smallGroup_\DA(\partial\HD_2)$
respectively.  For brevity, we will write $\tP_{\x}$ for $\tR(P_\x)$.
Thus we can write elements of 
$\smallGrSet_{\DA}(\HD_1,\x_1)
\times_{\smallGroup(\PMC_\md)}
\smallGrSet_{\DA}(\HD_2,\x_2)$
as
$(\tP_{\x_1}\cdot (n_1,\alpha_1,\beta_1))\times (\tP_{\x_2}\cdot
(n_2,\alpha_2,\beta_2))$ where $n_i\in\ZZ$, $\alpha_i\in
H_1(F(-\PMC_L(Y_i)))$
and $\beta_i\in H_1(F(\PMC_R(Y_i)))$.
We then define
$$p_\x ((\tP_{\x_1}\cdot(n_1,\alpha_1,\beta_1)) \times (\tP_{\x_2}\cdot(n_2,\alpha_2,\beta_2)) )
=i_*(\beta_1+\alpha_2)$$
where here $i_*$ is the inclusion map from $H_1(\partial_R
Y_1)=H_1(-\partial_L Y_2)$ to $H_1(Y;\ZZ)$.
With this definition, elements
\[s,t\in
\smallGrSet_{\DA}(\HD_1,\x_1)
\times_{\smallGroup(\PMC_\md)}
\smallGrSet_{\DA}(\HD_2,\x_2)\]
lie in the same $\smallGroup_{\DA}(\partial Y)$-orbit if and only if
$p_\x(s)=p_\x(t)$. 

The map $p_\x$ depends on the base generator~$\x$. However,
the map
\[
q=q_{\x}\co
\smallGrSet_{\DA}(\HD_1,\x_1)
\times_{\smallGroup(\PMC_\md)}
\smallGrSet_{\DA}(\HD_2,\x_2)
\to \SpinC(Y)
\]
defined by 
\[q_{\x}(s_1\times s_2)=\spinc(\x)+p_{\x}(s_1\times s_2)\] is
independent of~$\x$, in the following sense.
If $\x$ and $\y$ are two choices of base generators, then
$q_{\x}\circ \Phi^{\y}_{\x}=q_{\y}$,
where  here
$$\Phi^{\y}_{\x}\co 
\smallGrSet_{\DA}(\HD_1,\x_1)
\times_{\smallGroup(\PMC_\md)}
\smallGrSet_{\DA}(\HD_1,\x_2)
\to 
\smallGrSet_{\DA}(\HD_1,\y_1)
\times_{\smallGroup(\PMC_\md)}
\smallGrSet_{\DA}(\HD_1,\y_2)$$
given by $\Phi^{\y}_{\x}=(\Phi^{\y_1}_{\x_1}\times \Phi^{\y_2}_{\x_2})$ 
is the map gotten by putting together the two canonical identifications of 
grading sets (see Equation~\eqref{eq:IdentifyGradingSets}; see also
Lemma~\ref{prop:gr-prime-is-grading}).
More explicitly, if
there are $C_i\in\pi_2(\x_i,\y_i)$, then
$$\Phi^{\y}_{\x}(\tP_{\y_1}\cdot g_1\times \tP_{\y_2}\cdot g_2)
=(\tP_{\x_1}\cdot \tR(g(C_1))\cdot g_1\times \tP_{\x_2}\cdot \tR(g(C_2))\cdot g_2).$$
Now,
$i_*(\partial^\partial_R[C_1]+r_*(\partial^\partial_L[C_2]))=\epsilon (\x,
\y)$, where $\epsilon(\x,\y)$ is the map giving the difference in
$\SpinC$ structures between $\x$ and $\y$ as in the proof of Lemma~\ref{lem:same-spinc}. Thus, 
\begin{multline*}
  q_{\x} \circ \Phi^{\y}_{\x}(\tP_{\x_1}\cdot (n_1,\alpha_1,\beta_1)
  \times 
   \tP_{\x_2} \cdot (n_2,\alpha_2,\beta_2))   \\
   \begin{aligned}
  &= q_{\x}(\tP_{\x_1}\cdot \tR(g(C_1))\cdot (n_1,\alpha_1,\beta_1)
  \times \tP_{\x_2}\cdot \tR(g(C_2))\cdot (n_2,\alpha_2,\beta_2) )\\
  &=\spinc(\x)+ i_*(\beta_1+\partial^\partial_R[C_1]+r_*(\partial^\partial_L[C_2])
  +\alpha_2) \\
  &= \spinc(\x)+\epsilon(\x,\y)+i_*(\beta_1+\alpha_2) \\
  &=\spinc(\y)+i_*(\beta_1+\alpha_2) \\
  &=q_{\y}(\tP_{\y_1}\cdot (n_1,\alpha_1,\beta_1)\times \tP_{\y_2}
  \cdot (n_2,\alpha_2,\beta_2)),
   \end{aligned}
\end{multline*}
as claimed.

Thus, if we write $\spinc_i=\spinc|_{Y_i}$, the map $q$ defines an
identification of $\smallGroup_{\DA}(\partial Y)$-orbits in
$$\smallGrSet_{\DA}(Y_1,\spinc_1)
\times_{\smallGroup(\PMC_\md)}
\smallGrSet_{\DA}(Y_2,\spinc_2)$$
and those $\SpinC$ structures on $Y$ 
which are of the form $\spinc+\PD[i_*(h)]$ for some $h \in
H_1(\partial_R Y_1)$.

We refine this to a map of grading sets, as follows. Given orbits
$O_1$ and $O_2$ of $\smallGroup_{\DA}(\partial Y_1)$ and
$\smallGroup_{\DA}(\partial Y_2)$ respectively, 
fix a $\smallGroup_{\DA}(\partial Y)$-orbit
$O_{12}$ in $O_1\times O_2$, and suppose that
there is a generator $\x$ for $\HD$ which represents the corresponding 
$\SpinC$ structure. Without loss of generality, we can think of
$\smallGrSet_{\DA}(Y,\spinc)$ as the orbit of $\gr(\x)$,
and its components $\x_i$
as determining the grading sets $\smallGrSet_{\DA}(\HD_i,\x_i)$, so
that $O_{12}$ is contained in the orbit of $\gr(\x_1)\times \gr(\x_2)$.
Then define a map $O_{12}\to \smallGrSet_{\DA}(Y,\spinc)$ as
follows. Any element of $O_{12}$ can be represented as 
\[
\tP_{\x_1}\cdot (n_1,\alpha_1,\beta_1)\times \tP_{\x_2}\cdot (n_2,\alpha_2,\beta_2)
\]
where $\beta_1+\alpha_2=0$.
Then define a map $\phi\co O_{12}\to \smallGrSet_{\DA}(Y,\spinc)$ by
\[
\phi\bigl(\tP_{\x_1}\cdot (n_1,\alpha_1,\beta_1)\times  \tP_{\x_2}\cdot (n_2,\alpha_2,\beta_2)\bigr)
= \tP_{\x}\cdot (n_1+n_2,\alpha_1,\beta_2).
\]
It is clear that this defines a map of $\smallGroup_{\DA}(\partial \HD) =
\smallGroup(-\PMC_L)^\op\times_\lambda\smallGroup(\PMC_R)$ sets.

To verify that this map respects the gradings on
$\CFDAa(\HD_1)\DT\CFDAa(\HD_2)$ and $\CFDAa(\HD)$, suppose that
$\y=\y_1\otimes\y_2$ is another generator of $\CFDAa(\HD)$ and
$B\in\pi_2(\x,\y)$. We can decompose $B$ as $B_1\times B_2$ where
$B_i\in\pi_2(\x_i,\y_i)$ and $\bdy_R^\bdy B_1=\bdy_L^\bdy B_2$. Then
we have 
\begin{align*}
  \gr_\x(\y) &=\tP_\x\cdot(-e(B)-n_\x(B)-n_\y(B),r_*(\bdy_L^\bdy
  B_1),\bdy_R^\bdy B_2)\cdot \psi_\DA(\y)^{-1}\\
  \gr_{\x_1}(\y_1)&=\tP_{\x_1}\cdot(-e(B_1)-n_\x(B_1)-n_\y(B_1),r_*(\bdy_L^\bdy
  B_1),\bdy_R^\bdy B_1)\cdot\psi_\DA(\y_1)^{-1}\\
  \gr_{\x_2}(\y_2)&=\tP_{\x_2} \cdot (-e(B_2)-n_\x(B_2)-n_\y(B_2),r_*(\bdy_L^\bdy
  B_2),\bdy_R^\bdy B_2)\cdot\psi_\DA(\y_2)^{-1}.
\end{align*}
(Recall from Section~\ref{subsec:RefineCFDAGrading} that 
$\psi_\DA(\y) =
(\psi_{L,D}(I_{L,D}(\y))^{-1},\psi_{R,A}(I_{R,A}(\y)))$.) Thus,
\begin{multline*}
  \phi(\gr(\y_1)\times \gr(\y_2)) \\
  \begin{aligned}
  &=\phi\bigl(\tP_{\x_1}\cdot(-e(B_1)-n_{\x_1}(B_1)-n_{\y_1}(B_1),r_*(\bdy_L^\bdy
  B_1),\bdy_R^\bdy
  B_1)\cdot\psi_\DA(\y_1)^{-1}\times{} \\
  &\qquad \tP_{\x_2}(-e(B_2)-n_{\x_2}(B_2)-n_{\y_2}(B_2),r_*(\bdy_L^\bdy
  B_2),\bdy_R^\bdy B_2)\cdot\psi_\DA(\y_2)^{-1}\bigr)\\
  &=\phi\bigl(\tP_{\x_1}(-e(B_1)-n_{\x_1}(B_1)-n_{\y_1}(B_1),r_*(\bdy_L^\bdy
  B_1),\bdy_R^\bdy
  B_1)(\psi_{L,D}(I_{L,D}(\y_1)),0)\times{} \\
  &\qquad \tP_{\x_2}(-e(B_2)-n_{\x_2}(B_2)-n_{\y_2}(B_2),r_*(\bdy_L^\bdy
  B_2),\bdy_R^\bdy B_2)(0,\psi_{R,A}(I_{R,A}(\y_2))^{-1})\bigr)\\
  &= \tP_\x\cdot(-e(B_1)-n_{\x_1}(B_1)-n_{\y_1}(B_1)-e(B_2)-n_{\x_2}(B_2)-
  n_{\y_2}(B_2),r_*(\bdy_L^\bdy
  B_1),\bdy_R^\bdy B_2) \cdot {}\\
  &\qquad (\psi_{L,D}(I_{L,D}(\y_1)), \psi_{R,A}(I_{R,A}(\y_2))^{-1})\\
  &=(-e(B)-n_\x(B)-n_\y(B),r_*(\bdy_L^\bdy B),\bdy_RB)\cdot\psi_\DA(\y)^{-1}\\
  &=\gr(\y_1\otimes\y_2),
  \end{aligned}
\end{multline*}
as desired.
\end{proof}

\subsection{Hochschild homology and knot Floer homology}
\label{sec:self-pairing}
To give a precise statement of the self-pairing theorem, we will need
to discuss the relevant Alexander grading on knot Floer homology for
generalized open books.

Let $Y$ be a strongly bordered three-manifold with two boundary
components specified by $-\PMC$ and $\PMC$, and let $(Y^{\circ},K)$ be
its associated generalized open book, as in
Construction~\ref{construct:GeneralizedOpenBooks}.

Recall that $\PunctF(\PMC)=F(\PMC)\setminus \mathbb{D}^2$.  Then we can
think of $\PunctF\subset Y^{\circ}$ as an embedded Seifert surface for
$K$.  As such, it induces an integral grading on the knot Floer
homology $\CFKa(Y^{\circ},K)$. Specifically, thinking of knot
Floer homology as graded by relative $\SpinC$ structures
$\SpinC(Y^{\circ},K)$, the summand of $\HFKa(Y^{\circ},K)$ in Alexander
grading $i$ is the sum of knot Floer homology groups
over all relative $\SpinC$ structures $\spinc$ with 
$\OneHalf \langle c_1({\widehat \spinc}),[{\widehat F}]\rangle = i$,
where here ${\widehat F}$ is the surface gotten by capping off $\PunctF$
in the zero-surgery of $Y^{\circ}(F)$, and ${\widehat\spinc}$ is the extension
of the relative $\SpinC$ structure~$\spinc$ over the zero-surgery.

On the bordered side, the bimodule $\CFDAa(Y)$ splits according
to the strands grading $\CFDAa(Y)=\bigoplus_{i\in\ZZ}\CFDAa(Y,i)$
(see Equation~\eqref{eq:StrandSplittingDA});
and hence, so does its Hochschild homology.

The following is a generalization of Theorem~\ref{thm:Hochschild}:
\begin{theorem}
  \label{thm:DoublePairing}
  Let $Y$ be a strongly bordered three-manifold with two
  boundary components parameterized by $-\PMC$ and $\PMC$. 
  Let $(Y^\circ,K)$ be the open
  book obtained by gluing the boundary components of $Y$ together and
  performing $0$-surgery on $\gamma$. Then there is an identification
  between the knot Floer homology of the generalized
  open book and the Hochschild homology of the bimodule of $Y$
  \[
  \HFKa(Y^\circ,K)\cong\HH(\lsupv{\Alg(\PMC)}\CFDAa(Y)_{\AlgA{\PMC}}),
  \]
  which identifies the Alexander grading on knot Floer homology with the strands grading on the bimodule; i.e.,
  \[
  \HFKa(Y^\circ,K,i)\cong
  \HH(\lsupv{\AlgAS{\PMC}{i}}\CFDAa(Y,i)_{\AlgAS{\PMC}{i}}).
  \]
  Moreover, this isomorphism intertwines the $\ZZ$-set gradings on
  $\HH(\lsupv{\Alg(\PMC)}\CFDAa(Y)_{\AlgA{\PMC}})$
  (from Lemma~\ref{lem:Hoch-grading}) and on $\HFKa(Y^\circ,K)$.
\end{theorem}
(For the statement about gradings, we have chosen the same grading
refinement data $\psi$ (Definition~\ref{def:grading-refinement}) for the two
sides of $\HD$.)

We prove this theorem in two ways, first with nice diagrams and then
with deforming the diagonal.

\begin{proof}[Proof via nice diagrams.]
  As in the proof of Theorems \ref{thm:GenReparameterization}
  and~\ref{thm:GenComposition}, choose a nice diagram $\HD$ for $Y$.
  By \cite[Lemma~\ref*{LOT:lem:nice-admissible}]{LOT1}, the diagram
  $\HD$ is admissible. Hence,
  by Lemma~\ref{lem:admis-bd-DA}, the bimodule $\CFDAa(\HD)$ is
  bounded. Also, by Lemma~\ref{lem:dbl-bounded-self-glue}, the
  doubly-pointed Heegaard diagram $\HD^\circ$ is weakly admissible.

  Note that, by assumption, $\PMC_R(\HD)\cong-\PMC_L(\HD)$; denote
  $\PMC_R(\HD)$ simply by $\PMC$.

  View $\CFDAa(\HD)$ as a type \DA\  structure with structure maps
  $$\delta_{n+1}\co X(\HD)\otimes \AlgA{\PMC}^{\otimes n}\to
  \AlgA{\PMC}\otimes X(\HD).$$  By Proposition~\ref{prop:HochschildDA},
  the Hochschild homology of $\CFDAa(\HD)$ is computed as the homology
  of $(X(\HD)^\circ,\widetilde{\bdy})$ where $X(\HD)^\circ =
  X(\HD)/[\Idem(\PMC),X(\HD)]$ is the cyclicization of $X(\HD)$ and
  $\tilde{\bdy}$ is as in Formula~\eqref{eq:HochschildDiff}.

  Let $\HD^\circ$ be the doubly-pointed Heegaard diagram for
  $(Y^{\circ},K)$ gotten by self-gluing $\HD$, as in
  Construction~\ref{construct:GenOBHeeg}.  We will show that the
  complex $(X(\HD)^\circ,\widetilde{\bdy})$ is exactly
  $\CFKa(\HD^\circ)$.  First, as an $\Field$-vector space,
  $X(\HD)^\circ$ is isomorphic to $\CFKa(\HD^\circ)$:
  the $\Field$-vector space $X(\HD)^\circ$ has basis
  the generators $\x\in\S(\HD)$ such that $I_{L,D}(\x)=I_{R,A}(\x)$:
  there is a natural one-to-one correspondence between such generators
  in $\SA{\HD}$ and the generators of $\SA{\HD^{\circ}}$.

  Since $\HD$ is nice, the definition of $\widetilde{\bdy}$,
  Formula~\eqref{eq:HochschildDiff}, simplifies considerably. Indeed,
  for $n>1$,
  \[
  \pi\circ(R\circ\delta)^n\circ\iota=0.
  \]
  Consequently, the differential $\widetilde{\bdy}$ has two
  contributions, corresponding to the cases $n=0$ and $n=1$. The
  $n=0$ part of $\widetilde{\bdy}$ corresponds to provincial domains
  in the differential on $\CFDAa(\HD)$, i.e., rectangles and bigons in
  $\HD$. These also contribute in exactly the same way to the
  differential on $\CFKa(\HD^\circ)$.

  The $n=1$ part of $\widetilde{\bdy}$ corresponds to chains of the
  form
  \[
  \x\stackrel{\bdy}{\longrightarrow}
  \rho\y\stackrel{R}{\longrightarrow} \y\rho\stackrel{m_2}{\longrightarrow}\w.
  \]
  These correspond exactly to rectangles in $\HD^\circ$ which cross
  $\bdy_L\HD=\bdy_R\HD$: the first arrow comes from one half of the
  rectangle, which crosses the left boundary in a chord $\rho$, while
  the third arrow comes from the other half, crossing the boundary in
  the same $\rho$. In total,
  this rectangle contributes exactly as it would for $\CFKa(\HD)$.

  Certainly no bigons in $\HD^\circ$ cross through $\bdy_L\HD$, and no
  rectangle can cross $\bdy_L\HD$ twice. So, the differential on
  $\CFKa(\HD^\circ)$ is exactly the same as the differential
  $\widetilde{\bdy}$ on $X(\HD)^\circ$, proving
  the isomorphism.

  We turn next to the strands grading.
  Let $\SA{\HD}^\circ$ be the generators in $\SA{\HD}$ that survive in
  $X(\HD)^\circ$.  This set is naturally identified with $\SA{\HD^\circ}$.  To verify
  the statement about the Alexander and strand gradings, it suffices
  to show that generators $\x\in\SA{\HD}^\circ$ with
  $\#o_R(\x)=k+i$ are mapped
  under the natural one-to-one correspondence to generators
  $\x^{\circ}\in\SA{\HD^{\circ}}$ with Alexander grading equal to $i$.
  (Recall that $o_R(\x)$ denotes the set of $\alpha^R$-arcs
  which are occupied by the generator $\x$.)
 
  To verify the this assertion, observe that the surface $\PunctF$ is
  isotopic to the union of
  \begin{itemize}
  \item a regular neighborhood $N$ of
    $\cup_{i=1}^{2k}\alpha_i^{a,R}\cup\bdy_R\overline{\Sigma}$ and
  \item the descending disk of $\bdy N\setminus \bdy_R\overline{\Sigma}$.
  \end{itemize}
  So, it follows from the description of $\spinc(\x^\circ)$ from
  \cite[Section 2.6]{OS04:HolomorphicDisks} and \cite[Section
  2.3]{OS04:Knots} that
  \[
  \langle c_1({\widehat{ \spinc(\x^\circ)}}),{\widehat\PunctF}\rangle=-2k+2\#(\x^\circ\cap N)=2i.
  \]

  Finally we turn to the $\ZZ$-set gradings. Fix a generator
  $\x^\circ\in\Gen(\HD^\circ)$ for the self-glued Heegaard diagram,
  and let $\x\in\Gen(\HD)$ denote the corresponding generator for the
  bordered Heegaard diagram. Regard $\gr(\x)$ as an element of
  $S_\DA(\HD)/{\sim} = \tP_{\x_0}\backslash G_\DA(\bdy\HD) / {\sim}$,
  where $\sim$ is the equivalence relation
  $$(n,\alpha,\beta)\sim (n,\alpha+\beta,0),$$
  and otherwise the notation is as in Section~\ref{sec:cf-gradings}.
  (This is the same as the equivalence relation from
  Lemma~\ref{lem:Hoch-grading}.  For brevity, we denote
  $G_\DA(\bdy\HD)$ by~$G$.)

  We first show that the divisibility of $\gr(\x)$ and $\gr(\x^\circ)$
  are the same. Indeed, $n\cdot \gr(\x)=\gr(\x)$ means there is a
  periodic domain $P\in\pi_2(\x,\x)$ with
  $\tR(g(P))=(n,\alpha,-\alpha)$. But then $P$ closes up to give a periodic
  domain in $\pi_2(\x^\circ,\x^\circ)$ with $\ind(P)=n$.

  Next we identify the $\ZZ$-orbits in
  $S_\DA(\HD)/{\sim}$ 
  with the
  $\ZZ$-orbits of the grading set of $\HD^\circ$, i.e., the relative
  $\SpinC$-structures on $Y^\circ\setminus K$, as follows. Fix a
  $\SpinC$-structure $\spinc$ on $\HD$ and let $\x$ be a generator
  representing $\spinc$. Recall that the $\SpinC$-structure $\spinc$
  corresponds to the $\smallGroup$-orbit
  of $\gr(\x)$; we will use
  $\x$ as the base generator for this orbit. The sum map
  $G\to H_1(\bdy_L(Y(\HD)))$ given by
  \[
  (n,\alpha,\beta)\mapsto \alpha+\beta
  \]
  does not descend to the $G$-orbit of $\gr(\x)$, as there may be
  elements of $\tP_\x$ with non-trivial image under this map. However, 
  the inclusion $H_1(\bdy_L(Y(\HD)))\to
  H_1(Y^\circ\setminus K)$ kills the image of $\tP_\x$.
  If we further compose with the Poincar\'e duality isomorphism
  $H_1(Y^\circ\setminus K)\cong H^2(Y^\circ,K)$ we get a map
  \[
  p_\x \co \gr(\x)\cdot G \to H^2(Y^\circ,K).
  \]
  This map depends on the choice of $\x$; however, the map $q\co
  \gr(\x)\cdot G\to \SpinC(Y,K)$ defined by $q(\gr(\x)\cdot
  g)=\spinc(\x^\circ)+p_\x(g)$ is
  independent of the choice of $\x$. In fact, $q$ descends to
  an identification of $\ZZ$-orbits in
  $S_\DA(\HD)/{\sim}$ 
  with relative $\SpinC$-structures.

  Now, focus on the $\ZZ$-orbit in
  $S_\DA(\HD)/{\sim}$ 
  which
  contains the generator~$\x$.  The map of grading sets (on
  this $\ZZ$-orbit) is completely determined by the requirement that
  $\gr(\x)$ map to $\gr(\x^\circ)$. It remains to check that this map is
  compatible with the isomorphism
  \[
  \HH(\lsupv{\Alg(\PMC)}\CFDAa(Y)_{\Alg(\PMC)})\cong\HFKa(Y^\circ,K).
  \]
  Let
  $\y$ be in the $\ZZ$-orbit of $\x$. Then there is a domain
  $B\in\pi_2(\x,\y)$ so that $r_*(\bdy_L^\bdy B)=-\bdy_R^\bdy B$. Taking
  $\x$ as our base generator (i.e., setting
  $\gr(\x)=\tP_\x \subset\smallGroup_\DA(\bdy\HD)$), we have
  \begin{align*}
    \gr(\y)&=\tP_\x \cdot (-e(B)-n_\x(B)-n_\y(B),r_*(\bdy_L^\bdy
    B),\bdy_R^\bdy B) \cdot (\psi(I_{L,D}(\y)), \psi(I_{R,A}(\y))^{-1})\\
    \gr(\y)^\circ &= [(-e(B)-n_\x(B)-n_\y(B),0,0)]\\
    &=\lambda^{\gr(\x^\circ,\y^\circ)},
  \end{align*}  
  where $\gr(\y)^\circ$ is defined in
  Definition~\ref{def:Hoch-grading}; the factors $\psi(I_{L,D}(\y))$ and
  $\psi(I_{R,A}(\y))^{-1}$ cancel in $S(\HD)^\circ$ because $I_{L,D}(\y) = I_{R,A}(\y)$;
  $\gr(\x^\circ,\y^\circ)$ denotes the $(\ZZ/n)$-grading
  difference between $\x^\circ$ and $\y^\circ$; and
  the last equality
  follows from the fact that $B$ gives a domain $B^\circ$ in
  $\pi_2(\x^\circ,\y^\circ)$ with the same Euler measure and point
  measures as~$B$.
\end{proof}

\begin{figure}
  \centering
  %Font is 12 point.
  \includegraphics{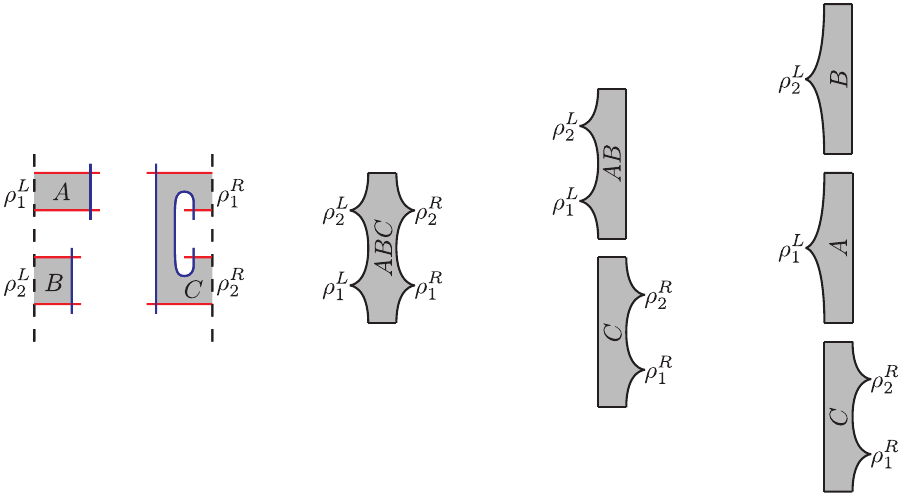}
  \caption{\textbf{Hochschild homology via deforming the diagonal.}
    Far left: a region in a Heegaard diagram $\HD$, contributing to
    the differential on $\CFKa(\HD^\circ)$. All of the $\alpha$-arcs
    shown are parts of different $\alpha$-curves. Center left: a
    schematic of the corresponding self-matched curve. Center right: a
    schematic of the corresponding $\infty$-shifted self-matched
    curve. Far right: a schematic of the corresponding
    $\infty$-shifted $\infty$-self-matched curve. Another interesting
    example can be obtained by reflecting the diagram horizontally.}
  \label{fig:deform-diag}
\end{figure}

\begin{proof}[Proof via deforming the diagonal (sketch).]
  Fix an admissible Heegaard diagram
  $\HD=(\Sigma,\allowbreak\alphas,\allowbreak\betas,\allowbreak z)$
  for $Y$.  Given a holomorphic map $u\co
  S\to\Sigma\times[0,1]\times\RR$ with right punctures
  $p_1^R,\dots,p_i^R$ and left punctures $p_1^L,\dots,p_j^L$ we have
  points
  \begin{align*}
    \ev_R(u)&=(t\circ u(p_1^R),\dots,t\circ u(p_i^R))\in\RR^i\\
    \ev_L(u)&=(t\circ u(p_1^L),\dots,t\circ u(p_i^L))\in\RR^j.
  \end{align*}
  By a \emph{self-matched curve} we mean a holomorphic map $u\co S\to
  \Sigma\times[0,1]\times\RR$ with the same number of right punctures
  as left punctures, labelled by the same Reeb chords in the same
  order, and such that $\ev_R(u)=\ev_L(u)$. Let $\cM^B_{\SM}$ denote
  the set of embedded, self-matched curves in the homology class
  $B$. One can show that, generically, $\cM^B_{\SM}$ is a manifold,
  transversely cut out and of dimension $e(B)+n_\x(B)+n_\y(B)-1$. 
  Moreover, for appropriate almost complex structures, the
  differential on $\CFKa(\HD^\circ)$ is given by
  \[
  \bdy \x=\sum_{\y}\sum_{{\substack{B\in\pi_2(\x,\y)\\
      e(B)+n_{\x}(B)+n_{\y}(B)=1}}}\!\!\!\!\!\!\!\left(\#\cM^B_{\SM}\right)\y.
  \]
  
  Next, we deform the condition of being self-matched in two
  stages. First, for $t\in[0,\infty)$, a \emph{$T$-shifted self-matched
    curve} is a curve $u$ with right punctures $p_1^R,\dots,p_j^R$ and
  left punctures $p_1^L,\dots,p_j^L$ so that for each $i=1,\dots,j$,
  \[
  t\circ u(p_i^R)+T=t\circ u(p_i^L)
  \]
  (and $u$ is asymptotic to the same Reeb chords at $p_i^R$ and
  $p_i^L$).  A $0$-shifted self-matched curve is just a self-matched
  curve as previously defined. Let $\cM^B_{T\Hyph S,\SM}$ denote the moduli
  space of $T$-shifted self-matched curves.  Defining $\bdy$ instead
  using $\cM^B_{T\Hyph S,\SM}$, we get a new chain complex which is
  homotopy equivalent to $\CFKa(\HD^\circ)$.

  Now, take $T\to\infty$. Sequences of $T_i$-shifted self-matched
  curves with $T_i\to\infty$ converge to many-story holomorphic
  combs $(u_1,u_2,\dots,u_l)$ where each
  $u_i\in\pi_2(\x_i,\x_{i+1})$ (with $\x_1=\x$ and $\x_{l+1}=\y$),
  subject to the following condition. Let
  $p_{i,1}^R,\dots,p_{i,j_i^R}^R$ denote the right punctures of $u_i$
  and $p_{i,1}^L,\dots,p_{i,j_i^L}^L$ denote the left punctures of
  $u_i$. Then:
  \begin{itemize} 
  \item for each $i=1,\dots,l-1$, $j_i^R=j_{i+1}^L$, and $j_l^R=j_1^L=0$;
  \item for each $i=1,\dots,l-1$ and $j=1,\dots,j_i$, $u$ is
    asymptotic to the same Reeb chords at $p_{i,j}^R$ and
    $p_{i+1,j}^L$; and
  \item for each $i=1,\dots,l-1$ and $j=1,\dots,j_i-1$,
    \begin{equation}
      t\circ u_i(p_{i,l+1}^R)-t\circ u_i(p_{i,l}^R)=t\circ
      u_{i+1}(p_{i+1,l+1}^L)-t\circ u_{i+1}(p_{i+1,l}^L).\label{eq:inf-match}
    \end{equation}
  \end{itemize}
  We call such combs \emph{$\infty$-shifted self-matched combs}, and
  let $\cM^B_{\infty\Hyph S,\SM}$ denote the moduli space of
  $\infty$-shifted self-matched combs. Again, using
  $\cM^B_{\infty\Hyph S,\SM}$ in place of $\cM^B_{\SM}$, we obtain a new
  chain complex, homotopy equivalent to $\CFKa(\HD^\circ)$.
  (See Figure~\ref{fig:deform-diag} for an example of this shifting
  and the further steps in the proof.)

  Next, we further deform the diagonal as follows. An
  \emph{$\infty$-shifted, $T$-self-matched holomorphic comb} is a
  holomorphic comb $(u_1,u_2,\dots,u_i)$ satisfying the same
  conditions as a $\infty$-shifted self-matched comb, except that
  Equation~\eqref{eq:inf-match} is replaced with the formula:
  \[
  T\cdot\left(t\circ u_i(p_{i,l+1}^R)-t\circ u_i(p_{i,l}^R)\right)=t\circ
  u_{i+1}(p_{i+1,l+1}^L)-t\circ u_{i+1}(p_{i+1,l}^L).
  \]
  An $\infty$-shifted, $1$-self-matched comb is the same as an
  $\infty$-shifted self-matched comb.

  Let $\cM^B_{\infty\Hyph S,T\Hyph \SM}$ denote the moduli space of
  $\infty$-shifted, $T$-self-matched combs.  Replacing $\cM^B_{\SM}$
  by $\cM^B_{\infty\Hyph S,T\Hyph \SM}$, we obtain yet another chain complex
  homotopy equivalent to $\CFKa(\HD^\circ)$. 

  Now, send $T\to\infty$. One can show that sequences of
  $\infty$-shifted, $T_i$-self-matched combs converge to holomorphic
  combs $(u_1,u_2,\dots,u_l)$ such that
  \begin{itemize}
  \item each $u_i$ is asymptotic to a sequence of sets of Reeb chords
    $\vec{\rhos}_i^R=(\rhos_{i,1}^R,\dots,\rhos_{i,j_i^R}^R)$ at
    $\bdy_R\Sigma$, and to a sequence of Reeb chords
    $\vec{\rho}_i^L=(\rho_{i,1}^L,\dots,\rho_{i,j_i^L}^L)$ at
    $\bdy_L\Sigma$;
  \item the sequence of algebra elements
    $a(\rhos_{1,1}^R),\dots,a(\rhos_{1,j_1^R}^R),a(\rhos_{2,1}^R),\dots$ and
    the sequence of algebra elements
    $a(-\vec{\rho}_1^L),a(-\vec{\rho}_2^L),\dots$ are the same; and
  \item each $u_i$ is rigid, as a one-story comb with the specified
    asymptotics.
  \end{itemize}
  Call such a comb a \emph{$\infty$-shifted, $\infty$-self-matched
    holomorphic comb}, and let $\cM^B_{\infty\Hyph S,\infty\Hyph \SM}$ denote
  the moduli space of such combs. Replacing $\cM^B_{\SM}$ by
  $\cM^B_{\infty\Hyph S,\infty\Hyph \SM}$, we obtain another chain complex
  homotopy equivalent to $\CFKa(\HD^\circ)$. But this chain complex
  also has an alternate description: it is the Hochschild complex for
  $\CFDAa(\HD)$ of Proposition~\ref{prop:HochschildDA}. This implies the result.

  The statements about gradings follows exactly as in the ``nice
  diagrams'' version of the proof.
\end{proof}

\begin{proof}[Proof of Theorem~\ref{thm:Hochschild}]
  This is immediate from Lemma~\ref{lem:open-book-is-open-book},
  Theorem~\ref{thm:DoublePairing} and the
  definition of $\CFDAa(\psi)$ in Section~\ref{sec:AutBimodules}.
\end{proof}

%%% Local Variables:
%%% mode: latex
%%% TeX-master: "Bimodules"
%%% End: 

\section{The mapping class group action}
\label{sec:mcg}

In this section, we show that the bimodules $\CFDAa(\phi)$ associated
to surface diffeomorphisms $\phi$ induce an action of the bordered
mapping class group on the derived category of
$\AlgA{\PMC}$-modules. A key step towards establishing this result is
that the bimodule associated to the identity surface diffeomorphism is
the identity map, i.e., $\CFDAa(\Id_{F(\PMC)})$ is homotopy equivalent
to $\lsupv{\AlgA{\PMC}}[\Id]_{\AlgA{\PMC}}$ (Definition~\ref{def:rank-1-DA-mods}), verifying
Theorem~\ref{thm:Id-is-Id}. This is done in
Section~\ref{sec:id-bim}.  The mapping class group action on the
derived categories is stated precisely in
Section~\ref{sec:GroupActionsOnCategories}, and verified in
Section~\ref{sec:MappingClassGroupAction}.

\subsection{Identity Bimodules}
\label{sec:id-bim}
We first prove that the identity map on $\PMC$ induces a bimodule which is quasi-isomorphic to the identity
bimodule on $\AlgA{\PMC}$, as stated in Theorem~\ref{thm:Id-is-Id}.

We make our notation slightly more precise than in the original
statement of the theorem, writing $\Alg(\PMC)$ instead of $\Alg(F)$,
so the desired result is $\CFDAa(\Id_{F(\PMC)})\simeq
\lsupv{\AlgA{\PMC}}[\Id_{\AlgA{\PMC}}]_{\AlgA{\PMC}}$ During the proof
we also make our notation also less precise, dropping the subscript
from the identity $\Id$ (which could be either the algebra of $\PMC$,
or the surface associated to $\PMC$): it should be clear from the
context.

\begin{proof}[Proof of Theorem~\ref{thm:Id-is-Id}]
  We start by arguing that $\lsupv{\Alg(\PMC)}\CFDAa(\Id)_{\Alg(\PMC)}$ is quasi-invertible in the
  sense of Definition~\ref{def:QuasiInvertible}.  Consider the
  canonical bordered Heegaard diagram for the identity diffeomorphism
  (Definition~\ref{def:ConstructHeegaardDiagram}), illustrated in
  Figure~\ref{fig:id-diagram-disk} (left). It is clear from inspection
  that $\lsupv{\Alg(\PMC)}\CFDAa(\Id)_{\Alg(\PMC)}$ is isomorphic to $\AlgA{F}$ as a left
  $\AlgA{F}$-module. Moreover, $\delta^1_1=0$: any non-trivial domain meets the type $A$ boundary
  in the canonical diagram. Thus, Lemma~\ref{lem:characterize-induced} applies, showing that
  $\CFDAa(\Id)_{\Alg(\PMC)}$ is
  isomorphic to $\lsup{\AlgA{F}}[\phi]_{\AlgA{F}}$ for some
  $\Ainf$-endomorphism $\phi$ of~$\AlgA{F}$.

  \begin{figure}
    \centering
    %Font is 12 point.
    \includegraphics{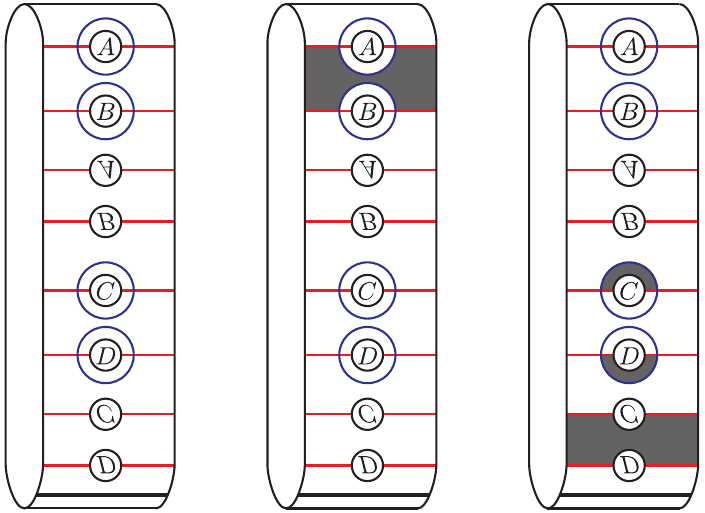}
    \caption{\textbf{The canonical Heegaard diagram for the identity
        map.} Left: the identity map of the split genus $2$
      surface. Center and right: two of the holomorphic disks implying
      $\phi_*(\rho)=\rho$ for any length $1$ chord $\rho$.}
    \label{fig:id-diagram-disk}
  \end{figure}

  For any Reeb chord~$\rho$ of length 1, $\phi_*(\rho) = \rho$, as
  there is an obvious holomorphic disk; see
  Figure~\ref{fig:id-diagram-disk} (center, right).  By
  Proposition~\ref{prop:Rigidity}, $\phi_1$ induces the identity map
  on the homology of~$\AlgA{F}$; in particular, $\phi$ is a
  quasi-isomorphism. By Proposition~\ref{prop:quasi-isom-invert-bimodule}
  $\lsupv{\Alg(\PMC)}\CFDAa(\Id)_{\Alg(\PMC)}$ is quasi-invertible.

  Next, according to Theorem~\ref{thm:Composition},
  $$\lsupv{\Alg(\PMC)}\CFDAa(\Id)_{\Alg(\PMC)} \DT
  \lsupv{\Alg(\PMC)}\CFDAa(\Id)_{\Alg(\PMC)} \simeq
  \lsupv{\Alg(\PMC)}\CFDAa(\Id)_{\Alg(\PMC)}.$$ Hence, applying $\DT$
  with the quasi-inverse to
  $\lsupv{\Alg(\PMC)}\CFDAa(\Id)_{\Alg(\PMC)}$ to both sides of the
  above quasi-isomorphism, we obtain the desired quasi-isomorphism
  \[
  \lsupv{\Alg(\PMC)}\CFDAa(\Id)_{\Alg(\PMC)}\simeq 
  \lsupv{\Alg(\PMC)}[\Id]_{\Alg(\PMC)}.\qedhere
  \]
\end{proof}

\begin{corollary}\label{Cor:MCG-Equivalences} Let
  $\phi\in\MCG_0(F(\PMC))$. Then the functors
  \begin{align*}
  \cdot\DT\lsupv{\AlgA{\PMC}}\CFDAa(\phi)_{\AlgA{\PMC}}&\co 
  \ModCat_{\AlgA{\PMC}}\to
  \ModCat_{\AlgA{\PMC}}\qquad\text{and}\\
  \lsupv{\AlgA{\PMC}}\CFDAa(\phi)_{\AlgA{\PMC}}\DT\cdot&\co\lsup{\AlgA{\PMC}}\ModCat\to \lsup{\AlgA{\PMC}}\ModCat
  \end{align*}
  are auto quasi-equivalences of $\ZZ$-set graded differential categories.
  
  More generally, if $\phi\in\MCG_0(F(\PMC),\allowbreak F(\PMC'))$ is
  in the mapping class
  groupoid, then the functors
  \begin{align*}
  \cdot\DT\lsupv{\AlgA{\PMC}}\CFDAa(\phi)_{\AlgA{\PMC'}}&\co \ModCat_{\AlgA{\PMC}}\to
  \ModCat_{\AlgA{\PMC'}},\\
  \cdot\DT\lsupv{\AlgA{\PMC},\AlgA{-\PMC'}}\CFDDa(\phi)&\co\ModCat_{\AlgA{\PMC}}\to
  \lsup{\AlgA{-\PMC'}}\ModCat,\\
  \lsupv{\AlgA{\PMC}}\CFDAa(\phi)_{\AlgA{\PMC'}}\DT\cdot&\co\lsup{\AlgA{\PMC'}}\ModCat\to \lsup{\AlgA{\PMC'}}\ModCat,
  \qquad\text{and}\\
  \CFAAa(\phi)_{\AlgA{-\PMC},\AlgA{\PMC'}}\DT\cdot&\co\lsup{\AlgA{\PMC'}}\ModCat\to \ModCat_{\AlgA{-\PMC'}}
  \end{align*}
  are quasi-equivalences of $\ZZ$-set graded differential categories.
\end{corollary}
\begin{proof}
  We will prove that 
  \[
  \cdot\DT\lsupv{\AlgA{\PMC}}\CFDAa(\phi)_{\AlgA{\PMC'}}\co \ModCat_{\AlgA{\PMC}}\to
  \ModCat_{\AlgA{\PMC'}}
  \]
  is a quasi-equivalence; the other cases are similar. By
  Lemma~\ref{lem:quasi-invert-gives-equiv-DA} it suffices to show that
  $\lsupv{\AlgA{\PMC}}\CFDAa(\phi)_{\AlgA{\PMC'}}$ is quasi-invertible.

  Fix any Heegaard diagrams $\HD$ for $\phi$ and $\HD'$ for $\phi^{-1}$. 
  It follows from Lemma~\ref{lemma:compose-MCG-diagrams} that
  $\HD{}_{\bdy_R}\cup_{\bdy_L}\HD'$ is a Heegaard diagram for the
  identity map $\Id_{F(\PMC)}\in \MCG_0(F(\PMC),F(\PMC))$ while
  $\HD'{}_{\bdy_R}\cup_{\bdy_L}\HD$ is a Heegaard diagram for $\Id_{F(\PMC')}\in
  \MCG_0(F(\PMC'),F(\PMC'))$. By Theorem~\ref{thm:GenComposition},
  \begin{align*}
    \lsupv{\Alg(\PMC)}\CFDAa(\HD)_{\Alg(\PMC')}\DT
    \lsupv{\Alg(\PMC')}\CFDAa(\HD')_{\Alg(\PMC)}&
    \simeq \lsupv{\Alg(\PMC)}\CFDAa(\HD\sos{\bdy_R}{\cup}{\bdy_L}\HD')_{\Alg(\PMC')}\\
    \lsupv{\Alg(\PMC')}\CFDAa(\HD')_{\Alg(\PMC)}
    \DT
    \lsupv{\Alg(\PMC)}\CFDAa(\HD)_{\Alg(\PMC')}&\simeq 
    \lsupv{\Alg(\PMC')}\CFDAa(\HD'{}_{\bdy_R}\cup_{\bdy_L}\HD)_{\Alg(\PMC')}.
  \end{align*}
  By Lemma~\ref{lemma:diagram-for-mcg-defined} and
  Proposition~\ref{prop:CFDA-inv}, 
  \begin{align*}
    \lsupv{\Alg(\PMC)}\CFDAa(\HD\sos{\bdy_R}{\cup}{\bdy_L}\HD')_{\Alg(\PMC)}&
    \simeq \lsupv{\Alg(\PMC)}\CFDAa(\HD(\Id_{F(\PMC)}))_{\Alg(\PMC)}\\
    \lsupv{\Alg(\PMC')}\CFDAa(\HD'\sos{\bdy_R}{\cup}{\bdy_L}\HD)_{\Alg(\PMC')}
    &\simeq \lsupv{\Alg(\PMC')}\CFDAa(\HD(\Id_{F(\PMC')}))_{\Alg(\PMC')}.
  \end{align*}
  So, by Theorem~\ref{thm:Id-is-Id},
  \begin{align*}
    \lsupv{\Alg(\PMC)}\CFDAa(\HD)_{\Alg(\PMC')}\DT\lsupv{\Alg(\PMC')}\CFDAa(\HD')_{\Alg(\PMC)}&\simeq \lsupv{\AlgA{\PMC}}[\Id]_{\AlgA{\PMC}}\\
    \lsupv{\Alg(\PMC')}\CFDAa(\HD')_{\Alg(\PMC)}
    \DT\lsupv{\AlgA{\PMC}}\CFDAa(\HD)_{\Alg(\PMC')}&\simeq 
    \lsupv{\AlgA{\PMC'}}[\Id]_{\AlgA{\PMC'}}.
  \end{align*}
  This proves the claim.
\end{proof}
Finally, as another corollary, we have
Theorem~\ref{thm:AlgebraDependsOnSurface}, the statement that
different pointed matched circles for a given surface have equivalent
derived categories.
\begin{proof}[Proof of Theorem~\ref{thm:AlgebraDependsOnSurface}]
  This follows from
  Corollary~\ref{Cor:MCG-Equivalences} by choosing any mapping class
  $\phi\in\MCG_0(\PtdMatchCirc_1,\PtdMatchCirc_2)$.
\end{proof}

\subsection{Group actions on categories}
\label{sec:GroupActionsOnCategories}
To state the mapping class group(oid) action precisely requires a
little categorical algebra, which we review in this
subsection.

\begin{definition}
  Let $G$ be a group and $\Cat$ a category. Let $\End(\Cat)$ denote
  the class of functors $F\co \Cat\to\Cat$.
  \begin{itemize}
  \item A \emph{strict action} of $G$ on $\Cat$ is a map
    $A\co G\to\End(\Cat)$ such that $A(\Id)$ is the identity functor
    and if $g,h\in G$ then $A(gh)=A(g)\circ A(h)$.
  \item A \emph{weak action} of $G$ on $\Cat$ is a map $A\co
    G\to\End(\Cat)$ together with a natural isomorphism $A_0$ of the
    identity functor $\Id_{\Cat}$ to $A(\Id)$; and for each $g,h\in
    G$, an isomorphism $A_2(g,h)$ from $A(g)\circ A(h)$ to
    $A(gh)$, so that the following diagrams commute:
    \[
    \begin{tikzpicture}
      \node at (0,0) (tl) {$A(g)\circ A(h)\circ A(k)$};
      \node at (5,0) (tr) {$ A(gh)\circ A(k)$};
      \node at (0,-2) (bl) {$A(g)\circ A(hk)$};
      \node at (5,-2) (br) {$A(ghk)$};
      \draw[->] (tl) to node[above]{\lab{A_2(g,h)\circ \Id_{A(k)}}} (tr);
      \draw[->] (tl) to node[left]{\lab{\Id_{A(g)}\circ A_2(h,k)}} (bl);
      \draw[->] (tr) to node[right]{\lab{A_2(gh,k)}} (br);
      \draw[->] (bl) to node[below]{\lab{A_2(g,hk)}} (br);
    \end{tikzpicture}
    \]
    \[
    \begin{tikzpicture}
      \node at (0,0) (tl) {$A(g)\circ \Id_{\Cat}$};
      \node at (3,0) (tr) {$ A(g)$};
      \node at (0,-2) (bl) {$A(g)\circ A(\Id)$};
      \node at (3,-2) (br) {$A(g\Id)$};
      \draw[->] (tl) to node[above]{\lab{=}} (tr);
      \draw[->] (tl) to node[left]{\lab{A(g)\circ A_0}} (bl);
      \draw[->] (bl) to node[below]{\lab{A_2(g,\Id)}} (br);
      \draw[->] (br) to node[right]{\lab{=}} (tr);
    \end{tikzpicture}    
    \qquad\qquad
    \begin{tikzpicture}
      \node at (0,0) (tl) {$\Id_{\Cat}\circ A(g)$};
      \node at (3,0) (tr) {$A(g)$};
      \node at (0,-2) (bl) {$A(\Id)\circ A(g)$};
      \node at (3,-2) (br) {$A(\Id g)$};
      \draw[->] (tl) to node[above]{\lab{=}} (tr);
      \draw[->] (tl) to node[left]{\lab{A_0\circ A(g) }} (bl);
      \draw[->] (bl) to node[below]{\lab{A_2(\Id,g)}} (br);
      \draw[->] (br) to node[right]{\lab{=}} (tr);
    \end{tikzpicture}    
    \]
    (Compare \cite[Section XI.2]{MacLane98:Categories}.)
  \end{itemize}
  If $\Cat$ has some extra structure (for example, if it is
  triangulated) then we replace $\End(\Cat)$ by the class of
  endofunctors preserving that structure.
\end{definition}
\begin{remark}
  This terminology is not entirely standard. In particular, the reader
  is cautioned that some sources call our weak action a strong action.
\end{remark}

We want to extend the notion of weak group actions to actions of
groupoids, so first we reinterpret it. Recall:
\begin{definition}\label{def:2-functor}
  If $\mathscr{D}$ and $\mathscr{E}$ are $2$-categories then a
  \emph{weak $2$-functor} from $\mathscr{D}$ to $\mathscr{E}$ consists of
  \begin{itemize} 
  \item a map $A\co \ob_{\mathscr{D}}\to \ob_{\mathscr{E}}$;
  \item for each $a,b\in \ob_{\mathscr{D}}$ a functor $A_{a,b}\co
    \Mor_{\mathscr{D}}(a,b)\to \Mor_{\mathscr{E}}(A(a),A(b))$;
  \item for each $a\in\ob_{\mathscr{D}}$ a 2-morphism $A_a\in\TMor(
    \Id_{A(a)},A_{a,a}(\Id_a))$; and
  \item for each $a,b,c\in\ob_{\mathscr{D}}$, $f \in
    \Mor_{\mathscr{D}}(b,c)$, and $g \in \Mor_{\mathscr{D}}(a,b)$,
    a 2-morphism
    \[
    A_{a,b,c}(f,g)\in  \TMor(A_{b,c}(f) \circ A_{a,b}(g), A_{a,c}(f\circ g)),
    \]
    forming a natural transformation of functors.  More precisely, as
    $f$ and $g$ vary, both $A_{b,c}(f) \circ A_{a,b}(g)$ and
    $A_{a,c}(f\circ g)$ give functors from $\Mor_{\mathscr{D}}(b,c)
    \times \Mor_{\mathscr{D}}(a,b)$ to $\Mor_{\mathscr{E}}(A(a),
    A(c))$.  Then $A_{a,b,c}$ is required to be a natural
    transformation between these two functors.
  \end{itemize}
  These data must satisfy:
  \begin{itemize}
  \item For any objects $a,b,c,d\in\ob_{\mathscr{D}}$ and morphisms
    $f\in\Mor_{\mathscr{D}}(c,d)$, $g\in\Mor_{\mathscr{D}}(b,c)$ and
    $h\in\Mor_{\mathscr{D}}(a,b)$ the diagram
    \[
    \begin{tikzpicture}
      \node at (0,0) (tl) {$A_{c,d}(f)\circ (A_{b,c}(g)\circ
        A_{a,b}(h))$};
      \node at (7,0) (tr) {$ (A_{c,d}(f)\circ A_{b,c}(g))\circ
        A_{a,b}(h)$};
    \node at (0,-2) (cl) {$A_{c,d}(f)\circ A_{a,c}(g\circ h)$};
    \node at (7,-2) (cr) {$A_{b,d}(f\circ g)\circ
      A_{a,b}(h)$};
    \node at (0,-4) (bl) {$A_{a,d}(f\circ (g\circ h))$};
    \node at (7,-4) (br) {$A_{a,d}((f\circ g)\circ h)$};
    \draw[->] (tl) to node[above]{\lab{=}} (tr);
    \draw[->] (bl) to node[below]{\lab{=}} (br);
    \draw[->] (tl) to node[left]{\lab{\Id_{A_{c,d}(f)}\circ A_{a,b,c}(g,h)}} (cl);
    \draw[->] (tr) to node[right]{\lab{ A_{b,c,d}(f,g)\circ \Id_{A_{a,b}(h)}}} (cr);
    \draw[->] (cl) to node[left]{\lab{A_{a,c,d}(f,g\circ h)}} (bl);
    \draw[->] (cr) to node[right]{\lab{A_{a,b,d}(f\circ g,h)}} (br);
    \end{tikzpicture}
    \]
    commutes.
  \item For any morphism $f\in\Mor(a,b)$, the diagrams
    \[
    \begin{tikzpicture}
      \node at (0,0) (tl) {$A_{a,b}(f)\circ \Id_{A(b)}$};
      \node at (4,0) (tr) {$A_{a,b}(f)$};
      \node at (0,-2) (bl) {$A_{a,b}(f)\circ A(\Id_b)$};
      \node at (4,-2) (br) {$A_{a,b}(f\circ \Id_b)$};
      \draw[->] (tl) to node[left]{\lab{A_{a,b}(f)\circ A_b}} (bl);
      \draw[->] (tl) to node[above]{\lab{=}} (tr);
      \draw[->] (bl) to node[below]{\lab{A_{a,b,b}(f,\Id)}} (br);
      \draw[->] (br) to node[right]{\lab{=}} (tr);
    \end{tikzpicture}\qquad
    \begin{tikzpicture}
      \node at (0,0) (tl) {$\Id_{A(a)}\circ A_{a,b}(f)$};
      \node at (4,0) (tr) {$A_{a,b}(f)$};
      \node at (0,-2) (bl) {$A(\Id_a)\circ A_{a,b}(f)$};
      \node at (4,-2) (br) {$A_{a,b}(\Id_a\circ f)$};
      \draw[->] (tl) to node[left]{\lab{A_a\circ A_{a,b}(f)}} (bl);
      \draw[->] (tl) to node[above]{\lab{=}} (tr);
      \draw[->] (bl) to node[below]{\lab{A_{a,a,b}(\Id,f)}} (br);
      \draw[->] (br) to node[right]{\lab{=}} (tr);
    \end{tikzpicture}    
    \]
    commute.
  \end{itemize}
\end{definition}
See \cite[Definition 4.1]{Benabou67:biCategories}, which defines the
notion more generally for weak $2$-categories (or bicategories),
although we choose to keep the standard convention for order of
composition.

We may view a group $G$ as a $2$-category $\Gamma$ with a single
object $\bullet$, $\Mor(\bullet,\bullet)=G$, and $\TMor(g,h)$ empty if
$g\neq h$ and consisting of the identity map if $g=h$. A category
$\Cat$ specifies a $2$-category $\CatEnd(\Cat)$ with a single object $\ast$,
$\Mor(\ast,\ast)=\End(\Cat)$, and $\TMor(F_1,F_2)$ the set of natural
transformations from $F_1$ to $F_2$.

\begin{lemma}
  With the above setup, a weak action of $G$ on $\Cat$ is a weak
  $2$-functor from $\Gamma$ to~$\CatEnd(\Cat)$.
\end{lemma}

\begin{proof}
  This is largely immediate from the definitions: given a weak
  action~$A$, we define a weak 2-functor~$B$ by
  \begin{itemize}
  \item $B(\bullet) = \ast$.
  \item $B_{\bullet,\bullet}(g) = A(g)$.  This function on the objects
    of $\Mor_\Gamma(\bullet,\bullet)$ extends trivially to the
    morphisms of $\Mor_\Gamma(\bullet,\bullet)$ (which are the
    2-morphisms of $\Gamma$).
  \item $B_{\bullet,\bullet,\bullet}(g,h) = A_2(g,h)$.  This map
    automatically defines a natural transformation, since
    $\Mor_\Gamma(\bullet, \bullet) \times \Mor_\Gamma(\bullet,
    \bullet)$ has only identity morphisms.
  \end{itemize}
  The diagrams that are required to commute are precisely the same in
  the two cases.
\end{proof}

This leads easily to the notion of a groupoid action.
\begin{definition}\label{def:groupoid-action}
  Let $\Cat_1,\dots,\Cat_n$ be categories and $\Gamma$ a groupoid. Make $\Gamma$
  into a $2$-category with only identity $2$-morphisms. Let
  $\CatEnd(\{\Cat_1,\dots,\Cat_n\})$ denote the full $2$-subcategory of
  $\CatCat$ generated by $\Cat_1,\dots,\Cat_n$. That is,
  \begin{align*}
    \ob(\Cat)&=\{\Cat_1,\dots,\Cat_n\}\\
    \Mor(\Cat_i,\Cat_j)&=\{\text{functors }F\co\Cat_i\to\Cat_j\}\\
    \TMor(F_1,F_2)&=\{\text{natural transformations from $F_1$ to $F_2$}\}.
  \end{align*}
  Then a \emph{weak action of $\Gamma$ on $\{\Cat_1,\dots,\Cat_n\}$} is a weak
  $2$-functor from $\Gamma$ to $\CatEnd(\Cat)$. (Again, if the categories
  $\Cat_1,\dots,\Cat_n$ are triangulated, say, then we restrict to
  triangulated functors.)
\end{definition}
(There is an obvious analogue when $\Cat$ has infinitely many
elements, but we shall not need this.)

Recall from Definition~\ref{def:cat-vary-gr-set} that
$\ModCat_{\AlgA{\PMC}}$ denotes the category of $\smallGroup$-set
graded right $\Ainf$-modules over $\AlgA{\PMC}$, and
$\HMod_*(\ModCat_{\AlgA{\PMC}})$ the category whose objects are
set-graded $\Ainf$-modules over $\AlgA{\PMC}$ and morphisms are
$\Ainf$-homotopy classes of $\Ainf$-module maps.  There is no natural
notion of degree~$0$ morphisms in $\ModCat_{\Alg(\PMC)}$, and so
$\HMod_*(\ModCat_{\Alg(\PMC)})$ is not a triangulated category in the
usual sense.  However, $\ModCat_{\Alg(\PMC)}$ does have a (non-full)
subcategory $\setModCat_{\Alg(\PMC),\sG(\PMC)}$ of modules graded by
the grading group $\sG(\PMC)$ (considered as a right $\sG(\PMC)$-set);
see Definition~\ref{def:cat-fixed-gr-set}.  On this subcategory there
is a notion of degree~$0$ morphisms and so one obtains a triangulated
category $\HMod(\setModCat_{\Alg(\PMC),\sG(\PMC)})$.

The goal of the rest of this section is to prove:
\begin{theorem}\label{thm:MCG-2-functoriality}
  The bimodules $\CFDAa(\phi)$ induce a weak action of the genus $k$
  bordered mapping class groupoid $\MCG_0(k)$ on %the $2$-category
  $\{\HMod_*(\ModCat_{\AlgA{\PMC}})\mid \mathrm{genus}(F(\PMC))=k\}$.
  This action preserves the subcategories %y
  $\{\HMod(\setModCat_{\Alg(\PMC),\sG(\PMC)})\}$ and acts by triangulated
  functors on $\{\HMod(\setModCat_{\Alg(\PMC),\sG(\PMC)})\}$.
\end{theorem}
In particular, for any pointed matched circle $\PMC$, the bimodules
$\CFDAa(\phi)$ for $\phi\co F(\PMC) \to F(\PMC)$ induce a weak action
of the genus $k$ bordered mapping
class group on $\HMod_*(\ModCat_{\AlgA{\PMC}})$, and for different
choices of $\PMC$ these actions are conjugate (in the obvious
sense).
One easily digestible piece of Theorem~\ref{thm:MCG-2-functoriality}
is that tensoring with the bimodules $\CFDAa(\phi)$ induces
equivalences of categories; this fact is
Corollary~\ref{Cor:MCG-Equivalences}.

\begin{remark}\label{rmk:bdy-dehn-twist}
  It is important to note that it is the {\em bordered} mapping class
  group which acts on the category, rather than the ordinary one.  For
  example, if $\PMC$ represents a surface of genus one, and $\phi$
  denotes the mapping class which is gotten by Dehn twist around the
  boundary of $\PMC$, then the tensor product with
  $\phi$ induces a non-trivial action on the category of
  $\Alg(\PMC)$-modules. For instance, there is a module with rank one
  over $\Field$ (and trivial differential) whose tensor product with
  $\CFDAa(\phi)$ has homology with rank $9$; see
  Figure~\ref{fig:UnbasedMappingClass}.
\begin{figure}
\input{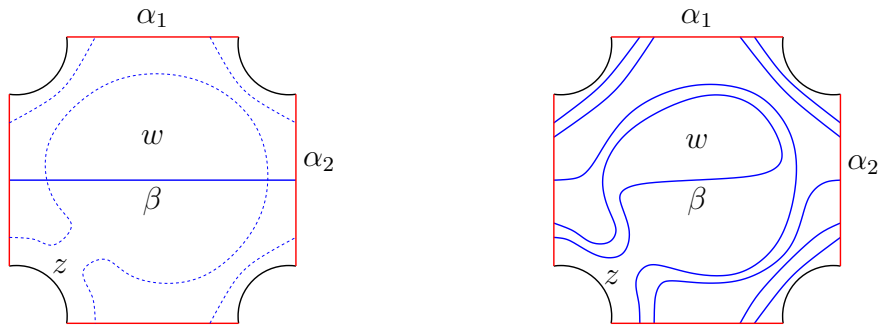}
\caption{\label{fig:UnbasedMappingClass} {\bf{Action of the
      unbased mapping class group.}}  The differential module
  corresponding to one idempotent can be thought of as the type $A$
  module associated to the doubly-pointed Heegaard diagram illustrated
  on the left. The Dehn twist along the dotted curve (which can be
  thought of as Dehn twist around the disk in the torus minus a disk)
  acts to give the diagram on the right. Since there are no
  differentials, the rank of the
  homology of the 
  resulting module is $9$.
  (The diagram depicts a left-handed Dehn twist applied to the $\beta$-circle, which
  corresponds to the action of a right-handed Dehn twist on the bordered three-manifold.)}
\end{figure}
\end{remark}

\begin{remark}
  It is possible to refine Theorem~\ref{thm:MCG-2-functoriality} by
  allowing triangulated actions for different grading sets beyond just
  $\sG(\PMC)$.  For
  instance, consider the category $\setModCat_{\Alg(\PMC),*}$ defined to
  be the disjoint union, over all right $\sG(\PMC)$-sets~$S$ on which
  $\lambda$ acts freely, of $\setModCat_{\Alg(\PMC),S}$.  Here the
  ``disjoint union'' of categories means the category whose objects
  are the union of the objects of the summands, with no morphisms
  between objects from different summands, and with the inherited
  morphisms between objects from the same summand.  Then each
  $\setModCat_{\Alg(\PMC),S}$ and therefore
  $\setModCat_{\Alg(\PMC),*}$ are \dg categories in the usual sense, and
  $\MCG_0(k)$ acts by triangulated functors on $\{\HMod(\setModCat_{\Alg(\PMC),*})\}$.
\end{remark}

\begin{remark}
  Group actions on algebraic categories have seen considerable
  interest recently; see for instance~\cite{GanterKapronov08:2reps},
  as well as the references in~\cite{KhTomas07:CobordismsCategories}.
  The reader might also wonder about group-ish structures with more
  interesting $2$-morphisms; such objects are called \emph{2-groups}
  and are studied in~\cite{BaezLauda04:2Groups}.
\end{remark}

\begin{remark}
  One could try to strengthen Theorem~\ref{thm:MCG-2-functoriality} as
  follows. Consider the $2$-category with objects pointed matched
  circles, $\Mor(\PMC,\PMC')$ the set of diffeomorphisms
  $F^\circ(\PMC)\to F^\circ(\PMC')$, and $2$-morphisms
  $\TMor(\phi,\psi)$ the set of isotopy classes of isotopies from $\phi$ to $\psi$. Then
  one could try to associate a weak $2$-functor from this category to
  the weak $2$-category of algebras, bimodules and homotopy classes of $\Ainf$-bimodule
  maps.  That is, one would associate a well-defined bimodule to each
  diffeomorphism and an $\Ainf$-homotopy equivalence of bimodules (well-defined up
  to homotopy) to
  each (isotopy class of) isotopy, satisfying appropriate coherence
  axioms. (One could, of course, imagine going farther, and looking for a functor
  between the $\infty$-category of surfaces, diffeomorphisms, paths of
  diffeomorphisms, \dots
  and the $\infty$-category of differential algebras, differential
  bimodules, differential bimodule homomorphisms, \dots.) Although
  this approach would lead to a slightly stronger result, it would
  require additional technicalities (for example, consistent choices
  of perturbation data).
\end{remark}

\subsection{Construction of the mapping class group action}
\label{sec:MappingClassGroupAction}

Before constructing the mapping class group action, we will need
to study the $DD$-identity bimodule.

We will need to consider simultaneously $\PMC$ and $-\PMC$.
To this end, recall from Section~\ref{sec:Gradings} that there is an
orientation-reversing map
$$r\co \PMC\to -\PMC,$$
which induces the map $a_i\mapsto a_{4k+1-i}$ on the points in~$\PMC$
(with respect to their orderings induced by the orientations on $Z$
and $-Z$). This induces a map from $[4k]/M$ to $[4k]/r(M)$.

\begin{definition}
  Let $\PMC$ be a pointed matched circle. If $\SetS$ and $r(\SetT)$
  form a partition of $[4k]/M$, we say that the idempotents
  $I_{\Alg(\PMC)}(\SetS)$ and $I_{\Alg(-\PMC)}(\SetT)$
  for $\Alg(\PMC)$ and $\Alg(-\PMC)$ are {\em complementary
    idempotents}.
\end{definition}

Complementary idempotents show up in the generating set for $\CFDDa$
for the standard Heegaard diagram for the identity map, as follows:

\begin{lemma}
  \label{lem:GeneratorsIdentityDD}
  Consider the standard Heegaard diagram $\HD$ for the identity map
  pictured in Figure~\ref{fig:id-diagram-disk}.  The generating set
  $\Gen(\HD)$ is in one-to-one correspondence with the set of
  idempotents $\Alg(\PMC)$: indeed, for each pair $(I,I')$ of
  complementary idempotents, there is a unique generator $\x=\x(I)$
  satisfying $(I\otimes I')\cdot \x=\x$. 
\end{lemma}

\begin{proof}
  This follows from a straightforward inspection of the diagram.
\end{proof}

We turn now to gradings on the identity type $DD$ bimodule.
The map 
$${R}\co \smallGroup(-\PMC)\to \smallGroup(\PMC)^{\op}$$
defined by ${R}(s,\eta)=(s,r_*(\eta))$ is a group
isomorphism.
Using this, we give $\smallGroup(\PMC)$ the structure of a left 
$\smallGroup(\PMC)\times_{\ZZ}\smallGroup(-\PMC)$-set by the rule
\[(g_1\times_{\ZZ} g_2)* h \coloneqq g_1\cdot h \cdot {R}(g_2)\]
(where the operation $\cdot$ on the right-hand-side refers
to multiplication in $\smallGroup(\PMC)$).
When referring to $\smallGroup(\PMC)$ as a
$\smallGroup(\PMC)\times_{\ZZ}\smallGroup(-\PMC)$-set in this way,
we denote it by $T$.

\begin{lemma}
  \label{lem:GradingIdentityDD}
  For $\HD$ the standard Heegaard diagram for the identity map,
  there is a natural identification of the grading set
  $\smallGrSet_\DD(\HD)$ of $\CFDDa(\HD)$ with the 
  $\smallGroup(\PMC)\times_{\ZZ}\smallGroup(-\PMC)$-set $T$.
\end{lemma}

\begin{proof}
  Recall that $P_{\x_0}$ is the image of the space of periodic domains under the map from
  Equation~\eqref{eq:DefineGPrimed}. 
  $\smallGrSet_\DD(\HD)$ can be defined as the quotient of
  $\smallGroup(\PMC)\times_{\ZZ}\smallGroup(-\PMC)$ by $(R\times R)(P_{\x_0})$.
  Another glance at the standard
  Heegaard diagram shows that $P_{\x_0}$ is generated by
  $(0,-r_*(m))\times_{\ZZ}(0,m)$, as $m$ runs over intervals in $\PMC$
  connecting matched pairs.  The map
  $$F\co \smallGroup(\PMC)\times_{\ZZ}\smallGroup(-\PMC)
  \to T$$
  defined by $F(g_1 \times g_2)=g_1\cdot R(g_2)$ induces an isomorphism of
  $\smallGroup(\PMC)\times_{\ZZ}\smallGroup(-\PMC)$-spaces
  \[
  f\co \bigl(\smallGroup(\PMC)\times_{\ZZ}\smallGroup(-\PMC)\bigr)\big/\bigl((R \times R)(P_{\x_0})\bigr)
  \to T.\qedhere
  \]
\end{proof}

\begin{lemma}
  \label{lem:RigidDD}
  Let $\PMC$ be a pointed matched circle. Any grading-preserving 
  homotopy auto-equivalence  of
  $\lsupv{\AlgA{\PMC},\AlgA{-\PMC}}\CFDDa(\Id)$
  as a type \DD\ bimodule is homotopic to the identity.
\end{lemma}

\begin{proof}
  It suffices to verify that for the standard Heegaard diagram $\HD$ of the
  identity map, pictured in Figure~\ref{fig:id-diagram-disk},
  the identity map is the only grading-preserving automorphism
  $\phi$ from  $\CFDDa(\HD)$ to $\CFDDa(\HD)$.

  It follows from Lemma~\ref{lem:GeneratorsIdentityDD} that our
  automorphism $\phi$ is specified by algebra elements
  $a(I,J)\in\AlgA{\PMC}$ and $b(I,J)\in\AlgA{-\PMC}$,
  indexed by pairs $I$ and $J$ of primitive idempotents for $\Alg(\PMC)$,
  whose compatibility with the idempotents is given by
  \[
  \xymatrix{I\cdot a(I,J)\cdot J=a(I,J) &
    I'\cdot b(I,J)\cdot J'=b(I,J)}
  \]
  (where $I'$ is complementary to~$I$ and $J'$ is complementary
  to~$J$), so that
  \[
  \phi(\x(I))=\sum_{J} (a(I,J)
  \otimes b(I,J))\otimes \x(J).
  \]

  To draw conclusions, we must turn to gradings.
  To this end, note that
  any two generators $\x(I)$ and $\x(J)$ can be connected by a domain $B$ with
  $e(B)+n_{\x(I)}(B)+n_{\x(J)}(B)=0$, and
  $\partial^{\partial_R}(B)=r_*(\partial^{\partial_L}(B))$. Thus, if we fix $I$,
  we can find one-chains $\alpha_J$ for each idempotent $J$, with 
  the property that
  $$\gr'(\x(J))=\bigl((0,\alpha_J)\times_{\ZZ}(0,r_*(\alpha_J))\bigr)*
    \gr'(\x(I)).$$
  Let $\psi$ and
  $\psi'$ be the grading refinement data for $\Alg(\PMC)$ and
  $\Alg(-\PMC)$, respectively.  Recall that
  \[
  \gr((a(I,J)\otimes b(I,J))\otimes \x(J))= (\psi(I)\times_{\ZZ}\psi'(I')) \cdot
  \gr'(a(I,J)\otimes b(I,J)) *\gr'(\x(J)).
  \]
  Thus, if $\gr((a(I,J)\otimes b(I,J))\otimes\x(J))=\gr(\x(I))$, then it follows that
  $$\gr'(a(I,J)\otimes b(I,J))\cdot \gr'(\x(J))=\gr'(\x(I)),$$
  so that
  \[
  \gr'(a(I,J)\otimes b(I,J))*((0,\alpha_J)\times_{\ZZ}(0,r_*(\alpha_J)))*\gr'(\x(I))=\gr'(\x(I)).
  \]
  Hence (according to Lemma~\ref{lem:GradingIdentityDD})
  there is another one-chain $\beta$ with the property that
  \[
  \gr'(a(I,J)\otimes b(I,J))\cdot \left((0,\alpha_J)\times_{\ZZ}(0,r_*(\alpha_J))\right) =
  (0,\beta)\times_{\ZZ}(0,r_*(\beta)).
  \]
  Using $m_{\PMC}(\alpha_J, \bdy \beta) =
  -m_{-\PMC}(r_*(\alpha_J), \bdy r_*(\beta))$, it follows that
  $\gr'(a(I,J)\otimes b(I,J))=(0,\beta-\alpha_J)\times_\ZZ\allowbreak(0,r_*(\beta-\alpha_J))$.
  From
  Lemma~\ref{lem:NegativeGradingsOnAlgebra}, it follows that each of $a(I,J)$ and
  $b(I,J)$ is an idempotent.

  Thus, $\phi$ has the form
  $$\phi(\x(I))=c(I,I)\otimes \x(I),$$
  where here $c(I,I)$ can be either $0$ or $1$. The fact that
  $\phi$ is an automorphism implies that each of these terms is $1$;
  i.e., $\phi$ is the identity map.
\end{proof}

\begin{lemma}\label{lem:bimodule-rigidity} Let $\HD$ and $\HD'$ be
  Heegaard diagrams representing
  $\phi\in\MCG_0(F(\PMC),\allowbreak F(\PMC'))$. Let
  $f,g\co\CFDAa(\HD)\to\CFDAa(\HD')$ be graded $\Ainf$-homotopy
  equivalences of bimodules. Then $f$ and $g$ are $\Ainf$-homotopic.
\end{lemma}

\begin{proof}
  We start by considering the identity map
  $\Id_{F(\PMC)}\in\MCG_0(F(\PMC),F(\PMC))$, with the Heegaard diagrams
  $\HD=\HD'$. By Lemma~\ref{lem:RigidDD}, $\CFDDa(\HD)$ has only one
  homotopy class of auto-equivalences. By
  Corollary~\ref{Cor:MCG-Equivalences}, $\CFAAa(\Id)\DT\cdot$ is an
  quasi-equivalence of categories from the category of type \DD\ 
  structures to the category of type DA structures, and by
  Theorem~\ref{thm:GenComposition} $\CFAAa(\Id)\DT\CFDDa(\HD)\simeq
  \CFDAa(\HD)$.  Thus, $\CFDAa(\HD)$ also has only one homotopy class
  of auto-equivalences.

  The case where $\HD\neq \HD'$ follows readily: $\CFDAa(\HD)$ and
  $\CFDAa(\HD')$ are homotopy equivalent, so the set of equivalences
  from $\CFDAa(\HD)$ to $\CFDAa(\HD')$ is identified with the set of
  homotopy auto-equivalences of $\CFDAa(\HD)$.  Finally, the case
  where $\phi\neq\Id$ follows, since $\CFDAa(\phi)$ is carried by the
  equivalence of categories given by $\CFDAa(\phi^{-1})\DT\cdot$ to
  $\CFDAa(\Id)$ (again, by Theorem~\ref{thm:GenComposition}).
\end{proof}

\begin{proof}[Proof of Theorem~\ref{thm:MCG-2-functoriality}]
For each mapping class $[\phi]\in\MCG_0(F(\PMC),F(\PMC'))$ choose a diffeomorphism
$\phi\co F(\PMC)\to F(\PMC')$ representing $[\phi]$. Recall from
Section~\ref{sec:AutBimodules} that there is a canonical bordered Heegaard diagram
$\HD(\phi)$ associated to $\phi$. For each  $\phi$
choose also a generic almost complex structure $J_\phi$ on
$\HD(\phi)\times[0,1]\times\RR$. Associated to these choices is a
well-defined bimodule $\CFDAa(\phi)=\CFDAa(\HD(\phi),J_\phi)$. Let
$\mathcal{F}_{\phi}$ denote the functor 
\[
\cdot\DT\CFDAa(\phi)\co \mathsf{H}(\ModCat_{\AlgA{\PMC}})\to \mathsf{H}(\ModCat_{\AlgA{\PMC'}}).
\]

Now, we define
an action~$A$ of $\MCG_0(k)$ on
$\{\mathsf{H}(\ModCat_{\AlgA{\PMC}})\mid \text{genus}(F(\PMC))=k\}$, as
in Definition~\ref{def:groupoid-action}, as follows. First, let
\begin{align*}
  A(\PMC)&= \HMod(\ModCat_{\AlgA{\PMC}}) \\
  A_{\PMC,\PMC'}([\phi])&= \mathcal{F}_\phi.
\end{align*}
It remains to define the correction terms $A_a$ and $A_{a,b,c}$ of
Definition~\ref{def:2-functor}.

By Theorem~\ref{thm:Id-is-Id} and Lemma~\ref{lem:bimodule-rigidity},
given a pointed matched circle $\PMC$ there is a unique graded
isomorphism $\CFDAa(\Id_{F(\PMC)})\to
\lsupv{\AlgA{\PMC}}[\Id]_{\AlgA{\PMC}}$.
Together with the canonical
natural isomorphism between the identity functor and
$\cdot\DT{}^{\AlgA{\PMC}}[\Id]_{\AlgA{\PMC}}$ 
(see Lemma~\ref{lem:IdentityDAIsIdentity})
this defines a
natural transformation
\[
A_\PMC\co \Id_{\mathsf{H}(\ModCat_{\AlgA{\PMC}})}\to \mathcal{F}_{\Id_{F(\PMC)}}.
\]

Similarly, given $[\phi_{12}]\in\MCG_0(\AlgA{\PMC_1},\AlgA{\PMC_2})$
and $[\phi_{23}]\in \MCG_0(\AlgA{\PMC_2},\AlgA{\PMC_3})$, let $\phi_{13}$
denote the chosen representative of
$[\phi_{23}\circ\phi_{12}]\in\MCG(\AlgA{\PMC_1},\AlgA{\PMC_3})$. Then,
by the pairing theorem (Theorem~\ref{thm:GenComposition}) and
Lemma~\ref{lem:bimodule-rigidity} there is a unique isomorphism
(in the derived category)
\[
\Phi_{123}\co \CFDAa(\phi_{12})\DT\CFDAa(\phi_{23})\to\CFDAa(\phi_{13}).
\]
This isomorphism induces a natural transformation
\[
A_{\PMC_1,\PMC_2,\PMC_3}([\phi_{23}], [\phi_{12}])\co\mathcal{F}_{\phi_{23}}\circ\mathcal{F}_{\phi_{12}}\to \mathcal{F}_{\phi_{13}}.
\]

This completes the definition of the weak $2$-functor $A$. It remains
to check that $A$ satisfies the three commutative diagrams of
Definition~\ref{def:2-functor}. But these follow trivially from
Lemma~\ref{lem:bimodule-rigidity}: both paths around the diagram
correspond to graded isomorphisms between a bimodule $M$ and a
bimodule of the form $\CFDAa(\phi)$, and that lemma guarantees that
there is a unique such isomorphism, so the diagrams must commute.

Finally, it is immediate from the definitions that these functors
preserve the subcategories %y
$\HMod(\setModCat_{\Alg(\PMC),\sG(\PMC)})$,
and act by triangulated functors on them.
\end{proof}

\begin{remark}\label{rmk:mcg-grading-action}
  We have constructed an action of the mapping class group action on
  the module category.  The reader might be interested in its
  behavior on grading sets.  Specifically, if
  $\phi\in\MCG_0(F(\PMC),F(\PMC))$ is a strongly based mapping class,
  then one can show that the grading set for $\CFDAa(\phi,\spinc)$ is a
  left-right $\smallGroup(\PMC)$-$\smallGroup(\PMC)$-set, for which
  both the left and the right actions are simply transitive. In turn, it is
  easy to see that isomorphism classes of such
  $\smallGroup(\PMC)$-$\smallGroup(\PMC)$-sets are in one-to-one
  correspondence with outer automorphisms of $\smallGroup(\PMC)$
  (fixing $\lambda$). Thus, we obtain a representation
  $\MCG_0(F(\PMC),F(\PMC))\to \Out(\smallGroup(\PMC))$.
  Projecting onto $H_1(F(\PMC))$, of course, we get the induced
  representation of the mapping class group on homology. But in fact,
  the entire representation
  is also determined by its behavior on homology. For more on this,
  see~\cite{LOT4}; see also Lemma~\ref{lem:GradingsDehnTwist} for an
  explicit action in the case where $g=1$.
\end{remark}

%%% Local Variables: 
%%% mode: latex
%%% TeX-master: "Bimodules"
%%% End: 

\section{Duality}
\label{sec:Duality}

In this section we will deduce Theorem~\ref{thm:Duality} from the fact
that $\CFAAa(\Id)$ and $\CFDDa(\Id)$ are quasi-inverses to each other
(Definition~\ref{def:QuasiInvertible}) (which in turn follows from the
pairing theorem (Theorem~\ref{thm:Composition}) and the fact that
$\CFDAa(\Id)\simeq \lsup{\Alg}[\Id]_\Alg$
(Theorem~\ref{thm:Id-is-Id})) and the algebraic results from
Section~\ref{sec:algebra-modules} (particularly
subsection~\ref{sec:mod-hom}). We give a genus $1$ illustration of
Theorem~\ref{thm:Duality} in Section~\ref{subsec:ExampleDuality}.

To keep notation simple, fix a pointed matched circle $\PMC$ and let
$\Alg=\Alg(\PMC)$ and $\Blg=\Alg(-\PMC)$.

Recall from Section~\ref{sec:bar-cobar-modules} that given a type $D$
structure $\lsup{\Alg}M$ we let $\lsub{\Alg}M=\lsub{\Alg}\Alg_\Alg\DT
\lsup{\Alg}M$ denote the corresponding type $A$ module; and similarly
for bimodules.
\begin{lemma}\label{lemma:AADDquasi-equiv}Let $\Clg$ and $\Elg$ be $\Ainf$ algebras. For any
  bimodules $\lsubsup{\Elg}{\Blg}M$ and $\lsubsup{\Clg}{\Blg}N$, the
  ``tensoring with the identity map'' morphism (see
  Proposition~\ref{prop:Hom-tensor-id})
  \[
  \Mor^{\Blg}(\lsubsup{\Elg}{\Blg}M,\lsubsup{\Clg}{\Blg}N)\to
  \Mor_{\Alg}(\CFAAa(\Id)_{\Alg,\Blg}\DT\lsubsup{\Elg}{\Blg}M,\CFAAa(\Id)_{\Alg,\Blg}\DT\lsubsup{\Clg}{\Blg}N)
  \]
  is a quasi-isomorphism of $(\Clg,\Elg)$-bimodules.
\end{lemma}
\begin{proof}
  The fact that this map respects the bimodule structure is immediate
  from Proposition~\ref{prop:Hom-tensor-id-D}.

  As noted above, it is immediate from Theorems~\ref{thm:Id-is-Id} and~\ref{thm:GenComposition}
  that the bimodules
  $\lsup{\Blg,\Alg}\CFDDa(\Id)$ and $\CFAAa(\Id)_{\Alg,\Blg}$
  are quasi-inverses to each other, as in
  Definition~\ref{def:QuasiInvertible}. It follows from
  Lemma~\ref{lem:quasi-invert-gives-equiv-DA} that the functors
  \begin{align*}
    \CFAAa(\Id)_{\Alg,\Blg}\DT\lsub{\Blg}\cdot &\co\HMod_*(\lsupv{\Blg}\ModCat)\to \HMod_*(\ModCat_\Alg)\\\shortintertext{and}
    \cdot_{\Alg}\DT\lsup{\Alg,\Blg}\CFDDa(\Id) &\co\HMod_*(\ModCat_{\Alg})\to \HMod_*(\lsupv{\Blg}\ModCat)
  \end{align*}
  are inverse equivalences of categories. But, on the level of
  morphisms, this says exactly that the ``tensoring with the
  identity'' map
  \[
  \Mor^{\Blg}(\lsup{\Blg}M,\lsup{\Blg}N)\to
  \Mor_{\Alg}(\CFAAa(\Id)_{\Alg,\Blg}\DT\lsup{\Blg}M,\CFAAa(\Id)_{\Alg,\Blg}\DT\lsup{\Blg}N)
  \]
  is a quasi-isomorphism.
\end{proof}

\begin{proposition}\label{prop:DDAA-duality}
  There is a quasi-isomorphism of bimodules
  \[
  \Mor^{\Blg}(\lsubsup{\Alg}{\Blg}{}\CFDDa(\Id),\lsupv{\Blg}[\Id]_{\Blg})\simeq
  \CFAAa(\Id)_{\Alg,\Blg}.
  \]
\end{proposition}
\begin{proof}
    There are quasi-isomorphisms of $\smallGroup$-set graded
    $(\Blg,\Alg)$-bimodules
  \begin{align*}
    \Mor^{\Blg}(\lsubsup{\Alg}{\Blg}{}\CFDDa(\Id),\lsupv{\Blg}[\Id]_{\Blg})&\to
    \Mor_{\Alg}(\CFAAa(\Id)_{\Alg,\Blg}\DT_\Blg\lsubsup{\Alg}{\Blg}{}\CFDDa(\Id),\CFAAa(\Id)_{\Alg,\Blg})\\
    &\to
    \Mor_{\Alg}(\lsub{\Alg}\Alg_{\Alg},\CFAAa(\Id)_{\Alg,\Blg})\\
    &\to \CFAAa(\Id)_{\Alg,\Blg}.
  \end{align*}
  Here, the first map is induced by tensoring with the identity, and
  is a quasi-isomorphism by Lemma~\ref{lemma:AADDquasi-equiv}.  The
  second map is induced by naturality of $\Mor$ and is a
  quasi-isomorphism because $\CFDDa$ and $\CFAAa$ are
  quasi-inverses. The third map comes from
  Lemma~\ref{lem:Tautology}. Note that all of these maps respect the
  $\smallGroup$-set grading.
\end{proof}

We are now in a position to prove Theorem~\ref{thm:Duality}, which we
restate and generalize as follows:
\begin{theorem}
  \label{thm:Duality-precise}
  Let $Y$ be a bordered three-manifold with boundary parameterized by $F(\PMC)$.  Then
  \begin{align*}
    \CFAa(Y)_{\Alg(\PMC)}&\simeq
    \Mor_{\Alg(-\PMC)}(\lsub{\Alg(\PMC),\Alg(-\PMC)}\CFDDa(\Id),\lsub{\Alg(-\PMC)}\CFDa(Y))\\
    &\simeq \Mor^{\Alg(-\PMC)}\bigl(\lsubsup{\Alg(\PMC)}{\Alg(-\PMC)}{}\CFDDa(\Id),\lsup{\Alg(-\PMC)}\CFDa(Y)\bigr).
  \end{align*}

  Similarly, suppose $Y_{12}$ is a strongly bordered three-manifold
  with boundary parameterized by $F(\PMC_1)$ and $F(\PMC_2)$. Then
  \[
  \lsup{\Alg(-\PMC_1)}\CFDAa(Y_{12})_{\Alg(\PMC_2)}\simeq\Mor^{\Alg(-\PMC_2)}\bigl(\lsubsup{\Alg(\PMC_2)}{\Alg(-\PMC_2)}{}\CFDDa(\Id_{\PMC_2}),\lsup{\Alg(-\PMC_1),\Alg(-\PMC_2)}\CFDDa(Y_{12})\bigr).
  \]
\end{theorem}
\begin{proof}
  As before, let $\Alg=\Alg(\PMC)$ and $\Blg=\Alg(-\PMC)$.
  By Proposition~\ref{prop:mor-iso-equiv},
  \[
  \Mor^{\Blg}\bigl(\lsubsup{\Alg}{\Blg}{}\CFDDa(\Id),\lsup{\Blg}\CFDa(Y)\bigr)\simeq
  \Mor_{\Blg}(\lsub{\Alg,\Blg}\CFDDa(\Id),\lsub{\Blg}\CFDa(Y)).
  \]
  By definition (Formula~\eqref{eq:MorDT}), Proposition~\ref{prop:DDAA-duality},
  and Theorem~\ref{thm:GenReparameterization} respectively,
  \begin{align*}
    \Mor^{\Blg}\bigl(\lsubsup{\Alg}{\Blg}{}\CFDDa(\Id),\lsup{\Blg}\CFDa(Y)\bigr)&\simeq
    \Mor^{\Blg}\bigl(\lsubsup{\Alg}{\Blg}{}\CFDDa(\Id),\lsupv{\Blg}[\Id]_{\Blg}\bigr)\DT\lsup{\Blg}\CFDa(Y)\\
    &\simeq \CFAAa(\Id)_{\Alg,\Blg}\DT\lsup{\Blg}\CFDa(Y)\\
    &\simeq \CFAa(Y)_{\Alg}.
  \end{align*}
  These isomorphisms are all maps of $\smallGroup$-set
  graded modules. 
  
  The bimodule case is similar.
\end{proof}

\begin{remark}
  One might imagine that
  \[
  \CFAa(Y)_{\Alg}\simeq
  \Mor^{\Blg}(\lsup{\Alg,\Blg}\CFDDa(\Id),\lsup{\Blg}\CFDa(Y))\DT\lsub{\Alg}\Alg_\Alg.
  \]
  This is in fact the case. However, the most obvious analogue for
  bimodules is false. See \cite{LOTHomPair}.
\end{remark}

\begin{remark}
  Like tensoring with $\CFDDa(\Id)$ or $\CFAAa(\Id)$, the
  duality exchanges the actions of $\Alg(F,i)$ and
  $\Alg(-F,-i)$. (Since modules associated to manifolds with
  connected boundary are supported in $i=0$, this reversal of
  $\SpinC$-structures is invisible in Theorem~\ref{thm:Duality}.)
\end{remark}

Of course, the above discussion
is much more useful once one calculates $\CFDDa(\Id)$.

\begin{remark}
\label{rmk:CFDD-id}
Consider the type \DD\ structure $B$ which is generated by elements of
the form $ I(r)\otimes I(s)$, where here $r,s\subset
[2k]$ are complementary sets of subsets of $[2k]=[4k]/M$, for a given
pointed matched circle. Let $R$ be the set of Reeb chords for $\PMC$. For
$\rho\in R$, let
$a(\rho)\in\Alg(\PtdMatchCirc)$ denote the algebra element
associated to $\rho$, and let $a(-\rho)\in\Alg(-\PtdMatchCirc)$
denote the corresponding algebra element for the pointed matched
circle with the opposite orientation.
We endow $B$ with a differential which is given
by 
\[
\delta^1(I(r)\otimes I(s))=\!\!\!\sum_{\substack{\rho\in R\\
    I(r)a(\rho)=a(\rho) I(r')\\ (r',s')\text{ complementary}}}
\!\!\!
(a(\rho)\otimes a(-\rho))\otimes (I(r')\otimes I(s')).
\]
It is shown in~\cite{LOT4} that the above explicitly-defined bimodule
$B$ is in fact quasi\hyp isomorphic to the bimodule $\CFDDa(\Id)$.
(See also Proposition~\ref{prop:TypeDDTorus} for a verification of this in the 
genus one case.)
\end{remark}

It is interesting to compare the bimodule of Remark~\ref{rmk:CFDD-id}
with the dualizing modules
in the theory of Koszul algebras. See for example~\cite{Priddy70:Koszul}.

%%% Local Variables: 
%%% mode: latex
%%% TeX-master: "Bimodules"
%%% End: 

\section{Bimodules for the torus}
\label{sec:torus-calc}

In \cite[Appendix~\ref*{LOT:app:Bimodules}]{LOT1}, we stated various
bimodules for
the torus. These included the bimodules $\CFDDa$ and $\CFAAa$ for the
identity cobordism, and also the type \DA\  bimodules for Dehn twists
along generators for the mapping class group. In this section we
verify those claims. Bimodules for generators of the
mapping class groupoid of a surface with arbitrary genus, given using
a different mechanism, are given in~\cite{LOT4}. 

We use notation for the torus algebra from
Section~\ref{subsec:GenusOneAlgebra} (mostly
from~\cite{LOT1}).  In Section~\ref{subsec:AAId1}, we calculate the
\AAm\ and \DD\ bimodules for the identity map. (The same techniques
can be used to calculate the \DA\ bimodule of the identity map: we do
not bother with this, in view of Theorem~\ref{thm:Id-is-Id}, which
ensures that it is simply the identity bimodule. It is worth noting,
though, that a direct calculation along the lines of
Section~\ref{subsec:AAId1} would result in a different, though
quasi-isomorphic, bimodule.)  In Section~\ref{subsec:DA-mcg} we
calculate the type \DA\ bimodules which represent the mapping class
group generators in genus one. We conclude by giving an illustration
of the duality theorem (Theorem~\ref{thm:Duality}) in
Section~\ref{subsec:ExampleDuality}.

Note that we focus here on the case $\Alg=\Alg(F,0)$. This is because
$\Alg(F,-1)\cong \Field$, and $\Alg(F,1)$ is quasi-isomorphic to
$\Field$, so all quasi-invertible bimodules over these algebras are
quasi-isomorphic to $\Field$.

\subsection{Bimodules for the identity map of the torus}
\label{subsec:AAId1}

\begin{figure}
\begin{center}
  \begin{tikzpicture}[y=54pt,x=1in]
    \node at (0,3) (w1) {${\mathbf w}_1$} ;
    \node at (2,3) (z1) {${\mathbf z}_1$} ;
    \node at (1,2) (y)  {${\mathbf y}$} ;
    \node at (1,1) (x)  {${\mathbf x}$} ;
    \node at (0,0) (w2) {${\mathbf w}_2$} ;
    \node at (2,0) (z2) {${\mathbf z}_2$} ;
    \draw[->] (w1) to node[above,sloped] {\lab{\sigma_1}} (y) ;
    \draw[->] (z1) to node[above] {\lab{\rho_1}} (y) ;
    \draw[->] (y)  to node[above,sloped] {\lab{(\sigma_2,\rho_2)}} (x) ;
    \draw[->] (x)  to node[below,sloped] {\lab{\rho_3}} (w2) ;
    \draw[->] (x)  to node[below] {\lab{\sigma_3}} (z2) ;
    \draw[->] (w1) to node[below,sloped] {\lab{1+(\sigma_{12},\rho_{23})}} (w2) ;
    \draw[->] (z1) to node[above,sloped] {\lab{1+(\sigma_{23},\rho_{12})}} (z2) ;
    \draw[->] (w1) to[pos=0.4] node[below,sloped] {\lab{(\sigma_{12},\rho_2)}} (x) ;
    \draw[->] (y)  to[pos=0.6] node[above,sloped] {\lab{(\sigma_2,\rho_{23})}} (w2) ;
    \draw[->] (z1) to[pos=0.4] node[below,sloped] {\lab{(\sigma_2,\rho_{12})}} (x) ;
    \draw[->] (y)  to[pos=0.6] node[above,sloped] {\lab{(\sigma_{23},\rho_2)}} (z2) ;
    \draw[->] (w1) to[out=-125,in=145] (-0.25,-0.25) to[out=-35,in=-150] node[pos=0,below,sloped]
       {\lab{(\sigma_{123},\rho_2)+(\sigma_3,\sigma_2,\sigma_1,\rho_2)}} (z2) ;
    \draw[->] (z1) to[out=-55,in=35] (2.25,-0.25) to[out=-145,in=-30] node[pos=0,below,sloped]
       {\lab{(\sigma_{2},\rho_{123})}} (w2) ;
  \end{tikzpicture}
\end{center}
\caption {The type \AAm\  bimodule $\CFAAa(\Id,0)$. 
  The labels on the arrows indicate the $\Ainf$
  operations; for instance, the label $(\sigma_{23},\rho_2)$ on the
  arrow from $\y$ to~$\mathbf{z}_2$ means that $m(\y,\sigma_{23},\rho_2)$
  contains a term~$\mathbf{z}_2$.}
\label{fig:AAmodAns}
\end{figure}
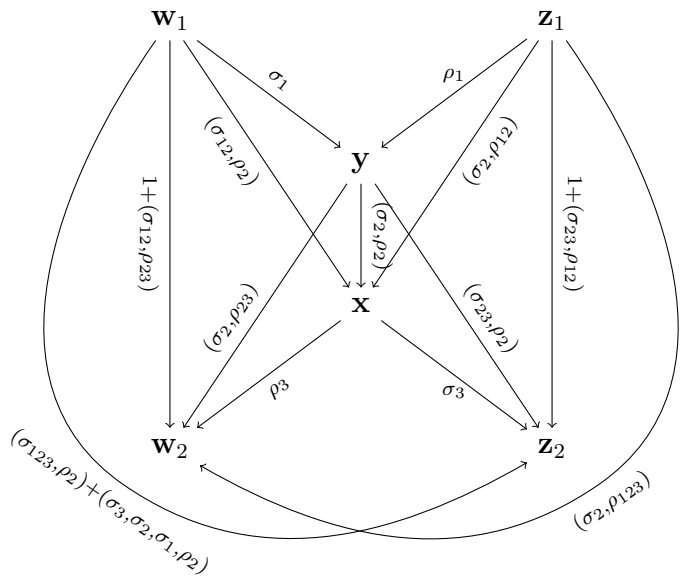

\begin{figure}
\begin{center}
\input{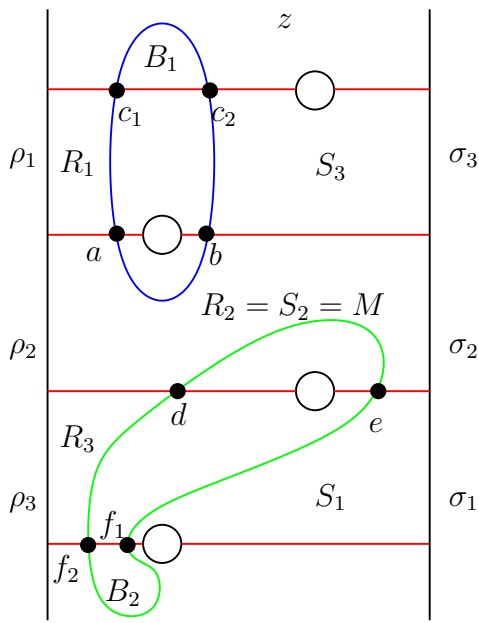}
\end{center}
\caption {Admissible Heegaard diagram for the identity map of the solid torus.
  There are two handles attached to the surface, connecting the
  circles that are vertically aligned with each other.}
\label{fig:HeegAA}
\end{figure}

Consider the unique pointed matched circle $\PMC$ for  a surface of
genus $1$. In the present section, we describe the bimodule
$\CFAAa(\Id,0)$, but first, we must set up some notation. Recall that
$\CFAAa(\Id,0)$ is a right-right
$\Alg(\PMC,0)$-$\Alg(-\PMC,0)$-bimodule. We write $\Alg=\Alg(\PMC,0)$,
and $\Blg=\Alg(-\PMC,0)$. Of course $\Alg\cong\Blg$, but we still find
it convenient to distinguish them, to help record which actions we are
using. Specifically, we think of $\Blg$ as having idempotents which we
denote $j_0$ and $j_1$ (corresponding to $\iota_0$ and $\iota_1$ in
$\Alg$), and generators $\sigma_i$ (corresponding to the $\rho_i$ in
$\Alg$).  With these conventions, then, we claim that the \AAm\
bimodule for the identity diffeomorphism of the torus (in the $i=0$
summand) is the \AAm\ bimodule illustrated in
Figure~\ref{fig:AAmodAns}.  The idempotent actions on the generators
can be easily determined from the actions on the algebra (or indeed
they can be read immediately off the Heegaard diagram in
Figure~\ref{fig:HeegAA}). We include only two of them here and leave
the others to the reader:
\[
\x\cdot (\iota_0\otimes j_0)=\x
\qquad
\y\cdot (\iota_1\otimes j_1)=\y
\]

Note that one can contract arrows to reduce to a quasi-isomorphic
bimodule with only two generators; however, in that model, the $\Ainf$
operations look rather complicated. (In particular, there are
infinitely many different ones, i.e., the module is not operationally
bounded in the sense of Definition~\ref{def:AA-bounded}.)

We verify our description of the \AAm\ bimodule for the identity by
drawing a suitable Heegaard diagram and analyzing the holomorphic
curves.  Unfortunately, the canonical Heegaard diagram for the
identity map given by Definition~\ref{def:ConstructHeegaardDiagram} is not
admissible.  This means that, although there are relatively few
generators, we could have infinitely many non-trivial $\Ainf$-products
(and hence infinitely many domains to consider). To simplify matters,
then, we apply some finger moves to get an admissible diagram $\HD$,
as illustrated in Figure~\ref{fig:HeegAA}.

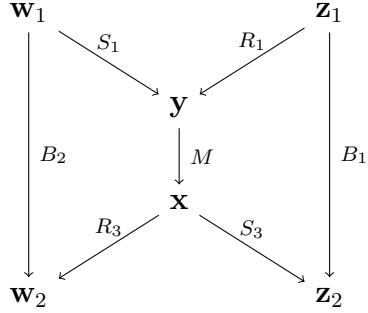
\begin{figure}
\begin{center}
  \begin{tikzpicture}[y=36pt,x=2cm]
    \node at (0,3) (BF1) {${\mathbf w}_1$} ;
    \node at (2,3) (C1E) {${\mathbf z}_1$} ;
    \node at (1,2) (AE)  {${\mathbf y}$} ;
    \node at (1,1) (BD)  {${\mathbf x}$} ;
    \node at (0,0) (BF2) {${\mathbf w}_2$} ;
    \node at (2,0) (C2E) {${\mathbf z}_2$} ;
    \draw [->] (BF1) to node[above] {\lab{S_1}} (AE) ;
    \draw [->] (C1E) to node[above] {\lab{R_1}} (AE) ;
    \draw [->] (AE) to node[auto] {\lab{M}} (BD) ;
    \draw [->] (BD) to node[above] {\lab{R_3}} (BF2) ;
    \draw [->] (BD) to node[above] {\lab{S_3}} (C2E) ;
    \draw [->] (BF1) to node[auto] {\lab{B_2}} (BF2) ;
    \draw [->] (C1E) to node[auto] {\lab{B_1}} (C2E) ;
  \end{tikzpicture}
\end{center}
\caption {Graph of domains.}
\label{fig:AAgraph}
\end{figure}

The diagram $\HD$ has six generators: 
${\mathbf y}=ae$, 
${\mathbf x}=bd$, 
${\mathbf w}_1=bf_1$,
${\mathbf w}_2=bf_2$, 
${\mathbf z}_1=c_1 e$, and 
${\mathbf z}_2=c_2 e$.  These can be connected by
the domains labelled in Figure~\ref{fig:HeegAA} to form the graph
shown in Figure~\ref{fig:AAgraph}.

This graph has the property that any homotopy class connecting two
generators can be realized as a path in the graph. Indeed,
any positive domain connecting two generators can be
represented by a connected path of downward-pointing edges. For
example, there are two positive domains from ${\mathbf w}_1$ to
${\mathbf w}_2$, and these are $B_2$ and the composite domain $S_1 M
R_3$.

The regions corresponding to the various domains from
Figure~\ref{fig:AAgraph} represent polygons. It follows that all the
$\Ainf$ operations from Figure~\ref{fig:AAmodAns} which drop height by
one are given by the corresponding labels; i.e.,
\begin{align*}
m({\mathbf w}_1,\sigma_1)&={\mathbf y} \\
m({\mathbf z}_1,\rho_1)&={\mathbf y} \\
m({\mathbf y},\sigma_2,\rho_2)&={\mathbf x} \\
m({\mathbf x},\rho_3)&={\mathbf w}_2 \\
m({\mathbf x},\sigma_3)&={\mathbf z}_2 \\
m({\mathbf w}_1)&={\mathbf w}_2 \\
m({\mathbf z}_1)&={\mathbf z}_2.
\end{align*}

The $\Ainfty$ relation for a type \AAm\ bimodule applied to the
${\mathbf x}$ components of $m^2({\mathbf
  w}_1\otimes\sigma_1\otimes\sigma_2\otimes\rho_2)$
and $m^2({\mathbf z}_1\otimes\sigma_2\otimes\rho_1\otimes\rho_2)$
forces the arrows from ${\mathbf w}_1$ and
${\mathbf z}_1$ to ${\mathbf x}$:
\begin{align*}
m({\mathbf w}_1,\sigma_{12},\rho_2)&={\mathbf x} \\
m({\mathbf z}_1,\sigma_2,\rho_{12})&={\mathbf x}.
\end{align*}
Similar considerations force the arrows from ${\mathbf y}$
to~${\mathbf w}_2$ and~${\mathbf z}_2$, from ${\mathbf z}_1$ to~${\mathbf z}_2$, and from ${\mathbf w}_1$ to~${\mathbf w}_2$:
\begin{align*}
  m({\mathbf y},\sigma_2,\rho_{23}) &= {\mathbf w}_2\\
  m({\mathbf y},\sigma_{23},\rho_2) &= {\mathbf z}_2\\
  m({\mathbf w}_1,\sigma_{12},\rho_{23}) &= {\mathbf w}_2\\
  m({\mathbf z}_1,\sigma_{23},\rho_{12}) &= {\mathbf z}_2.
\end{align*}
Considering the ${\mathbf w}_2$ and~${\mathbf z}_2$ components of $m^2({\mathbf z}_1\otimes\sigma_2\otimes\rho_1\otimes\rho_{23})$ and
$m^2({\mathbf w}_1\otimes \sigma_1\otimes\sigma_{23}\otimes\rho_2)$
respectively, we see that
\begin{align*}
m({\mathbf z}_1,\sigma_2,\rho_{123})&={\mathbf w}_2 \\
m({\mathbf w}_1,\sigma_{123},\rho_{2})&={\mathbf z}_2.
\end{align*}

We must finally consider the possibility of alternative higher
multiplications supported by the positive domains we have found.  For
example, the domain $S_1 M S_3$, as we have already seen,
represents $m({\mathbf w}_1,\sigma_{2},\rho_{123})={\mathbf z}_2$,
but has an alternative decomposition, where we cut to the right wherever
possible. This gives rise to a polygon representing
$$m({\mathbf w}_1,\sigma_3,\sigma_2,\sigma_1,\rho_2)={\mathbf z}_2.$$

Inspecting the existing positive domains, it is clear that no other
one can give rise to an alternate $\Ainf$ operation, as
corresponding cuts are not possible.

This concludes the verification that the identity map gives rise to
the bimodule pictured in Figure~\ref{fig:AAmodAns}. 
Having found all the
holomorphic curves for the \AAm\  bimodule, it is easy to write down
also the type \DD\ bimodule. To do this, we proceed as follows:
we relabel the algebra elements on the regions to be compatible
with type $D$ labelings (i.e., swap the order of the generators on the boundary),
then throw out some curves which cannot contribute for idempotent
reasons, and finally add up other curve counts. 

\begin{figure}
\begin{center}
  \begin{tikzpicture}[y=54pt,x=1in]
    \node at (0,3) (w1) {${\mathbf w}_1$} ;
    \node at (2,3) (z1) {${\mathbf z}_1$} ;
    \node at (1,2) (y)  {${\mathbf y}$} ;
    \node at (1,1) (x)  {${\mathbf x}$} ;
    \node at (0,0) (w2) {${\mathbf w}_2$} ;
    \node at (2,0) (z2) {${\mathbf z}_2$} ;
    \draw[->] (w1) to node[above] {\lab{\sigma_3}} (y) ;
    \draw[->] (z1) to node[above] {\lab{\rho_3}} (y) ;
    \draw[->] (y)  to node[right] {\lab{\sigma_2 \rho_2}} (x) ;
    \draw[->] (x)  to node[below] {\lab{\rho_1}} (w2) ;
    \draw[->] (x)  to node[below] {\lab{\sigma_1}} (z2) ;
    \draw[->] (w1) to node[left] {\lab{1}} (w2) ;
    \draw[->] (z1) to node[right] {\lab{1}} (z2) ;
    \draw[->] (z1) to[out=-55,in=35] (2.25,-0.25) to[out=-145,in=-30] node[pos=0,below,sloped]
       {\lab{\sigma_{2}\rho_{123}}} (w2) ;
  \end{tikzpicture}
\end{center}
\caption {The type \DD\ bimodule $\CFDDa(\Id,0)$. 
  The labels on the arrows indicate differentials; 
  for instance, the arrow from ${\mathbf y}$ to ${\mathbf x}$ signifies a term
  of $(\rho_2 \sigma_2)\otimes {\mathbf x}$ in $\partial {\mathbf y}$.}
\label{fig:DDmodAns}
\end{figure}
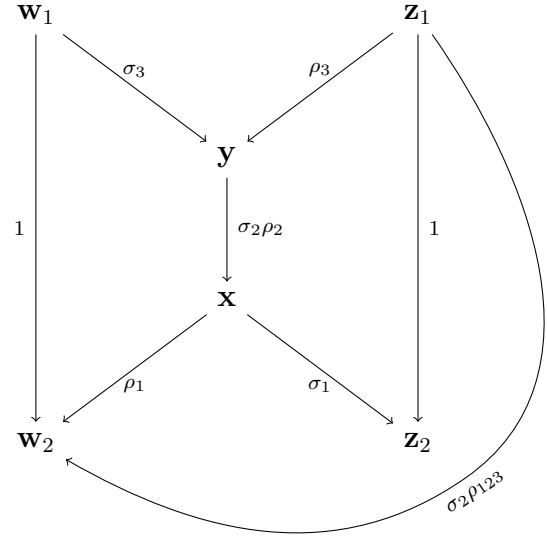

For instance, a curve which used to count as $m({\mathbf
  w}_1,\sigma_1)={\mathbf y}$ now counts as giving a term of
${\sigma_3}\otimes {\mathbf y}$ in $\partial {\mathbf w}_1$. Also, the
curve which used to count as $m({\mathbf
  w}_1,\sigma_{12},\rho_2)={\mathbf x}$ now does not count. Finally,
the contributions 
\[
m({\mathbf w}_1,\sigma_{123},\rho_2)={\mathbf
  z}_2 \qquad{\text{and}}\qquad m({\mathbf
  w}_1,\sigma_3,\sigma_2,\sigma_1,\rho_2)={\mathbf z}_2\]
both contribute to the coefficient of $(\sigma_{123} \rho_2)\otimes
{\mathbf z}_2$ in $\partial {\mathbf w}_1$ (and hence they cancel).
The results are summarized in Figure~\ref{fig:DDmodAns}.

We can simplify this further, to obtain the following:

\begin{proposition}
  \label{prop:TypeDDTorus}
  There is a type \DD\ identity bimodule for the torus
  with two generators $\mathbf p$ and $\mathbf q$, 
  satisfying 
  \[(\iota_0\otimes j_0)\cdot \mathbf p = \mathbf p
  \qquad 
  (\iota_1\otimes j_1)\cdot \mathbf q = \mathbf q \]
  and differential given by
  \[
  \partial \mathbf p = (\rho_1\sigma_3+ \rho_3\sigma_1 +
    \rho_{123}\sigma_{123}) \otimes \mathbf q \qquad
  \partial \mathbf q = (\rho_2\sigma_2)\otimes \mathbf p.
  \]
\end{proposition}
\begin{proof}
Substitute
\begin{align*}
\mathbf p&={\mathbf x}+\rho_1 \otimes {\mathbf w}_1+\sigma_1\otimes
  {\mathbf z}_1 + \sigma_{12}\rho_{123} \otimes {\mathbf w}_1
\\
\mathbf q&={\mathbf y}
\end{align*}
in the description from Figure~\ref{fig:DDmodAns}
to get a quasi-isomorphic submodule 
with the stated differential.
The idempotent actions follow immediately from the diagram.
(Note that the idempotents are the same as those arising from
the interpretation of the generators as generators for a type \AAm\ bimodule.)
\end{proof}

The type \DD\ module given in Figure~\ref{fig:DDmodAns} is
bounded in the sense of Definition~\ref{def:DD-bounded}, while the
module given in Proposition~\ref{prop:TypeDDTorus} is merely left and right
bounded.

Determining relative
gradings is also a straightforward matter, as in
Lemma~\ref{lem:GradingIdentityDD}.
%The octagon~$M$ connecting $\x$ and $\y$ demonstrates that 
%\[((0;0,1,0)\times_{\lambda}(0;0,1,0))\cdot \gr'(\mathbf p)=\gr'(\mathbf q).\]

\subsection{The \DA\ bimodules for the mapping class group of the torus}
\label{subsec:DA-mcg}

The mapping class group is generated by Dehn twists $\tau_\merid$ and
$\tau_\longitude$ along meridian and longitude respectively (i.e.,
$\tau_\merid$ takes an $n$-framed knot complement to an $n+1$-framed knot
complement). We describe the summand
$\CFDAa(\cdot,0)$ of the type \DA\  modules for Dehn twists about these
two curves and their inverses.  

\begin{figure}
\input{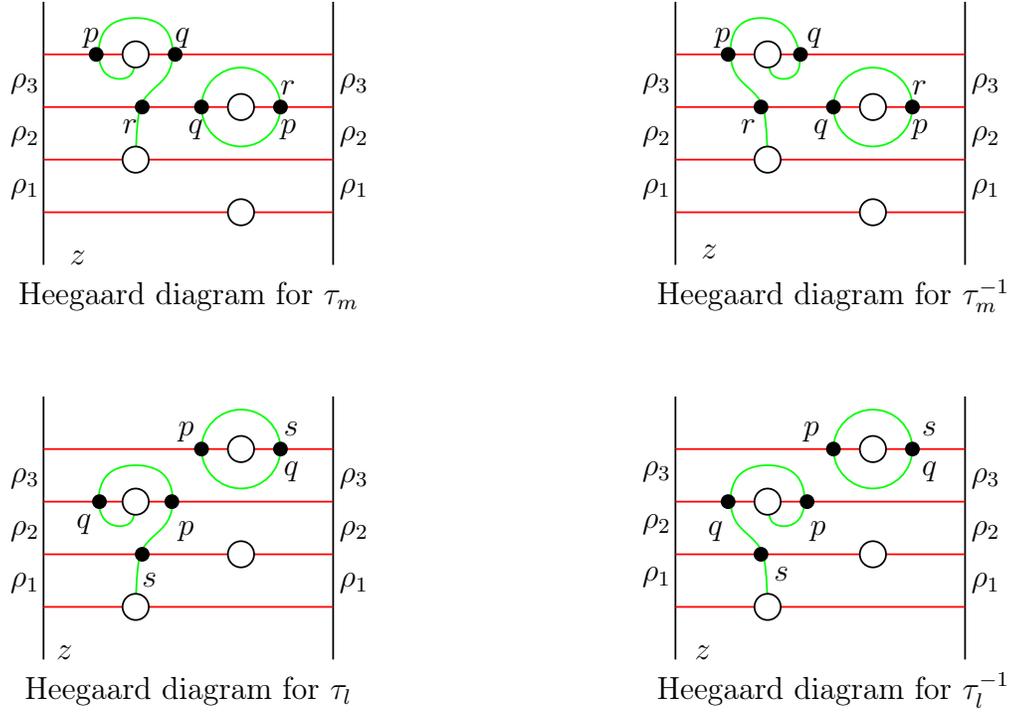}
\caption{\label{fig:DehnTwistsGenusOne} {\bf{Heegaard diagrams for
      mapping class group elements.}}  We have illustrated here genus
  two bordered diagrams for the generators of the genus one mapping
  class group and their inverses, as indicated.  In each of the four
  diagrams, there are three
  generators in the $i=0$ summand.}
\end{figure}

Heegaard diagrams for $\tau_\merid$, $\tau_\merid^{-1}$,
$\tau_\longitude$ and $\tau_\longitude^{-1}$ are illustrated in
Figure~\ref{fig:DehnTwistsGenusOne}.  Each of the four type \DA\
structures $\CFDAa(\tau_\merid,0)$, $\CFDAa(\tau_\merid^{-1},0)$,
$\CFDAa(\tau_\longitude,0)$, and $\CFDAa(\tau_\longitude^{-1},0)$ has
three generators, which are labelled by four possible letters
${\mathbf p}$, ${\mathbf q}$, ${\mathbf r}$, and ${\mathbf s}$.
($\CFDAa(\tau_\merid^{\pm 1})$ does not have ${\mathbf s}$ and
$\CFDAa(\tau_\longitude^{\pm 1})$ does not have ${\mathbf r}$.)  The
compatibility with the idempotents is given as follows:
\[
\iota_0\cdot {\mathbf p}\cdot \iota_0 = {\mathbf p} \qquad
\iota_1\cdot {\mathbf q}\cdot \iota_1 = {\mathbf q} \qquad
\iota_1\cdot {\mathbf r}\cdot \iota_0 = {\mathbf r} \qquad
\iota_0\cdot {\mathbf s}\cdot \iota_1 = {\mathbf s}.
\]

Next, we study the grading set for the modules. 
We will adopt the notation for elements of
$\bigGroup$ from \cite[Chapter~\ref*{LOT:chap:TorusBoundary}]{LOT1};
that is, elements of
$\bigGroup(T^2)$ are written as tuples $g=(m;i,j,k)$ where
$m\in\frac{1}{2}\ZZ$ is the Maslov component of $g$ and $i,j,k$
are the local multiplicities of $[g]\in H_1(\PMC\setminus
z,\mathbf{a})$ at the three relevant components of
$\PMC\setminus\mathbf{a}$. See
Section~\ref{subsec:GenusOneAlgebra}. Since $\smallGroup$ is a
subgroup of $\bigGroup$, this gives us notation for elements of
$\smallGroup$ as well.

\begin{definition}
  Let $f\co G\to G$ be a group homomorphism.
  Let $G_{f}$ be the associated left-right $G$-$G$-space
  whose underlying set is $G$, and whose action is given by
  $g_1 \star s * g_2 = g_1\cdot s\cdot f(g_2)$,
  where $\star$ denotes the left action on $G_f$,
  $*$ denotes the right action on $G_f$, and $\cdot$
  denotes multiplication in $G$.
\end{definition}

With this definition, if $G_{f}\cong G_g$
  as left-right $G$-$G$-spaces, then $f$ and $g$ are conjugate to
  one another.

Recall that the actions of the generating Dehn twists on homology is given by
\begin{align*}
  (\tau_\merid)_*(\merid)&=\merid & 
  (\tau_\merid)_*(\longitude)&=\merid+\longitude \\
  (\tau_\merid^{-1})_*(\merid)&=\merid & 
  (\tau_\merid^{-1})_*(\longitude)&=-\merid+\longitude \\
  (\tau_\longitude)_*(\merid)&=\merid-\longitude & 
  (\tau_\longitude)_*(\longitude)&=\longitude \\
  (\tau_\longitude^{-1})_*(\merid)&=\merid+\longitude & 
  (\tau_\longitude^{-1})_*(\longitude)&=\longitude.
\end{align*}
Here, we are thinking of the homology classes $\merid$ and $\longitude$
as represented by local multiplicities
\[\merid=(0,1,1) \qquad \longitude=(1,1,0).\]
These have canonical lifts to elements in $\smallGroup$, gotten by
$\Lmerid=\gr(\rho_{23})$ and $\Llongitude=\gr(\rho_{12})$; i.e.,
(using the grading refinement data from Equation~\eqref{eq:RefineTorus}),
we get
\[\Lmerid=(\OneHalf;0,1,1) \qquad \Llongitude=(-\OneHalf;1,1,0)\]
We define the following lifts of the action of homology to 
automorphisms of $\smallGroup(\PMC)$:
$f_\merid$, $f_\merid^{-1}$, $f_\longitude$, and
$f_{\longitude}^{-1}$.  These
are determined by the property that they fix $\lambda$, and
transform the other generators according to the following:
\begin{equation}
  \begin{aligned}
    \label{eq:LiftToGradingGroup}
  f^{\tau_\merid}(\Lmerid)&=\Lmerid & 
  f^{\tau_\merid}(\Llongitude)&=\lambda\cdot \Lmerid\cdot \Llongitude \\
  f^{\tau_\merid^{-1}}(\Lmerid)&=\Lmerid & 
  f^{\tau_\merid^{-1}}(\Llongitude)&=\lambda^{-1}\cdot \Lmerid^{-1}\cdot \Llongitude \\
  f^{\tau_\longitude}(\Lmerid)&=\lambda^{-1}\cdot\Lmerid\cdot \Llongitude^{-1} & 
  f^{\tau_\longitude}(\Llongitude)&=\Llongitude\\
  f^{\tau_\longitude^{-1}}(\Lmerid)&=\lambda\cdot \Lmerid\cdot \Llongitude& 
  f^{\tau_\longitude^{-1}}(\Llongitude)&=\Llongitude
  \end{aligned}
\end{equation}

\begin{lemma}
  \label{lem:GradingsDehnTwist}
  Let
  $\phi\in\{\tau_\merid,\tau_\merid^{-1},\tau_{\longitude},\tau_{\longitude}^{-1}\}$.
  For the grading refinement data of Equation~\eqref{eq:RefineTorus},
  if we base our grading sets around $\mathbf p$, the intersection point in the
  $\iota_0$ idempotent on both the left and the right,
  then the grading set for $\CFDAa(\phi)$ is identified with 
  $\smallGroup_{f^\phi}$, where $f^\phi$ is the homomorphism corresponding 
  to $\phi$ as described in Equations~\eqref{eq:LiftToGradingGroup}.

  With respect to these identifications of grading sets, we find
  \begin{align*}
    \gr_{\merid}({\mathbf p})&=(0;0,0,0)  
    & \gr_{\merid}({\mathbf q})&= (0;0,0,0)& 
    \gr_{\merid}({\mathbf r})&= (\OneHalf;1,1,0)\\
    \gr_{\merid^{-1}}({\mathbf p})&= (0;0,0,0)& 
    \gr_{\merid^{-1}}({\mathbf q})&= (0;0,0,0) &
    \gr_{\merid^{-1}}({\mathbf r})&=(-\OneHalf;-1,-1,0)\\
    \gr_{\longitude}({\mathbf p})&= (0;0,0,0) & 
    \gr_{\longitude}({\mathbf q}) &= (0;0,0,0) & 
    \gr_{\longitude}({\mathbf s}) &=(\OneHalf;1,1,0)\\
    \gr_{\longitude^{-1}}({\mathbf p})&= (0;0,0,0) & 
    \gr_{\longitude^{-1}}({\mathbf q}) &= (0;0,0,0) & 
    \gr_{\longitude^{-1}}({\mathbf s}) &=(-\OneHalf;-1,-1,0)\\
  \end{align*}
  where the subscript on $\gr$ indicates which diagram we are considering.
\end{lemma}

\begin{figure}
\input{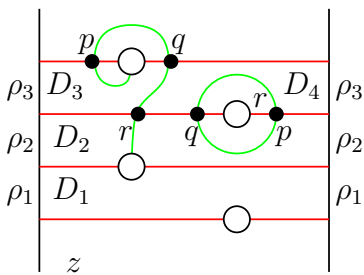}
\caption{\label{fig:DehnTwistMu} \textbf{Heegaard diagram for
      $\tau_\merid$, with labelled domains.}}
\end{figure}

\begin{proof}
  We describe in detail the calculations in the case where $\phi=\tau_\merid$
  (and consider $\gr=\gr_{\merid}$).  (The other cases are entirely parallel.)
  A basis for the space of periodic domains 
  $\pi_2({\mathbf r},{\mathbf r})$ is given
  by the domains $D_2+D_3+D_4$ and $D_1+D_2-D_4$, with diagrams labelled
  as in Figure~\ref{fig:DehnTwistMu}.

  Let $\PMC_L\amalg\PMC_R$ be the boundary of the Heegaard diagram for~$\phi$.
  The grading set for the larger ($\bigGroup$) grading 
  takes its values in the quotient of
  $\bigGroup(-\PMC_L)\times_{\lambda}\bigGroup(\PMC_R)$
  (viewed as a left-right $\bigGroup(-\PMC_L)$-$\bigGroup(\PMC_R)$-set)
  by relations coming from the periodic domains; specifically, 
  \begin{align*}
    \grb({\mathbf p})&=
    (\OneHalf;0,-1,-1)\cdot\grb({\mathbf p})\cdot (-\OneHalf;0,1,1)
    \\
    \grb({\mathbf p}) 
    &=
    (\OneHalf;-1,-1,0)\cdot\grb({\mathbf p})\cdot (-1;1,0,-1)
  \end{align*}
  The first of these equations comes from $D_2+D_3+D_4$, which has $e(D_2+D_3+D_4)=-2$ and
  $2n_{{\mathbf p}}(D_2+D_3+D_4)=2$, 
  $r_*(\bdy^{\bdy_L}(D_2+D_3+D_4))=(0,-1,-1)$,
  $\bdy^{\bdy_R}(D_2+D_3+D_4)=(0,1,1)$.
  The second equation comes from $D_1+D_2-D_4$,
  which has $e(D_1+D_2-D_4)=-\OneHalf$ and $2n_{{\mathbf p}}(D_1+D_2-D_4)=1$.)
  According to Equation~\eqref{eq:bimod-grading-refine}
  (see also Remark~\ref{rmk:RefineOtherModules}), we have
  $\gr({\mathbf p})=\psi(\iota_0)\cdot\grb({\mathbf p})\cdot
  \psi(\iota_0)^{-1}$;
  but since
  $\psi(\iota_0)$ is the identity, the above relations can be restated
  as:
  \begin{align*}
    \gr({\mathbf p})&=
    (\OneHalf;0,-1,-1)\cdot\gr({\mathbf p})\cdot (-\OneHalf;0,1,1)
    \\
    \gr({\mathbf p}) 
    &=
    (\OneHalf;-1,-1,0)\cdot\gr({\mathbf p})\cdot (-1;1,0,-1).
  \end{align*}
  It is now easy to see that the map from the grading
  set to $G_{f^{\tau_\merid}}$ given by
  \begin{equation}
    \label{eq:IdentifyGradingSet}
    g_1\cdot \gr({\mathbf p})\cdot g_2 \mapsto 
    g_1\cdot f^{\tau_\merid}(g_2)
  \end{equation}
  (which evidently sends $\gr({\mathbf p})$ to
  $(0;0,0,0)$)
  is an isomorphism of left-right $\bigGroup$-$\bigGroup$-sets.

  To calculate the grading of ${\mathbf r}$, for example, consider the
  domain $D_2$ from ${\mathbf r}$ to ${\mathbf p}$. 
  This domain has $e(D_2)=0$ and $n_{\mathbf r}+n_{\mathbf p}=1/2$, and hence
  it gives the relation 
  $\gr'({\mathbf p})=(-\OneHalf;0,-1,0)\star\gr'({\mathbf r})$. 
  Since $\gr({\mathbf r})=\psi(\iota_1)\cdot \gr'({\mathbf r})\cdot \psi(\iota_0)^{-1}$,
  it follows that 
  $$\gr({\mathbf r})=(\OneHalf;1,1,0).$$
  A similar calculation (only now using the domain $D_1$) shows
  $\gr({\mathbf p})=\gr({\mathbf q})$.
\end{proof}

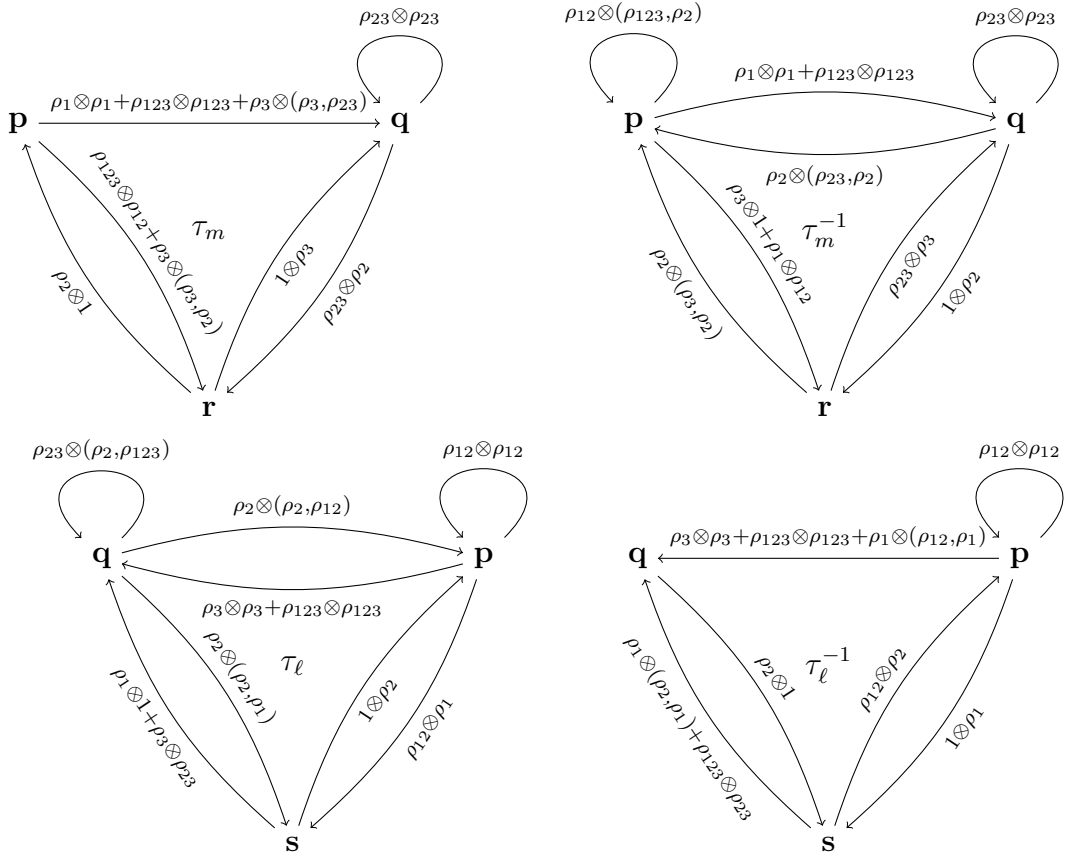
\begin{figure}
\begin{center}
  \begin{tikzpicture}[y=54pt,x=1in]
    \node at (0,2) (p) {${\mathbf p}$} ;
    \node at (2,2) (q) {${\mathbf q}$} ;
    \node at (1,0) (r) {${\mathbf r}$} ;
    \node at (1,1.25) (label) {$\tau_\merid$};
    \draw[->] (p) to node[above,sloped] {\lab{\rho_1\otimes\rho_1+\rho_{123}\otimes\rho_{123}+\rho_3\otimes (\rho_{3},\rho_{23})}}  (q)  ;
    \draw[->] (p) [bend left=15] to node[above,sloped] {\lab{\rho_{123}\otimes\rho_{12}+
\rho_3\otimes (\rho_{3},\rho_{2})}} (r) ;
    \draw[->] (q) [bend left=15] to node[below,sloped] {\lab{\rho_{23}\otimes \rho_{2}}} (r) ;
    \draw[->] (q) [loop] to node[above] {\lab{\rho_{23}\otimes\rho_{23}}}
                 (q) ;
    \draw[->] (r) [bend left=15] to node[below,sloped] {\lab{\rho_2\otimes 1}} (p) ;
    \draw[->] (r) [bend left=15] to node[below,sloped] {\lab{1\otimes \rho_3}} (q) ;
  \end{tikzpicture}
  \begin{tikzpicture}[y=54pt,x=1in]
    \node at (0,2) (p) {${\mathbf p}$} ;
    \node at (2,2) (q) {${\mathbf q}$} ;
    \node at (1,0) (r) {${\mathbf r}$} ;
    \node at (1,1.25) (label) {$\tau_\merid^{-1}$};
    \draw[->] (p) [bend left=15] to node[above,sloped] {\lab{\rho_1\otimes \rho_1+\rho_{123}\otimes\rho_{123}}}  (q)  ;
    \draw[->] (p) [loop] to node[above] {\lab{\rho_{12}\otimes(\rho_{123},\rho_2)}} (p) ;
    \draw[->] (p) [bend left=15] to node[above,sloped]  {\lab{\rho_3\otimes 1 + \rho_1\otimes \rho_{12}}} (r) ;
    \draw[->] (q) [bend left=15] to node[below,sloped] {\lab{\rho_2\otimes(\rho_{23},\rho_2)}} (p) ;
    \draw[->] (q) [bend left=15] to node[below,sloped] {\lab{1\otimes \rho_2}} 
                  (r) ;
    \draw[->] (q) [loop] to node[above] {\lab{\rho_{23}\otimes\rho_{23}}} (q) ;
    \draw[->] (r) [bend left=15] to node[below,sloped] {\lab{\rho_2\otimes (\rho_3,\rho_2)}} (p) ;
    \draw[->] (r) [bend left=15] to node[below,sloped] {\lab{\rho_{23}\otimes \rho_3}} (q) ;
  \end{tikzpicture} \\
  \begin{tikzpicture}[y=54pt,x=1in]
    \node at (0,2) (p) {${\mathbf q}$} ;
    \node at (2,2) (q) {${\mathbf p}$} ;
    \node at (1,0) (r) {${\mathbf s}$} ;
    \node at (1,1.25) (label) {$\tau_\longitude$};
    \draw[->] (p) [bend left=15] to node[above,sloped] {\lab{\rho_2\otimes(\rho_2,\rho_{12})}}  (q)  ;
    \draw[->] (p) [bend left=15] to node[above,sloped] {\lab{\rho_2\otimes (\rho_2,\rho_1)}} (r) ;
    \draw[->] (p) [loop] to node[above]  {\lab{\rho_{23}\otimes (\rho_2,\rho_{123})}}  (p) ;
    \draw[->] (q) [bend left=15] to node[below,sloped] {\lab{\rho_{12}\otimes\rho_1}} (r) ;
    \draw[->] (q) [loop] to node[above] {\lab{\rho_{12}\otimes\rho_{12}}} (q) ;
    \draw[->] (q) [bend left=15] to node[below,sloped] {\lab{\rho_3\otimes\rho_3+\rho_{123}\otimes\rho_{123}}} (p) ;
    \draw[->] (r) [bend left=15] to node[below,sloped] {\lab{\rho_1\otimes 1+\rho_3\otimes\rho_{23}}} (p) ;
    \draw[->] (r) [bend left=15] to node[below,sloped] {\lab{1\otimes \rho_2}} (q) ;
  \end{tikzpicture}
  \begin{tikzpicture}[y=54pt,x=1in]
    \node at (0,2) (p) {${\mathbf q}$} ;
    \node at (2,2) (q) {${\mathbf p}$} ;
    \node at (1,0) (r) {${\mathbf s}$} ;
    \node at (1,1.25) (label) {$\tau_\longitude^{-1}$};
    \draw[->] (p) [bend left=15] to node[above,sloped] {\lab{\rho_2\otimes 1}}  (r)  ;
    \draw[->] (q) [loop] to node[above] {\lab{\rho_{12}\otimes\rho_{12}}}
                (q) ;
    \draw[->] (q) [bend left=15] to node[below,sloped]  {\lab{1\otimes \rho_1}} (r) ;
    \draw[->] (q)  to node[above,sloped] {\lab{\rho_3\otimes\rho_3+\rho_{123}\otimes\rho_{123}+\rho_1\otimes(\rho_{12},\rho_1)}} (p) ;
    \draw[->] (r) [bend left=15] to node[below,sloped] {\lab{\rho_1\otimes (\rho_2,\rho_1)+\rho_{123}\otimes\rho_{23}}} (p) ;
    \draw[->] (r) [bend left=15] to node[above,sloped] {\lab{\rho_{12}\otimes\rho_2}} (q) ;
  \end{tikzpicture}
\end{center}
\caption {
\label{fig:DAbimodules}
Type \DA\ bimodules for torus mapping class group action.  These are
the module associated to $\tau_\merid$, $\tau_\merid^{-1}$,
$\tau_{\longitude}$, $\tau_{\longitude}^{_1}$ respectively. The
notation is as follows. Consider the module for $\tau_\merid$. The
label $\rho_1\otimes\rho_1$ on the horizontal arrow indicates that
$m(\mathbf{p},\rho_1)$ contains a term of the form
$\rho_1\mathbf{q}$. Similarly, the label
$\rho_3\otimes(\rho_3,\rho_{23})$ on that arrow indicates that
$m(\mathbf{p},\rho_3,\rho_{23})$ contains a term of the form
$\rho_3\mathbf{q}$.}
\end{figure}

Next we turn to computing the explicit bimodules.

\begin{proposition}
  The type \DA\ bimodules for $\tau_m$, $\tau_m^{-1}$, $\tau_\ell$,
  and $\tau_\ell^{-1}$ are as given in Figure~\ref{fig:DAbimodules}.
\end{proposition}

\begin{proof}
We give the proof in detail for $\CFDAa(\tau_m)$.
We enumerate domains which contribute to the type \DA\  actions,
organizing them by the algebra elements they contribute
on the type $D$ side. These algebra elements in turn
are determined by how the domain meets $\partial_L \HD$.

{\bf{Algebra element $1$.}} The only domain which is disjoint
from $\partial_L \HD$ is $D_4$. That represents a rectangle,
which therefore contributes to $\delta^1$. We denote this by writing
$$({\mathbf r},\rho_3)\stackrel{D_4}{\longrightarrow} {\mathbf q},$$
to mean that $\mathbf{q}$ occurs in $\delta^1(\mathbf{r},\rho_3)$.

{\bf{Algebra element $\rho_1$.}} There is only one valid domain which
could contribute $\rho_1$, namely the domain $D_1$ itself. Indeed, it
is a rectangle, starting at
${\mathbf p}$ and ending at ${\mathbf q}$. Therefore,
$\delta^1(\mathbf{p},\rho_1)$ contains $\rho_1\mathbf{q}$, or graphically
$$({\mathbf p},\rho_1)\stackrel{D_1}{\longrightarrow}\rho_1\otimes{\mathbf q}.$$

{\bf{Algebra element $\rho_2$.}} The only domain which could contribute
in this case is $D_2$; and that in turn is a bigon, giving
$$({\mathbf r})\stackrel{D_2}{\longrightarrow}\rho_2\otimes{\mathbf p}.$$

{\bf{Algebra element $\rho_3$.}} 
There are only two domains which can contribute $\rho_3$
on the type $D$ side, and these are $D_3 + D_4$ and $D_3+2 D_4$, which
we consider in turn.

{\it{$D_3+D_4$}.} 
This domain has two possible interpretations: either its input
consists of $({\mathbf p},\rho_2 \rho_3)$ or $({\mathbf p},\rho_3,
\rho_2)$. But the first interpretation is not valid: bilinearity over
the idempotents forces such a term to vanish. The second
interpretation leads to considering the domain $D_3 D_4$, thought of
as having a cut form ${\mathbf p}$ out to $\partial_R \HD$. As such,
it represents an annulus with a cut parameter at ${\mathbf r}$ (whose
other endpoints go out to the other boundary component). Such an 
annulus is always represented by a holomorphic curve. Hence we have:
$$({\mathbf p},\rho_3,\rho_2)\xrightarrow{D_3 + D_4}
\rho_3\otimes{\mathbf r}.$$

{\it{$D_3+2D_4$}.} There are three possible interpretations of
this domain: $({\mathbf p},\rho_3,\rho_{23})$, 
$({\mathbf p},\rho_{23},\rho_3)$, or $({\mathbf p},\rho_3,\rho_2,\rho_3)$.
The second interpretation is impossible because of idempotents.

The third interpretation gives a moduli space whose expected dimension
is non-zero; i.e., the Maslov index of the moduli space is wrong.  This
is neatly formulated in terms of the gradings calculated in 
Lemma~\ref{lem:GradingsDehnTwist}. Specifically,
substituting gradings calculated from that lemma, we see that
\begin{align*}
\lambda^{-1}\cdot \gr({\mathbf p}\otimes \rho_3[1]\otimes\rho_{2}[1]\otimes \rho_{3}[1])&=
\lambda^2 \cdot \gr({\mathbf p})\cdot f^\phi(\gr(\rho_3)\cdot \gr(\rho_{2}) \cdot \gr(\rho_3)) \\
&= (1;-1,0,1)
\end{align*}
(Note that the notation $\rho_i[1]$ means the element $\rho_i\in\Alg$ with 
a shift in its grading.)
which is different from
\[
\gr(\rho_3\otimes {\mathbf q})
= \gr(\rho_3)\cdot \gr({\mathbf q}) \\
= (0;-1,0,1).
\]
The reader might be concerned that this calculation depends
on several auxiliary choices, such as the choice of refinement data
$\psi$ and the identification of the grading set given in
Lemma~\ref{lem:GradingsDehnTwist}. However, the conclusion
that 
\[\lambda^{-1}\gr({\mathbf p}\otimes \rho_3[1]\otimes\rho_{2}[1]\otimes \rho_{3}[1])=
\lambda\cdot \gr(\rho_3\otimes {\mathbf q}), \] 
which rules out the possibility that $\rho_{3}\otimes {\mathbf q}$
appears in $m({\mathbf p},\rho_3,\rho_2,\rho_3)$,
is independent of this
choice (and indeed the exponent of $\lambda$ appearing on the
right-hand-side here gives the dimension of the relevant moduli space.)

This leaves only the first possible interpretation
$$Q_1\colon ({\mathbf p},\rho_3,\rho_{23})\xrightarrow{D_3+2D_4}
{\rho_3}\otimes {\mathbf q}.$$
(The label $Q$ here signifies that we have not (yet) determined
that this contribution is indeed $1\pmod{2}$.) Consider the $\Ainf$ 
relation with inputs $({\mathbf p}, \rho_3,\rho_2,\rho_3)$ and output
${\rho_3}\otimes{\mathbf q}$. The composite of 
$({\mathbf p},\rho_3,\rho_2)\stackrel{D_3 + D_4}{\longrightarrow}
\rho_3\otimes{\mathbf r}$ with 
$({\mathbf r},\rho_3)\stackrel{D_4}{\longrightarrow} {\mathbf q}$
gives on non-trivial term in this relation; the only possible
alternate contribution is gotten from the arrow $Q_1$ under consideration.
Thus, we have verified the existence of 
$$({\mathbf p},\rho_3,\rho_{23})\xrightarrow{D_3+2D_4}
{\rho_3}\otimes {\mathbf q}.$$

{\bf{Algebra element $\rho_{12}$.}} The possible domain is $D_1+D_2$,
thought of as a map
$$Q_2\co ({\mathbf r},\rho_1)\xrightarrow{D_1+D_2}\rho_{12} \cdot {\mathbf q}.$$
But this is incompatible with the idempotents of ${\mathbf q}$.

{\bf{Algebra element $\rho_{23}$.}} The domains are $D_2+D_3$,
$D_2+D_3+D_4$, and $D_2+D_3+2 D_4$. 

{\it{$D_2+D_3$}.} This domain must have a cut out to $\partial_L \HD$.
After this cut is made, the domain is a rectangle, giving a contribution
$$({\mathbf q}, \rho_{2})\xrightarrow{D_2+D_3}\rho_{23}\otimes {\mathbf r}.$$

{\it{$D_2+D_3+D_4$}.} This is a periodic domain, so could be
interpreted to have initial
and terminal generator either ${\mathbf p}$, ${\mathbf q}$, or
${\mathbf r}$. However, $\rho_{23}\otimes {\mathbf p}$ cannot appear
as the target of a type \DA\  action on ${\mathbf p}$: the left
idempotent of ${\mathbf p}$ is $\iota_0$, while the left idempotent of
$\rho_{23}$ is $\iota_1$.  Thus, we need consider only cases where the
initial and terminal point are ${\mathbf q}$ or ${\mathbf r}$. We
exclude first the latter case. Idempotents ensure that the only
possible interpretation of $D_2+D_3+D_4$ as a non-trivial contribution
to the type \DA\  module with initial generator ${\mathbf r}$ is to
think of it as a domain from $({\mathbf r},\rho_3,\rho_2)$ to
$\rho_{23}\otimes {\mathbf r}$.  But this is ruled out by gradings:
\begin{align*}
  \lambda^{-1}\cdot \gr({\mathbf r}\otimes \rho_3[1]\otimes \rho_2[1]) &=
  \lambda\cdot\gr({\mathbf r})\cdot f^\phi(\gr(\rho_3)\cdot \gr(\rho_2)) \\
  &= (1;1,2,1)),
\end{align*}
which is different from
\begin{align*}
  \gr(\rho_{23})\cdot \gr({\mathbf r}) 
  &= (0;1,2,1).
\end{align*}
Thus, the initial and terminal generator must be ${\mathbf q}$. Examining
idempotents once again, we see that there are only two possible interpretations
of this domain. One is as a domain from $({\mathbf q},\rho_2,\rho_3)$
to $\rho_{23}\otimes {\bf q}$. But this is ruled out by looking
at gradings. The only remaining interpretation 
of the domain is as a map
$$Q_3\colon ({\mathbf q},
\rho_{23})\xrightarrow{D_2+D_3+D_4}\rho_{23}\otimes
{\mathbf q}.$$
We will now establish the existence of the contribution by $Q_3$,
using the $\Ainf$ relation, together with the information about
the type \DA\  bimodule which we have so far collected. Specifically,
consider the $\rho_{23}\otimes {\mathbf q}$-coefficient of the
$\Ainf$ relation with inputs $({\mathbf q},\rho_2,\rho_3)$. We know
that $D_2+D_3$ and $D_{4}$ give a non-zero contribution of 
$m(m({\mathbf q},\rho_2),\rho_3)$. Other terms in this $\Ainf$ relation
include $m({\mathbf q},\rho_{23})$, which is counted by $Q_3$.
Possible other terms are: $m(m({\mathbf q}),\rho_2,\rho_3)$ and
$m(m({\mathbf q},\rho_2,\rho_3))$.  Looking back  at those terms which contribute
algebra elements $1$, $\rho_2$, $\rho_3$, and $\rho_{23}$, we see there are no
possible such terms.
This forces the existence of
$$({\mathbf q}, \rho_{23})\xrightarrow{D_2+D_3+D_4}\rho_{23}\otimes {\mathbf q}.$$

{\it{$D_2+D_3+2 D_4$.}} This domain starts at ${\mathbf r}$ and terminates
at ${\mathbf q}$. There are two possible interpretations of the domain
to $\rho_{23}\otimes {\mathbf q}$: one starts at 
$({\mathbf r}, \rho_3, \rho_{23})$
while the other starts at 
$({\mathbf r}, \rho_3, \rho_{2},\rho_{3})$.
Both possibilities are excluded by considering gradings:
\begin{align*}
\lambda^{-1}\gr({\mathbf r}\otimes\rho_3[1]\otimes\rho_{23}[1])
&=\lambda\cdot \gr(\rho_{23}\otimes {\mathbf q}) \\
\lambda^{-1}\gr({\mathbf r}\otimes\rho_3[1]\otimes\rho_{2}[1]\otimes\rho_{3}[1])
&=\lambda^2\cdot \gr(\rho_{23}\otimes {\mathbf q}).
\end{align*} 

{\bf{Algebra element $\rho_{123}$.}} 
The domains now are $D_1+D_2+D_3$ and $D_1+D_2+D_3+D_4$. 

{\it{$D_1+D_2+D_3+D_4$}; first visit.}  There are two conceivable
interpretations of this domain compatible with the idempotents: one
with input $({\mathbf p},\rho_{123})$, and the other with input
$({\mathbf p},\rho_3,\rho_2,\rho_1)$. The first is forced to exist by
the $\rho_{123}\otimes{\mathbf q}$-coefficient of the $\Ainf$ relation
with inputs $({\mathbf p}, \rho_{1},\rho_{23})$, since we have already
verified that $m_{2}(m_2({\mathbf
  p},\rho_1),\rho_{23})=\rho_{123}\otimes {\mathbf q}$.  We will
return to the second interpretation of $D_1+D_2+D_3+D_4$; but first we turn to the easier analysis of 
$D_1+D_2+D_3$.

{\it{$D_1+D_2+D_3$}.} The initial point is ${\mathbf p}$ and the 
terminal point ${\mathbf r}$, and hence it follows readily that the
only interpretation of this domain is as a map
$$Q_3\co m({\mathbf p},\rho_{12})\xrightarrow{D_1+D_2+D_3}
\rho_{123}\otimes {\mathbf r}.$$
We verify the existence of this map, by considering the
$\rho_{123}\otimes {\mathbf r}$-component of the 
$\Ainf$ relation with inputs
$({\mathbf p}, \rho_{12}, \rho_3)$. One contribution to this
is furnished by the domain $D_1+D_2+D_3+D_4$, interpretated as giving the operation
$m_2({\mathbf p},\rho_{123})=\rho_{123}\otimes {\mathbf r}$. 
The only conceivable alternate contribution is provided
by the composite of the presently considered domain ($Q_3$) with 
$({\mathbf r},\rho_3)\stackrel{D_4}{\longrightarrow}\rho_1\otimes{\mathbf q}$,
hence forcing the existence of 
$$m({\mathbf p},\rho_{12})\xrightarrow{D_1+D_2+D_3}
\rho_{123}\otimes {\mathbf r}.$$

{\it{$D_1+D_2+D_3+D_4$} revisited.} Recall that this domain had two
interpretations; one with one with input $({\mathbf p},\rho_{123})$,
and the other with input $({\mathbf p},\rho_3,\rho_2,\rho_1)$. We
already verified the existence of the curve with the first
interpretation.

We will see that $m_4({\mathbf p},\rho_3,\rho_2,\rho_1)$ vanishes (i.e., the
contribution of $D_1+D_2+D_3+D_4$ under this second interpretation is
zero), but this is surprisingly
subtle. The $\Ainfty$ relations on the type DA structure are not rich
enough to give information in this case: the inputs cannot be
factored, and the product of any non-idempotent element with
coefficient $\rho_{123}$ vanishes.
Moreover, there are two types of curves which can contribute to this
$m_4$: in one, there are two cuts going out to $\partial_L\HD$, and in 
the other, there are no such cuts. 
To calculate the $m_4$, we use some information which can be
extracted from the type \AAm\  bimodule, which has the advantage that the two 
kinds of curves are counted differently. To this
end, we label chords on the left of the diagram in a type $A$ manner,
i.e., placing algebra elements $\sigma_1$, $\sigma_2$, and $\sigma_3$
along $\partial_L \HD$ in the regions $D_3$, $D_2$, and $D_1$ respectively.

Interpret the domain $D_1+D_2+D_3+D_4$ as contributing
in the type \AAm\  bimodule. Cutting all the way out to the boundary at
each corner, we obtain a rectangle. Traversing its boundary, we 
find a contribution
$$X\co ({\mathbf p},\sigma_3,\sigma_2,\sigma_1,\rho_3,\rho_2,\rho_1)
\xrightarrow{D_1+D_2+D_3+D_4}{\mathbf q}.$$
We wish to argue that there is also a curve contributing
$$Y\co ({\mathbf p},\sigma_{123},\rho_3,\rho_2,\rho_1)\xrightarrow{D_1+D_2+D_3+D_4}
{\mathbf q}.$$
To see this, we consider the $\Ainf$ relation with inputs
$({\mathbf p}, \rho_3,\rho_2,\rho_1,\sigma_1,\sigma_{23})$
(recall that the positions of the $\sigma_i$ relative to the $\rho_i$ 
have no meaning). One contribution to this $\Ainf$ relation is the juxtaposition
$$({\mathbf p},\rho_{3},\rho_{2},\sigma_1)\xrightarrow{D_3+D_4}{\mathbf r}$$
(an annulus, which is easily seen to have a representative)
with
$$({\mathbf r},\rho_1,\sigma_{23})\xrightarrow{D_1+D_2}{\mathbf q}$$
(a rectangle).
The only term which can cancel this juxtaposition is 
$m_4({\mathbf p},\sigma_{123},\rho_3,\rho_2,\rho_1)$, forcing the existence
of the curve $Y$. 

Turning attention back to the type \DA\ bimodule, observe that both
the curves $X$ and $Y$ contribute to the $\rho_{123}\otimes {\mathbf
  q}$-coefficient of $m_4({\mathbf r},\rho_3,\rho_2,\rho_1)$. Thus,
taken together we see that the contribution of $D_1+D_2+D_3+D_4$ to
$m_4({\mathbf p},\rho_3,\rho_2,\rho_1)$ vanishes. (Recall that in the
calculation of the type \DD\ bimodule for the identity a similar
mechanism applied: the module was deduced from the \AAm\ identity bimodule,
and there were two distinct $\Ainfty$ operations in the \AAm\
bimodules from $\mathbf w_1$ to $\mathbf z_2$
which cancelled when interpreted as contributions in the type \DD\ identity bimodule.)

This completes the verification that $\CFDAa(\tau_\merid)$ has the form 
stated in Figure~\ref{fig:DAbimodules}. The other mapping class group
generators can be calculated in an entirely parallel way.
\end{proof}

\subsection{Duality for the genus one handlebody}
\label{subsec:ExampleDuality}

We illustrate now the duality theorem, Theorem~\ref{thm:Duality},
by explicitly verifying that the
description of the type $A$ module for the genus one handlebody gotten
by inspecting a Heegaard diagram with one generator is
quasi-isomorphic to the one gotten from the duality theorem.
This calculation will use the type \DD\ identity 
bimodule in the genus one
case.

Continuing notation from Section~\ref{subsec:AAId1}, we let
$\Alg=\Alg(\PMC,0)$ and $\Blg=\Alg(-\PMC,0)$. 
We abbreviate $\lsup{\Alg,\Blg}D$ for $\lsup{\Alg,\Blg}\CFDDa(\Id)$.

Our goal is to calculate the $\Blg$-module
$\Mor^{\Alg}(\lsub{\Blg}\Blg_{\Blg}\DT\lsup{\Alg,\Blg}D,\lsup{\Alg}N)$,
where
here $\lsup{\Alg}N$ is the type~$D$ module associated to a $0$-framed
solid torus from~\cite[Section~\ref*{LOT:sec:surg-exact-triangle}]{LOT1}. 
We will use the interpretation of this morphism
space provided by Corollary~\ref{cor:InterpretMor}.
$\lsup{\Alg}N$ has a single generator as a type $D$
module; and hence
$\lsub{\Alg}N=\lsub{\Alg}\Alg_{\Alg}\DT\lsup{\Alg}N$ is spanned by
three elements: $\x$, $\rho_2\otimes \x$ and $\rho_{12}\otimes
\x$. Note that $\rho_2\otimes \x$ is in the left $\iota_1$-idempotent,
while the two other elements are in the left $\iota_0$ idempotent. 

The differential on $\lsup{\Alg}N$ is given by
$\delta^1(\x)=\rho_{12}\otimes \x$, so the differential on
$\lsub{\Alg}N$ is given by $\bdy (\x)=\rho_{12}\otimes \x$ and
$\bdy(\rho_{2}\otimes \x)=\bdy(\rho_{12}\otimes \x)=0$.

\begin{figure}
\input{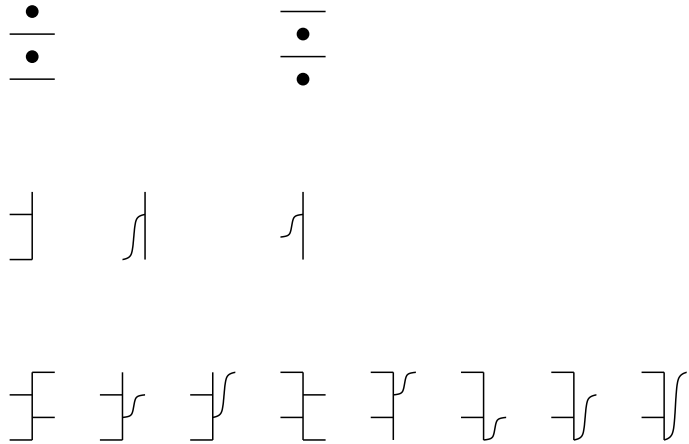}
\caption{\label{fig:ModuleGenerators} {\bf{Module generators.}} In the
  top row, we have displayed the two primitive idempotents ($\iota_0$
  and $\iota_1$ respectively) in the torus algebra. In the second row,
  we exhibit the three generators ($\x$, $\rho_{12}\otimes \x$, and
  $\rho_2\otimes \x$ respectively) of the module for a framed
  handlebody. In the third row, we have displayed the eight generators
  of $\lsubsup{\Blg}{\Alg}D$. They are
  $\iota_0\otimes j_0$, $\iota_0\otimes \sigma_2$, $\iota_0\otimes\sigma_{12}$,
  $\iota_1\otimes j_1$, $\iota_1\otimes\sigma_1$, $\iota_1\otimes
  \sigma_3$, $\iota_1\otimes\sigma_{23}$, and $\iota_1\otimes\sigma_{123}$
  respectively.}
\end{figure}

The module $\lsubsup{\Blg}{\Alg}{D}=
\lsub{\Blg}\Blg_{\Blg}\DT\lsup{\Alg\otimes\Blg}D$ has eight
generators, which are naturally labelled
$\iota_0\otimes j_0$, $\iota_0\otimes \sigma_2$, $\iota_0\otimes\sigma_{12}$, $\iota_1\otimes j_1$, $\iota_1\otimes\sigma_1$,
$\iota_1\otimes \sigma_3$, $\iota_1\otimes\sigma_{23}$, and
$\iota_1\otimes\sigma_{123}$. 
These are all enumerated in
Figure~\ref{fig:ModuleGenerators}. A homomorphism from $\lsubsup{\Blg}{\Alg}D$ to
$\lsup{\Alg}N$ is uniquely specified by where it takes each of those eight
generators. Of those eight generators, three are in the
$\iota_0$-idempotent, whereas five are in the $\iota_1$-idempotent.
Thus, the vector space 
$\Mor^{\Alg}(\lsubsup{\Blg}{\Alg}D,\lsup{\Alg}N)$
has $11$ basis
vectors, which are gotten by the six possible maps sending any of the
three elements $\{\iota_0\otimes j_0$, $\iota_0\otimes \sigma_2$
$\iota_0\otimes \sigma_{12}\}$ to any of the two elements
$\{\x,\rho_{12}\otimes \x\}$, or the five maps gotten by sending
any of the five elements
$\{\iota_1\otimes j_1, \iota_1\otimes\sigma_1, \iota_1\otimes \sigma_3,
\iota_1\otimes\sigma_{23}, \iota_1\otimes\sigma_{123}\}$
to
$\rho_2\otimes \x$.

These basis vectors are relabelled as indicated in
Figure~\ref{fig:ModuleDifferentials}. The relabeling has the convenient property
that for all $i=1,\dots,5$, the differential of $T_i$ contains a non-trivial
component in $H_i$; indeed, these are all the components of the differential coming from
$\lsubsup{\Blg}{\Alg}D$. There are three additional components in the differential,
taking $X$ to $H_3$, $T_1$ to $H_4$, and $T_2$ to $H_5$. (These are all
components coming from the differential in $\lsup{\Alg}N$).

\begin{figure}
\input{ModuleDifferentials}
\caption{\label{fig:ModuleDifferentials} {\bf{Differentials in $\Mor^{\Alg}(\lsubsup{\Blg}{\Alg}D,\lsup{\Alg}N)$.}}
  We have two charts, indicating the two types of elements (corresponding
  to the two left-idempotent), wherein we relabel the eleven basis vectors
  $\Mor^{\Alg}(\lsubsup{\Blg}{\Alg}D,\lsup{\Alg}N)$
  by $X$, $\{T_i\}_{i=1}^5$ and $\{H_i\}_{i=1}^S$.
  For instance, $T_1$ is the element which carries $\iota_0\otimes \sigma_2$
  to the element $\x$, while $H_4$ takes $\iota_0\otimes \sigma_2$
  to $\rho_{12}\otimes \x$. There are three additional components to the differential
  not captured by the convention that $\partial T_i = H_i + \cdots$: an $H_3$ term in $\partial X$,
an $H_4$ term in $\partial T_1$, and an $H_5$ term
in $\partial T_2$.}
\end{figure}

The homology of this complex is carried by $X+T_3$. 
It is
straightforward to verify that
\begin{align*}
  (X+T_3)\cdot \sigma_3 &= H_1 \\
  \partial (T_1+T_4) &= H_1 \\
  (T_1+T_4)\cdot \sigma_2  &= X+T_3 \\
  (T_1+T_4)\cdot \sigma_{23} &= H_1.
\end{align*}
Setting $\y=X+T_3$, we see readily that this module is $\Ainf$-homotopy
equivalent to the module generated by $\y$, with $\Ainf$ relations
indexed by integers $i\geq 0$ 
\begin{align*}
  m_{3+i}(\y,\sigma_3,\overbrace{\sigma_{23},\dots,\sigma_{23}}^{i}, \sigma_2)&=\y
\end{align*}
which is in fact the module obtained from the inspection of holomorphic disks
(see \cite[Lemma~\ref*{LOT:lemma:SolidTorus}]{LOT1}).

%%% Local Variables: 
%%% mode: latex
%%% TeX-master: "Bimodules"
%%% End: 

\bibliographystyle{hamsalpha}\bibliography{heegaardfloer}
\end{document}